\documentclass[12pt]{amsbook}

\usepackage{amssymb}
\usepackage{amsthm}
\usepackage{amscd}
\usepackage{pifont}
\usepackage{amsmath}
\usepackage{latexsym}
\usepackage{graphicx}
\usepackage{amsfonts}
\usepackage{cancel}
\usepackage{mathrsfs}
\usepackage[all]{xy}

\newtheorem{lemma}{Lemma}[chapter]
\newtheorem{prop}[lemma]{Proposition}
\newtheorem{thm}[lemma]{Theorem}
\newtheorem{cor}[lemma]{Corollary}
\newtheorem{aplemma}{Lemma~A.\hspace{-1.5mm}}
\newtheorem{approp}{Proposition~A.\hspace{-1.5mm}}
\newtheorem{apthm}{Theorem~A.\hspace{-1.5mm}}
\newtheorem{apcor}{Corollary~A.\hspace{-1.5mm}}
\newtheorem{intthm}{Theorem}

\newtheorem{prdef}[lemma]{Proposition-Definition}

\newtheorem*{notation}{Notation and Conventions}

\theoremstyle{definition}

\newtheorem{rem}{Remark}

\newtheorem{rema}[lemma]{Remark}

\newtheorem{remb}{Remark}

\newtheorem{defi}[lemma]{Definition}
\newtheorem{exa}[lemma]{Example}
\newtheorem{aprem}{Remark~A.\hspace{-1.5mm}}
\newtheorem{apdefi}{Definition~A.\hspace{-1.5mm}}
\newcommand{\bde}{\begin{defi}}
\newcommand{\ede}{\end{defi}\vspace{0mm}}
\newcommand{\ble}{\begin{lemma}}
\newcommand{\ele}{\end{lemma}}
\newcommand{\bpr}{\begin{prop}}
\newcommand{\epr}{\end{prop}}
\newcommand{\bt}{\begin{thm}}
\newcommand{\et}{\end{thm}}
\newcommand{\bco}{\begin{cor}}
\newcommand{\eco}{\end{cor}}
\newcommand{\bre}{\begin{rem}}
\newcommand{\ere}{\end{rem}}
\newcommand{\brea}{\begin{rema}}
\newcommand{\erea}{\end{rema}\vspace{1mm}}
\newcommand{\breb}{\begin{remb}}
\newcommand{\ereb}{\end{remb}\vspace{1mm}}
\newcommand{\bex}{\begin{exa}}
\newcommand{\eex}{\end{exa}}
\newcommand{\bpf}{\begin{proof}}
\newcommand{\epf}{\end{proof}\vspace{1mm}}

\newcommand{\bade}{\begin{apdefi}}
\newcommand{\eade}{\end{apdefi}}
\newcommand{\bale}{\begin{aplemma}}
\newcommand{\eale}{\end{aplemma}}
\newcommand{\bapr}{\begin{approp}}
\newcommand{\eapr}{\end{approp}}
\newcommand{\bat}{\begin{apthm}}
\newcommand{\eat}{\end{apthm}}
\newcommand{\baco}{\begin{apcor}}
\newcommand{\eaco}{\end{apcor}}
\newcommand{\bare}{\begin{aprem}}
\newcommand{\eare}{\end{aprem}}

\newcommand{\be}{\begin{enumerate}}
\newcommand{\ee}{\end{enumerate}}
\newcommand{\bcd}{\[\begin{CD}}
\newcommand{\ecd}{\end{CD}\]}
\newcommand{\bit}{\begin{itemize}}
\newcommand{\eit}{\end{itemize}}
\newcommand{\bq}{\begin{quote}}
\newcommand{\eq}{\end{quote}}
\newcommand{\ba}{\begin{array}}
\newcommand{\ea}{\end{array}}
\newcommand{\mcA}{\mathcal{A}}
\newcommand{\mcB}{\mathcal{B}}

\newcommand{\mcD}{\mathcal{D}}
\newcommand{\mcE}{\mathcal{E}}
\newcommand{\mcF}{\mathcal{F}}
\newcommand{\mcG}{\mathcal{G}}
\newcommand{\mcH}{\mathcal{H}}

\newcommand{\mcJ}{\mathcal{J}}
\newcommand{\mcK}{\mathcal{K}}
\newcommand{\mcL}{\mathcal{L}}

\newcommand{\mcN}{\mathcal{N}}
\newcommand{\mcO}{\mathcal{O}}
\newcommand{\mcP}{\mathcal{P}}

\newcommand{\mcR}{\mathcal{R}}
\newcommand{\mcS}{\mathcal{S}}
\newcommand{\mcT}{\mathcal{T}}
\newcommand{\mcU}{\mathcal{U}}
\newcommand{\mcV}{\mathcal{V}}
\newcommand{\mcW}{\mathcal{W}}

\newcommand{\mbA}{\mathbb{A}}
\newcommand{\mbB}{\mathbb{B}}
\newcommand{\mbC}{\mathbb{C}}

\newcommand{\mbF}{\mathbb{F}}
\newcommand{\mbG}{\mathbb{G}}
\newcommand{\mbH}{\mathbb{H}}

\newcommand{\mbJ}{\mathbb{J}}

\newcommand{\mbN}{\mathbb{N}}

\newcommand{\mbP}{\mathbb{P}}
\newcommand{\mbQ}{\mathbb{Q}}
\newcommand{\mbR}{\mathbb{R}}
\newcommand{\mbS}{\mathbb{S}}
\newcommand{\mbT}{\mathbb{T}}

\newcommand{\mbV}{\mathbb{V}}
\newcommand{\mbW}{\mathbb{W}}

\newcommand{\mbY}{\mathbb{Y}}
\newcommand{\mbZ}{\mathbb{Z}}

\newcommand{\mfA}{\mathfrak{A}}

\newcommand{\mfC}{\mathfrak{C}}
\newcommand{\mfD}{\mathfrak{D}}
\newcommand{\mfE}{\mathfrak{E}}

\newcommand{\mfM}{\mathfrak{M}}
\newcommand{\mfN}{\mathfrak{N}}
\newcommand{\mfO}{\mathfrak{O}}

\newcommand{\mfS}{\mathfrak{S}}

\newcommand{\mfU}{\mathfrak{U}}
\newcommand{\mfV}{\mathfrak{V}}

\newcommand{\mfb}{\mathfrak{b}}
\newcommand{\mfc}{\mathfrak{c}}

\newcommand{\mfe}{\mathfrak{e}}
\newcommand{\mff}{\mathfrak{f}}
\newcommand{\mfg}{\mathfrak{g}}
\newcommand{\mfh}{\mathfrak{h}}
\newcommand{\mfi}{\mathfrak{i}}

\newcommand{\mfl}{\mathfrak{l}}
\newcommand{\mfm}{\mathfrak{m}}
\newcommand{\mfn}{\mathfrak{n}}
\newcommand{\mfo}{\mathfrak{o}}
\newcommand{\mfp}{\mathfrak{p}}
\newcommand{\mfq}{\mathfrak{q}}
\newcommand{\mfr}{\mathfrak{r}}
\newcommand{\mfs}{\mathfrak{s}}
\newcommand{\mft}{\mathfrak{t}}

\newcommand{\msA}{\mathscr{A}}

\newcommand{\msC}{\mathscr{C}}

\newcommand{\msE}{\mathscr{E}}
\newcommand{\msF}{\mathscr{F}}
\newcommand{\msG}{\mathscr{G}}

\newcommand{\msL}{\mathscr{L}}

\newcommand{\msN}{\mathscr{N}}

\newcommand{\msP}{\mathscr{P}}

\newcommand{\msU}{\mathscr{U}}
\newcommand{\msV}{\mathscr{V}}

\newcommand{\msX}{\mathscr{X}}

\newcommand{\migi}{\rightarrow}
\newcommand{\longmigi}{\longrightarrow}

\newcommand{\isom}{\stackrel{\sim}{\migi}}

\newcommand{\migiincl}{\hookrightarrow}

\newcommand{\migisurj}{\twoheadrightarrow}

\newcommand{\triv}{\mcO}

\newcommand{\DD}{\mcD \text{{\it iff}}}
\newcommand{\JJ}{{^{\rotatebox{90}{{\tiny $\Join$}}}}}
\newcommand{\DDD}{\mcD\text{{\it iff}}}

\newcommand{\mr}{\mathrm}
\newcommand{\hidden}[1]{\,}

\newcommand{\ZZZ}{{^{\mr{Zzz...}}}}
\newcommand{\NILP}{{^{p\text{-}\mr{nilp}}}}
\newcommand{\EE}{\Box}

\newcommand{\BB}{\varTheta}
\newcommand{\bb}{\vartheta}

\newcommand{\A}{\mcF^\sharp}

\newcommand{\Spec}{\ooalign{$\hspace{0.5mm}\_$ \cr $\rotatebox[origin=c]{180}{$\mbV$}$}\hspace{-0.5mm}}

\newcommand{\DE}{{^\dagger}\mcE}
\newcommand{\Nab}{\nabla}

\newcommand{\DF}{{^\dagger}\mcF}
\newcommand{\Dif}{D\hspace{-3.0mm}{^{_{=}}}\hspace{1mm}}
\newcommand{\DV}{{^\dagger}\mcV}
\newcommand{\dV}{{^\dagger}\mcV}
\newcommand{\VV}{{\ooalign{$\vee$ \cr $\wedge$}}}
\newcommand{\TT}{\blacktriangle}
\newcommand{\Fus}{\rotatebox[origin=c]{180}{$\mbY$}}
\newcommand{\Om}{\Omega^\sharp}
\newcommand{\BL}{?}
\newcommand{\qq}{q_1\hspace{-5mm}\rotatebox[origin=c]{-10}{$\swarrow$}\hspace{0.5mm}}
\newcommand{\Dual}{\rotatebox[origin=c]{180}{$D$}\hspace{-0.5mm}}
\newcommand{\dual}{\blacktriangledown}
\newcommand{\GR}{{^{\rotatebox[origin=c]{-20}{$\vartriangle$}}}\hspace{-3mm}G}

\newcommand{\GRRR}{{^{\rotatebox[origin=c]{-20}{$\vartriangle$}}}\hspace{-2.2mm}G}
\newcommand{\flqq}{\fontencoding{T1}\fontfamily{ppl}\selectfont <<}
\newcommand{\frqq}{\fontencoding{T1}\fontfamily{ppl}\selectfont >>}

\makeindex

\begin{document}

\title[A theory of dormant opers on pointed stable curves]{A theory of dormant opers \\ on pointed stable curves
\\ ---a proof of Joshi's conjecture---}
\author{Yasuhiro Wakabayashi}
\date{}
\maketitle
\footnotetext{Y. Wakabayashi: Department of Mathematics, Tokyo Institute of Technology, 2-12-1 Ookayama, Meguro-ku, Tokyo 152-8551, JAPAN;}
\footnotetext{e-mail: {\tt wkbysh@math.titech.ac.jp};}
\begin{abstract}
This manuscript presents a detailed and original account of the theory of opers  defined on  pointed stable curves in arbitrary characteristic and their moduli. 
 In particular, it includes the development of the study of dormant opers, which are opers of a certain sort in 
 positive characteristic. 
 The theory of dormant opers (or more 
generally, opers in positive characteristic) on pointed stable curves, which has 
proved to be rather rich and deep, was born in the work of S. 
Mochizuki,  
 who developed the theory for $\mfs \mfl_2$-opers and used it to 
establish 
$p$-adic Teichm\"{u}ller theory. 
Some parts of Mochizuki's work were 
 later extended in the case of proper smooth curves by K. Joshi, C. Pauly, and other mathematicians. 
This manuscript represents an advance in the theory of opers 
 that  takes the subject beyond the work of Mochizuki, Joshi, and 
Pauly.
In particular, we provide general unified formulations and the 
basics of principal bundles and connections defined on families of  pointed stable curves. 
The notion of an oper is accordingly introduced 
 in the context of logarithmic algebraic geometry. 
 Some of the results 
 can be regarded as generalizations of results obtained in the 
 fundamental work on the geometric Langlands program developed by A. 
 Beilinson and V. Drinfeld. 
 We also describe various properties 
 and assertions about (dormant) opers,  such as  duality, comparison with differential operators, and compactification  of the moduli space. 
Our goal is to give an explicit  
 formula,  conjectured by Joshi, for the generic number of 
 dormant $\mathfrak{sl}_n$-opers.
  We do so by obtaining a detailed understanding of 
the 
 moduli space of dormant opers and computing the Gromov-Witten invariants for Quot-schemes in characteristic zero. This 
formula 
 reveals an interaction between studies in $p$-adic 
 Teichm\"{u}ller theory and certain areas of mathematics, including Gromov-Witten theory.
\end{abstract}
\setcounter{tocdepth}{1}
\setcounter{chapter}{-1}
\tableofcontents
\chapter{Introduction}




{\it \begin{flushright} 
\flqq Mon prince, il y a plus de cinquante ans 
que j'ay ou\"{i} dire  \\ 
\`{a} mon pere  qu'il y avoit dans ce chasteau une princesse, 
  \\ 
   la plus belle du monde; qu'elle y devoit dormir cent ans,   et qu'elle
 \\
   serait r\'{e}veill\'{e}e par le fils d'un roy, \`{a} qui elle estoit reserv\'{e}e.\frqq
   \end{flushright}}
 \begin{flushright}
 --- Charles Perrault, {\it La Belle au bois dormant}  (cf. ~\cite{Perra}).
 \end{flushright}


\hspace{3mm}

 This manuscript presents 
 the foundation 
 for a theory of opers  defined on  pointed stable curves in arbitrary characteristic and their moduli. 
 In particular, it includes the development of the study of dormant opers, which are opers of a certain sort in 
 positive characteristic. 
The final goal of this manuscript  is to prove 
the following conjectural formula: \vspace{5mm} \leavevmode\\
{\bf Joshi's conjecture} (cf. ~\cite{Jo14}, Conjecture 8.1)
\leavevmode\\
\vspace{0mm} 
 \begin{equation} \label{Joshi's conjecture}
 \mr{deg}_{}(\mfO \mfp^\ZZZ_{\mfs \mfl_n, X}/\overline{k})
 = \frac{ p^{(n-1)(g-1)-1}}{n!} \cdot  
 \sum_{\genfrac{.}{.}{0pt}{}{(\zeta_1, \cdots, \zeta_n) \in \mbC^{\times n} }{ \zeta_i^p=1, \ \zeta_i \neq \zeta_j (i\neq j)}}
 \frac{(\prod_{i=1}^n\zeta_i)^{(n-1)(g-1)}}{\prod_{i\neq j}(\zeta_i -\zeta_j)^{g-1}}.  \end{equation}
\vspace{3mm}

We shall explain what is displayed above. 
Let $\overline{k}$ be  an algebraically closed field  of characteristic $p> 0$ and
$X$  a  proper smooth curve of genus $g>1$ over  $\overline{k}$.
Assume  that $p> C$, where $C$ is  a specific constant  depending on both  $g$ and a fixed positive integer $n$.
 It is known  (cf. ~\cite{JP}, Theorem 5.4.1 and Corollary 6.1.6) that  the moduli functor $\mfO \mfp^\ZZZ_{\mfs \mfl_n, X}$ classifying {\it dormant $\mfs \mfl_n$-opers} on $X$ may be represented by  a nonempty finite scheme over $\overline{k}$.
Toward a further understanding of this moduli space,
 Kirti Joshi  proposed, in  ~\cite{Jo14},
the conjectural formula (\ref{Joshi's conjecture}) displayed above, 
computing the degree $\mr{deg}_{}(\mfO \mfp^\ZZZ_{\mfs \mfl_n, X}/\overline{k})$ of $\mfO \mfp^\ZZZ_{\mfs \mfl_n, X}$ over $\overline{k}$.
Here, the sum on the right-hand side of  (\ref{Joshi's conjecture}) is taken over the set of $n$-tuples $(\zeta_1, \cdots, \zeta_n) \in \mbC^{\times n}$ of $p$-th roots of unity in $\mbC$ (= the field of complex numbers)  that satisfy   $\zeta_i \neq \zeta_j$ if $i \neq j$.

\section{Review of opers on complex algebraic curves} \label{y01}

Before continuing our explanation of the present manuscript, let us briefly review what an {\it oper} is. 
Recall (cf. Definition \ref{y035}) that 
an oper, or  more precisely a $\mfg$-oper  for a semisimple Lie  algebra $\mfg$, 
 is defined as a principal bundle
 over an algebraic  curve equipped with an integrable connection satisfying certain conditions.
Opers on an algebraic  curve, or a formal disc,  over  $\mbC$ play a central role in some areas of mathematics.
For example, such objects  were studied   in the context of the geometric Langlands program (cf. ~\cite{BD1})   to construct Hecke eigensheaves on the moduli space of principal bundles by associating them with   the base space 
of Hitchin's integrable system.
Also, 
as discussed in  ~\cite{DS},  
 the moduli space of opers on a formal disc is realized as
the phase space of Drinfeld-Sokolov 
 hierarchies, i.e.,  certain infinite   families of nonlinear PDEs; 
 these hierarchies are certainly among the most studied examples of integrable systems because of their tau-symmetries.

Many mathematicians have studied various equivalent 
realizations of opers,
 including  certain kinds of differential operators (related to Schwarzian equations) between line bundles.
For example, 
$\mfs \mfl_2$-opers  may be interpreted as a   {\it projective structure}. 
A projective structure on a Riemann surface $X_\mbC$  is a maximal system of local coordinates on $X_\mbC$  in the analytic topology modeled on the complex projective line  $\mbP_\mbC^1 \left(= \mbC \sqcup \{ \infty \} \right)$ such  that 
on any two overlapping coordinate patches, the change of coordinates may be described   as a M\"{o}bius transformation.
On a Riemann surface  equipped with a projective  structure, there is  a {\it    projective geometry},  in the spirit of F. Klein's Erlangen program, that locally agrees with the geometry of $\mbP^1_\mbC$.
Projective  structures  
  arise in many areas of mathematics and are still actively being investigated.
Indeed,  such additional structures on Riemann surfaces 
play a major role in understanding the framework of uniformization theorem.
An important consequence of the uniformization theorem for Riemann surfaces is that any closed Riemann surface 
 is isomorphic either to   $\mbP^1_\mbC$, or to the quotient of the complex affine line $\mbA_\mbC^1  \left(=\mbC \subseteq \mbP_\mbC^1 \right)$ by a discrete group of translations (i.e.,  a  $1$-dimensional complex torus), or to the quotient of the $1$-dimensional complex hyperbolic space $\mbH^1_\mbC \ \left(\subseteq \mbP_\mbC^1 \right)$ by a torsion-free discrete subgroup of $\mr{SU}(1, 1)\cong \mr{SL}(2, \mbR)$.
This implies that by collecting various local inverses of a universal covering map,
   we have a canonical projective structure on each  closed Riemann surface.

Moreover, $\mfs \mfl_2$-opers were introduced and studied, under the name of {\it indigenous bundles}, in the work of R. C. Gunning (cf. ~\cite{G2}, \S\,2).
One may think of an indigenous bundle as an {\it algebraic} object encoding {\it analytic} (i.e., {\it non-algebraic})  uniformization data for  Riemann surfaces.
Also,  in the work of C. Teleman (cf.  ~\cite{Tele}), one may find 
the objects corresponding to $\mfs \mfl_n$-opers, under the name of {\it homographic structures}.
In  the case of characteristic zero, we refer the reader  to  ~\cite{Bis}, ~\cite{BeBi}, and ~\cite{Simp} for further studies of $\mfs \mfl_n$-opers, and to ~\cite{BF}, ~\cite{Fr3}, and ~\cite{Fr2} 
 for reviews and expositions  concerning $\mfg$-opers for  an arbitrary $\mfg$ over $\mbC$.

\section{What about opers in positive characteristic?} \label{y02}

In the present manuscript, we study $\mfg$-opers in {\it arbitrary characteristic}, i.e.,   ones including the case of {\it positive characteristic}.
Just as in the complex case mentioned above,
  one may define the notion of a $\mfg$-oper in characteristic $p>0$  because of the algebraic nature of its formulation.
For $\mfg = \mfs \mfl_2$, various properties of such objects and their  moduli space 
  were discussed in the context of the $p$-adic Teichm\"{u}ller  theory developed  by S. Mochizuki (cf. ~\cite{Mzk1}, ~\cite{Mzk2}).
From a different point of view, Y. Ihara developed,  in, e.g.,  ~\cite{Ihara1}, ~\cite{Ihara2},  a theory of Schwarzian  equations in an arithmetic context.
For the case where  $\mfg$ is more general (but the underlying curve is assumed to be  smooth over an algebraically closed field),
$\mfg$-opers in positive characteristic
    were studied   by K. Joshi, S. Ramanan, E. Z. Xia, J. K. Yu, C. Pauly, T. H. Chen, X. Zhu,  et al.
(cf. ~\cite{JRXY}, ~\cite{JP}, ~\cite{Jo14}).

One of the common key ingredients in the development of these works is 
the study of  the {\it $p$-curvature} of a connection.
Recall that the $p$-curvature of a connection is an invariant that measures   the obstruction to the  compatibility of $p$-power structures appearing  in certain associated spaces of infinitesimal (i.e., ``Lie'')  symmetries.
It is well-known that a connection on a vector bundle has vanishing $p$-curvature if and only if the corresponding homogeneous linear   differential equation  locally has a full set of  solutions  in the Zariski topology.
We say that a $\mfg$-oper is {\it dormant} (cf. Definition \ref{y089}) if the  $p$-curvature of a connection defining this oper vanishes  identically.

This class of opers enables us to 
 develop geometry in positive characteristic  from  the same viewpoint as the theory of projective structures on Riemann surfaces.
One thing to keep in mind is that due to its analytic formulation, the definition of  a projective  structure cannot be adopted in positive characteristic as it is.
Indeed, 
if  we had  defined it in a naive fashion,  it would be nothing  other than the  trivial example of a curve having such a structure.
Regarding this matter, Y. Hoshi (cf. ~\cite{Hos2}) introduced  {\it Frobenius-projective structures} (of finite level) to deal with  appropriate replacements.
It was shown (cf. ~\cite{Hos2}, Theorem A) that  
giving a Frobenius-projective structure of level $1$ is equivalent to
giving 
a dormant $\mfs \mfl_2$-oper;  this equivalence  is an analogue  of the relationship, mentioned before,  between projective structures on a Riemann surface and  $\mfs \mfl_2$-opers  on it. 
We refer the reader  to  ~\cite{Wak5}, ~\cite{Wak6} for  studies related to this topic.

There are  some other aspects of dormant opers.
For example, 
  dormant $\mfs \mfl_2$-opers  on a proper smooth curve correspond, via Cartier descent,  to a certain type of Frobenius-destabilized vector bundles of rank $2$ (cf. ~\cite{O1}, \S\,4, Proposition 4.2).  
This correspondence gives us an approach to understanding the generalized  Verschiebung map on the moduli space of rank $2$ semistable bundles.
Indeed, in the case of genus $2$ curves, we can compute (cf. ~\cite{O3}, \S\,1, Theorem 1.2) the degree of this rational map  by using the fact that the set of points classifying such bundles in this moduli space coincides with the undefined locus.
Also,  results  by  S. Mochizuki, F. Liu, and B. Osserman (cf. ~\cite{Mzk2}, ~\cite{LO}, ~\cite{Wak2}) 
provide 
an interesting connection 
 with the  combinatorics of convex  rational polytopes, as well as edge-colored graphs,   and  the geometry of moduli spaces classifying  dormant $\mfs \mfl_2$-opers.

As another example, 
the sheaf of locally exact $1$-forms on $X$ may be
 obtained (cf. ~\cite{Ray1}, Remark 4.1.2)  as the sheaf of horizontal sections of a certain dormant $\mfs \mfl_{(p-1)}$-oper twisted by a flat line bundle, 
 i.e., a certain $\mr{GL}_{(p-1)}$-oper (cf. Definition \ref{y0108}, (i)) with vanishing $p$-curvature. 
This specific object 
is of importance because
 the theta divisor associated to this bundle (cf. ~\cite{Ray1}, Theorem 4.1.1) gives information on the fundamental group of the underlying curve.
It leads us to 
the sort of phenomena studied in 
anabelian geometry for hyperbolic curves in  positive characteristic (cf. ~\cite{PS}, ~\cite{Ray2},  ~\cite{Tam1}, ~\cite{Tam2}). 
(But, we should note that the phenomena studied in the cited references are mostly related to the associated theta divisor  and not based on properties of the  oper structure on it.)

On the other hand,  in the context of   $p$-adic Teichm\"{u}ller  theory,  it is also  worth studying $\mfg$-opers with  nilpotent $p$-curvature, which will be called  {\it $p$-nilpotent $\mfg$-opers} (cf. Definition \ref{y089}).
I. I. Bouw and S. Wewers set up (cf. ~\cite{BW}) an equivalence  between 
$p$-nilpotent $\mfs \mfl_2$-opers
and certain  deformation data on the underlying curve.
 This equivalence allows us to  translate the problem of the existence  of such deformation data into one on  the existence of polynomial solutions of certain differential equations with additional properties.
 In  ~\cite{Hoshi}, Y. Hoshi  advanced the understanding of  indigenous bundles with nilpotent $p$-curvature   in {\it characteristic three} by creating a complete list of them (cf. ~\cite{Hoshi}, Theorem 6.1) for the case where the underlying curve is of genus $2$.
However, in the present manuscript, we will focus more on dormant $\mfg$-opers rather than  $p$-nilpotent ones.

\section{Counting problem of dormant opers} \label{yff03}

With the above-mentioned background   in mind, 
the study of dormant $\mfg$-opers  and their moduli could be placed within 
developments of certain  areas of mathematics, and hence, 
the following question regarding such objects may arise naturally from various points of view:
\begin{quote}\textit{Can one explicitly calculate how many 
   dormant $\mfg$-opers exist  on a fixed (general) curve?} \end{quote}

In what follows, we shall review 
previous results 
concerning 
  the answer to the  above question.
Let  $X/\overline{k}$ be as introduced at the beginning of this Introduction.
A theorem   in  $p$-adic Teichm\"{u}ller theory due to S. Mochizuki implies (cf. ~\cite{Mzk2}, Chap. II, \S\,2.3, Theorem 2.8)  that 
 if   $X$ is    sufficiently general  in the moduli stack $\mfM_{g}$ of proper smooth curves of genus $g$, then  $\mfO \mfp_{\mfs \mfl_2, X}^\ZZZ$ is  finite (as we mentioned above) and {\it \'{e}tale}  over $\overline{k}$. 
Hence,  the task of resolving the above question for the case where $\mfg = \mfs \mfl_2$ may be reduced to the explicit computation of  the degree $\mr{deg}(\mfO \mfp_{\mfs \mfl_2, X}^\ZZZ/\overline{k})$  of $\mfO \mfp_{\mfs \mfl_2, X}^\ZZZ$ over $\overline{k}$, that is, to proving  Joshi's conjecture for   $\mfg = \mfs \mfl_2$.

In the case of $g=2$,
S. Mochizuki (cf. ~\cite{Mzk2}, Chap.\,V, \S\,3.2,  Corollary 3.7), H. Lange-C. Pauly (cf. ~\cite{LP}, Theorem 2), and B. Osserman (cf. ~\cite{O4}, Theorem 1.2)
verified, by applying different methods,  the following  equality:
\begin{equation} \label{n=2}
  \mr{deg}(\mfO \mfp_{\mfs \mfl_2, X}^\ZZZ/\overline{k}) = \frac{1}{24}\cdot (p^3-p). 
  \end{equation}
Moreover, by extending the relevant formulations to the case where $X$ admits marked points and  nodal singularities,
  S. Mochizuki also gave
  an explicit description of dormant $\mfs \mfl_2$-opers on each {\it totally degenerate curve} (cf. Definition \ref{y0178}) in terms of  the {\it radii of atoms}    (cf. ~\cite{Mzk2}, Introduction, Theorem 1.3).
 This is essentially equivalent to  a previous result obtained by Y. Ihara (cf. ~\cite{Ihara1}, or \S\,\ref{QG1901} in this manuscript), who investigated the situation when a given hypergeometric differential equation in characteristic $p$ had a full set of solutions.
  As a result of the explicit description, we  obtained 
   a combinatorial procedure for explicitly computing   the value   $\mr{deg}(\mfO \mfp_{\mfs \mfl_2, X}^\ZZZ/\overline{k})$  according to a sort of fusion rule.
(These fusion rules will be generalized, in Chapter \ref{y0171} of the present manuscript,  to the case of an arbitrary $\mfg$).
The explicit description also   leads to a work by
  F. Liu and B. Osserman, who  have shown (cf. ~\cite{LO}, Theorem 2.1) that
the value $\mr{deg}(\mfO \mfp_{\mfs \mfl_2, X}^\ZZZ/\overline{k}) $
may be expressed as  a degree $3g-3$ polynomial with respect to ``$p$''.
They did so by applying
Ehrhart's theory, which concerns computing   the cardinality of the set of lattice points inside a polytope.

For an arbitrary $g$ and $\mfg = \mfs \mfl_n$,  K. Joshi conjectured, with  amazing insight, an explicit description 
of the value  $\mr{deg}(\mfO \mfp_{\mfs \mfl_n, X}^\ZZZ/\overline{k})$, as  displayed in (\ref{Joshi's conjecture}). 
In ~\cite{Wak}, the author proved,  by grace of that idea and  discussion due to K. Joshi et al. (cf. ~\cite{JP}, ~\cite{Jo14}), 
the following equality:
\begin{align} 
 \label{n=2g}
 \mr{deg}(\mfO \mfp_{\mfs \mfl_2, X}^\ZZZ/\overline{k})
 = \frac{ p^{g-1}}{2^{2g-1}} \cdot  \sum_{\theta =1}^{p-1}\frac{1}{\mr{sin}^{2g-2}(\frac{\pi   \theta}{p})}.  
 \end{align}
(Notice  that   the inequality   $p>2(g-1)$ was imposed   in the statement of ~\cite{Wak}, Corollary 5.4.
 But,  by combining 
 this result 
  with ~\cite{LO}, Theorem 2.1,  one may
  remove this condition; i.e., the asserted equality holds
    even if $g$ is  arbitrary. See also ~\cite{Wak2}, Corollary 8.11.)
 After easy calculations, one verifies that   (\ref{n=2}) and (\ref{n=2g}) are consistent  with   (\ref{Joshi's conjecture}).
 That is to say,  the conjecture of Joshi  has been proved affirmatively   for the case where $n =2$ and $g$ are  arbitrary.
A goal of the present manuscript is to 
prove  
Joshi's conjecture (cf. Theorem \ref{y020} asserted below),  i.e., to give an affirmative answer to the question above,   for the case where
$n$ is an arbitrary integer $>1$.
In ~\cite{Jo14},  the conjectural formula  was formulated under the condition that $p > C$, where $C = n(n-1)(n-2)(g-1)$. But, 
the left-hand side of this formula actually makes sense for every $n$,  and we shall  give its proof under the assumption that
  $p > n \cdot \mr{max} \{ g-1, 2 \}$.

\section{Comparison with arguments due to Mochizuki, Joshi, and Pauly} \label{y03}

The proof of Joshi's conjecture given in the present manuscript may be regarded as a simple generalization of the discussion of the  $n=2$ case.
As   in ~\cite{Wak}, 
our discussion 
 follows, to a substantial extent, the ideas discussed in ~\cite{JP}, and in ~\cite{Jo14}.
Indeed,  some of the results  are mild generalizations of  those obtained in ~\cite{JP} concerning $\mfs \mfl_n$-opers to the case of  {\it families of curves} over quite general base schemes.
(Such {\it relative} formulations allow us to discuss  {\it deformations} of various types of data.)

For example, Proposition \ref{y0202}  in the present manuscript corresponds to ~\cite{JP}, Theorem 3.1.6 (or ~\cite{JRXY}, \S\,5.3; ~\cite{SUN}, \S\,2, Lemma 2.1);
Proposition \ref{QG180}  corresponds to ~\cite{JP}, Theorem 5.4.1;
and Proposition \ref{y0204} corresponds to ~\cite{JP}, Proposition 5.4.2, which associates the moduli space $\mfO \mfp_{\mfs \mfl_n, X}^\ZZZ$ with a certain Quot-scheme in positive characteristic.
Also, as in   ~\cite{Wak},
we prove  the {\it vanishing of obstruction} to 
lifting  the  Quot-scheme to characteristic $0$
  (cf. 
  the proof of Theorem \ref{y0231}).
Thus, in order to complete the proof of Joshi's conjecture,  we relate the value $ \mr{deg}(\mfO \mfp_{\mfs \mfl_n, X}^\ZZZ/\overline{k})$
to the degree of the resulting  Quot-scheme over  $\mbC$ and  {\it directly} apply  
a formula of Holla (cf. Theorem \ref{y0230}), which is a special case of the Vafa-Intriligator formula.
Note that this  idea of  connecting  with the  formula of Y. Holla  is due to K. Joshi (cf. ~\cite{Jo14}).

Our discussion is founded on 
the generic (relative to $X$) \'{e}taleness of  $\mfO \mfp_{\mfs \mfl_n, X}^\ZZZ$ over $\overline{k}$.
As  mentioned above,  S. Mochizuki already proved this fact for $n=2$.
Unfortunately, however, the remaining case (i.e., $n >2$)
 remains unproven.
We will  prove  the generic  \'{e}taleness of  $\mfO \mfp_{\mfs \mfl_n, X}^\ZZZ/\overline{k}$ for an arbitrary $n$, which 
  is one of the main results of the present manuscript (cf. Theorem  \ref{3778567d}).
To this end, following the discussion in 
 ~\cite{Mzk2} (for the case of $n=2$), 
 we  start by  formulating the notion of  an $\mfs \mfl_n$-oper (or, more generally, a $\mfg$-oper for a semisimple Lie algebra $\mfg$) on a {\it family of pointed stable curves},  and then develop  a  theory of such objects.
Moreover, the steps along the path to our goal  (e.g., showing the nonemptyness of the moduli space of opers of prescribed radii and the finiteness of the moduli of dormant opers) will require us to study   a sort of  generalization  of  $\mfg$-opers (cf. Definition \ref{y035}, (i); ~\cite{BD1}, \S\,3.1.14; ~\cite{BD2}, \S\,5.2), which is   called a   {\it $(\mfg, \hslash)$-oper} (where $\hslash$  is  a parameter).

\section{Structure  and main theorems of the manuscript}\label{y04}

\subsection*{Chapter 1: Logarithmic connections over log schemes} \label{QG2100}

In the first half of Chapter 1, 
we discuss  
 the general theory of
 principal  $\mbG$-bundles
 and  $\hslash$-connections, where $\mbG$ is a smooth affine  algebraic group with Lie algebra $\mfg$  and 
$\hslash$ is  a parameter.
  The notion of an $\hslash$-connection was introduced by P. Deligne and used to realize Higgs bundles   as a degeneration of flat vector bundles, i.e., vector bundles equipped with a connection.
 Also, C. T. Simpson studied extensively  relevant  moduli spaces  in a series of papers ~\cite{Simp1}, ~\cite{Simp2}, and ~\cite{Simp}.
Here, we generalize the situation  
by dealing with $\hslash$-flat $\mbG$-bundles (which means   principal $\mbG$-bundles equipped with  a flat $\hslash$-connection) defined on 
a  family of smooth log schemes.
Related notions are reviewed in the context of logarithmic geometry, e.g.,  the Maurer-Cartan form,   Lie algebroids,  and curvature, etc.
One of the important steps for the subsequent  discussion is to describe  gauge transformations of an $\hslash$-connection using
 the exponential and logarithm maps. 
(For example,  the resulting descriptioin will be used to prove  that 
each $(\mfg, \hslash)$-oper has no nontrivial automorphisms.
This fact cannot be ignored from the viewpoint of representability of the relevant moduli categories.) 
If the base field is of characteristic $p>0$, then 
the prime $p$ has to be sufficiently large relative to (a certain invariant associated to) $\mbG$ in order to make that description possible.

In the second half of Chapter 1, we briefly recall the definition of a pointed stable curve and 
the stack $\overline{\mfM}_{g,r}$ over a field $k$ (where $g$ and $r$ are nonnegative integers with $2g-2+r>0$), i.e., the moduli stack classifying $r$-pointed stable curves of genus $g$ over $k$.
 By considering the natural log structure on a pointed stable curve, 
 we can  associate  a certain $\mfg$-valued   invariant to each $\hslash$-flat $\mbG$-bundle and a choice of marked points.
We call this invariant  ``{\it monodromy (operator)}''  according to the wording in ~\cite{Mzk1}, ~\cite{Mzk2}.
It determines  a sort of boundary condition to glue  together $\hslash$-flat $\mbG$-bundles on distinct pointed curves  along the fibers over the points of  attachment.

\subsection*{Chapter 2: Opers on a family of pointed stable curves} \label{QG2103}

In Chapter 2, we study the general theory  of opers and their moduli.
 Given 
 a semisimple algebraic $k$-group $\mbG$ of adjoint type (equipped with a split maximal torus $\mbT$ and a Borel subgroup $\mbB$ defined over $k$ containing $\mbT$) and $\hslash \in k$,
  we introduce the notion of a $(\mfg, \hslash)$-oper, where  $\mfg$ denotes  the Lie algebra of $\mbG$, defined on a family of pointed stable curves  (or more generally, a family of log curves) in  arbitrary characteristic.
 Compared with the previous studies on opers,   this would be  a formulation in a  more general situation because our discussion especially includes the case of singular curves.
 The goal of Chapter 2 is to prove the representability of the moduli category  classifying  pointed stable curves equipped with a $(\mfg, \hslash)$-oper.
 R. C. Gunning 
 proved (cf.  ~\cite{G2}, \S\,4, Corollary 2) that the set of projective structures (or equivalently, $\mfs \mfl_2$-opers) on a smooth curve over $\mbC$ forms an affine space modeled on the space of quadratic differentials. 
The generalization of this result for
  a general $\mfg$ was given by 
 A. Beilinson and V. Drinfeld   (cf. ~\cite{BD2}  \S\,3.4, Theorem; ~\cite{Bar}, Lemma 3.6).
 On the other hand,  S. Mochizuki  (cf. ~\cite{Mzk1}, Chap.\,I, Proposition 2.11) constructed the moduli stack of $\mfs \mfl_2$-opers on the universal family of pointed curves.
(Also, we can find  the construction of $\mr{SL}_n$-opers in positive characteristic in ~\cite{JP}, Proposition 3.2.3.)

There is one thing to keep in mind when generalizing  the construction of the moduli space.
That is, if the underlying curve has a singular point or a marked point,  there may be  a situation where  the moduli problem of $(\mfg, \hslash)$-opers cannot be solved    in terms of flat vector bundles  or differential operators because of   the lack of theta characteristics.
Taking this into consideration,
 we introduce a {\it local}   version of $(\mfg, \hslash)$-opers (cf. Definition \ref{QQ51}) and
  glue together  their moduli spaces defined locally on the underlying curve.
We apply this argument to the universal family of pointed curves $\msC_{g,r}$ parametrized by $\overline{\mfM}_{g,r}$.
Thus, we obtain 
 the  moduli stack 
  \index{$\mfO \mfp_{\mfg, \hslash, (\rho,) \msX}$, $\mfO \mfp_{\mfg, (\hslash,) (\rho,) g,r}$}
\begin{equation} \label{QH1066}
\mfO \mfp_{\mfg, \hslash, g,r} 
\end{equation}
(cf. (\ref{univ1})) representing the category of $r$-pointed stable curves of genus $g$ over $k$  together with a $(\mfg, \hslash)$-oper  on it.
It may be regarded as a stack over $\overline{\mfM}_{g,r}$ via  the projection $\mfO \mfp_{\mfg, \hslash, g,r} \migi \overline{\mfM}_{g,r}$ induced by forgetting the data of $(\mfg, \hslash)$-opers.
Notice that, in order to  construct $\mfO \mfp_{\mfg, \hslash, g,r}$, we apply in advance some of the  results 
 proved in Chapter 4.

A key ingredient in the study of  opers on a pointed curve is
 an invariant called {\it radius} (cf. Definition \ref{y056}).
 The radius of a $(\mfg, \hslash)$-oper at a marked point
 is an element  of the GIT quotient $\mfc := \mfc_{\mfg}$  of the  adjoint $\mbG$-action on $\mfg$   induced  by the monodromy.
It defines a sort of boundary condition
to glue together $(\mfg, \hslash)$-opers 
  according to the attachment of underlying  curves at this point.
 
 Let $\rho := (\rho_i)_{i=1}^r$ be an element of  $\mfc^{\times r}(k) \left(= \mfc (k)^{\times r} \right)$, i.e.,  a $k$-rational point of the $r$-th power product of $\mfc$  
 over $k$, where $\mfc^{\times 0} (k):= \{\emptyset \}$.
Each such element $\rho$ induces
 the substack
  \index{$\mfO \mfp_{\mfg, \hslash, (\rho,) \msX}$, $\mfO \mfp_{\mfg, (\hslash,) (\rho,) g,r}$}
\begin{equation}
 \mfO \mfp_{\mfg, \hslash, \rho, g,r} 
\end{equation}
of $ \mfO \mfp_{\mfg, \hslash, g,r}$ classifying $(\mfg, \hslash)$-oper {\it of radii $\rho$}, i.e., such that the radius  of this oper at the $i$-th marked point  coincides with $\rho_i$ for every $i$.
 In order to ensure 
 the non-emptiness of this stack,
 we observe the deformation family of the moduli spaces $\mfO \mfp_{\mfg, \hslash, g,r}$ parametrized by the value ``$\hslash$'' 
  and the degeneration at $\hslash =0$.
   Summarizing the results  in Chapter 2,  we obtain the following theorem (cf. Theorem \ref{QR90}).
(The case of $\mfg = \mfs \mfl_n$ with $1 < n <p$ will be proved in Chapter \ref{y0104}, Theorem \ref{y0129}.)

\begin{intthm}[Structure of  $\mfO \mfp_{\mfg, \hslash, g,r}$]
\label{y014}
Let $\mbG$ be a semisimple algebraic $k$-group of adjoint type (equipped with a split maximal torus $\mbT$ and a Borel subgroup $\mbB$ defined over $k$ containing $\mbT$.)
 In the case of $\mr{char}(k) =:p >0$, 
suppose further  that  $\mbG$ satisfies either 
``\,$\mbG = \mr{PGL}_n$ with $1<n <p$'' or  ``\,$2  h_\mbG <p$'' (cf. (\ref{QH9932}) for the definition of $h_\mbG$).
Denote by $\mfg$ the Lie algebra of $\mbG$.
Then, for any $\hslash \in k$ and $\rho \in \mfc^{\times r}(k)$,  $\mfO \mfp_{\mfg, \hslash,  g,r}$ (resp., $\mfO \mfp_{\mfg, \hslash, \rho, g,r}$) may be represented by a relative affine space over $\overline{\mfM}_{g,r}$,  and   the  relative dimension
is given by the value
\begin{align}
\aleph (\mfg) := 
\frac{2g-2+r}{2}\cdot \mr{dim}(\mfg) + \frac{r}{2}\cdot  \mr{rk}(\mfg)\hspace{4mm}
 \\
  \left(\text{resp.}, \,   
 {^c \aleph} (\mfg) := 
 \frac{2g-2+r}{2} \cdot \mr{dim}(\mfg) - \frac{r}{2}\cdot   \mr{rk}(\mfg)\right). \index{$\aleph (\mfg)$, ${^c \aleph} (\mfg)$}
 \notag 
 \end{align}
In particular, $\mfO \mfp_{\mfg, \hslash, g,r}$ (resp., $\mfO \mfp_{\mfg, \hslash, \rho, g,r}$) is a
nonempty, geometrically connected, and smooth Deligne-Mumford stack over $k$
 of dimension $3g-3+r + \aleph (\mfg)$ (resp., $3g-3+r + {^c \aleph} (\mfg)$).  
\end{intthm}

\subsection*{Chapter 3: Opers in positive characteristic}\label{y05}

In Chapter \ref{y064}, we study
$(\mfg, \hslash)$-opers in the case of $\mr{char}(k) =:p>0$.
To begin with, we introduce the logarithmic version of a restricted Lie algebroid and recall the $p$-curvature associated to a flat $\mbG$-bundle on  a log scheme.
For each $\rho \in \mfc^{\times r} (k)$, denote by
\index{$\mfO \mfp^\ZZZ_{\mfg, \hslash, (\rho,) \msX}$, $\mfO \mfp^\ZZZ_{\mfg, \hslash, (\rho,) g,r}$, $\mfO \mfp_{n, \hslash, \bb (, \rho)}^{\diamondsuit \ZZZ}$}
\begin{equation}\label{QH1100}
\mfO \mfp^\ZZZ_{\mfg, \hslash, g,r} \ \left(\text{resp}., \ \mfO \mfp^\ZZZ_{\mfg, \hslash, \rho, g,r}\right)
\end{equation}
the closed substack  of $\mfO \mfp_{\mfg, \hslash, g,r}$
classifying  {\it dormant} $(\mfg, \hslash)$-opers (resp., {\it dormant} $(\mfg, \hslash)$-opers of radii $\rho$).
One of the important consequences of this chapter is
the properness of  $\mfO \mfp^\ZZZ_{\mfg, \hslash, g,r} \ \left(\text{resp}., \ \mfO \mfp^\ZZZ_{\mfg, \hslash, \rho, g,r}\right)$ over $k$.
As an intermediate step in proving this fact, we consider a generalization of the  {\it Hitchin-Mochizuki morphism}
\index{$\rotatebox[origin=c]{58}{{\large $\kappa$}}^{\hspace{-0.0mm}}_{\hslash (, \rho)}$, $\rotatebox[origin=c]{58}{{\large $\kappa$}}^{\hspace{-0.0mm}\nabla}_{\hslash (, \rho)}$, $\rotatebox[origin=c]{58}{{\large $\kappa$}}^{\hspace{-0.0mm}\nabla}_{\hslash (, \rho) \mr{univ}}$, Hitchin-Mochizuki morphism}
\begin{equation} \label{QG25880}
 \rotatebox[origin=c]{58}{{\large $\kappa$}}^{\hspace{-0.0mm}\nabla}_{\hslash, \mr{univ}} :  \mfO \mfp_{\mfg, \hslash, g,r} \migi  \ooalign{$\bigoplus$ \cr $\bigotimes^\nabla_{\msC_{g,r}^{(1)}}$}  \   
 \left(\text{resp.}, \ \rotatebox[origin=c]{58}{{\large $\kappa$}}^{\hspace{-0.0mm}\nabla}_{\hslash, \rho, \mr{univ}} :  \mfO \mfp_{\mfg, \hslash, \rho, g,r} \migi  \ooalign{$\bigoplus$ \cr $\bigotimes^\nabla_{\Psi_\hslash (\rho), \msC_{g,r}^{(1)}}$}\right) 
 \end{equation}
(cf. (\ref{QG2500})).
This morphism is defined  as the restriction to the moduli stack of $(\mfg, \hslash)$-opers of a certain $p$-curvature analogue of the Hitchin fibration,  studied in ~\cite{JP}.
When specialized  at $\hslash =0$, it becomes the relative Frobenius morphism of  
$\ooalign{$\bigoplus$ \cr $\bigotimes^\nabla_{\msC_{g,r}^{(1)}}$}$ (resp., $\ooalign{$\bigoplus$ \cr $\bigotimes^\nabla_{\Psi_\hslash (\rho), \msC_{g,r}^{(1)}}$}$) over $\overline{\mfM}_{g,r}$ (cf. Proposition \ref{y093}).
This observation helps us to understand the structure of 
$ \rotatebox[origin=c]{58}{{\large $\kappa$}}^{\hspace{-0.0mm}\nabla}_{\hslash, \mr{univ}}$ (resp., $\rotatebox[origin=c]{58}{{\large $\kappa$}}^{\hspace{-0.0mm}\nabla}_{\hslash, \rho, \mr{univ}}$) for a general $\hslash$, as asserted below.
Here,  we denote the inverse image via
$ \rotatebox[origin=c]{58}{{\large $\kappa$}}^{\hspace{-0.0mm}\nabla}_\hslash$ (resp., $\rotatebox[origin=c]{58}{{\large $\kappa$}}^{\hspace{-0.0mm}\nabla}_{\hslash, \rho}$)
of the zero section  by 
\index{$\mfO \mfp^\NILP_{\mfg, \hslash, (\rho,)  \msX}$, $\mfO \mfp^\NILP_{\mfg, \hslash, (\rho,) g,r}$}
\begin{equation} \label{QH1101}
\mfO \mfp^\NILP_{\mfg, \hslash, g,r} \ \left(\text{resp}., \ \mfO \mfp^\NILP_{\mfg, \hslash, \rho, g,r}\right).
\end{equation}
Each $(\mfg, \hslash)$-oper classified by this stack is called  {\it $p$-nilpotent}.
The main results of Chapter 3 are  the following Theorems \ref{y015} and  \ref{y016} (cf. Theorems \ref{y0102}, \ref{QS1040}, and  \ref{y0101}   for the  full statements).
 Theorem \ref{y015} generalizes  results proved by S. Mochizuki (cf. ~\cite{Mzk1}, Chap.\,I, Theorem 2.3) and T.-H. Chen-X. Zhu
(cf. an unpublished version\footnote{This version is available at:  https://arxiv.org/pdf/1306.0299v1.pdf} of ~\cite{CZ}, \S\,1, Theorem 1.1).
Also,  note that the desired properness of $\mfO \mfp^\ZZZ_{g,\hslash, g,r}$ (resp., $\mfO \mfp^\ZZZ_{g,\hslash, \rho, g,r}$) follows directly from the finiteness of   $\mfO \mfp^\NILP_{g,\hslash, g,r}$ (resp., $\mfO \mfp^\NILP_{g,\hslash, \rho, g,r}$) asserted in Theorem \ref{y015};
this is  because $\mfO \mfp^\ZZZ_{\mfg, \hslash, g,r}$  (resp., $\mfO \mfp^\ZZZ_{\mfg, \hslash, \rho, g,r}$) may be regarded as a closed substack of $\mfO \mfp^\NILP_{\mfg, \hslash, g,r}$ (resp., $\mfO \mfp^\NILP_{\mfg, \hslash, \rho, g,r}$).

\begin{intthm}[Structure of $\mfO \mfp^\NILP_{\mfg, \hslash, g,r}$]  \label{y015}
Let $\mbG$ be as in Theorem \ref{y014} and 
suppose that 
$\mr{char}(k) =:p > 2  h_\mbG$.
Let us  fix an element $\hslash$ of $k^\times$ and an element
$\rho$ of  $\mfc^{\times r}(k)$ with $\Psi_\hslash (\rho) = [0]^{\times r}$ (where $\rho = \emptyset$ if $r =0$).
Then, the  morphism
$\rotatebox[origin=c]{58}{{\large $\kappa$}}^{\hspace{-0.0mm}\nabla}_{\hslash, \mr{univ}}$ (resp., $\rotatebox[origin=c]{58}{{\large $\kappa$}}^{\hspace{-0.0mm}\nabla}_{\hslash, \rho, \mr{univ}}$) is  finite and faithfully flat of  degree $p^{\aleph (\mfg)}$  (resp., $p^{^c \aleph (\mfg)}$).
In particular,
$\mfO \mfp_{\mfg, \hslash, g,r}^\NILP$ (resp., $\mfO \mfp_{\mfg, \hslash, \rho, g,r}^\NILP$) may be represented by a nonempty and proper Deligne-Mumford stack over $k$ of dimension $3g-3+r$, and the natural projection $\mfO \mfp_{\mfg, \hslash, g,r}^\NILP \migi \overline{\mfM}_{g,r}$ (resp., $\mfO \mfp_{\mfg, \hslash, \rho, g,r}^\NILP \migi \overline{\mfM}_{g,r}$) is finite and faithfully flat of degree $p^{\aleph (\mfg)}$  (resp., $p^{^c \aleph (\mfg)}$).
\end{intthm}

\begin{intthm}[Structure of $\mfO \mfp^\ZZZ_{\mfg, \hslash, g,r}$] \label{y016}
Let $\mbG$ be as in Theorem \ref{y014} and
suppose that $\mr{char}(k) = :p > 2  h_\mbG$.
Also,  let us fix  an element $\hslash$ of $k^\times$ and an element  $\rho$ of $ \mfc^{\times r}(k)$  (where  $\rho = \emptyset$ if $r =0$).
Then, 
 $\mfO \mfp^\ZZZ_{\mfg, \hslash, g,r}$
   may be represented  by a  nonempty  proper Deligne-Mumford stack over $k$ of dimension $3g-3+r$.
   Also, $\mfO \mfp^\ZZZ_{\mfg, \hslash, g,r}$   is finite over $\overline{\mfM}_{g,r}$ and has an irreducible component that dominates $\overline{\mfM}_{g,r}$.
   Finally, we have
\begin{equation} 
 \coprod_{\varrho \in \mfc^{\times r}(\mbF_p)} \mfO \mfp^\ZZZ_{\mfg, \hslash, \hslash \star \varrho, g,r} =   \mfO \mfp^\ZZZ_{\mfg, \hslash, g,r}.  \end{equation}
\end{intthm}

\begin{rema} \label{QG3048}
The  non-emptiness of $\mfO \mfp^\ZZZ_{\mfg, \hslash, g,r}$ asserted in  Theorem \ref{y016}
can be verified  by recalling the non-emptiness of $\mfO \mfp^\ZZZ_{\mfs \mfl_2, \hslash, g,r}$ and  constructing a closed immersion  $\mfO \mfp^\ZZZ_{\mfs \mfl_2, \hslash, g,r} \migiincl \mfO \mfp^\ZZZ_{\mfg, \hslash, g,r}$ by  using  a suitable $\mfs \mfl_2$-triple in $\mfg$.
In particular, 
since $\mfO \mfp^\ZZZ_{\mfs \mfl_2, \hslash, g,r}$ is  smooth and of dimension $3g-3+r$, as proved by S. Mochizuki,  the existence of such a closed  immersion implies that the dimension  of  $\mfO \mfp^\ZZZ_{\mfg, \hslash, g,r}$  is exactly   $3g-3+r$.
\end{rema}

\subsection*{Chapters 4-5: Comparison  and duality}\label{y06}

In Chapter \ref{y0104}, we consider the  description of    $(\mfs \mfl_n, \hslash)$-opers  in terms of 
 differential operators and vector bundles.
For a smooth curve  $X$ and a theta characteristic $\Theta$ 
on $X$ (i.e., a square root of the canonical bundle),
a projective connection on  $X$ in the sense of ~\cite{Fr}, \S\,4.1.1,   is a second-order differential operator  $\Theta^{\otimes (-1)} \migi \Theta^{\otimes 3}$  
such that the  principal symbol is $1$ and the subprincipal symbol is $0$.
In the usual manner, it 
 gives a rank $2$ flat vector bundle  equipped with  a line subbundle such that the Kodaira-Spencer map is an isomorphism. 
 Such an object is called a $\mr{GL}_2$-oper and moreover 
induces  
 an $\mfs \mfl_2$-oper via projectivization.
 The resulting correspondence from  projective connections to $\mr{GL}_2$-opers, as well as to   
   $\mfs \mfl_2$-opers,
   is  the simplest case of the main theorem of Chapter 4.
 We refer the reader to    ~\cite{BD2}, \S\,2, for  higher-order differential operators.
   To  generalize the situation (e.g., to families  of not necessarily smooth  stable curves),
we   introduce a generalized  theta characteristic, i.e.,
 an {\it $(n, \hslash)$-theta characteristic} (where $n$ is a positive integer and $\hslash$ is a parameter).
 After  choosing an  $(n, \hslash)$-theta characteristic $\bb$,
 we define the notions of an {\it $(n, \hslash, \bb)$-projective connection} and a {\it $(\mr{GL}_n, \hslash, \bb)$-oper}.
For example,  an $(n, \hslash, \bb)$-projective connection is defined as a certain $n$-th order  differential operator such that the principal symbol is $1$ and  the subprincipal symbol is characterized by $\bb$.
The main result of \S\,\ref{y0104} is as follows (cf. Theorem \ref{y0129} for the full statement):

\begin{intthm}  [{\bf Comparison of  moduli functors}]  \label{y017}
 Let $\msX$ be an $r$-pointed  stable curve of genus $g$ over a $k$-scheme $S$,  $\hslash$   an element of $\Gamma (S, \mcO_S)$, $n$  an integer with $1 < n$,    $\bb := (\BB, \nabla_\bb)$ an $(n, \hslash)$-theta characteristic    
of $\msX$, and $\rho$ an element of $\mfc_{\mfs \mfl_n}^{\times r} (S)$ (where $\rho = \emptyset$ if $r =0$).
Denote by   $\mfS \mfc \mfh_{/S}$ the category of $S$-schemes and by $\Box$ either the absence or presence of $\rho$.
 In the case of $\mr{char}(k) =:p>0$, we suppose further that $n < p$.
Then, there exists 
 a canonical sequence of isomorphisms
\begin{align} \label{QG3000}
\vcenter{\xymatrix@C=46pt@R=36pt{
 \mfD \mfi \mff \mff^\clubsuit_{n,  \hslash, \bb, \Box} 
 \ar[rd]^-{\sim}_-{\Lambda_{\hslash, \bb, \Box, \msX}^{\clubsuit \Rightarrow \diamondsuit}}
 &&  \overline{\mfO} \mfp_{n, \hslash, \Box}^\heartsuit
  \ar[rd]^{\sim}_-{\Lambda_{\Box, \msX}^{\heartsuit \Rightarrow \spadesuit}}& 
\\
&   \mfO \mfp_{n, \hslash, \bb, \Box}^\diamondsuit \ar[ur]_-{\sim}^-{ \Lambda_{\bb, \Box, \msX}^{\diamondsuit \Rightarrow \heartsuit}}& &
\mfO \mfp_{n, \hslash, \Box}^\spadesuit,
}}
\end{align}
where
\begin{itemize}
\item
$\mfO \mfp_{n, \hslash, \Box}^\spadesuit :=$ 
\index{$\mfO \mfp_{n, \hslash (, \rho)}^{\spadesuit}$, $\mfO \mfp_{n, \hslash (, \rho)}^{\spadesuit, \mr{O}}$, $\mfO \mfp_{n, \hslash (, \rho)}^{\spadesuit, \mr{Sp}}$}
the moduli functor  on $\mfS \mfc \mfh_{/S}$ classifying $(\mfs \mfl_n,\hslash)$-opers on $\msX$ of radii $\Box$, i.e., $\mfO \mfp_{n, \hslash, \Box}^\spadesuit := S \times_{s, \overline{\mfM}_{g,r}} \mfO \mfp_{\mfs \mfl_n, \hslash, \Box, g,r}$, where $s : S \migi \overline{\mfM}_{g,r}$ is the classifying morphism of $\msX$;
\item
$\overline{\mfO} \mfp_{n, \hslash, \Box}^\heartsuit :=$  \index{$\overline{\mfO} \mfp_{n, \hslash (,\rho)}^\heartsuit$}
the moduli functor   on $\mfS \mfc \mfh_{/S}$
classifying  {\it  equivalence classes}
of $(\mr{GL}_n, \hslash)$-opers on $\msX$ of radii $\Box$;
\item
$\mfO \mfp_{n, \hslash, \bb, \Box}^\diamondsuit :=$
 \index{$\mfO \mfp_{n, \hslash, \bb (, \rho)}^\diamondsuit$, $\mfO \mfp^{\diamondsuit, \mr{O}}_{n, \hslash}$, $\mfO \mfp^{\diamondsuit, \mr{Sp}}_{n, \hslash}$}
the moduli functor    on $\mfS \mfc \mfh_{/S}$
classifying   isomorphism classes of $(\mr{GL}_n, \hslash, \bb)$-opers on $\msX$ of radii $\Box$;
\item
$\mfD \mfi \mff \mff^\clubsuit_{n, \hslash, \bb, \Box} :=$ 
\index{$\DD^{\,\clubsuit}_{n, \hslash, \BB}$, $\DD^{\,\clubsuit}_{n, \hslash, \bb (, \rho)}$, $\mfD \mfi \mff \mff^{\,\clubsuit}_{n, \hslash, \bb (, \rho)}$}
the moduli functor     on $\mfS \mfc \mfh_{/S}$
 classifying 
$(n, \hslash, \bb)$-projective connections on $\msX$ of radii $\Box$.
\end{itemize}
In particular, all the functors appearing in this sequence may be represented by relative affine spaces over $S$.
\end{intthm}

We would like to make two comments on 
the above theorem in  characteristic $p>0$.
The first is that if
$0 < \mr{char}(k)=p \left(<n\right)$, 
then there always exists an $(n, \hslash)$-theta characteristic   even when it has a singular point or a marked point.
In particular, 
we can construct the moduli functors and   the sequence of isomorphisms   appearing in (\ref{QG3000}) for 
the universal  pointed curves $\msC_{g,r}$. 
Since $\mfD \mfi \mff \mff^\clubsuit_{n, \hslash, \bb}$ has a natural torsor structure modeled on a vector bundle,  it deduces the representability of $\mfO \mfp_{\mfs \mfl_n, g,r}$ (as well as $\mfO \mfp_{\mfs \mfl_n, g,r, \rho}$ for every $\rho$)  asserted in Theorem \ref{y014}.

The second comment concerns the equivalence of  the dormancy condition on opers.
Let ${^p}\mfO \mfp^{\spadesuit}_{n, \hslash, \Box}$
 \index{${^p}\mfO \mfp^{\spadesuit}_{n, \hslash (, \rho)}$, ${^p}\mfO \mfp^{\spadesuit, \mr{O}}_{n, \hslash (, \rho)}$, ${^p}\mfO \mfp^{\spadesuit, \mr{Sp}}_{n, \hslash (, \rho)}$} 
  (where $\Box$ denotes either the absence or presence of $\rho \in \mfc_{\mfs \mfl_n}^{\times r} (S)$) denote  the subscheme   of $\mfO \mfp^{\spadesuit}_{n, \hslash, \Box}$
classifying {\it dormant} $(\mfs \mfl_n, \hslash)$-opers,
  i.e., 
  ${^p}\mfO \mfp^{\spadesuit}_{n, \hslash, \Box} := \mfO \mfp^{\ZZZ}_{\mfs \mfl_n, \hslash, \rho, \msX}$ in the notation of (\ref{QS239}) and (\ref{QS240}).
Also, denote by $\mfD \mfi \mff \mff^{\clubsuit \mr{full}}_{n, 1, \bb, \Box}$ the subfunctor of $\mfD \mfi \mff \mff^{\clubsuit}_{n,  \hslash, \bb, \Box}$ classifying $(n, 1, \bb)$-projective connections {\it having a full set of root functions} in the sense of Definition \ref{D01fgh}.
Then, 
the composite isomorphism  $\mfD \mfi \mff \mff^{\clubsuit}_{n, 1, \bb, \Box} \isom \mfO \mfp^\spadesuit_{n, 1, \Box}$ resulting from  the above theorem  is restricted to an isomorphism
\begin{align}
\mfD \mfi \mff \mff^{\clubsuit \mr{full}}_{n, 1, \bb, \Box} \isom {^p}\mfO \mfp^{\spadesuit}_{n,1, \Box}
\end{align}
(cf. Propositions \ref{QW200} and \ref{QG1780}).
At the end of Chapter 4, 
we review  the previous study  on
 hypergeometric differential operators having a full set of root functions, or equivalently, dormant $\mfs \mfl_2$-oper   on the $3$-pointed projective line.
 This observation allows us to enumerate explicitly  such objects in terms of radii.

In Chapter \ref{p0139}, we discuss 
 $(\mfg, \hslash)$-opers  for some  classical types $\mfg$ other than $\mfs \mfl_n$, i.e., $\mfs \mfo_{2l+1}$ and $\mfs \mfp_{2m}$.
 One important concept related to these opers  is  the duality of $(\mr{GL}_n, \hslash, \bb)$-opers, which  may be interpreted  as  taking  the transposes 
 of the corresponding differential operators.
  In the cases of $\mfg = \mfs \mfo_{2l+1}$ with $2l+1 =n$ and $\mfg = \mfs \mfp_{2m}$ with $2m=n$,  
  we describe $(\mfg, \hslash)$-opers  as certain {\it self-dual}  $(\mr{GL}_n, \hslash, \bb)$-opers equipped with a nondegenerate bilinear form.
This description gives rise to  a canonical closed immersion
  $\mfO \mfp_{n, \hslash, \rho}^{\spadesuit,  \mr{O}} \migiincl \mfO \mfp_{n, \hslash, \rho}^\spadesuit$ (resp.,  $\mfO \mfp_{n, \hslash, \rho}^{\spadesuit,  \mr{Sp}} \migiincl \mfO \mfp_{n, \hslash, \rho}^\spadesuit$), where 
   $\mfO \mfp_{n, \hslash, \rho}^{\spadesuit, \mr{O}}$ 
   \index{$\mfO \mfp_{n, \hslash (, \rho)}^{\spadesuit}$, $\mfO \mfp_{n, \hslash (, \rho)}^{\spadesuit, \mr{O}}$, $\mfO \mfp_{n, \hslash (, \rho)}^{\spadesuit, \mr{Sp}}$}
    (resp., $\mfO \mfp_{n, \hslash, \rho}^{\spadesuit, \mr{Sp}}$) 
   \index{$\mfO \mfp_{n, \hslash (, \rho)}^{\spadesuit}$, $\mfO \mfp_{n, \hslash (, \rho)}^{\spadesuit, \mr{O}}$, $\mfO \mfp_{n, \hslash (, \rho)}^{\spadesuit, \mr{Sp}}$}
    denotes the moduli scheme  on $\mfS \mfc \mfh_{/S}$ classifying $(\mfs \mfo_{2l+1}, \hslash)$-opers (resp., $(\mfs \mfp_{2m}, \hslash)$-opers) on $\msX$ of radii $\rho \in \mfc_{\mfs \mfo_{2l+1}}^{\times r} (S)$ (resp., $\rho \in \mfc_{\mfs \mfp_{2m}}^{\times r} (S)$).
If $\mr{char}(k)=p>0$, 
it is restricted  to a closed immersion   ${^p}\mfO \mfp_{n, \hslash, \rho}^{\spadesuit,  \mr{O}} \migiincl {^p}\mfO \mfp_{n, \hslash, \rho}^\spadesuit$ (resp.,  ${^p}\mfO \mfp_{n, \hslash, \rho}^{\spadesuit,  \mr{Sp}} \migiincl {^p}\mfO \mfp_{n, \hslash, \rho}^\spadesuit$),
where 
   ${^p}\mfO \mfp_{n, \hslash, \rho}^{\spadesuit, \mr{O}}$ 
    \index{${^p}\mfO \mfp^{\spadesuit}_{n, \hslash (, \rho)}$, ${^p}\mfO \mfp^{\spadesuit, \mr{O}}_{n, \hslash (, \rho)}$, ${^p}\mfO \mfp^{\spadesuit, \mr{Sp}}_{n, \hslash (, \rho)}$} 
    (resp., ${^p}\mfO \mfp_{n, \hslash, \rho}^{\spadesuit, \mr{Sp}}$) 
   \index{${^p}\mfO \mfp^{\spadesuit}_{n, \hslash (, \rho)}$, ${^p}\mfO \mfp^{\spadesuit, \mr{O}}_{n, \hslash (, \rho)}$, ${^p}\mfO \mfp^{\spadesuit, \mr{Sp}}_{n, \hslash (, \rho)}$} 
  denotes  the 
closed subscheme of $\mfO \mfp_{n, \hslash, \rho}^{\spadesuit, \mr{O}}$ (resp., $\mfO \mfp_{n, \hslash, \rho}^{\spadesuit, \mr{Sp}}$)  classifying dormant ones.
The existence of such a closed immersion  deduces the
 generic unramifiedness 
  of $\mfO \mfp_{\mfg, g,r}^\ZZZ$  ($\mfg = \mfs \mfo_{2l+1}, \mfs \mfp_{2m}$)  from that of $\mfO \mfp_{\mfs \mfl_n, g,r}^\ZZZ$.
  Moreover, 
  the duality 
  yields a bijection  between (dormant) $(\mfg, \hslash)$-opers of radii $\rho$ and (dormant) $(\mfg, \hslash)$-opers of radii $(-1) \star \rho$ for $\mfg = \mfs \mfl_n, \mfs \mfo_{2l+1}, \mfs \mfp_{2m}$.
  These observations are  essential ingredients for constructing the pseudo-fusion ring  associated to  dormant $(\mfg, \hslash)$-opers for such $\mfg$'s, as discussed in Chapter 7.
In the following, we summarize the main  assertions proved  in Chapter 5  (cf. Theorems \ref{yp3128} and  \ref{ppp012} for the statements).

\begin{intthm}  [{\bf Duality of (dormant) opers}]  \label{ypp017}
Let $\msX$, $\hslash$, and $n$ be as in Theorem \ref{y017}.
Then, the following assertions hold:
\begin{itemize}
\item[(i)]
For each $\rho \in \mfc_{\mfs \mfl_n}^{\times r} (S)$, 
there exists a canonical $S$-isomorphism
 \index{$\Dual_{\bb (, \rho)}^{\, \diamondsuit}$, $\Dual_{(\rho)}^{\,\heartsuit}$, $\Dual_{(\rho)}^{\, \spadesuit}$, ${^p}\Dual_{(\rho)}^{\, \spadesuit}$}
\begin{align}
\Dual^{\,\spadesuit}_{\rho} : \mfO \mfp^{\spadesuit}_{n, \hslash, \rho}  \isom \mfO \mfp^{\spadesuit}_{n, \hslash, (-1) \star\rho}
\end{align}
satisfying the equalities  $\Dual_{(-1) \star \rho}^{\,\spadesuit} \circ \Dual^{\,\spadesuit}_{\rho} = \mr{id}$.
Moreover, if  $\mr{char}(k) =:p>0$,
then  $\Dual^{\,\spadesuit}_{\bb, \rho}$ restricts to an isomorphism
\begin{align}
{^p}\Dual_{\rho}^{\, \spadesuit}
: {^p}\mfO \mfp^{\spadesuit}_{n, \hslash, \rho}  \isom {^p}\mfO \mfp^{\spadesuit}_{n, \hslash, (-1) \star\rho}.
\end{align}
\item[(ii)]
Suppose that there exists an integer $l$ (resp., $m$) with $2l+1 =n$ (resp., $2m= n$) and $\rho$ lies in $\mfc_{\mfs \mfo_{2l+1}}^{\times r} (S)$
 (resp., $\mfc_{\mfs \mfp_{2m}}^{\times r} (S)$).
 (In particular, we have $\rho = (-1)\star \rho$).
 Then, 
 $\Dual^{\,\spadesuit}_{\rho}$ restricts to the identity morphism of $\mfO \mfp^{\spadesuit, \mr{O}}_{n, \hslash, \rho}$ (resp., $\mfO \mfp^{\spadesuit, \mr{Sp}}_{n, \hslash, \rho}$).
 If, moreover, $\mr{char}(k) =:p >0$,
 then
 ${^p}\Dual_{\rho}^{\, \spadesuit}$
 restricts to the identity morphism of ${^p}\mfO \mfp^{\spadesuit, \mr{O}}_{n, \hslash, \rho}$ (resp., ${^p}\mfO \mfp^{\spadesuit, \mr{Sp}}_{n, \hslash, \rho}$).
\end{itemize}
 \end{intthm}

\subsection*{Chapters 6-7: Deformation and degeneration of (dormant) opers}\label{y07}

In Chapter \ref{y0139}, 
we study  local deformations of (dormant) $(\mfg, \hslash)$-opers, where $\mfg$ is as before.
First, we review classical deformation theory of  pointed stable curves  and flat vector bundles in the context of logarithmic geometry.
After that, we describe  the deformation space of  a $(\mfg, \hslash)$-oper (of prescribed radii)  using  the de Rham cohomology of a certain complex     associated to that oper (cf. Propositions \ref{y0148},  \ref{y0151}).
A key observation concerning   flat bundles  in characteristic $p>0$ is  that 
the Cartier operator may be interpreted as
describing deformation of $p$-curvature (cf. Proposition \ref{y0157}).
By taking this into account, 
we achieve a detailed understanding of  the tangent space of the moduli stack classifying  dormant $(\mfg, \hslash)$-opers (cf. Proposition \ref{y0167}, Corollary \ref{y0169}).
One of the  results in Chapter \ref{y0139}  is  a criterion for the \'{e}taleness of that stack (cf. Corollary \ref{y0170}), which plays a crucial  role in  the proof of  Joshi's conjecture.

In Chapter \ref{y0171}, we examine  the behavior,  according to decomposition arising from various  clutching morphisms ,  
of   the generic degree 
 \begin{equation}
 \mr{deg}(\mfO\mfp^\ZZZ_{\mfg, \hslash, \hslash \star \rho, g, r}/\overline{\mfM}_{g,r }) \index{$\mr{deg}(\mfO\mfp^\ZZZ_{\mfg, \hslash, \hslash \star \rho, g, r}/\overline{\mfM}_{g,r })$} 
 \end{equation}
  of  $\mfO\mfp^\ZZZ_{\mfg, \hslash, \hslash \star \rho, g, r}/\overline{\mfM}_{g,r }$ for each $\rho := (\rho_i)_{i=1}^r \in \mfc^{\times r} (\mbF_p)$ and $\hslash \in k^\times$.
This value depends  neither on the choice of $\hslash \in k^\times$ nor  on the ordering of $\rho_1, \cdots, \rho_r$; this is why
there exists a well-defined function
\index{${^p}N_{\mfg,  0}$, ${^p}N_{\mfg,  g}$}
\begin{align} \label{QH2556}
{^p}N_{\mfg,  g} :  \mbN^{\mfc (\mbF_p)} 
\migi \mbZ; 
  \overline{\rho} \left(:=\textstyle{\sum_{i =1}^r \rho_i}\right) 
  \mapsto  \mr{deg}(\mfO\mfp^\ZZZ_{\mfg, \hslash, \hslash \star \rho, g, r}/\overline{\mfM}_{g,r }) 
\end{align}
 (cf. (\ref{f186})) on 
  the free commutative monoid $\mbN^{\mfc (\mbF_p)}$ with basis 
$\mfc (\mbF_p)$.
Under certain assumptions,
${^p}N_{\mfg, g} (\overline{\rho})$
coincides  precisely  with the number of  dormant $(\mfg,\hslash)$-opers  on a sufficiently general $r$-pointed stable curve of genus $g$ over an algebraically  closed field.
Also, as proved in  (cf. Proposition  \ref{c2}), 
 the function ${^p}N_{\mfg, 0}$
   satisfies 
a certain factorization property,  which will be referred to as  a  ``{\it pseudo-fusion rule}" on $\mfc (\mbF_p)$ (cf. Definition \ref{c34}).
The notion of a pseudo-fusion rule  may be regarded as a somewhat weaker variant of  a  fusion rule.
(For reviews and expositions concerning fusion rules, we refer the reader to ~\cite{Beau} or  ~\cite{Fuchs}.)
The latter half of Chapter \ref{y0171} is mainly devoted to developing  a general theory of pseudo-fusion rules.
By applying the results proved  there, we obtain a way to calculate 
the values ${^p}N_{\mfg, g} (\overline{\rho})$ (for various  $g \geq 0$) combinatorially  from the case of  $(g,r)= (0,3)$.
Moreover, 
the function ${^p}N_{\mfg, g}$ induces
a commutative  ring  
 \index{${^p}\Fus_{\mfg}$, pseudo-fusion ring for dormant $\mfg$-opers}
\begin{equation}  \label{QH27789}
 {^p}\Fus_{\mfg}, 
 \end{equation}
which will be  referred to as the {\it pseudo-fusion ring for dormant $\mfg$-opers} (cf. Definition \ref{c4}).
 At the end of Chapter \ref{y0171},
 we explicitly describe  the structure of ${^p}\Fus_{\mfs \mfl_2}$. 
In general,  the explicit description of 
 (the characters of) ${^p}\Fus_{\mfg}$ allows us to 
 perform a  computation for
  various values   ${^p}N_{\mfg,  g}(\overline{\rho})$, as 
asserted   in the following  (cf. Theorem \ref{c6} (ii)):

\begin{intthm}[Factorization property] \label{y018}
Suppose that  $\mr{char}(k) =:p >0$ and that  $\mfg$ is either $\mfs \mfl_n$ with $2n < p$,  $\mfs \mfo_{2l+1}$ with $4l+2<p$, or $\mfs \mfp_{2m}$ with $4m < p$.
Let $\rho = (\rho_i)_{i=1}^r \in \mfc^{\times r} (\mbF_p)$ (where  $\rho = \emptyset$ if $r =0$) and $\hslash \in k^\times$. 
Write  ${^p}\mfS_{\mfg}$ \index{${^p}\mfS_{\mfg}$}
 for the set of ring homomorphisms ${^p}\Fus_{\mfg} \migi \mbC$ and write ${^p}\mr{Cas}_\mfg := \sum_{\lambda \in \mfc (\mbF_p)} \lambda \ast \lambda^\dual \left( \in {^p}\Fus_{\mfg}\right)$. \index{${^p}\mr{Cas}_\mfg$}
Then, we have the following equality:
\begin{align}   \label{1060}
\mr{deg}(\mfO\mfp^\ZZZ_{\mfg, \hslash, \hslash \star \rho,  g, r}/\overline{\mfM}_{g,r})   \left(= 
{^p}N_{\mfg, g} (\sum_{i=1}^r \rho_i) \right)
  = \sum_{\chi \in {^p}\mfS_{\mfg}} \chi ({^p}\mr{Cas}_\mfg)^{g-1}\cdot \prod_{i=1}^r \chi(\rho_i).  \end{align}
If, moreover, $r  =0$ (hence $g >1$),  this equality means that
\begin{align}   \label{1dd060}
\mr{deg}(\mfO\mfp^\ZZZ_{\mfg, \hslash,  g, 0}/\overline{\mfM}_{g,0}) \ \left(  = 
{^p}N_{\mfg, g} (0) = {^p}N_{\mfg, g} (\overline{\emptyset})\right)  = \sum_{\chi \in {^p}\mfS_{\mfg}} \chi ({^p}\mr{Cas}_\mfg)^{g-1}. 
\end{align}

\end{intthm}

\subsection*{Chapters 8-9: A proof of Joshi's conjecture} \label{y08}

Chapter \ref{y0186} is devoted to proving the generic \'{e}taleness of
 $\mfO \mfp_{\mfg, \hslash, \hslash \star \rho, g,r}^\ZZZ$ over $\overline{\mfM}_{g,r}$ for  certain  classical $\mfg$'s.
For $\mfg = \mfs \mfl_2$,  this fact was already proved
 by S. Mochizuki (cf. ~\cite{Mzk2}, Chap.\,II, \S\,2.8, Theorem 2.8).
Thanks to the results  in
the previous chapters, 
  it suffices to prove that there is no nontrivial deformation of any dormant $\mfs \mfl_n$-oper on the $3$-pointed projective line  $\msP$.
We  give a proof of this claim (cf. Corollary \ref{002})  based on the idea that {\it infinitesimal deformations of a dormant oper on  $\msP$ are completely determined  by their radii and the sheaves of their horizontal sections} (as shown in ~\cite{Mzk2}).
Consequently,   we obtain the following theorem  (cf.  Theorem   \ref{3778567d}):

\begin{intthm}[Generic \'{e}taleness] \label{QG3457}
Let $k$ be a perfect field of characteristic $p>0$, $\hslash$ an element of $k^\times$, and $\rho$ an element of $\mfc^{\times r}_{\mfg} (\mbF_p)$ (where $\rho = \emptyset$ if $r =0$).
Suppose  that  $\mfg$ is either $\mfs \mfl_n$ with $2n < p$,  $\mfs \mfo_{2l+1}$ with $4l<p$, or $\mfs \mfp_{2m}$ with $4m < p$.
  Then, 
   $\mfO \mfp_{\mfg, \hslash, \hslash \star \rho, g,r}^\ZZZ$ is 
\'{e}tale over the points of $\overline{\mfM}_{g,r}$ classifying totally degenerate curves.
In particular, $\mfO \mfp_{\mfg, \hslash, \hslash \star \rho, g,r}^\ZZZ$ is 
 generically  \'{e}tale over $\overline{\mfM}_{g,r}$; i.e., any irreducible component that dominates  $\overline{\mfM}_{g,r}$ admits a dense open subscheme which is \'{e}tale over  $\overline{\mfM}_{g,r}$.
\end{intthm}

In Chapter \ref{y0198}, we give a proof of Joshi's conjecture by  investigating the relationship    between
 (the fiber over a point in $\overline{\mfM}_{g,0}$ of) the moduli stack  $\mfO \mfp_{\mfs \mfl_n, 1, g,0}^\ZZZ$ and  certain Quot-schemes  in positive characteristic,  denoted by $\mr{Quot}^{n, 0}_\BB$ and $\mr{Quot}^{n, \mcO}_\BB$.
These Quot-schemes are associated to
a smooth curve $X/S$ equipped with  a theta  characteristic $\BB$, and
 both classify rank $n$ subbundles of a rank $p$ vector bundle on $X$
   induced from $\BB$.
The  discussion in Chapter \ref{y0198} is 
parallel to the one  in ~\cite{Wak}, \S\S\,4-5, where the case of $\mfg = \mfs \mfl_2$ was studied.
Just as in {\it loc.\,cit.},
we construct   an isomorphism between the moduli space of dormant opers  and  $\mr{Quot}^{n, \mcO}_\BB$
 via  Cartier descent  of connections with vanishing $p$-curvature  (cf. Proposition \ref{y0204}).
In particular,  (under the assumption that $X/S$ is sufficiently general in $\overline{\mfM}_{g, 0}$) the problem can be reduced  to  a computation of the degree of
$\mr{Quot}^{n, \mcO}_\BB$, and moreover, to that  of the other Quot-scheme $\mr{Quot}^{n, 0}_\BB$.
By   the generic \'{e}taleness of $\mfO \mfp_{\mfs \mfl_n, 1, g, 0}^\ZZZ$ resulting from Theorem \ref{QG3457},
this value coincides with the degree of a certain Quot scheme {\it over $\mbC$} obtained by deforming the underlying curve and the rank $p$ vector bundle on it.
On the other hand, 
there is a formula  given by Y. Holla  for computing  the degree of this Quot-scheme, which  may be regarded as  a special case of the Vafa-Intriligator formula for explicitly computing  the Gromov-Witten invariants.
By applying this formula,
 we thus prove  the following theorem   (cf. Theorem \ref{y0251}).
 That is to say,  {\it Joshi's conjecture holds for sufficiently general curves}.

\begin{intthm}[{\bf Joshi's conjecture}]   \label{y020}
Suppose that $p>n \cdot \mr{max} \{g-1, 2 \}$ and $\hslash \in k^\times$.
  Then, the generic degree $\mr{deg}(\mfO \mfp_{\mfs \mfl_n, \hslash, g,0}^\ZZZ/\overline{\mfM}_{g, 0}) $
 of $\mfO \mfp_{\mfs \mfl_n, \hslash, g,0}^\ZZZ$ over $\overline{\mfM}_{g,0}$
 is given by the following formula:
 \vspace{3mm}
 \begin{equation} \mr{deg}(\mfO \mfp_{\mfs \mfl_n, \hslash, g,0}^\ZZZ/\overline{\mfM}_{g, 0})    =    \frac{ p^{(n-1)(g-1)-1}}{n!} \cdot  
 \sum_{\genfrac{.}{.}{0pt}{}{(\zeta_1, \cdots, \zeta_n) \in \mbC^{\times n} }{ \zeta_i^p=1, \ \zeta_i \neq \zeta_j (i\neq j)}}
\frac{(\prod_{i=1}^n\zeta_i)^{(n-1)(g-1)}}{\prod_{i\neq j}(\zeta_i -\zeta_j)^{g-1}}.  
 \end{equation}
\end{intthm}


\section{Acknowledgements} 
We cannot  express enough our sincere and deep gratitude to Professors Shinichi Mochizuki and Kirti Joshi.
Without their philosophical viewpoints, theoretical insights, and endless creativity, my study of mathematics would have remained ``{\it dormant}".
We deeply appreciate Professor Yuichiro Hoshi  
for his  helpful suggestions, as well as for reading preliminary versions of the present manuscript.
However, responsibility for any errors or misconceptions in the present manuscript  ultimately remains with the author. 
Also, 
we would like to
express our sincere appreciation to   the referee
 for  a careful reading of the manuscript and for various useful comments and suggestions.
Finally, we would  like to thank 
 Professors Go Yamashita and Akio Tamagawa (for their constructive comments concerning the discussion in Chapter \ref{y0171}),  and the various individuals (including {\it pointed stable curves and their moduli space}!) with whom the author became acquainted in Kyoto.
{\it I wrote  the present manuscript as  a gratitude letter to them}.

\section{Notation and Coventions} 

\subsection*{Categories and Schemes}
Denote by   $\mfS \mfe \mft$  \index{$\mfS \mfe \mft$, category of  sets} the category of (small) sets.
Throughout this manuscript, all schemes (and Deligne-Mumford stacks) are assumed to
be locally noetherian. 
 For a scheme $S$, we denote by 
  $\mfS \mfc \mfh_{/S}$ \index{$\mfS \mfc \mfh_{/S}$, category of  $S$-schemes} the category of  schemes over $S$ and by 
   $\mfE  \mft_{/S}$ 
   \index{$\mfE  \mft_{/S}$, small \'{e}tale site on $S$}
     the small \'{e}tale site on $S$.
  If $S = \mr{Spec}(k)$ for a field $k$, then we write
  $\mfS \mfc \mfh_{/k} := \mfS \mfc \mfh_{/\mr{Spec}(k)}$.

\subsection*{Fields} \label{QH8000}
Let $k$ be a field.
Then,  $\mr{char}(k)$ \index{$\mr{char}(k)$} denotes  the characteristic of $k$.
For a $k$-vector space $\mfh$, we denote by  $\mbS_k (\mfh)$ \index{$\mbS_k (\mfh)$} the symmetric algebra on $\mfh$ over $k$.
Occasionally, we will regard $\mfh$ as the affine  $k$-scheme 
$\mr{Spec}(\mbS_k (\mfh^\vee))$, where $\mfh^\vee$ is the dual space of $\mfh$.

\subsection*{Sheaves}\label{QH1120}
Let $S$ be a scheme.
 For a sheaf $\mcV$ on $S$, we use the notation ``\,$v \in \mcV$\,''
   for a local section $v$ of a sheaf $\mcV$.
If $\mcV$ is a $\mcO_S$-module,  then we denote by $\mcV^\vee$ \index{$\mcV^\vee$, dual of $\mcV$} the dual of $\mcV$, i.e., $\mcV^\vee := \mcH om_{\mcO_S}(\mcV, \mcO_S)$.
Next, by a {\bf vector bundle} \index{vector bundle} on $S$, we mean a locally free $\mcO_S$-module of finite rank. 
For a vector bundle $\mcV$ on $S$ and a nonnegative integer  $i$,
we denote  by $\bigwedge^i \mcV$ the $i$-th exterior power of $\mcV$.
Also, $\mbS_{\mcO_Y} (\mcV )$ \index{$\mbS_{\mcO_Y} (\mcV )$, $\mbS^j_{\mcO_Y} (\mcV )$} denotes   the symmetric algebra of $\mcV$ over $\mcO_S$; it decomposes into the direct sum 
$\mbS_{\mcO_S} (\mcV ) = \bigoplus_{j \geq 0}\mbS^j_{\mcO_S}(\mcV)$, where 
 $\mbS^j_{\mcO_S}(\mcV)$ ($j \geq 0$) is the  subsheaf consisting of sections of homogenous degree $j$.
 If $X$ is a scheme (or more generally, a Deligne-Mumford stack) over $S$,
 then we shall write $\mcT_{X/S}$ \index{$\mcT_{X/S}$, tangent bundle on $X/S$} (resp., $\Omega_{X/S}$) \index{$\Omega_{X/S}$, sheaf of $1$-forms on $X/S$} for the tangent bundle (resp., the sheaf of $1$-forms) on $X$ relative to $S$. 

\subsection*{Relative affine schemes}\label{QH2002}
Let  $s : S' \migi S$ be an affine morphism of schemes.
Then, we shall write $\mcO_S \langle S' \rangle$ \index{$\mcO_S \langle S' \rangle$, $\mcO_S$-algebra corresponding to $S'$} for the $\mcO_S$-algebra corresponding to $S'$, i.e., $\mcO_S \langle S' \rangle := s_*(\mcO_{S'})$ and $\mcS pec (\mcO_S \langle S' \rangle) \cong S'$.
Let  $S$ be  a scheme, or more generally a Deligne-Mumford stack,  and $\mcV$   a vector bundle on $S$.
 Denote by  
 $\mbV (\mcV)$ \index{$\mbV (\mcV)$} 
 the relative  affine scheme over $S$ associated to $\mcV$, i.e., the spectrum $\mcS pec(\mbS_{\mcO_T}(\mcV^\vee))$ of  $\mbS_{\mcO_T}(\mcV^\vee)$.
If $S$ is defined over a field $k$ and 
 $v: \mr{Spec}(R) \migi S$  is a $k$-morphism, where $R$ denotes a ring over $k$, then 
we  denote by $\mcT_v S \left(:= v^*(\mcT_{S/k}) \right)$ \index{$\mcT_v S$, tangent space of $S$ at $v$}  the $R$-module defined to be the tangent space of $S$ at $v$.

\subsection*{Complexes}\label{QH2001}

Let $\nabla : \mcK^0 \migi \mcK^1$ be a morphism of sheaves of abelian groups on a scheme $S$.
It may be regarded as a complex concentrated at degrees $0$ and $1$; we denote this complex by 
$\mcK^\bullet[\nabla]$. 
\index{$\mcK^\bullet[\nabla]$, complex defined by $\nabla$}
Next, let $f : X \migi S$ be  a morphism of schemes and $i$ a nonnegative integer.
 Then,
one may define the sheaf
  $\mbR^i f_*(\mcK^\bullet [\nabla])$ 
on $S$ obtained from $\mcK^\bullet [\nabla]$ by applying 
the $i$-th hyper-derived functor $ \mbR^i f_*(-)$  of $f_*(-)$ (cf. ~\cite{Kal}, (2.0)).
In particular, $\mbR^0f_*(\mcK^\bullet [\nabla]) = f_*(\mr{Ker}(\nabla))$.
If $S$ is affine, then $\Gamma (S,   \mbR^i f_*(\mcK^\bullet [\nabla]))$ may be identified with the $i$-th hypercohomology group 
$\mbH^i (S, \mcK^\bullet [\nabla])$.
Given an integer $n$ and a sheaf $\mcF$ on $S$,
  we define the complex
$\mcF [n]$ \index{$\mcF [n]$} to be  $\mcF$ (considered as a complex concentrated at degree $0$) shifted down by $n$, so that $\mcF[n]^{-n} = \mcF$ and $\mcF[n]^i = 0$ ($i \neq -n$).

\subsection*{Algebraic groups}\label{QH2000}
Suppose that 
we are working over a field $k$.
Then, denote by   $\mbG_m$ $\left(:= \mr{Spec}(k[z^{\pm 1}]) \right)$ 
\index{$\mbG$, algebraic group@--- $\mbG_m$, multiplicative group}
  the multiplicative group over $k$.
For a positive  integer $n$, 
we shall write $\mr{GL}_n$ 
\index{$\mbG$, algebraic group@--- $\mr{GL}_n$, general linear group}
 (resp., $\mr{PGL}_{n}$) 
 \index{$\mbG$, algebraic group@--- $\mr{PGL}_{n}$, projective linear group}
  for the general (resp.,  projective) linear group  of rank $n$ over $k$.
The Lie algebra $\mfg \mfl_n$ 
\index{$\mfg$, Lie algebra@--- $\mfg \mfl_n$, Lie algebra of $\mr{GL}_n$}
 of $\mr{GL}_n$ consists of $n \times n$ matrices over $k$.
Denote by $\mfs \mfl_n$ \index{$\mfg$, Lie algebra@--- $\mfs \mfl_n$, Lie algebra of $\mr{SL}_n$}
 the Lie subalgebra of $\mfg \mfl_n$ consisting of $n \times n$ matrices with vanishing trace.
If $n$ is invertible in $k$, then the Lie algebra $\mfp \mfg \mfl_n$
 \index{$\mfg$, Lie algebra@--- $\mfp \mfg \mfl_n$, Lie algebra of $\mr{PGL}_n$}
  of $\mr{PGL}_n$ is  isomorphic to $\mfs \mfl_n$  via the natural composite $\mfs \mfl_n \migiincl \mfg \mfl_n \migisurj \mfp \mfg \mfl_n$;  we will  occasionally identify $\mfs \mfl_n$ with $\mfp \mfg \mfl_n$ via this composite isomorphism.
Next, let $\mbG$ be an affine algebraic group over $k$ and
$\mfh$ a $k$-vector space equipped with a $\mbG$-action.
Then, the dual $\mfh^\vee$ of $\mfh$ is equipped with a $\mbG$-action given by $(h \cdot f)(v) := f (h^{-1} \cdot v)$ for $h \in \mbG$, $f \in \mfh^\vee$, and $v \in \mfh$.
In particular,  this $\mbG$-action induces a $\mbG$-action on the coordinate ring $\mbS_k (\mfh^\vee)$ of $\mfh$.



\chapter{Logarithmic connections over log schemes}  \label{y010}


In the first half of this chapter, we discuss a general theory of principal bundles and $\hslash$-connections, where $\hslash$ is 
a parameter,
  in the  context of logarithmic algebraic geometry.
 We shall refer the reader to, e.g., ~\cite{Simp1}, ~\cite{Simp2}, ~\cite{Simp},  ~\cite{Ar2}, ~\cite{Ar1}, and ~\cite{Tra} for the study of  $\hslash$-connections.
 One of the important results in this chapter  is the description of gauge transformations using the exponential and logarithm maps.
 This fact will be used to prove, e.g., that a $(\mfg, \hslash)$-oper has no nontrivial automorphisms.
 In the second half of this chapter, we briefly recall  pointed stable curves and the moduli stack  classifying them.
 By   considering the natural log structure on a pointed stable curve, we can associate a certain invariant, called ``{\it monodromy (operator)}'',  to each $\hslash$-flat principal bundle on the curve  and a choice of a marked point. 
 As we will discuss in Chapter \ref{y0171}, it determines a sort of boundary condition to glue together $\hslash$-flat principal bundles on pointed curves along the fibers over the points of attachment.

\begin{notation}
Let us fix a field $k$
\index{$k$, field} 
 and   a connected smooth affine algebraic group $\mbG$  
\index{$\mbG$, algebraic group}
over $k$.
 Denote by 
 $e$
  \index{$\mbG$, algebraic group@--- $e$, identity element of}
  the identity element of $\mbG$  and by 
 $\mfg$ 
 \index{$\mfg$, Lie algebra}
 the  Lie algebra of $\mbG$, i.e., the tangent space  $\mcT_e \mbG$
 of $\mbG$ at $e$.
Occasionally,  we regard  $\mfg$ as a $k$-scheme.
\end{notation}

\section{The adjoint representation and the Maurer-Cartan form} \label{y011}

\subsection{Derivations}\label{QS4868}
Let $T$ be a $k$-scheme and $f: Y \migi T$ a $T$-scheme.
Each  vector bundle  $\mcV$ on $Y$ defines   
 $\mcE nd_{f^{-1}(\mcO_T)}(\mcV)$ (resp., $\mcE nd_{\mcO_Y} (\mcV)$),   the sheaf  of $f^{-1}(\mcO_T)$-linear  (resp., $\mcO_Y$-linear) endomorphisms of $\mcV$.
Denote by 
\begin{equation} \label{QQ77}
 \Delta  : \mcO_Y \migi \mcE nd_{f^{-1}(\mcO_T)}(\mcV)
 \index{$\Delta$}
\end{equation}
 the diagonal embedding, i.e.,  the $f^{-1}(\mcO_T)$-linear morphism which, to any local section $a \in \mcO_Y$, assigns the locally defined endomorphism  of  $\mcV$ given by multiplication by $a$.
The local freeness of  $\mcV$ implies that $\Delta$ is  injective.
The sheaf $ \mcE nd_{f^{-1}(\mcO_T)}(\mcV)$ may be regarded as an
 $\mcO_Y$-module
in such a way that  $a \cdot D := \Delta (a) \circ D$ for any $a \in \mcO_Y$ and $D \in \mcE nd_{f^{-1}(\mcO_T)}(\mcV)$.
Also, $ \mcE nd_{f^{-1}(\mcO_T)}(\mcV)$ has  a structure of $\mcO_T$-Lie algebra on $Y$  (in the sense of ~\cite{Lan}, \S\,2.1), where 
 the Lie  bracket operation   is given by $[D_1, D_2] := D_1 \circ D_2 - D_2 \circ D_1$.
The sheaf $\mcE nd_{\mcO_Y}(\mcV)$ is   the $\mcO_T$-Lie subalgebra of $\mcE nd_{f^{-1}(\mcO_T)}(\mcV)$ consisting of local sections that commute with $\mr{Im}(\Delta)$.

Recall that
a {\bf derivation} \index{derivation} on $\mcO_Y$ relative to $T$ is an $f^{-1}(\mcO_T)$-linear endomorphism $D$ of $\mcO_Y$ (i.e., a global section of $\mcE nd_{f^{-1}(\mcO_T)}(\mcO_Y)$) that satisfies 
\begin{align}
D (a \cdot b) = a \cdot D (b)  + D(a) \cdot b
\end{align}
for any local sections $a, b\in \mcO_Y$. 
In particular,  we have  $D (a) =0$   for any $a \in f^{-1}(\mcO_T)$.
The  tangent bundle $\mcT_{Y/T} \left(:= \Omega_{Y/T}^\vee\right)$ of $Y/T$
may be identified, in the usual manner,  with the sheaf of locally defined derivations on $\mcO_Y$ relative to $T$;
 this sheaf  specifies an $\mcO_T$-Lie subalgebra of $\mcE nd_{f^{-1}(\mcO_T)} (\mcO_Y)$, meaning that $\mcT_{Y/T}$
   is closed under both the $\mcO_Y$-action and the Lie bracket operation    on  $\mcE nd_{f^{-1}(\mcO_T)} (\mcO_Y)$.

\subsection{Right/Left-translations}
Let $T$ be a scheme over $k$ and
$h : T \migi \mbG$  a morphism of $k$-schemes.
 Denote by $\Gamma_h \left(:= (\mr{id}_T, h) \right) : T \migi T \times_k \mbG$ \index{$\Gamma_h$, graph of $h$} the graph of $h$.
We write 
\begin{align} \label{QR180}
\mr{R}_h \ \left(\text{resp.,}  \ \mr{L}_h \right)
 \index{$\mr{R}_h$, right-translation} \index{ $\mr{L}_h$, left-translation}
 \end{align}
  for the right-translation  (resp., the left-translation) by $h$, i.e., the composite isomorphism
\begin{align} \label{translation}
& \hspace{12mm} \mr{R}_h : T \times_k \mbG \xrightarrow{(\mr{id}_T \times \mr{sw}) \circ (\Gamma_h \times \mr{id}_\mbG)} T \times_k  \mbG \times_k \mbG  \xrightarrow{\mr{id}_T \times \mr{mult}} T \times_k \mbG  \\
& \left(\text{resp.,} \ \mr{L}_h : T \times_k \mbG  \xrightarrow{\Gamma_h\times \mr{id}_\mbG} T \times_k \mbG \times_k \mbG  \xrightarrow{\mr{id}_T \times \mr{mult}} T \times_k \mbG \right),\notag
\end{align}
where ``$\mr{mult}$'' denotes the multiplication $\mbG \times_k \mbG \migi \mbG$ in $\mbG$ and ``$\mr{sw}$'' denotes  the automorphism  of $\mbG \times_k \mbG$ given by switching the two factors  $(h_1, h_2) \mapsto (h_2, h_1)$.
Note that 
\begin{align}
\mr{R}_{h_1} \circ \mr{R}_{h_2} = \mr{R}_{h_2 h_1}, 
 \  \ 
 \mr{L}_{h_1} \circ \mr{L}_{h_2} = \mr{L}_{h_1 h_2}, \ \  
 \text{and}   
 \ \
 \mr{R}_{h_1} \circ \mr{L}_{h_2} = \mr{L}_{h_2} \circ \mr{R}_{h_1}
 \end{align}
  for any $h_1, h_2 \in \mbG (T)$.

By considering the right $\mbG$-action  on $\mbG$ defined by $\mr{R}_{(-)}$, we often identify $\mfg$  with the space of right-invariant vector fields on $\mbG$.
That is to say, given any element $v \in \mfg = \mcT_e \mbG$, we can use right-translation to produce a  global vector field $v_\mbG$ with $v_\mbG |_h = d \mr{R}_h |_e (v)$ for every $h \in \mbG$, where $d \mr{R}_h |_e$ denotes the differential $\mcT_e \mbG \migi \mcT_h \mbG$ of $\mr{R}_h$ at $e$.

\subsection{The adjoint action}
Given a  morphism of $k$-schemes  $h : T \migi \mbG$, we shall write 
\begin{equation}  \label{QQ1}
\mr{Inn}_h : = \mr{L}_h \circ \mr{R}_{h^{-1}} \left(= \mr{R}_{h^{-1}}\circ \mr{L}_h \right) : T \times_k \mbG  \migi T \times_k \mbG,
\index{$\mr{Inn}_h$}
\end{equation}
which specifies an automorphism of the  $T$-scheme $T \times_k \mbG$.
Let  $e_T : T \migi \mbG$ denote  the identity in $\mbG (T)$, i.e., the composite $T \migi \mr{Spec}(k) \xrightarrow{e} \mbG$.
Denote by $\mr{GL}(\mfg)$ the automorphism group of the $k$-vector space $\mfg$. 
The {\bf adjoint representation}  of $\mbG$ is, by definition,   the morphism
\begin{equation} 
\label{adjrep} \mr{Ad}_\mbG :\mbG \migi \mr{GL}(\mfg) 
\index{$\mr{Ad}_\mbG$, adjoint representation}\end{equation}
  of algebraic $k$-groups which, to any 
$k$-morphism $h : T \migi \mbG$, assigns 
 the automorphism 
 $d\mr{Inn}_{h} |_{e_T}$
  of $T \times_k \mfg$ 
   defined as  the differential   of $\mr{Inn}_{h}$  at  $e_T$.
We occasionally   identify $\mr{Ad}_\mbG (h)$ with the corresponding $\mcO_T$-linear automorphism of $\mcO_T \otimes_k \mfg \left(= \Gamma_{e_T}^*(\mcT_{T \times_k \mbG/T})\right)$.
Write 
\begin{equation}\label{adjrep2}
\mr{ad} : = d\mr{Ad}_\mbG |_e : \mfg \migi \mr{End} (\mfg)  \left(= \mcT_{\mr{id}_\mfg} \mr{GL}(\mfg)\right),
\end{equation}
i.e., $\mr{ad}(v)(v') := [v, v']$ for $v$, $v' \in \mfg$.
The Lie algebra $\mfg$ has a left $\mbG$-action
 determined by $\mr{Ad}_\mbG$.
 Also, 
 the composition of  $\mr{Ad}_\mbG$ and $\mr{Ad}_{\mr{GL}(\mfg)} : \mr{GL}(\mfg) \migi \mr{GL}(\mr{End}(\mfg))$
  determines a left $\mbG$-action on $\mr{End}(\mfg)$; more concretely, it  is given by  
 $(h, a) \mapsto 
  \mr{Ad}_\mbG (h) \circ a \circ \mr{Ad}_{\mbG}(h^{-1})$
  for every $h \in \mbG$, $a \in \mr{End}(\mfg)$.
This $\mbG$-action on $\mr{End} (\mfg)$  is compatible with the adjoint $\mbG$-action on $\mfg$ via 
 $\mr{ad} : \mfg \migi \mr{End} (\mfg)$.

\subsection{The Maurer-Cartan form} \label{QR193}
In the following, let us recall the Maurer-Cartan form on the algebraic group $\mbG$  (cf. ~\cite{Sha}, Chap.\,3, \S\,1, Definition 1.3, or ~\cite{GS}, \S\,8.6).

\bde \label{QQ4}
 The {\bf (left-invariant) Maurer-Cartan form} on $\mbG$  is the global section
\begin{equation} \label{MC}
\Theta_\mbG \in \Gamma (\mbG, \Omega_{\mbG/k}\otimes_k \mfg)  \index{$\Theta_\mbG$, Maurer-Cartan form on $\mbG$} 
\end{equation}
 corresponding to the differential
\begin{align}
 d\mr{L}_{\mr{inv}} |_{\Gamma_{\mr{id}_\mbG}} : & \left(\mcT_{\mbG/k} = \right)  \Gamma_{\mr{id}_\mbG}^*(\mcT_{\mbG \times_k \mbG/\mbG}) 
 \migi \Gamma_{e_\mbG}^*(\mcT_{\mbG \times_k \mbG/\mbG})  \left(
 = \mcO_\mbG \otimes_k \mfg\right) 
\end{align}
of the left-multiplication $\mr{L}_{\mr{inv}} : \mbG \times_k \mbG \isom \mbG \times_k \mbG_m$ by 
 the inversion $\mr{inv} : \mbG \isom \mbG$ in $\mbG$,  where 
  $\mcT_{\mbG \times_k \mbG/\mbG}$ denotes the    tangent bundle of  $\mbG \times_k \mbG$ relative to  the first projection $\mbG \times_k \mbG \migi \mbG$.
(In particular, $\Theta_\mbG$ is  a $\mfg$-valued $1$-form on $\mbG$ in the terminology  below.)
 Occasionally, it will be convenient to  identify $\Theta_\mbG$ with the corresponding  $\mcO_\mbG$-linear morphism
  $\mcT_{\mbG/k} \migi \mcO_\mbG \otimes_k \mfg$.
 \ede

\begin{rema} \label{QQ19} 
 The Maurer-Cartan form $\Theta_\mbG$ can be characterized by  the condition that,
 for any $k$-morphism $h: T \migi \mbG$, the pull-back
\begin{equation} \label{MCcondition}
h^*(\Theta_\mbG)  \in  \Gamma (T, h^*(\Omega_{\mbG/k}) \otimes_k \mfg) \ \left(\cong \mr{Hom}_{\mcO_T}(h^*(\mcT_{\mbG/k}), \mcO_T\otimes \mfg)\right)
\end{equation}
corresponds to  the differential
\begin{align}
 d\mr{L}_{h^{-1}} |_{\Gamma_h}  : & \left(h^*(\mcT_{\mbG/k}) = \right)  \Gamma_h^*(\mcT_{T \times_k \mbG/T}) 
 \migi \Gamma_{e_T}^*(\mcT_{T \times_k \mbG/T})  \left(
 = \mcO_T \otimes_k \mfg\right) 
\end{align}
of $\mr{L}_{h^{-1}}$ at  $\Gamma_h \in (T \times_k \mbG) (T)$.
This implies that if $\partial \in \Gamma (\mbG, \mcT_{\mbG/k})$ is a left-invariant vector field, then the equality $\Theta_\mbG (\partial) = 1 \otimes \partial |_e$ holds. 
 \end{rema}

Let $T$ be a smooth $k$-scheme (hence, $\mcT_{T/k}$ is a vector bundle).
For a nonnegative integer $q$, a {\bf $\mfg$-valued $q$-form on  $T$}
\index{$\mfg$-valued $q$-form}
 is   a global  section of the vector bundle $(\bigwedge^q \Omega_{T/k}) \otimes_k \mfg$.

Let $q_1$, $q_2$ be nonnegative integers such that both $q_1!$ and $q_2 !$  are  invertible in $k$.
Also, let us take local sections $\delta_1 \in (\bigwedge^{q_1} \Omega_{T/k}) \otimes_k \mfg$,
$\delta_2 \in (\bigwedge^{q_2} \Omega_{T/k}) \otimes_k \mfg$ defined on  a  common open subset of $T$.
Then,  we define
\begin{align} \label{QQ20}
[\delta_1, \delta_2] \in  (\bigwedge^{q_1+q_2} \Omega_{T/k}) \otimes_k \mfg
\end{align}
to be the the $\mfg$-valued $(q_1 + q_2)$-form given by
\begin{align}
& \ [\delta_1, \delta_2] (\partial_1, \cdots, \partial_{q_1+q_2}) \\
 := & \ \frac{1}{q_1! q_2!} \cdot \sum_{\sigma \in \mfS_{q_1 + q_2}} (-1)^{\mr{sgn}(\sigma)} [\delta_1 (\partial_{\sigma (1)}, \cdots, \partial_{\sigma (q_1)}), \delta_2 (\partial_{\sigma (q_1 +1)}, \cdots, \partial_{\sigma (q_1 +q_2)})], \notag
\end{align}
where $\mfS_{q_1 + q_2}$ denotes the symmetric group of $q_1 +q_2$ letters.
For example,
for $\mfg$-valued $1$-forms $\delta_1$, $\delta_2$, 
we have
\begin{align} \label{QQ21}
[\delta_1, \delta_2] (\partial_1, \partial_2) = [\delta_1 (\partial_1), \delta_2 (\partial_2)] - [\delta_1 (\partial_2), \delta_2 (\partial_1)].
\end{align}
In particular, any $\mfg$-valued $1$-form $\delta$ satisfies
\begin{align} \label{QQ22}
[\delta, \delta] (\partial_1, \partial_2) = 2 \cdot [\delta (\partial_1), \delta (\partial_2)].
\end{align}

Finally,   
 the {\bf  exterior differential} \index{$d$, exterior differential, universal derivation} is the morphism
\begin{align} \label{QQ23}
d : (\bigwedge^q \Omega_{T/k}) \otimes_k \mfg \migi 
(\bigwedge^{q+1} \Omega_{T/k}) \otimes_k \mfg
\end{align}
($q \geq 0$)  determined  by  $d (\overline{\delta} \otimes v) = d \overline{\delta} \otimes v$
for any 
$\overline{\delta} \in \bigwedge^q \Omega_{T/k}$, $v \in \mfg$.
One verifies immediately that 
\begin{align} \label{QQ24}
d [\delta_1, \delta_2] = [d \delta_1, \delta_2] + (-1)^{q_1} \cdot [\delta_1, d \delta_2]
\end{align}
for every $\delta_i \in (\bigwedge^{q_i} \Omega_{T/k})\otimes_k \mfg$ ($i=1,2$).

\bpr[The Maurer-Cartan equation] \label{QQ16}
Suppose that  $2$ is invertible in $k$.
Then, the Maurer-Cartan form $\Theta_\mbG$ satisfies an equality
  \begin{align} \label{QQ2}
 d \Theta_\mbG + \frac{1}{2} \cdot [\Theta_\mbG,  \Theta_\mbG] = 0
 \end{align}
of elements in $\Gamma (\mbG, (\bigwedge^2 \Omega_{\mbG/k})\otimes_k \mfg)$.
\epr
\begin{proof}
Let us take  two local sections $\partial_1, \partial_2 \in \mcT_{\mbG/k}$.
Since $\Theta_\mbG$ is a vector-valued $1$-form on the  smooth $k$-scheme $\mbG$,   the equality
\begin{align} \label{QR200}
d \Theta_\mbG (\partial_1, \partial_2) = \partial_1 (\Theta_\mbG (\partial_2)) - \partial_2 (\Theta_\mbG (\partial_1)) - \Theta_\mbG ([\partial_1, \partial_2])
\end{align}
holds because of   an argument  similar to the proof of   ~\cite{Sha}, Chap.\,1, \S\,5, Lemma 5.15.
If  both  $\partial_1$ and $\partial_2$ are  left-invariant, then
$[\partial_1, \partial_2]$  is left-invariant.
Moreover, thanks to  the final comment in Remark \ref{QQ19},
we obtain the equalities
 \begin{align} \label{QG5000}
 \partial_1 (\Theta_\mbG (\partial_2)) = \partial_2 (\Theta_\mbG (\partial_1)) =0
 \end{align}
   and  obtain the following sequence of equalities:
\begin{align} \label{QG5001}
\Theta_\mbG ([\partial_1, \partial_2]) &= 1 \otimes [\partial_1, \partial_2] |_e \\
&= 1 \otimes [\partial_1 |_e, \partial_2 |_e] \notag \\
&=  [1 \otimes \partial_1 |_e,  1 \otimes \partial_2 |_e]  \notag \\
&= [\Theta_\mbG (\partial_1), \Theta_\mbG (\partial_2)]  \notag \\
& \hspace{-1mm}\stackrel{(\ref{QQ22})}{=} \frac{1}{2} \cdot [\Theta_\mbG, \Theta_\mbG] (\partial_1, \partial_2). \notag
\end{align}
By (\ref{QR200}), (\ref{QG5000}), and (\ref{QG5001}), 
 the equality 
\begin{align} \label{QQ201}
\left(d \Theta_\mbG + \frac{1}{2} \cdot [\Theta_\mbG,  \Theta_\mbG]\right) (\partial_1, \partial_2)= 0
\end{align}
holds for any left-invariant vector fields $\partial_1$, $\partial_2$.
Since the left-invariant vector fields on $\mbG$ generates $\mcT_{\mbG/k}$,
(\ref{QQ201}) implies  the desired equality  (\ref{QQ2}).
\end{proof}

\subsection{Change of structure group}

Let $\mbG'$ be another connected smooth affine algebraic group over $k$ with  Lie algebra $\mfg'$ and $w : \mbG \migi \mbG'$ a $k$-morphism of algebraic groups. 
The differential $d w : \mcT_{\mbG/k} \migi w^*(\mcT_{\mbG'/k})$ of $w$ 
restricts to
a $k$-linear morphism $dw |_{e} : \mfg \migi \mfg'$.

\vspace{3mm}
\bpr \label{QQ25}
The Maurer-Cartan form $\Theta_{\mbG'}$ on $\mbG'$ makes the following square diagram commute:
\begin{align}   \label{MCchange}
\vcenter{\xymatrix@C=46pt@R=36pt{
\mcT_{\mbG/k} \ar[r]^-{\Theta_{\mbG}} \ar[d]_-{dw}& \mcO_{\mbG} \otimes_k \mfg\ar[d]^-{\mr{id} \otimes dw |_{e}}\\
w^*(\mcT_{\mbG'/k}) \ar[r]_-{w^*(\Theta_{\mbG'})}& 
\mcO_{\mbG} \otimes_k \mfg'
}}
\end{align}
\epr
\begin{proof}
The assertion follows from the definitions of $\Theta_{\mbG}$ and $\Theta_{\mbG'}$. 
\end{proof}

\begin{exa} \label{QQ29}
Let us consider the case where
  $\mbG = \mr{GL}_n$ for a positive integer $n$.
    By identifying  $\mr{GL}_n$ with  the  open subscheme $\left\{ R \in \mfg \mfl_n \, | \, \mr{det}(R) \in k^\times \right\}$ of $\mfg \mfl_n$, we have an open immersion
    $A : \mr{GL}_n \migiincl \mfg \mfl_n$.
 The  differential of $A$ defines an $\mcO_{\mr{GL}_n}$-linear morphism
 \begin{align}
 dA : \mcT_{\mr{GL}_n} \migi \left( A^*(\mcT_{\mfg \mfl_n})  = A^*(\mcO_{\mfg \mfl_n} \otimes_k \mfg \mfl_n) \right) = \mcO_{\mr{GL}_n} \otimes_k \mfg,
 \end{align}
 where we also  use the notation $dA$ to denote the corresponding $\mfg \mfl_n$-valued $1$-form on $\mr{GL}_n$.
  Then,  for any $k$-morphism $h : T \migi \mr{GL}_n$ and  a local section $\partial \in \mcT_{T/k}$, the Maurer-Cartan form $\Theta_{\mr{GL}_n}$ on $\mr{GL}_n$ satisfies 
 \begin{align} \label{EE339}
h^* (\Theta_{\mr{GL}_n}) (\partial)
 & = d \mr{L}_{h^{-1}} |_{\Gamma_h} (\partial)     \\
  & = d \mr{L}_{h^{-1}} (d h ( \partial))    \notag \\
 & = (A \circ h)^{-1}  (d (A \circ h) (\partial)) \notag \\
 & = h^*(A^{-1}dA) (\partial). \notag
 \end{align}
This implies  an equality 
\begin{align} \label{QQ27}
\Theta_{\mr{GL}_n} = A^{-1} d A
\end{align}
of $\mfg \mfl_n$-valued $1$-forms on $\mr{GL}_n$.
  Moreover, 
 if  $\mbG$ is  an algebraic  subgroup  of $\mr{GL}_n$, then it follows from Proposition \ref{QQ25} that
 \begin{align} \label{QR220}
 \Theta_\mbG = (A^{-1} dA )|_\mbG.
 \end{align}
\end{exa}

\section{Principal $\mbG$-bundles and their associated Lie algebroids}\label{y012}

Next, we
discuss
 principal bundles on {\it log schemes}.
We refer the reader to  ~\cite{KaKa}, 
~\cite{ILL}, and ~\cite{KaFu} for basic properties on log schemes.

\subsection{Principal $\mbG$-bundles}

To begin with, we recall  the definition of a principal $\mbG$-bundle.
Let $Y$ be a $k$-scheme.

\bde
 \label{QQ3}
\begin{itemize}
\item[(i)]
 A {\bf (principal) $\mbG$-bundle}  \index{$\mcE, $(principal) $\mbG$-bundle} (in the \'{e}tale topology)
  on $Y$ is a scheme $\mcE$ over $Y$ equipped with a right $\mbG$-action $\mcE \times_k \mbG \migi \mcE$ which is \'{e}tale locally trivial in the following sense:
 there exist an \'{e}tale surjective $k$-morphism $Y' \migi Y$ and a $\mbG$-equivariant isomorphism $Y' \times_Y \mcE \isom Y' \times_k \mbG$ over $Y'$, where both $Y' \times_Y \mcE$ and $Y' \times_k \mbG$ are equipped with  a right $\mbG$-action  induced from the $\mbG$-action on the  second factor of the product.
 \item[(ii)]
 Let $\mcE$, $\mcE'$ be $\mbG$-bundles on $Y$.
Then,  an {\bf isomorphism of $\mbG$-bundles} from  $\mcE$ to $\mcE'$
 is defined as a $\mbG$-equivariant isomorphism $\mcE  \isom \mcE'$ over $Y$.
 \end{itemize}
 \ede

\begin{exa}[Trivial bundle] \label{QQ9999} 
The {\bf trivial $\mbG$-bundle} 
\index{$\mcE, $(principal) $\mbG$-bundle@--- trivial} 
on $Y$ is  the product $Y \times_k \mbG$, where $\mbG$ acts only on the second factor,   considered as a scheme over $Y$ by the first projection $Y \times_k \mbG \migi Y$.
Any $\mbG$-bundle is, by definition,  \'{e}tale locally  isomorphic to the trivial $\mbG$-bundle.

The left-translation $\mr{L}_h$  by each element $h \in \mbG (Y)$ specifies an automorphism of the trivial $\mbG$-bundle $Y \times_k \mbG$, i.e., preserves   the right $\mbG$-action.
Conversely,  each automorphism $\eta$ of  the trivial $\mbG$-bundle $Y \times_k \mbG$ can be expressed as  a   left-multiplication.
Indeed,
let us take $h$ to be the composite $Y \xrightarrow{(\mr{id}_Y, e_Y)}Y \times_k \mbG \xrightarrow{\eta} Y \times_k \mbG \xrightarrow{\mr{pr}_2} \mbG$.
The two automorphisms $\eta$, $\mr{L}_h$ of $Y \times_k \mbG$ coincide when 
restricted to the closed subscheme $\mr{Im}(\Gamma_{e_Y}) \left(\subseteq Y \times_k \mbG \right)$.
Hence, for any $y \in Y$ and $h \in \mbG$, we have 
\begin{align}
\eta (y, h) = \eta (y, e)\cdot h = \mr{L}_h (y, e) \cdot h = \mr{L}_h (y, h).
\end{align}
This implies  $\eta = \mr{L}_\eta$, as desired. 
One verifies that  such an element $h$ is uniquely determined.
 \end{exa}

\subsection{Pull-back and base-change}
 \label{QQ40} 
 Let $\mcE$ be a $\mbG$-bundle on $Y$ and  $u: Y' \migi Y$
 a scheme   over $Y$.
Then,  the fiber product
 \index{$\mcE, $(principal) $\mbG$-bundle@--- $u^*(\mcE)$, pull-back of} 
 \begin{align} \label{QR240}
 u^*(\mcE) := Y' \times_Y \mcE
 \end{align} 
forms a $\mbG$-bundle on $Y'$ by considering  a right $\mbG$-action induced from the $\mbG$-action on the second factor $\mcE$ of the product.
 We shall refer to $u^*(\mcE)$ as the {\bf pull-back} of $\mcE$ to $Y'$.

Suppose that $Y'$ comes from  base-change via a $k$-morphism $t : T' \migi T$, i.e., $Y' := T' \times_T Y$.
Then, 
 we denote the pull-back of $\mcE$ to $Y'$  by 
 \index{$\mcE, $(principal) $\mbG$-bundle@--- $t^*(\mcE)$, base-change of}
 \begin{align} \label{QR241}
 t^*(\mcE) \left(:= Y' \times_Y \mcE = T' \times_T \mcE \right)
 \end{align}
  and
 refer to it 
  as the {\bf base-change} of $\mcE$ to $T'$.

\subsection{Change of structure group} \label{QQ119}
Let $\mbG'$ be another connected smooth affine algebraic group over $k$ and $w : \mbG \migi \mbG'$ a $k$-morphism of algebraic groups.
Then, each $\mbG$-bundle $\mcE$ on $Y$ induces a $\mbG'$-bundle  via  change of structure group by $w$; we denote this  $\mbG'$-bundle  by
\index{$\mcE, $(principal) $\mbG$-bundle@--- $\mcE \times^\mbG \mbG'$, $\mcE \times^{\mbG, w }\mbG'$}
\begin{align} \label{QQ120}
\mcE \times^\mbG \mbG' \left(\text{or} \  \mcE \times^{\mbG, w }\mbG' \right).
\end{align}
To be precise, let us set $\mcE \times^\mbG \mbG'$ to be the quotient of $\mcE \times_k \mbG'$ by the  right $\mbG'$-action given by $(v, h')\cdot h := (v \cdot h, h^{-1}\cdot h')$, where $v \in \mcE$, $h \in \mbG$, and $h' \in \mbG$.
 The natural right $\mbG'$-action on $\mcE \times_k \mbG'$ induces a well-defined  right $\mbG'$-action on 
 $\mcE \times^\mbG \mbG'$, by which  $\mcE \times^\mbG \mbG'$ forms a $\mbG'$-bundle on $Y$.

\bde
\label{QQ123}
Let $\mbB$ be a subgroup of $\mbG$ and $\mcE$ a $\mbG$-bundle on $Y$.
By a {\bf $\mbB$-reduction} \index{$\mbB$-reduction} of $\mcE$, we mean  a $\mbB$-bundle $\mcE_\mbB$ on $Y$ together with an isomorphism $\mcE_B \times^\mbB \mbG \isom \mcE$ of $\mbG$-bundles.
 \ede

Let $Y$ be a $k$-scheme and $\mcE$ a $\mbG$-bundle on $Y$. 
Given  a $k$-vector space $\mfh$ equipped with a left $\mbG$-action, we  write 
\begin{align} \label{QR243}
\mfh_{\mcE} 
\end{align}
for the vector bundle on $Y$ associated to the  affine space $\mcE \times^\mbG \mfh \left(:= (\mcE \times_k \mfh) /\mbG \right)$  over $Y$, i.e., satisfying that
 $\mbV (\mfh_{\mcE}) \cong \mcE \times^\mbG \mfh$.
 If $\mfh = \mfg$, then
  we refer to
 $\mfg_\mcE$ \index{$\mfg_\mcE$, adjoint  bundle associated to $\mcE$}
 as the {\bf adjoint (vector) bundle associated to $\mcE$}.
 This is  the vector bundle corresponding to the $\mr{GL}(\mfg) \left(\cong \mr{GL}_{\mr{dim}(\mfg)} \right)$-bundle $\mcE \times^{\mbG, \mr{Ad}_\mbG} \mr{GL}(\mfg)$ and equipped with an $\mcO_Y$-Lie algebra structure induced   from the Lie algebra structure on $\mfg$.

\subsection{Logarithmic derivations} \label{QR250}
Given a log scheme (or more generally, a log stack) indicated, say,  by  \index{$Y^\mr{log}$} $Y^\mr{log}$, we shall write $Y$ for the underlying scheme (stack) of $Y^\mr{log}$, and 
$\alpha_Y : M_Y \migi \mcO_Y$ \index{$\alpha_Y$} for the morphism of sheaves of monoids defining the log structure of $Y^\mr{log}$.

We shall recall the notion of a logarithmic derivation.
Let $T^\mr{log} := (T, \alpha_T : M_T \migi \mcO_T)$, $Y^\mr{log} := (Y, \alpha_Y : M_Y \migi \mcO_Y)$ be  fine log schemes over $k$ (cf. ~\cite{KaKa}, (2.3)) and $f^\mr{log} : Y^\mr{log} \migi T^\mr{log} $ a  $k$-morphism of log schemes.
By definition,  a {\bf    logarithmic derivation}  index{logarithmic derivation} $\partial$ 
on $(\mcO_Y, M_Y)$ relative to $(\mcO_T, M_T)$  (cf. ~\cite{Og}, \S\,1, Definition 1.1.2)
is given as a pair 
\begin{equation} \label{QR249}
\partial =(D : \mcO_Y \migi \mcO_Y, \delta : M_Y \migi \mcO_Y),
\end{equation}
where $D$ denotes  a (non-logarithmic) derivation on $\mcO_Y$ over $T$ and $\delta$ denotes  a monoid homomorphism from $M_Y$ to the underlying additive monoid of $\mcO_Y$ satisfying the following two conditions:
\begin{itemize}
\item
$D (\alpha_Y (m)) = \alpha_Y (m) \cdot \delta (m)$ for any local section $m \in M_Y$;
\vspace{1.5mm}
\item
$\delta (f^*(m))=0$ for any local section
 $m \in M_T$.
\end{itemize}

Let us write 
\begin{align} \label{QR244}
\mcT_{Y^\mr{log}/T^\mr{log}}
 \index{$\mcT_{Y^\mr{log}/T^\mr{log}}$, logarithmic tangent bundle on $Y^\mr{log}/T^\mr{log}$}
\end{align}
 for the sheaf of  locally defined logarithmic derivations of $(\mcO_Y, M_Y)$ relative to $(\mcO_T, M_T)$.
 It forms an $\mcO_T$-Lie algebras on $Y$  with Lie bracket operation defined by
 \begin{align} \label{QR500}
 [(D_1, \delta_1), (D_2, \delta_2)] := ([D_1, D_2], D_1 \circ \delta_2 - D_2 \circ \delta_1)
 \end{align}
 (cf. ~\cite{Og}, \S\,1, Proposition 1.1.7).
  Also,  write 
  \begin{align} \label{QR245}
  \Omega_{Y^\mr{log}/T^\mr{log}} 
  \index{$\Omega_{Y^\mr{log}/T^\mr{log}}$, sheaf of logarithmic $1$-forms on $Y^\mr{log}/T^\mr{log}$}
   \end{align}
   for   the sheaf of logarithmic $1$-forms on $Y^\mr{log}$ relative to  $T^\mr{log}$, as defined in ~\cite{KaKa}, (1.7).
 In the usual manner,  the sheaf $\mcT_{Y^\mr{log}/T^\mr{log}}$ may be identified with the dual of $\Omega_{Y^\mr{log}/T^\mr{log}}$.
 (The formations of such sheaves may be extended naturally to the case where both $T^\mr{log}$ and $Y^\mr{log}$ are log stacks.)
If
 $Y^\mr{log}/T^\mr{log}$ is log smooth   (i.e., ``smooth" in the sense of ~\cite{KaKa}, (3.3)), then  both $\mcT_{Y^\mr{log}/T^\mr{log}}$ and $\Omega_{Y^\mr{log}/T^\mr{log}}$ are  vector bundles (cf ~\cite{KaKa}, Proposition (3.10)).
The universal logarithmic derivation on $Y^\mr{log}$ (relative to $T^\mr{log}$) will be denoted by
 \index{$d$, exterior differential, universal derivation}
\begin{align} \label{QR280}
d : \mcO_Y \migi \Omega_{Y^\mr{log}/T^\mr{log}}. 
\end{align}

\begin{exa}[Trivial log strucutre] \label{QQ102} 
 Suppose that $f^\mr{log} : Y^\mr{log} \migi T^\mr{log}$ is strict in the sense of ~\cite{ILL}, \S\,1.2.
 Then, the assignment given by $\partial  \mapsto D$ for every $\partial$ as in  (\ref{QR249})
   defines an isomorphism $\mcT_{Y^\mr{log}/T^\mr{log}} \isom \mcT_{Y/T}$.
 Indeed, if $D$ is   a local section  of $\mcT_{Y/T}$,
 then the morphism $\mcO_Y^\times \migi \mcO_Y$ determined by  $m \mapsto \frac{D (m)}{m}$ extends to a well-defined monoid homomorphism $\delta_D: M_Y \migi \mcO_Y$.
 One verifies that  the pair $(D, \delta_D)$ specifies a section of $\mcT_{Y^\mr{log}/T^\mr{log}}$ and 
 the assignment $D \mapsto (D, \delta_D)$ gives the inverse to 
 the assignment $\partial \mapsto D$.
  \end{exa}

\begin{exa}[Normal crossing divisor]\label{QQ82}
Let  $T$ be  a $k$-scheme,  $Y$ a smooth $T$-scheme,  and $D$ the relative reduced divisor on $Y$ (relative to $T$) with normal crossings.
If $j : Y\setminus D \migiincl Y$  denotes the open immersion, then
the inclusion $M_Y = \mcO_Y \cap j_*(\mcO^\times_{Y \setminus D}) \migiincl \mcO_Y$ defines  a fine log structure on $Y$.
Choose   a system of local coordinates   $(x_1, \cdots, x_n)$  on $Y$ (relative to $T$) such that $D$ is defined by the equation $\prod_{i=1}^r x_i^{a_i} =0$ (where $0\leq r \leq n$ and $a_1, \cdots, a_r \geq 1$).
Then, the sheaf of monoids $M_Y$ can be expressed locally as $M_Y = \mcO_Y^\times \cdot  \prod_{i=1}^r x_i^{\mbZ_{\geq 0}}$.

For each $i=1, \cdots, n$,  let  $\partial_i$ denote the non-logarithmic derivation $\frac{\partial}{\partial x_i}$ on $\mcO_Y$.
In addition, for each $i =1, \cdots, r$,
denote by $\partial^\mr{log}_i$ the logarithmic derivation
$(x_i \cdot \partial_i, \delta_i)$, where  $\delta_i$ is the monoid homomorphism $\mcO_Y^\times \cdot  \prod_{i=1}^r x_i^{\mbZ_{\geq 0}} \migi \mcO_Y$ given by $u \cdot \prod_{i=1}^r x_i^{l_i} \mapsto x_i \cdot \frac{\partial_i u}{u} + l_i$.
Then, $\mcT_{Y^\mr{log}/T}$ decomposes locally into the direct sum
\begin{align}
\mcT_{Y^\mr{log}/T} = \left(\bigoplus_{i=1}^r \mcO_Y \cdot \partial^\mr{log}_{i}\right) \oplus \left(\bigoplus_{i=r+1}^n  \mcO_Y \cdot \partial_{i} \right).
\end{align}
  \end{exa}

\subsection{Lie algebroids} \label{QR300}
 We shall discuss  the logarithmic version of  Lie algebroids  (cf. ~\cite{Lan}, \S\,2.1, Definition 2.1 for  the non-logarithmic version) and the Lie algebroid associated to a $\mbG$-bundle.
 
\bde \label{QS39}
A {\bf Lie  algebroid} \index{Lie algebroid} on $Y^\mr{log}/T^\mr{log}$ is a pair 
 \begin{align}
 (\mcV, \alpha)
 \end{align}
  consisting of an  $\mcO_T$-Lie algebras $\mcV$ on $Y$ and an $\mcO_Y$-linear morphism $\alpha : \mcV \migi \mcT_{Y^\mr{log}/T^\mr{log}}$, which commutes  with the  Lie bracket operations and satisfies the  equality
\begin{align} \label{QG430}
[\partial_1, a \cdot \partial_2] = \alpha (\partial_1) (a) \cdot \partial_2 + a \cdot [\partial_1, \partial_2]
\end{align}
for any local sections $a \in \mcO_Y$ and $\partial_1, \partial_2 \in \mcV$.
The morphism $\alpha$ is called the {\bf anchor (map)} 
\index{Lie algebroid@--- anchor (map) of}
 of $\mcV$. 
If there is no fear of confusion, we only use the symbol $\mcV$ to indicate the  Lie algebroid  $(\mcV,  \alpha)$.
 \ede

In what follows, let us construct a Lie algebroid associated to a $\mbG$-bundle on $Y$ in the context of logarithmic geometry.
Let $\mcE$ be a $\mbG$-bundle on $Y$ and denote by $\pi : \mcE \migi Y$ its structure morphism.
By pulling-back the log structure of $Y^\mr{log}$  via $\pi$, one may obtain a log structure on $\mcE$; we denote  the resulting log scheme by   $\mcE^\mr{log}$.
The right $\mbG$-action on $\mcE$ induces a right $\mbG$-action on  the direct image   $\pi_*(\mcT_{\mcE^\mr{log}/T^\mr{log}})$ (resp.,  $\pi_*(\mcT_{\mcE^\mr{log}/Y^\mr{log}})$) of  the sheaf   of logarithmic derivations $\mcT_{\mcE^\mr{log}/T^\mr{log}}$ (resp., $\mcT_{\mcE^\mr{log}/Y^\mr{log}}$).
 We shall set
 \begin{equation}
 \widetilde{\mcT}_{\mcE^\mr{log}/T^\mr{log}} := \pi_*(\mcT_{\mcE^\mr{log}/T^\mr{log}})^\mbG, \hspace{10mm} \widetilde{\mcT}_{\mcE^\mr{log}/Y^\mr{log}}   := \pi_*(\mcT_{\mcE^\mr{log}/Y^\mr{log}})^\mbG,
 \index{$\widetilde{\mcT}_{\mcE^\mr{log}/T^\mr{log}}$}
 \end{equation}
 i.e., $\widetilde{\mcT}_{\mcE^\mr{log}/T^\mr{log}}$ and $\widetilde{\mcT}_{\mcE^\mr{log}/Y^\mr{log}}$ are   the subsheaves of $\mbG$-invariant sections  of $\pi_*(\mcT_{\mcE^\mr{log}/T^\mr{log}})$ and $\pi_*(\mcT_{\mcE^\mr{log}/Y^\mr{log}})$ respectively.
Under the natural identification of 
$\mfg$ with the space of right-invariant vector fields on $\mbG$, 
 $\widetilde{\mcT}_{\mcE^\mr{log}/Y^\mr{log}}$ may be  identified with  the adjoint vector bundle $\mfg_{\mcE}$  associated to $\mcE$.
Differentiating $\pi$ yields a sequence of $\mcO_{\mcE}$-modules
\begin{align} \label{QR301}
0 \longmigi \mcT_{\mcE^\mr{log}/Y^\mr{log}} \longmigi \mcT_{\mcE^\mr{log}/T^\mr{log}} \longmigi \pi^*(\mcT_{Y^\mr{log}/T^\mr{log}}) \longmigi 0
\end{align}
(i.e., the dual of the  sequence displayed  in ~\cite{KaKa}, Proposition (3.12)),   which is exact because of the local triviality of the $\mbG$-bundle $\mcE$.
By applying the functor $\pi_*(-)^\mbG$ to this sequence, we obtain
a short exact sequence of $\mcO_Y$-modules
\begin{equation}  0 \longmigi  \mfg_\mcE \longmigi  \widetilde{\mcT}_{\mcE^\mr{log}/T^\mr{log}} \stackrel{d^\mr{log}_\mcE}{\longmigi}  \mcT_{Y^\mr{log}/T^\mr{log}} \longmigi 0.
\label{Ex0}
\index{$d^\mr{log}_\mcE$}
\end{equation} 
We shall regard $\mfg_\mcE$ as an $\mcO_Y$-submodule of  $\widetilde{\mcT}_{\mcE^\mr{log}/T^\mr{log}}$ via the second arrow $\mfg_\mcE \migiincl  \widetilde{\mcT}_{\mcE^\mr{log}/T^\mr{log}}$.
The sheaf $\widetilde{\mcT}_{\mcE^\mr{log}/T^\mr{log}}$ admits a structure of Lie bracket $[-,-]$ arising from that of $\mcT_{\mcE^\mr{log}/T^\mr{log}}$, and $d_\mcE^\mr{log}$ commutes with the respective Lie bracket operations  on
$\widetilde{\mcT}_{\mcE^\mr{log}/T^\mr{log}}$ and $\mcT_{Y^\mr{log}/T^\mr{log}}$.
Moreover,  by means of the  local description of $[-,-]$ (cf. Example \ref{QQ131} below)  the equality 
\begin{align} \label{QS41}
[\partial_1, a \cdot \partial_2] = d_\mcE^\mr{log} (\partial_1)(a) \cdot \partial_2 + a \cdot [\partial_1, \partial_2]
\end{align}
can be verified for any local sections $a \in \mcO_Y$ and $\partial_1, \partial_2 \in \widetilde{\mcT}_{\mcE^\mr{log}/T^\mr{log}}$.
That is to say, 
the pair 
\begin{align} \label{QQ170}
( \widetilde{\mcT}_{\mcE^\mr{log}/T^\mr{log}},  d^\mr{log}_\mcE)
\end{align}
 specifies   a Lie algebroid on $Y^\mr{log}/T^\mr{log}$.
We   call it   the {\bf  Lie algebroid (on $Y^\mr{log}/T^\mr{log}$) associated to  $\mcE$}.

\begin{exa}[Trivial bundle] \label{QQ131} 
Let us consider the Lie  algebroid associated to the trivial $\mbG$-bundle.
Denote by $\mr{pr}_1$ (resp., $\mr{pr}_2$) the first projection $Y \times_k \mbG \migi Y$ (resp., the second projection $Y \times_k \mbG \migi \mbG$).
In particular, $\mr{pr}_1$ is nothing but the structure morphism of the trivial $\mbG$-bundle $Y \times_k \mbG$ on $Y$.
The direct sum of the differentials of $\mr{pr}_1$ and $\mr{pr}_2$ gives an isomorphism
\begin{align} \label{QQ13}
\left(\mcT_{(Y \times_k \mbG)^\mr{log}/T^\mr{log}} =\right)\mcT_{Y^\mr{log} \times_k \mbG/T^\mr{log}} \isom \mr{pr}_1^*(\mcT_{Y^\mr{log}/T^\mr{log}})\oplus \mr{pr}_2^*(\mcT_{\mbG/k}).
\end{align}
By applying the functor $\mr{pr}_{1*}(-)^\mbG$ to this isomorphism, we obtain 
a decomposition 
\begin{align} \label{QQ136}
\tau : \widetilde{\mcT}_{Y^\mr{log}\times_k \mbG/T^\mr{log}} \isom \mcT_{Y^\mr{log}/T^\mr{log}} \oplus \mcO_Y \otimes_k \mfg
\index{$\tau$}
\end{align}
of $\widetilde{\mcT}_{Y^\mr{log} \times_k \mbG/T^\mr{log}}$, where the second factor  $\mcO_Y \otimes_k \mfg$ of  the codomain is naturally isomorphic to $\mfg_{Y \times_k \mbG}$.
The morphism $d_{Y \times_k \mbG}^\mr{log}$ coincides, via $\tau$, with the first projection $\mcT_{Y^\mr{log}/T^\mr{log}} \oplus \mcO_Y \otimes_k \mfg \migi \mcT_{Y^\mr{log}/T^\mr{log}}$.
Under the decomposition $\tau$, the Lie bracket operation on $\widetilde{\mcT}_{Y^\mr{log}\times_k \mbG/T^\mr{log}}$ is given by 
\begin{align} \label{QR922}
& [(\partial_1, a_1 \otimes v_1), (\partial_2, a_2 \otimes v_2)]  \\
= & ([\partial_1, \partial_2], (a_1 \cdot a_2) \otimes [v_1, v_2] + \partial_1 (a_2)\otimes v_2-\partial_2(a_1)\otimes v_1)\notag  
\end{align}
for any local sections $\partial_i \in \mcT_{Y^\mr{log}/T^\mr{log}}$, $a_i \in \mcO_Y$, and any  $v_i \in \mfg$ ($i=1,2$).
By means of this description,    (\ref{QS41}) can be immediately verified  for the trivial $\mbG$-bundle. 
\end{exa}

\subsection{Morphisms of Lie algebroids}
Next, we define  a morphism of Lie algebroids  and exhibit  two  examples  arising  from $\mbG$-bundles.

\bde \label{QG200}
Let $\mcV$, $\mcV'$
 be Lie algebroids on $Y^\mr{log}/T^\mr{log}$.
 Then, a {\bf morphism of Lie algebroids} from $\mcV$ to $\mcV'$
 is an $\mcO_Y$-linear morphism $\mcV \migi \mcV'$ that commutes with the respective Lie bracket operations  and  with the respective anchor maps.
\ede

\begin{exa} \label{QG401}
 An automorphism of $\mbG$-bundles induces  an automorphism between their associated  Lie algebroids. 
Indeed, let $\pi: \mcE \migi Y$ and  $\pi' : \mcE' \migi Y$ be $\mbG$-bundles  on $Y$ and $\eta : \mcE \isom\mcE'$  an isomorphism of $\mbG$-bundles.
The differential of $\eta$ gives a $\mbG$-equivariant   isomorphism
$d \eta : \pi_*(\mcT_{\mcE^{\mr{log}}/T^\mr{log}}) \isom \pi'_*(\mcT_{\mcE'^{^\mr{log}}/T^\mr{log}})$ of $\mcO_Y$-modules.
Then, the restriction 
 \begin{align} \label{QQ10}
 d\eta^\mbG :\widetilde{\mcT}_{\mcE^\mr{log}/T^\mr{log}} \isom  \widetilde{\mcT}_{\mcE'^{\mr{log}}/T^\mr{log}}\index{$ d\eta^\mbG$}
 \end{align} 
of this isomorphism specifies an isomorphism of Lie algebroids.
\end{exa}

\begin{exa} \label{QG402}
The change of structure group  for a  $\mbG$-bundle induces a morphism  of  Lie algebroids.
In fact,  let $\mcE$ be a $\mbG$-bundle on $Y$, and 
let $\mbG'$ and $w$ be as in \S\,\ref{QQ119}.
Then, one verifies that the differential of the natural morphism $\mcE \migi \mcE \times^\mbG \mbG'$ (i.e., the morphism given by $v \mapsto (v, e)$ for every  $v \in \mcE$) yields
a morphism of Lie algebroids
\begin{align} \label{QQ810}
(dw)_\mcE : \widetilde{\mcT}_{\mcE^\mr{log}/T^\mr{log}} \migi  \widetilde{\mcT}_{(\mcE \times^\mbG \mbG')^\mr{log}/T^\mr{log}}.
\index{$(dw)_\mcE$}
\end{align}
\end{exa}

\begin{rema} \label{QQ130} 
The Lie algebroid $\widetilde{\mcT}_{\mcE^\mr{log}/T^\mr{log}}$ can be described in a different way after choosing a section of $\mcE$.
To see this, suppose that there is 
  a global section $\sigma : Y \migi \mcE$ of the $\mbG$-bundle $\pi : \mcE \migi Y$.
  This section
   induces an  isomorphism  of $\mbG$-bundles
$\eta_\sigma : Y \times_k \mbG \isom \mcE$
 given by $(y, g) \mapsto \sigma (y) \cdot g$ for $y \in Y$ and  $g \in \mbG$.
Hence, we have the
following composite isomorphism:
\begin{align} \label{QQ133}
\mcT_{Y^\mr{log}/T^\mr{log}} \oplus \mcO_Y \otimes \mfg 
&\isom \Gamma_{e_Y}^*(\mr{pr}_1^*(\mcT_{T^\mr{log}/T^\mr{log}}))\oplus \Gamma_{e_Y}^*(\mr{pr}_2^*(\mcT_{\mbG/k})) \\
& \isom \Gamma_{e_Y}^*(\mcT_{Y^\mr{log}\times_k \mbG/T^\mr{log}}) \notag \\
& \isom \Gamma_{e_Y}^*(\eta_\sigma^*(\mcT_{\mcE^\mr{log}/T^\mr{log}})) \notag \\
& \isom \sigma^*(\mcT_{\mcE^\mr{log}/T^\mr{log}}), \notag
\end{align}
where the second arrow  arises from (\ref{QQ13}).
By composing this  isomorphism with  $\tau \circ d (\eta_\sigma^{-1})^\mbG \circ$ (cf. (\ref{QQ136}) for the definition of $\tau$), we obtain an isomorphism of $\mcO_Y$-modules
\begin{align} \label{QQ139}
\widetilde{\mcT}_{\mcE^\mr{log}/T^\mr{log}} \isom \sigma^*(\mcT_{\mcE^\mr{log}/T^\mr{log}}).
\end{align}
This description will be used in the proof of Lemma \ref{QR349}.
 \end{rema}

\section{Logarithmic connections and curvature} \label{QQ129}

\subsection{Logarithmic connections} \label{QG405}
We shall introduce   
 an $\hslash$-twisted  version of a logarithmic connection defined  on a principal bundle.
Let $Y^\mr{log}/T^\mr{log}$ be as above.
Also, let $\mcE$ be a $\mbG$-bundle on $Y$ and    $\hslash$  an element of $\Gamma(T, \mcO_T)$.

\bde \label{y021}
An {\bf $\hslash$-$T^\mr{log}$-connection} \index{$\nabla$, $\hslash$-$T^\mr{log}$-connection}
on  $\mcE$
  is an $\mcO_Y$-linear morphism
\begin{align} \label{QQ75}
\nabla : \mcT_{Y^\mr{log}/T^\mr{log}} \migi \widetilde{\mcT}_{\mcE^\mr{log}/T^\mr{log}}
\end{align}
 satisfying 
$d^\mr{log}_\mcE \circ \nabla = \hslash \cdot \mr{id}_{\mcT_{Y^\mr{log}/T^\mr{log}}}$.
For simplicity, a $1$-$T^\mr{log}$-connection will be referred to as a {\bf $T^\mr{log}$-connection}.
  \ede

\begin{exa}[The case of $\hslash =0$] \label{QQ8} 
Note that 
 giving a $0$-$T^\mr{log}$-connection on $\mcE$ is equivalent to giving an $\mcO_Y$-linear morphism
$\mcT_{Y^\mr{log}/^\mr{log}} \migi \mfg_{\mcE}$.
Indeed, the image of a $0$-$T^\mr{log}$-connection $\nabla : \mcT_{Y^\mr{log}/T^\mr{log}} \migi \widetilde{\mcT}_{\mcE^\mr{log}/T^\mr{log}}$ is contained in $\mr{Ker}(d_\mcE^\mr{log}) = \mfg_\mcE$ because of the equality $d_\mcE^\mr{log} \circ \nabla = 0$.
 \end{exa}

\begin{rema}[The space  of connections]
 \label{QQ0995}
Let $\nabla$ be an $\hslash$-$T^\mr{log}$-connection on $\mcE$.
Then, for an $\mcO_Y$-linear morphism $\nabla' : \mcT_{Y^\mr{log}/T^\mr{log}} \migi \widetilde{\mcT}_{\mcE^\mr{log}/T^\mr{log}}$,
it defines an $\hslash$-$T^\mr{log}$-connection on $\mcE$ if and only if
 $\nabla' = \nabla + R$ for some $R \in \mr{Hom}_{\mcO_Y}(\mcT_{Y^\mr{log}/T^\mr{log}}, \mfg_\mcE)$ $\left(\cong \Gamma (Y, \Omega_{Y^\mr{log}/T^\mr{log}} \otimes \mfg_\mcE) \right)$.
It follows that the set  of $\hslash$-$T^\mr{log}$-connections on $\mcE$ forms, if it is nonempty,  an affine space  modeled on the $k$-vector space $\Gamma (Y, \Omega_{Y^\mr{log}/T^\mr{log}} \otimes \mfg_\mcE)$.
\end{rema}

\subsection{The trivial bundle} \label{QR304}
We here consider  the case of the trivial  $\mbG$-bundle $\mcE = Y \times_k \mbG$.  
Given an $\hslash$-$T^\mr{log}$-connection $\nabla$ on $Y \times_k \mbG$,
we shall write $\nabla^\mfg$ for the composite
\index{$\nabla$, $\hslash$-$T^\mr{log}$-connection@--- $\nabla^\mfg$}
\begin{align} \label{QQ12}
\nabla^\mfg : \mcT_{Y^\mr{log}/T^\mr{log}} \xrightarrow{\nabla}  \widetilde{\mcT}_{Y^\mr{log}\times_k \mbG/T^\mr{log}} \xrightarrow{\tau}
\mcT_{Y^\mr{log}/T^\mr{log}} \oplus \mcO_Y \otimes_k \mfg \xrightarrow{\mr{pr}_2}
\mcO_Y \otimes_k \mfg.
\end{align}
This morphism   characterizes the connection $\nabla$
because the 
 equality
$\nabla = (\hslash \cdot \mr{id}_{\mcT_{Y^\mr{log}/T^\mr{log}}}, \nabla^\mfg)$  holds under
the identification $\widetilde{\mcT}_{Y^\mr{log}\times_k \mbG/T^\mr{log}} = \mcT_{Y^\mr{log}/T^\mr{log}}\oplus \mcO_Y \otimes_k \mfg$ given by
  $\tau$.
Then, 
 $\nabla \mapsto \nabla^\mfg$ defines a bijective correspondence between the set of $\hslash$-$T^\mr{log}$-connections on $Y \times_k \mbG$ and the set of $\mcO_Y$-linear morphisms $\mcT_{Y^\mr{log}/T^\mr{log}} \migi \mcO_Y \otimes_k \mfg$.
 In particular, if $\mbG = \mbG_m$ (hence, $\mfg = k$),
then any $\hslash$-$T^\mr{log}$-connection on the trivial $\mbG_m$-bundle $Y \times_k \mbG_m$ can be given by 
a global section of $\Omega_{Y^\mr{log}/T^\mr{log}}$
$\left(\cong \mcH om_{\mcO_Y}(\mcT_{Y^\mr{log}/T^\mr{log}}, \mcO_Y) \right)$.

Also,  there is  a trivial example of an $\hslash$-$T^\mr{log}$-connection
\index{$\nabla$, $\hslash$-$T^\mr{log}$-connection@--- $\nabla^\mr{triv}$, trivial}
\begin{align} \label{QR340}
\nabla^\mr{triv} : \mcT_{Y^\mr{log}/T^\mr{log}}\migi   \widetilde{\mcT}_{Y^\mr{log} \times_k \mbG/T^\mr{log}}
\end{align}
corresponding to the zero map $\mcT_{Y^\mr{log}/T^\mr{log}} \migi \mcO_Y \otimes_k \mfg$.
We  refer to $\nabla^\mr{triv}$ as the {\bf  trivial connection} on $Y \times_k \mbG$.


\subsection{Pull-back}\label{QQ41} 
Let $\nabla$ be an $\hslash$-$T^\mr{log}$-connection on $\mcE$.
Also, let $Y'^{\mr{log}}$ be a fine log scheme equipped with 
a {\it  log \'{e}tale} morphism $u^\mr{log} : Y'^{\mr{log}} \migi Y^\mr{log}$.
Then, we can  obtain an $\hslash$-$T^\mr{log}$-connection on  the pull-back $\mcE' := u^*(\mcE)$ (cf. (\ref{QR240})) arising from  $\nabla$.
Indeed, 
 because of the log \'{e}taleness of $u^\mr{log}$, the differential of the natural  projection $u_\mcE : \mcE' \left(=Y' \times_Y \mcE \right) \migi \mcE$ yields the following  isomorphism of short exact sequences:
\begin{align} \label{QR310}
\vcenter{\xymatrix@C=36pt@R=36pt{
0 \ar[r] 
 &\mfg_{\mcE'} \ar[r] \ar[d]^-{\wr}
 &\widetilde{\mcT}_{\mcE'^{\mr{log}}/T^\mr{log}} \ar[r]^-{d_{\mcE'}^\mr{log}} \ar[d]^-{\wr}
 &\mcT_{Y'^{\mr{log}}/T^\mr{log}} \ar[d]^-{\wr} \ar[r] & 0 \\
0\ar[r]
&u^*(\mfg_{\mcE}) \ar[r]
&u^*(\widetilde{\mcT}_{\mcE^\mr{log}/T^\mr{log}})  \ar[r]_-{u^*(d_\mcE^\mr{log})}
&u^*(\mcT_{Y^{\mr{log}}/T^\mr{log}}) \ar[r] & 0.
}}
\end{align}
Then, the pull-back $u^*(\mcT_{Y^{\mr{log}}/T^\mr{log}}) \migi u^*(\widetilde{\mcT}_{\mcE^\mr{log}/T^\mr{log}})$ of $\nabla$ by $u$ specifies, via the above isomorphism of sequences, an $\mcO_{Y'}$-linear morphism
$\mcT_{Y'^{\mr{log}}/T^\mr{log}} \migi \widetilde{\mcT}_{\mcE'^{\mr{log}}/T^\mr{log}}$ forming an
 $\hslash$-$T^\mr{log}$-connection on $\mcE'$.
 We 
 shall write 
 \index{$\nabla$, $\hslash$-$T^\mr{log}$-connection@--- $u^*(\nabla)$, pull-back of $\nabla$}
  \begin{align} \label{QR312}
 u^*(\nabla) \ \left(\text{or} \ u_\mcE^*(\nabla)\right)
 \end{align}
for  this connection and refer to it as the {\bf pull-back} of $\nabla$ to  $Y'^{\mr{log}}$.
Conversely, one verifies that the formation of  an $\hslash$-$T^\mr{log}$-connection on $Y^\mr{log}/T^\mr{log}$ has  descent with respect to the \'{e}tale topology on $Y$.

\subsection{Base-change}\label{QQ98} 
Let $T'^{\mr{log}}$ be a fine log scheme over $k$ equipped with a $k$-morphism $t^\mr{log} : T'^{\mr{log}} \migi T^\mr{log}$, and write $Y'^{\mr{log}}:= T'^{\mr{log}} \times_{T^\mr{log}} Y^\mr{log}$.
Then,  
we can construct  the base-change of an $\hslash$-$T^\mr{log}$-connection  over   $T'^{\mr{log}}$.
Indeed,   the differential of the natural projection $t^*(\mcE) \migi \mcE$
 yields an isomorphism 
\begin{align} \label{QR337}
\widetilde{\mcT}_{t^*(\mcE)^\mr{log}/T'^{\mr{log}}}\isom t^*(\widetilde{\mcT}_{\mcE^\mr{log}/T^\mr{log}})
\end{align}
 compatible with the respective surjections onto $\mcT_{Y'^{\mr{log}}/T'^{\mr{log}}} \left(\cong t^*(\mcT_{Y^\mr{log}/T^\mr{log}}) \right)$, where we  use the notation $t^*(-)$ to denote pull-back   by  the projection $Y' \migi Y$ (cf. (\ref{QR241})).
Hence, for an $\hslash$-$T^\mr{log}$-connection $\nabla$ on $\mcE$, 
its pull-back $t^*(\mcT_{Y^\mr{log}/T^\mr{log}}) \migi t^*(\widetilde{\mcT}_{\mcE^\mr{log}/T^\mr{log}})$ by  the projection  $Y' \migi Y$ specifies, via (\ref{QR337}),
 a $t^*(\hslash)$-$T'^{\mr{log}}$-connection
 \index{$\nabla$, $\hslash$-$T^\mr{log}$-connection@--- $t^*(\nabla)$, base-change}
\begin{align} \label{QQ100}
t^*(\nabla) : \mcT_{Y'^{\mr{log}}/T'^{\mr{log}}} \migi \widetilde{\mcT}_{t(\mcE)^\mr{log}/T'^{\mr{log}}} 
\end{align}
on  $t^*(\mcE)$,  where $t^*(\hslash)$ denotes the image of $\hslash$ in $\Gamma (T', \mcO_{T'})$.


\subsection{Gauge transformation} \label{QR331}
In what follows, we define the notion of  gauge transformation  and examine  how the local description of an $\hslash$-$T^\mr{log}$-connection   changes after applying  a gauge transformation. 

Let $\pi': \mcE' \migi Y$ be another $\mbG$-bundle on $Y$ and $\eta : \mcE \isom \mcE'$ an isomorphism of $\mbG$-bundles.
Then, for an $\hslash$-$T^\mr{log}$-connection
$\nabla$
  on $\mcE$,  the composite
\begin{equation} \label{comps}
\eta_*(\nabla) := d\eta^\mbG \circ \nabla : \mcT_{Y^\mr{log}/T^\mr{log}} \migi  \widetilde{\mcT}_{\mcE'^{\mr{log}}/T^\mr{log}} 
\end{equation}
(cf. (\ref{QQ10}) for the definition of $d\eta^\mbG$)
 forms an $\hslash$-$T^\mr{log}$-connection on $\mcE'$.

\bde
 \label{QQ90}
Suppose that $\mcE' = \mcE$, i.e., $\eta$ is an automorphism of $\mcE$.
Then, 
we  shall refer to $\eta_*(\nabla)$ as the {\bf gauge transformation}  of $\nabla$ by $\eta$. 
\index{$\nabla$, $\hslash$-$T^\mr{log}$-connection@--- $\eta_* (\nabla)$, gauge transformation of $\nabla$ by $\eta$}
 \ede


Let us observe the following lemma.
\ble \label{QR349}
We shall identify $\widetilde{\mcT}_{Y^\mr{log}\times_k \mbG/T^\mr{log}}$ 
with $\mcT_{Y^\mr{log}/T^\mr{log}} \oplus \mcO_Y \otimes_k \mfg$ by means of $\tau$  (cf. (\ref{QQ136})).
Then, for 
a morphism of $k$-schemes  $h: Y \migi \mbG$, 
the automorphism  $d \mr{L}_h^\mbG$ (cf. (\ref{QQ10})) of 
$\mcT_{Y^\mr{log}/T^\mr{log}} \oplus \mcO_Y \otimes_k \mfg$
 can be expressed  as 
\begin{equation}
d \mr{L}_h^\mbG = (\mr{id}_{\mcT_{Y^\mr{log}/T^\mr{log}}}, (h^{-1})^*(\Theta_\mbG) + \mr{Ad}_\mbG (h)), 
\end{equation}
where
$\Theta_\mbG$ and  $\mr{Ad}_\mbG (h)$ are  considered  as   $\mcO_Y$-linear  morphisms $\mcT_{Y^\mr{log}/T^\mr{log}} \migi \mcO_Y \otimes_k \mfg$ and  $\mcO_Y \otimes_k \mfg \isom \mcO_Y \otimes_k \mfg$ respectively.
\ele
\begin{proof}
The automorphism $d \mr{L}_h^\mbG$
 is compatible with the  surjection
$d^\mr{log}_{Y \times_k \mbG}$, i.e., the first projection 
$\mcT_{Y^\mr{log}/T^\mr{log}} \oplus \mcO_Y \otimes_k \mfg \migi \mcT_{Y^\mr{log}/T^\mr{log}}$.
This implies that there exist
$\mcO_Y$-linear morphisms
\begin{equation}
D^1 : \mcT_{Y^\mr{log}/T^\mr{log}} \migi \mcO_Y \otimes_k \mfg \ \  \text{and} \ \ D^2 : \mcO_Y \otimes_k \mfg \migi \mcO_Y \otimes_k \mfg
\end{equation}
satisfying   $d \mr{L}_h^\mbG = (\mr{id}_{\mcT_{Y^\mr{log}/T^\mr{log}}}, D^1 + D^2)$.
In what follows, let us find out what  $D_1$ and $D_2$ are.

First, observe that, 
 for any local section $\partial \in h^*(\mcT_{\mbG/k})$,  the following sequence of equalities holds:
\begin{align}
D^1 (\partial) & =  d\mr{L}_h^\mbG ((\partial, 0)) 
   = d \mr{L}_{h} |_{\Gamma_{h^{-1}}} (\partial) 
    \stackrel{}{=} (h^{-1})^*(\Theta_\mbG) (\partial),
    \end{align}
where the second equality follows from  (\ref{QQ139}) and the last equality follows from Remark \ref{QQ19}.
This implies    $D^1 = (h^{-1})^*(\Theta_\mbG)$.

On the other hand, recall that each element of $\mfg$,
 considered as a vector field on $\mbG$, is invariant under the right-translation by any element of $\mbG$.
Hence, the automorphism of $\mr{pr}_{1*}(\mcT_{Y^\mr{log} \times_k \mbG/T^\mr{log}})$ induced by  the differential of $\mr{Inn}_{h} \left(= \mr{L}_{h}\circ \mr{R}_{h^{-1}} \right)$
coincides with that of $\mr{L}_{h}$.
It follows that  $d \mr{L}_h^\mbG$ restricts to 
the differential of $\mr{Inn}_h$ at $e_Y$, which coincides with  $\mr{Ad}_\mbG(h)$.
That is to say,  we have $D^2 = \mr{Ad}_\mbG(h)$, completing the proof.
\end{proof}

The following proposition follows directly from the above lemma.

\vspace{3mm}
\bpr \label{y024}
Let $h$ be a morphism of $k$-schemes  $Y \migi \mbG$.
Then, 
 for an $\hslash$-$T^\mr{log}$-connection $\nabla$ on the trivial $\mbG$-bundle $Y \times_k \mbG$,  the gauge transformation $\mr{L}_{h*}(\nabla)$ 
  of $\nabla$ by $\mr{L}_h$ satisfies  the following equality of morphisms  $\mcT_{Y^\mr{log}/T^\mr{log}} \migi \mcO_Y \otimes_k \mfg$:
\begin{equation} \label{QQ14}
\mr{L}_{h*}(\nabla)^\mfg =  \hslash \cdot (h^{-1})^*(\Theta_\mbG)  + \mr{Ad}_\mbG (h)\circ \nabla^\mfg.
\end{equation}
(cf. (\ref{QQ12}) for the definition of $\nabla^\mfg$).
\epr

\subsection{Curvature}\label{QG420}

We shall  define the notion of curvature.
Roughly speaking,  it measures the failure of an $\hslash$-$T^\mr{log}$-connection  to commute with the Lie bracket operation $\hslash \cdot [-, -]$ on $\mcT_{Y^\mr{log}/T^\mr{log}}$ and the usual Lie bracket operation $[-, -]$ on   $\widetilde{\mcT}_{\mcE^\mr{log}/T^\mr{log}}$.
 Let  $\nabla :  \mcT_{Y^\mr{log}/T^\mr{log}} \migi \widetilde{\mcT}_{\mcE^\mr{log}/T^\mr{log}}$ be  an $\hslash$-$T^\mr{log}$-connection on $\mcE$.
For 
 local sections $\partial_1$, $\partial_2 \in \mcT_{Y^\mr{log}/T^\mr{log}}$,
 we write 
  \begin{equation} \label{QR601}
 \psi^{\nabla} (\partial_1 \wedge \partial_2)   := [\nabla (\partial_1), \nabla (\partial_2)] -\hslash \cdot \nabla ([\partial_1, \partial_2]).  \end{equation}
If $a$ is a local section of $\mcO_Y$, then
the following equalities hold:
\begin{align}
\psi^{\nabla} (\partial_1 \wedge \partial_2) &= -\psi^{\nabla} (\partial_2 \wedge \partial_1), \\
\psi^{\nabla} ((a \cdot \partial_1) \wedge \partial_2) &= \psi^{\nabla} (\partial_1 \wedge (a \cdot \partial_2)) = a \cdot \psi^{\nabla} (\partial_1 \wedge \partial_2), \notag
\end{align}
where the lower sequence of equalities follows from (\ref{QG430}).
So  the assignment  $\partial_1 \wedge \partial_2 \mapsto  \psi^{\nabla} (\partial_1 \wedge \partial_2)$ defines  an {\it $\mcO_Y$-linear} morphism 
 \begin{equation} \label{QR620}
  \psi^{\nabla} : \bigwedge^2  \mcT_{Y^\mr{log}/T^\mr{log}}  \migi  \widetilde{\mcT}_{\mcE^\mr{log}/T^\mr{log}}. 
  \end{equation}

\bde \label{QQ880}
We shall refer to $\psi^{\nabla}$  as  the {\bf curvature}  of $\nabla$.
\index{$\nabla$, $\hslash$-$T^\mr{log}$-connection@--- $\psi^\nabla$, curvature of}
Also, we shall say that $\nabla$ is ${\bf flat}$
\index{$\nabla$, $\hslash$-$T^\mr{log}$-connection@--- flat (= integrable)}
(or {\bf  integrable}) if  $ \psi^{\nabla} =0$.
\ede

\begin{rema} \label{QQ881}
Since $d_\mcE^\mr{log}$ preserves the Lie bracket operations,  
we have
\begin{align}
d_\mcE^\mr{log} \circ \psi^\nabla (\partial_1 \wedge \partial_2) &= 
d_\mcE^\mr{log} ([\nabla (\partial_1), \nabla (\partial_2)]-\hslash \cdot \nabla ([\partial_1, \partial_2])) \\
& = [d_\mcE^\mr{log} \circ \nabla (\partial_1), d_\mcE^\mr{log} \circ \nabla (\partial_2)] - \hslash \cdot d_\mcE^\mr{log} \circ \nabla ([\partial_1, \partial_2]) \notag \\
&= [\hslash \cdot \partial_1, \hslash \cdot \partial_2] - \hslash \cdot \hslash \cdot [\partial_1, \partial_2] \notag \\
& = 0. \notag
\end{align}
It follows that  
  $d_\mcE^\mr{log} \circ \psi^{\nabla} =0$, or equivalently,
 that
   the image of $\psi^{\nabla}$ is contained in $\mfg_\mcE \left(= \mr{Ker}(d_\mcE^\mr{log}) \right)$.
Hence,  $\psi^{\nabla}$ may be regarded   as an $\mcO_Y$-linear morphism $\bigwedge^2 \mcT_{Y^\mr{log}/T^\mr{log}} \migi \mfg_\mcE$.
\end{rema}

\begin{rema}[Log curve] \label{QQffg14} 
  Suppose that  $Y^\mr{log}/T^\mr{log}$ is a log curve in the sense of Definition \ref{pY5019} described later.
  Then,  $\mcT_{Y^\mr{log}/T^\mr{log}}$ is a line bundle, so  $\bigwedge^2  \mcT_{Y^\mr{log}/T^\mr{log}}=0$.
It follows that  any $\hslash$-$T^\mr{log}$-connection is automatically  flat. 
Since the main part of the subsequent discussions in this manuscript  deals with connections on  log curves,  
the concept of flatness may not be so important to
 us.
 \end{rema}

\begin{exa}[Trivial bundle] \label{QQ9} 
Suppose that $2$ is invertible in $k$.
Let $\nabla$ be an  $\hslash$-$T^\mr{log}$-connection on the trivial $\mbG$-bundle $Y \times_k \mbG$.
In the natural manner, we shall  regard  $\nabla^\mfg$  (cf. (\ref{QQ12})) as
an element of $\Gamma (Y, \Omega_{Y^\mr{log}/T^\mr{log}}\otimes_k \mfg)$
and $\psi^\nabla$ as an element of $\Gamma (Y, (\bigwedge^2\Omega_{Y^\mr{log}/T^\mr{log}})\otimes_k \mfg)$.
Then,  the following equality holds:
\begin{align} \label{QR990}
\psi^\nabla = d\nabla^\mfg + \frac{\hslash}{2} \cdot [\nabla^\mfg, \nabla^\mfg],
\end{align}
where the operators ``$[-,-]$'' and ``$d$'' are defined
 in the same fashions as (\ref{QQ20}) and  (\ref{QQ23}) respectively.
In particular, if $\omega \left(\in \Gamma (Y, \Omega_{Y^\mr{log}/T^\mr{log}}) \right)$ is the $1$-form corresponding to an $\hslash$-$T^\mr{log}$-connection $\nabla$ on the trivial $\mbG_m$-bundle (cf. \S\,\ref{QR304}),
then $\nabla$ is flat if and only if $\omega$ is closed.
\end{exa}

\bpr  \label{y025}
Let 
$\mcE$ be a $\mbG$-bundle on $Y$, $\nabla$ an $\hslash$-$T^\mr{log}$-connection on $\mcE$, and $\eta$ an automorphism of $\mcE$.
Denote by $\mr{Ad}_\mbG (\eta)$ the automorphism of $\mfg_\mcE$ induced from  $\eta$ via  change of structure group by $\mr{Ad}_\mbG$. 
Then,
the equality
\begin{equation}
 \psi^{\eta_*(\nabla)} = \mr{Ad}_\mbG(\eta) \circ \psi^{\nabla}
\end{equation}
of morphisms $\bigwedge^2\mcT_{Y^\mr{log}/T^\mr{log}} \migi \mfg_\mcE$
 holds.
 In particular, the flatness of an $\hslash$-$T^\mr{log}$-connection is closed under gauge transformations.
\epr
\begin{proof}
The assertion follows immediately  from 
 the fact that the automorphism $d \eta^\mbG$ is  compatible with the Lie bracket operation  on $\widetilde{\mcT}_{\mcE^\mr{log}/T^\mr{log}}$ and  restricts to $\mr{Ad}_\mbG (\eta)$.
\end{proof}


\bde
 \label{QQ6}
  \begin{itemize}
  \item[(i)]
  An  {\bf  $\hslash$-flat $\mbG$-bundle} \index{$\hslash$-flat $\mbG$-bundle}
 on $Y^\mr{log}/T^\mr{log}$ is a pair
 \begin{align}
 (\mcE, \nabla)
 \end{align}
  consisting of a $\mbG$-bundle  $\mcE$ on $Y$ and a flat  $\hslash$-$T^\mr{log}$-connection $\nabla$ on $\mcE$.
  For simplicity, a $1$-flat  $\mbG$-bundle will  be referred to as a {\bf flat   $\mbG$-bundle}.
  \item[(ii)]
  Let $(\mcE, \nabla)$, $(\mcE', \nabla')$ be $\hslash$-flat $\mbG$-bundles on $Y^\mr{log}/T^\mr{log}$.
  Then, an {\bf isomorphism of $\hslash$-flat $\mbG$-bundles} from $(\mcE, \nabla)$ to $(\mcE', \nabla')$ is defined as an isomorphism of $\mbG$-bundles $\eta : \mcE \isom \mcE'$ 
  with $\eta_*(\nabla) = \nabla'$.
  \end{itemize}
  \ede
\subsection{Change of structure group} \label{QQ809}
Let $\mcE$ and $\hslash$ be as before.
Also, let $\mbG'$ and $w$ be as in \S\,\ref{QQ119}.
If $\nabla$ is  an $\hslash$-$T^\mr{log}$-connection  on $\mcE$,
then the composite
\index{$\nabla$, $\hslash$-$T^\mr{log}$-connection@--- $w_*(\nabla)$}
\begin{align} \label{QQ813}
w_*(\nabla) : \mcT_{Y^\mr{log}/T^\mr{log}} \xrightarrow{\nabla} \widetilde{\mcT}_{\mcE^\mr{log}/T^\mr{log}} \xrightarrow{(dw)_\mcE} \widetilde{\mcT}_{(\mcE \times^\mbG \mbG')^\mr{log}/T^\mr{log}}
\end{align}
specifies an $\hslash$-$T^\mr{log}$-connection on $\mcE \times^\mbG \mbG'$.
This connection  will be denoted by $\nabla \times^{\mbG, w} \mbG'$ or $\nabla \times^{\mbG} \mbG'$  depending on  the situation.
Since  
$(dw)_\mcE$ preserves  the Lie bracket operations,  
we have
\begin{align} \label{QS56}
\psi^{w_*(\nabla)} = (dw)_\mcE \circ \psi^\nabla.
\end{align}
It follows that the flatness of $\psi^\nabla$ implies that of $\psi^{w_*(\nabla)}$, and vice versa if $w$ is injective.
In particular,   an $\hslash$-flat $\mbG$-bundle $(\mcE, \nabla)$ yields 
 an $\hslash$-flat $\mbG'$-bundle  
$(\mcE\times^\mbG \mbG', w_*(\nabla))$.

\subsection{Change of parameter} \label{QQ1030}
Let  $\nabla$ be an $\hslash$-$T^\mr{log}$-connection on $\mcE$.
Also, let us take $c \in \Gamma (T, \mcO_T)$.
Then, the $\mcO_T$-linear morphism
\index{$\nabla$, $\hslash$-$T^\mr{log}$-connection@--- $c  \cdot \nabla$}
\begin{align} \label{QQ1032}
c  \cdot \nabla : \mcT_{Y^\mr{log}/T^\mr{log}} \migi \widetilde{\mcT}_{\mcE^\mr{log}/T^\mr{log}},
\end{align}
i.e., the morphism $\nabla$ multiplied by $c$, forms a $(c \cdot \hslash)$-$T^\mr{log}$-connection on $\mcE$.
One verifies that the curvature of  $c \cdot \nabla$ is given by
\begin{align} \label{EE1033}
\psi^{c \cdot \nabla} = c^2 \cdot \psi^{\nabla}.
\end{align}
In particular, $c \cdot \nabla$ is flat if $\nabla$ is flat,  and the inverse is true when $c \in \Gamma (T, \mcO^\times_T)$.

\subsection{Products} \label{QS7}
Let $\mbG'$ be another connected smooth affine algebraic group over $k$ with Lie algebra $\mfg'$.
Then, the Lie algebra of the product $\mbG \times_k \mbG'$ is  the direct sum $\mfg \oplus \mfg'$.
Next, let $\mcE$ (resp., $\mcE'$) be a $\mbG$-bundle (resp., a $\mbG'$-bundle) on $Y$.
The fiber product
\index{$\mcE, $(principal) $\mbG$-bundle@--- $\mcE \boxtimes \mcE'$}
\begin{align} \label{QS8}
\mcE \boxtimes \mcE' := \mcE \times_Y \mcE'
\end{align}
admits, in the natural manner, a structure of $\mbG \times_k \mbG'$-bundle.
There exists a canonical isomorphism $(\mfg \oplus \mfg')_{\mcE \boxtimes \mcE'} \isom \mfg_\mcE \oplus \mfg'_{\mcE'}$, and 
 the following sequence is exact:
\begin{align} \label{QS9}
0 \longmigi \widetilde{\mcT}_{(\mcE \boxtimes \mcE')^\mr{log}/T^\mr{log}} \longmigi \widetilde{\mcT}_{\mcE^\mr{log}/T^\mr{log}} \oplus 
\widetilde{\mcT}_{\mcE'^{\mr{log}}/T^\mr{log}} \xrightarrow{(d_\mcE^\mr{log}, - d_{\mcE'}^\mr{log})} \mcT_{Y^\mr{log}/T^\mr{log}} \longmigi 0,
\end{align}
where the second arrow is obtained by differentiating the projections $\mcE \boxtimes \mcE' \migi \mcE$, $\mcE \boxtimes \mcE' \migi \mcE'$.
By using the second arrow, we shall regard  $\widetilde{\mcT}_{(\mcE \boxtimes \mcE')^\mr{log}/T^\mr{log}}$ as an $\mcO_Y$-submodule of 
 $\widetilde{\mcT}_{\mcE^\mr{log}/T^\mr{log}} \oplus 
\widetilde{\mcT}_{\mcE'^{\mr{log}}/T^\mr{log}}$.

Now, let $\nabla$ (resp., $\nabla'$) be an $\hslash$-$T^\mr{log}$-connection on $\mcE$ (resp., $\mcE'$).
Then,  the image of the morphism $(\nabla, \nabla') : \mcT_{Y^\mr{log}/T^\mr{log}} \migi \widetilde{\mcT}_{\mcE^\mr{log}/T^\mr{log}} \oplus 
\widetilde{\mcT}_{\mcE'^{\mr{log}}/T^\mr{log}}$ is contained in 
$\widetilde{\mcT}_{(\mcE \boxtimes \mcE')^\mr{log}/T^\mr{log}}$.
The resulting $\mcO_Y$-linear morphism
\index{$\nabla$, $\hslash$-$T^\mr{log}$-connection@--- $\nabla \boxtimes \nabla'$}
\begin{align} \label{QS11}
\nabla \boxtimes \nabla' : \mcT_{Y^\mr{log}/T^\mr{log}} \migi \widetilde{\mcT}_{(\mcE \boxtimes \mcE')^\mr{log}/T^\mr{log}}
\end{align}
specifies an $\hslash$-$T^\mr{log}$-connection on $\mcE \boxtimes \mcE'$.
The curvature $\psi^{\nabla \boxtimes \nabla'}$ of $\nabla \boxtimes \nabla'$ satisfies 
\begin{align} \label{QS10}
\psi^{\nabla \boxtimes \nabla'} = (\psi^\nabla, \psi^{\nabla'})
\end{align}
(as morphisms $\mcT_{Y^\mr{log}/T^\mr{log}} \migi (\mfg \oplus \mfg')_{\mcE \boxtimes \mcE'} \cong \mfg_\mcE \oplus \mfg'_{\mcE'}$).
In particular, if  $(\mcE, \nabla)$ and $(\mcE', \nabla')$  form an $\hslash$-flat $\mbG$-bundle  and  an $\hslash$-flat $\mbG'$-bundle respectively, then the pair $(\mcE \boxtimes \mcE', \nabla \boxtimes \nabla')$ forms an $\hslash$-flat $\mbG \times_k \mbG'$-bundle; we refer to it  as the {\bf product} of $(\mcE, \nabla)$ and $(\mcE', \nabla')$.

\section{Explicit descriptions of gauge transformations}\label{y026}

In what follows, we describe gauge transformations by using 
the exponential and logarithm maps, which connect  unipotent elements in $\mbG$ with  nilpotent elements in $\mfg$. 

\subsection{$p$-unipotent/$p$-nilpotent elements} \label{QR401}
Denote by $p$ the characteristic $\mr{char}(k)$ of $k$.
In this section, we assume   (except for the discussion in Remark \ref{ppp037})  that  $k$ is {\it perfect} and that 
 $(\mbG, \mbT)$ is a  {\it  split  semisimple algebraic group over $k$ of adjoint type}, where $\mbT$ 
 \index{$\mbG$, algebraic group@— $\mbT$, maximal torus of}
 denotes a maximal torus of $\mbG$.
 In particular, the adjoint representation $\mr{Ad}_\mbG : \mbG \migi \mr{GL}(\mfg)$ becomes a closed immersion;
 by using this, we regard  $\mbG$ as  a closed subgroup of $\mr{GL} (\mfg)$.
 Fix a Borel subgroup $\mbB$
 \index{$\mbG$, algebraic group@— $\mbB$, Borel subgroup of}
   of $\mbG$ defined over $k$ containing $\mbT$, and set $\mbN := [\mbB, \mbB]$, i.e., the unipotent radical of $\mbB$. 
   \index{$\mbG$, algebraic group@— $\mbN$}
Denote by $\mfn$ the Lie algebra of $\mbN$.
Here, note that each root $\beta$ in $\mbG$  can be expressed uniquely as $\beta = \sum_{\alpha \in \Gamma} n_\alpha \cdot \alpha$ with $n_\alpha \in \mbZ$, where $\Gamma$ denotes the set of simple roots in $\mbB$ with respect to $\mbT$;
we call $\sum_{\alpha \in \Gamma} n_\alpha$ the {\bf height} \index{height}
of $\beta$.
We shall set
\index{$h_\mbG$}
 \begin{align} \label{QH9932}
 h_\mbG := 1 + (\text{the height of a root with maximum height}).
\end{align}
If $\mbG$ is simple, then $h_\mbG$ coincides with the Coxeter number of $\mbG$.
We  assume further that   {\it $p=0$ or  $p> 2  h_\mbG$}.
In particular, we will not deal with the case where $p=2, 3$.

\begin{rema} \label{QR1000}
By the semisimplicity  of  $\mbG$ and the assumption ``$p=0$ or  $p> 2 h_\mbG$'' imposed above,  the Killing form on
$\mfg$ is nondegenerate; this fact follows from a formula for the discriminant  of the Killing form proved in ~\cite{SpSt}, \S\,4, (4.8) (cf. ~\cite{Ngo}, \S\,1.1).
Hence, since $k$ is a perfect field  of characteristic $\neq 2,3$  and $(\mbG, \mbT)$ is split over $k$,
the semisimple Lie algebra $\mfg$ is classical in the sense of ~\cite{Sel}, Chap.\,II, \S\,3.
It follows (cf. ~\cite{Sel}, Chap.\,III, \S\,5, Theorem III.5.1) that 
the image of the adjoint representation $\mr{Ad}_\mbG$ 
coincides with the identity component  $\mr{Aut}(\mfg)^0 \left(\subseteq \mr{GL}(\mfg) \right)$  of the group of Lie algebra automorphisms of $\mfg$.
In particular,  
 $\mr{Ad}_\mbG$   induces an isomorphism
$\mbG \isom \mr{Aut}(\mfg)^0$.

Also, according to ~\cite{Sel}, Chap.\,I, \S\,7,  Corollary,  any semisimple Lie algebra with nondegenerate Killing form can be decomposed  into a direct sum of simple Lie algebras. 
 Thus, the problem of understanding the structures of such  semisimple Lie algebras can be largely reduced to those of simple Lie algebras.
 For example, by applying  \cite{Car},  Chap.\,5, \S\,5.5, Proposition 5.5.2, we see that
if $\left(\mr{char}(k) =: \right) p>2 h_\mbG$,
  then  the equality $\mr{ad}(v)^p =0$ holds for any nilpotent element $v$ of $\mfg$.

\end{rema}

Denote by $\mbG^\text{unip}$
 \index{$\mbG$, algebraic group@--- $\mbG^\text{unip}$, unipotent variety of} 
 (resp., $\mfg^\text{nilp}$) 
 \index{$\mfg$, Lie algebra@--- $\mfg^\text{nilp}$, nilpotent variety of} 
the unipotent variety of $\mbG$ (resp., the nilpotent variety of $\mfg$), which  
is  the Zariski closed subset of $\mbG$ (resp., $\mfg$)  consisting of all unipotent elements of $\mbG$   (resp., all  nilpotent elements of $\mfg$).
If $p =0$ (resp., $p >0$), then
we shall write 
\begin{equation}
\mr{GL}(\mfg)^{\text{$p$-unip}}
\index{$\mr{GL}(\mfg)^{\text{$p$-unip}}$}
\end{equation}
for  the Zariski closed subset  in  $\mr{GL}(\mfg)$  consisting of all elements  $h$  with  $(h -\mr{id}_\mfg)^{m} = 0$ for some $m \in \mbZ_{>0}$ (resp., $(h -\mr{id}_\mfg)^{p} = 0$). 
Also, write 
\begin{equation}
\mr{End}(\mfg)^{\text{$p$-nilp}}
\index{$\mr{End}(\mfg)^{\text{$p$-nilp}}$}
\end{equation}
 for the Zariski closed subset in  $\mr{End}(\mfg)$ consisting  of  all elements $v$  satisfying that $v^m =0$ for some $m \in \mbZ_{>0}$ (resp., $v^p =0$).
We endow $\mbG^\text{unip}$ (resp., $\mfg^\text{nilp}$; resp.,  $\mr{GL}(\mfg)^{\text{$p$-unip}}$; resp.,  $\mr{End}(\mfg)^{\text{$p$-nilp}}$) with the reduced subscheme structure of  
$\mbG$ (resp., $\mfg$; resp., $\mr{GL}(\mfg)$; resp., $\mr{End}(\mfg)$).
In particular, there exist natural closed immersions $\mbN \migiincl \mbG^\text{unip}$,  $\mfn \migiincl \mfg^\text{nilp}$.


\subsection{Exponentials and logarithms} \label{QR400}

For  each  $v \in \mr{End}(\mfg)^{\text{$p$-nilp}}$,  
the {\bf  exponential}  of $v$ is defined to be 
\begin{equation} \label{QQ30}
\mr{Exp} (v) := \sum_{s =0}^\infty \frac{1}{s!} \cdot v^s \ \left( = \sum_{s =0}^{p-1} \frac{1}{s!} \cdot v^s,  \ \text{if $p>0$}\right).  \index{$\mr{Exp}$, exponential map}
\end{equation} 
Also, for $h \in \mr{GL} (\mfg)^{\text{$p$-unip}}$,
the {\bf logarithm} of $h$ is defined to be
 \begin{equation} \label{QQ31}
 \mr{Log} (h)  := \sum_{s =1}^\infty \frac{(-1)^{s-1}}{s} \cdot (h-\mr{id}_\mfg)^s \ \left( = \sum_{s =1}^{p-1} \frac{(-1)^{s-1}}{s} \cdot (h-\mr{id}_\mfg)^s,  \ \text{if $p>0$} \right). \index{$\mr{Log}$, logarithm map}
 \end{equation}
The assignments $v \mapsto \mr{Exp} (v)$ and $h \mapsto \mr{log}(h)$ determine  isomorphisms of $k$-schemes
\begin{equation} \label{QQ32}
\mr{Exp} : \mr{End}(\mfg)^{\text{$p$-nilp}} \isom \mr{GL}(\mfg)^{\text{$p$-unip}}, \hspace{5mm} \mr{Log} :  \mr{GL}(\mfg)^{\text{$p$-unip}}\isom 
\mr{End}(\mfg)^{\text{$p$-nilp}}, 
\end{equation}
 which are inverses of each other.

\bpr \label{QR703}
Let  $v: T \migi \mr{End}(\mfg)^{p\text{-}\mr{nilp}}$ (resp., $h : T \migi \mr{GL}(\mfg)^{p\text{-}\mr{unip}}$) be a morphism of $k$-schemes.
Then, the following  equality  holds:
\begin{align} \label{pp100}
\mr{Exp}(-v) = \mr{Exp}(v)^{-1} \ \left(\text{resp.,} \   \mr{Log} (h^{-1}) = - \mr{Log} (h)\right).
\end{align}
\epr
\begin{proof}
The non-resp'd assertion follows from the observation  
that 
\begin{align}
\mr{Exp}(v) \cdot \mr{Exp}(-v) = \mr{Exp}(v + (-v)) = \mr{Exp}(0) =\mr{id}_{\mfg}.
\end{align}
The resp'd assertion is a direct consequence of the  non-resp'd assertion together with   the fact that $\mr{Exp}$ and $\mr{Log}$ are inverses of each other.
\end{proof}

The group homomorphism  $\mr{Ad}_{\mr{GL}(\mfg)} : \mr{GL} (\mfg) \migi \mr{GL}(\mr{End}(\mfg))$ (resp., $\mr{Inn}_{(-)} : \mr{GL}(\mfg) \migi \mr{Aut}(\mr{GL}(\mfg))$)
  defines  a left $\mr{GL}(\mfg)$-action  on $\mr{End}(\mfg)^{\text{$p$-nilp}}$ (resp., $\mr{GL}(\mfg)^{\text{$p$-unip}}$).
  Moreover,  it induces a left $\mbG$-action on  $\mr{End}(\mfg)^{\text{$p$-nilp}}$ (resp., $\mr{GL}(\mfg)^{\text{$p$-unip}}$) via composition  with    $\mr{Ad}_\mbG : \mbG \migi \mr{GL}(\mfg)$.
Both $\mr{Exp}$ and $\mr{Log}$ are compatible with  the $\mr{GL}(\mfg)$-actions, and hence with  the induced $\mbG$-actions.
(That is to say,  if $p$ is sufficiently large relative to $h_\mbG$, then  $\mr{Exp}$ defines a Springer isomorphism for $\mr{GL}(\mfg)$.)


According to the final comment in Remark \ref{QR1000}, 
the morphisms $\mr{ad} : \mfg \migiincl \mr{End}(\mfg)$ and  $\mr{Ad}_\mbG : \mbG \migiincl \mr{GL} (\mfg)$
 restrict to  closed immersions 
$\mfg^\text{nilp} \migiincl \mr{End}(\mfg)^{\text{$p$-nilp}}$ and 
$\mbG^\text{unip} \migiincl  \mr{GL}(\mfg)^{\text{$p$-unip}}$ respectively.
Moreover,  by  the definition of $\mr{Exp}$, 
the image $\mr{Exp} (\mr{ad}(v)) \in \mr{GL}(\mfg)^{\text{$p$-unip}}$ of each $v \in \mfg^\text{nilp}$
defines a  Lie algebra automorphism of $\mfg$.
This implies  (since $\mbG \cong \mr{Aut}(\mfg)^0$ as mentioned in Remark \ref{QR1000})  that 
$\mr{Exp}$ restricts to a closed immersion
\begin{align} \label{QR610}
\mr{Exp}: \mfg^\text{nilp} \migiincl  \mbG^\text{unip}.
\index{$\mr{Exp}$, exponential map}
\end{align}
This morphism in fact is an isomorphism because $\mfg^\text{nilp}$ and $\mbG^\text{unip}$ are irreducible and of the same dimension.

\bpr \label{QR602}
The isomorphism  $\mr{Exp}$ restricts to an isomorphism  of $k$-schemes
\begin{align} \label{QR603}
\mr{Exp} :  \mfn \isom \mbN.
\end{align}
\epr
\begin{proof}
Let us fix an identification  of (the underlying $k$-vector space of) $\mfg$ with $k^{\oplus l}$, where $l : =\mr{dim}(\mfg)$.
Under this identification,  we have $\mr{GL}(\mfg) =\mr{GL}_l$ and $\mr{End}(\mfg) = \mfg \mfl_l$.
Denote by   $\mbN_{\mr{GL}_l}$ (resp., $\mfn_{\mfg \mfl_l}$)
the subgroup of $\mr{GL}_l$ (resp., $\mfg \mfl_l$) consisting of upper triangular matrices with diagonal entries $1$ (resp., $0$).
The isomorphism $\mr{Exp} : \mr{End}(\mfg)^{p\text{-}\mr{nilp}}\isom \mr{GL}(\mfg)^{p \text{-}\mr{unip}}$ restricts, by construction,  to an isomorphism
$\mr{End}(\mfg)^{p\text{-}\mr{nilp}} \cap \mfn_{\mfg \mfl_l} \isom \mr{GL}(\mfg)^{p \text{-}\mr{unip}} \cap \mbN_{\mr{GL}_l}$.
On the other hand,  $\mbN$  is  contained, via $\mr{Ad}_\mbG$,  in a maximal closed connected unipotent subgroup  of $\mr{GL}_l$, which is conjugate to $\mbN_{\mr{GL}_l}$ (cf. ~\cite{Malle}, \S\,6, Corollary 6.11, (c)).
Hence, since $\mr{Exp}$ is $\mr{GL}(\mfg)$-equivariant, we may assume, 
 after possibly replacing the fixed identification $\mfg \cong k^{\oplus l}$ with another,
 that $\mbN \subseteq  \mbN_{\mr{GL}_l}$.
 This assumption implies that  $\mbN = \mbG^{\mr{unip}} \cap \mbN_{\mr{GL}_l}$ and  $\mfn = \mfg^{\mr{nilp}} \cap \mfn_{\mfg \mfl_l}$.
Thus,  $\mr{Exp}$  restricts to $\mfn \left(=\mfg^{\mr{nilp}} \cap \mfn_{\mfg \mfl_l} \right) \isom  \mbN \left(=\mbG^{\mr{unip}} \cap \mbN_{\mr{GL}_l}  \right)$ (cf. (\ref{QR610})), 
 completing  the proof of the assertion.
\end{proof}

By the above proposition, we obtain the inverse
\begin{align} \label{QR630}
\mr{Log} :  \mbN \isom \mfn
\index{$\mr{Log}$, logarithm map}
\end{align}
to $\mr{Exp}: \mbN \isom \mfn$, which makes the following square diagram commute:
\begin{align}\label{comm331}
\vcenter{\xymatrix@C=46pt@R=36pt{
 \mbN \ar[r]^-{\mr{Ad}_\mbG |_\mbN} \ar[d]^-{\wr}_-{\mr{Log}}& \mr{GL}(\mfg)^{\text{$p$-unip}}  \ar[d]_-{\wr}^-{\mr{Log}} \\
\mfn\ar[r]_-{\mr{ad} |_\mfn}&  \mr{End}(\mfg)^{\text{$p$-nilp}}.
}}
\end{align}
The commutativity of  this diagram  implies the following proposition.

\bpr \label{y027}
 For a $k$-scheme $T$ and a $T$-rational point $h : T \migi \mbN$ of $\mbN$,
we have an equality
\begin{equation}
\mr{Ad}_\mbG (h) = \mr{Exp}(\mr{ad}(\mr{Log} (h)))  \left( = \sum_{s =0}^\infty \frac{1}{s !} \cdot \mr{ad}(\mr{Log} (h))^s\right)
\end{equation}
of $\mcO_T$-linear automorphisms of $\mcO_T \otimes_k \mfg$.
\epr


\subsection{Gauge transformation using the logarithm map} \label{QR650}

Since the tangent bundle $\mcT_{\mfg/k}$ of the $k$-scheme  $\mfg$ is canonically isomorphic to $\mcO_\mfg \otimes_k \mfg$, the differential of $\mr{Log} : \mbN \isom \mfn$
defines an $\mcO_\mbN$-linear morphism
\begin{equation} \label{QR800}
d\mr{Log} : \mcT_{\mbN/k} \migi \mcO_\mbN \otimes_k \mfn,
\end{equation}
or equivalently, an element of $\Gamma (\mbN, \Omega_{\mbN/k} \otimes_k \mfn)$.

\bpr \label{QR410}
The following equality holds:
\begin{align} \label{QR411}
\Theta_\mbN = d \mr{Log}. 
\end{align}
In particular, for a $k$-scheme $T$ and a $T$-rational point $h : T \migi \mbN$ of $\mbN$, the  equality
\begin{align} \label{QR801}
h^*(\Theta_\mbN) = d \mr{Log}(h) \left(:= d (\mr{Log} \circ h) \right)
\end{align}
of morphisms $\mcT_{T/k} \migi \mcO_T \otimes \mfn$ holds.
\epr
\begin{proof}
After possibly changing the base
 by a field extension of $k$,
we may assume that $k$ is algebraically closed.
Denote by $A: \mr{GL}(\mfg) \migiincl \mr{End}(\mfg)$ the natural open immersion, and 
consider  $\mbN$ as a closed subscheme of $\mr{End} (\mfg)$ via the composite of $\mr{Ad}_\mbG$ and $A$.
In particular, (since $\mbN \subseteq \mr{GL}(\mfg)$) we have  
\begin{align} \label{QR655}
\Theta_\mbN = \Theta_{\mr{GL}(\mfg)} |_{\mbN} = (A^{-1} dA) |_\mbN
\end{align}
(cf. Example \ref{QQ29}).
Now,  
let $h$ be a $k$-rational point of $\mbN \left(\subseteq \mr{End}(\mfg) \right)$ and 
 take  an element $R \in \mcT_h \mbN$, which classifies the element $h + \epsilon \cdot R$ of $\mbN (k [\epsilon]/(\epsilon^2)) \subseteq \mr{End}(\mfg) (k [\epsilon]/(\epsilon^2)) = \mr{End}(\mfg) (k) \oplus \epsilon \cdot \mr{End} (\mfg) (k)$.
The following sequence of equalities holds:
\begin{align} \label{QR659}
\mr{Log} (h + \epsilon \cdot  R) &= \sum_{s=1}^\infty \frac{(-1)^{s-1}}{s} \cdot (h - \mr{id}_\mfg + \epsilon  \cdot R)^{s} \\
&= \sum_{s=1}^\infty \frac{(-1)^{s-1}}{s} \cdot (h - \mr{id}_\mfg)^s + \epsilon \cdot \sum_{s=1}^\infty (-1)^{s-1} \cdot (h-\mr{id}_\mfg)^{s-1} \cdot R \notag \\
& = \mr{Log}(h) + \epsilon \cdot h^{-1}  \cdot R. \notag
\end{align}
Hence,  we have
\begin{align}
d\mr{Log} |_h (R) \stackrel{(\ref{QR659})}{=} h^{-1} \cdot R 
\stackrel{(\ref{QR655})}{=} \Theta_\mbG |_h (R).
\end{align}
By applying  this fact to various  $h$'s and $R$'s,
we have  $\Theta_\mbN = d \mr{Log}$, as desired.
\end{proof}

\bco \label{y028}
(Let us keep the notation and  assumptions described  at the beginning of \S\,\ref{QR401}.)
Let $Y^\mr{log}/T^\mr{log}$ be as in \S\,\ref{QR250}.
Also, let  $\hslash$ be an element of $\Gamma (T, \mcO_T)$, $\nabla$  an $\hslash$-$T^\mr{log}$-connection on the trivial $\mbG$-bundle $Y \times_k \mbG$,
 and $h :  Y \migi \mbN$ a $Y$-rational point of $\mbN$.
Then, we have an equality
\begin{align} \label{QG5023}
\mr{L}_{h*} (\nabla)^\mfg =  - \hslash  \cdot d \mr{Log}(h) +    \sum_{s =0}^\infty \frac{1}{s !} \cdot \mr{ad}(\mr{Log}(h))^s \circ \nabla^\mfg
\end{align}
of morphisms $\mcT_{Y^\mr{log}/T^\mr{log}}\migi \mcO_Y \otimes_k \mfg$.
\eco
\begin{proof}
The assertion follows from Propositions \ref{y024},  \ref{QR703}, \ref{y027}, and  \ref{QR410}.
\end{proof}

\begin{rema} \label{ppp037} 
Let $n$ be an integer with $n>1$.
If $p=0$ or  $n<p$, then 
we can obtain the description of a gauge transformation  (\ref{QG5023})  for  $\mbG = \mr{GL}_n$  ($n>1$) 
although it is not of adjoint type.
Indeed, let  $\mbN$ (resp., $\mfn$) denote 
the subgroup of $\mr{GL}_n$ (resp., $\mfg \mfl_n$) consisting of upper triangular matrices with diagonal entries $1$ (resp., $0$).
If  $p>0$, then 
it follows from  the assumption $n<p$ that 
 $v^p= (h-1)^p =0$ for any $v \in \mfn$ and $h \in \mbN$.
Then,  the exponential $\mr{Exp}(v)$ of $v$  and  the logarithm $\mr{Log}(h)$ of $h$  can be  defined as in (\ref{QQ30}) and (\ref{QQ31}) respectively.
The assignments $v \mapsto \mr{Exp}(v)$  and $h \mapsto \mr{Log}(h)$ specify
isomorphisms $\mr{Exp}: \mfn \isom \mbN$ and $\mr{Log}: \mbN \isom \mfn$ respectively,
 which are inverses of each other and satisfy the same assertion as Proposition \ref{QR703}.
For each  $h \in \mbN (T) \subseteq  \Gamma (T, \mcO_T \otimes_k \mfg \mfl_n)$  and $R \in \Gamma (T,  \mcO_T \otimes_k \mfg \mfl_n)$, 
the equality $\mr{Ad}_{\mr{GL}_n} (R) = h \cdot R \cdot h^{-1}$ holds.
On the other hand, as in the proof of Proposition \ref{QR410},
we have $\mr{\Theta}_\mbN = d \mr{Log}$.
Hence, keeping the notation in Corollary \ref{y028}, we obtain from Proposition \ref{y024} an equality
\begin{align}
\mr{L}_{h*}(\nabla)^{\mfg \mfl_n} (\partial)  =  - \hslash \cdot d \mr{Log}(h) + h \cdot \nabla^{\mfg \mfl_n}(\partial) \cdot h^{-1}
\end{align}
for any local section $\partial \in \mcT_{Y^\mr{log}/T^\mr{log}}$.

 \end{rema}


\section{The moduli stack of pointed stable curves} \label{y029}

We shall recall pointed stable curves and 
their canonical log structures pulled-back from  the moduli stack.
Let  us fix a pair of nonnegative integers $(g, r)$.

\subsection{Pointed curves} \label{QR779}
Let $S$ be a $k$-scheme.
By a {\bf curve} 
 \index{curve}  over  $S$, we mean   
a proper flat scheme $f : X \migi S$ over $S$ all of whose geometric fibers are connected, reduced, and of dimension $1$. 

\bde \label{QR860}
\begin{itemize}
\item[(i)]
 An {\bf $r$-pointed curve} \index{curve@--- $r$-pointed} over $S$ is defined as a collection of data
 \index{$\msX$, pointed (stable) curve}
\begin{align} \label{QR833}
\msX := (f: X \migi S, \{ \sigma_i\}_{i=1}^r), 
\end{align}
consisting of   a curve $f : X \migi S$ over $S$ and an ordered set of $r$ sections   $\sigma_i : S \migi X$ ($i = 1, \cdots, r$)   of $f$ whose images are disjoint, i.e.,
$\mr{Im}(\sigma_i) \cap \mr{Im}(\sigma_j) = \emptyset$ if $i \neq j$,  and land in the smooth locus of $f$. 
We refer to $\sigma_i$ as the {\bf $i$-th marked point} of $\msX$.
A $0$-pointed curve is nothing but a curve  in  the above sense and  will be  often referred to as an {\bf unpointed curve}.
\item[(ii)]
Suppose that, for each $j \in \{1,2\}$, we are given a $k$-scheme $S_j$ and an $r$-pointed curve $\msX_j := (f_j: X_j \migi S_j, \{ \sigma_{ji} \}_{i=1}^r)$ over $S_j$.
Then, a {\bf morphism of $r$-pointed curves} from 
$\msX_1$ to $\msX_2$ is a pair of $k$-morphisms
\begin{align} \label{QQ09}
(\phi : S_1 \migi S_2, \Phi : X_1 \migi X_2)
\end{align}
 such that the following squares
   form cartesian commutative diagrams:
\begin{align}
\vcenter{\xymatrix@C=46pt@R=36pt{
X_1 \ar[r]^-{\Phi} \ar[d]_-{f_1}& X_2\ar[d]^-{f_2}\\
S_1 \ar[r]_-{\phi}& S_2,
}}
\hspace{10mm}
\vcenter{\xymatrix@C=46pt@R=36pt{
X_1 \ar[r]^-{\Phi} & X_2\\
S_1 \ar[u]^-{\sigma_{1i}}  \ar[r]_-{\phi}& S_2 \ar[u]_-{\sigma_{2i}}
}} \hspace{3mm}\ (i=1, \cdots, r). \hspace{-15mm}
\end{align}
\end{itemize}
\ede

\bde \label{QR863}
\begin{itemize}
\item[(i)]
Let $\overline{k}$ be an algebraically closed field over $k$ and $X$  a curve over  $\overline{k}$.
A $\overline{k}$-rational point $x$ of $X$ is called
 a {\bf nodal point}, or an  {\bf ordinary double point},
   \index{curve@--- nodal point (= ordinary double point) of}
  if the complete local ring $\widehat{\mcO}_{X, x}$ at $x$ is isomorphic to $\overline{k} [[u, v]]/ (uv)$.
\item[(ii)]
An {\bf $r$-pointed smooth} (resp., {\bf prestable})\index{curve@--- smooth} \index{curve@--- prestable}
{\bf curve} over  $S$ is an $r$-pointed  curve $\msX := (f : X \migi S, \{ \sigma_i \}_i)$ over $S$
such that, for each  geometric point $s : \mr{Spec}(\overline{k}) \migi S$ of $S$,
any point of the fiber  $X_{s} := f^{-1}(s)$ is a smooth point (resp.,  either a smooth point or a nodal point).
\item[(iii)]
Let  $\msX:= (f :X \migi S, \{ \sigma_i\}_{i})$ be  an $r$-pointed  curve  over $S$.
We shall say that
$\msX$  is {\bf of (arithmetic) genus $g$} \index{curve@--- of (arithmetic) genus $g$} if, 
for any geometric point
$s : \mr{Spec}(\overline{k}) \migi S$ of $S$,   the equality $\mr{dim}_{\overline{k}} (H^1 (X_{s}, \mcO_{X_{s}})) = g$ holds.
\end{itemize}
\ede

Let $\msX := (X/S, \{ \sigma_i \}_{i=1}^r)$ be an $r$-pointed  curve  over $S$.
Then, for an $S$-scheme $S'$,
the collection of data
\index{$\msX$, pointed (stable) curve@--- $S' \times_S\msX$, base-change of}
\begin{align} \label{QQ02}
S' \times_S\msX := ((S' \times_S X)/S', \{(\mr{id}_{S'}, \sigma_i) : S' \migi S' \times_S X \}_{i=1}^r)
\end{align}
forms an $r$-pointed  curve over $S'$.
We refer to it as the {\bf base-change} of $\msX$ to $S'$.
If $\msX$ is prestable (resp.,  of genus $g$), then so is $S' \times_S \msX$.

\subsection{Pointed stable curves} \label{QQ03}
Let
 $\msX:= (f :X \migi S, \{ \sigma_i\}_{i=1}^r)$ be   an $r$-pointed prestable curve of genus $g$ over $S$. Since 
 $f$ is a local complete
intersection morphism, it follows from the theory of duality that the dualizing sheaf $\omega_{X/S}$ on $X$ relative to $S$ exists.
We shall set
\begin{align} \label{QQ3000}
\omega_{X/S}^\mr{log} : = \omega_{X/S}(\sum_{i=1}^r [\sigma_i]),
\index{$\omega_{X/S}^\mr{log}$}
\end{align}
 where $[\sigma_i]$ ($i =1, \cdots, r$) denotes the relative effective reduced divisor on $X$ determined by  the image of $\sigma_i$.
Note that $\omega_{X/S}^\mr{log}$ is a line bundle and its  fiber over any  point of $S$ is  of degree $2g-2+r$.
If $S = \mr{Spec}(\overline{k})$ for an algebraically closed field $\overline{k}$ over $k$, then the following three conditions are equivalent:
 \begin{itemize}
 \item[(a)]
The line bundle  $\omega^\mr{log}_{X/\overline{k}}$ is ample.
\item[(b)]
The automorphism group $\mr{Aut}(\msX)$ of $\msX$ over $\overline{k}$ (i.e., over  the identity morphism of $\mr{Spec}(\overline{k})$)  is finite.
 \item[(c)]
The inequality  $2g-2+r>0$ holds and, for 
every nonsingular rational component $E$ of $X$, 
 the number of points in $E$ contained  in either  the rest components of $X$ or the image of marked points is at least $3$.
   \end{itemize}

\bde
 \label{QQ35} 
  An {\bf $r$-pointed stable curve} ({\bf of genus $g$}) \index{curve@--- stable} over $S$ is an $r$-pointed prestable curve $\msX$ (of genus $g$) over $S$ such that 
  the fiber  $\msX_s := s \times_S \msX$ of $\msX$ over any geometric point $s$  of $S$  
  satisfies one of the  equivalent three conditions (a)-(c) described  above.
   \ede

\subsection{The moduli stack of pointed stable curves} \label{QR780}
Suppose that
$2g-2+r>0$.
Then, $r$-pointed stable curves of genus $g$ over a $k$-scheme and morphisms  between them form a category
\begin{align}
\overline{\mfM}_{g, r}  \migi \mfS \mfc \mfh_{/k}
\index{$\overline{\mfM}_{g,r}$, moduli stack of $r$-pointed stable curves of genus $g$}
\end{align}
 (cf. ~\cite{Kn2}, Definition 1.1) over the category $\mfS \mfc \mfh_{/k}$, where the projection $\overline{\mfM}_{g,r} \migi \mfS \mfc \mfh_{/k}$ is given by $(X/S,  \{\sigma_i\}_{i}) \mapsto S$.
One verifies that $\overline{\mfM}_{g,r}$ forms a category fibered in groupoids over $\mfS \mfc \mfh_{/k}$.
Moreover,  
it may be represented by a geometrically connected, proper, and smooth Deligne-Mumford stack over $k$ of dimension $3g-3+r$ (cf. ~\cite{Kn2}, Corollary 2.6 and Theorem 2.7; ~\cite{DM}, \S\,5).
The stack $\overline{\mfM}_{g,r}$ contains  a dense open substack 
\begin{align}
\mfM_{g,r}  \left(\subseteq \overline{\mfM}_{g,r} \right)
\index{$\mfM_{g,r}$, moduli stack of $r$-pointed smooth curves of genus $g$}
\end{align}
classifying $r$-pointed {\it smooth} curves of genus $g$.
The complement $\overline{\mfM}_{g,r} \setminus \mfM_{g,r}$ defines a divisor with normal crossings.
Denote by  
\begin{align} \label{pp101}
\msC_{g,r} := (f_{\mr{univ}}: \mfC_{g,r} \migi \overline{\mfM}_{g,r}, \{  \sigma_{\mr{univ}, i} : \overline{\mfM}_{g,r} \migi \mfC_{g,r}\}_{i=1}^r)
\index{$\msC_{g,r}$, universal family of curves over $\overline{\mfM}_{g,r}$}
\end{align}
the universal family of $r$-pointed stable curves over $\overline{\mfM}_{g,r}$.
Hence, for each $r$-pointed stable curve $\msX$ of genus $g$ over a $k$-scheme $S$,
there exists a unique (up to isomorphism) morphism $(\phi, \Phi) : \msX \migi \msC_{g,r}$ (cf. (\ref{QQ09})), inducing  an isomorphism $(\mr{id}_S, \Phi) :  \msX \isom S \times_{\phi, \overline{\mfM}_{g,r}} \msC_{g,r}$ over $S$.

\begin{exa}[The case of $g=0$]
Denote by 
$\mbP^1$ \index{$\mbP^1$, projective line} the projective line   over $k$ and by 
$[0]$, $[1]$, and $[\infty]$ the $k$-rational points of   $\mbP^1$ determined by the values $0$, $1$, and $\infty$ respectively.
After ordering  the points $ [0], [1], [\infty]$ suitably, we obtain 
 a unique (up to isomorphism) $3$-pointed stable curve 
\begin{equation} \label{1051}
\msP := (\mbP^1/k, \{ [0], [1], [\infty] \})
\index{$\msP$, $3$-pointed projective line}
\end{equation}
 of genus $0$ over $k$.
 One verifies that $\msP$ does not have nontrivial automorphisms and that  any $3$-pointed stable curve of genus $0$ over a $k$-scheme $S$  is isomorphic to the trivial deformation of $\msP$ over $S$.
That is to say, we have  $\overline{\mfM}_{0, 3} \cong \mfM_{0, 3} \cong \mr{Spec}(k)$.

Moreover, assigning  the fourth point of each $4$-pointed projective line gives
a $k$-isomorphism $\mfM_{0,4} \isom \mbP^1 \setminus \{[0], [1], [\infty] \}$, extending to  $\overline{\mfM}_{0, 4} \isom \mbP^1$.
If $r$ is an integer with $r \geq 4$, then $\overline{\mfM}_{0, r+1}$ can be  constructed by blowing-up a certain locus in $\overline{\mfM}_{0, r} \times_{\overline{\mfM}_{0, r-1}} \overline{\mfM}_{0, r}$.
In particular, $\overline{\mfM}_{0, r}$ ($r \geq 3$) can be represented by a  $k$-scheme  of dimension $r$. 
\end{exa}

\subsection{Log curves} \label{QR781}
Let us introduce the definition of a log curve, which is slightly  different  from    ~\cite{KaFu}, Definition 1.1 (and ~\cite{ACGH},  Definition 4.5) because our definition does not impose the connectedness condition.

\bde \label{pY5019} 
Let $T^\mr{log}$ be  an  fs log scheme 
 (cf. Remark \ref{QQ110} described below). 
A {\bf log curve} \index{log curve} over $T^\mr{log}$ is a  log smooth  integral morphism 
\begin{align} \label{QQ011}
f^\mr{log} : Y^\mr{log} \migi T^\mr{log}
\end{align}
between fs log schemes  (cf. ~\cite{KaKa}, \S\,3)
 such that the geometric fibers of the underlying morphism $f : Y \migi T$  are reduced  $1$-dimensional schemes.
 (In particular, both $\mcT_{Y^\mr{log}/T^\mr{log}}$ and $\Omega_{Y^\mr{log}/T^\mr{log}}$ are line bundle, and  the underlying morphism $f : Y \migi T$ is flat by ~\cite{KaKa}, Corollary 4.5.)
   \ede

\begin{rema} \label{QQ110} 
A log scheme is called {\bf fs}, i.e.,  {\bf fine and saturated} (cf. ~\cite{ILL}, \S\,1.3) \index{fs (= fine and saturated)}
 if it has \'{e}tale locally a chart modeled on a fine (= integral and finitely generated) and saturated monoid.
 The reason for the ``fs'' assumption in Definition \ref{pY5019} is to avoid cusps or worse singularities (cf. ~\cite{ACGH}, Example 4.6).
  \end{rema}

\bde \label{ppy021}
Let $T^\mr{log}$ be an fs log scheme and $Y^\mr{log}$ a log curve over $T^\mr{log}$.
 A {\bf logarithmic chart}, or a {\bf log chart} \index{logarithmic chart (= log chart)} for short,  on $Y^\mr{log}/T^\mr{log}$
is a pair 
\begin{align} \label{QQ010}
(U, \partial)
\end{align}
 consisting of 
an \'{e}tale $Y$-scheme $U$
(where we equip $U$ with a log structure  pulled-back from   $Y^\mr{log}$ and denote the resulting log scheme  by $U^{\mr{log}}$) and an element  $\partial$ of $\Gamma (U, \mcT_{U^\mr{log}/T^\mr{log}})$ 
that generates $\mcT_{U^\mr{log}/T^\mr{log}}$.
We shall say that a log chart $(U, \partial)$ is {\bf global} \index{logarithmic chart (= log chart)@--- global} if $U = Y$.
   \ede

\subsection{The logarithmic structure of a pointed stable curve} \label{QR781}

Recall from  ~\cite{KaFu}, Theorem 4.5, that $\overline{\mfM}_{g,r}$ has a natural log structure given by the divisor at infinity, i.e., the divisor $\overline{\mfM}_{g,r} \setminus \mfM_{g,r}$;
 we shall denote the resulting fs log stack  by 
 \begin{align} \label{QQ013}
 \overline{\mfM}_{g,r}^{\mr{log}}.
 \end{align}
Moreover, by taking the divisor which is the union of the $\sigma_{\mr{univ}, i}$'s and the pull-back of the divisor at infinity of $\overline{\mfM}_{g,r}$, we obtain a log structure on $\mfC_{g,r}$; we denote  the resulting log stack by $\mfC^{\mr{log}}_{g,r}$.
Both $ \overline{\mfM}_{g,r}^{\mr{log}}$ and 
$\mfC^\mr{log}_{g,r}$ are log smooth over $k$.
The structure morphism $f_\mr{univ} : \mfC_{g,r} \migi \overline{\mfM}_{g,r}$ of $\mfC_{g,r}$
    extends naturally to a morphism 
$f^{\mr{log}}_\mr{univ} : \mfC^\mr{log}_{g,r} \migi \overline{\mfM}^\mr{log}_{g,r}$ of log stacks, forming a family of   log curves.

Now, let $S$ be a $k$-scheme  and $\msX := (f: X \migi S, \{ \sigma_i \}_{i=1}^r)$ an $r$-pointed stable curve of genus $g$ over $S$.
Then, $\msX$ determines its classifying morphism $\phi : S \migi \overline{\mfM}_{g,r}$ and  a morphism 
$\Phi : X \migi \mfC_{g,r}$ over $\phi$.
By pulling-back the  log structures of $\overline{\mfM}^{\mr{log}}_{g,r}$ and $\mfC^{\mr{log}}_{g,r}$,
we obtain log structures on $S$ and $X$ respectively; we denote  the resulting schemes  by $S^{\mr{log}}$ and $X^{\mr{log}}$ respectively. 
Both $S^\mr{log}$ and $X^\mr{log}$ are  fs log schemes, and  the morphism $f : X \migi S$ extends to a log curve $f^\mr{log} : X^{\mr{log}} \migi S^{\mr{log}}$
over  $S^\mr{log}$ (cf. ~\cite{KaFu}, Theorem 2.6). 
The line bundle $\omega^\mr{log}_{X/S}$ (cf. (\ref{QQ3000})) is naturally isomorphic to $\Omega_{X^\mr{log}/S^\mr{log}}$.

\begin{rema} \label{QG19}
Recall here the characterization of pointed stable curves in terms of log structures.
Let  $S^\mr{log} := (S, \alpha_S: M_S \migi \mcO_S)$  be an fs log scheme over $k$.
Recall from ~\cite{KaFu}, \S\,1, Definition 1.2, that a log curve $X^\mr{log} := (X, \alpha_X: M_X \migi \mcO_X)$ over $S^\mr{log}$ 
is called a {\bf stable log curve of type $(g,r)$} \index{log curve@--- stable log curve of type $(g,r)$} if, for any log scheme $S_0^{\mr{log}}$ over $S^\mr{log}$ whose underlying scheme $S_0$ is the spectrum of a separably closed field,
the following two conditions are fulfilled:

\begin{itemize} 
\item
The underlying curve $X_0$ of  the fiber $X_0^\mr{log} := S_0^\mr{log} \times_{S^\mr{log}} X^\mr{log}$  is connected and of genus $g$,  and there exists  $r$ pairwise distinct closed points $s_1, \cdots, s_r$ in the smooth locus of $X_0$ such that the relative characteristic $\overline{M}_{X_0^\mr{log}/S^\mr{log}_0}$ (cf. ~\cite{KaFu}, Introduction)
 of $X^\mr{log}_0/S^\mr{log}_0$ satisfies 
\begin{align} \label{QG25}
\overline{M}_{X_0^\mr{log}/S_0^\mr{log}} \cong  \mbN_{s_1} \oplus \cdots \oplus \mbN_{s_r} \oplus \mbZ_{t_1} \oplus \cdots \oplus \mbZ_{t_l}.
\end{align}
Here,  $\{ t_1, \cdots, t_l \}$ is the set of all nodal points on $X_0$, and for a monoid $M$ and a point $x \in X_0$, we denote by $M_x$ the skyscraper sheaf supported only on $x$ by  $M$.
\item
The equality $H^0 (X_0, \mcT_{X_0^\mr{log}/S_0^\mr{log}}) =0$ holds.
\end{itemize}
Suppose that $2g-2+r >0$.
As discussed in ~\cite{KaFu}, the category of $r$-pointed stable curves of genus $g$ is equivalent to the category of stable log curves of type $(g,r)$.
Indeed, 
 let $\overline{\mfN}^\mr{log}$ denote the stack of stable log curves of type $(g,r)$ over $\mfS \mfc \mfh_{/k}^{\mr{log}}$, where  $\mfS \mfc \mfh_{/k}^{\mr{log}}$ denotes  the category of fs log schemes  over $k$.
On the other hand,  $\overline{\mfM}_{g,r}^\mr{log}$ may be considered as a stack over $\mfS \mfc \mfh_{/k}^{\mr{log}}$ in the natural manner.
If we are given
an  fs log scheme $S^\mr{log}$ over $k$ and an $S^\mr{log}$-rational point $s_\mfm^\mr{log} : S^\mr{log} \migi \overline{\mfM}_{g,r}^\mr{log}$,
then the log curve $S^\mr{log} \times_{\overline{\mfM}^\mr{log}_{g,r}} \mfC_{g,r}^\mr{log}$ forms a stable log curve of type $(g,r)$, which determines an $S^\mr{log}$-rational point $s^\mr{log}_\mfn$ of $\overline{\mfN}_{g,r}$.
According to ~\cite{KaFu}, \S\,4, Theorem 4.1, the resulting assignment $s^\mr{log}_\mfm \mapsto s^\mr{log}_\mfn$ defined for  various $S^\mr{log}$'s determines
 an isomorphism
\begin{align} \label{QG24}
\overline{\mfM}_{g,r}^\mr{log} \isom \overline{\mfN}_{g,r}^\mr{log}
\end{align}
of stacks over $\mfS \mfc \mfh_{/k}^\mr{log}$.

\end{rema}

\section{Local pointed curves and monodromy} \label{y030}

\subsection{Local pointed curves} \label{QQ015}
Let $S$ be a $k$-scheme and $r$ a nonnegative integer.
To begin with, let us introduce a local version  of pointed stable curve, as follows. 

\bde \label{QQ970}
\begin{itemize}
\item[(i)]
A {\bf local $r$-pointed  curve} 
\index{$\msU$, local $r$-pointed curve} 
over $S$ is a collection of data
\begin{align} \label{QQ971}
\msU := (U/S, \{ \sigma^U_i :S_i \migi U \}_{i =1}^r),
\end{align}
where $U$ denotes a curve over $S$ and each $\sigma_i^U$ ($i=1, \cdots, r$) denotes a $k$-morphism
from a possibly empty \'{e}tale $S$-scheme  $S_i$ to $U$ satisfying the following condition: 
 there exist   an $r$-pointed  curve $\msX:= (X/S, \{ \sigma_i \}_{i=1}^r)$ over $S$ and an \'{e}tale morphism $u: U \migi X$
such that, for each $i =1, \cdots, r$,  the equality  $S_i = U \times_{u, X, \sigma_i} S$  holds and $\sigma_i^U$ coincides with the base-change of $\sigma_i$ by  $u$.
We shall refer to $\msX$ as a {\bf base} \index{local $r$-pointed  curve@--- base for} for  $\msU$.

\item[(ii)]
We shall say that a local $r$-pointed curve $\msU$ is {\bf stable} \index{local $r$-pointed  curve@--- stable}  if
there exists a base $\msX$ for  $\msU$ that is stable of constant genus.
\end{itemize}
\ede

\begin{rema} \label{QH6}
We regard an $r$-pointed stable curve $\msX$  as a local $r$-pointed stable curve, where a base is taken to be  $\msX$  itself.
\end{rema}

Let $\msU:= (U/S, \{ \sigma^U_i : S_i \migi U \}_{i=1}^r)$ be a local $r$-pointed curve over $S$.
Suppose that $\msU$ is obtained by restricting  an $r$-pointed curve  $\msX$  in the manner of  Definition \ref{QQ970}.
For each $i=1, \cdots, r$, we shall denote by  
$[\sigma^U_i] \left(\subseteq U \right)$
 the scheme-theoretic image of $\sigma^U_i$, which is  possibly empty.
 Then, we obtain   the relative effective divisor 
 \begin{align} \label{QQ979}
 D_\msU \index{$D_\msU$, $D_\msX$}
 \end{align}
on $U$ defined as the union of $[\sigma_i^U]$'s.
 The smooth locus  $U^\mr{sm}$ \index{$U^\mr{sm}$}  of $U \setminus \mr{Supp}(D_\msU)$ relative to  $S$ is a scheme-theoretically dense open subscheme of $U$.

Now, suppose  further that $\msX$ is stable, in particular, both $X$ and $S$ are equipped with the natural log structures as  discussed before.
By pulling-back the log structure of $X^\mr{log}$, one may equip $U$ with a log structure;
if $U^\mr{log}$ denotes the resulting log scheme, then
 we obtain a log curve $U^\mr{log}/S^\mr{log}$.
 The log structures of $U^\mr{log}$ and $S^\mr{log}$  depend on the choice of $\msX$, but the isomorphism class of $\Omega_{U^\mr{log}/S^\mr{log}}$, as well as of $\mcT_{U^\mr{log}/S^\mr{log}}$, is independent of that choice.
Moreover, the restriction of the log structure of $U^\mr{log}$ to  $U^\mr{sm}$ coincides with the log structure  pulled-back from $S^\mr{log}$.

Note that, for each $i=1, \cdots, r$,  there exists a canonical isomorphism
\begin{align} \label{QQ040}
\mr{Res}_{\msU, i} : \sigma_i^{U*}(\Omega_{U^\mr{log}/S^\mr{log}}) \isom \mcO_{S_i} \index{$\mr{Res}_{\msU, i}$, residue map}
\end{align}
i.e., the  {\bf residue map}, which maps any local section of the  form $d\mr{log} (x) \left(:=\frac{1}{x}  \cdot dx\right) \in \sigma^{U*}_i (\Omega_{U^\mr{log}/S^\mr{log}})$ to $1 \in \mcO_{S_i}$, where $x$ denotes a local function defining the closed subscheme $[\sigma_i^U]$ of $U$.

\subsection{Pull-back} \label{Q989}
Let $\msU := (U/S, \{ \sigma_i^U : S_i \migi U \}_{i=1}^r)$ be a local $r$-pointed  (stable) curve over $S$.
Then, for an \'{e}tale $U$-scheme $U'$, the collection of data
\index{$\msU$, local $r$-pointed curve@--- $U' \times_U \msU$, pull-back of $\msU$}
\begin{align} \label{QQ034}
U' \times_U \msU := (U'/S, \{ \sigma^U_i \times \mr{id}_{U'} :  S_i \times_U U' \migi U'\}_{i=1}^r) 
\end{align}
forms a local $r$-pointed (stable) curve over $S$.
We refer to it as the {\bf pull-back} of $\msU$ to $U'$.

\subsection{Monodromy of a connection} \label{QQ020}
Let $\msU := (U/S, \{ \sigma_i : S_i \migi U \}_{i=1}^r)$ be a local $r$-pointed  stable curve over $S$.
After choosing a base for $\msU$, we obtain a log curve $U^\mr{log}/S^\mr{log}$ extending $U/S$.
 In what follows, we suppose that $r \geq 1$ and fix  an element $i \in \{1, \cdots, r \}$ with $[\sigma_i^U] \neq \emptyset$.
 Let $\pi : \mcE \migi U$  be a $\mbG$-bundle   on $U$.
 As discussed in \S\,\ref{QR300}, we obtain a log scheme 
  $\mcE^\mr{log}$ defined as $\mcE$ equipped with the log structure pulled-back from $U^\mr{log}$.
Both $\mcT_{U^\mr{log}/S^{\mr{log}}}$ and $\mcT_{\mcE^\mr{log}/S^{\mr{log}}}$ are locally free, and 
the  natural morphisms
\begin{align} \label{EEE01} \iota_U : \mcT_{U^\mr{log}/S^{\mr{log}}} \migi \mcT_{U/S}, \ \ \ \iota_\mcE : \mcT_{\mcE^\mr{log}/S^{\mr{log}}} \migi \mcT_{\mcE/S}.
\end{align}
 are isomorphisms over the  scheme-theoretically dense open subscheme $U^\mr{sm}$ of $U$.
It follows that  
these morphisms are injective over the entire space $U$.
The direct image  $\pi_*(\mcT_{\mcE^\mr{log}/S^{\mr{log}}}) \migiincl \pi_*(\mcT_{\mcE/S})$ of  $\iota_\mcE$ is compatible with the   $\mbG$-actions.
Hence, it  restricts to an injection 
\begin{equation} \widetilde{\iota}_\mcE :\widetilde{\mcT}_{\mcE^\mr{log}/S^\mr{log}} \migiincl \widetilde{\mcT}_{\mcE/S} \  \left(: = \pi_*(\mcT_{\mcE/S})^\mbG \right),\end{equation}
fitting  into the following
 morphism of short exact sequences:
\begin{align} \label{QW6200}
\vcenter{\xymatrix@C=26pt@R=36pt{
0 \ar[r] 
&\mfg_\mcE  \ar[r]\ar[d]^-{ \mr{id}}
&\widetilde{\mcT}_{\mcE^\mr{log}/S^\mr{log}}   \ar[r]^-{d_\mcE^\mr{log}}\ar[d]^-{\widetilde{\iota}_\mcE}
&  \mcT_{U^\mr{log}/S^\mr{log}}  \ar[r] \ar[d]^-{\iota_U} & 0 \\
0 \ar[r]&\mfg_\mcE  \ar[r]
&\widetilde{\mcT}_{\mcE/S}  \ar[r]^-{d_\mcE} 
& \mcT_{U/S} \ar[r] & 0, 
}}
\end{align}
where the lower horizontal sequence is   (\ref{Ex0}) for  the non-logarithmic schemes $S$, $U$, and $\mcE$.
Also,  by applying, to the lower horizontal sequence, the natural transformation 
 $(-) \migi \sigma^U_{i*}(\sigma^{U*}_i(-))$ arising from  the adjunction relation ``$\sigma_i^{U*}(-) \dashv \sigma^U_{i*}(-)$''
(i.e., ``the functor $\sigma^{U*}_i(-)$ is left adjoint to the functor $\sigma^{U}_{i*}(-)$''), we obtain  a morphism of sequences
 \begin{align} \label{QW6201}
\vcenter{\xymatrix@C=16pt@R=36pt{
0 \ar[r] & \mfg_\mcE  \ar[r]\ar[d]
&\widetilde{\mcT}_{\mcE/S}  \ar[r]^-{d_\mcE}\ar[d]
& \mcT_{U/S} \ar[r] \ar[d] & 0 \\
0 \ar[r]& \sigma^U_{i*}(\sigma_i^{U*}(\mfg_\mcE))   \ar[r]
& \sigma^U_{i*}(\sigma_i^{U*}(\widetilde{\mcT}_{\mcE/S})) \ar[r] 
& \sigma^U_{i*}(\sigma_i^{U*}(\mcT_{U/S})) \ar[r] & 0.
}}
\end{align}
The above  two horizontal sequences  are    exact because of the local triviality of the  $\mbG$-bundle $\mcE$ with respect to  the \'{e}tale topology.
Composing (\ref{QW6200}) and (\ref{QW6201}) gives  the following  morphism of sequences: 
 \begin{align} \label{diag81}
\vcenter{\xymatrix@C=16pt@R=36pt{
0 \ar[r] & \mfg_\mcE \ar[r]\ar[d]
& \widetilde{\mcT}_{\mcE^\mr{log}/S^\mr{log}}  \ar[r]^-{d_\mcE^\mr{log}}\ar[d]^-{\widetilde{\iota}'_\mcE}
&\mcT_{U^\mr{log}/S^\mr{log}}  \ar[r] \ar[d]^-{\iota'_{U} } & 0 \\
0 \ar[r]&\sigma^U_{i*}(\sigma_i^{U*}(\mfg_\mcE))  \ar[r]
& \sigma^U_{i*}(\sigma_i^{U*}(\widetilde{\mcT}_{\mcE/S})) \ar[r] 
& \sigma^U_{i*}(\sigma_i^{U*}(\mcT_{U/S})) \ar[r] & 0. 
}}
\end{align}
Since the natural surjection $\mcT_{U/S} \migisurj \mr{Coker}(\iota_U)$ induces an isomorphism
\begin{equation} \sigma^U_{i*}(\sigma_i^{U*}(\mcT_{U/S})) \cong \sigma^U_{i*}(\sigma_i^{U*}( \mr{Coker}(\iota_U))),
\end{equation}
the right-hand vertical arrow $\iota'_U$ in (\ref{diag81}) turns out to be the zero map.

Now, let us take  an $\hslash$-$S^\mr{log}$-connection $\nabla$ on $\mcE$.
It follows from the above discussion that the composite 
\begin{equation}    
 \widetilde{\iota}'_\mcE   \circ  \nabla   :  \mcT_{U^\mr{log}/S^\mr{log}} \migi \sigma^U_{i*}(\sigma_i^{U*}(\widetilde{\mcT}_{\mcE/S})) \end{equation}
factors through the injection $\sigma^U_{i*}(\sigma_i^{U*}(\mfg_\mcE)) \migiincl \sigma^U_{i*}(\sigma_i^{U*}(\widetilde{\mcT}_{\mcE/S}))$.
The resulting morphism
$\mcT_{U^\mr{log}/S^\mr{log}} \migi \sigma^U_{i*}\sigma_i^{U*}(\mfg_\mcE)$ corresponds, via the adjunction relation ``$\sigma_i^{U*}(-) \dashv \sigma^U_{i*}(-)$", to a morphism 
\begin{equation} \label{mono}  \sigma^{U*}_i(\mcT_{U^\mr{log}/S^\mr{log}}) \migi \sigma_i^{U*}(\mfg_\mcE). \end{equation}
The image of  $1 \in \Gamma(S_i, \mcO_{S_i})$ via the composite
 \begin{equation} \label{QQ066}
  \mcO_{S_i} \xrightarrow{\mr{Res}_{\msU, i}^\vee} \sigma_i^{U*}(\mcT_{U^\mr{log}/S^\mr{log}}) \xrightarrow{(\ref{mono})}  \sigma_i^{U*}(\mfg_\mcE)
  \end{equation}
determines  a global section 
\begin{equation} \label{monodromy}
 \mu_i^{\nabla} \in \Gamma(S_i, \sigma_i^{U*}(\mfg_\mcE)) \ \left( = \Gamma (S_i, \mfg_{\sigma_i^{U*}(\mcE)})\right).
\end{equation}

\bde \label{y031} 
 We shall refer to $\mu_i^{\nabla}$ as the {\bf monodromy (operator)}  
 \index{$\nabla$, $\hslash$-$T^\mr{log}$-connection@--- $\mu^\nabla_i$, monodromy (operator)}
 of $\nabla$ at  $\sigma_i^U$.
\ede


\begin{exa}\label{QQ1011}
Let us consider the case where $\mcE = U \times_k \mbG$, i.e., $\mcE$ is the trivial $\mbG$-bundle.
Suppose that there exists a  local function $x \in \mcO_U$ defining $[\sigma_i^U] \left(\subseteq U \right)$.
If  $\partial_x$  denotes the local  section of $\mcT_{U^\mr{log}/S^\mr{log}}$
defined as the dual of $dx \in \Omega_{U^\mr{log}/S^\mr{log}}$, then
the section $x \cdot \partial_x$ locally  generates  $\mcT_{U^\mr{log}/S^\mr{log}}$.
For each $S^\mr{log}$-connection $\nabla$ on $U \times_k \mbG$,  we have 
\begin{align} \label{QQ1012}
\mu_i^\nabla = \nabla^\mfg (x \cdot \partial_x) |_{\sigma_i}. 
\end{align}
\end{exa}

\subsection{Change of structure group} \label{QQ070}
Let $\mbG'$ and $w$ be as in \S\,\ref{QQ119}, and denote by $\mfg'$ the Lie algebra of $\mbG'$.
Also, let $(\mcE, \nabla)$ be an $\hslash$-flat $\mbG$-bundle on $U^\mr{log}/S^\mr{log}$.  
As discussed in \S\,\ref{QQ809}, $\nabla$ induces  an $\hslash$-$S^\mr{log}$-connection $w_*(\nabla)$  on the $\mbG'$-bundle $\mcE \times^\mbG \mbG'$.
Then, the monodromy $\mu_i^{w_*(\nabla)}$ of $w_*(\nabla)$ at $\sigma_i^U$  ($i=1, \cdots, r$) satisfies 
\begin{align}
\mu_i^{w_*(\nabla)} = \mr{Ad}_w (\mu_i^\nabla),
\end{align}
where $\mr{Ad}_w$ denotes the $\mcO_U$-linear morphism $\mfg_{\sigma^{U*}_i(\mcE)} \migi \mfg'_{\sigma^{U*}_i(\mcE)} \left(= \mfg'_{\sigma^{U*}_i(\mcE \times^\mbG \mbG')} \right)$ induced naturally by $w$.

\subsection{Change of parameter} \label{QQ1034}
Let 
$(\mcE, \nabla)$ be as above 
 and 
 $c$  an element of $\Gamma (S, \mcO_S)$.
As discussed in \S\,\ref{QQ1030}, $\nabla$ induces a $(c \cdot \hslash)$-$S^\mr{log}$-connection $c \cdot \nabla$ on $\mcE$.
Then, the monodromy $\mu_i^{c \cdot \nabla}$ of $c \cdot \nabla$ at $\sigma_i^U$ ($i=1, \cdots, r$) satisfies 
\begin{align} \label{QQ1039}
\mu_i^{c\cdot \nabla} =  c \cdot \mu_i^\nabla.
 \end{align}

\subsection{Gauge transformation} \label{QQ073}
Also, let $\eta$ be an automorphism of the $\mbG$-bundle $\mcE$.
Then, 
the monodromy $\mu_i^{\eta_{*} (\nabla)}$  at $\sigma_i^U$ ($i=1, \cdots, r$) of  the gauge transformation $\eta_{*} (\nabla)$ by $\eta$ (cf. Definition \ref{QQ90}) 
satisfies 
\begin{align} \label{QR0111}
\mu_i^{\eta_*(\nabla)} =  d\eta^\mbG |_{\sigma_i^U} (\mu^{\nabla}_i)
\end{align}
(cf. (\ref{QQ10}) for the definition of $d\eta^\mbG$).
In particular, if $\mcE = U \times_k \mbG$ and $\eta = \mr{L}_h$ for a morphism of $k$-schemes $h: U \migi \mbG$,
then it follows from Lemma \ref{QR349} that
\begin{align} \label{QQ1013}
\mu_i^{\mr{L}_{h*}(\nabla)} = \mr{Ad}_\mbG (h \circ \sigma_i^U) (\mu_i^\nabla) \left(= \mr{Ad}_\mbG (h) |_{\sigma^U_i} (\mu_i^\nabla) \right).
\end{align}



\chapter{Opers on a family of pointed stable curves}  \label{y032}

In this chapter, we study a general theory of $(\mfg, \hslash)$-opers on pointed stable curves and their moduli for  a semisimple Lie algebra $\mfg$.
Compared to previous studies on opers, our  formulation would be  fairly general   because  it especially includes the case of nodal curves.
The goal of this chapter is to prove the representability of the moduli category classifying  $(\mfg, \hslash)$-opers.
This result generalizes results by R. C. Gunning (cf. ~\cite{G2}), A. Beilinson-V. Drinfeld (cf. ~\cite{BD2}), and S. Mochizuki (cf. ~\cite{Mzk1}). 
Also, a key ingredient in studying opers on a pointed curve is an invariant called {\it radius} (cf. Definition \ref{y056}).
The radius of a given $(\mfg, \hslash)$-oper at a marked point is defined as an element of the adjoint quotient of $\mfg$ and  induced by the monodromy.
We will also prove the representability of the moduli category of $(\mfg, \hslash)$-opers of prescribed radii.

\begin{notation}
Let $k$
\index{$k$, field} 
 be a perfect field and $(\mbG, \mbT)$ 
\index{$\mbG$, algebraic group}
\index{$\mbG$, algebraic group@— $\mbT$, maximal torus of}
 a split semisimple algebraic group over $k$ of adjoint type, where $\mbT$ denotes a maximal torus of $\mbG$.
Let us fix 
a Borel subgroup $\mbB$ 
\index{$\mbG$, algebraic group@— $\mbB$, Borel subgroup of}
defined over $k$  containing $\mbT$.
Write $\mbN := [\mbB, \mbB]$, i.e., the unipotent radical of $\mbB$.
\index{$\mbG$, algebraic group@— $\mbN$}
Also, write $\mfg$, 
\index{$\mfg$, Lie algebra}
$\mfb$, $\mfn$, and $\mft$ for    the  Lie algebras of $\mbG$,  $\mbB$, $\mbN$, and $\mbT$ respectively (hence $\mft \subseteq \mfn \subseteq \mfb \subseteq \mfg$).

\end{notation}

\section{$(\mfg, \hslash)$-opers on a log curve} \label{y034}

\subsection{The Cartan decomposition}\label{QR21}
Let  \index{$\Gamma$, set of simple roots}$\Gamma \subseteq \mr{Hom}(\mbT, \mbG_m)$  (:= the additive group of all characters of $\mbT$)
 denote the set of  simple  roots in $\mbB$ with respect to $\mbT$.
For each $\beta \in \mr{Hom}(\mbT, \mbG_m)$, we set
\index{$\mfg$, Lie algebra@--- $\mfg^\beta$}
\begin{equation} 
\mfg^\beta := \big\{ x \in \mfg \ \big| \ \text{$\mr{Ad}_\mbG(t)(x) =  \beta (t) \cdot x$ for all  $t\in \mbT$}  \big\}.
\end{equation}
Recall  that there exists a unique  Lie algebra grading 
\index{$\mfg$, Lie algebra@--- $\mfg_j$, $\mfg^j$}
\begin{equation} \mfg = \bigoplus_{ j \in \mbZ}\mfg_j
\label{grading}
\end{equation}
(i.e.,  $[\mfg_{j_1}, \mfg_{j_2}] \subseteq \mfg_{j_1 +j_2}$ for $j_1$, $j_2 \in \mbZ$)
satisfying the following conditions: 
\begin{itemize}
\item
$\mfg_j = \mfg_{-j} = 0$ for all $j > \mr{rk}(\mfg)$, where $\mr{rk}(\mfg)$ denotes the rank of $\mfg$; 
\index{$\mfg$, Lie algebra@--- $\mr{rk}(\mfg)$,  rank of}
\item
$\mfg_0 = \mft$, $\mfg_1=\bigoplus_{\alpha \in \Gamma}\mfg^\alpha$ and $\mfg_{-1}= \bigoplus_{\alpha \in \Gamma}\mfg^{-\alpha}$.
\end{itemize}
The associated decreasing filtration  $\{ \mfg^j \}_{j \in \mbZ}$ on $\mfg$ is  defined by 
\index{$\mfg$, Lie algebra@--- $\mfg_j$, $\mfg^j$}
\begin{equation} \mfg^j:= \bigoplus_{l \geq j}\mfg_l,  
\label{filtration}
\end{equation}
which  is closed under the adjoint action of $\mbB$
 and satisfies the following conditions: 
 \begin{itemize}
\item
$\mfg^{-j} = \mfg$ and $\mfg^{j} = 0$ for all $j > \mr{rk}(\mfg)$;
\item
$\mfg^0 = \mfb$, $\mfg^1= \mfn$, and $[\mfg^{j_1}, \mfg^{j_2}] \subseteq \mfg^{j_1 +j_2}$ for $j_1$, $j_2 \in \mbZ$.
\end{itemize}

\subsection{A filtration on a Lie algebroid} 
\label{QR29}
By means of the filtration $\{ \mfg^j \}_{j \in \mbZ}$ defined above, we  can construct a filtration on the Lie algebroid associated to a $\mbG$-bundle  equipped with a   $\mbB$-reduction.

Let 
$S^\mr{log}$ be an fs log scheme 
 over $k$,  $U^\mr{log}$ a log curve over $S^\mr{log}$,  
and
$\pi_\mbB : \mcE_\mbB  \migi U$  a  $\mbB$-bundle on  $U$.
Denote by
$\pi_\mbG : (\mcE_\mbB \times^\mbB \mbG =:) \ \mcE_\mbG \migi U$ the  $\mbG$-bundle obtained via 
 change of structure group by  the inclusion $\mbB \migiincl \mbG$.
 Hence, $\mcE_\mbB$ specifies a $\mbB$-reduction of  $\mcE_\mbG$.
The natural morphism $\mcE_\mbB \migi \mcE_\mbG$ yields 
an isomorphism $\mfg_{\mcE_\mbB} \isom \mfg_{\mcE_\mbG }$.
Moreover, it gives 
the following  injective morphism of short exact sequences:
\index{$\iota_{\mfg/\mfb}$, $\widetilde{\iota}_{\mfg/\mfb}$}
\begin{align} \label{extension} 
\vcenter{\xymatrix@C=26pt@R=36pt{
0 \ar[r] &\mfb_{ \mcE_\mbB}  \ar[r]\ar[d]^-{\iota_{\mfg/\mfb}}
&  \widetilde{\mcT}_{\mcE_\mbB^\mr{log}/S^\mr{log}} \ar[r]^-{d_{\mcE_\mbB}^\mr{log}}\ar[d]^-{\widetilde{\iota}_{\mfg/\mfb}}
& \mcT_{U^\mr{log}/S^\mr{log}}  \ar[r] \ar[d]_-{\wr}^-{\mr{id} } & 0 \\
0 \ar[r]& \mfg_{\mcE_\mbG} \ar[r]
& \widetilde{\mcT}_{\mcE_\mbG^\mr{log}/S^\mr{log}} \ar[r]_-{d_{\mcE_\mbG}^\mr{log}} 
&\mcT_{U^\mr{log}/S^\mr{log}}  \ar[r] & 0. 
}}
\end{align}
Each  $\mfg^j \ (\subseteq \mfg)$  is closed under the adjoint action of $\mbB$, so  induces 
a subbundle $\mfg^j_{\mcE_\mbB}$ of $\mfg_{\mcE_\mbB} \left(= \mfg_{\mcE_\mbG} \right)$.
The collection 
$\{ \mfg^j_{\mcE_\mbB} \}_{j \in \mbZ}$
 forms a decreasing filtration on $\mfg_{\mcE_\mbB}$.
Moreover, for each $j \in \mbZ$, we write
\begin{align} \label{QQ0100}
\widetilde{\mcT}^j_{\mcE^\mr{log}_\mbG/S^\mr{log}} := \mr{Im}(\widetilde{\iota}_{\mfg/\mfb}) + \mfg_{\mcE_\mbB}^j \left(\subseteq \widetilde{\mcT}_{\mcE^\mr{log}_\mbG/S^\mr{log}} \right).
\end{align}
Then, the collection  $\{ \widetilde{\mcT}^j_{\mcE_\mbG^\mr{log}/S^\mr{log}}\}_{j \in \mbZ}$ forms a  decreasing filtration on the vector bundle  $\widetilde{\mcT}_{\mcE_\mbG^\mr{log}/S^\mr{log}}$.
The equality  $\widetilde{\mcT}^j_{\mcE^\mr{log}/S^\mr{log}} = \widetilde{\mcT}^0_{\mcE^\mr{log}/S^\mr{log}}$ holds for every $j \geq 0$,  and   the natural inclusion $\mfg_{\mcE_\mbB} \left(=\mfg_{\mcE_\mbG}\right) \migiincl \widetilde{\mcT}_{\mcE_\mbG^\mr{log}/S^\mr{log}}$ induces an isomorphism
\begin{align}  \label{grading2} 
 \mfg^{j-1}_{\mcE_\mbB}/ \mfg^{j}_{\mcE_\mbB} \isom \widetilde{\mcT}_{\mcE_\mbG^\mr{log}/S^\mr{log}}^{j-1}/ \widetilde{\mcT}_{\mcE_\mbG^\mr{log}/S^\mr{log}}^{j}
\end{align}
for every  $j \leq 0$.
Note that  each $\mfg^{-\alpha}$ ($\alpha \in \Gamma$) is closed under the $\mbB$-action  defined by  the composite $\mbB \migisurj \left(\mbB/\mbN = \right)\mbT \xrightarrow{\mr{Ad}_\mbG} \mr{Aut} (\mfg^{-\alpha})$.
Hence,   the  decomposition 
 ($\mfg_{-1} =$) $\mfg^{-1}/\mfg^0 = \bigoplus_{\alpha \in \Gamma} \mfg^{-\alpha}$
 gives  a canonical decomposition 
\begin{equation} \label{decom97} \widetilde{\mcT}_{\mcE_\mbG^\mr{log}/S^\mr{log}}^{-1} / \widetilde{\mcT}^0_{\mcE_\mbG^\mr{log}/S^\mr{log}} = \bigoplus_{\alpha \in \Gamma}    \mfg^{-\alpha}_{\mcE_\mbB} .\end{equation}


\subsection{The definition of a $(\mfg, \hslash)$-oper}

We now turn to define a $(\mfg, \hslash)$-oper on a log curve.
Let $\hslash$ be an element of $\Gamma (S, \mcO_S)$.

\bde \label{y035} 
\begin{itemize}
\item[(i)]
 Let  us consider a pair 
 \index{$\msE^\spadesuit$, $(\mfg, \hslash)$-oper}
 \begin{equation} \msE^\spadesuit := 
 (\mcE_\mbB, \nabla)
 \end{equation}
consisting of a  $\mbB$-bundle  $\mcE_\mbB$ on $U$ and an $\hslash$-$S^\mr{log}$-connection $\nabla$
 on the $\mbG$-bundle 
 $\mcE_\mbG \left(:= \mcE_\mbB \times^\mbB \mbG \right)$
   induced by $\mcE_\mbB$.
Then, we shall say that 
$\msE^\spadesuit$
  is 
a {\bf  $(\mfg, \hslash)$-oper} 
on  $U^\mr{log}/S^\mr{log}$
  if it satisfies the following two conditions:
\begin{itemize}
\item
 $\nabla (\mcT_{U^\mr{log}/S^\mr{log}}) \subseteq  \widetilde{\mcT}_{\mcE_\mbG^\mr{log}/S^\mr{log}}^{-1}$;
\item
For any $\alpha \in \Gamma$, the composite
\begin{equation}
\label{isomoper}
 \mcT_{U^\mr{log}/S^\mr{log}} \xrightarrow{\nabla} \widetilde{\mcT}_{\mcE_\mbG^\mr{log}/S^\mr{log}}^{-1} \migisurj \widetilde{\mcT}_{\mcE^\mr{log}_\mbG/S^\mr{log}}^{-1} /\widetilde{\mcT}^0_{\mcE_\mbG^\mr{log}/S^\mr{log}} \migisurj   \mfg^{-\alpha}_{\mcE_\mbB}
  \end{equation}
is an isomorphism, where the third arrow denotes the natural projection with respect to the   decomposition (\ref{decom97}).
\end{itemize}
For simplicity, a $(\mfg, 1)$-oper will be  referred to as a {\bf $\mfg$-oper}.
Next, suppose further that $U^\mr{log}/S^\mr{log}$ arises from a local pointed stable curve $\msU$ together with a choice of its base  (e.g., $U^\mr{log}= X^\mr{log}$ for a pointed stable curve $\msX := (X/S, \{ \sigma_i \}_i)$).
Then,
 we shall refer to any  $(\mfg, \hslash)$-oper on $U^\mr{log}/S^\mr{log}$ as a {\bf $(\mfg, \hslash)$-oper on $\msU$} if there is no fear of confusion.

\item[(ii)]
Let $\msE^\spadesuit := (\mcE_\mbB, \nabla)$, $\msE'^{\spadesuit} := (\mcE'_\mbB, \nabla')$ be $(\mfg, \hslash)$-opers on  $U^\mr{log}/S^\mr{log}$.
Then, an {\bf isomorphism of $(\mfg, \hslash)$-opers} 
 from $\msE^\spadesuit$ to $\msE'^\spadesuit$ is  defined as an isomorphism $\mcE_\mbB \isom \mcE'_\mbB$ of  $\mbB$-bundles  that induces, via change of structure group by $\mbB \migiincl \mbG$,
   an isomorphism of $\hslash$-flat $\mbG$-bundles $(\mcE_\mbG, \nabla) \isom (\mcE'_\mbG, \nabla')$.
\end{itemize}
  \ede

\begin{exa}[The case of $\mfg = \mfs \mfl_2$] \label{QR0203}
Suppose that  $2$ is invertible in $k$.
In the case of $\mfg = \mfs \mfl_2 =\mfp \mfg \mfl_2$,
we can describe the  definition of an $(\mfs \mfl_2, \hslash)$-oper  a little simpler.
As we will define later  (cf. (\ref{BorelB})), 
let $T$ (resp., $B$)  denote  the maximal torus (resp., the Borel subgroup) of $\mr{PGL}_2$ consisting of  the images, via $\mr{GL}_2 \migisurj \mr{PGL}_2$, of invertible diagonal matrices (resp., invertible upper triangular matrices).
Then,
an $(\mfs \mfl_2, \hslash)$-oper on $U^\mr{log}/S^\mr{log}$ with respect to the pair of choices $(T, B)$ is given  as  
an $\hslash$-flat $\mr{PGL}_2$-bundle $(\mcE_{\mr{PGL}_2}, \nabla)$ on  $U^\mr{log}/S^\mr{log}$ equipped with a  $B$-reduction $\mcE_B$ on $\mcE_{\mr{PGL}_2}$ such that
the composite
\begin{align} \label{QR0205}
\mcT_{U^\mr{log}/S^\mr{log}} \xrightarrow{\nabla} \widetilde{\mcT}_{\mcE_{\mr{PGL}_2}^\mr{log}/S^\mr{log}} \migisurj 
\widetilde{\mcT}_{\mcE_{\mr{PGL}_2}^\mr{log}/S^\mr{log}}/
\widetilde{\mcT}_{\mcE_{B}^\mr{log}/S^\mr{log}}
\end{align}
is an isomorphism.
If we choose a log chart $(V, \partial)$ on $U^\mr{log}/S^\mr{log}$ and a local trivialization $V \times_k \mr{PGL}_2 \isom \mcE_{\mr{PGL}_2} |_V$ of $\mcE_{\mr{PGL}_2}$, inducing $V \times_k B \isom \mcE_B$,
the element $\nabla^{\mfg} (\partial) \in \mfp \mfg \mfl_n (V)$ can be represented by a $2 \times 2$ matrix of the form 
$\nabla^{\mfg} (\partial) = \begin{pmatrix}  * & *  \\ u & * \end{pmatrix}$, where the $*$'s indicate arbitrary functions on $V$ and $u$ is a nowhere vanishing function, i.e., an element of $\mbG_m (V)$.
We refer to  Example \ref{RRRH} for the local description of an $(\mfs \mfl_n, \hslash)$-oper for a general $n$.
\end{exa}

\subsection{A grading of the adjoint bundle}
\label{y037} 
Let $\msE^\spadesuit := (\mcE_\mbB, \nabla_\mcE)$ be a $(\mfg, \hslash)$-oper on  $U^\mr{log}/S^\mr{log}$.
Consider the composite isomorphism $\mcT_{U^\mr{log}/S^\mr{log}}\isom \mfg^{-\alpha}_{\mcE_\mbB}$ (for each  $\alpha \in \Gamma$) displayed in   (\ref{isomoper}) and the decomposition  $\mfg = \bigoplus_{j \in \mbZ} \mfg_j$ (cf. (\ref{grading})). 
Then,  each subquotient $\mfg^{j}_{\mcE_{\mbB}}/\mfg^{j+1}_{\mcE_{\mbB}}$ 
of the filtration $\{\mfg^j_{\mcE_{\mbB}} \}_{j \in \mbZ}$
 is isomorphic to a dirent sum of finite copies of $\Omega_{U^\mr{log}/S^\mr{log}}^{\otimes j}$, i.e.,
 \begin{align} \label{QQ0104}
 \mfg^{j}_{\mcE_{\mbB}}/\mfg^{j+1}_{\mcE_{\mbB}} \cong \Omega_{U^\mr{log}/S^\mr{log}}^{\otimes j} \otimes_k \mfg_j.
 \end{align}
In particular, we have $\mfg^{0}_{\mcE_{\mbB}}/\mfg^{1}_{\mcE_{\mbB}} \cong \mcO_{U}^{\oplus \mr{rk}(\mfg)}$.


\subsection{Pull-back and base-change} 
\label{QQ46}
Let $\msE^\spadesuit := (\mcE_\mbB, \nabla_\mcE)$ be a $(\mfg, \hslash)$-oper on $U^\mr{log}/S^\mr{log}$, and 
suppose that we are given  
an \'{e}tale $U$-scheme $u: U' \migi U$.
Then, $U'$  is equipped with a log structure pulled-back from $U^\mr{log}$;
if we denote the resulting log scheme by $U'^{\mr{log}}$, then it defines a log curve over $S^\mr{log}$.
One verifies that
the pair of pull-backs 
\index{$\msE^\spadesuit$, $(\mfg, \hslash)$-oper@--- $u^*(\msE^\spadesuit)$, pull-back of}
\begin{align} \label{QQ0141}
u^*(\msE^\spadesuit) := (u^*(\mcE_\mbB), u^*(\nabla))
\end{align}
 (cf. (\ref{QR240}), (\ref{QR312}))
 specifies a $(\mfg, \hslash)$-oper on $U'^{\mr{log}}/S^\mr{log}$.
 Conversely,  the formation of a $(\mfg, \hslash)$-oper  has descent with respect to the \'{e}tale topology on $U$.

Also, let $S'^{\mr{log}}$ be an fs log scheme and 
$s^\mr{log} : S'^{\mr{log}} \migi S^\mr{log}$  a morphism of log schemes.
Then, the pair of base-changes 
\index{$\msE^\spadesuit$, $(\mfg, \hslash)$-oper@--- $s^*(\msE^\spadesuit)$, base-change}
\begin{align} \label{QQ866}
s^*(\msE^\spadesuit) := (s^*(\mcE), s^*(\nabla))
\end{align}
(cf. (\ref{QR241}), (\ref{QQ100})) forms a $(\mfg, s^*(\hslash))$-oper on  the log curve $(S'^{\mr{log}} \times_{S^\mr{log}}U^\mr{log})/S'^{\mr{log}}$.


\subsection{Change of parameter}
\label{QQ1070}
Let $\msE^\spadesuit := (\mcE_\mbB, \nabla)$ be as above
 and let $c \in \Gamma (S, \mcO_S^\times)$.
Then, the pair 
\index{$\msE^\spadesuit$, $(\mfg, \hslash)$-oper@--- $c \cdot \msE^\spadesuit$}
\begin{align} \label{QQ1078}
c \cdot \msE^\spadesuit := (\mcE_\mbB, c \cdot \nabla) 
\end{align}
(cf. (\ref{QQ1032})) specifies a $(\mfg,  c\cdot \hslash)$-oper on $U^\mr{log}/S^\mr{log}$.

\section{Local description  and automorphisms of a $(\mfg, \hslash)$-oper} \label{QQ50}

Hereafter, let  $S^\mr{log}$ be an fs log scheme over $k$,  $U^\mr{log}$ a log curve  over $S^\mr{log}$, and $\hslash$ an element of $\Gamma (S, \mcO_S)$.
Also, let us assume in this section 
that 
{\it  either ``\,$\mr{char}(k)=0$'' or ``\,$\mr{char}(k)> 2  h_\mbG$'' or  ``\,$\mbG = \mr{PGL}_n$ with $1 < n < \mr{char}(k)$'' is fulfilled}.

\subsection{Generators of root spaces} \label{QS8375}

We shall fix a collection of data 
\index{$\ \qq\,$}
\begin{align} \label{QG5697}
\qq := \{ x_\alpha \}_{\alpha\in \Gamma}, 
\end{align}
 where each $x_\alpha$ denotes a generator of $\mfg^\alpha$.
(If we take two such collections, then they   are, in an evident sense,  conjugate by an element of $\mbT (k)$.
For that reason, the results obtained in the following  discussion are  essentially independent of the choice of $\ \qq$.)
In particular, the pair $(\mbB, \ \qq)$ specifies a {\it pinning} of the split semisimple group $(\mbG, \mbT)$.
We  set
\begin{equation}  \label{sss01}  q_{1} := \sum_{\alpha \in \Gamma} x_\alpha, \ \ \ \ \check{\rho} := \sum_{\alpha \in \Gamma} \check{\omega}_\alpha, 
\index{$q_{1}$, $q_{-1}$} \index{$\check{\rho}$} \end{equation}
where $\check{\omega}_\alpha$ denotes the fundamental coweight of $\alpha$, considered as an element of $\mft$ via differentiation.
Then,  there is  a unique  collection $\{ y_\alpha \}_{\alpha \in \Gamma}$, where $y_\alpha$ is a generator of $\mfg^{-\alpha}$, such that if we write
\index{$q_{1}$, $q_{-1}$}
\begin{equation} \label{q_{-1}}
q_{-1}  := \sum_{\alpha \in \Gamma} y_\alpha \in \mfg_{-1},
\end{equation} 
then the set 
 $\{  q_{-1}, 2 \check{\rho}, q_1  \}$ forms an $\mfs \mfl_2$-triple.

\begin{exa}[Case of $\mbG = \mr{PGL}_n$] \label{QS8376}
Let us consider the case where $\mbG = \mr{PGL}_n$ (hence $\mfg = \mfp \mfg \mfl_n = \mfs \mfl_n$ because of the assumption $n <p$).
We set
\begin{equation} \label{BorelB}
T \left(\text{resp.,} \ B\right)
\index{$T$} \index{$B$}
\end{equation}
to be the Borel subgroup (resp., the maximal torus) of $\mr{PGL}_n$ 
defined as  the image, via the natural quotient $\mr{GL}_n \migisurj \mr{PGL}_n$,
 of  the group of 
 invertible  upper-triangular matrices (resp., invertible diagonal matrices).
If $\Gamma$ denotes, as before, the set of simple roots in $B$ with respect to $T$,
then we can find a unique collection 
$\ \qq := \{ x_\alpha \}_{\alpha \in \Gamma}$ of generators in  the root spaces $\mfs \mfl_n^\alpha$ so that the resulting 
   $\mfs \mfl_2$-triple  $\{  q_{-1}, 2  \check{\rho}, q_1  \}$  can be given by 
\begin{align} \label{QQ0230}
& q_{-1} = \begin{pmatrix} 0 & 0 & \cdots  & 0 & 0 \\  1 & 0 & \cdots  & 0 & 0 \\  0 & 1 & \cdots  & 0 & 0 \\ \vdots  & \vdots  & \ddots  & \vdots  & \vdots  \\  0 & 0 & \cdots  & 1 & 0    \end{pmatrix}, \  
2 \check{\rho} = \begin{pmatrix} n-1 & 0 & 0 & \cdots   & 0 \\  0 & n-3 & 0 &  \cdots   & 0 \\  0 & 0 & n-5 & \cdots   & 0 \\ \vdots  & \vdots  & \vdots  & \ddots  & \vdots  \\  0 & 0 & 0 &  \cdots   & - (n-1)
\end{pmatrix}, 
\end{align} 
\begin{align}
 & 
  q_1 = \begin{pmatrix} 0 & n-1 & 0 & 0 & \cdots   & 0 \\  0 & 0 & 2(n-2) & 0 &  \cdots   & 0 \\  0 & 0 & 0  & 3(n-2) & \cdots   & 0 \\ \vdots  & \vdots  & \vdots  & \vdots &  \ddots  & \vdots  \\  0 & 0 & 0 & 0&  \cdots   & n-1   \\ 0 & 0 & 0 & 0&  \cdots   &0    \end{pmatrix}. 
  \notag
\end{align} 
We suppose that  $T$ (resp., $B$; resp., $\ \qq$)   is  selected   as the  preferred  choice among all maximal torus (resp., all Borel subgroups; resp., all collections of generators in $\mfg^\alpha$'s)
whenever we are  dealing with the algebraic group $\mr{PGL}_n$.

\end{exa}

\subsection{Local $(\mfg, \hslash)$-opers} \label{QG445}

In what follows, let us introduce the local version of a $(\mfg, \hslash)$-oper 
 and
observe (cf. Proposition \ref{y039}) that, after  application of  gauge  transformation,  any such object may be ``normalized'' with respect to 
the $\mfg_{-1}$-components in the matrix form of that connection.

\bde \label{QQ51}
A {\bf local $(\mfg, \hslash)$-oper}  \index{local $(\mfg, \hslash)$-oper} on $U^\mr{log}/S^\mr{log}$ is an $\hslash$-$S^\mr{log}$-connection $\nabla$ on the trivial $\mbG$-bundle $U \times_k \mbG$ such that the pair $(U \times_k \mbB, \nabla)$ forms a $(\mfg, \hslash)$-oper on $U^\mr{log}/S^\mr{log}$.
  \ede

\begin{rema} \label{QG7000}
It is clear  that any $(\mfg, \hslash)$-oper is obtained  \'{e}tale locally from a local $(\mfg, \hslash)$-oper because the underlying $\mbB$-bundle is locally trivial.
\end{rema}

Let $\nabla$ be a local $(\mfg, \hslash)$-oper on $U^\mr{log}/S^\mr{log}$.
Then, by  the definition of a $(\mfg, \hslash)$-oper, we have 
$\mr{Im}(\nabla^\mfg) \subseteq \mcO_U \otimes_k \mfg^{-1}$ (cf. (\ref{QQ12}) for the definition of $\nabla^\mfg$) and, for each $\alpha \in \Gamma$,
the composite 
\begin{equation}  \label{QQ49}
 \mcT_{U^\mr{log}/S^\mr{log}} \xrightarrow{\nabla^\mfg} \mcO_U \otimes_k \mfg^{-1} \migisurj  \left(\mcO_U \otimes_k (\mfg^{-1}/\mfg^0)=\right) \mcO_U \otimes_k \mfg_{-1} \migisurj \mcO_U \otimes_k \mfg^{-\alpha}
 \end{equation}
is an isomorphism.

\bde \label{y038}
Let us keep the above notation.
Suppose further  that there exists 
 a global  log chart $(U, \partial)$  on $U^\mr{log}/S^\mr{log}$  (cf. Definition \ref{ppy021}).
Then, we shall say that  $\nabla$ is {\bf  $\ \qq$-prenormal  relative to  $(U, \partial)$}  \index{local $(\mfg, \hslash)$-oper@--- $\, \qq$-prenormal  relative to  $(U, \partial)$}
 if, for any $\alpha \in \Gamma$,  
 the image of $\partial$ via (\ref{QQ49}) coincides with  $1 \otimes y_\alpha$, or equivalently, 
 the composite 
 \begin{equation} 
 \mcT_{U^\mr{log}/S^\mr{log}} \xrightarrow{\nabla^\mfg} \mcO_U \otimes_k \mfg^{-1} \migisurj   \mcO_U \otimes_k \mfg_{-1} 
 \end{equation}
 sends $\partial$ to $1 \otimes q_{-1}$.
  \ede
\begin{exa} \label{RRRH}
Let us consider the  case where  $\mbG = \mr{PGL}_n$ (together with $T$, $B$,  and   $\, \qq$ defined   in Example \ref{QS8376}).
In particular, $\mfg = \mfp \mfg \mfl_n = \mfs \mfl_n$.
Suppose that  there exists a global  log chart $(U, \partial)$ of $U^\mr{log}/S^\mr{log}$.
 Then,  
an $\hslash$-$S^\mr{log}$-connection $\nabla$ on $U \times_k \mr{PGL}_n$ forms a local $(\mfg, \hslash)$-oper   (resp., a  local  $(\mfg, \hslash)$-oper that is $\ \qq$-prenormal  relative to $(U, \partial)$)  if and only if 
the element $\nabla^\mfg (\partial)$ of  $\mfp \mfg \mfl_n (U) \left(= \mfs \mfl_n (U) \right)$ can be represented by  an $n \times n$ matrix of the form
\begin{align} \label{QQ191}
 \begin{pmatrix} * & *& *& \cdots & *& * \\ u_1& * & *& \cdots & * & * \\ 0 & u_2  & *& \cdots & * & * \\ 0 & 0 &u_3& \cdots & *& * \\ \vdots & \vdots & \vdots & \ddots& \vdots & \vdots  \\ 0 & 0 & 0 & \cdots & u_{n-1} & *  \end{pmatrix}
 \ \left(\text{resp.,} \ 
 \begin{pmatrix} * & *& *& \cdots & *& * \\ 1& * & *& \cdots & * & * \\ 0 & 1 & *& \cdots & * & * \\ 0 & 0 &1& \cdots & *& * \\ \vdots & \vdots & \vdots & \ddots& \vdots & \vdots  \\ 0 & 0 & 0 & \cdots & 1 & *  \end{pmatrix} \right),
\end{align}
 where the $*$'s indicate arbitrary   functions on $U$ and   $u_1, \cdots, u_{n-1}$ in the non-resp'd portion are nowhere vanishing functions, i.e.,  elements in  $\mbG_m (U)$.
\end{exa}

\bpr  \label{y039}
Let $\nabla$ be a local $(\mfg, \hslash)$-oper on $U^\mr{log}/S^\mr{log}$.
Suppose that  there exists 
 a global  log chart $(U, \partial)$  on $U^\mr{log}/S^\mr{log}$.
 Then, there exists a unique $U$-rational point  $h : U \migi \mbT$  of $\mbT$ 
 such that $\mr{L}_{h*}(\nabla)$ is a $\,\qq$-prenormal  local $(\mfg, \hslash)$-oper   relative to $(U, \partial)$.
\epr
\begin{proof}
For each $\alpha \in \Gamma$,  let $u_\alpha$ be the element of $\Gamma (U, \mcO^\times_U) \left(=\mbG_m (U) \right)$ such that the image of $\partial$ via (\ref{QQ49}) coincides with $u_\alpha \otimes y_\alpha$. 
Let $h : U \migi \mbT$ be an arbitrary  $U$-rational point of $\mbT \left(\subseteq \mbG\right)$.
The automorphism $d\mr{L}_{h}^\mbG$ (cf.  (\ref{QQ10})) of $\widetilde{\mcT}_{U^\mr{log} \times_k \mbG/S^\mr{log}}$ obtained from the differential of 
the left-translation $\mr{L}_{h}$
 is compatible with    the filtration $\{ \widetilde{\mcT}^j_{U^\mr{log} \times_k\mbG/S^\mr{log}}\}_{j \in \mbZ}$.
  Hence, it 
 induces an automorphism
 of $\bigoplus_{\alpha \in \Gamma}    \mfg^{-\alpha}_{U \times_k \mbB}$ via 
   the decomposition
\begin{align} \label{QQ67}
\widetilde{\mcT}_{U^\mr{log} \times_k \mbG/S^\mr{log}}^{-1} / \widetilde{\mcT}^0_{U^\mr{log} \times_k \mbG/S^\mr{log}} = \bigoplus_{\alpha \in \Gamma}    \mfg^{-\alpha}_{U \times_k \mbB}
\end{align}
 (cf. (\ref{decom97})).
This automorphism
acts on each component $ \mfg^{-\alpha}_{U \times_k \mbB} \left(\cong \mcO_U \otimes_k \mfg^{-\alpha}\right)$ as the adjoint action by 
$h$; this action  coincides 
 with the automorphism $h^{-1}_\alpha \in \mr{GL}(\mfg^{-\alpha}_{U \times_k \mbB}) \left(=\mbG_m (U) \right)$ given by multiplication by $(-\alpha) \circ h \left(=\alpha\circ (h^{-1})\right) : U \migi \mbG_m$.
The assignment $h \mapsto (h^{-1}_\alpha)_{\alpha \in \Gamma}$ determines an isomorphism 
\begin{equation} \label{torus}
\mbT (U)\isom \prod_{\alpha \in \Gamma} \mr{GL}(\mfg^{-\alpha}_{U \times_k \mbB}).
\end{equation}
It follows that  we can find a  $U$-rational point $h$ of $T$ mapped to $(u_\alpha^{-1})_{\alpha \in \Gamma}$ by  this isomorphism.
Then,  the composite  (\ref{QQ49}) with $\nabla$ replaced by $\mr{L}_{h*}(\nabla) \left(= d\mr{L}_h^\mbG \circ \nabla \right)$ sends $\partial$ to $u_\alpha^{-1} \cdot (u_\alpha \otimes y_\alpha) = 1 \otimes y_\alpha$.
That is to say, this $U$-rational point  $h$  fulfills the requirement.
This  completes the existence assertion.
Since the uniqueness portion follows immediately  from the injectivity of (\ref{torus}),
we finish the proof of the assertion.
\end{proof}



\subsection{Automorphisms of a $(\mfg, \hslash)$-oper} 
\label{QR26}

The following proposition was already proved in ~\cite{BD2}, \S\,1.3, Proposition (or ~\cite{Bar}, Proposition 2.1),  if $U^\mr{log}/S^\mr{log}$  is an unpointed smooth curve over the field of complex numbers $\mbC$.
Even if $k$ has positive  characteristic and $U$ is an arbitrary log curve, 
 we can also prove it by  the same argument  as in {\it loc.\,cit.} thanks to the condition imposed at the beginning of this section.


\bpr \label{y040}
(Recall  that  either  ``\,$\mr{char}(k)=0$'',  ``\,$\mr{char}(k)> 2  h_\mbG$'',  or  ``\,$\mbG = \mr{PGL}_n$ with $1 < n < p$'' is fulfilled.)
 Any $(\mfg, \hslash)$-oper on $U^\mr{log}/S^\mr{log}$ does not have nontrivial automorphisms.
 \epr
\begin{proof}
We only consider the case where either  ``\,$\mr{char}(k)=0$'' or  ``\,$\mr{char}(k)> 2 h_\mbG$''  is fulfilled because the remaining  case follows from an entirely similar discussion (using  the result in Remark \ref{ppp037}).
Let $\msE^\spadesuit$ be a $(\mfg, \hslash)$-oper on $U^\mr{log}/S^\mr{log}$.
 We may assume, after possibly replacing $U$ with its covering in the \'{e}tale topology,   that  there exist
a global  log chart $(U, \partial)$ and  a local $(\mfg, \hslash)$-oper $\nabla$ on $U^\mr{log}/S^\mr{log}$ with $\msE^\spadesuit =(U \times_k \mbB, \nabla)$.
Moreover, by Proposition  \ref{y039} above, it suffices to consider the case where  
$\nabla$
 is $\ \qq$-prenormal relative to  $(U, \partial)$.
Let us take an automorphism $\eta$ of $\msE^\spadesuit$, i.e., an automorphism $\eta$
   of the  trivial $\mbG$-bundle  $U \times_k \mbG$ which preserves the  $\mbB$-reduction $U \times_k \mbB$ and satisfies 
   $\eta_*(\nabla) = \nabla$.
  Then, we can find a $U$-rational point $h: U \migi \mbB$ of $\mbB$ with $\eta = \mr{L}_h$ (cf. Example \ref{QQ9999}).
Let $\overline{h} : U \migi \mbT$ denote the composite 
\begin{equation}
\overline{h} : U \xrightarrow{h} \mbB \migisurj \left( \mbB/\mbN =\right) \mbT.
\end{equation}
Notice that  
$d \mr{L}_{h}^\mbG$ and $d \mr{L}_{\overline{h}}^\mbG$
 act on  $\bigoplus_{\alpha \in \Gamma}    \mfg^{-\alpha}_{U \times_k \mbB}$  in the same manner.
Hence,  $\mr{L}_{\overline{h}*}(\nabla)$ turns out to be $\ \qq$-prenormal  relative to $(U, \partial)$ because of the assumption   $d \mr{L}_{h}^\mbG \circ \nabla   \left(=\mr{L}_{h*}(\nabla) =  \eta_*(\nabla) \right)=\nabla$.
It follows from the uniqueness portion of Proposition \ref{y039} that
$\overline{h}=e$,
 i.e.,
 $h$
   factors through the inclusion $\mbN \migiincl \mbB$; we abuse notation by using $h$ for 
  the resulting morphism $U \migi \mbN$. 
It induces a $U$-rational point $\mr{Log} (h) : U \migi \mfn \left(=\mfg^1 \right)$ (cf. (\ref{QR630})). 
By Corollary \ref{y028}, the equality   $\mr{L}_{h*}(\nabla)  = \nabla$ implies  
\begin{equation} \label{above}
\nabla^\mfg (\partial) = -\hslash \cdot  d \mr{Log} (h) (\partial) + \sum_{s =0}^\infty \frac{1}{s !} \cdot \mr{ad}(\mr{Log} (h))^s(\nabla^\mfg (\partial)).
  \end{equation}
Now, we shall assume  that $\mr{Log}(h) \neq 0$.
Then, the positive integer  
\begin{equation}
j_0 := \mr{max}\left\{ j \in \mbZ_{>0}  \ \big| \  \mr{Log} (h) \in \mfg^j (U) \right\}
\end{equation}
 is well-defined.
For each $j \in \mbZ$,  let 
\begin{equation} \mr{Log} (h)_j : U \migi \mfg_j  \ \ \text{and} \ \ \nabla^\mfg (\partial)_j \in \Gamma (U, \mcO_U \otimes_k \mfg_j)
\end{equation}
 be the $j$-th components of, respectively,  $\mr{Log} (h) : U \migi \left(\mfg^1 \subseteq \right) \mfg$ and  $\nabla^\mfg (\partial) \in \Gamma (U, \mcO_U \otimes_k \mfg)$ with respect to the decomposition  $\mfg = \bigoplus_{j \in \mbZ} \mfg_j$.
By the definition of $j_0$, we see that 
 $d \mr{Log} (h) (\partial) \in \Gamma (U, \mcO_U \otimes_k \mfg^{j_0})$.
Since $\nabla^\mfg (\partial)_{-1} = 1 \otimes q_{-1}$ and $\nabla^\mfg (\partial)_j= 0$ for  $j < -1$, 
the section  $\mr{ad}(\mr{Log} (h)_j)^s(\nabla^\mfg (\partial)_{j'})$ lies in $\Gamma (U, \mcO_U \otimes_k \mfg^{j_0})$ if either ``$s \geq 2$" or ``$s =1$ and $j + j' \geq j_0$" is satisfied. 
Hence, 
 (\ref{above}) implies   that,  in the component  $\Gamma (U, \mcO_U \otimes_k \mfg_{j_0 -1})$ of $\Gamma (U, \mcO_U \otimes_k \mfg)$, we have
\begin{equation}
\mr{ad}(\mr{Log} (h)_{j_0})(\nabla^\mfg (\partial)_{-1}) \ \left(= - \mr{ad}(1 \otimes q_{-1}) (\mr{Log} (h)_{j_0})\right) = 0.
\end{equation}
But, since $j_0>0$,  the morphism
$\mr{ad}(1 \otimes q_{-1}) : \mcO_U \otimes_k \mfg_{j_0} \migi \mcO_U \otimes_k \mfg_{j_0-1}$ is injective.
It follows that   $\mr{Log} (h)_{j_0} =0$, which 
  contradicts  the definition of $j_0$.
Thus, $\mr{Log}(h)$ must be zero.
Because of the bijectivity of $\mr{Log} : \mbN \isom \mfn$,
we have  $h = e$, i.e., 
  $\eta$ coincides with  the identity morphism of $\msE^\spadesuit$.
This completes the proof of the proposition.
\end{proof}

\subsection{The sheaf of $(\mfg, \hslash)$-opers} \label{QR30}
By taking account of the \'{e}tale descent property of $(\mfg, \hslash)$-opers,
one may define the stack in groupoids
\begin{equation} \label{opersite}
\mcO p_{\mfg, \hslash, U^\mr{log}/S^\mr{log}} \migi  \mfE  \mft_{/U}
\index{$\mcO p_{\mfg, \hslash, U^\mr{log}/S^\mr{log}}$, $\mcO p_{\mfg, \hslash}$}
\end{equation}
over  the small \'{e}tale site  $\mfE  \mft_{/U}$ on $U$ whose category of sections over an
\'{e}tale $U$-scheme $U'$
is the groupoid of 
 $(\mfg, \hslash)$-opers on $U'^{\mr{log}} /S^\mr{log}$, where the log structure on $U'$ is obtained by pulling-back that on  $U^\mr{log}$.
It follows from Proposition  \ref{y040} that $\mcO p_{\mfg, \hslash, U^\mr{log}/S^\mr{log}} $ is fibered in equivalence relations (cf. ~\cite{FGA}, Definition 3.39), i.e., specifies a set-valued sheaf on $\mfE  \mft_{/U}$. 
If there is no fear of confusion, we shall write
\begin{align} \label{QQ0300}
\mcO p_{\mfg, \hslash} := \mcO p_{\mfg, \hslash, U^\mr{log}/S^\mr{log}}.
\end{align}

\section{The case of $\mfg = \mfs \mfl_2$} \label{y042}

We shall collect some facts about $(\mfs \mfl_2, \hslash)$-opers  used in the subsequent discussion.
These facts will be proved in Chapter 4 later, so  it should be noted that we are in advance referring  to them without giving their   proofs.
In this section, we suppose that  {\it  $2$ is invertible in $k$}.

\subsection{}\label{QR40}
Let $\mbG^\odot$ be the projective linear group over $k$ of rank $2$, i.e., $\mbG^\odot := \mr{PGL}_2$, and  
$\mbB^\odot$ (resp., $\mbT^\odot$)  the Borel subgroup (resp., the maximal torus) of $\mbG^\odot$ defined as the image, via the  natural quotient $\mr{GL}_2 \migisurj \mbG^\odot$,  of invertible  upper-triangular matrices (resp., invertible diagonal matrices).
Let us write $\mbN^\odot := [\mbB^\odot, \mbB^\odot]$.
Also, write $\mfg^\odot$, $\mfb^\odot$,  $\mfn^\odot$, and $\mft^\odot$ for the Lie algebras of
$\mbG^\odot$, $\mbB^\odot$,   $\mbN^\odot$, and $\mbT^\odot$ respectively.
 In particular,  we have 
 \index{$\mfg$, Lie algebra@--- $\mfg^\odot$}
 $ \mfg^\odot = \mfs \mfl_2$, $\mfg^\odot_{1} = \mfn^\odot$, 
 $\mfg^\odot_0 = \mft^\odot  =  \mfb^\odot/\mfn^\odot$,  and $\mfg^\odot_{-1} = \mfg^\odot / \mfb^\odot$.
The set of simple roots $\Gamma$ in $\mbB^\odot$ with respect to $\mbT^\odot$ consists exactly  of one element $\alpha^\odot$. 
We can choose the  $\mfs \mfl_2$-triple $\{q_{-1}^\odot, 2\check{\rho}^\odot,   q_{1}^\odot\}$ given by 
\begin{align} \label{QS20}
q_1^\odot \left(=x_{\alpha^\odot}\right)  := \overline{\begin{pmatrix}  0& 1 \\  0 & 0\end{pmatrix}},
\hspace{2mm}
 \check{\rho}^\odot \left(=\check{\omega}_{\alpha^\odot}\right): = \overline{\begin{pmatrix}  1/2& 0 \\  0 & -1/2\end{pmatrix}},
\hspace{2mm}
q_{-1}^\odot \left(=y_{\alpha^\odot}\right): = \overline{\begin{pmatrix}  0& 0 \\  1 & 0\end{pmatrix}}.
\end{align}
In particular, we obtain $\ \qq^\odot := \{ x_{\alpha^\odot}\}$ as in   (\ref{QG5697}) for $\mfg^\odot$.


\subsection{Results in Chap.\,4: Part I}\label{QQ0402}

 Let  us write 
\index{$\DE_{\mbB, \hslash, U^\mr{log}/S^\mr{log}}$, $\DE_{\mbB}$, $\DE_{\mbG, \hslash, U^\mr{log}/S^\mr{log}}$, $\DE_{\mbG}$}
 \begin{align}  \label{QQ0406}
 \DE_{\mbB^\odot, \hslash, U^\mr{log}/S^\mr{log}}, \ \text{or simply} \ 
 \DE_{\mbB^\odot},
 \end{align}
  for  the 
 $\mbB^\odot$-bundle  on $U$ constructed  in (\ref{QQ0430}) for  $n=2$.
 The $\mbG^\odot$-bundle induced by $\DE_{\mbB^\odot}$ via change of structure group by the inclusion $\mbB^\odot \migiincl \mbG^\odot$ will be denoted by
 \index{$\DE_{\mbB, \hslash, U^\mr{log}/S^\mr{log}}$, $\DE_{\mbB}$, $\DE_{\mbG, \hslash, U^\mr{log}/S^\mr{log}}$, $\DE_{\mbG}$}
 \begin{align} \label{QQ0408}
 \DE_{\mbG^\odot, \hslash, U^\mr{log}/S^\mr{log}}, \ \text{or simply} \ 
 \DE_{\mbG^\odot}.
 \index{$\DE_{\mbG^\odot, \hslash, U^\mr{log}/S^\mr{log}}$, $\DE_{\mbG^\odot}$}
 \end{align}
For each log chart $(V, \partial)$ on $U^\mr{log}/S^\mr{log}$, there exists a specific  trivialization 
\index{$\mr{triv}_{\mbG, (V,\partial)}$}
\begin{align} \label{QR48}
\mr{triv}_{\mbB^\odot, (V, \partial)} : V \times_k \mbB^\odot  \isom \DE_{\mbB^\odot} |_V
\end{align}
 (cf. (\ref{QQ0431})) of  $\DE_{\mbB^\odot}$ restricted to $V$.
 The formations of both  $ \DE_{\mbB^\odot}$ and $\mr{triv}_{\mbB^\odot, (V, \partial)}$ are  functorial with respect to pull-back to \'{e}tale $U$-schemes, as well as base-change to fs $S^\mr{log}$-log schemes.
 
 Also, let us suppose that $U^\mr{log}/S^\mr{log}$ comes from  a local $r$-pointed stable curve $\msU :=(U/S, \{ \sigma^{U}_i :S_i \migi U \})$ over $S$ with $r >0$ together with  a choice of its base $\msX$.
 Then, for each $i=1, \cdots, r$, there exists a canonical trivialization
 \index{$\mr{triv}_{\mbG, \sigma_i^U}$}
 \begin{align} \label{QW4599}
 \mr{triv}_{\mbB^\odot, \sigma^U_i} : S_i \times_k \mbB^\odot \isom \sigma^{U*}_i (\DE_{\mbB^\odot})
 \end{align}
 of the $\mbB^\odot$-bundle  $\sigma^{U*}_i (\DE_{\mbB^\odot})$ on $S_i$ (cf. Remark \ref{QW3982}).

\subsection{Results in Chap.\,4: Part II} \label{QQ0450}
 In Definition \ref{y0132}, we introduce the notion of a {\it normal}  $(\mfg^\odot, \hslash)$-oper, which is a $(\mfg, \hslash)$-oper  of the form $(\DE_{\mbB^\odot}, \nabla)$ for a suitable $\hslash$-$S^\mr{log}$-connection $\nabla$ on
 $\DE_{\mbG^\odot}$.
 The formation of a normal $(\mfg^\odot, \hslash)$-oper  
   is   closed 
 under pull-back to  \'{e}tale $U$-schemes, as well as base-change to fs $S^\mr{log}$-log schemes.
Also, it follows from  Proposition \ref{QQ0470} that
 for any $(\mfg^\odot, \hslash)$-oper  $\msE^\spadesuit$ on $U^\mr{log}/S^\mr{log}$,
there exists  uniquely a pair
\begin{align} \label{QQ0475}
({^\dagger}\msE^\spadesuit, \mr{nor}_{\msE^\spadesuit})
\end{align}
consisting of a normal $(\mfg^\odot, \hslash)$-oper on $U^\mr{log}/S^\mr{log}$ and an isomorphism $\mr{nor}_{\msE^\spadesuit} : {^\dagger}\msE^\spadesuit \isom \msE^\spadesuit$ of $(\mfg^\odot, \hslash)$-opers.

\subsection{Results in Chap.\,4: Part III} \label{QQ0490}
For $j = -1, 0, 1$,  
there exists  a canonical isomorphism
$(\mfg^\odot_j)_{\DE_{\mbB^\odot}}  \isom \Omega_{U^\mr{log}/S^\mr{log}}^{\otimes j}$
 of $\mcO_U$-modules (cf.  (\ref{QW56})).
 In particular,   the case of $j=1$ gives an isomorphism
 \begin{align}
 \label{QR0104}
 \mcH om_{\mcO_U}(\mcT_{U^\mr{log}/S^\mr{log}}, (\mfn^\odot)_{\DE_{\mbB^\odot}})  \isom \Omega_{U^\mr{log}/S^\mr{log}}^{\otimes 2}.
 \end{align}

Now, let  $\msE^\spadesuit := (\DE_{\mbB^\odot}, \nabla)$  be a normal $(\mfg^\odot, \hslash)$-oper on $U^\mr{log}/S^\mr{log}$  and 
let $R \in \Gamma (U, \Omega_{U^\mr{log}/S^\mr{log}}^{\otimes 2})$.
By identifying $R$ with an $\mcO_U$-linear morphism $\mcT_{U^\mr{log}/S^\mr{log}} \migi  (\mfn^\odot)_{\DE_{\mbB^\odot}}$ $\left(\subseteq \widetilde{\mcT}_{\DE^{\mr{log}}_{\mbG^\odot}/S^\mr{log}} \right)$ via (\ref{QR0104}), 
we obtain  the sum 
 \begin{equation} \label{QQ0499}
 \nabla + R : \mcT_{U^\mr{log}/S^\mr{log}} \migi \widetilde{\mcT}_{\DE^{\mr{log}}_{\mbG^\odot}/S^\mr{log}},
 \end{equation}
which  defines   an $\hslash$-$S^\mr{log}$-connection  on $\DE_{\mbG^\odot}$.
 According to the discussion in \S\,\ref{QW6}, 
the resulting pair
\index{$\msE^\spadesuit$, $(\mfg, \hslash)$-oper@--- $\msE^{\spadesuit}_{+ R}$}
\begin{equation} \label{plus1}
\msE^{\spadesuit}_{+ R} := (\DE_{\mbB^\odot},  \nabla + R)
\end{equation}
forms a normal $(\mfg^\odot, \hslash)$-oper on $U^\mr{log}/S^\mr{log}$, and 
the assignment $(\msE^{\spadesuit}, R) \mapsto \msE^{\spadesuit}_{+ R}$ is functorial 
with respect to pull-back to  \'{e}tale $U$-schemes, as well as base-change to   fs $S^\mr{log}$-log schemes.
Moreover, 
 this assignment determines  a structure of $\Omega^{\otimes 2}_{U^\mr{log}/S^\mr{log}}$-torsor
\begin{align} \label{QQ0500}
\mcO p_{\mfg^\odot, \hslash} \times \Omega_{U^\mr{log}/S^\mr{log}}^{\otimes 2} \migi  \mcO p_{\mfg^\odot, \hslash}
\end{align}
 on $\mcO p_{\mfg^\odot, \hslash}$, where $\Omega_{U^\mr{log}/S^\mr{log}}^{\otimes 2}$ are considered as a sheaf of abelian groups on $\mfE \mft_{/U}$ in an evident manner.
In particular,  
the stalk of $\mcO p_{\mfg^\odot, \hslash}$ at any geometric point of $U$ is nonempty.

\section{The  $\mbG$-bundle $\DE_\mbG$ and the vector bundle $\DV_\mfg$} \label{y044}

Next, let us go back to the case of 
a general  $\mfg$.
In the rest of this chapter,
 we assume, unless otherwise stated,   
 that  {\it either ``\,$\mr{char}(k)=0$'' or ``\,$\mr{char}(k)  >2  h_\mbG$'' is fulfilled}.

\subsection{Change of structure group arising from an $\mfs \mfl_2$-triple}\label{QQ0533}
 
The $\mfs \mfl_2$-triple $\{  q_{-1}, 2  \check{\rho}, q_1  \}$
 associated to  the fixed collection $\ \qq$
   determines an injection
    \index{$\iota_\mfg$, $\iota_\mbB$, $\iota_\mfc$}
\begin{equation} \label{injg} \iota_\mfg : \mfg^\odot \migiincl \mfg,  
\end{equation}
i.e., $\iota_\mfg$ satisfies $\iota_\mfg (q_1^\odot) = q_1$, $\iota_\mfg (2  \check{\rho}^\odot) = 2  \check{\rho}$, and $\iota_\mfg (q_{-1}^\odot) = q_{-1}$.
In particular, we have
 $\iota_\mfg (\mfb^\odot) \subseteq \mfb$.

\bpr \label{py040}
 There exists 
 a unique  injective morphism
   \index{$\iota_\mfg$, $\iota_\mbB$, $\iota_\mfc$}
  \begin{align} \label{rrrr}
   \iota_\mbB : \mbB^\odot \migiincl \mbB 
  \end{align}
of algebraic $k$-groups which satisfies $\iota_\mbB (\mbT^\odot) \subseteq \mbT$, $\iota_\mbB (\mbN^\odot) \subseteq \mbN$ and 
 whose differential at the identity element $e \in \mbB^\odot$ coincides with the restriction $\iota_\mfg |_{ \mfb^\odot} : \mfb^\odot \migiincl \mfb$ of $\iota_\mfg$.
 \epr
\begin{proof}
Note that  the composite
\begin{align} \label{rrrr2}
\iota_\mbN : \mbN^\odot \xrightarrow{\mr{Log}} \mfn^\odot \xrightarrow{\iota_\mfg |_{ \mfn^\odot}} \mfn \xrightarrow{\mr{Exp}} \mbN
\end{align}
 becomes  a morphism of algebraic groups over $k$.
Also,  the following composite forms  a morphism of algebraic groups over $k$:
\begin{align} \label{rrrr1}
\iota_\mbT : \mbT^\odot \isom \mr{GL} (\mfn^\odot) \migiincl \prod_{\alpha \in \Gamma} \mr{GL} (\mfg^{\alpha}) \isom \mbT, 
\end{align}
where
\begin{itemize}
\item
the first arrow denotes the isomorphism which, to any $h \in \mbT^\odot$, assigns the automorphism of $\mfn^\odot$ given by multiplication by $\alpha (h) \in \mbG_m$ (where $\alpha$ denotes the unique simple root in $\mbB^\odot$);
\item
the second arrow denotes the diagonal embedding under the natural identifications $\mr{GL} (\mfn^\odot) = \mr{GL} (\mfg^{\alpha}) \left(= \mbG_m \right)$;
\item
the  third arrow denotes the inverse of the isomorphism  $\mbT \isom \prod_{\alpha \in \Gamma} \mr{GL} (\mfg^{\alpha})$
which, to any $h \in \mbT$, assigns  the collection  $(h_\alpha)_{\alpha \in \Gamma} \in \prod_{\alpha \in \Gamma} \mr{GL} (\mfg^{\alpha})$ consisting of  the automorphisms $h_\alpha$ ($\alpha \in \Gamma$) given by  multiplication by $\alpha \circ h \in \mbG_m$.
\end{itemize}
In the following, we shall  construct a morphism $\iota_\mbB : \mbB^\odot \migi \mbB$ extending 
both the morphisms $\iota_\mbN$ and $\iota_\mbT$.
Let us take $t \in \mbT^\odot$ and $n \in \mbN^\odot$,  and write $n^t := t^{-1}\cdot   n \cdot  t \in \mbN^\odot$.
Also, let  us choose $2 \times 2$ matrices   of the forms  $\begin{pmatrix} a & 0 \\ 0 & 1 \end{pmatrix}$ and  $\begin{pmatrix} 1 & b \\ 0 & 1 \end{pmatrix}$  representing $t$ and $n$ respectively.
Then, $n^t$ may be represented by the matrix $\begin{pmatrix} 1 & a^{-1} \cdot   b \\ 0 & 1 \end{pmatrix}$, and 
the images $\mr{Log} (n)$ and $\mr{Log} (n^t)$ in $\mfn^\odot$ may be represented by 
 $\begin{pmatrix} 0 & b \\ 0 & 0 \end{pmatrix}$ and $\begin{pmatrix} 0 & a^{-1} \cdot   b \\ 0 & 0 \end{pmatrix}$ respectively.
Now, let $\alpha \in \mr{Hom} (\mbT, \mbG_m)$ be a root   with $\mfg^\alpha \in \mfg_j$ (for some $j \in \mbZ$)
 and let $v \in \mfg^\alpha$.
By identifying elements of $\mbN$ with automorphisms of  $\mfg$ via the inclusion $\mbN \migiincl (\mbG \cong) \ \mr{Aut}^0 (\mfg)$,
 we have a sequence of equalities
\begin{align}
(\iota_\mbT (t) \cdot \iota_\mbN ( n^t)) (v)  
&=  \iota_\mbT (t) \left(\sum_{s =0}^\infty \frac{1}{s !} \cdot  \mr{ad}(a^{-1} \cdot   b \cdot    q_1)^s (v)\right)  \\
& = \sum_{s =0}^\infty \frac{1}{s !} \cdot a^{s + j} \cdot  \mr{ad}(a^{-1} \cdot   b  \cdot  q_1)^s (v)  \notag \\
& = \sum_{s =0}^\infty \frac{1}{s !} \cdot  \mr{ad}(b \cdot    q_1)^s (a^j \cdot   v) \notag \\
& = \iota_\mbN (n)  (a^j  \cdot  v) \notag \\
&  =  \iota_\mbN (n) \cdot \iota_\mbT (t) (v), \notag
\end{align}
where the second equality follows from the fact that  $\mr{ad}(a^{-1}\cdot  b    \cdot  q_1)^s (v) \in \mfg_{s + j}$.
This implies that
 $\iota_\mbT (t) \cdot \iota_\mbN ( n^t) = \iota_\mbN (n) \cdot \iota_\mbT (t)$
 for any $t \in \mbT^\odot$, $n \in \mbN^\odot$.
Thus, 
the assignment $(n, t) \mapsto \iota_\mbN (n) \cdot \iota_\mbT (t)$ determines  a well-defined morphism
$\iota_\mbB : (\mbN^\odot \rtimes \mbT^\odot =) \ \mbB^\odot \migi \mbB$ of algebraic groups, extending  both $\iota_\mbT$ and $\iota_\mbN$.
It is clear that  $\iota_\mbB$ is a unique morphism  satisfying  the required conditions.
Consequently, 
 we complete  the proof of Proposition \ref{py040}.
\end{proof}

\subsection{A morphism from $\mcO p_{\mfg^\odot, \hslash}$ to $\mcO p_{\mfg, \hslash}$}\label{QQ0644}
Let  $\msE_\odot^\spadesuit := (\mcE_{\mbB^\odot}, \nabla)$ be a $(\mfg^\odot, \hslash)$-oper   on $U^\mr{log}/S^\mr{log}$.
 We shall write $\mcE_{\mbG^\odot} := \mcE_{\mbB^\odot} \times^{\mbB^\odot} \mbG^\odot$,
$\mcE_{\mbB} := \mcE_{\mbB^\odot} \times^{\mbB^\odot, \iota_\mbB} \mbB$, and $\mcE_\mbG := \mcE_\mbB \times^{\mbB} \mbG$, where $\iota_\mbB$ denotes the morphism $\mbB^\odot \migiincl \mbB$ resulting from Proposition \ref{py040}.
Observe that the following diagram is commutative:
\begin{align} \label{Efrtyui}
\vcenter{\xymatrix@C=26pt@R=16pt{ 
& \mfb^\odot_{\mcE_{\mbB^\odot}} \ar[rr]^-{\text{inclusion}} 
\ar[dl]_{\text{(b)}} \ar[dd]|(.50)\hole^(.70){\iota_{\mfg^\odot/\mfb^\odot}}
& &  \widetilde{\mcT}_{\mcE^\mr{log}_{\mbB^\odot}/S^\mr{log}}
 \ar[dl]^{\text{(a)}} \\
 \mfb_{\mcE_{\mbB}} (= \mfb_{\mcE_{\mbB^\odot}})
  \ar[rr]^{\hspace{17mm}\text{inclusion}} \ar[dd]_-{\iota_{\mfg/\mfb}} 
  & &
  \widetilde{\mcT}_{\mcE^\mr{log}_{\mbB}/S^\mr{log}}
  & \\
& \mfg^\odot_{\mcE_{\mbB^\odot}}
 \ar[dl]^-{\text{(c)}} 
 & & 
  \\
\mfg_{\mcE_{\mbG}},  
& &  & 
}}
\end{align}
where
\begin{itemize}
\item
(a) is obtained by differentiating the inclusion $\mcE^\mr{log}_{\mbB^\odot} \migiincl \mcE^\mr{log}_{\mbB}$ induced by  
$\iota_\mbB$;
\item
(b) denotes the restriction of $(a)$, which  coincides with the restriction  $\iota_\mfg |_{\mfb^\odot} : \mfb^\odot  \migiincl \mfb$ of $\iota_\mfg$ twisted by $\mcE_{\mbB^\odot}$;
\item
(c) denotes the injection $\iota_\mfg : \mfg^\odot \migiincl \mfg$ twisted by $\mcE_{\mbB^\odot}$.
\end{itemize}
On the other hand, the following two commutative square diagrams are  cocartesian:
\begin{align} \label{QQ0800}
\vcenter{\xymatrix@C=46pt@R=36pt{
 \mfb^\odot_{\mcE_{\mbB^\odot}} \ar[r]^-{\text{inclusion}} \ar[d]_-{\iota_{\mfg^\odot/b^\odot}}& \widetilde{\mcT}_{\mcE^\mr{log}_{\mbB^\odot}/S^\mr{log}} \ar[d]^-{\widetilde{\iota}_{\mfg^\odot/b^\odot}} \\
 \mfg^\odot_{\mcE_{\mbG^\odot}}  \ar[r]_-{\text{inclusion}}& \widetilde{\mcT}_{\mcE^\mr{log}_{\mbG^\odot}/S^\mr{log}},
}} \hspace{10mm}
\vcenter{\xymatrix@C=46pt@R=36pt{
 \mfb_{\mcE_{\mbB}} \ar[r]^-{\text{inclusion}} \ar[d]_-{\iota_{\mfg/b}}& \widetilde{\mcT}_{\mcE^\mr{log}_{\mbB}/S^\mr{log}} \ar[d]^-{\widetilde{\iota}_{\mfg/b}} \\
 \mfg_{\mcE_{\mbG}}  \ar[r]_-{\text{inclusion}}& \widetilde{\mcT}_{\mcE^\mr{log}_{\mbG}/S^\mr{log}}.
}}
\end{align}
Hence, the diagram displayed in (\ref{Efrtyui}) induces an $\mcO_U$-linear injection 
\begin{align} \label{QQ0911}
d \iota_{\mcE_{\mbB^\odot}} : \widetilde{\mcT}_{\mcE^\mr{log}_{\mbG^\odot}/S^\mr{log}} \migi
\widetilde{\mcT}_{\mcE^\mr{log}_{\mbG}/S^\mr{log}}.
\end{align}
By taking account of the local description of Lie bracket operations (cf.   (\ref{QR922})), we see that $d \iota_{\mcE_{\mbB^\odot}}$ becomes a morphism of Lie algebroids. 
It follows that the composite
\begin{align} \label{QQ0923}
\iota_{\mbG *}(\nabla) : \mcT_{U^\mr{log}/S^\mr{log}}  \xrightarrow{\nabla} \widetilde{\mcT}_{\mcE^\mr{log}_{\mbG^\odot}/S^\mr{log}} \xrightarrow{d \iota_{\mcE_{\mbB^\odot}}}
\widetilde{\mcT}_{\mcE^\mr{log}_{\mbG}/S^\mr{log}}
\end{align}
specifies an $\hslash$-$S^\mr{log}$-connection on $\mcE_\mbG$ and 
 the pair 
\begin{equation} \label{associated} 
\iota_{\mbG *}(\msE^\spadesuit_\odot) :=
 (\mcE_{\mbB},  \iota_{\mbG*}(\nabla))
 \end{equation}
  forms a $(\mfg, \hslash)$-oper on $U^\mr{log}/S^\mr{log}$.
The assignment $\msE_\odot^\spadesuit \mapsto \iota_{\mbG *}(\msE_\odot^\spadesuit)$ is functorial with respect to pull-back to \'{e}tale $U$-schemes,  so it determines a morphism of sheaves 
 \index{$\iota_{\mfg}^{\mcO p}$, $\iota_{\mfg}^{\mfO \mfp}$, ${^p}\iota_{\mfg}^{\mfO \mfp}$}
\begin{equation}  
\label{opincl}
 \iota_{\mfg}^{\mcO p}  : \mcO  p_{\mfg^\odot, \hslash} \migi \mcO p_{\mfg, \hslash}.
   \end{equation}
In particular,  
the stalk of $\mcO p_{\mfg, \hslash}$ at any geometric  point of $U$ is nonempty because this is true for  $\mfg = \mfg^\odot$, as mentioned before.

\subsection{The  $\mbG$-bundle $\DE_\mbG$}\label{QQ0534}
By
    change of structure group via the morphism  $\iota_\mbB$ (resp., the composite of $\iota_\mbB$ and the inclusion $\mbB \migiincl \mbG$),
    the $\mbB^\odot$-bundle $\DE_{\mbB^\odot}$ induces a $\mbB$-bundle (resp., a $\mbG$-bundle)  on $U$, which will be denoted by
    \index{$\DE_{\mbB, \hslash, U^\mr{log}/S^\mr{log}}$, $\DE_{\mbB}$, $\DE_{\mbG, \hslash, U^\mr{log}/S^\mr{log}}$, $\DE_{\mbG}$}
\begin{equation} \label{canonicalEg2}
\DE_{\mbB, \hslash, U^\mr{log}/S^\mr{log}} \ \text{or simply} \ \DE_{\mbB}  \ \left(\text{resp.,} \  \DE_{\mbG, \hslash, U^\mr{log}/S^\mr{log}} \ \text{or simply} \ \DE_{\mbG}\right).
\end{equation}
If $(V, \partial)$ is a log chart on $U^\mr{log}/S^\mr{log}$, then the trivialization $\mr{triv}_{\mbB^\odot,  (V,\partial)}$  induces a trivialization 
\index{$\mr{triv}_{\mbG, (V,\partial)}$}
\begin{equation} \label{canonicalTg}
\mr{triv}_{\mbG, (V,\partial)} : V \times_k \mbG \isom \DE_\mbG |_V
\end{equation}
of the $\mbG$-bundle  $\DE_\mbG |_V$ on $V$.
This isomorphism is, by construction,  compatible with the respective $\mbB$-reductions $V \times_k \mbB$ and  $\DE_{\mbB}|_V$. 

Also, let us suppose that  $U^\mr{log}/S^\mr{log}$ arises from a local $r (>0)$-pointed stable curve $\msU := (U/S, \{ \sigma_i \}_{i=1}^r)$ over $S$ together with  a choice of its base $\msX$.
Then, for each $i=1, \cdots, r$, we have a canonical trivialization
\index{$\mr{triv}_{\mbG, \sigma_i^U}$}
\begin{align} \label{QW4421}
\mr{triv}_{\mbG, \sigma_i^U} : S_i \times_k \mbG \isom \sigma^{U*}_i (\DE_\mbG)
\end{align}
induced by $\mr{triv}_{\mbB^\odot, \sigma_i^U}$.

\bpr \label{y0130hh2} 
Let $(V, \partial)$ be a log chart on $U^\mr{log}/S^\mr{log}$ and let $a \in \Gamma (V, \mcO_V^\times)$ (hence, the pair $(V, a \cdot  \partial)$ forms a log chart on $U^\mr{log}/S^\mr{log}$).
Then, the composite automorphism 
\begin{align}
\mr{triv}_{\mbG, (V, \partial)}^{-1} \circ \mr{triv}_{\mbG, (V, a \cdot \partial)}
\end{align}
of  the trivial $\mbG$-bundle  $V \times_k \mbG$ coincides with  
the  left-translation $\mr{L}_h$ by
$h = \iota_\mbB (\overline{\begin{pmatrix} 1 & 0 \\ 0 & a \end{pmatrix}}) \in \mbG (V)$.
In particular, the $\mcO_V$-linear automorphism of $\mcO_V \otimes_k \mfg$ arising from  differentiating $\mr{triv}_{\mbG, (V, \partial)}^{-1} \circ \mr{triv}_{\mbG, (V, a \cdot \partial)}$  coincides with the automorphism given by assigning  $1 \otimes v \mapsto a^{-j}  \otimes v$ for any   $v \in \mfg_j$ with $j \in \mbZ$.
\epr
\begin{proof}
The  assertion follows immediately from the definitions of $\DE_{\mbG^\odot}$ and $\mr{triv}_{\mbG, (V, \partial)}$. 
\end{proof}

\subsection{The space $\mfg^{\mr{ad}(q_1)}$}\label{QQ0560}
Let us consider the space $\mfg^{\mr{ad}(q_1)}$ 
\index{$\mfg$, Lie algebra@--- $\mfg^{\mr{ad}(q_1)}$, $\mfg^{\mr{ad}(q_{-1})}$, $\mfg_{\mr{ad}(q_1)}$, $\mfg_{\mr{ad}(q_{-1})}$}
of $\mr{ad}(q_{1})$-invariants, i.e., 
\begin{equation} \label{QQ800}
 \mfg^{\mr{ad}(q_1)} := \left\{ x \in \mfg \, \big| \, \mr{ad}(q_1)(x) =0\right\}.
\end{equation}
The $k$-vector space 
$\mfg^{\mr{ad}(q_1)}$ has dimension equal to $\mr{rk}(\mfg)$ (cf. ~\cite{BD2}, \S\,3.4).
The decomposition  $\mfg = \bigoplus_{j \in \mbZ} \mfg_j$ (cf. (\ref{grading})) yields a decomposition  
\begin{equation}
  \mfg^{\mr{ad}(q_1)} = \bigoplus_{j \in \mbZ}  \mfg^{\mr{ad}(q_1)}_j
  \label{ad; grading}
  \end{equation}
on $\mfg^{\mr{ad}(q_1)}$.
Since $\{ q_{-1}, 2 \check{\rho}, q_1 \}$ forms an $\mfs \mfl_2$-triple and either ``\,$\mr{char}(k) =0$'' or ``\,$\mr{char}(k)  > 2 h_\mbG$'' is fulfilled, 
the $j$-th component $\mfg_j$ of $\mfg$ (for each $j \in \mbZ$) can be decomposed into the direct sum
\begin{align} \label{QR02}
 \mfg_j = \mr{ad}(q_{-1})(\mfg_{j+1}) \oplus \mfg^{\mr{ad}(q_1)}_j
\end{align}
(cf. ~\cite{Ngo}, the comment after Chap.\,1,  \S\,1.2, Lemme 1.2.1).
The subspace $\mfg^{\mr{ad}(q_1)}$ of $\mfg$ and 
its  components
 $ \mfg^{\mr{ad}(q_1)}_j$  for all $j$ are closed under the adjoint action  of $\mbB^\odot$ via $\iota_\mbB$ and satisfies that $ \mfg^{\mr{ad}(q_1)}_j = 0$ if $j \leq 0$ (cf. ~\cite{BD2}, \S\,3.4).
 
 If $\mfh$ is a $k$-vector space equipped with a $\mbG_m$-action $(h, v) \mapsto {^h}v$ ($h \in \mbG_m$, $v \in \mfh$),
 then we shall write $\mfh (1)$ \index{$\mfh (1)$} for the $k$-vector space $\mfh$ equipped with the  $\mbG_m$-action given by $(h, v) \mapsto h \cdot {^h}v$.
Let us identify $\mbB^\odot / \mbN^\odot$ with $\mbG_m$ via the adjoint action $\mbB^\odot / \mbN^\odot \isom \mr{GL}(\mfn^\odot) = \mbG_m$.
The $\mbB^\odot$-action on $\mfg^{\mr{ad}(q_1)}$ factors through the surjection $\mbB^\odot \migisurj \mbG_m$, so it induces a $\mbG_m$-action on $\mfg^{\mr{ad}(q_1)}$.
In particular, we obtain a $k$-vector space $\mfg^{\mr{ad}(q_1)}(1)$ equipped with a $\mbG_m$-action.

\subsection{The  vector bundle $\DV_{\mfg}$}\label{QG9000}

By the above discussion, we obtain a vector bundle $\mfg^{\mr{ad}(q_1)}_{\DE_{\mbB^\odot}}$ on $U$ of rank equal to $\mr{rk}(\mfg)$.
Write
\begin{equation}   \label{QQ801}
\DV_{\mfg, U^\mr{log}/S^\mr{log}}  \left(\text{or simply} \  \DV_{\mfg}\right)
:=  \Omega_{U^\mr{log}/S^\mr{log}}\otimes  \mfg^{\mr{ad}(q_1)}_{\DE_{\mbB^\odot}}.
\index{$\DV_{\mfg, U^\mr{log}/S^\mr{log}}$, $\DV_{\mfg}$}
\end{equation}
If $\Omega_{U^\mr{log}/S^\mr{log}}^\times$ denotes the $\mbG_m$-bundle corresponding to $\Omega_{U^\mr{log}/S^\mr{log}}$,
 then $\DV_{\mfg}$ is isomorphic to the vector bundle  $(\Omega_{U^\mr{log}/S^\mr{log}}^\times) \times^{\mbG_m}\mfg^{\mr{ad}(q_1)}(1)$.
For example, $\DV_{\mfg^\odot}$ is canonically isomorphic to  $\Omega_{U^\mr{log}/S^\mr{log}}^{\otimes 2}$ (cf. \S\,\ref{QQ0490}).

The decomposition   (\ref{ad; grading})
induces  
a decomposition 
$\mcV_{\mfg} = \bigoplus_{j \in \mbZ}  \mcV_{\mfg, j}$, which gives 
a  decreasing filtration $\{ \mcV^j_{\mfg} \}_{j \in \mbZ}$ on $\mcV_\mfg$
by putting  $\mcV^j_{\mfg} := \bigoplus_{j' \geq j}  \mcV_{\mfg, j'}$.
According to  the discussion in \S\, \ref{y037}, we have 
\begin{align} \label{QR020}
\DV_{\mfg, j} \cong \Omega_{U^\mr{log}/S^\mr{log}}^{\otimes (j+1)} \otimes_k \mfg^{\mr{ad}(q_1)}_j \left(=(\Omega_{U^\mr{log}/S^\mr{log}}^\times) \times^{\mbG_m}\mfg_j^{\mr{ad}(q_1)} (1)\right)
\end{align}
 for every $j \in \mbZ$.
 
Next, the $\mbB^\odot$-equivariant inclusion $\mfg^{\mr{ad}(q_1)} \migiincl \mfg$ yields
an $\mcO_U$-linear injection
\begin{equation}   
\varsigma   : \DV_{\mfg} \migiincl \left( \Omega_{U^\mr{log}/S^\mr{log}}\otimes   \mfg_{\DE_{\mbB^\odot}}= \right)\Omega_{U^\mr{log}/S^\mr{log}}\otimes   \mfg_{\DE_{\mbB}}. \label{inclV}
\index{$\varsigma$}
\end{equation}
Also,
since $q_1 \left(= \iota_\mfg (q_1^\odot) \right)$ belongs to $\mfg^{\mr{ad}(q_1)}$, the injection $\iota_\mfg$ restricts to  
an injection  $(\mfg^\odot)^{\mr{ad}(q_1^\odot)} \migiincl \mfg^{\mr{ad}(q_1)}$, and hence,  induces 
 an $\mcO_U$-linear injection
\begin{equation}  
\iota_{\mfg}^\mcV :  \Omega_{U^\mr{log}/S^\mr{log}}^{\otimes 2}   \migiincl \DV_{\mfg}  \label{inclOmega}
\index{$\iota_{\mfg}^\mcV$}
\end{equation}
under the natural identification $\Omega_{U^\mr{log}/S^\mr{log}}^{\otimes 2} = \DV_{\mfg^\odot}$.
By using $\varsigma$ (resp., $\iota_{\mfg}^\mcV$), we shall regard $\DV_{\mfg}$ (resp., $\Omega_{U^\mr{log}/S^\mr{log}}^{\otimes 2}$) as an $\mcO_U$-submodule of $\Omega_{U^\mr{log}/S^\mr{log}}\otimes   \mfg_{\DE_{\mbB}}$ (resp., $\DV_{\mfg}$).

\begin{rema} \label{QR69}
Note that the $\mbB^\odot$-action on $\mfg^{\mr{ad}(q_1)}$ factors through the quotient $\mbB^\odot \migisurj \left(\mbB^\odot/\mbN^\odot = \right) \mbT^\odot$.
Hence, we have $\mfg^{\mr{ad}(q_1)}_{\DE_{\mbB^\odot}} = \mfg^{\mr{ad}(q_1)}_{\DE_{\mbB^\odot} \times^{\mbB^\odot} \mbT^\odot}$. 
On the other hand, as mentioned in Remark \ref{QQ0588},
the formation of the $T^\odot$-bundle $\DE_{\mbB^\odot} \times^{\mbB^\odot} \mbT^\odot$ does not depend on the choice of $\hslash$.
It follows that
the formation of the vector bundle $\DV_\mfg$ does not depend on the choice of $\hslash$.
To be precise, if $\DV'_\mfg$ denotes the vector bundle ``$\DV_\mfg$'' defined by using another $\hslash'$ instead of $\hslash$, then (by means of (\ref{QW29}) for  $n=2$) we have a canonical isomorphism $\DV_\mfg \isom \DV'_\mfg$.
\end{rema}

\section{$\ \qq$-normality and a torsor structure on the sheaf of $(\mfg, \hslash)$-opers} \label{y047}

\subsection{$\qq$-normal $(\mfg, \hslash)$-opers} \label{QG460}

We shall define the notions of a $\ \qq$-normal $(\mfg, \hslash)$-oper  formulated   in  the local and global situations respectively.

\bde[Local case] \label{y0252}
Suppose that there exists a global log chart $(U, \partial)$ on $U^\mr{log}/S^\mr{log}$.
Let $\nabla$ be a local $(\mfg, \hslash)$-oper on $U^\mr{log}/S^\mr{log}$.
Then, we shall say that  $\nabla$ is {\bf $\ \qq$-normal relative to $(U, \partial)$} \index{local $(\mfg, \hslash)$-oper@--- $\ \qq$-normal relative to $(U, \partial)$}
 if the equality 
\begin{align} \label{QQ0601}
\nabla^\mfg (\partial) = v + 1 \otimes q_{-1}
\end{align}
holds  for some $v \in \Gamma (U, \mcO_U \otimes_k \mfg^{\mr{ad}(q_1)})$.
(In particular, a $\ \qq$-normal local $(\mfg, \hslash)$-oper  relative to $(U, \partial)$ is $\ \qq$-prenormal relative to the same log chart).
  \ede


\bde[Global case] \label{QQ0998}
Let $\msE^\spadesuit  := (\mcE_\mbB, \nabla)$ be a $(\mfg, \hslash)$-oper on $U^\mr{log}/S^\mr{log}$.
Then, we shall say that $\msE^\spadesuit$ is 
\index{$\msE^\spadesuit$, $(\mfg, \hslash)$-oper@--- $\, \qq$-normal}
 {\bf  $\, \qq$-normal} if $\mcE_\mbB = \DE_\mbB$ and,
for any log chart $(V, \partial)$ on $U^\mr{log}/S^\mr{log}$,
the pull-back $\mr{triv}_{\mbG, (V, \partial)}^*(\nabla)$, forming a local $(\mfg, \hslash)$-oper on $V^\mr{log}/S^\mr{log}$,  is  $\ \qq$-normal relative to $(V, \partial)$.
\ede

\begin{rema}[A different choice of $\ \qq$] \label{QSfff17}
Let $\ \qq':= \{ x'_\alpha \}_{\alpha \in \Gamma}$ be another choice of generators of $\mfg^\alpha$'s.
Note that there exists an element $h \in \mbT$ with $\mr{Ad}_\mbG (h) (x_\alpha) = x'_\alpha$ for every $\alpha \in \Gamma$.
Let  $\iota'_\mbB$ denote the morphism resulting from Proposition \ref{py040} where  ``$\ \qq$'' is  replaced by  $\ \qq'$.
Then,  the equality $\mr{Inn}_h \circ \iota_\mbB = \iota'_\mbB$ holds.
This implies that
if $\DE'_\mbB$ denotes the $\mbB$-bundle obtained by using $\iota'_\mbG$ in the same manner as $\DE_\mbB$, then it is canonically isomorphic to  $\DE_\mbB \times^{\mbB, \mr{Inn}_h}\mbB$.  
Moreover, for a $\ \qq$-normal  $(\mfg, \hslash)$-oper $\msE^\spadesuit := (\DE_\mbB, \nabla)$, the  $(\mfg, \hslash)$-oper $\mr{Inn}_{h*}(\msE^\spadesuit) := (\DE'_\mbB, \mr{Inn}_{h*}(\nabla))$ turns out to be  $\,  \qq'$-normal.
Thus, 
the definition of $\ \qq$-normality is independent of the choice of $\ \qq$ up to gauge transformation by an element of  $\mbT$.
\end{rema}

\begin{rema}[The case of $\mfg = \mfs \mfl_2$]
It is clear that a $(\mfg^\odot, \hslash)$-oper is $\ \qq$-normal if and only if it is normal in the sense of Definition \ref{y0132}.
\end{rema}

\begin{rema} \label{QR01}
Suppose that there exists a global log chart $(U, \partial)$ on $U^\mr{log}/S^\mr{log}$.
Then, a local $(\mfg, \hslash)$-oper $\nabla$ is $\, \qq$-normal relative to $(U, \partial)$ (in the sense of Definition \ref{y0252}) if and only if  the pair   $(U \times_k \mbB, \nabla)$ becomes a  $\ \qq$-normal  $(\mfg, \hslash)$-oper   (in the sense of Definition \ref{QQ0998}) under the identification $U \times_k \mbB = \DE_\mbB$ given by $\mr{triv}_{\mbG, (U, \partial)}$.
\end{rema}
\bpr \label{QR60}
Let $\msE^\spadesuit_\odot$ be a $(\mfg^\odot, \hslash)$-oper on $U^\mr{log}/S^\mr{log}$.
If $\msE^\spadesuit$ is $\ \qq^\odot$-normal, then the $(\mfg, \hslash)$-oper
$\iota_{\mbG*}(\msE^\spadesuit)$ (cf. (\ref{associated})) is $\ \qq$-normal. 
\epr
\begin{proof}
The assertion follows from the definitions of $\iota_\mbB$ and $\, \qq$-normality. 
\end{proof}

\subsection{$\, \qq$-normalization} \label{QR060}

The following proposition asserts that every isomorphism class of a $(\mfg, \hslash)$-oper  contains exactly one $\, \qq$-normal $(\mfg, \hslash)$-oper.

\bpr  \label{y049} 
Let   $\msE^\spadesuit$
  be  a $(\mfg, \hslash)$-oper
  on $U^\mr{log}/S^\mr{log}$.
  Then, there exists uniquely 
  a pair 
  \begin{align} \label{QQ830}
  ({^\dagger}\msE^{\spadesuit}, \mr{nor}_{\msE^\spadesuit})
  \end{align}
   consisting of 
  a $\, \qq$-normal  $(\mfg, \hslash)$-oper 
${^\dagger}\msE^{\spadesuit}$
 on $U^\mr{log}/S^\mr{log}$ 
and an isomorphism $\mr{nor}_{\msE^\spadesuit} : {^\dagger}\msE^{\spadesuit} \isom \msE^\spadesuit$
 of $(\mfg, \hslash)$-opers.
 Moreover, the well-defined assignment $\msE^\spadesuit \mapsto ({^\dagger}\msE^{\spadesuit}, \mr{nor}_{\msE^\spadesuit})$  is functorial with respect to
   pull-back to   \'{e}tale $U$-schemes, as well as base-change to fs $S^\mr{log}$-log schemes.
\epr
\begin{proof}
Let us write $\msE^\spadesuit := (\mcE_\mbB, \nabla)$.
Since $\mcO p_{\mfg, \hslash}$ is a sheaf (cf. \S\,\ref{QR30}), 
the statement is  of local nature.
Thus,  after possibly pulling-back $\msE^\spadesuit$ to an \'{e}tale $U$-scheme,   we may assume that $\mcE_\mbB =  U \times_k \mbB$ and there exists a global  log chart $(U,\partial)$ on $U^\mr{log}/S^\mr{log}$ (hence we can identify $\DE_\mbG$ with $U \times_k \mbG$ by means of $\mr{triv}_{\mbG, (U, \partial)}$).
Moreover, by Proposition  \ref{y039},
 it suffices to prove the assertion of the case where
 the local $(\mfg, \hslash)$-oper $\nabla$
  is $\, \qq$-prenormal relative to
  $(U, \partial)$.
 
 For each $\hslash$-$S^\mr{log}$-connection $\nabla_0$ on $U \times_k \mbG$, we set
 \begin{align}\label{QR07}
 N (\nabla_0) := \mr{max}\left\{ j \in \mbZ \, \Big| \, \nabla_0^{\mfg}(\partial) -1 \otimes q_{-1} \in \Gamma(U,  \mcO_U \otimes_k (\mfg^{\mr{ad}(q_1)} +\mfg^j))  
  \right\},
\end{align}
 where $N (\nabla_0) := \infty$ if $\nabla_0^{\mfg}(\partial )-1 \otimes q_{-1} \in \Gamma(U,  \mcO_U \otimes_k \mfg^{\mr{ad}(q_1)})$.
Assume  now  that 
$j_0 := N (\nabla) < \infty$.
The integer $j_0$ is positive and 
 satisfies the equalities $\mfg_{j_0} = \mr{ad}(q_{-1})(\mfg_{j_0+1}) \oplus \mfg^{\mr{ad}(q_1)}_{j_0}$ (cf. (\ref{QR02})) and $\mfg_{j_0+1} \cap \mr{Ker}(\mr{ad}(q_{-1}))=0$.
 Hence,
the composite  
\begin{align} \label{QR09}
\mfg_{j_0 +1} \xrightarrow{\mr{ad}(q_{-1})}  \mfg_{j_0} \xrightarrow{\text{inclusion}}  \mfg^{j_0}  \migi 
(\mfg^{\mr{ad}(q_1)} + \mfg^{j_0}) /(\mfg^{\mr{ad}(q_1)} + \mfg^{j_0+1})
\end{align}
turns out to be an isomorphism.
  Let $v \in \Gamma (U, \mcO_U \otimes_k \mfg_{j_0 +1}) \left(= \mfg_{j_0 +1}(U) \right)$ be the image of $\nabla^\mfg (\partial)-1 \otimes q_{-1}$ via the composite
  \begin{align}
   \mcO_U \otimes_k (\mfg^{\mr{ad}(q_1)}+ \mfg^{j_0}) 
  &
   \migisurj \mcO_U \otimes_k (\mfg^{\mr{ad}(q_1)} + \mfg^{j_0}) /(\mfg^{\mr{ad}(q_1)} + \mfg^{j_0+1})  \\
  &\isom \mcO_U \otimes_k \mfg_{j_0 +1}, \notag
  \end{align}
  where the second arrow is the isomorphism induced  by (the inverse of) (\ref{QR09}).
Since $\mfg_{j_0+1}$ is contained in $\mfn$, the $U$-rational point   $h   := \mr{Exp}(v) :   U \migi \mbN$ (cf. (\ref{QR603}) for the definition of $\mr{Exp}$) can  be defined.
 By Corollary \ref{y028}, we have 
 \begin{align} \label{equi46}
\mr{L}_{h*}(\nabla)^\mfg (\partial) 
& = -\hslash \cdot d \mr{Log}(h) (\partial) + \sum_{s =0}^\infty \frac{1}{s !} \cdot \mr{ad} (\mr{Log}(h))^s (\nabla^\mfg (\partial)) \\
  & = -\hslash \cdot d v (\partial) + \sum_{s =0}^\infty \frac{1}{s !} \cdot \mr{ad} (v)^s (\nabla^\mfg (\partial)). \notag
 \end{align}
 Both the sections $-\hslash \cdot d v (\partial)$ and $\mr{ad}(v)^s (\nabla^\mfg (\partial))$ ($s \geq 2$) belong to  $\Gamma (U, \mcO_U \otimes_k \mfg^{j_0+1})$. 
 Hence, modulo $\Gamma (U, \mcO_U \otimes_k (\mfg^{\mr{ad}(q_1)}+\mfg^{j_0+1}))$, we obtain from (\ref{equi46}) the following  equalities: 
 \begin{align}
\mr{L}_{h*}(\nabla)^\mfg  (\partial) & \equiv \nabla^\mfg (\partial) + \mr{ad}(v)(\nabla^\mfg (\partial)) \\
  & \equiv   \nabla^\mfg (\partial) - \mr{ad}(\nabla^\mfg (\partial)) (v) \notag \\
  & \equiv  \nabla^\mfg (\partial) - \mr{ad}(1 \otimes q_{-1}) (v) \notag  \\
 & \equiv \nabla^\mfg (\partial) - (\nabla^\mfg (\partial) - 1\otimes q_{-1}) \notag\\
 & \equiv 1 \otimes q_{-1}, \notag
 \end{align}
 where the third ``$\equiv$'' follows from the condition that $\nabla$ is $\, \qq$-prenormal  relative to $(U, \partial)$, and  the fourth ``$\equiv$'' follows from the definition of $v$.
This implies that $\nabla' := \mr{L}_{h_*}(\nabla)$ is a local $(\mfg, \hslash)$-oper (which is $\, \qq$-renormal relative to $(U, \partial)$) satisfying the inequality
$\left(j_0=\right) N (\nabla) < N (\nabla')$.
Moreover, the left-translation $\mr{L}_{h^{-1}}$ by $h^{-1}$ specifies an isomorphism  $(U \times_k \mbB, \nabla') \isom (U \times_k \mbB, \nabla)$ of $(\mfg, \hslash)$-opers.
By applying  successively the above construction of  $\nabla'$ starting from  $\nabla$, one obtains
a local $(\mfg, \hslash)$-oper ${^\dagger}\nabla$ on $U^\mr{log}/S^\mr{log}$ with
$N ({^\dagger}\nabla) = \infty$ (i.e., being $\, \qq$-normal) together with an isomorphism   $(U \times_k \mbB, {^\dagger}\nabla) \isom (U \times_k \mbB, \nabla)$ of $(\mfg, \hslash)$-opers.
 This proves the existence assertion.

The  uniqueness and the functoriality follow from the above discussion.
This completes the proof of Proposition \ref{y049}.
\end{proof}

\bco \label{QQ820}
 The underlying $\mbB$-bundle of any $(\mfg, \hslash)$-oper on $U^\mr{log}/S^\mr{log}$ may be  trivialized  Zariski locally on $U$.
\eco
\begin{proof}
The assertion follows from Proposition \ref{y049} and  
the existence of trivializations $\mr{triv}_{\mbG, (V, \partial)} |_{U \times_k \mbB} : U \times_k \mbB \isom \DE_\mbB |_V$ defined for various  log charts $(V, \partial)$ such that $V$ are open subschemes of $U$. 
\end{proof}

\bco \label{QS480}
Suppose that there exists a global log chart $(U, \partial)$ on $U^\mr{log}/S^\mr{log}$.
Let us consider, in the natural manner,   $\mcO_U$ as  a sheaf on $\mfE \mft_{/U}$.
Given an \'{e}tale $U$-scheme $V$  and  an element  $v \in \Gamma (V, \mcO_U \otimes_k \mfg^{\mr{ad}(q_1)})$,
we shall write $\nabla_v$ for the $\hslash$-$S^\mr{log}$-connection on $V \times_k \mbG$ determined  by the condition that  $\nabla_v^\mfg (\partial) =1 \otimes q_{-1} + v$.
Then, the assignment $v \mapsto (V \times_k \mbG, \nabla_v)$ defines an isomorphism
\begin{align} \label{QS700}
\mr{triv}_{(U, \partial)}^{\mcO p} : 
\mcO_U \otimes_k \mfg^{\mr{ad}(q_1)} \isom \mcO p_{\mfg, \hslash}
\end{align}
of \'{e}tale sheaves on $U$.
\eco
\begin{proof}
The assertion is a direct consequence of Proposition \ref{y049}.
\end{proof}

\subsection{A $\DV_\mfg$-torsor structure on $\mcO p_{\mfg, \hslash}$}
Let  us take a global section  $R$ of $\Omega_{U^\mr{log}/S^\mr{log}}\otimes   \mfb_{\DE_{\mbB}}$ and regard it as an
$\mcO_U$-linear morphism 
$\mcT_{U^\mr{log}/S^\mr{log}} \migi  \mfb_{\DE_{\mbB}} \left(\subseteq \widetilde{\mcT}_{\DE^{\mr{log}}_{\mbG}/S^\mr{log}} \right)$.
Then,  for 
a $\, \qq$-normal $(\mfg, \hslash)$-oper   $\msE^\spadesuit := (\DE_\mbB, \nabla)$ on $U^\mr{log}/S^\mr{log}$,
the pair
\index{$\msE^\spadesuit$, $(\mfg, \hslash)$-oper@--- $\msE^\spadesuit_{+R}$}
\begin{align} \label{QR036}
\msE^\spadesuit_{+R} := (\DE_\mbB, \nabla+R)
\end{align}
forms a $(\mfg, \hslash)$-oper.
One verifies that $\msE^\spadesuit_{+R}$ is $\ \qq$-normal if and only if $R$ lies in $\Gamma (U, \DV_\mfg)$.
Also, any  $\, \qq$-normal  $(\mfg, \hslash)$-oper $\msE^\spadesuit$ on $U^\mr{log}/S^\mr{log}$  may be obtained in this way, i.e.,  coincides with  $\msE^\spadesuit_{+R}$ for some $R \in \Gamma (U, \DV_\mfg)$.
It follows  that the set of  $\, \qq$-normal  $(\mfg, \hslash)$-opers on $U^\mr{log}/S^\mr{log}$, identified  with the set $\mcO p_{\mfg, \hslash} (U)$ because of   Proposition   \ref{y049},  forms a torsor modeled on $\Gamma (U, \DV_\mfg)$.
This torsor structure is compatible with pull-back to  \'{e}tale $U$-schemes, as well as base-change to fs $S^\mr{log}$-log schemes.
In particular,  by regarding $\DV_\mfg$ as  a sheaf on $\mfE \mft_{/U}$, 
we obtain a structure of $\DV_\mfg$-torsor
\begin{align} \label{QR043}
\mcO p_{\mfg, \hslash} \times \DV_\mfg \migi \mcO p_{\mfg, \hslash}
\end{align}
on the sheaf $\mcO p_{\mfg, \hslash}$.
By construction, the following diagram is commutative:
\begin{align} \label{QR050}
\vcenter{\xymatrix@C=46pt@R=36pt{
\mcO p_{\mfg^\odot, \hslash} \times \Omega_{U^\mr{log}/S^\mr{log}}^{\otimes 2}
  \ar[r]^-{(\ref{QQ0500})} \ar[d]_-{\iota_{\mfg}^{\mcO p} \times \iota_{\mfg}^\mcV}& \mcO p_{\mfg^\odot, \hslash} \ar[d]^-{\iota_{\mfg}^{\mcO p}}\\
\mcO p_{\mfg, \hslash} \times \DV_\mfg  \ar[r]_-{(\ref{QR043})} & \mcO p_{\mfg, \hslash}.
}}
\end{align}
The following  assertion follows from the above discussion together with the fact that $\mcO p_{\mfg^\odot, \hslash}$ is an $\Omega_{U^\mr{log}/S^\mr{log}}^{\otimes 2}$-torsor with respect to the action given by (\ref{QQ0500}).

\bpr \label{QQ1089}
The sheaf $\mcO p_{\mfg, \hslash}$  may be represented by a relative affine space over $U$ modeled on $\mbV (\DV_\mfg)$.
Moreover, the morphism of sheaves
\begin{align}
  \mcO p_{\mfg_\odot, \hslash} \times^{\Omega_{U^\mr{log}/S^\mr{log}}^{\otimes 2}} \DV_{\mfg} \migi \mcO p_{\mfg, \hslash}.
 \end{align}
 arising from $\iota_\mfg^{\mcO p}$  and   (\ref{QR043}) together with the commutativity of (\ref{QR050}) is an isomorphism.
\epr

\section{$(\mfg, \hslash)$-opers on pointed stable curves} \label{QQ1059}

\subsection{}
Let $r$ be  a nonnegative integer,  $S$  a $k$-scheme, and $\msU = (f: U \migi S, \{ \sigma_i^U: S_i \migi U \}_{i=1}^r)$ a local $r$-pointed stable curve over $S$.
By choosing a base for $\mcU$,  we obtain log structures on $U$ and $S$, by which $U/S$ extends to a log curve $U^\mr{log}/S^\mr{log}$.
In this situation,  we occasionally use the notation $\DV_{\mfg, \msU}$ to  denote the vector bundle $\DV_{\mfg, U^\mr{log}/S^\mr{log}}$ if there is no fear of confusion.
 Since $\mr{dim}(\mfg^{\mr{ad}(q_1)}) \left(= \sum_{j \in \mbZ} \mr{dim}(\mfg^{\mr{ad}(q_1)}_j) \right) = \mr{rk}(\mfg)$, we have
\begin{align}  \label{Vtriv}
 \sigma^{U^*}_i (\DV_\mfg) 
& = \bigoplus_{j \in \mbZ} \sigma^{U*}_i (\DV_{\mfg, j}) \\
& \hspace{-2mm}\stackrel{(\ref{QR020})}{\cong}
\bigoplus_{i=1}^r \bigoplus_{j \in \mbZ}\sigma^{U*}_i(\Omega_{U^\mr{log}/S^\mr{log}}^{\otimes (j+1)}) \otimes_k \mfg^{\mr{ad}(q_1)}_j \notag \\
& \cong  \bigoplus_{j \in \mbZ}\mcO_{S_i} \otimes_k  \mfg^{\mr{ad}(q_1)}_j \notag \\
&\cong  \mcO_{S_i}^{\oplus \mr{rk}(\mfg)}, \notag
 \end{align}
where the penultimate ``$\cong$'' follows from the residue maps $\mr{Res}_{\msU, i}: \sigma^{U*}_i (\Omega_{U^\mr{log}/S^\mr{log}}) \isom \mcO_{S_i}$ (cf. (\ref{QQ040})). 

\subsection{Cohomology of $\DV_\mfg$}
 \label{QR024}
Let us consider the case of  a (non-local) pointed stable curve.
Let $(g,r)$ be a pair of nonnegative integers with $2g-2 +r >0$, 
$S$  a $k$-scheme,  and $\msX := (f: X \migi S, \{ \sigma_i \}_{i=1}^r)$
an $r$-pointed stable curve of genus $g$ over $S$.
We shall write
 \begin{align}
 &  \hspace{14mm} \label{aleph}\aleph (\mfg)  := \frac{2g-2+r}{2}\cdot \mr{dim}(\mfg) + \frac{r}{2}\cdot \mr{rk}(\mfg) 
 \index{$\aleph (\mfg)$, ${^c \aleph} (\mfg)$}
 \\
& \left(\text{resp.,} \ {^c \aleph} (\mfg) := \frac{2g-2+r}{2}\cdot \mr{dim}(\mfg) - \frac{r}{2}\cdot \mr{rk}(\mfg) \right). \notag \end{align}
Both $\aleph (\mfg)$ and  ${^c \aleph} (\mfg)$  are verified to be  integers.

\bpr \label{y046}
Recall (cf. (\ref{QQ979})) that
$D_\msX$ denotes the relative effective divisor  on $X$ defined as the union of the images of the $\sigma_i$'s.
Then, the following assertions hold:
\begin{itemize}
\item[(i)]
$\mbR^1 f_*(\DV_{\mfg, \msX}) = \mbR^1 f_*(\DV_{\mfg, \msX}(-D_\msX)) = 0$.
\item[(ii)]
 $f_*(\DV_{\mfg, \msX})$ (resp., $f_*(\DV_{\mfg, \msX} (-D_\msX))$)  is a vector bundle on $S$ of rank $\aleph (\mfg)$ (resp., ${^c \aleph} (\mfg)$).
 Moreover, if $s : S' \migi S$ is a morphism of $k$-schemes, then the natural  $S'$-morphism 
\begin{align} \label{bc}
 & \hspace{13mm} S' \times_S\mbV (f_*(\DV_{\mfg, \msX}))   \migi \mbV ((\mr{id}_{S'} \times f)_*(\DV_{\mfg, \msX'})) \\
& \hspace{-1mm} \left(\text{resp}.,   \ S' \times_S\mbV(f_*(\DV_{\mfg, \msX} (-D_\msX)))  \migi \mbV ((\mr{id}_{S'} \times f)_*(\DV_{\mfg, \msX'}(-D_{\msX'})))\right) \notag
\end{align}
 is an isomorphism, where $\msX' := S' \times_S \msX$.
\end{itemize}
 \epr
\begin{proof}
 Because of (\ref{QR020}) and the fact that $\mfg_j^{\mr{ad}(q_1)} =0$ for $j \leq 0$, 
 we see that $\DV_{\mfg}$ and $\DV_{\mfg} (-D_\msX)$ are isomorphic to direct sums of $\Omega_{X^\mr{log}/S^\mr{log}}^{\otimes l}$'s and $\Omega_{X^\mr{log}/S^\mr{log}}^{\otimes l}(- D_{\msX})$'s (for various $l \geq 2$) respectively.
 But, the equalities 
\begin{equation}\mbR^1f_* (\Omega_{X^\mr{log}/S^\mr{log}}^{\otimes l})= \mbR^1f_* (\Omega^{\otimes l}_{X^\mr{log}/S^\mr{log}}(-D_{\msX})) =0\end{equation}
 hold for any $l \geq2$.
It follows that 
\begin{equation} \label{zero1}
 \mbR^1 f_*(\DV_{\mfg, \msX}) = \mbR^1 f_*(\DV_{\mfg, \msX} (-D_\msX)) = 0, \end{equation}
thus completing   the proof of assertion (i).
Moreover,  (\ref{zero1})  implies (cf. ~\cite{Har}, Chap.\,III, Theorem 12.11 (b)) that 
both $f_*(\DV_{\mfg, \msX})$ and $f_*(\DV_{\mfg, \msX} (-D_\msX))$ are   vector bundles on $S$ and  that the morphisms displayed in  (\ref{bc}) are  isomorphisms.

In what follows, let us   compute the ranks of $f_*(\DV_{\mfg, \msX})$ and  $f_*(\DV_{\mfg, \msX} (-D_\msX))$.
First, we
refer to  the fact  that 
the isomorphism $\mr{Kos}_{\mfg} : \mfg^{\mr{ad}(q_1)} \isom \mfc$ defined in \S\,\ref{y054} (cf. (\ref{Kosg})) is $\mbG_m$-equivariant  with respect to  the  $\mbG_m$-action  on $\mfg^{\mr{ad}(q_1)}(1)$ and the $\mbG_m$-action on $\mfc$ that comes from the homotheties on $\mfg$ (i.e., the $\mbG_m$-action $\star$).
Hence,  if $e_l$ ($l = 1, \cdots, \mr{rk}(\mfg)$) denote
the degrees of homogeneous generators $u_l$ in $\mbS_k (\mfg^\vee)^\mbG$,
then $ \DV_{\mfg, \msX}$ is isomorphic to $\bigoplus_{i =1}^{\mr{rk}(\mfg)} \Omega_{X^\mr{log}/S^\mr{log}}^{\otimes e_i}$.
The rank of $f_*( \DV_{\mfg, \msX} )$ can be computed by the following sequence of equalities:
 \begin{align} \label{eq34}
 \mr{rk}(f_*( \DV_{\mfg, \msX} )) & = \sum_{i=1}^{\mr{rk}(\mfg)} \mr{rk}(f_*(\Omega_{X^\mr{log}/S^\mr{log}}^{\otimes e_i}))  \\
 & =  \sum_{i=1}^{\mr{rk}(\mfg)} ((g-1)  (2  e_i -1) + r \cdot e_i) \notag  \\
& = (g-1)\cdot \mr{dim}(\mfg) + \frac{r}{2} \cdot (\mr{dim}(\mfg) + \mr{rk}(\mfg)) \notag \\
&= \aleph (\mfg),
 \notag \end{align}
where the third equality follows from  the equality
 $\sum_{i=1}^{\mr{rk}(\mfg)} (2  e_i -1) = \mr{dim}(\mfg)$ (cf. ~\cite{BD1}, 2.2.4, Remark (iv)).
On the other hand,   
by applying the functor $\mbR^1f_*(-)$ to the short exact sequence
\begin{equation} 
0 \longmigi  \DV_{\mfg, \msX} (-D_\msX) \longmigi \DV_{\mfg, \msX} \longmigi \bigoplus_{i=1}^r \mcO_{S_i}^{\oplus \mr{rk}(\mfg)} \longmigi 0 \end{equation}
(cf. (\ref{Vtriv})), we obtain (by taking account of  (\ref{zero1})) an equality
\begin{equation}  \label{eq33}
 \mr{rk}(f_*( \DV_{\mfg, \msX} )) = \mr{rk}(f_*(\DV_{\mfg, \msX} (-D_\msX))) + r \cdot \mr{rk}(\mfg).\end{equation}
Thus, (\ref{eq33}) and (\ref{eq34}) together  imply
  the required computation for the rank of $f_*(\DV_\mfg (-D_\msX))$.
 This completes the proof of  assertion (ii).
 \end{proof}

\subsection{The moduli stack of $(\mfg, \hslash)$-opers} \label{QG9003}
Since the formation of a $(\mfg, \hslash)$-oper is  closed under base-change (cf. \S\,\ref{QQ46}), 
one may define the $\mfS \mfe \mft$-valued contravariant functor 
 \index{$\mfO \mfp_{\mfg, \hslash, (\rho,) \msX}$, $\mfO \mfp_{\mfg, (\hslash,) (\rho,) g,r}$}
\begin{equation}
\label{opermoduli}
 \mfO \mfp_{\mfg, \hslash, \msX} : \mfS \mfc \mfh_{/S} \migi \mfS \mfe \mft 
 \end{equation}
on $\mfS \mfc \mfh_{/S}$ which, to any $S$-scheme $S'$, assigns the set of isomorphism classes of $(\mfg, \hslash)$-opers on $S' \times_S \msX$. 
In particular, we have $\mcO p_{\mfg, \hslash, X^\mr{log}/S^\mr{log}}(X) = \mfO \mfp_{\mfg, \hslash, \msX} (S)$. 
The torsor structure defined in (\ref{QR043}) induces 
a $\mbV (f_* (\DV_{\mfg, \msX}))$-action
\begin{align} \label{QR094}
\mr{act}_{\mfg} : \mfO \mfp_{\mfg, \hslash, \msX} \times_S \mbV (f_* (\DV_{\mfg, \msX})) \migi 
\mfO \mfp_{\mfg, \hslash, \msX}
\end{align}
on $\mfO \mfp_{\mfg, \hslash, \msX}$ compatible with base-change to $S$-schemes.
Also, $\iota_\mfg^{\mcO p}$ (cf. (\ref{opincl})) gives rise to a morphism
 \index{$\iota_{\mfg}^{\mcO p}$, $\iota_{\mfg}^{\mfO \mfp}$, ${^p}\iota_{\mfg}^{\mfO \mfp}$}
\begin{align} \label{QR066}
\iota_\mfg^{\mfO \mfp} :\mfO \mfp_{\mfg^\odot, \hslash, \msX} \migi  
\mfO \mfp_{\mfg, \hslash, \msX}
\end{align}
of functors, which makes the following square diagram commute:
\begin{align} \label{QR069}
\vcenter{\xymatrix@C=46pt@R=36pt{
 \mfO \mfp_{\mfg^\odot, \hslash, \msX} \times_S \mbV (f_* (\Omega_{X^\mr{log}/S^\mr{log}}^{\otimes 2})) 
 \ar[r]^-{\mr{act}_{\mfg^\odot}} \ar[d]_-{\iota_\mfg^{\mfO \mfp} \times \iota_\mfg^{\mbV}}& \mfO \mfp_{\mfg^\odot, \hslash, \msX} \ar[d]^-{\iota_\mfg^{\mfO \mfp}}\\
 \mfO \mfp_{\mfg, \hslash, \msX} \times_S \mbV (f_* (\DV_{\mfg, \msX})) \ar[r]_-{\mr{act}_{\mfg}} & \mfO \mfp_{\mfg, \hslash, \msX},
}}
\end{align}
where $\iota_\mfg^{\mbV}$ denotes the closed immersion $\left(\mbV (f_*(\DV_{\mfg^\odot, \msX})) =\right)\mbV (f_* (\Omega_{X^\mr{log}/S^\mr{log}}^{\otimes 2}))  \migiincl \mbV (f_* (\DV_{\mfg, \msX})) $ induced by $\iota_\mfg^\mcV$ (cf. (\ref{inclOmega})).

 The following proposition was proved in ~\cite{BD2}  \S\,3.4, Theorem (cf.  ~\cite{Bar}, Lemma 3.6), for the case of smooth curves over $\mbC$,  and  in  ~\cite{Mzk1}, Chap.\,I, Proposition 2.11,  for the case of (arbitrary pointed stable curves but) $\mfg = \mfs \mfl_2$.


\bt \label{y051} 
(Recall    
 that  either ``\,$\mr{char}(k)=0$'' or ``\,$\mr{char}(k) >2  h_\mbG$'' is fulfilled.) The following assertions hold:
\begin{itemize}
\item[(i)]
$\mfO \mfp_{\mfg, \hslash, \msX}$ forms  a relative 
affine space over $S$ modeled on  $\mbV (f_*(\DV_{\mfg, \msX}))$ with respect to the action $\mr{act}_\mfg$ given by (\ref{QR094}), in particular, may be  represented  
by  a relative affine space of relative dimension $\aleph (\mfg)$ (cf. Proposition \ref{y046}). 
\item[(ii)]
The $S$-morphism
\begin{align} \label{QR070}
\mfO \mfp_{\mfg^\odot, \hslash, \msX} \times^{ \mbV (f_*(\Omega_{X^\mr{log}/S^\mr{log}}^{\otimes 2}))}  \mbV (f_*(\DV_{\mfg, \msX})) \migi \mfO \mfp_{\mfg, \hslash, \msX}
\end{align}
arising from $\iota_\mfg^{\mfO \mfp}$ and $\mr{act}_\mfg$ (together with the commutativity of (\ref{QR069})) is an isomorphism.
In particular, $\iota_\mfg^{\mfO \mfp}$ is a closed immersion.
\end{itemize}
\et
\begin{proof}
It follows from  Proposition \ref{y046}, (i),  that \'{e}tale locally on $S$ there always exists a $(\mfg, \hslash)$-oper on $\msX$.
Hence,  
the assertions follow from Propositions \ref{QQ1089}.
\end{proof}

\begin{rema} \label{QQ891} 
 If $\mfg$ is simple, then 
 the morphism $\iota_\mfg^\mcV$ restricts to  an isomorphism 
 $\Omega_{X^\mr{log}/S^\mr{log}}^{\otimes 2}\isom \DV_{\mfg, 1}$.
 Hence, 
  (\ref{QR070}) becomes
an isomorphism
\begin{equation}
\mfO \mfp_{\mfg^\odot, \hslash, \msX} \times_S \mbV (f_*( \DV^{2}_{\mfg})) \isom \mfO \mfp_{\mfg, \hslash, \msX}
\label{gissimple}
\end{equation}
(cf.  ~\cite{BD2}, \S\,3.4, Remark).
 \end{rema}

\begin{rema} \label{QQ892} 
Let $(\mbG', \mbT')$ be another split semisimple algebraic group over $k$ of adjoint type.
If $k$ is of positive characteristic, then we suppose  that $\mr{char}(k) > 2  h_{\mbG'}$. 
Denote by $\mfg'$ the Lie algebra of $\mbG'$,  and
 fix a Borel subgroup $\mbB'$ of $\mbG'$ defined over $k$ containing $\mbT'$.
Also, let  $\msE^\spadesuit := (\mcE_\mbB, \nabla)$ be a $(\mfg, \hslash)$-oper on $\msX$ and $\msE'^{\spadesuit} := (\mcE_{\mbB'}, \nabla')$  a $(\mfg', \hslash)$-oper on $\msX$.
Then, one verifies that  $\mbB \times_k \mbB'$ is a Borel subgroup of $\mbG \times_k \mbG'$ and  the pair
\index{$\msE^\spadesuit$, $(\mfg, \hslash)$-oper@--- $\msE^\spadesuit \boxtimes \msE'^\spadesuit$}
\begin{align} \label{QS14}
\msE^\spadesuit \boxtimes \msE'^\spadesuit := (\mcE_\mbB \boxtimes \mcE_{\mbB'}, \nabla \boxtimes \nabla')
\end{align} 
(cf. (\ref{QS8}), (\ref{QS11})) forms a $(\mfg \oplus \mfg', \hslash)$-oper on $\msX$.
The assignment $(\msE^\spadesuit, \msE'^\spadesuit) \mapsto \msE^\spadesuit \boxtimes \msE'^\spadesuit$ is functorial with respect to $S$ and the resulting $S$-morphism
\begin{equation} \label{product}
\mfO \mfp_{\mfg \oplus \mfg', \hslash, \msX}  \migi \mfO \mfp_{\mfg, \hslash, \msX} \times_S \mfO \mfp_{\mfg', \hslash, \msX}
\end{equation}
is an isomorphism (cf.  ~\cite{BD2}, \S\,3.4, Remark).
 \end{rema}


\subsection{The moduli stack of $(\mfg, 0)$-opers} \label{QR94}
We shall  specialize to  the case of $\hslash =0$.
Let us refer to a result in Example \ref{y0134} described later, which says   that if $\hslash =0$, then there exists a canonical decomposition 
\begin{align} \label{QR0108}
\mfg^\odot_{\DE_{\mbB^\odot}} = (\mfg^\odot_{-1})_{\DE_{\mbB^\odot}} \oplus (\mfg^\odot_{0})_{\DE_{\mbB^\odot}} \oplus (\mfg^\odot_{1})_{\DE_{\mbB^\odot}}.
\end{align}
Let $\nabla^{\blacktriangleright}$ denote the composite
\begin{align}
\nabla^{\blacktriangleright} 
: \mcT_{X^\mr{log}/S^\mr{log}} \isom (\mfg^\odot_{-1})_{\DE_{\mbB^\odot}} \migiincl 
\mfg^\odot_{\DE_{\mbB^\odot}} \migiincl \mfg_{\DE_{\mbB}},
\end{align}
where the first arrow is the canonical isomorphism mentioned in \S\,\ref{QQ0490},
 the second arrow  denotes the  inclusion to the first factor with respect to the decomposition   (\ref{QR0108}), and the last arrow denotes the injection induced by $\iota_\mfg$. 
Then, 
it specifies
 a $0$-$S^\mr{log}$-connection on $\DE_{\mbG}$ (cf. Example \ref{QQ8}),  and 
 the pair
 \begin{align} \label{QR99}
{^\dagger}\msE^{\blacktriangleright} :=  (\DE_{\mbB}, \nabla^{\blacktriangleright})
 \end{align}
 forms  a $\, \qq$-normal  $(\mfg, 0)$-oper on $U^\mr{log}/S^\mr{log}$.
This $(\mfg, 0)$-oper determines   a global section $S \migi \mfO \mfp_{\mfg, 0, \msX}$, which gives
 a trivialization
 \begin{align} \label{QR0304}
\mbV (f_*(\DV_{\mfg, \msX})) \isom \mfO \mfp_{\mfg, 0, \msX}
\end{align}
of the relative affine space  $\mfO \mfp_{\mfg, 0, \msX}$.
This trivialization
sends each $R \in \Gamma (S, f_*(\DV_{\mfg, \msX}))$ $\left(= \Gamma (X, \DV_{\mfg, \msX})\right)$ to 
  ${^\dagger}\msE^{\blacktriangleright}_{+R}$.
Now, observe that   the $\mbG_m$-action   on $\mfg^{\mr{ad}(q_1)}(1)$ (cf. \S\,\ref{QQ0560})  induces  a $\mbG_m$-action on  $\mbV (f_*(\DV_{\mfg, \msX}))$; the opposite of 
this action corresponds to  a $\mbZ_{\geq 0}$-grading of the $\mcO_S$-algebra $\mcO_S \langle \mbV (f_*(\DV_{\mfg, \msX})) \rangle$.
Hence, it defines a $\mbZ_{\geq 0}$-grading on the $\mcO_S$-algebra  $\mcO_S \langle  \mfO \mfp_{\mfg, 0, \msX}\rangle$ via the $\mcO_S$-algebra isomorphism 
\begin{align} \label{QHH3}
\mcO_S \langle\mfO \mfp_{\mfg, 0, \msX}\rangle \isom \mcO_S \langle \mbV (f_*(\DV_{\mfg, \msX})) \rangle
\end{align}
induced by  (\ref{QR0304}),

\subsection{A filtration  on  $\mcO_S \langle \mfO \mfp_{\mfg, \hslash, \msX}\rangle$} \label{QS491}

Let us go back to the case of a general $\hslash$ and 
 construct an increasing filtration 
on the $\mcO_S$-algebra $\mcO_S \langle \mfO \mfp_{\mfg, \hslash, \msX}\rangle$.
Recall that   $\mfO \mfp_{\mfg, \hslash, \msX}$ admits   canonically a structure of relative affine space over $S$ modeled on $\mbV (f_*(\DV_{\mfg, \msX}))$.
Hence, each point of $S$ 
has an open neighborhood $S'$ on which there exists an isomorphism
\begin{align} \label{QS4066}
\left(S' \times_S \mbV (f_*(\DV_{\mfg, \msX})) =\right) \mbV (f'_*(\DV_{\mfg, \msX'}))\isom \mfO \mfp_{\mfg, \hslash,  \msX'} \left(= S' \times_S \mfO \mfp_{\mfg, \hslash, \msX}  \right),
\end{align}
where $\msX' := S' \times_S \msX$ and    $f'$ denotes the projection $S' \times_S X \migi S'$;
 it corresponds to an isomorphism
of  $\mcO_{S'}$-algebras
\begin{align} \label{QS503}
\mcO_{S'} \langle \mfO \mfp_{\mfg, \hslash, \msX'}\rangle  \isom \mcO_{S'} \langle  \mbV (f'_*(\DV_{\mfg, \msX'})) \rangle.
\end{align}
The $\mbZ_{\geq 0}$-grading on $\mcO_{S'} \langle \mbV (f'_*(\DV_{\mfg, \msX'})) \rangle$
 defines, via this isomorphism,  an increasing  filtration 
 \begin{align} \label{QS502}
 \left\{ \mcO_{S'} \langle \mfO \mfp_{\mfg, \hslash, \msX'} \rangle^{j}\right\}_{j \in \mbZ_{\geq 0}}
 \end{align}
  on $\mcO_{S'} \langle \mfO \mfp_{\mfg, \hslash, \msX'} \rangle$.
  Moreover,   
  (\ref{QS503})
   induces
    an isomorphism
  \begin{align} \label{QS501}
\mr{gr}(\mcO_{S'} \langle \mfO \mfp_{\mfg, \hslash, \msX'}\rangle)  \isom \mcO_{S'} \langle  \mbV (f'_*(\mcV_{\mfg, \msX'})) \rangle
  \end{align} 
  between the respective  $\mbZ_{>0}$-graded $\mcO_{S'}$-algebras associated to the domain and codomain. 
Both 
the formations of  (\ref{QS502}) and (\ref{QS501})
  do not depend on the choice of  (\ref{QS4066}).
 Hence,
   the filtrations $ \left\{ \mcO_{S'} \langle \mfO \mfp_{\mfg, \hslash, \msX'} \rangle^{j}\right\}_j$  for various
   open subschemes $S' \left(\subseteq S \right)$ 
     may be glued together to obtain  an increasing filtration
 \begin{align} \label{QS504}
\left\{ \mcO_S \langle \mfO \mfp_{\mfg, \hslash, \msX}\rangle^j \right\}_{j \in \mbZ_{\geq 0}}
\end{align}
 on the $\mcO_S$-algebra $\mcO_S \langle \mfO \mfp_{\mfg, \hslash, \msX}\rangle$ which has
   an isomorphism
 \begin{align} \label{QS496}
 \mr{gr} (\mcO_S \langle \mfO \mfp_{\mfg, \hslash, \msX}\rangle) \isom \mcO_S \langle \mbV (f_*(\DV_{\mfg, \msX}))\rangle 
 \end{align}
 of $\mbZ_{\geq 0}$-graded $\mcO_S$-algebras.
Composing (\ref{QHH3}) and  (\ref{QS496})  gives  the following proposition, by which we can think of $\mcO_S \langle \mfO \mfp_{\mfg, \hslash, \msX}\rangle$ as a {\it quantization} of $\mcO_S \langle \mfO \mfp_{\mfg, 0, \msX}\rangle$ in some sense:

\bpr \label{QS497}
There exists a canonical isomorphism of $\mbZ_{\geq 0}$-graded  $\mcO_S$-algebras
\begin{align} \label{QS498}
\mr{gr} (\mcO_S \langle \mfO \mfp_{\mfg, \hslash, \msX}\rangle) \isom \mcO_S \langle \mfO \mfp_{\mfg, 0, \msX}\rangle. 
\end{align}
\epr

\section{The adjoint quotient of $\mfg$}\label{y054}

In the subsequent   discussion, 
we  
  introduce the definition of {\it radius} (cf. Definition \ref{y056}) and study the moduli space of $(\mfg, \hslash)$-opers of prescribed radii.
To do this, we  first review some basic facts concerning the quotient of $\mfg$ by the adjoint $\mbG$-action.

Denote by $\mbW$ \index{$\mbW$, Weyl group of $(\mbG, \mbT)$} the Weyl group of $(\mbG, \mbT)$, and
consider the GIT quotient $\mfg \ooalign{$/$ \cr $\,/$}\hspace{-0.5mm}\mbG$  (resp., $\mft \ooalign{$/$ \cr $\,/$}\hspace{-0.5mm}\mbW$) of 
   $\mfg$ (resp., $\mft$) by 
the adjoint action of $\mbG$ (resp.,  $\mbW$), i.e., the spectrum of the ring of polynomial invariants $\mbS_k (\mfg^\vee)^\mbG$ (resp., $\mbS_k (\mft^\vee)^\mbW$) on  $\mfg$ (resp., $\mft$).
Let us use the notation 
 \begin{equation} 
 \mfc_\mfg, \ \text{or simply}  \ \mfc, \index{$\mfc$, $\mfc_\mfg$}
 \label{c}
 \end{equation}
to denote the $k$-scheme $\mft \ooalign{$/$ \cr $\,/$}\hspace{-0.5mm}\mbW$.
Also, denote by $[0]$ \index{$[0]$, $[0]^{\times r}$, $[0]_U$} the $k$-rational point of $\mfc$ defined as the image of the zero element $0 \in \mft$ via the  quotient $\mft \migisurj \mfc$.
 A Chevalley's theorem asserts (cf. ~\cite{Ngo}, Theorem 1.1.1; ~\cite{KW}, Chap.\,VI, Theorem 8.2) that
the $k$-morphism $\mfc\migi \mfg \ooalign{$/$ \cr $\,/$}\hspace{-0.5mm}\mbG$ induced by 
the natural $k$-algebra morphism $\mbS_k (\mfg^\vee)^\mbG \migi \mbS_k (\mft^\vee)^\mbW$ is an isomorphism
 because we have assumed that  either ``\,$\mr{char} (k) =0$'' or ``\,$\mr{char}(k) > 2  h_\mbG$'' is   satisfied.
 Also,  
in  $\mbS_k (\mfg^\vee)^\mbG \left(\cong \mbS_k (\mft^\vee)^\mbW \right)$, there are $\mr{rk}(\mfg)$ algebraically independent (relative to $k$) homogeneous elements, which generate  $\mbS_k (\mfg^\vee)^\mbG$.
 
Denote by 
\begin{equation} \chi :  \mfg \migi \mfc
\label{chi}\index{$\chi$, adjoint quotient of $\mfg$}
\end{equation}
the morphism of $k$-schemes defines as  the composite of the natural quotient $\mfg \migisurj \mfg \ooalign{$/$ \cr $\,/$}\hspace{-0.5mm}\mbG$ and the inverse of the resulting  isomorphism $\mfc\isom \mfg \ooalign{$/$ \cr $\,/$}\hspace{-0.5mm}\mbG$.
Here, recall  the quotient stack $[\mfg /\mbG]$  of $\mfg$  by the adjoint $\mbG$-action;  it represents  the category fibered in groupoids over  $\mfS \mfc \mfh_{/k}$ whose category  of sections over  a $k$-scheme $T$ is  the groupoid of pairs $(\mcE, v)$ consisting of a $\mbG$-bundle  $\mcE$ on $T$ and an element $v \in \Gamma (T, \mfg_\mcE)$.
Then, $\chi$  factors through the  natural quotient  $\mfg \migi [\mfg /\mbG]$.
 We shall denote by 
 \begin{equation}    [\chi] :[\mfg/\mbG] \migi \mfc \label{[chi]}
 \index{$[\chi]$, adjoint quotient of $[\mfg/\mbG]$}
 \end{equation}
 the resulting morphism.

Next, note that  $\mfc$ ($\cong \mr{Spec}(\mbS_k (\mfg^\vee)^\mbG)$) admits a canonical $\mbG_m$-action
\begin{align} \label{action}
\star : \mbG_m \times_k \mfc   \migi  \mfc;  (h,    a)  \mapsto h \star a
\index{$\star$}
    \end{align}
    coming  from the homotheties on $\mfg$.
The $k$-algebra   $\mbS_k (\mfg^\vee)$ is equipped with the grading corresponding to 
  the opposite of this $\mbG_m$-action. 
Let us select 
a choice of homogeneous generators $u_l \in \mbS_k (\mfg^\vee)^\mbG$ 
 ($l = 1, 2, \cdots, \mr{rk}(\mfg)$), which 
determines  an isomorphism of $k$-algebras  $k[u_1, \cdots, u_{\mr{rk}(\mfg)}] \isom  \mbS_k (\mfg^\vee)^\mbG$.
Under the  identification  of  $\mfc$  with  the  $\mr{rk}(\mfg)$-dimensional  affine space $\mbA^{\mr{rk}(\mfg)}$ using  the resulting  $k$-isomorphism
\begin{equation} \label{isomCA}
 \mfc \isom \mbA^{\mr{rk}(\mfg)} \left(:=\mr{Spec}(k[u_1, \cdots, u_{\mr{rk}(\mfg)}])\right),
 \end{equation} 
the action $\star$ may be expressed as 
\begin{equation} h \star (a_1, \cdots, a_{\mr{rk}(\mfg)}) = (h^{e_1} \cdot a_1, \cdots, h^{e_{\mr{rk}(\mfg)}} \cdot a_{\mr{rk}(\mfg)}),
\label{action2}
\end{equation}
where  each $e_l$ denotes the degree of $u_l$.
The action $\star$  extends uniquely    to a   $k$-morphism  $\mbA^1 \times_k \mfc \migi \mfc$; we also use the notation $\star$ to denote this morphism.
Let  $r$ be a nonnegative integer and $S$ a $k$-scheme.
 Denote by $\mfc^{\times r}$  the product of $r$ copies of $\mfc$ over $k$.
Then, for each $\hslash \in \Gamma (S, \mcO_S) \left(= \mbA^1 (S) \right)$ and  
$\rho := (\rho_i)_{i=1}^r \in \mfc^{\times r}(S)$, we write $\hslash \star \rho := (\hslash \star \rho_i)_{i=1}^r \left(\in \mfc^{\times r} (S) \right)$.
(In the case of  $r=0$,  we set  $\mfc^{\times r}(S) := \{ \emptyset \}$  and  $\rho = \hslash \star \rho = \emptyset$.)

 Finally, let
\begin{equation}
  \mr{Kos}_\mfg : \mfg^{\mr{ad}(q_1)} \migiincl \mfg \xrightarrow{\chi} \mfc   \label{Kosg}
  \index{$\mr{Kos}_\mfg$}
\end{equation}
be the composite of $\chi$ and the inclusion $\mfg^{\mr{ad}(q_1)}  \migiincl \mfg$  given by  $v \mapsto v + q_{-1}$ for any $v \in  \mfg^{\mr{ad}(q_1)}$.
Since   either ``\,$\mr{char} (k) =0$'' or ``\,$\mr{char}(k) > 2  h_\mbG$'' is   satisfied,  
 $\mr{Kos}_\mfg$ turns out to be  an isomorphism of $k$-schemes (cf. ~\cite{Ngo}, Lemma 1.2.1).
In particular, we have $\mr{dim}(\mfc) =  \mr{rk} (\mfg) = \mr{dim}(\mfg^{\mr{ad}(q_1)})$.
Moreover,   $\mr{Kos}_{\mfg}$ may be regarded as a {\it $\mbG_m$-equivariant} isomorphism  $\mfg^{\mr{ad}(q_1)}(1) \isom \mfc$.

\begin{exa} \label{QQ1290}
Let us describe the adjoint quotient  for
 $\mbG = \mr{PGL}_n$  (with $n<p$) by using  the natural identification $\mfp \mfg \mfl_n = \mfs \mfl_n$.
Under this identification, the Lie algebra $\mft$ of $T$ (cf. Example \ref{QS8376})  is given as
\begin{align} \label{QR1}
\mft := \Big\{ \mr{diag} (a_1, \cdots, a_n) \in \mfs \mfl_n\, \Big| \, a_1, \cdots, a_n \in k, \sum_{j=1}^n a_j =0\Big\},
\end{align}
where $\mr{diag} (a_1, \cdots, a_n)$ denotes the $n \times n$ diagonal  matrix with  diagonal entries $a_1, \cdots, a_n$.
For $i=1, \cdots, n$, denote by  $\lambda_i$ the element of $\mft^\vee$ defined  by  $\lambda_i (\mr{diag} (a_1, \cdots, a_n))$
$= a_i$.
Then, we have $\mbS_k (\mft^\vee) = k [\lambda_1, \cdots, \lambda_n]/(\lambda_1 + \cdots + \lambda_n)$.
The action of $\mbW = \mfS_n$ ($:=$ the symmetric group of $n$ letters) on $\mft$ is given by $\sigma \cdot \lambda_i := \lambda_{\sigma (i)}$ for each $\sigma \in \mfS_n$.
It follows that $\mbS_k (\mft^\vee)^\mbW$ coincides with the ring of symmetric polynomials in $\lambda_1, \cdots, \lambda_n$.
Now,  consider the following elementary symmetric polynomials of $\lambda_1, \cdots, \lambda_n$:
\begin{align} \label{QW305}
s_1^\mft := \sum_{j=1}^n \lambda_j, \hspace{5mm}
s^\mft_2 := \sum_{j<j'} \lambda_{j} \cdot  \lambda_{j'}, \hspace{3mm}
\cdots,   \hspace{3mm}
s^\mft_n := \prod_{j=1}^n \lambda_j.
\end{align}
The equality $s^\mft_1 =0$ holds in $\mbS_k (\mft^\vee)^\mbW$
 and    the remaining  polynomials  $s^\mft_2, \cdots, s^\mft_n$ form  algebraically independent generators of $\mbS_k (\mfh^\vee)^\mbW$, i.e., $\mbS_k (\mfh^\vee)^\mbW = k[s^\mft_1, \cdots, s^\mft_n]$. 
Let us define 
 $s^{\mfs \mfl}_j \in \mbS_k (\mfs \mfl_n^\vee)$ ($j=1, \cdots, n$) to be the polynomial functions determined  by 
 \begin{align} \label{QW304}
 \mr{det}(t\cdot E-x) = t^n + \sum_{j=1}^n (-1)^j \cdot  s^{\mfs \mfl}_j (x) \cdot t^{n-j}
 \end{align}
  for any $x \in \mfs \mfl_n$, where  $E$ denotes the unit matrix of size $n$.
These polynomials belong to  $\mbS_k (\mfs \mfl_n^\vee)^{\mr{PGL}_n}$  and satisfy that
  $s^{\mfs \mfl}_1 =0$ and $\mbS_k (\mfs \mfl_n^\vee)^{\mr{PGL}_n} = k [s^{\mfs \mfl}_2, \cdots, s^{\mfs \mfl}_n]$. 
 Moreover,  the morphism $\chi : \mfs \mfl_n \left(= \mfp \mfg \mfl_n \right)\migi \mfc_{\mfs \mfl_n} \left(=\mfc_{\mfp \mfg \mfl_n} \right)$ is given by   $s^\mft_j \mapsto s^{\mfs \mfl}_j$ ($j=2, \cdots, n$).
\end{exa}

\section{Radii of $(\mfg, \hslash)$-opers}\label{y055}

In this section, 
we   introduce the notion of radius and 
discuss   $(\mfg, \hslash)$-opers 
 of prescribed radii.

\subsection{The definition of radius} \label{QR0110}
 Let $S$ be a $k$-scheme and  
 $\msU := (U/S, \{ \sigma^U_i : S_i \migi U \}_{i=1}^r)$ a local $r$-pointed stable curve over $S$, where $r$ is assumed to be positive unless otherwise stated.
 After choosing a base for $\msU$, we 
  obtain a log curve $U^\mr{log}/S^\mr{log}$ extending  $U/S$.
  Also, let $\hslash$ be an element of $\Gamma (S, \mcO_S)$.
  
\bde\label{y056} 
Let $\mcE$ be a $\mbG$-bundle on $U$ and $\nabla$ an $\hslash$-$S^\mr{log}$-connection on $\mcE$.
For each 
$i =1, \cdots, r$, we 
write
 \index{$\nabla$, $\hslash$-$T^\mr{log}$-connection@--- $\rho_i^{\nabla}$, radius}
\begin{equation}  \label{QQ1000}
 \rho_i^{\nabla} \in \mfc (S_i) 
 \end{equation}
for the image, via $[\chi] : [\mfg /\mbG] \migi \mfc$, of the $S_i$-rational point of $[\mfg /\mbG]$ classifying the pair $(\sigma_i^{U*}(\mcE), \mu_i^{\nabla})$ (cf. Definition \ref{y031} for the definition of $\mu_i^{\nabla}$).
We shall refer to $ \rho_i^{\nabla}$ as the {\bf radius}
 \index{$\nabla$, $\hslash$-$T^\mr{log}$-connection@--- radius of}
 of $\nabla$ at $\sigma_i^U$.
 \ede

\begin{rema} \label{QQ1017}
It follows from  (\ref{QR0111}) that
the radius of an $\hslash$-$S^\mr{log}$-connection 
does not change under carrying out a  gauge transformation.
That is to say, 
  $(\mcE, \nabla) \cong (\mcE', \nabla')$ implies $\rho_i^\nabla = \rho_i^{\nabla'}$ for every $i =1, \cdots, r$.
\end{rema}

\begin{exa}[Change of parameter]\label{QQ1040}
Let $(\mcE, \nabla)$ be an $\hslash$-flat $\mbG$-bundle on $U^\mr{log}/S^\mr{log}$ and let $c \in \Gamma (S, \mcO_S) \left(= \mbA^1 (S) \right)$.
Then,  we obtain a $(c \cdot \hslash)$-flat bundle $(\mcE, c \cdot \nabla)$ (cf. (\ref{QQ1032})).
By   the definition of  the action $\star$ on $\mfc$,  the equality
\begin{align} \label{QQ1041}
\rho_i^{c \cdot \nabla} = c \star \rho_i^\nabla
\end{align}
holds for every $i =1, \cdots, r$.
\end{exa}

 \subsection{Radii of a $(\mfg, \hslash)$-oper} \label{QR0200}
Let us fix $\rho := (\rho_i)_{i=1}^r \in \prod_{i=1}^r \mfc (S_i)$, where 
$\prod_{i=1}^r \mfc (S_i) := \{ \emptyset \}$ and $\rho := \emptyset$ when  $r =0$. 
 
\bde \label{y057}
Let $\msE^\spadesuit := (\mcE_\mbB, \nabla)$ be a $(\mfg, \hslash)$-oper on $\msU$.
 Then, we shall say that $\msE^\spadesuit$ is {\bf of radii $\rho$}
 \index{$\msE^\spadesuit$, $(\mfg, \hslash)$-oper@--- of radii $\rho$}
  if  $\rho^{\nabla}_i =\rho_i$ for every $i = 1, \cdots, r$.
When  $r =0$,
 we refer to any $(\mfg, \hslash)$-oper   as being  {\bf  of radii $\emptyset$.}
 \ede

\begin{rema} \label{QQ1001}
Let us briefly mention the relationship between $\mfg^\odot \left(=\mfs \mfl_2\right)$-opers (of prescribed radii) and {\it torally indigenous bundles} (of prescribed radii). \index{torally indigenous bundle} We shall refer the reader to  ~\cite{Mzk1} and ~\cite{Mzk2} for the definition of a  torally indigenous bundle.

Let  $\msX := (X/S, \{ \sigma_i \}_{i=1}^r)$  be an $r$-pointed stable curve 
and take a collection
 $\msE^\spadesuit := (\mcE_{\mbG^\odot}, \nabla, \mcE_{\mbB^\odot})$ defining a   $\mfg^\odot$-oper on $\msX$ (cf. Example \ref{QR0203}).
If  $\mcP$ denotes the $\mbP^1$-bundle on $X$ corresponding to the $\mbG^\odot$-bundle  $\mcE_{\mbG^\odot}$, then 
the $\mbB^\odot$-reduction $\mcE_{\mbB^\odot}$  determines a global  section $X \migi \mcP$ of the natural projection $\mcP \migi X$.
Indeed, recall that 
 the projective line  $\mbP^1$ may be identified  with the flag variety $\mbG^\odot /\mbB^\odot$ classifying Borel subgroups in $\mbG^\odot$.
Hence, the $\mbB^\odot$-reduction $\mcE_{\mbB^\odot}$ (which is, in some sense, viewed as  a family of Borel subgroups in  $\mbG^\odot$ twisted by $\mcE_{\mbG^\odot}$) specifies a section 
\begin{align}
\sigma  : X \migi  \left(\mcE_{\mbG^\odot} \times^{\mbG^\odot}(\mbG^\odot /\mbB^\odot) \cong \mcE_{\mbG^\odot} \times^{\mbG^\odot} \mbP^1 =\right) \mcP.
\end{align}
On the other hand, $\nabla$ induces a logarithmic connection $\nabla_\mcP$ on $\mcP$ (in the sense of Grothendieck formulated  in terms of descent data from 
the PD-envelope  of 
  the diagonal in $X \times_S X$).
 The bijectivity of   (\ref{QR0205}) may be interpreted as  the condition that the 
 Kodaira-Spencer morphism $\mcT_{X^\mr{log}/S^\mr{log}} \migi \sigma^*(\mcT_{\mcP/X})$ at $\sigma$ relative to $\nabla_\mcP$  (cf. ~\cite{Mzk1}, Chap.\,I, \S\,2) is an isomorphism.
Thus, the pair $(\mcP, \nabla_\mcP)$ forms a torally indigenous bundle on $X^\mr{log}/S^\mr{log}$  (cf. ~\cite{Mzk2}, Chap.\,I, \S\,4, Definition 4.1).
One verifies that {\it the assignment $\msE^\spadesuit \mapsto (\mcP, \nabla_\mcP)$ defines a bijective correspondence between the set of isomorphism classes of $\mfg^\odot$-opers on $\msX$ and the set of isomorphism classes of torally indigenous bundles on $X^\mr{log}/S^\mr{log}$ in the sense of Mochizuki}.

Next, let us consider  the composite
\begin{align}
\xi : \Gamma (S, \mcO_S) \migi \mfg^\odot (S) \xrightarrow{\chi} \mfc_{\mfg^\odot} (S),
\end{align}
where the first arrow is given by  
$a \mapsto \begin{pmatrix} a & 0 \\ 0 & -a \end{pmatrix}$.
For each $i = 1, \cdots, r$, we fix an element  $\lambda_i$ of $\{ 0 \} \cup \Gamma (S, \mcO^\times_S) \subseteq \Gamma (S, \mcO_S)$.
Then, 
one verifies that   {\it  
a $\mfg^\odot$-oper $\msE^\spadesuit$  is of radii $(\xi (\lambda_i) )_{i=1}^r$ if and only if  the corresponding  torally indigenous bundle $(\mcP, \nabla_\mcP)$    is of radii $(\lambda_i)_{i=1}^r$ in the sense of  ~\cite{Mzk2}, Chap.\,I, \S\,4, Definition 4.1.}
\end{rema}

\bpr \label{QQ1020}
Let 
$\msE^\spadesuit := ({^\dagger}\mcE_\mbB, \nabla)$ be a $\, \qq$-normal $(\mfg, \hslash)$-oper on $\msU$   and $R \in \Gamma (U, \DV_{\mfg, \msU})$.
Suppose that $\msE^\spadesuit$ is of radii $\rho : = (\rho_i)_{i=1}^r \in \prod_{i=1}^r \mfc (S_i)$.
Then, the $(\mfg, \hslash)$-oper $\msE_{+R}^\spadesuit := ({^\dagger}\mcE_\mbB, \nabla + R)$ (cf. (\ref{plus1})) is of  radii $\rho$ if and only if 
$R$ lies in $\Gamma (U, \DV_{\mfg, \msU}(-D_{\msU}))$.
\epr
\begin{proof}
Let us fix $i \in \{1, \cdots, r \}$.
After possibly replacing $U$ with its  open subscheme,
we may assume that 
there exists a global log chart
 $(U, \partial)$ on $U^\mr{log}/S^\mr{log}$ and $D_j^U = \emptyset$ for every $j \neq i$.
The differential of  $\mr{triv}_{\mbG, (U, \partial)}: U \times_k \mbG \isom \DE_\mbG$ (cf. (\ref{canonicalTg})) determines  an isomorphism
$\mcO_U \otimes_k\mfg \isom \mfg_{\DE_\mbG}$;
it induces   an  isomorphism
\begin{align}\label{QR0257}
\tau_i : \mfg (S_i) \left(=  \Gamma (S_i, \mcO_{S_i} \otimes_k \mfg)\right) \isom
\Gamma (S_i, \sigma^{U*}_i(\mfg_{\DE_\mbG})). 
\end{align}
We shall regard $\Gamma (S_i, \sigma^{U*}_i (\DV_{\mfg, \msU}))$ as a submodule of $\Gamma (S_i, \sigma^{U*}_i (\mfg_{\DE_\mbG}))$ via the following composite injection:
\begin{align}
\Gamma (S_i, \sigma^{U*}_i (\DV_{\mfg, \msU})) &\migiincl \Gamma (S_i, \sigma^{U*}_i (\Omega_{U^\mr{log}/S^\mr{log}}\otimes \mfg_{\DE_\mbG})) \\
&\isom \Gamma (S_i, \sigma^{U*}_i (\Omega_{U^\mr{log}/S^\mr{log}})\otimes \sigma^{U*}_i (\mfg_{\DE_\mbG})) \notag \\
& \isom \Gamma (S_i, \sigma^{U*}_i (\mfg_{\DE_\mbG})), \notag
\end{align}
where  the first arrow arises from $\varsigma$ (cf. (\ref{inclV})) and the last arrow follows from  the residue map $\mr{Res}_{\msU, i}$ (cf. (\ref{QQ040})).
Then, $\tau_i$ restricts to an isomorphism
$\mfg^{\mr{ad}(q_1)} (S_i) \isom \Gamma (S_i, \sigma^{U*}_i (\DV_{\mfg, \msU}))$.
Since $\msE^\spadesuit$ is $\ \qq$-normal,   $\tau_i^{-1}(\mu_i^\nabla)$ can be expressed as $1 \otimes q_{-1} + v$ for some $v \in \mfg^{\mr{ad}(q_1)} (S_i)$.
Hence,
for each $u \in \Gamma (S_i, \sigma^{U*}_i (\DV_{\mfg, \msU}))$,  the following equalities  hold  in $\mfc (S_i)$:
\begin{align}
[\chi] (\sigma^{U*}_i(\DE_\mbG), \mu^\nabla_i + u) = \chi (\tau_i^{-1} (\mu^\nabla_i + u)) = \chi (1 \otimes q_{-1} + (v + \tau_i^{-1}(u))).
\end{align}
This implies (from the bijectivity of $\mr{Kos}_\mfg$) that
the map $\Gamma (S_i, \sigma^{U*}_i (\DV_{\mfg, \msU}))  \migi \mfc (S_i)$ given by $u \mapsto [\chi] (\sigma^{U*}_i(\DE_\mbG), \mu_i^\nabla + u)$ is bijective.
Thus, the equality $\rho_i^{\nabla + R} = \rho_i^\nabla$ holds if and only if the equality  $R = 0$ holds  in $\Gamma (S_i, \sigma^*_i (\DV_{\mfg, \msU}))$, or equivalently, $R \in \Gamma (U, \DV_{\mfg, \msU}(-D_\msU))$.
This completes the proof of the assertion.
\end{proof}

The following assertion is a direct consequence of the above proposition.

\bco \label{QR0901}
Let $\rho := (\rho_i)_{i=1}^r \in \prod_{i=1}^r  \mfc (S_i)$, and denote by
$O_\rho$
 the subset of $\mcO p_{\mfg, \hslash, U^\mr{log}/S^\mr{log}}(U)$ consisting of
$(\mfg, \hslash)$-opers  of radii $\rho$.
Also, suppose that $O_\rho \neq \emptyset$.
Then, the $\Gamma (U, \DV_{\mfg, \msU})$-torsor  structure on $\mcO p_{\mfg, \hslash, U^\mr{log}/S^\mr{log}}(U)$ defined by  (\ref{QR043}) induces  a $\Gamma (U, \DV_{\mfg, \msU}(-D_\msU))$-torsor  structure  on $O_\rho$.
\eco

\subsection{The radius map $\mr{Rad}_\hslash$} \label{QQ1008}

Let $(g,r)$ be a pair of nonnegative integers with $2g-2+r>0$, $S$ a $k$-scheme, and  $\msX:= (f: X \migi S, \{ \sigma_i \}_{i=1}^r)$  an $r$-pointed stable curve of genus $g$ over $S$.
Also,  write $\mfc_S := S \times_k \mfc$.
For each $i \in \{1, \cdots, r \}$,
we shall set
\begin{align} \label{QQ1027}
\mr{Rad}_{\hslash, i} : \mfO \mfp_{\mfg, \hslash, \msX} \migi \mfc_S
\index{$\mr{Rad}_{\hslash, i}$, $\mr{Rad}_\hslash$, radius map}
\end{align}
to be the $S$-morphism given by assigning $\msE^\spadesuit :=(\mcE_\mbB, \nabla) \mapsto \rho_i^\nabla$.
Hence, we obtain an $S$-morphism
\begin{align} \label{QQ1026}
\mr{Rad}_\hslash := \prod_{i=1}^r \mr{Rad}_{\hslash, i} : \mfO \mfp_{\mfg, \hslash, \msX} \migi \mfc^{\times r}_S,
\end{align} 
where $\mfc^{\times r}_S$ denotes the product of $r$ copies of $\mfc_S$ over $S$.
Given  $\rho := (\rho_i)_{i=1}^r \in \mfc^{\times r}(S) \left(=\mfc_S^{\times r} (S) \right)$,
 we shall write
 \index{$\mfO \mfp_{\mfg, \hslash, (\rho,) \msX}$, $\mfO \mfp_{\mfg, (\hslash,) (\rho,) g,r}$}
\begin{align} \label{QQ1025}
\mfO \mfp_{\mfg, \hslash, \rho, \msX} := \mr{Rad}_{\hslash}^{-1}(\rho), 
\end{align}
i.e., the closed subscheme of $\mfO \mfp_{\mfg, \hslash\msX}$ defined as 
the scheme-theoretic  inverse image of $\rho$ via $\mr{Rad}_\hslash$.
That is to say, $\mfO \mfp_{\mfg, \hslash, \rho, \msX}$ is the $S$-scheme representing 
the $\mfS \mfe \mft$-valued contravariant  functor on $\mfS \mfc \mfh_{/S}$ which, to any $S$-scheme $s: S' \migi S$, assigns the set of isomorphism classes of $(\mfg, \hslash)$-opers on $S' \times_S \msX$ of radii $\rho \circ s := (\rho_i \circ s)_{i=1}^r \in \mfc^{\times r} (S')$. 
In particular, if $r=0$, then we have 
 $\mfO \mfp_{\mfg, \hslash, \emptyset, \msX} = \mfO \mfp_{\mfg, \hslash, \msX}$.

\section{The moduli space of $(\mfg, \hslash)$-opers of prescribed radii} \label{y059}

\subsection{The action $\star$ on the set of  radii} \label{QR0800}
Let us write $\mbA^1 = \mr{Spec}(k [\hat{\hslash}])$, \index{$\mbA^1$, affine line} i.e.,  we use the symbol $\hat{\hslash}$ to denote the  parameter of the affine line $\mbA^1$. 
Also, write $\mbA_S^1 := \mbA^1 \times_k S$.
For $\rho : = (\rho_i)_{i=1}^r \in \mfc^{\times r} (S)$, we shall set
\begin{equation}
\hat{\rho} := ( \hat{\rho}_i)_{i =1}^r \in \mfc^{\times r} (\mbA^1_S),
\end{equation}
 where $\hat{\rho}_i$ ($i = 1, \cdots, r$)  denotes  the composite 
 \begin{equation}
 \hat{\rho}_i : \mbA^1_S \left(= \mbA^1 \times_k S\right) \xrightarrow{\mr{id}_{\mbA^1} \times \rho_i} \mbA^1 \times_k \mfc \xrightarrow{\star} \mfc.
 \end{equation}
(We set $\hat{\rho} := \emptyset$ if $r =0$ and $\rho = \emptyset$.)
Then, the   equality $ \hat{\rho} \circ (\hslash, \mr{id}_S) = \hslash \star \rho$ holds
  for any $\hslash \in \Gamma (S, \mcO_S)\left(=\mbA^1 (S) \right)$.

\subsection{$(\mfg, \hat{\hslash})$-opers with  formal parameter $\hat{\hslash}$} \label{QR0801}
Write  $\hat{\msX} := \mbA^1_S \times_S \msX$, i.e.,
the base-change  of $\msX$ to  $\mbA^1_S$.
Let us consider 
the $\mbA^1_S$-scheme
\begin{align}
\mfO \mfp_{\mfg, \hat{\hslash},  \hat{\msX}} \ \left(\text{resp.,} \ \mfO \mfp_{\mfg, \hat{\hslash}, \hat{\rho}, \hat{\msX}}\right).
\end{align}
 If we regard it as an $S$-scheme, then it represents the $\mfS \mfe \mft$-valued contravariant  functor on $\mfS \mfc \mfh_{/S}$ assigning, to any 
 $S$-scheme $s : S' \migi S$,  the set of pairs $(\hslash, [\msE^\spadesuit])$ consisting of  $\hslash \in \Gamma (S', \mcO_{S'}) = \mbA^1 (S')$
 and the isomorphism class $[\msE^\spadesuit]$ of  a $(\mfg, \hslash)$-oper $\msE^\spadesuit$ on $S' \times_S \msX$  (resp., a $(\mfg, \hslash)$-oper $\msE^\spadesuit$ on $S' \times_S \msX$ of radii $\hslash \star \rho \in \mfc^{\times r}(S')$).
Denote by 
\begin{equation}   \label{QQ1056}
\Spec : \mfO \mfp_{\mfg, \hat{\hslash}, \hat{\msX}} \migi \mbA^1_S
\ \left(\text{resp.,} \ \Spec_\rho : \mfO \mfp_{\mfg, \hat{\hslash}, \hat{\rho}, \hat{\msX}} \migi \mbA^1_S \right) 
\index{$\Spec$, $\Spec_\rho$}
\end{equation}
the  $S$-morphism 
given by assigning $(\hslash, [\msE^\spadesuit]) \mapsto \hslash$.
By applying  Proposition \ref{y046}, (ii), and Theorem \ref{y051}, (i),  to $\hat{\msX}$, 
we see that the morphism $\Spec$
  forms a relative affine space over $\mbA^1_S$ modeled on
$\mbA_S^1 \times_S \mbV (f_* (\DV_{\mfg, \msX})) \left(\cong \mbV ((\mr{id}_{\mbA^1}\times f)_*(\DV_{\mfg, \hat{\msX}})) \right)$.
For each $\hslash \in \Gamma (S, \mcO_S) \left(= \mbA^1_S (S) \right)$, 
 the natural $S$-morphism
\begin{equation} \label{QQ1066}
\mfO \mfp_{\mfg, \hslash, \msX} \migi 
S \times_{\hslash,  \mbA^1_S, \,  \scalebox{0.7}{$\Spec$}} \mfO \mfp_{\mfg, \hat{\hslash}, \hat{\msX}}   
\left(\text{resp.,} \ \mfO \mfp_{\mfg, \hslash, \hslash \star \rho, \msX} \migi 
S \times_{\hslash,  \mbA^1_S, \,  \scalebox{0.7}{$\Spec_\rho$}} \mfO \mfp_{\mfg, \hat{\hslash}, \hat{\rho}, \hat{\msX}}\right)
\end{equation}
is an isomorphism.
Moreover, $\mfO \mfp_{\mfg, \hat{\hslash},  \hat{\msX}}$ has a $\mbG_m$($= \mr{Spec}(k [\hat{\hslash}^{\pm 1}])$)-action 
\begin{equation} \label{action61}
\mbG_m \times_k  \mfO \mfp_{\mfg, \hat{\hslash}, \hat{\msX}}  \migi \mfO \mfp_{\mfg, \hat{\hslash}, \hat{\msX}}
\end{equation}
 determined by  $(c, (\hslash, [\msE^\spadesuit])) \mapsto (c \cdot \hslash, [c \cdot \msE^\spadesuit])$ (cf. (\ref{QQ1078}) for the definition of $c \cdot \msE^\spadesuit$), which makes the following square diagram commute:
 \begin{align} \label{QR0455}
\vcenter{\xymatrix@C=76pt@R=36pt{
\mbG_m \times_k  \mfO \mfp_{\mfg, \hat{\hslash}, \hat{\msX}} \ar[d]_-{\mr{id}_{\mbG_m} \times (\scalebox{0.7}{$\Spec$},   \mr{Rad}_{\hat{\hslash}})} \ar[r]^-{(\ref{action61})}&\mfO \mfp_{\mfg, \hat{\hslash}, \hat{\msX}} \ar[d]^-{(\scalebox{0.7}{$\Spec$}, \mr{Rad}_{\hat{\hslash}})}
\\
\mbG_m \times_k \mbA^1_S \times_S \mfc_S^{\times r} \ar[r]_-{(c, \hslash, \varrho) \mapsto (c\cdot \hslash, c \star \varrho)}&\mbA^1_S \times_S \mfc_S^{\times r}.
}}
\end{align}
The commutativity of this diagram implies that (\ref{action61}) restricts to a $\mbG_m$-action
\begin{equation} \label{QR0506}
\mbG_m \times_k  \mfO \mfp_{\mfg, \hat{\hslash}, \hat{\rho}, \hat{\msX}}  \migi \mfO \mfp_{\mfg, \hat{\hslash}, \hat{\rho}, \hat{\msX}}.
\end{equation}
on $\mfO \mfp_{\mfg, \hat{\hslash}, \hat{\rho}, \hat{\msX}}$.
If $\hslash \in \Gamma (S, \mcO^\times_S) \left(= \mbG_m (S) \right)$, then the composite
\begin{align} \label{QR0589}
\mbG_m \times_k  \mfO \mfp_{\mfg, \hslash, \msX} 
  \xrightarrow{\mr{inclusion}} \mbG_m \times_k 
 \Spec^{-1}((\mbG_m)_S)
    \xrightarrow{(\ref{action61})} 
  \Spec^{-1}((\mbG_m)_S),
\end{align}
where $(\mbG_m)_S :=  \mbG_m \times_k S \left(\subseteq \mbA^1_S \right)$, is an isomorphism.
Moreover,
the following diagram is commutative:
\begin{align} \label{QR0600}
\vcenter{\xymatrix@C=56pt@R=36pt{
\mbG_m \times_k  \mfO \mfp_{\mfg, \hslash, \msX} \left(= (\mbG_m)_S \times_S \mfO \mfp_{\mfg, \hslash, \msX}  \right) \ar[r]^-{(\ref{QR0589})} \ar[d]_-{\mr{id}_{\mbG_m} \times \mr{Rad}_\hslash} & \Spec^{-1}((\mbG_m)_S) \ar[d]^-{(\scalebox{0.7}{$\Spec$},  \mr{Rad}_{\hat{\hslash}})} \\
(\mbG_m)_S \times_S  \mfc_S^{\times r}\ar[r]_-{(c, \varrho) \mapsto (c \cdot \hslash, c \star \varrho)} & (\mbG_m)_S \times_S \mfc_S^{\times r}.
}}
\end{align}
One verifies that (\ref{QR0589}) 
 restricts to an isomorphism
\begin{align} \label{QR0588}
\mbG_m \times_k  \mfO \mfp_{\mfg, \hslash, \hslash \star \rho,  \msX}  \isom
\Spec^{-1}_\rho ((\mbG_m)_S).
\end{align}
In particular, the assignment $\msE^\spadesuit \mapsto \hslash \cdot \msE^\spadesuit$ determines an isomorphism
\begin{align} \label{QS1002}
\mfO \mfp_{\mfg, 1, \msX} \isom  \mfO \mfp_{\mfg, \hslash, \msX} \ \left( \text{resp.,} \ \mfO \mfp_{\mfg, 1, \rho, \msX} \isom  \mfO \mfp_{\mfg, \hslash, \hslash \star \rho,  \msX} \right).
\end{align}

\subsection{The representability of $\mfO \mfp_{\mfg, \hslash, \hslash \star \rho, \msX}$} \label{QS1005}

In the following, we prove the representability of $\mfO \mfp_{\mfg, \hslash, \hslash \star\rho, \msX}$ by observing the radius map $\mr{Rad}_\hslash$ specialized at $\hslash =0$.

\bt\label{QQ1060}
(Recall that   either ``\,$\mr{char}(k)=0$'' or ``\,$\mr{char}(k) >2  h_\mbG$'' is fulfilled.) 
The morphism 
\begin{align} \label{QR0624}
(\Spec, \mr{Rad}_{\hat{\hslash}}) : \mfO \mfp_{\mfg, \hat{\hslash}, \hat{\msX}} \migi \mbA^1_S \times_S \mfc_S^{\times r}
\end{align}
forms a relative affine space over $\mbA^1_S \times_S \mfc_S^{\times r}$ modeled on $(\mbA^1_S \times_S \mfc_S^{\times r}) \times_S \mbV (f_*(\DV_{\mfg, \msX}))$.
In particular, for  any  $\hslash \in \Gamma (S, \mcO_S)$ and $\rho \in \mfc^{\times r}(S)$, 
 the $S$-scheme $\mfO \mfp_{\mfg, \hslash, \hslash \star \rho, \msX}$
 $\left(=(\Spec, \mr{Rad}_{\hat{\hslash}})^{-1}((\hslash, \rho))  \right)$ forms a relative affine space   modeled on $\mbV (f_*(\DV_{\mfg, \msX}(-D_\msX)))$.
 \et
\begin{proof}
By Corollary \ref{QR0901}, the problem is reduced to proving the claim that
$(\Spec, \mr{Rad}_{\hat{\hslash}})$ is smooth and surjective.
In order to  prove this, we may assume, after possibly changing the base,
that  $S = \mr{Spec}(k)$ and $k$ is algebraically closed.

First, we shall prove the smoothness of $(\Spec, \mr{Rad}_{\hat{\hslash}})$.
According to  Corollary  \ref{QR0901},
the fiber $(\Spec, \mr{Rad}_{\hat{\hslash}})^{-1}((\hslash, \rho)) \left( = \mfO \mfp_{\mfg, \hslash, \hslash \star \rho, \msX}\right)$ of  $(\Spec, \mr{Rad}_{\hat{\hslash}})$ over  each point  $(\hslash, \rho) \in (\mbA^1 \times_k \mfc^{\times r})(k)$ forms, if it is nonempty, an affine space modeled on $\Gamma (X, \DV_{\mfg, \msX}(-D_\msX))$.
In particular,  $(\Spec, \mr{Rad}_{\hat{\hslash}})^{-1}((\hslash, \rho))$ is smooth and has dimension  equal to   
the value
\begin{align}
{^c \aleph} (\mfg) =  (\aleph (\mfg) +1) - (1+r \cdot \mr{rk}(\mfg)) = \mr{dim}(\mfO \mfp_{\mfg, \hat{\hslash}, \hat{\msX}}) - \mr{dim}(\mbA^1 \times_k\mfc^{\times r})
\end{align}
 (cf. Proposition \ref{y046}, (ii)).
Hence, since both $\mfO \mfp_{\mfg, \hat{\hslash}, \hat{\msX}}$ and $\mbA^1 \times_k \mfc^{\times r}$ are smooth over $k$, the smoothness of $(\Spec, \mr{Rad}_{\hat{\hslash}})$   follows from  ~\cite{AlKl},  Chap.\,V, Proposition (3.5) (or ~\cite{Har}, Chap.\,III, Exercise 10.9).

The remaining portion is  the surjectivity of  $(\Spec, \mr{Rad}_{\hat{\hslash}})$.
To begin with,  we shall prove  the  surjectivity  of  
$\mr{Res}_0 : \mfO \mfp_{\mfg, 0, \msX} \migi \mfc^{\times r}$.
Given  a global section $s$ of $\Omega_{X^\mr{log}/k^\mr{log}}^{\otimes l}$  for  $l \in \mbZ_{>0}$, we shall  denote by  $s (\sigma_i)$ 
the image of $\sigma^*_i (s)$   via $\mr{Res}^{\otimes l}_{\msX, i} : \sigma^*_i (\Omega^{\otimes l}_{X^\mr{log}/k^\mr{log}}) \isom k$.
Recall that there exists a composite isomorphism
\begin{align} \label{QR0345}
\mfO \mfp_{\mfg, 0, \msX} \isom  \Gamma (X, \DV_{\mfg}) \isom \bigoplus_{j\in \mbZ_{>0}} \Gamma (X, \DV_{\mfg,  j}) \isom \bigoplus_{j \in \mbZ_{>0}} \Gamma (X, \Omega_{X^\mr{log}/k^\mr{log}}^{\otimes (j+1)})\otimes_k \mfg_j^{\mr{ad}(q_1)}, \hspace{-3mm}
\end{align}
where  the first arrow denotes the inverse to (\ref{QR0304}) and the third arrow is obtained by (\ref{QR020}).
Under the identifications given by this composite  and $\mr{Kos}_\mfg$ respectively,
the morphism $\mr{Rad}_0$ may be described  as a $k$-linear morphism
\begin{align}  \label{QR0388}
\bigoplus_{j \in \mbZ_{>0}} \Gamma (X, \Omega_{X^\mr{log}/k^\mr{log}}^{\otimes (j+1)})\otimes_k \mfg_j^{\mr{ad}(q_1)} \migi \bigoplus_{j\in \mbZ_{>0}}\mfg^{\mr{ad}(q_{1})}_j
\end{align}
given by $\sum_{j} s \otimes v \mapsto \sum_j s(\sigma_i) \cdot v$ for $s \in \Gamma (X, \Omega_{X^\mr{log}/k^\mr{log}}^{\otimes (j+1)})$ and $v \in \mfg_j^{\mr{ad}(q_1)}$.
Since the degree   of $\Omega_{X^\mr{log}/k^\mr{log}}^{\otimes (j+1)}$ for every $j \in \mbZ_{>0}$ is greater than $2g-2+r$,
   (\ref{QR0388}) turns out to be surjective (cf. ~\cite{Har}, Chap.\,IV, \S\,4.3, Corollary 3.2, (a)).
   That is to say, we obtain  the surjectivity of $\mr{Rad}_0$.

We turn to  finish the proof of the surjectivity.
Suppose, on the contrary,  that $(\Spec, \mr{Rad}_{\hat{\hslash}})$ is not surjective, i.e.,
there exists a $k$-rational point $(\hslash_0, \rho_0)$ of  $\mbA^1 \times_k \mfc^{\times r}$ outside  the set-theoretic image  $\mr{Im}((\Spec,\mr{Rad}_{\hat{\hslash}}))$  of $(\Spec, \mr{Rad}_{\hat{\hslash}})$.
By the surjectivity of $\mr{Rad}_0$ proved above, $\hslash_0$ is contained in $\mbG_m \left(\subseteq \mbA^1 \right)$.
Denote by $\rho_0^{\mbG_m}$ the image of $\mbG_m \times \{ \rho_0 \} \left(\subseteq \mbG_m \times \mfc^{\times r} \right)$ via the lower horizontal arrow in (\ref{QR0600}).
The set $\overline{\rho}_0^{\mbG_m} := \left\{ (c \cdot \hslash_0, c \star \rho_0) \, | \, c \in \mbA^1\right\}$ forms a closed subset of $\mbA^1 \times_k \mfc^{\times r}$ containing $\rho_0^{\mbG_m}$.
The commutativity of  (\ref{QR0600}) implies that 
 $\rho_0^{\mbG_m}$ is contained in $(\mbA^1 \times_k \mfc^{\times r}) \setminus \mr{Im} ((\Spec,\mr{Rad}_{\hat{\hslash}}))$.
Since  $(\Spec,\mr{Rad}_{\hat{\hslash}})$ is flat as proved already, 
$\mr{Im} ((\Spec,\mr{Rad}_{\hat{\hslash}}))$ is open in  $\mbA^1 \times_k \mfc^{\times r}$.
It follows that  $\overline{\rho}_0^{\mbG_m} \subseteq (\mbA^1 \times \mfc^{\times r}) \setminus \mr{Im} ((\Spec,\mr{Rad}_{\hat{\hslash}}))$.
But,  $(0, 0 \star \rho_0)$ belongs to $\overline{\rho}_0^{\mbG_m}$,
 so this contradicts the surjectivity  of  $\mr{Rad}_0$.
Consequently,  $(\Spec,\mr{Rad}_{\hat{\hslash}})$ must be  surjective.
 This completes the proof of this proposition.
\end{proof}


\begin{rema} \label{QQ1003}
The injection
$\iota_\mfg$ restricts to an injection $(\mfg^\odot)^{\mr{ad}(q_1^\odot)} \migiincl \mfg^{\mr{ad}(q_1)}$.
Under the identifications $(\mfg^\odot)^{\mr{ad}(q_1^\odot)} = \mfc_{\mfg^\odot}$, $\mfg^{\mr{ad}(q_1)} = \mfc_\mfg$ given by $\mr{Kos}_{\mfg^\odot}$, $\mr{Kos}_\mfg$ respectively,
this injection becomes an injection
$\iota_\mfc : \mfc_{\mfg^\odot} \migiincl \mfc_\mfg$.
\index{$\iota_\mfg$, $\iota_\mbB$, $\iota_\mfc$}
For $\rho^\odot := (\rho^\odot_i )_{i=1}^r \in \mfc_{\mfg^\odot}^{\times r} (S)$, we shall write 
\begin{align} \label{QS18}
\iota_{\mfc*}(\rho^\odot) := (\iota_\mfc \circ \rho^\odot_i)_{i=1}^r \in \mfc^{\times r}_\mfg (S).
\end{align}
If $\msE_\odot^\spadesuit$ is a $(\mfg^\odot, \hslash)$-oper of radii $\rho^\odot$, then  the induced  $(\mfg, \hslash)$-oper  
$\iota_{\mbG*}(\msE^\spadesuit_\odot)$ (cf. (\ref{associated}))
is of radii  $\iota_{\mfc*}(\rho^\odot)$.
Hence,  $\iota_\mfg^{\mfO \mfp}$ (cf. (\ref{QR066})) restricts to
 a closed immersion
\begin{align} \label{QS19}
\iota_{\mfg, \rho^\odot}^{\mfO \mfp} : \mfO \mfp_{\mfg^\odot, \hslash, \rho^\odot, \msX} \migi \mfO \mfp_{\mfg, \hslash, \iota_{\mfc*}(\rho^\odot), \msX}.
\end{align}
Also, $\mfO \mfp_{\mfg^\odot, \hslash, \rho^\odot, \msX}$ forms a relative affine space over $S$ modeled on $\mbV (f_*(\Omega_{X^\mr{log}/S^\mr{log}}^{\otimes 2} (-D_\msX)))$ $\left(= \mbV (f_*(\DV_{\mfg^\odot, \msX}(-D_\msX))) \right)$ and $\iota_{\mfg, \rho^\odot}^{\mfO \mfp}$
 induces  an $S$-isomorphism
\begin{align} \label{QS17}
 \mfO \mfp_{\mfg^\odot, \hslash, \rho^\odot, \msX} \times^{\mbV (f_*(\Omega_{X^\mr{log}/S^\mr{log}}^{\otimes 2} (-D_\msX)))}  \mbV (f_*(\DV_{\mfg, \msX}(-D_\msX))) \isom \mfO \mfp_{\mfg, \hslash, \iota_{\mfc*}(\rho^\odot), \msX}.
\end{align}
If, moreover, $\mfg$ is simple, then (\ref{QS17}) becomes an $S$-isomorphism
\begin{align} \label{QHH4}
 \mfO \mfp_{\mfg^\odot, \hslash, \rho^\odot, \msX} \times_S
 \mbV (f_*(\DV^2_{\mfg}(-D_\msX))) \isom \mfO \mfp_{\mfg, \hslash, \iota_{\mfc*}(\rho^\odot), \msX}
\end{align}
(cf. (\ref{gissimple})).
\end{rema}

\begin{rema} \label{QQ1004}
Let us keep the notation in  Remark \ref{QQ892}.
There exists a canonical identification  $\mfc_{\mfg \oplus \mfg'} = \mfc_\mfg \times_k \mfc_{\mfg'}$.
Given  $\rho := (\rho_i)_{i=1}^r \in \mfc_\mfg^{\times r} (S)$ and  
$\rho' := (\rho'_i)_{i=1}^r \in \mfc_{\mfg'}^{\times r} (S)$,
we obtain 
\begin{align}
(\rho, \rho') := ((\rho_i, \rho'_i))_{i=1}^r \in  \mfc_{\mfg \oplus \mfg'}^{\times r} (S).
\end{align}
If $\msE^\spadesuit$ is a $(\mfg, \hslash)$-oper on $\msX$ of radii $\rho$ and 
$\msE'^\spadesuit$ is  a $(\mfg', \hslash)$-oper on $\msX$ of radii $\rho'$,
then the $(\mfg \oplus \mfg', \hslash)$-oper 
$\msE^\spadesuit \boxtimes \msE'^\spadesuit$ (cf. (\ref{QS14}))
 is of  radii $(\rho, \rho')$.
Hence,   (\ref{product}) restricts to  an $S$-isomorphism
\begin{equation} \label{productrad}
 \mfO \mfp_{\mfg \oplus \mfg', \hslash, (\rho, \rho'), \msX}  \isom \mfO \mfp_{\mfg, \hslash, \rho, \msX} \times_S \mfO \mfp_{\mfg', \hslash, \rho', \msX}.
\end{equation}

\end{rema}

\section{The universal moduli stack} \label{QR70}

In this final section of Chapter 2, we deal with   the  universal moduli stack of $(\mfg, \hslash)$-opers, i.e., $\mfO \mfp_{\mfg, \hslash, (\rho,) \msX}$ in the case of $\msX = \msC_{g,r}$. 
Let $(g,r)$ be a pair of nonnegative integers with $2g-2+r>0$, $\hslash$ an element of $k$, and $\rho \in \mfc^{\times r}(k)$ (where $\rho := \emptyset$ if $r =0$).
Denote by
 \index{$\mfO \mfp_{\mfg, \hslash, (\rho,) \msX}$, $\mfO \mfp_{\mfg, (\hslash,) (\rho,) g,r}$}
\begin{align} \label{univ1}
 \mfO \mfp_{\mfg, \hslash, g,r}  \left(\text{resp.,} \ \mfO \mfp_{\mfg, \hslash, \rho,  g,r} \right)
\end{align}
the category over $\mfS \mfc \mfh_{/k}$ defined as follows:
\begin{itemize}
\item
The objects are the pairs  $(\msX, \msE^\spadesuit)$, where $\msX$ is an $r$-pointed stable curve of genus $g$ over a $k$-scheme  $S$ and $\msE^\spadesuit$ is a $(\mfg, \hslash)$-oper on $\msX$ (resp.,  a $(\mfg, \hslash)$-oper on $\msX$ of radii $\rho$);
\item
The morphisms from $(\msX, \msE^\spadesuit)$ to $(\msX', \msE'^\spadesuit)$ are 
 morphisms $(\phi, \Phi) : \msX \migi \msX'$ of $r$-pointed stable curves  with $\msE^\spadesuit \cong \phi^*(\msE'^\spadesuit)$.
\item
The projection $\mfO \mfp_{\mfg, \hslash, g,r} \migi \mfS \mfc \mfh_{/k}$ (resp., $\mfO \mfp_{\mfg, \hslash, \rho,  g,r} \migi \mfS \mfc \mfh_{/k}$) is given by assigning, to each pair $(\msX, \msE^\spadesuit)$ as above, the base $k$-scheme $S$ of $\msX$. 
\end{itemize}
In particular,  we have $\mfO \mfp_{\mfg, \hslash, \emptyset,  g, 0} = \mfO \mfp_{\mfg, \hslash, g, 0}$.
The assignment $(\msX, \msE^\spadesuit) \mapsto \msX$ defines  a functor  $\mfO \mfp_{\mfg, \hslash, g,r} \migi \overline{\mfM}_{g,r}$ (resp., $\mfO \mfp_{\mfg, \hslash, \rho,  g,r}  \migi \overline{\mfM}_{g,r}$).
If we take a $k$-scheme $S$ and a $k$-morphism  $s : S \migi \overline{\mfM}_{g,r}$, which classifies  an $r$-pointed stable curve  $\msX$ over $S$, then
there exists a canonical $S$-isomorphism 
\begin{align}
S \times_{s, \overline{\mfM}_{g,r}} \mfO \mfp_{\mfg, \hslash, (\rho,) g,r}  \isom \mfO \mfp_{\mfg, \hslash, (\rho,) \msX}.
 \end{align}
For simplicity, we write
 \index{$\mfO \mfp_{\mfg, \hslash, (\rho,) \msX}$, $\mfO \mfp_{\mfg, (\hslash,) (\rho,) g,r}$}
\begin{align} \label{univ1gg}
 \mfO \mfp_{\mfg, g,r} :=\mfO \mfp_{\mfg, 1, g,r}   \left(\text{resp.,} \ \mfO \mfp_{\mfg, \rho,  g,r} := \mfO \mfp_{\mfg, \hslash, \rho,  g,r} \right).
\end{align}
By applying  results obtained so far, we obtain the following theorem concerning the geometric structure of this moduli category  of $(\mfg, \hslash)$-opers; it is a generalization of a result for $\mfg = \mfs \mfl_2$ proved by S. Mochizuki (cf. ~\cite{Mzk1}, Chap.\,I, Theorem 2.3).

\bt[cf. Theorem \ref{y014}]  \label{QR90}
Assume  that
either ``\,$\mr{char}(k) =0$'' or ``\,$\mr{char}(k) > 2  h_\mbG$'' is fulfilled.
Let $\hslash \in k$ and $\rho \in \mfc^{\times r}(k)$ (where  $\rho := \emptyset$ if $r=0$).
Then, $\mfO \mfp_{\mfg, \hslash, g,r}$ (resp., $\mfO \mfp_{\mfg, \hslash, \rho, g,r}$) may be represented by a relative affine space over $\overline{\mfM}_{g,r}$ of relative dimension $\aleph (\mfg)$ (resp., $^c \aleph (\mfg)$) (cf. (\ref{aleph})).
In particular, $\mfO \mfp_{\mfg, \hslash, g,r}$ (resp., $\mfO \mfp_{\mfg, \hslash, \rho, g,r}$) is a
nonempty, geometrically connected, and smooth Deligne-Mumford stack over $k$ of dimension $3g-3+r + \aleph (\mfg)$ (resp., $3g-3+r + {^c \aleph} (\mfg)$).  
  \et
\begin{proof}
The assertion follows from Theorems  \ref{y051}, (i), and    \ref{QQ1060} applied to  various $r$-pointed stable curves $\msX$ of genus $g$.
\end{proof}

Finally, 
 note that the  Deligne-Mumford stack $\mfO \mfp_{\mfg, \hslash,  g,r}$ (resp., $\mfO \mfp_{\mfg, \hslash, \rho, g,r}$)   
  admits a  log structure
 pulled-back from the log structure of $\overline{\mfM}_{g,r}^\mr{log}$;
we denote the resulting log stack by
\begin{align}\label{univ1p}
 \mfO \mfp_{\mfg, \hslash,  g,r}^\mr{log} \ \left(\text{resp.,} \  \mfO \mfp_{\mfg, \hslash, \rho, g,r}^\mr{log}\right).
\end{align}
By the above theorem, $\mfO \mfp_{\mfg, \hslash,  g,r}^\mr{log}$  (resp., $\mfO \mfp_{\mfg, \hslash, \rho, g,r}^\mr{log}$) turns out to be  log smooth over $k$.

\chapter{Opers in positive characteristic} \label{y064} 

The aim of this chapter is to study  $(\mfg, \hslash)$-opers in positive characteristic.
To begin with, we introduce the logarithmic version of a restricted Lie algebroid and recall the $p$-curvature associated to a flat principal bundle on a log scheme.
After that,  
we define the notions of a {\it dormant} oper and a {\it $p$-nilpotent} oper, as opers satisfying certain vanishing conditions on $p$-curvature.
The central characters in our discussion are 
 the moduli spaces classifying these objects.
 For example,  the moduli space of $p$-nilpotent opers can be obtained as the kernel of the {\it Hitchin-Mochizuki morphism}.
As  an important consequence of this chapter, we conclude the properness of these moduli spaces by proving the finiteness of the Hitchin-Mochizuki morphism.

\begin{notation}
Let $k$ 
\index{$k$, field} 
be a perfect field of characteristic $p >0$ and 
$\mbG$
\index{$\mbG$, algebraic group}
 a connected smooth affine  algebraic group over $k$ with Lie algebra $\mfg$.
\index{$\mfg$, Lie algebra}
\end{notation}

\section{Frobenius twists and relative Frobenius morphisms}\label{y065}

\subsection{Frobenius morphisms}
Let $T$ be a scheme  over $k$ ($\supseteq \mbF_p =:\mbZ/p  \mbZ$) and $f : Y \migi T$ a scheme over $T$.
Denote by $F_T : T \migi T$ (resp., $F_Y : Y \migi Y$) \index{$F_T$, absolute Frobenius endomorphism} the absolute (i.e., $p$-th power) Frobenius endomorphism of $T$ (resp., $Y$).
The  {\bf Frobenius twist} \index{$Y^{(1)}$, Frobenius twist of $Y$} of $Y$ over $T$ is, by definition,  the base-change 
\begin{align}
Y^{(1)} \left(:=  T \times_{T, F_T}Y \right)
\end{align}
 of $Y$ via
 $F_T : T \migi T$.
Denote by $f^{(1)} : Y^{(1)} \migi T$ the structure morphism of $Y^{(1)}$.
The {\bf relative Frobenius morphism} \index{$F_{Y/T}$, relative Frobenius morphism} of $Y$ over $T$ is, by definition,   the unique morphism $F_{Y/T} : Y \migi Y^{(1)}$  over $T$ that fits into the following commutative diagram: 
\begin{align} \label{QS3360}
\vcenter{\xymatrix@C=36pt@R=36pt{
Y \ar@/^10pt/[rrrrd]^{F_Y}\ar@/_10pt/[ddrr]_{f} \ar[rrd]_{  \ \ \ \  \  \   F_{Y/T}} & & &   &   \\
& & Y^{(1)}  \ar[rr]_{F_T \times \mr{id}_Y} \ar[d]^-{f^{(1)}}  \ar@{}[rrd]|{\Box}  &  &  Y \ar[d]^-{f} \\
&  &T \ar[rr]_{F_T} & &  T.
}}
\end{align}

\subsection{The Frobenius twist of a pointed stable curve} \label{QS87}
Let $r$ be a nonnegative integer, $S$ a scheme over $k$, and 
 $\msU := (U/S, \{ \sigma^U_i : S_i \migi U \}_{i=1}^r)$  a local $r$-pointed stable curve  over  $S$.
For $i = 1, \cdots, r$, we write 
\begin{equation} \label{QS1401}
\sigma_i^{U(1)} := (\mr{id}_S |_{S_i}, \sigma^U_i \circ F_{S^i}) : S_i \migi U^{(1)}  \left(= S \times_{S, F_S}U \right).
\end{equation}
Then, the collection of data 
\index{$\msU$, local $r$-pointed curve@--- $\msU^{(1)}$, Frobenius twist of $\msU$}
\begin{equation} \label{QS99}
\msU^{(1)} := (U^{(1)}/S, \{ \sigma_i^{U(1)} \}_{i=1}^r)
\end{equation}
 specifies a local $r$-pointed stable curve over $S$.
If $\msU$ is taken to be an $r$-pointed stable curve (of genus $g \geq 0$), then the resulting collection $\msU^{(1)}$ becomes an $r$-pointed stable curve (of genus $g$).

\subsection{Relatively torsion-free sheaves} \label{QS88}
Let $S$ be a $k$-scheme and $U$ a prestable curve over $S$.
If $U/S$ is smooth,
then the relative Frobenius morphism $F_{U/S} : U \migi U^{(1)}$ is finite and faithfully flat of degree $p$, so
 the direct image $F_{U/S*}(\mcO_U)$ of $\mcO_U$ is a vector bundle on $U^{(1)}$ of rank $p$.
This is not true in general  because  $F_{U/S} : U \migi U^{(1)}$ is not flat  at a non-smooth point of $U$.
But, as is shown in Proposition \ref{torsionE} below,   $F_{U/S*}(\mcO_U)$ satisfies a certain local property even  at a non-smooth point, i.e., $F_{U/S*}(\mcO_U)$ is  {\it relatively torsion-free} (cf. ~\cite{Ses}, Septi\`{e}me partie, D\'{e}finition 1 or   ~\cite{Ja}, Definition 1.0.1).
Any vector bundle  on $U$ is relatively torsion-free, and vice versa 
on  the smooth locus of $U$ relative to $S$.
We here recall the definition.

\bde \label{QS1406}
Let  $n$ be a positive integer and $\mcF$  a coherent $\mcO_{U}$-module.
Then, $\mcF$ is called
{\bf relatively torsion-free of rank $n$}
\index{relatively torsion-free sheaf}
 if it is  flat over $S$ and, for each geometric  point $s \in S$,
the $\mcO_{U_s}$-module $\mcF |_{U_s}$ on the fiber $U_s := s \times_S U$   is of rank $n$ (cf. ~\cite{HL}, \S\,1.2, Definition 1.2.2) and has no associated primes of height one.

\ede

The following proposition will be used in the discussion of Chapter \ref{y0139}.

\bpr  \label{torsionE} 
 For any vector bundle $\mcV$ on $U$ of rank $n>0$,  the $\mcO_{U^{(1)}}$-module $F_{U/S*}(\mcV)$ is relatively torsion-free of rank $n \cdot p$.
  \epr
\begin{proof}
Since $\mcV$ is locally isomorphic to $\mcO_U^{\oplus n}$ and 
  $F_{U/S*}(\mcV)$ is  flat over $S$, 
we may assume, without loss of generality, that $S = \mr{Spec}(k)$, $\mcV = \mcO_U$, and $k$ is algebraically closed.
 As  explained above, $F_{U/k*}(\mcO_U)$ is locally free of rank $p$ over the smooth locus of $U^{(1)}$.
Hence, it suffices to  prove that $F_{U/k*}(\mcO_U)$ is (relatively) torsion-free of rank $p$ at each  nodal point of $U$.
Let us take 
 a nodal point $q$ of  $U$, and write $q^{(1)} := F_{U/k}(q) \in U^{(1)}(k)$.
Denote by $\widehat{\mcO}_{U^{(1)}, q^{(1)}}$  the completions of the stalk of $\mcO_{U^{(1)}}$ at $q^{(1)}$ with respect to its maximal ideal $\mfm_{U^{(1)}, q^{(1)}}$  and  by $\widehat{\mfm}_{U^{(1)},q^{(1)}}$  the maximal ideal of $\widehat{\mcO}_{U^{(1)}, q^{(1)}}$.
The problem is to prove the torsion-freeness of  
 the $\mfm_{U^{(1)}, q^{(1)}}$-adic  completion $\widehat{\mcO}_{U, \mfq}$  of
 the stalk $\mcO_{U, q}$ of $\mcO_U$ at $q$.
Recall that  $\widehat{\mcO}_{U, q} \cong k [[x,y]]/(xy)$.
By using  the injection  $F_{U/k}^* : \widehat{\mcO}_{U^{(1)}, q^{(1)}} \migi \widehat{\mcO}_{U, q}$ induced by $F_{U/k}$, we shall identify $\widehat{\mcO}_{U^{(1)}, q^{(1)}}$ with the subring $k [[x^p,y^p]]/(x^py^p)$ of $k [[x,y]]/(xy)$.
In particular,  we have $\widehat{\mfm}_{U^{(1)}, q^{(1)}} \cong x^p \cdot k[[x^p]] \oplus y^p \cdot k[[y^p]]$).
Now, let us define a $k [[x^p,y^p]]/(x^py^p)$-linear  morphism 
\begin{equation}
(k [[x^p,y^p]]/(x^py^p)) \oplus (x^p \cdot k[[x^p]] \oplus y^p \cdot k[[y^p]])^{\oplus (p-1)} \migi k [[x,y]]/(xy)
\end{equation}
given by 
\begin{equation}
(A, (x^p \cdot B_i + y^p \cdot C_i)_{i=1}^{p-1}) \mapsto A+ \sum_{i =1}^{p-1}(x^i \cdot B_i+ y^i \cdot C_i)
 \end{equation}
for $A \in k [[x^p,y^p]]/(x^py^p)$, and $B_i \in k[[x^p]]$, $C_i \in k[[y^p]]$ ($i = 1, \cdots, p-1$).
This morphism is verified to be bijective.
Thus, $\widehat{\mcO}_{U, q} $ is isomorphic to $\widehat{\mcO}_{U^{(1)}, q^{(1)}} \oplus \widehat{\mfm}_{U^{(1)},q^{(1)}}^{\oplus (p-1)}$, which is torsion-free  (cf. ~\cite{Ja}, the discussion at the beginning of \S\,1.2).
This completes the proof of Proposition \ref{torsionE}.
\end{proof}
\vspace{5mm}

Also, the following proposition will be used in the proof of Proposition-Definition  \ref{y088}. 

\bpr  \label{torsion}
 Let $\mcV$ be a vector bundle on $U^{(1)}$ of rank $n>0$.
 Denote by
 \begin{equation}
 w_{\mcV} : \mcV \migi F_{U/S*}(F^*_{U/S}(\mcV))
 \end{equation}
the $\mcO_{U^{(1)}}$-linear morphism
 corresponding, 
via  the adjunction relation 
 ``$F^*_{U/S} (-) \dashv F_{U/S*}(-)$",
  to the identity morphism of $F^*_{U/S}(\mcV)$.
 Then, $w_{\mcV}$ is injective and its cokernel $\mr{Coker}(w_{\mcV})$ is relatively torsion-free of rank $n \cdot (p-1)$.
   \epr
\begin{proof}
As the statement is  of local nature, we may assume (since $\mcV$ is locally isomorphic to $\mcO_{U^{(1)}}^{\oplus n}$) that  
$\mcV = \mcO_{U^{(1)}} =: \mcO$.
It is verified that  $w_{\mcO}$ is injective over
the smooth  locus $U^{\mr{sm}}$ of $U$ relative to $S$.
Since the open 
 subscheme $U^{\mr{sm}}$ of $U$   is scheme-theoretically dense,
 $w_{\mcO}$ is  injective over the entire space $U$.

Next, note that the formation of   $w_{\mcO}$ is functorial with respect to base-change to  any $S$-scheme $s : S' \migi S$ in the following sense:
  if we write $U' : = S' \times_S U$ and $\mcO' := \mcO_{U'^{(1)}}$, then the square diagram 
\begin{align} \label{QS101}
\vcenter{\xymatrix@C=66pt@R=36pt{
(s \times \mr{id}_U)^*(\mcO)  \ar[r]^-{(s \times \mr{id}_U)^*(w_{\mcO})} \ar[d]_-{\wr}&(s \times \mr{id}_U)^*( F_{U/S*}(F^*_{U/S}(\mcO))) \ar[d]^-{\wr}\\
\mcO' \ar[r]_-{w_{\mcO'}}& F_{U'/S*}(F^*_{U'/S} (\mcO'))
}}
\end{align}
 is commutative, where the right-hand vertical arrow arises from the base-change map $(s \times \mr{id}_S)^*(F_{U/S*}(\mcO_U)) \isom F_{U'/S*}(s^*(\mcO_U))$.
 The above discussion applied  to $\mcO'$ implies the injectivity of  
 $w_{\mcO'}$.
 Hence,  
 $w_{\mcO}$ is   universally injective with respect to base-change to $S$-schemes.
By ~\cite{MAT}, p.\,17, Theorem 1,   
the cokernel $\mr{Coker}(w_{\mcO})$ turns out to be flat over $S$.
Thus,  in order to finish the proof,
it suffices to consider the case where  $S = \mr{Spec} (k)$  and $k$ is algebraically closed.
We shall prove  the torsion-freeness of  $\mr{Coker}(w_{\mcO})$ at a {\it nodal} point 
(resp., a {\it smooth} point) $q^{(1)}$ of $U^{(1)}$.
Let us keep the notation in the proof of Proposition \ref{torsionE}.
According to  the discussion there,
the $\mfm_{U^{(1)}, q^{(1)}}$-adic completion of $w_{\mcO}$
  may be identified with
the inclusion $\widehat{\mcO}_{U^{(1)}, q^{(1)}} \migiincl \widehat{\mcO}_{U^{(1)}, q^{(1)}} \oplus \widehat{\mfm}_{U^{(1)},q^{(1)}}^{\oplus (p-1)}$
(resp., $\widehat{\mcO}_{U^{(1)}, q^{(1)}} \migiincl \widehat{\mcO}_{U^{(1)}, q^{(1)}}^{\oplus p}$)
 into the first factor.
Hence, the completion of  $\mr{Coker}(w_{\mcO})$ at $q^{(1)}$ is isomorphic to 
$\widehat{\mfm}_{U^{(1)},q^{(1)}}^{\oplus (p-1)}$ (resp.,  $\widehat{\mcO}_{U^{(1)}, q^{(1)}}^{\oplus (p-1)}$), which is torsion-free of rank $p-1$.
This completes the proof of Proposition \ref{torsion}.
\end{proof}

\section{Lie algebroids in positive characteristic} \label{y066}

In this section, we   recall the definition of a restricted Lie algebra.
The  $p$-curvature of a flat $\mbG$-bundle will be defined by means of the structure of restricted Lie algebra on the associated Lie algebroid.

\subsection{Restricted Lie algebras} \label{QG500}

First, let us recall the notion of a restricted Lie algebra. 

\bde \label{QS1499}
Let $R$ be a (unital, associative, and commutative) $k$-algebra.
A {\bf restricted Lie algebra} (or a {\bf $p$-Lie algebra}) 
\index{restricted Lie algebra (= $p$-Lie algebra)}
over $R$ is a Lie algebra $\mfh$ over $R$ together with a map $(-)^{[p]} :\mfh \migi \mfh$ satisfying the following conditions:
\begin{itemize}
\item
$(\mr{ad}(v))^p = \mr{ad}(v^{[p]})$ for all $v \in \mfh$;
\item
$(a\cdot v)^{[p]} = a^p \cdot v^{[p]}$ for all $a \in R$, $v \in \mfh$;
\item
$(v + u)^{[p]} = v^{[p]}+ u^{[p]} + \sum_{i=1}^{p-1} \frac{s_i (v, u)}{i}$ for all $v, u \in \mfh$, where  $s_i (v, u)$ denotes the coefficient of $t^{i-1}$ in the  expression $\mr{ad}(vt + u)^{p-1}(v) \in \mfh [t] \left(:=\mfh \otimes_R R[t]\right)$.
\end{itemize}
We shall refer to the operation $(-)^{[p]}$ as the  {\bf $p$-power operation} on $\mfh$. \index{$(-)^{[p]}$, $p$-power operation}
\ede

\begin{exa} \label{QR78}
Let $\mfh$ be a unital associative $R$-algebra.
Then, the bracket operation $[v, u] := vu -uv$ and the $p$-power operation $v^{[p]} := v^p$ make $\mfh$ into a restricted Lie algebra over $R$.
For example, 
the space of $R$-linear endomorphisms $\mr{End}_R(V)$ of $R$-module $V$
 gives a typical example of  a restricted Lie algebra over $R$.
\end{exa}

\subsection{The $p$-power operation on the tangent bundle} \label{QQ080}
Let $f : Y \migi T$ be a scheme over a $k$-scheme $T$.
In the natural fashion, we can define the notion of 
a {\it restricted $\mcO_T$-Lie algebra} on $Y$ (i.e., a sheaf of restricted $\mcO_T$-Lie algebra in the sense of ~\cite{Lan}, \S\,4.1).
Just as in the case of $\mr{End}_R(V)$ (cf. Example \ref{QR78} above),
the sheaf   $\mcE nd_{f^{-1}(\mcO_T)} (\mcV)$  for each vector bundle $\mcV$ (cf. \S\,\ref{QS4868}) has a natural structure of restricted $\mcO_T$-Lie algebra on $Y$ with 
  $p$-power operation $v^{[p]} := v^p$ (= the $p$-th iteration of $v$) for any local section $v \in \mcE nd_{f^{-1}(\mcO_T)} (\mcV)$.
Recall that the tangent  bundle $\mcT_{Y/T}$ may be identified with the subsheaf of 
$\mcE nd_{f^{-1}(\mcO_T)} (\mcO_Y)$ consisting of 
locally defined derivations on $\mcO_Y$ relative to $T$.
The $p$-power operation on $\mcE nd_{f^{-1}(\mcO_T)} (\mcO_Y)$ restricts to  a structure of restricted $\mcO_T$-Lie algebra
\begin{align} \label{QS1410}
(-)^{[p]} : \mcT_{Y/T} \migi \mcT_{Y/T}
\end{align}
on $\mcT_{Y/T}$.
Indeed, let us take   a local section $\partial$ of $\mcT_{Y/T}$, which corresponds to a locally defined derivation.
The $p$-th iterate $\partial^{[p]} := \partial^p$ of the derivation $\partial$ is verified to form a derivation, i.e., specifies a local section of $\mcT_{Y/T}$.
This implies that $\mcT_{Y/T}$ is closed under  the $p$-power operation on $\mcE nd_{f^{-1}(\mcO_T)} (\mcO_Y)$.

More generally, this discussion can be generalized to log schemes. 
Let $f : Y^\mr{log} \migi T^\mr{log}$ be a morphism of fine log schemes over $k$.
If $\partial := (D, \delta)$ (cf. (\ref{QR249})) is the logarithmic derivation corresponding to a local section of $\mcT_{Y^\mr{log}/T^\mr{log}}$, 
then 
the pair  
\begin{align} \label{QR98}
\partial^{[p]} \left(=(D, \delta)^{[p]} \right) := (D^p, F_Y^*\circ \delta + D^{p-1}\circ \delta)
\end{align}
specifies a logarithmic derivation.
Moreover, by  ~\cite{Og},  Proposition 1.2.1 (or, ~\cite{Ogu3}, Chap.\,V, \S\,4, Theorem 4.1.4), the equality 
\begin{align} \label{QW87}
(a \cdot \partial)^{[p]} = a^p \cdot \partial^{[p]} + (a \cdot \partial)^{p-1}(a) \cdot \partial \left(= a^p \cdot \partial^{[p]} - a \cdot \partial^{p-1}(a^{p-1}) \cdot \partial \right)
\end{align}
holds for any local section $a \in \mcO_Y$.
The sheaf $\mcT_{Y^\mr{log}/T^\mr{log}}$ together with the resulting operation
\begin{align} \label{QS1466}
 (-)^{[p]} : \mcT_{Y^\mr{log}/T^\mr{log}} \migi \mcT_{Y^\mr{log}/T^\mr{log}}
 \end{align}
 forms a restricted $\mcO_T$-Lie algebra on $Y$.
If the log structures on  $Y^\mr{log}$ and $T^\mr{log}$ are trivial, the  $p$-power structure on $\mcT_{Y^\mr{log}/T^\mr{log}}$ coincides with  (\ref{QS1410}).

\subsection{The restricted Lie algebra of an algebraic group} \label{QS32}

The Lie algebra of any smooth affine algebraic group admits a structure of restricted Lie algebra.
  Indeed,  the $p$-power operation 
$\partial \mapsto \partial^{[p]}$
on $\mcT_{\mbG/k}$
 defines  a structure of restricted Lie algebra 
\begin{align} \label{QS28}
(-)^{[p]} : \mfg \migi \mfg
\end{align}
 on the Lie algebra $\mfg$ via  the  identification of  $\mfg$ with the space of  right-invariant vector fields on $\mbG$.
Occasionally,  this $p$-power operation will be regarded as an endomorphism of the $k$-scheme $\mfg$.
According to ~\cite{Sel}, Chap.\,I, \S\,8, the second Corollary of Theorem I.8.1,
a $p$-power operation on $\mfg$
making it into 
  a restricted Lie algebra is uniquely determined if the Killing form on $\mfg$ is nondegenerate.

Let   $\mbG'$ be  another connected smooth affine algebraic group over $k$ with Lie algebra $\mfg'$ and $w : \mbG \migi \mbG'$ a morphism of algebraic groups over $k$.
Then, the differential $d w : \mfg \migi \mfg'$ preserves the $p$-power operations  (cf. ~\cite{Spr}, 4.4.9).
For example, if $\mbG' = \mr{GL}(\mfg)$ and $w = \mr{Ad}_\mbG$, then this fact just means
the first condition appearing in Definition \ref{QS1499}.
Also, if $\mbG' = \mbG$ and $w = \mr{Inn}_h$ for $h \in \mbG$,
then  it says that the automorphism $d \mr{Inn}_h |_e$ of $\mfg$    preserves the $p$-power operations.
 This implies that, for a $\mbG$-bundle  $\mcE$  on $Y$, 
the $p$-power operation  on $\mfg$ induces   a structure of   restricted $\mcO_Y$-Lie algebra
\begin{align} \label{QS1469}
(-)^{[p]} : \mfg_\mcE \migi \mfg_\mcE
\end{align}
  on the adjoint bundle $\mfg_\mcE$.
  In the case where  $\mcE = Y \times_k \mbG$, 
   the  $p$-power operation on $\mcO_Y \otimes_k \mfg \left(=\mfg_\mcE\right)$
   is given by $a \otimes v \mapsto a^p \otimes v^{[p]}$ for any $a \in \mcO_Y$ and $v \in \mfg$.

Moreover, the assignment 
$(\mcE, R) \mapsto  (\mcE, R)^{[p]} := (\mcE, R^{[p]})$ (for a $\mbG$-bundle $\mcE$ on a $k$-scheme and a global section $R$ of $\mfg_\mcE$)
defines a $k$-endomorphism 
\begin{align} \label{QS26}
(-)^{[p]} :  [\mfg/\mbG] \migi  [\mfg/\mbG]
\end{align}
of the quotient stack $[\mfg /\mbG]$, which is compatible with the $p$-power operation on $\mfg$ via the natural quotient $\mfg \migi [\mfg/\mbG]$.

\subsection{Restricted Lie algebroids} \label{QG567}

In the following, we introduce  the notion  of a restricted Lie algebroid defined on a log scheme.
The non-logarithmic version can be found   in  ~\cite{Lan}, \S\,4.2, Definition 4.2.



\bde[] \label{QS46}
\begin{itemize}
\item[(i)]
A {\bf restricted Lie algebroid} \index{restricted Lie algebroid} on $Y^\mr{log}/T^\mr{log}$ is a triple 
\begin{align} \label{QS47}
(\mcV,  \alpha, (-)^{[p]})
\end{align}
 consisting of a Lie algebroid $(\mcV, \alpha)$ on $Y^\mr{log}/T^\mr{log}$ (cf. Definition \ref{QS39}) and a map $(-)^{[p]} : \mcV \migi \mcV$
 satisfying the following conditions:
 \begin{itemize}
 \item
 $(\mcV, (-)^{[p]})$ forms a restricted $\mcO_T$-Lie algebra on $Y$;
 \item
 The anchor map $\alpha$ commutes with the respective $p$-power operations on $\mcV$ and $\mcT_{Y^\mr{log}/T^\mr{log}}$;
 \item
 The equality
\begin{align} \label{QS48}
(a \cdot \partial)^{[p]} = a^p \cdot \partial^{[p]} + \alpha (a \cdot \partial)^{p-1}(a) \cdot \partial
\end{align}
holds for any local sections $a \in \mcO_Y$, $\partial \in \mcV$.
 \end{itemize}
 If there is no fear of confusion, we  only use the symbol $\mcV$  to indicate the restricted  $\mcO_T$-Lie algebroid $(\mcV, \alpha, (-)^{[p]})$.
 \item[(ii)]
 Let $\mcV$, $\mcV'$ be restricted Lie algebroids on $Y^\mr{log}/T^\mr{log}$.
 Then, a {\bf morphism of restricted Lie algebroids} from $\mcV$ to $\mcV'$ is defined as a morphism of Lie algebroids $\mcV \migi \mcV'$ (cf. Definition \ref{QG200}) compatible with the respective  $p$-power operations on $\mcV$ and $\mcV'$.
 \end{itemize}
 \ede

\subsection{Restricted Lie algebroids associated to  $\mbG$-bundles} \label{QS1433}

Let $\pi : \mcE \migi Y$ be a $\mbG$-bundle on $Y$.
The $p$-power operation on  $\mcT_{\mcE^\mr{log}/T^\mr{log}}$  restricts to 
a $p$-power operation $(-)^{[p]}$ on 
the Lie algebroid $\widetilde{\mcT}_{\mcE^\mr{log}/T^\mr{log}} \left(= (\pi_*(\mcT_{\mcE^\mr{log}/T^\mr{log}}))^\mbG\right)$.
By construction, 
it is consistent, via restriction,  with  the $p$-power operation on the adjoint bundle $\mfg_\mcE \left(\subseteq \widetilde{\mcT}_{\mcE^\mr{log}/T^\mr{log}} \right)$ constructed in \S\,\ref{QS32}, and 
the anchor map $d_\mcE^\mr{log}$ preserves the $p$-power operations.
Moreover,
by considering the local situation  (i.e., the case where $\mcE = Y \times_k \mbG$), we can verify
  (\ref{QS48})  for   $\widetilde{\mcT}_{\mcE^\mr{log}/T^\mr{log}}$.
That is to say,
the triple 
\begin{align} \label{QS1500}
 (\widetilde{\mcT}_{\mcE^\mr{log}/T^\mr{log}}, d_\mcE^\mr{log}, (-)^{[p]})
 \end{align}
 forms  a restricted Lie algebroid  on $Y^\mr{log}/T^\mr{log}$ associated to $\mcE$.


\section{$p$-curvature  on an $\hslash$-flat $\mbG$-bundle} \label{QS3338}

Next, we  shall recall the notion of $p$-curvature and discuss basic properties concerning $p$-curvature.
Roughly speaking, the $p$-curvature measures 
the obstruction to the compatibility, with respect to an $\hslash$-$T^\mr{log}$-connection,
of 
 the $p$-power operation $\hslash^{p-1} \cdot (-)^{[p]}$ on $\mcT_{Y^\mr{log}/T^\mr{log}}$ and the  $p$-power operation on $\widetilde{\mcT}_{\mcE^\mr{log}/T^\mr{log}}$ defined above.

\subsection{$p$-curvature}

Let $\mcE$ be as above.
Also, let $\hslash$ be an element of $\Gamma (T, \mcO_T)$ and $\nabla$  a flat $\hslash$-$T^\mr{log}$-connection on $\mcE$.
For each local section $\partial \in \mcT_{Y^\mr{log}/T^\mr{log}}$,
we write 
\begin{align} \label{QG601}
{^p}\psi^\nabla (\partial) := \nabla (\partial)^{[p]} -\hslash^{p-1} \cdot \nabla (\partial^{[p]}).
\end{align}
Then, the following assertion holds:

\ble \label{QW82}
 Let $\partial, \partial_1, \partial_2$ be local sections of $\mcT_{Y^\mr{log}/T^\mr{log}}$,  and  let  $a$ be  a  local section of  $\mcO_Y$.
(Suppose that these sections are defined on a common open subset of $Y$.)
Then, the following equalities hold:
\begin{align} \label{QG600}
{^p}\psi^{\nabla}(a\cdot \partial) = a^p \cdot {^p}\psi^\nabla (\partial),
\hspace{5mm}
{^p}\psi^\nabla (\partial_1 + \partial_2) = {^p}\psi^\nabla (\partial_1) + {^p}\psi^\nabla (\partial_2).
\end{align}
\ele
\begin{proof}
The first equality  follows from the following sequence of equalities:
\begin{align}
{^p}\psi^\nabla (a \cdot \partial) = & \ \nabla (a \cdot \partial)^{[p]}-\hslash^{p-1}  \cdot \nabla ((a \cdot \partial)^{[p]})
\\
= &  \  (a \cdot \nabla (\partial))^{[p]} - \hslash^{p-1} \cdot \nabla (a^p \cdot \partial^{[p]} +  (a \cdot \partial)^{p-1}(a) \cdot \partial)\notag \\
= & \ a^p \cdot \nabla (\partial)^{[p]} + d_\mcE^\mr{log} (a \cdot \nabla (\partial))^{p-1} (a) \cdot \nabla (\partial) \notag
\\
& \ - \hslash^{p-1} \cdot \nabla (a^p \cdot \partial^{[p]})  -\hslash^{p-1} \cdot \nabla ( (a \cdot \partial)^{p-1}(a) \cdot \partial)\notag \\
= & \ a^p \cdot \nabla (\partial)^{[p]} + \cancel{(a \cdot  \hslash \cdot \partial)^{p-1}(a) \cdot \nabla (\partial)}  \notag \\
& \ - \hslash^{p-1} \cdot \nabla (a^p \cdot \partial^{[p]})  - \cancel{\nabla ( (\hslash \cdot a \cdot \partial)^{p-1}(a) \cdot \partial)}\notag \\
= & \ a^p \cdot ((\nabla (\partial))^{[p]}-\hslash^{p-1}  \cdot \nabla (\partial^{[p]})) \notag \\
= & \ a^p \cdot {^p}\psi^\nabla (\partial), \notag
\end{align}
where the first and second equalities follow from (\ref{QW87}) and (\ref{QS48}) respectively.

To verify the second equality, we use the flatness of $\nabla$.
Indeed, since $\nabla$ preserves  the Lie bracket operations,
 the equality $s_i (\nabla (\partial_1), \nabla (\partial_2)) = h^{p-1} \cdot \nabla (s_i (\partial_1, \partial_2))$ holds for every $i$ (cf. Definition \ref{QS1499} for the definition of $s_i$).
 Hence, we have 
 \begin{align}
 {^p}\psi^\nabla (\partial_1 + \partial_2) = & \ \nabla (\partial_1 + \partial_2)^{[p]}-\hslash^{p-1} \cdot \nabla ((\partial_1 + \partial_2)^{[p]}) \\
 = &\ (\nabla (\partial_1)+ \nabla (\partial_2))^{[p]} -\hslash^{p-1} \cdot \nabla \left(\partial^{[p]}_1 + \partial_2^{[p]} + \sum_{i=1}^{p-1} \frac{s_i (\partial_1, \partial_2)}{i}\right) \notag \\
 = & \ \nabla (\partial_1)^{[p]} + \nabla (\partial_2)^{[p]} + \cancel{\sum_{i=1}^{p-1} \frac{s_i (\nabla (\partial_1), \nabla (\partial_2))}{i}}  \notag \\
 &\  - \hslash^{p-1} \cdot \nabla (\partial_1^{[p]}) -\hslash^{p-1} \cdot \nabla (\partial_2^{[p]}) - \cancel{\sum_{i=1}^{p-1} \frac{\hslash^{p-1} \cdot \nabla (s_i (\partial_1, \partial_2))}{i}} \notag  \\
 = & \ \left(\nabla (\partial_1)^{[p]} -\hslash^{p-1} \cdot \nabla (\partial_1^{[p]})  \right) + \left(\nabla (\partial_2)^{[p]} -\hslash^{p-1} \cdot \nabla (\partial_2^{[p]})  \right) \notag \\
 = & \ {^p}\psi^\nabla (\partial_1) + {^p}\psi^\nabla (\partial_2), \notag
 \end{align} 
 thus completing  the proof of the second equality.
\end{proof}


The above lemma implies  that
 the assignment $F^{-1}_Y(\partial) \mapsto {^p}\psi^\nabla (\partial)$ (for any $\partial \in \mcT_{Y^\mr{log}/T^\mr{log}}$) defines 
 an {\it  $\mcO_Y$-linear} morphism
\begin{align}   \label{QR100}
 {^p}\psi^\nabla : F_Y^*(\mcT_{Y^\mr{log}/T^\mr{log}}) \migi  \widetilde{\mcT}_{\mcE^\mr{log}/T^\mr{log}}. 
  \end{align}

\bde \label{y068} 
We shall refer to
 ${^p \psi}^{\nabla}$  as the {\bf $p$-curvature}  of $\nabla$.
 \index{$\nabla$, $\hslash$-$T^\mr{log}$-connection@--- ${^p}\psi^\nabla$, $p$-curvature of}
 Also, $\nabla$ is called {\bf $p$-flat} (or {\bf $p$-integrable})
 \index{$\nabla$, $\hslash$-$T^\mr{log}$-connection@--- $p$-flat (= $p$-integrable)}
   if   ${^p}\psi^\nabla =0$.
 \ede

\begin{rema} \label{QG665}
Since $d_\mcE^\mr{log}$ satisfies  $d^\mr{log}_\mcE \circ \nabla = \hslash \cdot \mr{id}_{\mcT_{Y^\mr{log}/T^\mr{log}}}$ and  
 preserves the  $p$-power operations, 
the image of ${^p}\psi^\nabla$ is contained in $\mfg_\mcE \left(=\mr{Ker}(d^\mr{log}_\mcE)\right)$.
Thus,  the morphism  ${^p}\psi^\nabla$ may be regarded as 
an $\mcO_Y$-linear morphism $F^*_Y(\mcT_{Y^\mr{log}/T^\mr{log}}) \migi \mfg_\mcE$,  or equivalently 
 a global section of $F^*_Y(\Omega_{Y^\mr{log}/T^\mr{log}})\otimes \mfg_\mcE$. 
\end{rema}

\bpr \label{QS82}
Let us keep the above notation.
Also, let $\eta$ be an automorphism of the $\mbG$-bundle   $\mcE$ and denote by $\mr{Ad}_\mbG (\eta)$ the automorphism of $\mfg_\mcE$ induced by $\eta$ via change of structure group by $\mr{Ad}_\mbG$.
Then, the equality
\begin{align}
{^p}\psi^{\eta_*(\nabla)} = \mr{Ad}_\mbG (\eta) \circ  {^p}\psi^\nabla
\end{align}
(cf. (\ref{comps}) for the definition of $\eta_*(\nabla)$) of morphisms $F_Y^*(\mcT_{Y^\mr{log}/T^\mr{log}}) \migi  \mfg_\mcE$ holds.
In particular, the $p$-flatness of an $\hslash$-$T^\mr{log}$-connection is closed under gauge transformations.
\epr
\begin{proof}
The assertion follows immediately from the definition of $p$-curvature together with   the fact that 
the automorphism $d \eta^\mbG : \widetilde{\mcT}_{\mcE^\mr{log}/T^\mr{log}} \isom  \widetilde{\mcT}_{\mcE^\mr{log}/T^\mr{log}}$ (cf. (\ref{QQ10})) preserves  the $p$-power operations and  restricts to $\mr{Ad}_\mbG (\eta)$.
\end{proof}

\subsection{Change of structure group}\label{QS60}
Let  $\mbG'$ be  another connected smooth affine algebraic group over $k$ with  Lie algebra $\mfg'$ and $w : \mbG \migi \mbG'$  a morphism of algebraic groups over $k$.
Also, let $(\mcE, \nabla)$ be  an $\hslash$-flat $\mbG$-bundle on  $Y^\mr{log}/T^\mr{log}$.
Since 
the differential  $dw : \mfg \migi \mfg'$ of $w$ preserves  the $p$-power operations,
 the morphism  of Lie algebroids $(dw)_\mcE : \widetilde{\mcT}_{\mcE^\mr{log}/T^\mr{log}} \migi \widetilde{\mcT}_{(\mcE \times^\mbG \mbG')^\mr{log}/T^\mr{log}}$ (cf. (\ref{QQ810})) preserves   the $p$-power operations, i.e., defines a morphism of restricted Lie algebroids.
Hence, the $p$-curvature ${^p}\psi^{w_*(\nabla)}$ of the $\hslash$-$T^\mr{log}$-connection $w_*(\nabla)$ (cf. (\ref{QQ813})) on $\mcE \times^\mbG \mbG'$ satisfies 
\begin{align} \label{QS61}
{^p}\psi^{w_*(\nabla)} = (d w)_\mcE \circ {^p}\psi^{\nabla}.
\end{align}
In particular, 
the $p$-flatness of $\nabla$
 implies 
 that of $w_*(\nabla)$,
 and vice versa if $dw$ is injective. 


\subsection{Change of parameter} \label{QS62}
Let $(\mcE, \nabla)$  be  as above and 
 $c$  an element of $\Gamma (T, \mcO_T)$.
Then, the $p$-curvature ${^p}\psi^{c \cdot \nabla}$ of the flat  $(c \cdot \hslash)$-$T^\mr{log}$-connection $c \cdot \nabla$ (cf. (\ref{QQ1032})) is given by 
\begin{align} \label{QR101}
{^p}\psi^{c\cdot \nabla} = c^p \cdot {^p}\psi^\nabla.
\end{align}
In particular, 
the $p$-flatness of $\nabla$
implies
that of $c \cdot \nabla$,
 and vice versa if $c \in \Gamma (T, \mcO_T^\times)$.

\subsection{Products} \label{QS63}
Let $\mbG'$ be as in \S\,\ref{QS60}, and 
 let $(\mcE, \nabla)$ (resp., $(\mcE', \nabla')$) be an $\hslash$-flat $\mbG$-bundle (resp., an $\hslash$-flat $\mbG'$-bundle) on $Y^\mr{log}/T^\mr{log}$.
As discussed in \S\,\ref{QS7},  the product $(\mcE \boxtimes \mcE', \nabla \boxtimes \nabla')$, forming an $\hslash$-flat $\mbG \times_k \mbG'$-bundle, can be constructed.
Note that the inclusion $\widetilde{\mcT}_{(\mcE \boxtimes \mcE')^\mr{log}/T^\mr{log}} \migiincl \widetilde{\mcT}_{\mcE^\mr{log}/T^\mr{log}} \oplus \widetilde{\mcT}_{\mcE'^{\mr{log}}/T^\mr{log}}$ (cf. (\ref{QS9})) preserves the  $p$-power  operations, where the $p$-power operation on the direct sum $\widetilde{\mcT}_{\mcE^\mr{log}/T^\mr{log}} \oplus \widetilde{\mcT}_{\mcE'^{\mr{log}}/T^\mr{log}}$ is the operation obtained  naturally from both factors.
Hence, the $p$-curvature ${^p}\psi^{\nabla \boxtimes \nabla'}$ of $\nabla \boxtimes \nabla'$ satisfies, via this inclusion,  the following equality:
\begin{align}
\label{QS68}
{^p}\psi^{\nabla \boxtimes \nabla'} = {^p}\psi^\nabla \oplus  {^p}\psi^{\nabla'}.
\end{align}
In particular, 
$\nabla \boxtimes \nabla'$ is $p$-flat 
 if and only if  both $\nabla$ and $\nabla'$ are $p$-flat.

\subsection{Canonical connections} \label{y071}

We shall construct  a canonical  
connection on the Frobenius pull-back of a $\mbG$-bundle.
Let $\pi : \mcG \migi Y^{(1)}$ be a $\mbG$-bundle on  $Y^{(1)}$.
Denote by 
\begin{equation} \label{QR102}
 \pi_F: F^*_{Y/T}(\mcG) \migi Y
\end{equation}
 the $\mbG$-bundle on  $Y$ 
defined as  the base-change of 
 $\mcG$
   via  $F_{Y/T} :Y \migi Y^{(1)}$.
Also, denote by 
 \begin{equation} 
  F_{\mcG/Y/T} : F^*_{Y/T}(\mcG) \ \left(:= Y \times_{Y^{(1)}, \pi} \mcG\right) \migi \mcG\end{equation}
the projection to the second factor (i.e., $F_{\mcG/Y/T} := F_{Y/T}\times \mr{id}_\mcG$), which commutes  with the respective $\mbG$-actions of $F^*_{Y/T}(\mcG)$ and $\mcG$.
Differentiating this morphism yields an $\mcO_{F^*_{Y/T}(\mcG)}$-linear  isomorphism
\begin{align} \label{QR105}
 \mcT_{F^*_{Y/T}(\mcG)^\mr{log}/Y^\mr{log}} \isom 
 F_{\mcG/Y/T}^*(\mcT_{\mcG^\mr{log}/Y^{(1) \mr{log}}}). 
\end{align}
Consider a morphism of short exact sequences  of $\mcO_{F^*_{Y/T}(\mcG)}$-modules
\begin{align} \label{canconn}
\vcenter{\xymatrix@C=9pt@R=36pt{
0 \ar[r] &\mcT_{F^*_{Y/T}(\mcG)^\mr{log}/Y^\mr{log}}   \ar[r]\ar[d]_-{\wr}^-{(\ref{QR105})}
&  \mcT_{F^*_{Y/T}(\mcG)^\mr{log}/T^\mr{log}} \ar[r]\ar[d]^-{d F_{\mcG/Y/T}}
&  \pi_F^{*}(\mcT_{Y^\mr{log}/T^\mr{log}}) \ar[r] \ar[d]^-{dF_{Y/T}} 
& 0 \\
0 \ar[r]&F_{\mcG/Y/T}^*(\mcT_{\mcG^\mr{log}/Y^{(1) \mr{log}}})   \ar[r]
&F^*_{\mcG/ Y/T}(\mcT_{\mcG^\mr{log}/T^\mr{log}})  \ar[r] 
&\pi_F^{*} (F_{Y/T}^*(\mcT_{Y^{(1) \mr{log}}/T^\mr{log}}))  \ar[r] & 0, 
}}\hspace{-6mm}
\end{align}
 where 
\begin{itemize}
\item
the upper horizontal sequence is obtained by differentiating $\pi_F$;
\item
the lower horizontal  sequence is obtained by differentiating $\pi$ and successively  pulling-back via $F_{\mcG/Y/T}$ under the natural isomorphism
\begin{equation}
 \pi_F^{*} (F_{Y/T}^*(\mcT_{Y^{(1) \mr{log}}/T^\mr{log}}))  \isom F^*_{\mcG/Y/T}(\pi^*(\mcT_{Y^{(1) \mr{log}}/T^\mr{log}}));
  \end{equation}
\item
the middle and right-hand vertical  arrows are obtained by   differentiating  the morphisms $F_{\mcG/Y/T}$ and $F_{Y/T}$ respectively. 
\end{itemize}
By applying the functor $\pi_{F*}(-)^\mbG$ to this diagram, we obtain
a morphism of short exact sequence of 
 $\mcO_Y$-modules
\begin{align} \label{canconn2}
\vcenter{\xymatrix@C=16pt@R=36pt{
0 \ar[r] &\mfg_{F^*_{Y/T}(\mcG)}  \ar[r]\ar[d]_-{\wr}
&\widetilde{\mcT}_{F^*_{Y/T}(\mcG)^\mr{log}/T^\mr{log}}   \ar[r]^-{d_{F^*_{Y/T}(\mcG)}^\mr{log}}\ar[d]^-{\pi_{F*}(dF_{\mcG/Y/T})^\mbG}
& \mcT_{Y^\mr{log}/T^\mr{log}}  \ar[r] \ar[d]^-{\pi_{F*}(dF_{Y/T})^\mbG} 
& 0 \\
0 \ar[r]&F_{Y/T}^*(\mfg_\mcG) \ar[r]
& \pi_{F*}(F^*_{\mcG/Y/T}(\mcT_{\mcG^\mr{log}/T^\mr{log}}))^\mbG \ar[r] 
& F_{Y/T}^*(\mcT_{Y^{(1) \mr{log}}/T^\mr{log}}) \ar[r] & 0. 
}} \hspace{-5mm}
\end{align}
Since  $d (F^*_{Y/T}(a)) = 0$  for any $a \in \mcO_{Y^{(1)}}$,
the  right-hand vertical  arrow $\pi_{F*}(dF_{Y/T})^\mbG$ vanishes identically on $Y$.
Hence, the middle vertical  arrow $\pi_{F*}(dF_{\mcG/Y/T})^\mbG$
  factors through the injection
$F_{Y/T}^*(\mfg_\mcG)  \migiincl  \pi_{F*}(F^*_{\mcG/Y/T}(\mcT_{\mcG^\mr{log}/T^\mr{log}}))^\mbG$.
The resulting morphism 
$\widetilde{\mcT}_{F^*_{Y/T}(\mcG)^\mr{log}/T^\mr{log}}  \migi  F_{Y/T}^*(\mfg_\mcG)$
 determines, via the left-hand vertical  arrow $\mfg_{F^*_{Y/T}(\mcG)} \isom F_{Y/T}^*(\mfg_\mcG)$, a split surjection of the upper horizontal  sequence in (\ref{canconn2}).
Thus, 
the corresponding split injection
\index{$\nabla$, $\hslash$-$T^\mr{log}$-connection@--- $\nabla^\mr{can}_{\mcG}$, $\nabla^{\mr{can}}_{\mcG, \hslash}$,  $\nabla_\mcW^\mr{can}$, canonical}
\begin{equation}   \label{QR111}
\nabla^\mr{can}_{\mcG} :  \mcT_{Y^\mr{log}/T^\mr{log}} \migi   \widetilde{\mcT}_{F^*_{Y/T}(\mcG)^\mr{log}/T^\mr{log}}
 \end{equation}
of that sequence 
specifies  a $T^\mr{log}$-connection on $F^*_{Y/T}(\mcG)$.
We refer to it as the {\bf  canonical $T^\mr{log}$-connection} on $F^*_{Y/T}(\mcG)$  (cf. ~\cite{Kal}, Theorem (5.1), for the non-logarithmic version).
Moreover, given  $\hslash \in \Gamma (T, \mcO_T)$, 
we obtain an $\hslash$-$T^\mr{log}$-connection
\begin{equation} \label{QR110}
 \nabla^{\mr{can}}_{\mcG, \hslash} := \hslash \cdot \nabla_\mcG^\mr{can}\end{equation}
and  refer to it as the {\bf canonical $\hslash$-$T^\mr{log}$-connection}.

\bpr \label{y073}
Let us keep the above notation.
Then, the canonical $\hslash$-$T^\mr{log}$-connection $ \nabla^{\mr{can}}_{\mcG, \hslash}$ is flat and $p$-flat.
  \epr
\begin{proof}
By (\ref{EE1033}) and (\ref{QR101}),
 it suffices to consider  the case where $\hslash =1$.
Moreover, since the problem is of local nature,
we can assume that $\mcG = Y^{(1)} \times_k \mbG$, i.e., the trivial $\mbG$-bundle.
Then, since $F_{(Y^{(1)} \times_k \mbG)/Y/T} = F_{Y/T} \times \mr{id}_{\mbG}$,  the $T^\mr{log}$-connection  $\nabla^{\mr{can}}_{\mcG}$ coincides with the trivial connection $\nabla^\mr{triv}$ (cf. (\ref{QR340})) on $Y \times_k \mbG \left(= F^*_{Y/T}(Y^{(1)} \times_k \mbG)\right)$.
The flatness and   $p$-flatness of $\nabla^\mr{triv}$ can be verified immediately by considering the various definitions involved, so we  complete the proof of the proposition.
\end{proof}


\section{Dormant/$p$-nilpotent $(\mfg, \hslash)$-opers and their moduli spaces} \label{QS1712} 

In this section, we introduce dormant $(\mfg, \hslash)$-opers and $p$-nilpotent $(\mfg, \hslash)$-opers, which are certain $(\mfg, \hslash)$-opers characterized in terms of $p$-curvature.
Moreover, the moduli spaces classifying these opers  will be defined.

Hereafter,  suppose  that
$(\mbG, \mbT)$
\index{$\mbG$, algebraic group}\index{$\mbG$, algebraic group@— $\mbT$, maximal torus of}
is a split semisimple algebraic group over $k$ of adjoint type (where $\mbT$ denotes a maximal torus $\mbT$ of $\mbG$) 
with $p> 2  h_\mbG$.
Let us fix 
a Borel subgroup $\mbB$ 
\index{$\mbG$, algebraic group@— $\mbB$, Borel subgroup of}
defined over $k$ containing $\mbT$. 

\subsection{The adjoint quotient of the $p$-power operation}  \label{y075}

Since $\mbG$ is a Chevalley group, 
all the groups  
$\mbG$, $\mbT$,  and $\mbB$
 can 
 be defined over the prime field $\mbF_p := \mbZ/p\mbZ$.
It follows that the induced  $k$-scheme $\mfc$ (cf. (\ref{c}))  arises from an $\mbF_p$-scheme $\mfc_{\mbF_p}$, which means that  $\mfc = \mfc_{\mbF_p} \times_{\mbF_p} k$.
 In particular,  the  Frobenius twist    $\mfc^{(1)}$ is isomorphic to $\mfc$ itself, and hence, 
 the relative Frobenius morphism  $F_{\mfc /k}$ may be regarded as   a $k$-endomorphism of $\mfc$.
We denote by $\mfc (\mbF_p)$ the (finite) set of $\mbF_p$-rational points of $\mfc_{\mbF_p}$, i.e.,
the set of elements in $\mfc (k)$ invariant under  $F_{\mfc/k}$.
The element $[0] \in \mfc (k)$ lies in $\mfc (\mbF_p)$. For each nonnegative  integer $r$, we write 
\begin{align}\label{QS3001}
[0]^{\times r} := ([0], [0], \cdots, [0]) \in \mfc(\mbF_p)^{\times r},
\index{$[0]$, $[0]^{\times r}$, $[0]_U$}
\end{align}
where $[0]^{\times 0} := \emptyset$.

The following proposition will be used to prove Proposition \ref{y077} described  below.

\bpr \label{y076}
The square diagrams
\begin{align} \label{comm53}
\vcenter{\xymatrix@C=46pt@R=36pt{
\mfg \ar[r]^-{\chi} \ar[d]_-{(-)^{[p]}}& \mfc \ar[d]^-{F_{\mfc/k}}\\
\mfg \ar[r]_-{\chi}& \mfc
}}
\hspace{10mm}
\vcenter{\xymatrix@C=46pt@R=36pt{
[\mfg /\mbG]\ar[r]^-{[\chi]} \ar[d]_-{(-)^{[p]}}& \mfc\ar[d]^-{F_{\mfc/k}}\\
[\mfg/\mbG]\ar[r]_-{[\chi]}&\mfc
}}
\end{align}
are  commutative (cf. (\ref{QS26}) for the definition of the endomorphism $(-)^{[p]}$ of $[\mfg/\mbG]$).
  \epr
\begin{proof}
The $p$-power operation $(-)^{[p]}$ on $\mfg$ is compatible with the adjoint $\mbG$-action, so
it induces  an endomorphism $F'$ of $\mfc$ via  the quotient $\chi : \mfg \migi \mfc$.
On the other hand, the $p$-power  operation $(-)^{[p]}$ on $\mfg$  restricts to 
that on the Lie algebra $\mft$ of $\mbT$, which  coincides with the  Frobenius morphism $F_{\mft/k} : \mft \migi \mft \left(=\mft^{(1)}\right)$ (cf. ~\cite{Spr}, 4.4.10. Example (2)). 
Recall that  $\mfc$ can be  obtained as the quotient of $\mft$ by the $\mbW$-action and that the quotient $\mft \migisurj \mfc$ is compatible with both  the  Frobenius morphisms and the $p$-power operations.
Hence, $F'$ must  coincide with $F_{\mfc/k}$.
This completes 
the commutativity of the left-hand square in (\ref{comm53}), which also implies the  commutativity of the right-hand square. 
\end{proof}


\subsection{Radius of a $p$-flat connection} \label{QS1256}
Let $r$ be a positive  integer, $S$ a $k$-scheme,   and $\msU := (U/S, \{ \sigma^U_i : S_i \migi U \}_{i=1}^r)$ a local $r$-pointed stable curve over $S$.
After choosing a base for $\msU$, we obtain a log curve $U^\mr{log}/S^\mr{log}$ extending $U/S$.
Also, let $\hslash$ be an element of $\Gamma (S, \mcO_S)$ and $(\mcE, \nabla)$  an $\hslash$-flat $\mbG$-bundle on $U^\mr{log}/S^\mr{log}$.
For each  $i \in \{1, \cdots, r \}$, 
the restriction  $\sigma^{U*}_i ({^p}\psi^\nabla)$
 of ${^p}\psi^\nabla$ via  $\sigma^U_i$ may be regarded as an element of $\Gamma (S_i, \mfg_{\sigma_i^{U*}(\mcE)})$ under the
 identification $\mfg_{\sigma_i^{U*}(\mcE)} = \sigma^{U*}_i (F_U^*(\Omega_{U^\mr{log}/S^\mr{log}})\otimes \mfg_\mcE)$ given by 
 the  isomorphism
\begin{align}
F^*_{S_i}(\mr{Res}_{\msU, i}) :  \sigma^{U*}_i (F_{U}^*(\Omega_{U^\mr{log}/S^\mr{log}})) \left(= F_{S_i}^*(\sigma^{U*}_i (\Omega_{U^\mr{log}/S^\mr{log}})) \right) \isom \mcO_{S_i}
\end{align}
 (cf. (\ref{QQ040})).
By the definition of $p$-curvature, 
the following equality  
 holds: 
\begin{equation}
 \label{pcurvM} 
\sigma^{U*}_i ({^p}\psi^\nabla)  =  (\mu_i^{\nabla})^{[p]} - \hslash^{p-1} \cdot \mu_i^{\nabla}.
 \end{equation}

\bpr \label{y077} 
Let $(\mcE, \nabla)$   be as above, and 
suppose further that $\nabla$ is $p$-flat.
Moreover, we shall fix an element $i \in \{1, \cdots, r \}$.
Then,  the equality
\begin{equation} \label{Frho} 
F_{\mfc/k} \circ \rho_i^{\nabla} = \hslash^{p-1} \star  \rho_i^{\nabla}  \left(= \rho_i^{\hslash^{p-1} \cdot \nabla}\right)\end{equation}
of elements of $\mfc (S_i)$ holds (cf.   (\ref{action}) for the definition of $\star$).
In particular, the following assertions hold:
\begin{itemize}
\item[(i)]
If, moreover, $\hslash \in \Gamma (S, \mcO_S^\times)$, then $ \hslash^{-1} \star \rho_i^{\nabla} \in \mfc (\mbF_p)$;
\item[(ii)]
If, moreover,  $\hslash = 0$ and $S_i$ is reduced, then $\rho_i^{ \nabla} =[0] \in \mfc (\mbF_p)$.
\end{itemize}
  \epr
\begin{proof}
Since $\nabla$ is $p$-flat,
  (\ref{pcurvM}) is equivalent to the equality  $(\mu_i^{\nabla})^{[p]} = \hslash^{p-1} \cdot \mu_i^{\nabla}$.
By Proposition \ref{y076}, we have a sequence of equalities
  \begin{align}  
 F_{\mfc/k} \circ \rho_{i}^{\nabla}  
   &=  F_{\mfc/k} \circ [\chi] ((\sigma^{U*}_i(\mcE), \mu_{i}^{\nabla}))  \\
   & = [\chi]( (\sigma^{U*}_i(\mcE), \mu_{i}^{\nabla})^{[p]}) \notag  \\
& = [\chi]( (\sigma^{U*}_i(\mcE), (\mu_{i}^{\nabla})^{[p]})) \notag \\
& = [\chi] ((\sigma^{U*}_i(\mcE), \hslash^{p-1} \cdot \mu_i^{\nabla})) \notag \\
& =  \hslash^{p-1} \star \rho_i^{\nabla}, \notag
\end{align} 
thus completing the proof of  (\ref{Frho}).

Next, we consider assertions (i) and (ii).
Let us fix an isomorphism $\mfc \isom \mbA^{\mr{rk}(\mfg)} := \mr{Spec}(k [u_1, \cdots, u_{\mr{rk}(\mfg)}])$ as in (\ref{isomCA}).  
 The element $\rho_i^{\nabla}$  can be written, via this isomorphism,  as
\begin{equation} 
 \rho_i^{\nabla} = (a_1, \cdots, a_{\mr{rk}(\mfg)})  \left(=: (a_j)_{j=1}^{\mr{rk}(\mfg)}\right) \in \mbA^{\mr{rk}(\mfg)} ( S_i).
\end{equation}
for some $a_j \in \Gamma (S_i, \mcO_{S_i})$ ($j = 1, \cdots, \mr{rk}(\mfg)$).
Under this equality, 
we can describe $\hslash^{-1} \star \rho_i^\nabla$  as 
\begin{align}
\hslash^{-1} \star \rho_i^\nabla = (\hslash^{-e_1} \cdot a_1, \cdots, \hslash^{-e_{\mr{rk}(\mfg)}} \cdot a_{\mr{rk}(\mfg)}).
\end{align}
 Moreover,  (\ref{Frho}) reads 
 $(a_j^p)_{j=1}^{\mr{rk}(\mfg)} = (\hslash^{(p-1) \cdot e_j} \cdot a_j)_{j=1}^{\mr{rk}(\mfg)}$.
If $\hslash \in \Gamma (S, \mcO_S^\times)$ (i.e., the assumption in  (i) is satisfied), then,  for every $j =1, \cdots, \mr{rk}(\mfg)$,  we have 
\begin{equation} a_j^p = \hslash^{(p-1) \cdot e_j} \cdot a_j  \Longleftrightarrow (\hslash^{-e_j} \cdot a_j)^{p} = \hslash^{-e_j} \cdot a_j  \Longleftrightarrow \hslash^{-e_j} \cdot a_j \in \mbF_p.
\end{equation}
This implies that  $\hslash^{-1} \star \rho_i^{\nabla} \in \mfc (\mbF_p)$.
On the other hand, if $\hslash = 0$ and $S_i$ is reduced (i.e., the assumption in  (ii) is satisfied), then 
\begin{equation}  a_j^p = \hslash^{(p-1) \cdot e_j} \cdot a_j  \ \left(=0\right)\Longleftrightarrow a_j =0.\end{equation}
That is to say, the equality $\rho_i^\nabla =0$ holds.
This completes the proof of Proposition \ref{y077}.
\end{proof}

\subsection{Dormant/$p$-nilpotent $(\mfg, \hslash)$-opers}\label{QS130}
Let $S^\mr{log}$ be an fs log scheme over $k$ and 
$U^\mr{log}$ a log curve over $S^\mr{log}$.
For a line bundle $\mcL$ on $U$, we shall write $\mcL^\times$
for the $\mbG_m$-bundle  associated to $\mcL$.
Now, let $(\mcE, \nabla)$ be an $\hslash$-flat $\mbG$-bundle on $U^\mr{log}/S^\mr{log}$.
Then, the pair  $(\mcE, {^p}\psi^\nabla)$  specifies a global section  of the quotient stack $[(F_U^*(\Omega_{U^\mr{log}/S^\mr{log}})^\times \times^{\mbG_m} \mfg)/\mbG]$ over $U$; the composite  of  this section  with  the morphism
\begin{align}
[(F_X^*(\Omega_{U^\mr{log}/S^\mr{log}})^\times \times^{\mbG_m} \mfg)/\mbG] \migi F^*_U(\Omega_{U^\mr{log}/S^\mr{log}})^\times \times^{\mbG_m} \mfc 
\end{align}
induced by $\chi : \mfg \migi \mfc$ specifies a $U$-rational point 
\begin{align} \label{QS105}
{^p}\psi^\nabla_{/\hspace{-0.9mm}/\mbG} \in (F^*_U(\Omega_{U^\mr{log}/S^\mr{log}})^\times \times^{\mbG_m} \mfc) (U)
\index{${^p}\psi^\nabla_{/\hspace{-0.9mm}/\mbG}$}
\end{align}
of $F^*_U(\Omega_{U^\mr{log}/S^\mr{log}})^\times \times^{\mbG_m} \mfc$.
By Proposition \ref{QS82}, ${^p}\psi^\nabla_{/\hspace{-0.9mm}/\mbG}$  depends only on the isomorphism class  of $(\mcE, \nabla)$.
Here, denote by 
\begin{align} \label{QS3000}
[0]_U  \in (F^*_U(\Omega_{U^\mr{log}/S^\mr{log}})^\times \times^{\mbG_m} \mfc) (U)
\index{$[0]$, $[0]^{\times r}$, $[0]_U$}
\end{align}
the $U$-rational point  determined by $[0] \in \mfc (k)$.

\begin{rema}[Change of structure group]
Let  $\mbG'$, $\mfg'$, and $w$ be as in \S\,\ref{QS60}.
Suppose further that $\mbG'$ is a split  semisimple algebraic group over $k$  of adjoint type with $p > 2  h_{\mbG'}$.
The adjoint $\mbG$-action on $\mfg$ is compatible with the adjoint $\mbG'$-action on $\mfg'$ via $w$ and  its  differential $dw : \mfg \migi \mfg'$.
Hence, $d w$ induces a $k$-morphism $\mfc_\mfg \migi \mfc_{\mfg'}$ via taking the respective quotients of $\mfg$ and $\mfg'$.
This morphism 
 preserves the $\mbG_m$-actions, so
  induces a $U$-morphism
\begin{align}
(d w)_U : F^*_U(\Omega_{U^\mr{log}/S^\mr{log}})^\times \times^{\mbG_m} \mfc_{\mfg} \migi F^*_U(\Omega_{U^\mr{log}/S^\mr{log}})^\times \times^{\mbG_m} \mfc_{\mfg'}.
\end{align}
One verifies that (\ref{QS61}) deduces   the following equality:
\begin{align} \label{QH23}
{^p}\psi^{w_*(\nabla)}_{/\hspace{-0.9mm}/\mbG} =(d w)_U \circ  {^p}\psi^\nabla_{/\hspace{-0.9mm}/\mbG}.
\end{align}
In particular, ${^p}\psi^\nabla_{/\hspace{-0.9mm}/\mbG} =[0]$ implies ${^p}\psi^{w_*(\nabla)}_{/\hspace{-0.9mm}/\mbG} =[0]_U$.
\end{rema}

\bde \label{y089}
Let $\msE^\spadesuit := (\mcE_\mbB, \nabla)$ be a $(\mfg, \hslash)$-oper  on $U^\mr{log}/S^\mr{log}$.
Then,  we shall say that 
   $\msE^\spadesuit$
is {\bf dormant}  (resp., {\bf $p$-nilpotent}) 
\index{$\msE^\spadesuit$, $(\mfg, \hslash)$-oper@--- dormant} 
\index{$\msE^\spadesuit$, $(\mfg, \hslash)$-oper@--- $p$-nilpotent}   if 
${^p}\psi^\nabla =0$  (resp., ${^p}\psi^\nabla_{/\hspace{-0.9mm}/\mbG} = [0]_U$).
 \ede

\subsection{The moduli spaces $\mfO \mfp^\ZZZ_{\mfg, \hslash, \msX}$ and $\mfO \mfp^\NILP_{\mfg, \hslash, \msX}$}
Let $(g,r)$ be a pair of nonnegative integers with $2g-2+r >0$ and $\msX := (f: X \migi S, \{ \sigma_i \}_{i=1}^r)$ an $r$-pointed stable curve of genus $g$ over $S$.
Denote by 
\index{$\mfO \mfp^\ZZZ_{\mfg, \hslash, (\rho,) \msX}$, $\mfO \mfp^\ZZZ_{\mfg, \hslash, (\rho,) g,r}$, $\mfO \mfp_{n, \hslash, \bb (, \rho)}^{\diamondsuit \ZZZ}$}
\index{$\mfO \mfp^\NILP_{\mfg, \hslash, (\rho,)  \msX}$, $\mfO \mfp^\NILP_{\mfg, \hslash, (\rho,) g,r}$}
\begin{align} \label{QS239}
\mfO \mfp^\ZZZ_{\mfg, \hslash, \msX} \ \left(\text{resp.,} \  \mfO \mfp^\NILP_{\mfg, \hslash, \msX} \right)
\end{align}
 (cf. ~\cite{Perra} for the notation ``Zzz..."!) the subfunctor  of $\mfO \mfp_{\mfg, \hslash, \msX}$ classifying dormant $(\mfg, \hslash)$-opers (resp., $p$-nilpotent $(\mfg, \hslash)$-opers).
For each $\rho \in \mfc^{\times r} (S)$,  we obtain  the subfunctor 
\index{$\mfO \mfp^\ZZZ_{\mfg, \hslash, (\rho,) \msX}$, $\mfO \mfp^\ZZZ_{\mfg, \hslash, (\rho,) g,r}$, $\mfO \mfp_{n, \hslash, \bb (, \rho)}^{\diamondsuit \ZZZ}$}
\index{$\mfO \mfp^\NILP_{\mfg, \hslash, (\rho,)  \msX}$, $\mfO \mfp^\NILP_{\mfg, \hslash, (\rho,) g,r}$}
\begin{align} \label{QS240}
\mfO \mfp^\ZZZ_{\mfg, \hslash, \rho,  \msX} \ \left(\text{resp.,} \  \mfO \mfp^\NILP_{\mfg, \hslash,  \rho, \msX} \right)
\end{align}
 of $\mfO \mfp_{\mfg, \hslash, \rho, \msX}$ classifying dormant $(\mfg, \hslash)$-opers (resp., $p$-nilpotent $(\mfg, \hslash)$-opers) of radii $\rho$.
These functors 
may be represented by $S$-schemes and
we have a commutative diagram  of $S$-schemes
\begin{align} \label{dormant10}
\vcenter{\xymatrix@C=46pt@R=36pt{
\mfO \mfp^\ZZZ_{\mfg, \hslash, \rho, \msX} 
 \ar[r] \ar[d]&\mfO \mfp^\NILP_{\mfg, \hslash, \rho, \msX}   \ar[r] \ar[d]&  \mfO \mfp_{\mfg, \hslash, \rho, \msX} \ar[d]\\
\mfO \mfp^\ZZZ_{\mfg, \hslash, \msX} 
 \ar[r]&\mfO \mfp^\NILP_{\mfg, \hslash, \msX}  \ar[r]&  \mfO \mfp_{\mfg, \hslash, \msX},
}}
\end{align}
where  the arrows are all closed immersions and 
the both sides of
 squares  are  cartesian.


\bpr \label{QS241}
Suppose that  $r>0$ and let $\rho$ be an element of $\mfc^{\times r}(S)$.
Then, the following assertions hold:
\begin{itemize}
\item[(i)]
Suppose further  that $\hslash \in \Gamma (S, \mcO_S^\times)$.
Then, 
$\mfO \mfp^\ZZZ_{\mfg, \hslash, \rho,  \msX} = \emptyset$  unless $\hslash^{-1} \star \rho \in \mfc (\mbF_p)^{\times r}$. 
In particular, we have
\begin{align}
\mfO \mfp^\ZZZ_{\mfg, \hslash, \msX} = \coprod_{\varrho \in \mfc (\mbF_p)^{\times r}} \mfO \mfp^\ZZZ_{\mfg, \hslash, \hslash \star \varrho,  \msX}.
\end{align}
\item[(ii)]
Suppose further that  $\hslash =0$ and $S$ is reduced.
 Then,  $\mfO \mfp^\ZZZ_{\mfg, \hslash, \rho,  \msX} =\emptyset$ unless $\rho =  [0]^{\times r} \in \mfc^{\times r}(k) \left(\subseteq \mfc^{\times r}(S) \right)$.
 In particular,  we have
 \begin{align}
 \mfO \mfp^\ZZZ_{\mfg, 0, \msX} =  \mfO \mfp^\ZZZ_{\mfg, 0, [0]^{\times r},  \msX}. 
 \end{align}
\end{itemize}
\epr
\begin{proof}
Assertion (i) (resp., (ii)) follows directly from Proposition \ref{y077},  (i) (resp., (ii)).
\end{proof}

\bpr \label{y081} 
\begin{itemize}
\item[(i)]
Let $\msE^\spadesuit_\odot := (\mcE_{\mbB^\odot}, \nabla)$ be a $(\mfg^\odot, \hslash)$-oper on $\msX$. Denote by $\iota_{\mbG*}(\msE^\spadesuit_\odot) := (\mcE_\mbB, \iota_{\mbG*}(\nabla))$ the $(\mfg, \hslash)$-oper
associated to  $\msE^\spadesuit_\odot$ (cf.  (\ref{associated})).
Then, the equality
\begin{align}
{^p}\psi^{\iota_{\mbG*}(\nabla)} = d \iota_{\mcE_{\mbB^\odot}} \circ {^p}\psi^\nabla
\end{align}
holds (cf. (\ref{QQ0911}) for the definition of $d \iota_{\mcE_{\mbB^\odot}}$).
In particular, $\msE^\spadesuit_\odot$ is dormant if and only if $\iota_{\mbG*}(\msE^\spadesuit_\odot)$ is dormant.
\item[(ii)]
The  closed  immersion  $\iota_{\mfg}^{\mfO \mfp} : \mfO \mfp_{\mfg^\odot, \hslash, \msX} \migiincl \mfO \mfp_{\mfg, \hslash, \msX}$ (cf. (\ref{QR066}))
restricts to  a closed immersion 
\index{$\iota_{\mfg}^{\mcO p}$, $\iota_{\mfg}^{\mfO \mfp}$, ${^p}\iota_{\mfg}^{\mfO \mfp}$}
\begin{equation} \label{QR130}
  {^p}\iota^{\mfO \mfp}_{\mfg} : \mfO \mfp^\ZZZ_{\mfg^\odot, \hslash, \msX} \migi \mfO \mfp^\ZZZ_{\mfg, \hslash, \msX},
   \end{equation}
which makes  the following commutative  diagram cartesian:
\begin{align} \label{QW9000}
\vcenter{\xymatrix@C=46pt@R=36pt{
  \mfO \mfp^\ZZZ_{\mfg^\odot, \hslash, \msX} \ar[r]^-{{^p}\iota^{\mfO \mfp}_{\mfg}}  \ar[d]_-{\mr{inclusion}}& \mfO \mfp^\ZZZ_{\mfg, \hslash, \msX}   \ar[d]^-{\mr{inclusion}} \\
  \mfO \mfp_{\mfg^\odot, \hslash, \msX}  \ar[r]_-{\iota_{\mfg}^{\mfO \mfp}}&  \mfO \mfp_{\mfg, \hslash, \msX}.
}}
\end{align}
\end{itemize}
\epr
\begin{proof}
The  former assertion of (i) follows from the fact that $d \iota_{\mcE_{\mbB^\odot}}$ preserves the  $p$-power operations.
The latter assertion of (i) follows from the injectivity of $d \iota_{\mcE_{\mbB^\odot}}$.
Assertion (ii) follows directly from assertion (i).
\end{proof}

\begin{rema}  \label{QR131}
Let us keep the notation in Remarks \ref{QQ892} and \ref{QQ1004} (under the assumption $p > 2  h_\mbG$).
Then, 
it follows from (\ref{QS68}) that (\ref{product}) restricts to an $S$-isomorphism
\begin{align} \label{QS190}
\mfO \mfp_{\mfg \oplus \mfg', \hslash, \msX}^\ZZZ \isom \mfO \mfp_{\mfg, \hslash, \msX}^\ZZZ \times_S \mfO \mfp_{\mfg',\hslash, \msX}^\ZZZ.
\end{align}
Moreover, if $r>0$, then for  $\rho \in \mfc_\mfg (\mbF_p)^{\times r}$ and $\rho' \in \mfc_{\mfg'}(\mbF_p)^{\times r}$, (\ref{productrad}) restricts to an $S$-isomorphism
\begin{align} \label{QS3004}
\mfO \mfp_{\mfg \oplus \mfg', \hslash, (\rho, \rho'), \msX}^\ZZZ \isom \mfO \mfp_{\mfg, \hslash, \rho, \msX}^\ZZZ \times_S  \mfO \mfp_{\mfg', \hslash, \rho', \msX}^\ZZZ. 
\end{align}
The same assertions hold for the moduli space of $p$-nilpotent $(\mfg \oplus \mfg', \hslash)$-opers (of radii $(\rho,\rho')$).

\end{rema}

\section{The Hitchin-Mochizuki morphism} \label{QS230}

\subsection{} \label{QS1299}
Let us fix  a collection $\ \qq := \{ x_\alpha \}_{\alpha}$
\index{$\ \qq\,$}
 of generators of $\mfg^\alpha$'s 
  (for simple positive roots $\alpha$ in $\mbB$) as in (\ref{QG5697}).
This collection  associates elements $q_1$ (cf. (\ref{sss01}))  and $q_{-1}$ (cf. (\ref{q_{-1}})), and hence we obtain an isomorphism
  $\mr{Kos}_\mfg : \mfg^{\mr{ad} (q_1)} \isom \mfc$ (cf. (\ref{Kosg})).
We shall denote by
\begin{align} \label{QS3008}
\Psi_\hslash : \mfc_S \migi \mfc_S
\index{$\Psi_\hslash$}
\end{align}
 the $S$-endomorphism of $\mfc_S := S \times_k \mfc$  given by assigning
\begin{equation} \label{QS355}
\lambda \mapsto \Psi_\hslash (\lambda) := \chi \circ ((q_{-1} + \mr{Kos}_\mfg^{-1} (\lambda))^{[p]} - \hslash^{p-1} \cdot  (q_{-1} +\mr{Kos}_\mfg^{-1}(\lambda)))
\end{equation}
for any $\lambda \in \mfc_S(S')$, where $S'$ is an $S$-scheme.
If $r>0$, then  for each 
 $\rho := (\rho_i)_{i=1}^r  \in \mfc^{\times r} (S) \left(= \mfc_S(S)^{\times r} \right)$, we shall   write 
\begin{equation} \label{Frob6}
\Psi_\hslash  (\rho) :=   (\Psi_\hslash (\rho_i))_{i =1}^r \in \mfc^{\times r} (S).
\end{equation}

\begin{rema} \label{QG3021}
Let 
$\msE^\spadesuit := (\mcE_\mbB, \nabla)$ be  a $(\mfg, \hslash)$-oper on $\msX$.
Then, for every $i =1, \cdots, r$,  we have  
\begin{align} \label{QG3020}
\chi (\sigma_i^*({^p}\psi^\nabla)) = \Psi_\hslash (\rho^\nabla_i).
\end{align}
  In particular,  $\Psi_\hslash (\rho_i^\nabla) = [0]$  if $\msE^\spadesuit$ is $p$-nilpotent. 
To prove  (\ref{QG3020}), we may assume that $\msE^\spadesuit$ is $\ \qq$-normal (cf. Propositions \ref{y049} and  \ref{QS82}).
Under this assumption, the monodromy $\mu^\nabla_i$ of $\nabla$ satisfies   $\sigma_i^*(\mu^\nabla_i) = q_{-1} + \mr{Kos}^{-1}_\mfg (\rho_i^\nabla)$.
Hence,   (\ref{pcurvM}) deduces  (\ref{QG3020}).

\end{rema}

\begin{exa} \label{QS388}
If $\hslash =0$, then  it follows from Proposition \ref{y076} that
 the equality
\begin{equation} \label{rho0}
\Psi_0 (\lambda) = F_{\mfc/k} \circ \lambda
\end{equation}
holds for $\lambda \in \mfc_S (S')$, where $S'$ is any $S$-scheme. 
That is to say,   $\Psi_{0}$ coincides with the relative Frobenius morphism $F_{\mfc_S/S} : \mfc_S \migi \mfc_S \left(= \mfc_S^{(1)} \right)$ of $\mfc_S$
 over $S$.
\end{exa}

\subsection{A grading on the coordinate ring of $\mfc_S$} \label{QS409}
The $\mcO_S$-algebra $\mcO_S \langle\mfc_S\rangle$ corresponding to the affine  $S$-scheme $\mfc_S$ has a $\mbZ$-grading $\mcO_S \langle\mfc_S\rangle = \bigoplus_{j \in \mbZ} \mcO_S \langle\mfc_S\rangle_j$
 induced by  the $\mbG_m$-action on $\mfc$ opposite to  the action $\star$.
 This grading  is positive,  meaning  that   $\mcO_S \langle\mfc_S\rangle_j = 0$ if $j <0$.
 By letting 
 \begin{align} \label{QS3033}
 \mcO_S \langle\mfc_S\rangle^j := \bigoplus_{j' \leq j} \mcO_S \langle\mfc_S\rangle_{j'},
 \end{align}
  we obtain an increasing  filtration
 $\{ \mcO_S \langle\mfc_S\rangle^j \}_{j \in \mbZ_{\geq 0}}$ on $\mcO_S \langle\mfc_S\rangle$.


\bpr \label{y087} 
\begin{itemize}
\item[(i)]
The $\mcO_S$-algebra endomorphism  $\Psi_{\hslash}^*$ of $\mcO_S \langle\mfc_S\rangle$ corresponding to $\Psi_{\hslash}$
satisfies that
\begin{align}
\Psi_{\hslash}^*(\mcO_S \langle\mfc_S\rangle^j) \subseteq  \mcO_S \langle\mfc_S\rangle^{p\cdot j}
\end{align}
 for each $j \in \mbZ_{\geq 0}$.
\item[(ii)]
For each $j \in \mbZ_{\geq 0}$,  we shall denote by 
$\mr{gr}^j(\Psi_{\hslash}^*)$ the $\mcO_S$-linear  morphism  
\begin{align}
\mcO_S \langle\mfc_S\rangle_j \left(= \mcO_S \langle\mfc_S\rangle^j /\mcO_S \langle\mfc_S\rangle^{j-1} \right) \migi \mcO_S \langle\mfc_S\rangle_{p\cdot j} \left(=   \mcO_S \langle\mfc_S\rangle^{p \cdot j} /\mcO_S \langle\mfc_S\rangle^{p \cdot j-1}\right)
\end{align}
induced by $\Psi_\hslash^* |_{ \mcO_S \langle\mfc_S\rangle^j}$.
Then, the $\mcO_S$-algebra endomorphism
\begin{align}
\mr{gr}(\Psi_{\hslash}^*) := \bigoplus_{j \in \mbZ_{\geq 0}} \mr{gr}^j(\Psi_{\hslash}^*) : \mcO_S \langle\mfc_S\rangle \migi \left( \bigoplus_{j \in \mbZ_{\geq 0}}  \mcO_S \langle\mfc_S\rangle_{p \cdot j}\subseteq\right) \mcO_S \langle\mfc_S\rangle
\end{align}
 of $\mcO_S \langle\mfc_S\rangle$ coincides with the 
$\mcO_S$-algebra endomorphism $F^*_{\mfc_S/S} : \mcO_S \langle\mfc_S\rangle \migi \mcO_S \langle\mfc_S\rangle$ induced by $F_{\mfc_S/S}$.
\item[(iii)]
 $\Psi_{\hslash}$ is finite and faithfully flat of degree $p^{\mr{rk}(\mfg)}$.
\end{itemize}
\epr
\begin{proof}
To prove the assertions, it suffices to consider the case where  $S = \mbA^1 := \mr{Spec}(k[\hat{\hslash}])$ and $\hslash = \hat{\hslash}$.
Indeed,   a general case can be deduced from this case  by   base-change via  the morphism $\hslash : S \migi \mbA^1$.
Also, we may assume, without loss of generality,  that $k$ is algebraically closed.

First, we shall  consider assertion (i).
Let $\alpha : \mbG \times_k \mfg^{\mr{ad}(q_1)} \migi \mfg$ be the morphism given by $(h, v) \mapsto \mr{Ad}_\mbG  (h)(q_{-1} + v)$.
Also, let $\Psi_{\mfg, \hat{\hslash}}$ be the $\mbA^1$-endomorphism of 
$\mbA^1 \times_k \mfg$ given by $v \mapsto v^{[p]} - \hat{\hslash}^{p-1} \cdot v$.
Then, for each $(h, v) \in  \mbG \times_k \mfg^{\mr{ad}(q_1)}$, we have
\begin{align}
& \ \ \  \ \chi \circ \Psi_{\mfg, \hat{\hslash}} \circ \alpha ((h, v))   \\
&= 
\chi \circ \Psi_{\mfg, \hat{\hslash}} (\mr{Ad}_\mbG (h) (q_{-1}+ v)) \notag  \\
&= \chi ((\mr{Ad}_\mbG (h) (q_{-1}+ v))^{[p]}- \hat{\hslash}^{p-1} \cdot \mr{Ad}_\mbG (h) (q_{-1}+ v)) \notag  \\
& = \chi (\mr{Ad}_\mbG (h) ((q_{-1}+v)^{[p]}-\hat{\hslash}^{p-1} \cdot (q_{-1}+v))) \notag \\
& =   \chi ((q_{-1}+v)^{[p]}-\hat{\hslash}^{p-1} \cdot(q_{-1}+v))
\notag \\
& =   \chi ((q_{-1}+\mr{Kos}^{-1}_\mfg (\chi (v)))^{[p]}-\hat{\hslash}^{p-1} \cdot(q_{-1} + \mr{Kos}_\mfg^{-1}(\chi (v))))
\notag \\
& = \Psi_{\hat{\hslash}} \circ \chi (v)  \notag \\
&= \Psi_{\hat{\hslash}} \circ \chi (\mr{Ad}_\mbG (h) (v)) \notag \\
&= \Psi_{\hat{\hslash}} \circ \chi \circ \alpha ((h, v)). \notag
\end{align}
 This implies  $\chi \circ \Psi_{\mfg, \hat{\hslash}} \circ \alpha = \Psi_{\hat{\hslash}} \circ \chi \circ \alpha$.
 The elements  $q_{-1}+v$ for various $v \in \mfg^{\mr{ad}(q_1)}$ are regular (cf. ~\cite{Ngo},   
 the comment preceding Lemma 1.2.3), 
 so the morphism $\alpha$ is dominant.
 Hence, we have $\chi \circ \Psi_{\mfg, \hat{\hslash}}  = \Psi_{\hat{\hslash}} \circ \chi$.
 If 
  \begin{align}
 \Psi^*_{\mfg, \hat{\hslash}} : k[\hat{\hslash}] \otimes_k \mbS_k (\mfg^\vee) \migi k[\hat{\hslash}] \otimes_k \mbS_k (\mfg^\vee)
 \end{align}
denotes   the $k[\hat{\hslash}]$-algebra endomorphism  corresponding to $\Psi_{\mfg, \hat{\hslash}}$, then the equality $\chi \circ \Psi_{\mfg, \hat{\hslash}}  = \Psi_{\hat{\hslash}} \circ \chi$ induces 
  a commutative square diagram
\begin{align} \label{fff0108gg}
\vcenter{\xymatrix@C=50pt@R=36pt{
\hspace{0mm}\Gamma (\mbA^1, \mcO_{\mbA^1} \langle \mfc_{\mbA^1}\rangle) \left(= k[\hat{\hslash}]\otimes_k\mbS_k (\mfg^\vee)^\mbG\right) \ar[r]^-{\mr{inclusion}} \ar[d]_-{\Psi^*_{\hat{\hslash}}} & k[\hat{\hslash}] \otimes_k \mbS_k (\mfg^\vee)  \ar[d]^-{\Psi^*_{\mfg, \hat{\hslash}}}\\
\hspace{0mm} \Gamma (\mbA^1, \mcO_{\mbA^1} \langle \mfc_{\mbA^1}\rangle) \left(=  k[\hat{\hslash}]\otimes_k\mbS_k (\mfg^\vee)^\mbG\right)\ar[r]_-{\mr{inclusion}}& k[\hat{\hslash}] \otimes_k \mbS_k (\mfg^\vee).  
}} \hspace{0mm}
\end{align}
 By  the definition of $\Psi_{\mfg, \hat{\hslash}}$,  we obtain the inclusion relation 
\begin{align} \label{QS400}
\Psi^*_{\mfg, \hat{\hslash}} (k[\hat{\hslash}]\otimes_k\mbS_k^{\leq j} (\mfg^\vee)) \subseteq k[\hat{\hslash}]\otimes_k\mbS_k^{\leq p\cdot j} (\mfg^\vee)
\end{align}
for each $j \in \mbZ_{\geq 0}$, where $\mbS_k^{\leq j'} (\mfg^\vee)$ ($j' \in \mbZ_{\geq 0}$) denotes the $k$-subspace of $\mbS_k (\mfg^\vee)$ consisting of elements of degree $\leq j'$.
 Since the inclusion  
\begin{align} \label{QS404}
\Gamma (\mbA^1, \mcO_{\mbA^1} \langle \mfc_{\mbA^1}\rangle) \left(=k[\hat{\hslash}]\otimes_k\mbS_k (\mfg^\vee)^\mbG\right) \migiincl k[\hat{\hslash}]\otimes_k\mbS_k (\mfg^\vee)
\end{align}
 preserves the   gradings,
  (\ref{QS400}) and the commutativity of (\ref{fff0108gg}) together 
 imply that 
 \begin{align}
 \Psi^*_{\hat{\hslash}}(\Gamma (\mbA^1, \mcO_{\mbA^1} \langle \mfc_{\mbA^1}\rangle^j)) \subseteq \Gamma (\mbA^1, \mcO_{\mbA^1} \langle \mfc_{\mbA^1}\rangle^{p \cdot j}).
 \end{align}
This completes the proof of assertion (i).

Next,  we shall consider assertion (ii).
Denote by $\Psi^*_{\mfg, 0}$ the $k[\hat{\hslash}]$-algebra endomorphism of $k [\hat{\hslash}] \otimes_k \mbS_k (\mfg^\vee)$ corresponding to
the endomorphism $\mbA^1 \times_k \mfg \migi \mbA^1 \times_k \mfg$  given   by
the $p$-power operation $(-)^{[p]}$ on $\mfg$.
Also, let  $\mr{gr} (\Psi^*_{\mfg, \hat{\hslash}})$ and $\mr{gr}(\Psi^*_{\mfg, 0})$ denote the endomorphisms of 
$k[\hat{\hslash}] \otimes_k \mbS_k (\mfg^\vee) \left(= \mr{gr}(k[\hat{\hslash}]\otimes_k \mbS_k (\mfg^\vee)) \right)$ induced,  via taking gradings,  from $\Psi^*_{\mfg, \hat{\hslash}}$ and $\Psi^*_{\mfg, 0}$ respectively.
By  definition,  the equality $\mr{gr} (\Psi^*_{\mfg, \hat{\hslash}})= \mr{gr}(\Psi^*_{\mfg, 0})$ holds.
Hence, the commutativity of (\ref{fff0108gg}) implies the commutativity  of the following square diagram:
\begin{align} \label{QS401}
\vcenter{\xymatrix@C=46pt@R=36pt{
\Gamma (\mbA^1, \mcO_{\mbA^1}\langle \mfc_{\mbA^1} \rangle)
 \ar[r]^-{\mr{inclusion}} \ar[d]_-{\mr{gr}(\Psi^*_{\hat{\hslash}})}& k [\hat{\hslash}] \otimes_k \mbS_k (\mfg^\vee)  \ar[d]^-{\mr{gr}(\Psi^*_{\mfg, 0})}\\
\Gamma (\mbA^1, \mcO_{\mbA^1}\langle \mfc_{\mbA^1} \rangle) 
 \ar[r]_-{\mr{inclusion}}& k[\hat{\hslash}] \otimes_k \mbS_k (\mfg^\vee).
}}
\end{align}
On the other hand,
 by Proposition \ref{y076},
the inclusion    (\ref{QS404}) is compatible with the endomorphisms $F^*_{\mfc_{\mbA^1}/\mbA^1} \left(= \mr{gr}(F^*_{\mfc_{\mbA^1}/\mbA^1}) \right)$, $\mr{gr} (\Psi^*_{\mfg, 0})$.
Thus,  we have  $\mr{gr} (\Psi^*_{\hat{\hslash}}) = F^*_{\mfc_{\mbA^1}/\mbA^1}$,  as desired.

Finally,  we shall consider assertion (iii).
By assertion (ii),  
$\Gamma (\mbA^1, \mcO_{\mbA^1} \langle \mfc_{\mbA^1} \rangle)$
 turns out to be  a finite $\mr{gr}(\Psi_{\hat{\hslash}}^*) (\Gamma (\mbA^1, \mcO_{\mbA^1} \langle \mfc_{\mbA^1} \rangle))$-module.
 It follows from a routine argument that
 $\Gamma (\mbA^1, \mcO_{\mbA^1} \langle \mfc_{\mbA^1} \rangle)$ is a finite  $\Psi_{\hat{\hslash}}^*(\Gamma (\mbA^1, \mcO_{\mbA^1} \langle \mfc_{\mbA^1} \rangle))$-module, i.e., 
 $\Psi_{\hat{\hslash}}$ is finite.
 Since the domain and codomain of $\Psi_{\hat{\hslash}}$  are irreducible, smooth, and of the same dimension, one verifies immediately that  $\Psi_{\hat{\hslash}}$ is also faithfully flat.
 The computation of its degree may be accomplished by computing the degree of its fiber  $\Psi_{\hat{\hslash}} |_{\hat{\hslash} =0}$  over the point ``$\hat{\hslash} =0$''  in $\mbA^1$.
 But, $\Psi_{\hat{\hslash}} |_{\hat{\hslash} =0}$ coincides with $F_{\mfc/k}$  (cf. Example \ref{QS388}), so 
 its degree is given by $\left(p^{\mr{dim}(\mfc)} =\right) p^{\mr{rk}(\mfg)}$. 
This completes the proof of assertion (iii).
\end{proof}

\subsection{The spindles $\bigoplus$\hspace{-4.6mm}$\bigotimes_{\msX}$ and   $\bigoplus$\hspace{-4.6mm}$\bigotimes_{\msX}^\nabla$}\label{y083}

Denote by
\index{$\ooalign{$\bigoplus$ \cr $\bigotimes_{\msX (, \rho)}$}$, $\ooalign{$\bigoplus$ \cr $\bigotimes_{\msX (, \rho)}^\nabla$}$}
\begin{equation}  \label{spindle}
\ooalign{$\bigoplus$ \cr $\bigotimes_{\msX}$}
 \  \left(\text{resp.,} \  \ooalign{$\bigoplus$ \cr $\bigotimes_{\msX}^\nabla$} \right)
 \end{equation}
(cf. ~\cite{Perra} for the notation ``$\ooalign{$\bigoplus$ \cr $\bigotimes$}$\hspace{-1mm}" that represents a {\it spindle})
 the $\mfS \mfe \mft$-valued contravariant functor on $\mfS \mfc \mfh_{/S}$  which, to any $S$-scheme $s : S' \migi S$, assigns the set of  morphisms  $S' \times_S  X\migi (F_X^*(\Omega_{X^\mr{log}/S^\mr{log}}))^\times \times^{\mbG_m}\mfc$  (resp., 
$S' \times_S X  \migi \Omega_{X^{\mr{log}}/S^\mr{log}}^\times \times^{\mbG_m}\mfc$) over $X$. 
Here, for each line bundle $\mcL$, we denote by $\mcL^\times$ the $\mbG_m$-torsor corresponding to $\mcL$.
One verifies that both $\ooalign{$\bigoplus$ \cr $\bigotimes_{\msX}$}$ and 
$\ooalign{$\bigoplus$ \cr $\bigotimes_{\msX}^\nabla$}$
may be represented by a finite type $S$-scheme
(cf. (\ref{fff0105}) in the proof of Proposition \ref{y085}).
There exists a canonical identification  $\ooalign{$\bigoplus$ \cr $\bigotimes_{\msX^{(1)}}^\nabla$} = (\ooalign{$\bigoplus$ \cr $\bigotimes_{\msX}^\nabla$})^{(1)}$ between $S$-schemes.

Next, let us suppose that $r >0$.
Given an $S$-scheme $s : S' \migi S$  and  an element $w$ of  $\ooalign{$\bigoplus$ \cr $\bigotimes_{\msX}$} (S')$ (resp., $\ooalign{$\bigoplus$ \cr $\bigotimes^\nabla_{\msX}$} (S')$), we obtain, for each $i \in \{1, \cdots, r \}$,  
the composite $S$-morphism
\begin{align} \label{QS1018} 
w |_{\sigma_i} :  S'
\isom S' \times_S S
 \xrightarrow{ \sigma^*_i (w)} \sigma^*_{i}(\Omega)^\times \times^{\mbG_m} \mfc 
\isom \mcO_S^\times \times^{\mbG_m} \mfc 
\isom \mfc_S \left(= S \times_k \mfc \right),  
\end{align}
where $\Omega : = F_X^*(\Omega_{X^{\mr{log}}/S^\mr{log}})$ (resp., $\Omega := \Omega_{X^{\mr{log}}/S^\mr{log}}$) and  the third  arrow  arises from 
the isomorphism $\sigma^*_i(\Omega) \isom \mcO_S$ induced by  $\mr{Res}_{\msX, i}$ (cf. (\ref{QQ040})).
The assignment $w \mapsto (w |_{\sigma_i})_{i=1}^r$ defines an $S$-morphism 
\index{$\mr{Rad}^{\ooalign{$\oplus$ \cr $\otimes$}}$, $\mr{Rad}^{\ooalign{$\oplus$ \cr $\otimes$}\hspace{-1mm}, \nabla}$}
\begin{equation} \label{spindle441}
\mr{Rad}^{\ooalign{$\oplus$ \cr $\otimes$}}  \hspace{-1mm} : \ooalign{$\bigoplus$ \cr $\bigotimes_{\msX}$} \migi \mfc_S^{\times r} \ \left(\text{resp.,} \ \mr{Rad}^{\ooalign{$\oplus$ \cr $\otimes$}\hspace{-1mm}, \nabla}  : \ooalign{$\bigoplus$ \cr $\bigotimes^\nabla_{\msX}$} \migi \mfc_S^{\times r}  \right),
\end{equation}
where $\mfc_S^{\times r}$ denotes the product of $r$ copies of $\mfc_S$ over $S$.


\subsection{Sleeping Beauty} \label{QS321}

Denote by  
\index{$\rotatebox[origin=c]{58}{{\large $\kappa$}}^{\hspace{-0.0mm}}_{\hslash (, \rho)}$, $\rotatebox[origin=c]{58}{{\large $\kappa$}}^{\hspace{-0.0mm}\nabla}_{\hslash (, \rho)}$, $\rotatebox[origin=c]{58}{{\large $\kappa$}}^{\hspace{-0.0mm}\nabla}_{\hslash (, \rho) \mr{univ}}$, Hitchin-Mochizuki morphism}
\begin{equation}  \label{dormant8} 
\rotatebox[origin=c]{58}{{\large $\kappa$}}_{\hslash} :  \mfO \mfp_{\mfg, \hslash, \msX} \migi  \ooalign{$\bigoplus$ \cr $\bigotimes_{\msX}$} 
  \end{equation}
(cf. ~\cite{Perra} for the notation ``$\rotatebox[origin=c]{58}{{\large $\kappa$}}$" that represents an {\it old fairy}) 
 the $S$-morphism   which, to each $(\mfg, \hslash)$-oper $\msE^\spadesuit := (\mcE_\mbB, \nabla)$, assigns ${^p}\psi^\nabla_{/ \hspace{-0.9mm}/\mbG}$. 
The morphism
\begin{align} \label{QS116}
\text{\Pisymbol{pzd}{42}}_{\hspace{-0.5mm}\msX}  :S \migi \ooalign{$\bigoplus$ \cr $\bigotimes_{\msX}$}
\index{$\text{\Pisymbol{pzd}{42}}_{\hspace{-0.5mm}\msX}$}
\end{align}
(cf. ~\cite{Perra} for the notation ``$\text{\Pisymbol{pzd}{42}}$" that represents  a {\it princess's finger})
determined by  $[0]_X : X \migi (F_X^*(\Omega_{X^\mr{log}/S^\mr{log}}))^\times \times^{\mbG_m} \mfc$ (cf. (\ref{QS3000})) makes the following square diagram cartesian:
 \begin{align} \label{QR112}
\hspace{0mm}\vcenter{\xymatrix@C=46pt@R=36pt{
\hspace{-20mm}\left(\mfO \mfp^\ZZZ_{\mfg, \hslash,  \msX} \subseteq \right)\mfO \mfp^\NILP_{\mfg, \hslash,   \msX} \ar[r]^{\mr{inclusion}} \ar[d]& \mfO \mfp_{\mfg, \hslash,   \msX}  \ar[d]^-{\rotatebox[origin=c]{58}{{\large $\kappa$}}_{\hslash}}\\
S \ar[r]_-{\text{\Pisymbol{pzd}{42}}_{\hspace{-0.5mm}\msX}}& \ooalign{$\bigoplus$ \cr $\bigotimes_{\msX}$}.
}}
\end{align}

\ble \label{QS250}
Suppose that $r>0$, and
let $\msE^\spadesuit := (\mcE_\mbB, \nabla)$ be a $(\mfg, \hslash)$-oper on $\msX$.
Then,  for each $i =1, \cdots, r$,   we have the following  equality  of  elements in $\mfc (S)$:
 \begin{align} \label{QS268}
 {^p}\psi^\nabla_{/\hspace{-0.9mm}/\mbG} |_{\sigma_i} = \Psi_{\hslash}(\rho_i^\nabla).
 \end{align}
(cf. (\ref{QS3008}) for the definition of $\Psi_{\hslash}$).
\ele
\begin{proof}
The formations of both  sides of (\ref{QS268}) depend only on the isomorphism class of $\msE^\spadesuit$ and the problem is of local nature,
Hence, after possibly replacing $S$ with its  covering in the Zariski  topology,  one may assume 
that there exists  a  log chart $(U, \partial)$ on $X^\mr{log}/S^\mr{log}$ with $\mr{Im}(\sigma_i) \subseteq U$  and  $\msE^\spadesuit$ is  $\ \qq$-normal (cf. Proposition \ref{y049}).
 If
 \begin{align}
 \tau_i : \mfg (S) \left(=\Gamma (S, \mcO_S \otimes_k \mfg) \right) \isom \Gamma (S, \mfg_{\sigma^*_i (\DE_\mbG)})
 \end{align}
   denotes the isomorphism   induced naturally by $\mr{triv}_{\mbG, (U, \partial)} : U \times_k \mbG \isom \DE_\mbG |_U$ (cf. (\ref{canonicalTg}) and  (\ref{QR0257})), then
 the equality $\mu_i^\nabla = \tau_i(q_{-1} + \mr{Kos}_\mfg^{-1}(\rho_i^\nabla))$ holds.
 Hence, we have
 \begin{align}
&  \ \ \  \  {^p}\psi^\nabla_{/\hspace{-0.9mm}/\mbG} |_{\sigma_i}  \\
&= [\chi] ( (\sigma^*_i (\DE_\mbG), \sigma_i^*({^p}\psi^\nabla)) ) \notag \\
 & =   [\chi] ( (\sigma^*_i (\DE_\mbG), (\mu^\nabla_i)^{[p]}-\hslash^{p-1} \cdot \mu_i^\nabla) ) \notag \\
 & = [\chi] ( (\sigma^*_i (\DE_\mbG), \tau_i(q_{-1} + \mr{Kos}_\mfg^{-1}(\rho_i^\nabla))^{[p]}-\hslash^{p-1} \cdot \tau_i(q_{-1} + \mr{Kos}_\mfg^{-1}(\rho_i^\nabla))) ) \notag \\
 & =  [\chi] ( (\sigma^*_i (\DE_\mbG), \tau_i ((q_{-1} + \mr{Kos}_\mfg^{-1}(\rho_i^\nabla))^{[p]}-\hslash^{p-1} \cdot (q_{-1} + \mr{Kos}_\mfg^{-1}(\rho_i^\nabla))) )) \notag \\
 & = [\chi ] ((S \times_k \mbG,  (q_{-1} + \mr{Kos}_\mfg^{-1}(\rho_i^\nabla))^{[p]}-\hslash^{p-1} \cdot (q_{-1} + \mr{Kos}_\mfg^{-1}(\rho_i^\nabla)))) \notag \\
 & = \chi \circ  ((q_{-1} + \mr{Kos}_\mfg^{-1}(\rho_i^\nabla))^{[p]}-\hslash^{p-1} \cdot (q_{-1} + \mr{Kos}_\mfg^{-1}(\rho_i^\nabla))) \notag \\
 & = \Psi_\hslash (\rho_i^\nabla), \notag
 \end{align}
where the  second equality follows from (\ref{pcurvM}).
This completes the proof of the assertion.
\end{proof}

The above lemma implies the commutativity of 
the following diagram:
\begin{align} \label{QS1020}
\vcenter{\xymatrix@C=46pt@R=36pt{
\mfO \mfp_{\mfg, \hslash, \msX}  \ar[r]^-{\rotatebox[origin=c]{58}{{\large $\kappa$}}_{\hslash}} \ar[d]_-{\mr{Rad}_\hslash}& \ooalign{$\bigoplus$ \cr $\bigotimes_{\msX}$}  \ar[d]^-{\mr{Rad}^{\tiny{\ooalign{$\oplus$ \cr $\otimes$}}}}\\
 \mfc_S^{\times r} \ar[r]_-{\Psi_\hslash^{\times r}}& \mfc_S^{\times r},
}}
\end{align}
where $\Psi_\hslash^{\times r}$ denotes the $r$-th power of $\Psi_\hslash$ and the both sides of vertical arrows are taken to be  the natural morphisms onto $S$. 
(If $r =0$, then we set $\mfc_S^{\times 0} := S$  and $\Psi_\hslash^{\times r} := \mr{id}_S$.) 
Given an element $\rho \in \mfc^{\times r} (S) = \mfc_S^{\times r} (S)$, where $\mfc^{\times 0}(S) := \{ \emptyset \}$ and  $\rho := \emptyset$ if $r =0$, we shall write
\index{$\ooalign{$\bigoplus$ \cr $\bigotimes_{\msX (, \rho)}$}$, $\ooalign{$\bigoplus$ \cr $\bigotimes_{\msX (, \rho)}^\nabla$}$}
\begin{equation}  \label{spindle31}\ooalign{$\bigoplus$ \cr $\bigotimes_{\rho, \msX}$} := (\mr{Rad}^{\ooalign{$\oplus$ \cr $\otimes$}}\hspace{-1mm})^{-1} (\rho) \ \left(\text{resp.,} \ \ooalign{$\bigoplus$ \cr $\bigotimes^\nabla_{\rho, \msX}$}  := (\mr{Rad}^{\ooalign{$\oplus$ \cr $\otimes$}\hspace{-1mm}, \nabla})^{-1} (\rho)\right),
\end{equation}
i.e.,  the closed subscheme of  $\ooalign{$\bigoplus$ \cr $\bigotimes_{\msX}$}$ (resp., $\ooalign{$\bigoplus$ \cr $\bigotimes^\nabla_{\msX}$}$)  defined to be the scheme-theretic inverse image of $\rho$ via $\mr{Rad}^{\tiny{\ooalign{$\oplus$ \cr $\otimes$}}}$\hspace{-1mm} (resp., $\mr{Rad}^{\ooalign{$\oplus$ \cr $\otimes$}\hspace{-1mm}, \nabla}$).
By the commutativity of (\ref{QS1020}),  the morphism $\rotatebox[origin=c]{58}{{\large $\kappa$}}_{\hslash}$ restricts to an $S$-morphism
\index{$\rotatebox[origin=c]{58}{{\large $\kappa$}}^{\hspace{-0.0mm}}_{\hslash (, \rho)}$, $\rotatebox[origin=c]{58}{{\large $\kappa$}}^{\hspace{-0.0mm}\nabla}_{\hslash (, \rho)}$, $\rotatebox[origin=c]{58}{{\large $\kappa$}}^{\hspace{-0.0mm}\nabla}_{\hslash (, \rho) \mr{univ}}$, Hitchin-Mochizuki morphism} 
\begin{align} \label{QR11HH0}
\rotatebox[origin=c]{58}{{\large $\kappa$}}_{\hslash, \rho} :  \mfO \mfp_{\mfg, \hslash, \rho, \msX} \migi  \ooalign{$\bigoplus$ \cr $\bigotimes_{\Psi_\hslash (\rho), \msX}$}.
\end{align}


\subsection{The Hitchin-Mochizuki morphism} \label{QS200}
Given an element $\rho$ of $\mfc^{\times r}(S)$,
we shall  denote by 
\begin{equation} \label{spindle18}
F_{X/S}^{\ooalign{$\oplus$ \cr $\otimes$}}
 :   {\ooalign{$\bigoplus$ \cr $\bigotimes^\nabla_{\msX^{(1)}}$}}  \migi
{\ooalign{$\bigoplus$ \cr $\bigotimes_{\msX}$}} \ \ \ 
\left(\text{resp.}, F_{X/S, \rho}^{\ooalign{$\oplus$ \cr $\otimes$}} : {\ooalign{$\bigoplus$ \cr $\bigotimes^\nabla_{\rho, \msX^{(1)}}$}}  \migi {\ooalign{$\bigoplus$ \cr $\bigotimes_{\rho, \msX}$}}\right)
 \index{$F_{X/S (, \rho)}^{\ooalign{$\oplus$ \cr $\otimes$}}$}
  \end{equation}
  the $S$-morphism given
by pull-back via $F_{X/S}$.  

\bpr \label{y085}
The morphism $F_{X/S}^{\ooalign{$\oplus$ \cr $\otimes$}}$ (resp., $F_{X/S, \rho}^{\ooalign{$\oplus$ \cr $\otimes$}}$)  is a closed immersion.
 \epr
\begin{proof}
Since $F_{X/S, \rho}^{\ooalign{$\oplus$ \cr $\otimes$}}$ may be obtained as 
the restriction of  $F_{X/S}^{\ooalign{$\oplus$ \cr $\otimes$}}$,
 it suffices to prove the non-resp'd assertion.

Let us fix an isomorphism  $\mfc \isom  \mr{Spec}(k[u_1, \cdots, u_\mr{rk}(\mfg)])$ as in (\ref{isomCA}).
If $\mcL$ is a line bundle  on either $X$ or $X^{(1)}$, 
then the isomorphism fixed right now  allows us to identify
$\mcL^\times \times^{\mbG_m} \mfc$
with
  the relative affine space  associated to $\bigoplus_{l=1}^{\mr{rk}(\mfg)} \mcL^{\otimes e_l}$,
  i.e., we have $\mcL^\times \times^{\mbG_m} \mfc \cong \mbV (\bigoplus_{l=1}^{\mr{rk}(\mfg)} \mcL^{\otimes e_l})$.
 Hence, there exist isomorphisms of $S$-schemes
 \begin{equation} \label{fff0105}
 {\ooalign{$\bigoplus$ \cr $\bigotimes^\nabla_{\msX^{(1)}}$}} \isom \mbV (f^{(1)}_*(\bigoplus_{l=1}^{\mr{rk}(\mfg)}  \Omega^{\otimes e_l}_{X^{(1) \mr{log}}/S^\mr{log}})), \hspace{3mm}
  {\ooalign{$\bigoplus$ \cr $\bigotimes_{\msX}$}} \isom \mbV (f_*(\bigoplus_{l=1}^{\mr{rk}(\mfg)} F^*_X(\Omega_{X^\mr{log}/S^\mr{log}}^{\otimes e_l}))).
 \end{equation}
  For $l = 1, \cdots, \mr{rk}(\mfg)$, we shall write
\begin{equation}   \label{QS3056}
\omega_l :  \Omega^{\otimes e_l}_{X^{(1) \mr{log}}/S^\mr{log}} \migi F_{X/S*}(F_{X/S}^*(\Omega^{\otimes e_l}_{X^{(1) \mr{log}}/S^\mr{log}})) \ \left(\cong F_{X/S *}(F_X^*(\Omega_{X^\mr{log}/S^\mr{log}}^{\otimes  e_l}))\right)\end{equation}
for the 
morphism of $\mcO_{X^{(1)}}$-modules corresponding, 
via  the adjunction relation 
 ``$F^*_{X/S} (-) \dashv F_{X/S*}(-)$",
   to 
 the identity morphism of $F^*_{X/S}(\Omega^{\otimes e_l}_{X^{(1)\mr{log}}/S^\mr{log}})$.
The $\mcO_{X^{(1)}}$-module $F_{X/S *}(F_{X}^*(\Omega_{X^\mr{log}/S^\mr{log}}^{\otimes  e_l}))$ is relatively torsion-free of rank $p$ (cf. Proposition \ref{torsionE}) and $\omega_l$ is easily verified to be injective over the scheme-theoretically dense open subscheme  $X^{(1)\mr{sm}}$  of $X^{(1)}\setminus \mr{Supp}(D_\msX)$ (cf. \S\,\ref{y030}).
It follows that the morphism  $\omega_l$, as well as its direct image $f^{(1)}(\omega_l)$, is injective.
Since $e_l \geq 2$, 
the following equalities hold:
\begin{align}
& \mbR^1f^{(1)}_*(\Omega^{e_l}_{X^{(1)\mr{log}}/S^\mr{log}}) = 0,  \\
&  \mbR^1f^{(1)}_*(F_{X/S *}(F_{X}^*(\Omega_{X^{\mr{log}}/S^\mr{log}}^{\otimes  e_l}))) \ \left(= \mbR^1f_*(\Omega_{X^\mr{log}/S^\mr{log}}^{\otimes p \cdot e_l})\right)=0.   \notag
\end{align}
Hence, both the domain and codomain of  $f_*^{(1)}(\omega_l)$ are locally free (cf. ~\cite{Har}, Chap. III, Theorem 12.11 (b)),  and
the formation of  $f_*^{(1)}(\omega_l)$ is
 functorial with respect to  $S$.
By considering  
the above argument with $S$ replaced by  various $S$-schemes,
 one concludes that $f^{(1)}_*(\omega_l)$ is {\it universally injective}
with respect to base-change to $S$-schemes. 
In particular,   $\mr{Coker}(f^{(1)}_*(\omega_l))$ turns out to be flat over $S$ (cf. ~\cite{MAT}, p.\,17, Theorem 1), and 
 the dual 
\begin{equation}
f_*^{(1)}(\omega_l)^\vee : f^{(1)}_*(F_{X/S *}(F_X^*(\Omega_{X^\mr{log}/S^\mr{log}}^{\otimes  e_l})))^\vee \migi f^{(1)}_*(\Omega^{\otimes e_l}_{X^{(1) \mr{log}}/S^\mr{log}})^\vee
\end{equation}
of $f_*^{(1)}(\omega_l)$  is surjective for all $l$.
Since 
the morphism  $F_{X/S}^{\ooalign{$\oplus$ \cr $\otimes$}}$  is obtained by (taking the spectrum of)  the surjection
\begin{equation}
 \mbS_{\mcO_S}(\bigoplus_{l =1}^{\mr{rk}(\mfg)}f^{(1)}_*(F_{X/S *}(F_X^*(\Omega_{X^\mr{log}/S^\mr{log}}^{\otimes  e_l})))^\vee) \migisurj \mbS_{\mcO_S}(\bigoplus_{l =1}^{\mr{rk}(\mfg)}f^{(1)}_*(\Omega^{\otimes e_l}_{X^{(1) \mr{log}}/S^\mr{log}})^\vee)
\end{equation}
 induced by $\bigoplus_{l =1}^{\mr{rk}(\mfg)}f_*^{(1)}(\omega_l)^\vee$,
we see that  $F_{X/S}^{\ooalign{$\oplus$ \cr $\otimes$}}$ defines  a closed immersion.
This completes the proof of the assertion.
\end{proof}

\begin{prdef}\label{y088}
Let  $\rho$ be an element of $\mfc^{\times r} (S)$.
Then, the morphism $\rotatebox[origin=c]{58}{{\large $\kappa$}}_\hslash$ (resp., $\rotatebox[origin=c]{58}{{\large $\kappa$}}_{\hslash, \rho}$)
 factors through the closed immersion 
 $F^{\ooalign{$\oplus$ \cr $\otimes$}}_{X/S}$ (resp., $F^{\ooalign{$\oplus$ \cr $\otimes$}}_{X/S, \Psi_\hslash (\rho)}$).
 That is to say, 
 there exists a unique  morphism
\index{$\rotatebox[origin=c]{58}{{\large $\kappa$}}^{\hspace{-0.0mm}}_{\hslash (, \rho)}$, $\rotatebox[origin=c]{58}{{\large $\kappa$}}^{\hspace{-0.0mm}\nabla}_{\hslash (, \rho)}$, $\rotatebox[origin=c]{58}{{\large $\kappa$}}^{\hspace{-0.0mm}\nabla}_{\hslash (, \rho) \mr{univ}}$, Hitchin-Mochizuki morphism}
\begin{equation} \label{QG2500}
 \rotatebox[origin=c]{58}{{\large $\kappa$}}^{\hspace{-0.0mm}\nabla}_\hslash :  \mfO \mfp_{\mfg, \hslash, \msX} \migi  \ooalign{$\bigoplus$ \cr $\bigotimes^\nabla_{\msX^{(1)}}$}  \   
 \left(\text{resp.}, \ \rotatebox[origin=c]{58}{{\large $\kappa$}}^{\hspace{-0.0mm}\nabla}_{\hslash, \rho} :  \mfO \mfp_{\mfg, \hslash, \rho, \msX} \migi  \ooalign{$\bigoplus$ \cr $\bigotimes^\nabla_{\Psi_\hslash (\rho), \msX^{(1)}}$}\right) 
 \end{equation}
 satisfying  
 \begin{equation} \label{QHH11}
 F^{\ooalign{$\oplus$ \cr $\otimes$}}_{X/S}
  \circ
   \rotatebox[origin=c]{58}{{\large $\kappa$}}^{\hspace{-0.0mm}\nabla}_{\hslash} =  \rotatebox[origin=c]{58}{{\large $\kappa$}}_{\hslash} \ \left(\text{resp.},  \ 
   F^{\ooalign{$\oplus$ \cr $\otimes$}}_{X/S, \Psi_\hslash (\rho)}
    \circ \rotatebox[origin=c]{58}{{\large $\kappa$}}^{\hspace{-0.0mm}\nabla}_{\hslash, \rho} =  \rotatebox[origin=c]{58}{{\large $\kappa$}}_{\hslash, \rho}\right).
 \end{equation}
 We shall refer to $ \rotatebox[origin=c]{58}{{\large $\kappa$}}^{\hspace{-0.0mm}\nabla}_\hslash$ (resp., $\rotatebox[origin=c]{58}{{\large $\kappa$}}^{\hspace{-0.0mm}\nabla}_{\hslash, \rho}$) as the {\bf Hitchin-Mochizuki morphism}.
 \end{prdef}
\begin{proof}
It suffices to consider 
the non-resp'd assertion because the resp'd assertion is entirely similar.
In what follows, let 
us keep the notation in
the proof of Proposition \ref{y085}.
 
Let $\msE^\spadesuit = (\mcE_\mbB, \nabla_\mcE)$  be  a $(\mfg, \hslash)$-oper on $\msX$, which determines an $S$-rational point $s_{\msE^\spadesuit}$ of $\mfO \mfp_{\mfg, \hslash, \msX}$;
its image $\rotatebox[origin=c]{58}{{\large $\kappa$}}_\hslash (s_{\msE^\spadesuit})$
 corresponds, via the second  isomorphism in (\ref{fff0105}), an element  
\begin{equation}
R_{\msE^\spadesuit} \in \Gamma (X^{(1)}, \bigoplus_{l=1}^{\mr{rk}(\mfg)} F_{X/S*}(F_X^*(\Omega_{X^{\mr{log}}/S^{\mr{log}}}^{\otimes e_l}))).
\end{equation}
By  an argument similar to the argument discussed in the proof  of ~\cite{CZ}, Proposition  3.1 (or ~\cite{LP},  Proposition 3.2),
we see that the restriction $R_{\msE^\spadesuit} |_{X^{\mr{sm}}}$ of  $R_{\msE^\spadesuit}$ to the scheme-theoretically dense open subscheme $X^{\mr{sm}}$ of $X\setminus \mr{Supp}(D_\msX)$ lies   in the image of 
$\bigoplus_{l =1}^{\mr{rk} (\mfg)} w_l  |_{X^{\mr{sm}}}$.
Since 
the cokernel of $\bigoplus_l \omega_l$ is  relatively torsion-free (cf. Proposition \ref{torsion}), 
 $R_{\msE^\spadesuit}$ itself lies in   the image of $\bigoplus_{l =1}^{\mr{rk} (\mfg)} w_l$.
By considering  this fact with  $\msE^\spadesuit$ replaced with 
the base-change to 
every   $S$-scheme,
we conclude  that $\rotatebox[origin=c]{58}{{\large $\kappa$}}_\hslash$ factors through  $F^{\ooalign{$\oplus$ \cr $\otimes$}}_{X/S}$.
This completes the proof of the assertion.
\end{proof}

According to   Proposition-Definition \ref{y088}, the commutative diagram  (\ref{QS1020}) induces
a commutative square diagram
\begin{align} \label{FGHJ}
\vcenter{\xymatrix@C=46pt@R=36pt{
  \mfO \mfp_{\mfg, \hslash, \msX} 
  \ar[r]^-{\rotatebox[origin=c]{58}{{\large $\kappa$}}^{\hspace{-0.0mm}\nabla}_{\hslash}} \ar[d]_-{\mr{Rad}_\hslash}& \ooalign{$\bigoplus$ \cr $\bigotimes_{\msX^{(1)}}^\nabla$}  
   \ar[d]^-{\mr{Rad}^{\ooalign{$\oplus$ \cr $\otimes$}\hspace{-1mm}, \nabla}}
  \\
\mfc_S^{\times r}\ar[r]_-{\Psi_\hslash^{\times r}}& \mfc_S^{\times r}.
}}
\end{align}
In particular, 
since  both $\Psi_\hslash^{\times r}$ and $\mr{Rad}_\hslash$ are surjective (cf. Proposition \ref{y087}, (iii), and Theorem \ref{QQ1060}), 
we see that
$\mr{Rad}^{\ooalign{$\oplus$ \cr $\otimes$}\hspace{-1mm}, \nabla}$
  is surjective.

\section{Varying the parameter $\hslash$} \label{y090}

We shall discuss dormant and $p$-nilpotent $(\mfg, \hslash)$-opers with  universal parameter $\hslash$.
After that, we describe the degeneration of the Hitchin-Mochizuki morphism  $\rotatebox[origin=c]{58}{{\large $\kappa$}}_\hslash$ at $\hslash =0$.

\subsection{The $\mbG_m$-action on $\mfO \mfp^\ZZZ_{\mfg, \hat{\hslash}, \hat{\msX}}$}
 \label{QS1101}
Let us keep the notation  in \S\,\ref{y059} (e.g., $\hat{\hslash}$, $\hat{\msX}$, and $\Spec$, etc.).
Also, let us take an element $\rho \in \mfc^{\times r} (S)$.
It follows from  (\ref{QR101}) that  the $\mbG_m$-action on $\mfO \mfp_{\mfg, \hat{\hslash}, \hat{\msX}}$ (resp., $\mfO \mfp_{\mfg, \hat{\hslash}, \hat{\rho}, \hat{\msX}}$) defined in (\ref{action61}) (cf. (\ref{QR0506}))  restricts to    a $\mbG_m$-action
\begin{align} \label{QR151}
\mbG_m \times \mfO \mfp_{\mfg, \hat{\hslash},  \hat{\msX}}^\ZZZ \migi  \mfO \mfp_{\mfg, \hat{\hslash}, \hat{ \msX}}^\ZZZ \ 
\left(\text{resp.,} \ \mbG_m \times \mfO \mfp_{\mfg, \hat{\hslash}, \hat{\rho},  \hat{\msX}}^\ZZZ \migi  \mfO \mfp_{\mfg, \hat{\hslash}, \hat{\rho}, \hat{ \msX}}^\ZZZ \right)
\end{align}
on $\mfO \mfp_{\mfg, \hat{\hslash},  \hat{\msX}}^\ZZZ$ (resp., $\mfO \mfp_{\mfg, \hat{\hslash},  \hat{\rho},  \hat{\msX}}^\ZZZ$).
If $\hslash$  is an element of $\Gamma (S, \mcO_S^\times)$, then
this morphism restricts to an isomorphism
\begin{align} \label{QS3201}
\mbG_m \times_k \mfO \mfp^\ZZZ_{\mfg, \hslash, \msX} \isom 
\Spec^{-1}((\mbG_m)_S)) \cap  \mfO \mfp_{\mfg, \hat{\hslash},  \hat{\msX}}^\ZZZ
\hspace{9mm}
\\
\left(\text{resp.,} \  \mbG_m \times_k \mfO \mfp^\ZZZ_{\mfg, \hslash, \hslash \star \rho, \msX} \isom 
\Spec^{-1}_{\rho} ((\mbG_m)_S)) \cap  \mfO \mfp_{\mfg, \hat{\hslash}, \hat{\rho},  \hat{\msX}}^\ZZZ
  \right). \notag
\end{align}
In particular,  the assignment $\msE^\spadesuit \mapsto \hslash \cdot \msE^\spadesuit$  determines an isomorphism
\begin{align} \label{QS1100}
\mfO \mfp^\ZZZ_{\mfg, 1, \msX} \isom \mfO \mfp_{\mfg, \hslash, \msX}^\ZZZ \ 
\left(\text{resp.,} \ \mfO \mfp^\ZZZ_{\mfg, 1, \rho, \msX} \isom \mfO \mfp_{\mfg, \hslash, \hslash \star \rho, \msX}^\ZZZ\right)
\end{align}
over $S$.
(\ref{QS3201}) and (\ref{QS1100}) are compatible with (\ref{QR0588}) and (\ref{QS1002}) respecitvely.


\subsection{The $\mbG_m$-action on $\mfO \mfp^\NILP_{\mfg, \hat{\hslash}, \hat{\msX}}$} 
\label{QS1102}
On the other hand, the $\mbG_m$-action on $\mfc$ defined by $(c, \varrho) \mapsto c^p \star \varrho$ (for any $c \in \mbG_m$, $\varrho \in \mfc$) induces  a $\mbG_m$-action
\begin{align} \label{action098}
\mbG_m \times_k \ooalign{$\bigoplus$ \cr $\bigotimes_{\hat{\msX}}$}   \migi  \ooalign{$\bigoplus$ \cr $\bigotimes_{\hat{\msX}}$}.
\end{align}
It follows from (\ref{QR101}) that the following square diagram  is commutative:
\begin{align} \label{QR144}
\vcenter{\xymatrix@C=46pt@R=36pt{
\mbG_m \times_k \mfO \mfp_{\mfg, \hat{\hslash}, \hat{\msX}}\ar[r]^-{(\ref{QR0588})} \ar[d]_-{\mr{id}_{\mbG_m} \times \rotatebox[origin=c]{58}{{\large $\kappa$}}_{\hat{\hslash}}}& \mfO \mfp_{\mfg, \hat{\hslash}, \hat{\msX}} \ar[d]^-{\rotatebox[origin=c]{58}{{\large $\kappa$}}_{\hat{\hslash}}} \\
\mbG_m \times_k \ooalign{$\bigoplus$ \cr $\bigotimes_{\hat{\msX}}$} \ar[r]_-{(\ref{action098})}& \ooalign{$\bigoplus$ \cr $\bigotimes_{\hat{\msX}}$}.
}}
\end{align}
Since the inverse image of 
$\mr{Im}(\text{\Pisymbol{pzd}{42}}_{\hspace{-0.5mm}\hat{\msX}})$ 
 via (\ref{action098}) coincides with $\mbG_m \times_k \mr{Im}(\text{\Pisymbol{pzd}{42}}_{\hspace{-0.5mm}\hat{\msX}})$, the above diagram induces
a $\mbG_m$-action
\begin{align} \label{QR150}
\mbG_m \times \mfO \mfp_{\mfg, \hat{\hslash},  \hat{\msX}}^\NILP \migi  \mfO \mfp_{\mfg, \hat{\hslash}, \hat{ \msX}}^\NILP \
\left(\text{resp.,} \  \mbG_m \times \mfO \mfp_{\mfg, \hat{\hslash}, \hat{\rho},  \hat{\msX}}^\NILP \migi  \mfO \mfp_{\mfg, \hat{\hslash}, \hat{\rho}, \hat{ \msX}}^\NILP\right).
\end{align}
 on $\mfO \mfp_{\mfg, \hat{\hslash},   \hat{\msX}}^\NILP$ (resp., $\mfO \mfp_{\mfg, \hat{\hslash}, \hat{\rho},  \hat{\msX}}^\NILP$) compatible with  that on $\mfO \mfp_{\mfg, \hat{\hslash},  \hat{\msX}}^\ZZZ$ (resp., $\mfO \mfp_{\mfg, \hat{\hslash}, \hat{\rho}, \hat{\msX}}^\ZZZ$).
Moreover,  if $\hslash \in \Gamma (S, \mcO_S^\times)$, then  (\ref{QR0588}) restricts to 
\begin{align} \label{QS1007}
\mbG_m \times_k \mfO \mfp^\NILP_{\mfg, \hslash, \hslash \star \rho, \msX} \isom 
\Spec^{-1} ((\mbG_m)_S)) \cap  \mfO \mfp_{\mfg, \hat{\hslash},  \hat{\msX}}^\NILP  \hspace{4mm}\\
\left(\text{resp.,} \ \mbG_m \times_k \mfO \mfp^\NILP_{\mfg, \hslash, \hslash \star \rho, \msX} \isom \Spec^{-1}_{\rho} ((\mbG_m)_S)) \cap  \mfO \mfp_{\mfg,  \hat{\hslash},  \hat{\rho}, \hat{\msX}}^\NILP   \right) \notag
\end{align}
which is compatible with (\ref{QS3201}).
In particular, the assignment $\msE^\spadesuit \mapsto \hslash \cdot \msE^\spadesuit$  determines an $S$-isomorphism
\begin{align} \label{QS1104}
\mfO \mfp^\NILP_{\mfg, 1, \rho, \msX} \isom \mfO \mfp_{\mfg, \hslash, \hslash \star \rho, \msX}^\NILP.
\end{align}


\subsection{Degeneration  at $\hslash =0$} \label{y092}

We specialize to 
the case of  $\hslash = 0$.
Since
the isomorphism $\mr{Kos}_\mfg : \mfg^{\mr{ad}(p_1)} \isom\mfc$ (cf. (\ref{Kosg})) is $\mbG_m$-equivariant  with respect to  the $\mbG_m$-action  on 
$\mfg^{\mr{ad}(p_1)}(1)$ and 
 the $\mbG_m$-action $\star$ on $\mfc$, it  induces  an isomorphism
\begin{equation} \label{KosgVV}
\mbV (f_*(\DV_{\mfg,  \msX})) \isom  \ooalign{$\bigoplus$ \cr $\bigotimes_{\msX}^\nabla$}
\end{equation}
of $S$-schemes  (cf.  the proof of Proposition \ref{y046}).
 In particular,  $\ooalign{$\bigoplus$ \cr $\bigotimes_{\msX}^\nabla$}$ 
is a relative affine spaces over $S$ of relative dimension  ${\aleph}(\mfg)$ (cf. Proposition \ref{y046}, (ii)).
For  an element  $\rho := (\rho_i)_{i=1}^r \in \mfc^{\times r} (S)$ (where $\mfc^{\times r}(S) := \{ \emptyset \}$ and  $\rho := \emptyset$ if $r=0$), 
denote by 
\begin{align} \label{QS199}
\mbV (f_*(\DV_{\mfg,  \msX}))_\rho
\end{align}
the closed subscheme of $\mbV (f_*(\DV_{\mfg,  \msX}))$ corresponding to $\ooalign{$\bigoplus$ \cr $\bigotimes_{\rho, \msX}^\nabla$} \left(\subseteq  \ooalign{$\bigoplus$ \cr $\bigotimes_{\msX}^\nabla$} \right)$ via (\ref{KosgVV}).

\bpr \label{y093} 
\begin{itemize}
\item[(i)]
The square diagram
\begin{align} \label{QS160}
\vcenter{\xymatrix@C=46pt@R=36pt{
 \mbV (f_*(\DV_{\mfg, \msX}))\ar[r]_-{\sim}^-{(\ref{QR0304})} \ar[d]^-{\wr}_-{(\ref{KosgVV})} & \mfO \mfp_{\mfg, 0, \msX} \ar[d]^-{\rotatebox[origin=c]{58}{{\large $\kappa$}}^{\hspace{-0.0mm}\nabla}_{\hspace{-0.0mm}0}}\\
 \ooalign{$\bigoplus$ \cr $\bigotimes^\nabla_{\msX}$}  \ar[r]_-{F_{\tiny{\ooalign{$\oplus$ \cr $\otimes^\nabla_{\msX}$}}/S}}&  \ooalign{$\bigoplus$ \cr $\bigotimes^\nabla_{\msX^{(1)}}$} \left(=  ({\ooalign{$\bigoplus$ \cr $\bigotimes^\nabla_{\msX}$}})^{(1)} \right)
}}
\end{align}
is commutative.
In particular, $\rotatebox[origin=c]{58}{{\large $\kappa$}}_{0}^{\hspace{-0.0mm}\nabla}$  is finite and faithfully flat of degree $p^{\aleph(\mfg)}$. 
\item[(ii)]
Let $\rho$ be an element of $\rho \in \mfc^{\times r} (S)$.
Then,  (\ref{QR0304}) restricts to an $S$-isomorphism
\begin{align} \label{QS3442}
\mbV (f_*(\DV_{\mfg, \msX}))_\rho \isom  
\mfO \mfp_{\mfg, 0, \rho, \msX},
\end{align}
and (\ref{QS160}) restricts to a  commutative diagram
\begin{align} \label{QS163}
\vcenter{\xymatrix@C=46pt@R=36pt{
 \mbV (f_*(\DV_{\mfg, \msX}))_\rho\ar[r]_-{\sim}^-{(\ref{QS3442})} \ar[d]^-{\wr} & \mfO \mfp_{\mfg, 0, \rho, \msX} \ar[d]^-{\rotatebox[origin=c]{58}{{\large $\kappa$}}^{\hspace{0mm}\nabla}_{\hspace{-0.0mm}0, \rho}}\\
 \ooalign{$\bigoplus$ \cr $\bigotimes^\nabla_{\rho, \msX}$}  \ar[r]_-{F_{{\tiny\ooalign{$\oplus$ \cr $\otimes^\nabla_{\rho, \msX}$}}/S}}&  \ooalign{$\bigoplus$ \cr $\bigotimes^\nabla_{\Psi_\hslash (\rho),  \msX^{(1)}}$} \left(=  ({\ooalign{$\bigoplus$ \cr $\bigotimes^\nabla_{\rho, \msX}$}})^{(1)} \right),
}}
\end{align}
where the left-hand vertical arrow denotes the restriction of (\ref{KosgVV}) to $\mbV (f_*(\DV_{\mfg, \msX}))_\rho$.
In particular,  $\rotatebox[origin=c]{58}{{\large $\kappa$}}^{\hspace{-0.0mm}\nabla}_{\hspace{0mm} 0, \rho} $ is finite and faithfully flat of degree $p^{{^c}\aleph(\mfg)}$. 
 \end{itemize}
 \epr
\begin{proof}
Assertion (i) follows from  the definitions  of the various morphisms in (\ref{QS160}) and Proposition \ref{y076}.
Assertion (ii) follows from assertion (i) together with the result of Example \ref{QS388}.
\end{proof}

\bco  \label{y094}
Suppose that  we are given $\rho \in \mfc^{\times r} (S)$ with  $\Psi_\hslash (\rho) = [0]^{\times r}$. 
Then, $\mfO \mfp_{\mfg, 0, \msX}^\NILP$
(resp., $\mfO \mfp_{\mfg, 0, \rho, \msX}^\NILP$)
is finite and faithfully flat over $S$ of degree $p^{{\aleph} (\mfg)}$ (resp., $p^{{^c\aleph} (\mfg)}$).
Moreover, $\mfO \mfp_{\mfg, 0, \msX}^\NILP$ (resp., $\mfO \mfp_{\mfg, 0, \rho, \msX}^\NILP$) is, Zariski locally on $S$,  isomorphic to 
the $S$-scheme 
\begin{align} \label{QS300ghj}
& \hspace{12mm}S \times_k \mr{Spec}(k[\varepsilon_1, \cdots, \varepsilon_{\aleph (\mfg)}]/(\varepsilon^p_1, \cdots, \varepsilon^p_{\aleph (\mfg)}))
 \\
&\left(\text{resp.,} \  S \times_k \mr{Spec}(k[\varepsilon_1, \cdots, \varepsilon_{{^c}\aleph (\mfg)}]/(\varepsilon^p_1, \cdots, \varepsilon^p_{{^c}\aleph (\mfg)})) \right) \notag
\end{align}
(cf. the following Remark for a more detailed description of this $S$-scheme).
\eco

 \begin{rema} \label{QS1200}
For a vector bundle $\mcV$ on $S$, we write
$\mbJ_p (\mcV)$ for the relative affine scheme over $S$ determined by the $\mcO_S$-algebra $\mbS_{\mcO_S} (\mcV^\vee)/\langle\mbS^p_{\mcO_S} (\mcV^\vee)\rangle$, i.e., the quotient of the symmetric $\mcO_S$-algebra $\mbS_{\mcO_S} (\mcV^\vee)$  by the ideal sheaf generated by   the homogenous local sections of degree $p$.
Then,  by taking account of  (\ref{QR020}),   (\ref{QR112}),  and Proposition  \ref{y093}, we  can describe $\mfO \mfp^\NILP_{\mfg, 0, \msX}$ (resp.,  $\mfO \mfp^\NILP_{\mfg, 0, \rho, \msX}$)  as 
 \begin{align} \label{QS1201}
 \mfO \mfp^\NILP_{\mfg, 0, \msX} \cong \mbJ_p (f_*(\DV_{\mfg}))  \cong  \prod_{j =1}^{\mr{rk}(\mfg)} \mbJ_p (f_*(\Omega_{X^\mr{log}/S^\mr{log}}^{\otimes (j+1)}))^{\times \mr{dim}(\mfg_i^{\mr{ad}(q_1)})} \hspace{18mm}\\
 \left(\text{resp.,} \   \mfO \mfp^\NILP_{\mfg, 0, \rho, \msX} \cong \mbJ_p (f_*(\DV_{\mfg}(-D_\msX)))  \cong  \prod_{j =1}^{\mr{rk}(\mfg)} \mbJ_p (f_*(\Omega_{X^\mr{log}/S^\mr{log}}^{\otimes (j+1)}(-D_\msX)))^{\times \mr{dim}(\mfg_i^{\mr{ad}(q_1)})}  \right), \notag
 \end{align}
 where ``$\prod$'' means  product  over $S$ and $(-)^{\times d}$ means $d$-th power  product over $S$.
 \end{rema}
 
 \bco \label{QS3332}
 $\mfO \mfp^\ZZZ_{\mfg, 0, \msX}$ is finite over $S$.
 Moreover, if 
 $S$ is reduced, then
  we have 
 \begin{align} \label{QS3333}
 (\mfO \mfp^\ZZZ_{\mfg, 0, \msX})_{\mr{red}} \cong S,
 \end{align}
 where  the left-hand side denotes the reduced scheme associated to $\mfO \mfp_{\mfg, 0, \msX}$.
 \eco

\section{The finiteness of the Hitchin-Mochizuki morphism} \label{y095}

The goal of this section is to 
obtain a detailed understanding of  the morphism $\rotatebox[origin=c]{58}{{\large $\kappa$}}^{\hspace{-0.0mm}\nabla}_{\hspace{0mm} \hslash, (\rho)}$ for an arbitrary $\hslash$.
As a corollary, 
 the finiteness of 
$\mfO \mfp_{\mfg, \hslash, \msX}^\ZZZ$ will be proved.

\subsection{The filtration on $\mcO_S \langle \bigoplus$\hspace{-3.9mm}$\bigotimes^\nabla_{\msX} \rangle$} \label{QS3334}

As constructed  in (\ref{QS504}),
there exists a canonical filtration $\{ \mcO_S \langle \mfO \mfp_{\mfg, \hslash, \msX} \rangle^j \}_{j \geq 0}$ on the $\mcO_S$-algebra $\mcO_S \langle \mfO \mfp_{\mfg, \hslash, \msX} \rangle$.
On the other hand,
the  $\mbG_m$-action   on $\mfc$  opposite to ``$\star$'' (cf. (\ref{action}))
 induces  a  $\mbZ_{\geq 0}$-grading 
 $\mcO_S \langle \ooalign{$\bigoplus$ \cr $\bigotimes^\nabla_{\msX}$} \rangle = \bigoplus_{j \in \mbZ_{\geq 0}} \mcO_S \langle \ooalign{$\bigoplus$ \cr $\bigotimes^\nabla_{\msX}$} \rangle_j$
 on 
  $\mcO_S \langle \ooalign{$\bigoplus$ \cr $\bigotimes^\nabla_{\msX}$} \rangle$. By letting
\begin{align}
\mcO_S \langle \ooalign{$\bigoplus$ \cr $\bigotimes^\nabla_{\msX}$} \rangle^j := \bigoplus_{j' \leq j} \mcO_S \langle \ooalign{$\bigoplus$ \cr $\bigotimes^\nabla_{\msX}$} \rangle_{j'},
\end{align}
we obtain an increasing filtration $\{ \mcO_S \langle \ooalign{$\bigoplus$ \cr $\bigotimes^\nabla_{\msX}$} \rangle^j \}_{j \in \mbZ_{\geq 0}}$.
Let    us denote by
\begin{align} \label{QS4019}
\rotatebox[origin=c]{58}{{\large $\kappa$}}^{\nabla*}_{\hslash} : 
\mcO_S \langle \ooalign{$\bigoplus$ \cr $\bigotimes^\nabla_{\msX^{(1)}}$} \rangle \migi \mcO_S\langle\mfO \mfp_{\mfg, \hslash, \msX}\rangle.
\end{align}
the $\mcO_S$-algebra morphism induced by $\rotatebox[origin=c]{58}{{\large $\kappa$}}^{\nabla}_{\hslash}$ .
Then, we obtain the following assertion:

 \ble \label{QS4012}%
 For any  $j \in \mbZ_{\geq 0}$, the inclusion relation
 \begin{equation} \label{QS943}
\rotatebox[origin=c]{58}{{\large $\kappa$}}^{\nabla*}_{\hslash} (\mcO_S \langle \ooalign{$\bigoplus$ \cr $\bigotimes^\nabla_{\msX^{(1)}}$} \rangle^j) \subseteq \mcO_S\langle\mfO \mfp_{\mfg, \hslash, \msX}\rangle^{p \cdot j}
 \end{equation}
 is fulfilled.
 \ele
 \begin{proof}
 To prove this lemma,  one may assume, after possibly changing the base $S$ to its Zariski open covering,
    that there exists a $(\mfg, \hslash)$-oper on $\msX$.
 This implies that there exists an isomorphism $\mbV (f_{*}(\DV_{\mfg, \msX})) \isom \mfO \mfp_{\mfg, \hslash, \msX}$ over $S$.
 To be precise, 
 if  we fix a $\, \qq$-normal  $(\mfg, \hslash)$-oper  $\msE^\spadesuit := (\DE_\mbB, \nabla)$  on $\msX$  (cf. Proposition \ref{y049}), 
 then  the assignment $R \mapsto \msE^\spadesuit_{+ R} :=  (\DE_\mbB, \nabla + R)$ (cf. (\ref{QR036}))  for each $R \in f_{*}(\DV_{\mfg, \msX})$ determines an isomorphism $\tau_\msX : \mbV (f_* (\DV_{\mfg, \msX})) \isom \mfO \mfp_{\mfg, \hslash, \msX}$.
The $\mcO_S$-algebra isomorphism 
\begin{align} \label{QS4002}
\tau_\msX^* : \mcO_S \langle \mfO \mfp_{\mfg, \hslash, \msX} \rangle \isom \mcO_S \langle \mbV (f_* (\mcV_{\mfg, \msX})) \rangle
\end{align}
corresponding to $\tau_\msX$ preserves the filtrations because of
the definition of $\{ \mcO_S \langle \mfO \mfp_{\mfg, \hslash, \msX} \rangle^j\}_{j}$.

Next, let us consider the $\mbZ_{\geq 0}$-grading on $\mcO_S \langle \mbV (f_* (F^*_X(\Omega_{X^\mr{log}/S^\mr{log}}) \otimes \mfg_{\DE_\mbG})) \rangle$ induced by  the $\mbG_m$-action  on $\mfg$ opposite to the action coming from the  homotheties on $\mfg$;
 this grading defines in the natural manner an increasing filtration 
 \begin{align}
 \{ \mcO_S \langle\mbV (f_* (F^*_X(\Omega_{X^\mr{log}/S^\mr{log}})  \otimes \mfg_{\DE_\mbG})) \rangle^j \}_{j \in \mbZ_{\geq 0}}.
 \end{align}
 Denote by 
 $\psi_\msX : \mbV (f_*(\DV_{\mfg, \msX})) \migi \mbV (f_*(F^*_X(\Omega_{X^\mr{log}/S^\mr{log}}) \otimes\mfg_{\DE_\mbG}))$ the $S$-morphism
 given by $R \mapsto {^p}\psi^{\nabla + R}$.
By  the definition of $p$-curvature,
the $\mcO_S$-algebra morphism
\begin{align} \label{QS4001}
\psi^*_\msX :  \mcO_S \langle \mbV (f_*(F^*_X(\Omega_{X^\mr{log}/S^\mr{log}}) \otimes\mfg_{\DE_\mbG}))\rangle
\migi \mcO_S \langle \mbV (f_*(\DV_{\mfg, \msX}))\rangle
\end{align}
corresponding to $\psi_\msX$ satisfies that
\begin{align} \label{QS4000}
\psi^*_\msX (\mcO_S \langle \mbV (f_*(F^*_X(\Omega_{X^\mr{log}/S^\mr{log}}) \otimes\mfg_{\DE_\mbG}))\rangle^j) \subseteq
\mcO_S \langle \mbV (f_*(\DV_{\mfg, \msX}))\rangle^{p \cdot j}
\end{align}
for every $j \in \mbZ_{\geq 0}$.

Moreover, 
let $\chi_\msX :  \mbV (f_* (F^*_X(\Omega_{X^\mr{log}/S^\mr{log}}) \otimes \mfg_{\DE_\mbG})) \migi \ooalign{$\bigoplus$ \cr $\bigotimes_{\msX}$}$ denote the $S$-morphism induced by $\chi : \mfg \migi \mfc$.
The morphism  $\chi$ is $\mbG_m$-equivariant,  so
 the $\mcO_S$-algebra morphism 
 \begin{align} \label{QS4004}
 \chi_\msX^* : \mcO_S \langle \ooalign{$\bigoplus$ \cr $\bigotimes_{\msX}$} \rangle \migi \mcO_S \langle \mbV (f_* (F^*_X(\Omega_{X^\mr{log}/S^\mr{log}}) \otimes \mfg_{\DE_\mbG})) \rangle
 \end{align}
 corresponding to $\chi_\msX$
 preserves the respective filtrations.
 
 Since 
 the composite $(\tau_\msX^*)^{-1}\circ \psi^*_\msX\circ \chi_\msX^*$ coincides with
 the $\mcO_S$-algebra morphism  $\rotatebox[origin=c]{58}{{\large $\kappa$}}^{*}_{\hslash} :\mcO_S \langle \ooalign{$\bigoplus$ \cr $\bigotimes_{\msX}$} \rangle \migi \mcO_S\langle\mfO \mfp_{\mfg, \hslash, \msX}\rangle$ induced by $\rotatebox[origin=c]{58}{{\large $\kappa$}}_{\hslash}$,
 the desired inclusion relation  follows from
   (\ref{QHH11}) and 
 the various compatibilities of  filtrations obtained above.
 \end{proof}
 
 \subsection{The finiteness of the Hitchin-Mochizuki morphism} \label{QS4010}
 For each $j \in \mbZ_{\geq 0}$, denote by $\mr{gr}^j(\rotatebox[origin=c]{58}{{\large $\kappa$}}^{\nabla*}_{\hslash})$ the morphism
\begin{align}
\mcO_S \langle \ooalign{$\bigoplus$ \cr $\bigotimes^\nabla_{\msX}$}  \rangle_j
\left(= \mcO_S \langle \ooalign{$\bigoplus$ \cr $\bigotimes^\nabla_{\msX}$}  \rangle^j / \mcO_S \langle \ooalign{$\bigoplus$ \cr $\bigotimes^\nabla_{\msX}$}  \rangle^{j-1} \right) \migi  \mr{gr}^{p\cdot j} (\mcO_S \langle \mfO \mfp_{\mfg, \hslash, \msX}\rangle)
\end{align}
induced by $\rotatebox[origin=c]{58}{{\large $\kappa$}}^{\nabla*}_{\hslash} |_{\mcO_S \langle \tiny{\ooalign{$\bigoplus$ \cr $\bigotimes^\nabla_{\msX}$} } \rangle^j}$ (cf.  Lemma \ref{QS4012}).
Then, 
the direct sum $\mr{gr}(\rotatebox[origin=c]{58}{{\large $\kappa$}}^{\nabla*}_{\hslash}) := \bigoplus_{j \in \mbZ_{\geq 0}} \mr{gr}^j(\rotatebox[origin=c]{58}{{\large $\kappa$}}^{\nabla*}_{\hslash})$ determines 
an $\mcO_S$-algebra morphism
\begin{align} \label{QS4014}
\mr{gr}(\rotatebox[origin=c]{58}{{\large $\kappa$}}^{\nabla*}_{\hslash}) 
:  \mcO_S \langle \ooalign{$\bigoplus$ \cr $\bigotimes^\nabla_{\msX}$}  \rangle
 \migi 
\left(\bigoplus_{j \in \mbZ_{\geq 0}}  \mr{gr}^{p \cdot j}(\mcO_S \langle \mfO \mfp_{\mfg, \hslash, \msX}\rangle)\subseteq  \right)
 \mr{gr} (\mcO_S \langle \mfO \mfp_{\mfg, \hslash, \msX}\rangle).
\end{align}
 
 
\bpr \label{y096}
Let us keep the notation in \S\,\ref{y059}.
We shall consider  the $\mcO_{\mbA^1_S}$-algebra  morphism 
\begin{align}
&\hspace{12mm}\rotatebox[origin=c]{58}{{\large $\kappa$}}^{\nabla*}_{\hat{\hslash}}   : \mcO_{\mbA^1_S} \langle \ooalign{$\bigoplus$ \cr $\bigotimes^\nabla_{\hat{\msX}^{(1)}}$} \rangle \migi \mcO_{\mbA^1_S}\langle\mfO \mfp_{\mfg, \hat{\hslash}, \hat{\msX}}\rangle
\\
 & \left(\text{resp.,} \  \rotatebox[origin=c]{58}{{\large $\kappa$}}^{\nabla*}_{\hat{0}}  : \mcO_{\mbA^1_S} \langle \ooalign{$\bigoplus$ \cr $\bigotimes^\nabla_{\hat{\msX}^{(1)}}$} \rangle \migi \mcO_{\mbA^1_S}\langle\mfO \mfp_{\mfg, 0, \hat{\msX}}\rangle \right),  \notag
\end{align}
 defined to be 
the morphism $\rotatebox[origin=c]{58}{{\large $\kappa$}}^{\nabla*}_{\hslash}$ (cf. (\ref{QS4019}))  where
the   triple ``$(S, \msX, \hslash)$'' is replaced by  $(\mbA^1_S, \hat{\msX}, \hat{\hslash})$ (resp., $(\mbA^1_S, \hat{\msX}, 0)$).
Then, the following diagram is commutative:
\begin{align} \label{QS505}
\vcenter{\xymatrix@C=26pt@R=46pt{
 \mr{gr} (\mcO_{\mbA^1_S} \langle \mfO \mfp_{\mfg, \hat{\hslash}, \hat{\msX}}\rangle) \ar[rr]_-{\sim}^-{(\ref{QS498})} &&  \mcO_{\mbA^1_S} \langle \mfO \mfp_{\mfg, 0, \hat{\msX}}\rangle  \\
& \mcO_{\mbA^1_S} \langle \ooalign{$\bigoplus$ \cr $\bigotimes^\nabla_{\hat{\msX}}$}  \rangle \ar[ul]^-{\mr{gr}(\rotatebox[origin=c]{58}{{\large $\kappa$}}^{\nabla*}_{\hat{\hslash}})} \ar[ur]_-{\rotatebox[origin=c]{58}{{\large $\kappa$}}^{\nabla*}_{\hat{0}}}.  &
}}
\end{align}
In particular,  the morphism $\mr{gr}(\rotatebox[origin=c]{58}{{\large $\kappa$}}^{\nabla*}_{\hslash})$  does not depend on the choice of $\hslash \in \Gamma (S, \mcO_S)$  and coincides with $\rotatebox[origin=c]{58}{{\large $\kappa$}}^{\nabla*}_{0}$ under the identification  $\mr{gr} (\mcO_S \langle \mfO \mfp_{\mfg, \hslash, \msX} \rangle) = \mcO_S \langle \mfO \mfp_{\mfg, 0, \msX} \rangle$ given by (\ref{QS498}).
 \epr
\begin{proof}
Let  $\psi_{\hat{\msX}, \hat{\hslash}}^*$ (resp.,  $\psi_{\hat{\msX}, \hat{0}}^*$) denote the morphism
$\psi^*_\msX$  (cf. (\ref{QS4001})) with $\msX$ and $\hslash$ replaced by  $\hat{\msX}$ and $\hat{\hslash}$ (resp., $0$) respectively.
Then, 
by considering  the definition of $p$-curvature, we see that 
 the morphism $\mr{gr}(\psi_{\hat{\msX}, \hat{\hslash}}^*)$, i.e., the graded morphism   associated to  $\psi_{\hat{\msX}, \hat{\hslash}}^*$, coincides with $\psi_{\hat{\msX}, 0}^*$.
This implies (from the discussion in the proof of Lemma \ref{QS4012}) the desired commutativity of (\ref{QS505}).
The latter assertion follows from this result by  taking the fiber over the point ``$\hat{\hslash} = \hslash$'' in $\mbA^1_S$.
\end{proof}

By using the above proposition, we obtain  the following assertion,
generalizing ~\cite{Mzk1}, Chap.\,I, Theorem 2.3, and  the first unpublished version of ~\cite{CZ}, \S\,1, Theorem 1.1.

\bco \label{QS300}
\begin{itemize}
\item[(i)]
The 
morphism 
$\rotatebox[origin=c]{58}{{\large $\kappa$}}^{\hspace{-0.0mm}\nabla}_{\hat{\hslash}} : \mfO \mfp_{\mfg, \hat{\hslash}, \hat{\msX}} \migi \mbA^1_S \times_S  \ooalign{$\bigoplus$ \cr $\bigotimes^\nabla_{\msX}$}  \left(= \ooalign{$\bigoplus$ \cr $\bigotimes^\nabla_{\hat{\msX}}$} \right)$
is finite and faithfully flat of degree $p^{\aleph (\mfg)}$.
In particular, for each $\hslash \in \Gamma (S, \mcO_S) \left(= \mbA^1 (S)\right)$,
the structure morphism $\mfO \mfp^\NILP_{\mfg, \hslash, \msX} \migi S$ of  the $S$-scheme $\mfO \mfp^\NILP_{\mfg, \hslash, \msX}$ (resp., the Hitchin-Mochizuki  morphism $\rotatebox[origin=c]{58}{{\large $\kappa$}}^{\hspace{-0.0mm}\nabla}_{\hslash} : \mfO \mfp_{\mfg, \hslash, \msX} \migi \ooalign{$\bigoplus$ \cr $\bigotimes^\nabla_{\msX}$}$)  is finite and faithfully flat  of degree  $p^{\aleph (\mfg)}$.

\item[(ii)]
Let $\rho$ be  an element of  $\mfc^{\times r} (S)$ (where $\rho := \emptyset$ if $r =0$). 
Then, the morphism $\rotatebox[origin=c]{58}{{\large $\kappa$}}^{\hspace{-0.0mm}\nabla}_{\hat{\hslash}, \rho}  : \mfO \mfp_{\mfg, \hat{\hslash}, \rho, \hat{\msX}} \migi \ooalign{$\bigoplus$ \cr $\bigotimes^\nabla_{\Psi_\hslash (\rho), \hat{\msX}^{(1)}}$}$  is finite and faithfully flat of degree $p^{^c \aleph (\mfg)}$.
In particular,
 if $\hslash$ is an element of $\Gamma (S, \mcO_S) \left(= \mbA^1 (S) \right)$ with $\Psi_\hslash (\rho) = [0]^{\times r}$, then  the projection  $\mfO \mfp^\NILP_{\mfg, \hslash, [0]^{\times r}, \msX} \migi S$ is  finite and faithfully flat   of degree $p^{^c \aleph (\mfg)}$.
\end{itemize}
\eco
\begin{proof}
We only consider assertion (ii) because the proof of assertion (i) follows from a similar (and relatively  simple) argument using Proposition \ref{y093}, (i),   instead   of (ii) in the  same  proposition.
First, let us prove the former assertion of (ii).
Both $\mfO \mfp_{\mfg, \hat{\hslash}, \rho,  \hat{\msX}}$ and $\ooalign{$\bigoplus$ \cr $\bigotimes^\nabla_{\Psi_{\hat{\hslash}}(\rho), \hat{\msX}}$}$  are flat over $S$, 
so we may assume, after taking the respective fibers, that
 $S = \mr{Spec}(k)$.
Because of  Proposition \ref{y093}, (i),   and the commutativity of  (\ref{QS505}) proved in Proposition \ref{y096},
$\mr{gr} (\mcO_{\mbA^1} \langle \mfO \mfp_{\mfg, \hat{\hslash}, \hat{\msX}}\rangle)$ is a finite
 $\mr{gr} (\mcO_{\mbA^1}\langle \ooalign{$\bigoplus$ \cr $\bigotimes^\nabla_{\hat{\msX}}$} \rangle)$-module.
 Hence, $\mcO_{\mbA^1} \langle \mfO \mfp_{\mfg, \hat{\hslash}, \hat{\msX}}\rangle$ is a finite $\mcO_{\mbA^1}\langle \ooalign{$\bigoplus$ \cr $\bigotimes^\nabla_{\hat{\msX}}$} \rangle$-module, or equivalently, 
 the morphism  $\rotatebox[origin=c]{58}{{\large $\kappa$}}^{\hspace{-0.0mm}\nabla}_{\hat{\hslash}}$ is finite. 
 Note that $\rotatebox[origin=c]{58}{{\large $\kappa$}}^{\hspace{-0.0mm}\nabla}_{\rho, \hat{\hslash}}$ may be 
 obtained from  the finite morphism $\rotatebox[origin=c]{58}{{\large $\kappa$}}^{\hspace{-0.0mm}\nabla}_{\hat{\hslash}}$ by restricting the domain and codomain to their closed schemes.
 This implies the finiteness of  $\rotatebox[origin=c]{58}{{\large $\kappa$}}^{\hspace{-0.0mm}\nabla}_{\rho, \hat{\hslash}}$.
 Moreover, 
 since the domain and codomain of  $\rotatebox[origin=c]{58}{{\large $\kappa$}}^{\hspace{-0.0mm}\nabla}_{\rho, \hat{\hslash}}$ are irreducible, smooth, and of the same dimension, one verifies immediately that   $\rotatebox[origin=c]{58}{{\large $\kappa$}}^{\hspace{-0.0mm}\nabla}_{\rho, \hat{\hslash}}$ is also faithfully flat.
Now we compute its degree.
To this end, 
it suffices to consider the fiber over $0 \in \mbA^1$, i.e., compute the degree of 
$\rotatebox[origin=c]{58}{{\large $\kappa$}}_{0, \rho}^{\hspace{-0.0mm} \nabla}$.
But,  the computation in this case has already finished  in Proposition \ref{y093}, (ii), so  we have
$\mr{deg}(\rotatebox[origin=c]{58}{{\large $\kappa$}}^{\hspace{-0.0mm}\nabla}_{\rho, \hat{\hslash}}) = \mr{deg}(\rotatebox[origin=c]{58}{{\large $\kappa$}}^{\hspace{-0.0mm}\nabla}_{\rho, 0}) = p^{{^c}\aleph (\mfg)}$.
This completes the proof of the former assertion.
The latter assertion  follows directly from the former assertion. 
\end{proof}

\bco \label{QS1010}
 For any $\hslash \in \Gamma (S, \mcO_S)$, $\mfO \mfp_{\mfg, \hslash, \msX}^\ZZZ$ is  finite  over $S$.
\eco
\begin{proof}
The assertion follows from the finiteness of $\mfO \mfp_{\mfg, \hslash, \msX}^\NILP$ proved in Corollary \ref{QS300}, (i),  and the fact that $\mfO \mfp_{\mfg, \hslash, \msX}^\ZZZ$ is a closed subscheme of $\mfO \mfp_{\mfg, \hslash, \msX}^\NILP$.
\end{proof}

\section{The universal moduli stacks of  dormant/$p$-nilpotent $(\mfg, \hslash)$-opers} \label{y099}

This final  section deals with  the   moduli stacks of dormant/$p$-nilpotent  $(\mfg, \hslash)$-opers 
on the universal family 
 of pointed stable curves.
 We  describe several facts concerning the geometric structures of these stacks  deduced from  results obtained so far.

Let $(\mbG, \mbT)$ be, as before, a split semisimple algebraic group of adjoint type over $k$.
We suppose that $p > 2  h_\mbG$.
Also, let $(g,r)$ be a pair of nonnegative integers with $2g-2+r>0$, $\hslash$ an element of $k$,  and $\rho$   an element of $\mfc^{\times r} (k)$ (where $\mfc^{\times r} (k) := \{ \emptyset \}$ and $\rho := \emptyset$ if $r=0$).

\subsection{Dormant case} \label{QS1031}
First, let us consider  the moduli stack of   dormant $(\mfg, \hslash)$-opers.
For each $\hslash \in k$, we shall
denote by 
\index{$\mfO \mfp^\ZZZ_{\mfg, \hslash, (\rho,) \msX}$, $\mfO \mfp^\ZZZ_{\mfg, \hslash, (\rho,) g,r}$, $\mfO \mfp_{n, \hslash, \bb (, \rho)}^{\diamondsuit \ZZZ}$}
\begin{align} \label{QR139}
\mfO \mfp^\ZZZ_{\mfg, \hslash, g,r} \ \left(\text{resp.,} \  \mfO \mfp^\ZZZ_{\mfg, \hslash, \rho,  g,r} \right)
\end{align}
the closed substack of $\mfO \mfp_{\mfg, \hslash, g,r}$ (resp., $\mfO \mfp_{\mfg, \hslash, \rho,  g,r}$) (cf. (\ref{univ1})) classifying dormant  $(\mfg, \hslash)$-opers.
For simplicity, we write 
$\mfO \mfp^\ZZZ_{\mfg, g,r} := \mfO \mfp^\ZZZ_{\mfg, 1, g,r}$ (resp., $ \mfO \mfp^\ZZZ_{\mfg, \rho,  g,r} :=  \mfO \mfp^\ZZZ_{\mfg, 1, \rho,  g,r}$).
The Deligne-Mumford stack $\mfO \mfp^\ZZZ_{\mfg, \hslash,  g,r}$ (resp.,  $\mfO \mfp^\ZZZ_{\mfg, \hslash, \rho, g,r}$)
  admits a  log structure
 pulled-back from  $\overline{\mfM}_{g,r}^\mr{log}$;
we denote the resulting log stack by
\begin{align}\label{AS1055}
  \mfO \mfp^{^\mr{Zzz...} \mr{log}}_{\mfg, \hslash,  g,r} \ \left(\text{resp.,} \  \mfO \mfp^{^\mr{Zzz...} \mr{log}}_{\mfg, \hslash, \rho, g,r} \right).
\end{align}
 
\bpr[cf. Theorem \ref{y016}]
 \label{y0102}
\begin{itemize}
\item[(i)]
If $\hslash \in k^\times$, then 
we have 
\begin{equation} \label{decompZzz}     \coprod_{\varrho \in  \mfc^{\times r}(\mbF_p)} \mfO \mfp^\ZZZ_{\mfg, \hslash, \hslash \star \varrho, g,r} =   \mfO \mfp^\ZZZ_{\mfg, \hslash, g,r}.
  \end{equation}
In particular,  $\mfO \mfp^\ZZZ_{\mfg, \hslash, \hslash \star \rho, g,r}$ (for each $\rho \in \mfc^{\times r}(\mbF_p)$) may be represented by  a (possibly empty)  proper Deligne-Mumford stack over $k$ which is finite over $\overline{\mfM}_{g,r}$.
\item[(ii)]
If $\hslash =0$, then we have
\begin{align}
\mfO \mfp^\ZZZ_{\mfg, 0, [0]^{\times r}, g,r} = \mfO \mfp^\ZZZ_{\mfg, 0, g,r}.
\end{align}
Moreover, 
 the natural projection
\begin{equation}
(\mfO \mfp^\ZZZ_{\mfg, 0, g,r})_\mr{red} \migi
  \overline{\mfM}_{g,r}
\end{equation}
 is an isomorphism of $k$-stacks, where $(\mfO \mfp^\ZZZ_{\mfg, 0,  g,r})_\mr{red}$ denotes the reduced stack associated to $\mfO \mfp^\ZZZ_{\mfg, 0, g,r}$.
In particular,   $\mfO \mfp^\ZZZ_{\mfg, 0, g,r}$ is nonempty and geometrically irreducible of dimension $3g-3+r$.
\end{itemize}
  \epr
\begin{proof}
Assertions (i) and (ii) follows
from Corollary \ref{QS241}, (i) and (ii),  respectively.
\end{proof}

Here, we shall refer to
the following 
Mochizuki's result;  it concerns the moduli stack of dormant $\mfg^\odot (= \mfs \mfl_2)$-opers and   was proved in the work of $p$-adic Teichm\"{u}ller theory (cf. Remark \ref{QQ1001} for the correspondence between $\mfg^\odot$-opers and torally indigenous bundles).

\bt \label{QS1067}
Suppose that $p >2$.
Then, 
$\mfO \mfp^\ZZZ_{\mfg^\odot, g,r}$
 is  a nonempty,  smooth, and  proper
 Deligne-Mumford stack over $k$ of dimension $3g-3+r$.
 Also, 
the projection $\mfO \mfp^\ZZZ_{\mfg^\odot,  g, r} \migi \overline{\mfM}_{g, r}$ (resp., $\mfO \mfp^\ZZZ_{\mfg^\odot,  \rho, g, r} \migi \overline{\mfM}_{g, r}$) is finite, faithfully flat, and generically \'{e}tale.
Finally,  
 $\mfO \mfp^\ZZZ_{\mfg^\odot, g,0} $ (resp., $\mfO \mfp^\ZZZ_{\mfg^\odot, \rho, g,r}$ for $r>0$) is  geometrically connected  if it is nonempty.
 \et
\begin{proof}
See  ~\cite{Mzk2}, Chap.\,II, \S\,2.3, Theorem 2.8.
\end{proof}

 \begin{rema} \label{QS1088}
~\cite{Mzk2}, Chap.\,II, \S\,2.3, Theorem 2.8, does not 
 mention the non-emptiness of  $\mfO \mfp_{\mfg^\odot, g, r}^\ZZZ$, i.e., the moduli stack of dormant torally  indigenous bundles; 
the statement 
 is only concerned with its substack $\mfO \mfp_{\mfg^\odot, \rho, g, r}^\ZZZ$ ($=\overline{\mcN}^{\rho, \mr{dor}}_{g, r}$ in the notation of {\it loc.\,cit.}), and whether it  is  empty or not depends  completely  on  the data of  radii.
But, 
 the non-emptiness of $\mfO \mfp_{\mfg^\odot, g, r}^\ZZZ$  can be immediately verified.
Indeed, to verify this property,  let us consider 
 a totally degenerate curve, i.e., a pointed stable curve obtained  by gluing together $3$-pointed projective lines at their marked points  (cf. Definition \ref{y0178} for the precise definition).
We  can find   the molecule structure on such a curve 
such that its restriction to (the normalization of) each rational component forms
the atom   each of whose radius is $1$.
Then, it follows from ~\cite{Mzk2}, Introduction, Theorem 1.3, (2), that  the indigenous bundle  corresponding  to this molecule is dormant.
It follows that the fibers of the projection $\mfO \mfp_{\mfg^\odot,  g, r}^\ZZZ \migi \overline{\mfM}_{g,r}$ over the points classifying totally degenerate curves are always nonempty, in particular, $\mfO \mfp_{\mfg^\odot,  g, r}^\ZZZ \neq \emptyset$.  
\end{rema}
 
 \begin{rema} \label{QS1089}
Suppose that  $\hslash \in k^{\times}$.
Then, since $\mfO \mfp^\ZZZ_{\mfg^\odot, \hslash, g,r} \cong \mfO \mfp^\ZZZ_{\mfg^\odot,  g,r}$ and $\mfO \mfp^\ZZZ_{\mfg^\odot, \hslash, \hslash \star \rho,  g,r} \cong \mfO \mfp^\ZZZ_{\mfg^\odot, \rho, g,r}$   (cf. (\ref{QS1100})),
 the same assertion as Theorem \ref{QS1067} holds if $\mfO \mfp^\ZZZ_{\mfg^\odot, g,r}$ (resp., $\mfO \mfp^\ZZZ_{\mfg^\odot, \rho,  g,r}$) is replaced with
$\mfO \mfp^\ZZZ_{\mfg^\odot, \hslash, g,r}$ (resp., $\mfO \mfp^\ZZZ_{\mfg^\odot, \hslash, \hslash \star \rho, g,r}$).
 \end{rema}

 \begin{rema} \label{QG37}
 It is an important problem to determine for what elements $\rho \in \mfc^{\times r} (\mbF_p)$ the associated moduli stack $\mfO \mfp^{\ZZZ}_{\mfg, \hslash, \hslash \star \rho, g,r}$ with $\hslash \in k^\times$ is nonempty.  
 We can give  a complete  answer to this problem  for $\mfg = \mfg^\odot$ 
 (cf. ~\cite{Mzk2}, Introduction, \S\,1.2, Theorem 1.3, or the discussion of  \S\,\ref{y0772} in the present manuscript).
 But,  
 little is known about the general case at the time of writing the present manuscript.

 One of the facts that can be easily verified is that  $\mfO \mfp^{\ZZZ}_{\mfg, \hslash, \hslash \star \rho, g,r} = \emptyset$ unless $\rho \in  (\mfc (\mbF_p) \setminus \{ [0] \})^{\times r}$.
 Indeed, suppose, on the contrary, that there exist an $r$-pointed stable curve $\msX$ of genus $g$ and a dormant $(\mfg, \hslash)$-oper $\msE^\spadesuit := (\DE_\mbB, \nabla)$  on $\msX$ of radii $\rho := (\rho_i)_i$ with $\rho_{i_0} = [0]$ for some $i_0 \in \{1, \cdots, r \}$.
 After possibly applying Proposition \ref{y049} and changing the base space of $\msX$,
 let us  assume that $k$ is algebraically closed, $\msX$ is defined over $k$, and $\msE^\spadesuit$ is $\, \qq$-normal.
  The monodromy $\mu_{i_0}^\nabla$ of $\nabla$ at the $i_0$-th marked point  $\sigma_{i_0}$ of $\msX$ may be expressed as $\mu_{i_0}^\nabla = 1 \otimes q_{-1} + R$ for some $R \in \mfg^{\mr{ad}(q_1)}$. 
  But,   $R$ must be zero because  $\mr{Kos}_\mfg$ is an isomorphism and $\rho_{i_0} = [0]$.
  It follows that 
    \begin{align}
 \sigma^*_{i_0} ({^p}\psi^\nabla) &= (\mu_{i_0}^\nabla)^{[p]} -\hslash^{p-1} \cdot \mu_{i_0}^\nabla \\
 &= 1 \otimes q_{-1}^{[p]} -\hslash^{p-1} \cdot (1 \otimes q_{-1}) \notag \\
 & = - \hslash^{p-1} \otimes q_{-1}  \notag \\
 &\neq 0.\notag
  \end{align}
  This contradicts the dormancy condition on $\msE^\spadesuit$.
  Consequently, we have proved the desired property.
\end{rema}

For the moduli stack  $\mfO \mfp^\ZZZ_{\mfg, \hslash, g,r}$ with $\mfg$ general, 
 the following assertion can be deduced from results proved so far together with the above  theorem due to S. Mochizuki:

\bt[cf. Theorem \ref{y016}] \label{QS1040}
\begin{itemize}
\item[(i)]
$\mfO \mfp^\ZZZ_{\mfg, \hslash, g,r}$  may be represented   by   a    proper Deligne-Mumford stacks over $k$, and the natural projections  $\mfO \mfp^\ZZZ_{\mfg, \hslash, g,r}$ $\migi \overline{\mfM}_{g,r}$ is   finite.
\item[(ii)]
The assignment $(\msX, \msE^\spadesuit) \mapsto (\msX, \iota_{\mbG*}(\msE^\spadesuit))$ (cf. (\ref{associated}))  determines  a closed immersion
\begin{align} \label{QS4567}
\mfO \mfp^\ZZZ_{\mfg^\odot, \hslash, g,r}  \migiincl \mfO \mfp^\ZZZ_{\mfg, \hslash, g,r}
\end{align}
over $\overline{\mfM}_{g,r}$.
In particular, $\mfO \mfp^\ZZZ_{\mfg, \hslash, g,r}$ is nonempty, of dimension $3g-3+r$,  and has an irreducible component that dominates $\overline{\mfM}_{g,r}$.
\end{itemize}
\et
\begin{proof}
The assertion follows from Proposition \ref{y081}, (ii), Corollary \ref{QS1010}, and Theorem \ref{QS1067}.
\end{proof}

\subsection{$p$-nilpotent case} 
\label{QS1030}
Next, let us describe an assertion concerning the geometric structure of the universal moduli stack classifying $p$-nilpotent $(\mfg, \hslash)$-opers.
Suppose that $\Psi_\hslash (\rho) = [0]^{\times r}$ (where $\Psi_\hslash (\rho) := \emptyset$ if $\rho = \emptyset$).
As in the case of a pointed stable curve over a scheme, we can construct $\ooalign{$\bigoplus$ \cr $\bigotimes^\nabla_{\msC_{g,r}^{(1)}}$}$  (resp., $\ooalign{$\bigoplus$ \cr $\bigotimes^\nabla_{\Psi_\hslash (\rho), \msC_{g,r}^{(1)}}$}$), i.e., 
 the relative affine space $ \ooalign{$\bigoplus$ \cr $\bigotimes^\nabla_{\msX^{(1)}}$}$ (resp., $\ooalign{$\bigoplus$ \cr $\bigotimes^\nabla_{\Psi_\hslash (\rho), \msX^{(1)}}$}$) for the case where  $\msX$ is taken to be  the universal family of pointed curves $\msC_{g,r}$. 
Then,  the Hitchin-Mochizuki morphism (cf. (\ref{QG2500}))  for $\msC_{g,r}$ defines a morphism 
\index{$\rotatebox[origin=c]{58}{{\large $\kappa$}}^{\hspace{-0.0mm}}_{\hslash (, \rho)}$, $\rotatebox[origin=c]{58}{{\large $\kappa$}}^{\hspace{-0.0mm}\nabla}_{\hslash (, \rho)}$, $\rotatebox[origin=c]{58}{{\large $\kappa$}}^{\hspace{-0.0mm}\nabla}_{\hslash (, \rho) \mr{univ}}$, Hitchin-Mochizuki morphism}
\begin{equation} \label{QG3041}
 \rotatebox[origin=c]{58}{{\large $\kappa$}}^{\hspace{-0.0mm}\nabla}_{\hslash, \mr{univ}} :  \mfO \mfp_{\mfg, \hslash, g,r} \migi  \ooalign{$\bigoplus$ \cr $\bigotimes^\nabla_{\msC_{g,r}^{(1)}}$}  \   
 \left(\text{resp.}, \ \rotatebox[origin=c]{58}{{\large $\kappa$}}^{\hspace{-0.0mm}\nabla}_{\hslash, \rho, \mr{univ}} :  \mfO \mfp_{\mfg, \hslash, \rho, g,r} \migi  \ooalign{$\bigoplus$ \cr $\bigotimes^\nabla_{\Psi_\hslash (\rho), \msC_{g,r}^{(1)}}$}\right) 
 \end{equation}
over $\overline{\mfM}_{g,r}$.
Also, denote by 
\index{$\mfO \mfp^\NILP_{\mfg, \hslash, (\rho,)  \msX}$, $\mfO \mfp^\NILP_{\mfg, \hslash, (\rho,) g,r}$}
\begin{align} \label{QR13J9}
\mfO \mfp^\NILP_{\mfg, \hslash, g,r} \ \left(\text{resp.,} \  \mfO \mfp^\NILP_{\mfg, \hslash, \rho, g,r} \right)
\end{align}
the inverse image of the zero section of $\ooalign{$\bigoplus$ \cr $\bigotimes^\nabla_{\msC_{g,r}^{(1)}}$}$ (resp., $\ooalign{$\bigoplus$ \cr $\bigotimes^\nabla_{\Psi_\hslash (\rho), \msC_{g,r}^{(1)}}$}$)
 via 
$\rotatebox[origin=c]{58}{{\large $\kappa$}}^{\hspace{-0.0mm}\nabla}_{\hslash, \mr{univ}}$ (resp., $\rotatebox[origin=c]{58}{{\large $\kappa$}}^{\hspace{-0.0mm}\nabla}_{\hslash, \rho, \mr{univ}}$), i.e., 
the closed substack of $\mfO \mfp_{\mfg, \hslash, g,r}$ (resp., $\mfO \mfp_{\mfg, \hslash, \rho, g,r}$) classifying $p$-nilpotent $(\mfg, \hslash)$-opers.
The following assertion is   a generalization of a result for $\mfg = \mfs \mfl_n$ proved by S. Mochizuki (cf. ~\cite{Mzk1}, Chap.\,I, Theorem 2.3), as well as of a result for  general $\mfg$'s  by  T.-H. Chen-X. Zhu (cf. the first unpublished version of ~\cite{CZ}, \S\,1, Theorem 1.1).


\bt[cf. Theorem \ref{y015}] 
 \label{y0101} 
\begin{itemize}
\item[(i)]
The morphism
$\rotatebox[origin=c]{58}{{\large $\kappa$}}^{\hspace{-0.0mm}\nabla}_{\hslash, \mr{univ}}$ (resp., $\rotatebox[origin=c]{58}{{\large $\kappa$}}^{\hspace{-0.0mm}\nabla}_{\hslash, \rho, \mr{univ}}$) is  finite and faithfully flat of  degree $p^{\aleph (\mfg)}$  (resp., $p^{^c \aleph (\mfg)}$).
In particular,
$\mfO \mfp_{\mfg, \hslash, g,r}^\NILP$ (resp., $\mfO \mfp_{\mfg, \hslash, \rho, g,r}^\NILP$) may be represented by a nonempty and proper Deligne-Mumford stack over $k$ of dimension $3g-3+r$, and the natural projection $\mfO \mfp_{\mfg, \hslash, g,r}^\NILP \migi \overline{\mfM}_{g,r}$ (resp., $\mfO \mfp_{\mfg, \hslash, \rho, g,r}^\NILP \migi \overline{\mfM}_{g,r}$) is finite and faithfully flat of degree $p^{\aleph (\mfg)}$  (resp., $p^{^c \aleph (\mfg)}$).
\item[(ii)]
If, moreover, $\hslash = 0$, then
the natural projection
\begin{equation}
(\mfO \mfp^\NILP_{\mfg, 0, g,r})_\mr{red} 
 \migi \overline{\mfM}_{g,r}
 \end{equation}
 is an isomorphism of $k$-stacks, where $(\mfO \mfp^\NILP_{\mfg, 0, g,r})_\mr{red}$
    denotes the  reduced stack associated to $\mfO \mfp^\NILP_{\mfg, 0, g,r}$. 
In particular, $\mfO \mfp^\NILP_{\mfg, 0, g,r}$ is geometrically irreducible.
\end{itemize}
  \et
\begin{proof}
The assertions follow from Corollaries \ref{y094} and  \ref{QS300}, (i), (ii).
\end{proof}


\chapter{Flat vector bundles and differential operators} \label{y0104}

In this chapter, we consider the description of $(\mfs \mfl_n, \hslash)$-opers in terms of differential operators and vector bundles.
It  is well-known that
 $\mfs \mfl_n$-opers on  a smooth curve correspond bijectively to  certain differential operators between  line bundles on that curve.
To extend the situation (e.g., to the case of a  family of possibly singular  stable curves), we introduce a generalized theta characteristic, i.e., an $(n, \hslash)$-theta characteristic.
After choosing an $(n, \hslash)$-theta characteristic $\bb$,
we define the notion of an {\it $(n, \hslash, \bb)$-projective connection} and a {\it $(\mr{GL}_n, \hslash, \bb)$-oper}.
For example, an $(n, \hslash, \bb)$-projective connection is defined as a certain $n$-th order differential operator such that the principal symbol is $1$ and the subprincipal symbol is characterized by $\bb$.
The main result of this chapter  provides a generalization of the correspondence mentioned above, connecting  $\mfs \mfl_n$-opers, $(\mr{GL}_n, \hslash, \bb)$-opers, and $(n, \hslash, \bb)$-projective connections. 
Moreover, we discuss the equivalence of the dormancy condition on opers in positive characteristic and the fullness of root functions of the corresponding projective connections.

\begin{notation}
Let  us fix a perfect field $k$
\index{$k$, field} 
 and a positive integer $n$.
\end{notation}

\section{Logarithmic $\hslash$-connections on vector bundles} \label{y0105}

In Definition \ref{y021}, we already defined  a logarithmic  connection on a principal bundle.
In what follows, we formulate and study the version  defined on a vector bundle.
As proved in \S\,\ref{y0106} later, these  definitions of a connection  are equivalent under the identification between  rank $n$ vector bundles and   $\mr{GL}_n$-bundles.

\subsection{Logarithmic $\hslash$-connections}\label{QS4800}
Let $T^\mr{log} := (T, \alpha_T : M_T \migi \mcO_T)$, 
$Y^{\mr{log}} := (Y, \alpha_Y : M_Y \migi \mcO_Y)$ be fine log schemes over $k$  and 
$f^\mr{log} : Y^\mr{log} \migi T^\mr{log}$ a  $k$-morphism of log schemes.
Also, let 
 $\hslash$  be an element of $\Gamma (T, \mcO_T)$  and $\mcV$   a rank $n$   vector bundle  on $Y$  (or more generally, an $\mcO_Y$-module).
 To simplify the notation, we write $\mcT := \mcT_{Y^\mr{log}/T^\mr{log}}$ and $\Omega := \Omega_{Y^\mr{log}/T^\mr{log}}$.

\bde 
 \label{QQ65}
An {\bf $\hslash$-$T^\mr{log}$-connection}
\index{$\nabla$, $\hslash$-$T^\mr{log}$-connection}
 on  $\mcV$ is  an
$f^{-1}(\mcO_T)$-linear morphism
\begin{equation}
 \Nab : \mcV \migi \Omega \otimes \mcV 
  \end{equation}
 satisfying  
 $\Nab (a \cdot v) = \hslash \cdot da \otimes v + a \cdot \Nab (v)$
 for any  local sections $a \in \mcO_Y$, $v \in \mcV$.
 Each local section in  the kernel $\mr{Ker}(\nabla)$ is called {\bf horizontal}.
 \index{$\nabla$, $\hslash$-$T^\mr{log}$-connection@--- horizontal (section) of}
 We shall refer to any $1$-$T^\mr{log}$-connection  as a {\bf $T^\mr{log}$-connection}.
  \ede

\begin{exa}[The case of $\hslash =0$] \label{QS9900}
A $0$-$T^\mr{log}$-connection on  $\mcV$ is nothing but an {\it $\mcO_Y$-linear} morphism $\mcV \migi \Omega \otimes \mcV$, or equivalently an element of $\Gamma (Y, \Omega \otimes \mcE nd_{\mcO_Y}(\mcV))$.
In particular, there always exists a $0$-$T^\mr{log}$-connection on any vector bundle.
A $0$-$T^\mr{log}$-connection are often  called a {\it Higgs field}.
\end{exa}

\subsection{Curvature} \label{QS5006} 
Now, let $\Nab: \mcV \migi \Omega \otimes \mcV$ be an $\hslash$-$T^\mr{log}$-connection on $\mcV$.
Given  a local section $\partial$ of $\mcT \left(= \Omega^\vee \right)$,
we denote by  $\Nab_\partial$  the $f^{-1}(\mcO_T)$-linear endomorphism of $\mcV$ defined as $(\partial \otimes \mr{id}_\mcV) \circ \Nab$.
Also, for  local sections $\partial_1, \partial_2 \in \mcT$,
we write 
\begin{align}
\psi^\nabla (\partial_1 \wedge \partial_2) :=  [\Nab_{\partial_1}, \Nab_{\partial_2}] - \Nab_{[\partial_1, \partial_2]}.
\end{align}
It is verified that
$\psi^\nabla (\partial_1 \wedge \partial_2) = - \psi^\nabla (\partial_2 \wedge \partial_1)$ and $\psi^\nabla ((a \cdot \partial)\wedge \partial_2) = \psi^\nabla (\partial_1 \wedge (a \cdot \partial_2)) = a\cdot \psi^\nabla (\partial_1 \wedge \partial_2)$ for any local sections $\partial_1, \partial_2 \in \mcT$ and $a \in \mcO_Y$.
Hence,
the assignment  $\partial_1 \wedge \partial_2 \mapsto \psi^\nabla (\partial_1 \wedge \partial_2)$
determines  an {\it $\mcO_Y$-linear} morphism
\index{$\nabla$, $\hslash$-$T^\mr{log}$-connection@--- $\psi^\nabla$, curvature of}
\begin{align} \label{QQ78}
\psi^{\Nab} : \bigwedge^2 \mcT \migi \mcE nd_{f^{-1}(\mcO_T)} (\mcV).
\end{align}

\bde[] \label{QQ1099}
We shall refer to $\psi^\nabla$ as the {\bf curvature}  of $\Nab$.
Also,  we shall say that $\Nab$ is {\bf flat} (or {\bf integrable}) \index{$\nabla$, $\hslash$-$T^\mr{log}$-connection@--- flat (= integrable)} if $\psi^{\Nab}=0$.
  (In particular, if $Y^\mr{log}/T^\mr{log}$ is a log curve, which implies that $\mcT$ is a line bundle, then any $\hslash$-$S^\mr{log}$-connection is automatically flat.)
\ede

\begin{rema}
It is immediately verified that 
the image of $\psi^{\Nab}$ lies  in the subsheaf $\mcE nd_{\mcO_Y}(\mcV)$ of $\mcE nd_{f^{-1}(\mcO_T)} (\mcV)$.
Hence, $\psi^{\Nab}$ may be regarded as an $\mcO_Y$-linear morphism $\bigwedge^2 \mcT \migi \mcE nd_{\mcO_Y}(\mcV)$, or equivalently  a global section of $\left(\bigwedge^2 \Omega\right) \otimes \mcE nd_{\mcO_Y}(\mcV)$.
\end{rema}

\bde \label{QQ76}
\begin{itemize}
\item[(i)]
An {\bf $\hslash$-flat (vector) bundle}
\index{$\hslash$-flat (vector) bundle}
 on $Y^\mr{log}/T^\mr{log}$ (of rank $n$) is a pair 
\begin{align} \label{QS53278}
\msV := (\mcV, \Nab)
\end{align}
 consisting of a vector bundle $\mcV$ on $Y$ (of rank $n$) and  a flat  $\hslash$-$T^\mr{log}$-connection $\Nab$ on $\mcV$. If, moreover, $\mcV$ is of rank $1$, then  $\msV$ is called  an {\bf $\hslash$-flat line bundle} on $Y^\mr{log}/T^\mr{log}$.  
For simplicity, a $1$-flat bundle is called  a {\bf flat (vector) bundle}.
\item[(ii)]
Let $\msV := (\mcV, \Nab)$ and $\msV' :=(\mcV', \Nab')$ be $\hslash$-flat  bundles on $Y^\mr{log}/T^\mr{log}$.
A {\bf morphism of $\hslash$-flat bundles} from $\msV$ to $\msV'$ is an $\mcO_Y$-linear morphism 
$\eta : \mcV \migi \mcV'$ compatible  with the  $\hslash$-$T^\mr{log}$-connections, i.e.,  satisfying   $\nabla' \circ \eta = (\mr{id}_\Omega \otimes \eta)\circ \nabla$.
\end{itemize}
 \ede

\begin{exa}[Trivial bundle] \label{QQ80}
A trivial example of an $\hslash$-flat bundle can be  provided  by the pair 
\begin{align} \label{QS6000}
(\mcO_Y, \hslash \cdot d)
\end{align}
  consisting of   the structure sheaf $\mcO_Y$ and $\hslash \cdot d$, i.e.,  the universal derivation $\mcO_Y \migi \Omega$ multiplied by $\hslash$;
we refer to it  as the {\bf trivial $\hslash$-flat (line) bundle}.

Also,  we can describe each  $\hslash$-$T^\mr{log}$-connection $\Nab$ on $\mcO_Y$   as $\nabla = \hslash \cdot d + \nabla^{\mr{GL}_1}$ for a unique $1$-form   $\nabla^{\mr{GL}_1} \in \Gamma (Y, \Omega) \left(= \mr{Hom}_{\mcO_Y}(\mcO_Y, \Omega\otimes \mcO_Y) \right)$.
It is immediately  verified  that $\Nab$ is flat if and only if $\nabla^{\mr{GL}_1}$ is closed, i.e., $d\nabla^{\mr{GL}_1}=0$.

In general, for any  $\hslash$-$T^\mr{log}$-connection $\nabla$ on $\mcO_Y^{\oplus n}$,
there exists a unique element
\index{$\nabla$, $\hslash$-$T^\mr{log}$-connection@--- $\nabla^{\mr{GL}_n}$}
\begin{align} \label{QS7144}
\nabla^{\mr{GL}_n}  \in \mr{Hom}_{\mcO_Y} (\mcO_Y^{\oplus n}, \Omega \otimes (\mcO_Y^{\oplus n})) \left(= \Gamma (Y, \Omega \otimes \mcE nd_{\mcO_Y}(\mcO_Y^{\oplus n})) \right)
\end{align}
satisfying  $\nabla  = \hslash \cdot d^{\oplus n} + \nabla^{\mr{GL}_n}$. \end{exa}

\subsection{Pull-back and base-change} \label{QS8100}
Let $\mcV$ be a vector bundle on $Y$ and $\nabla$ an $\hslash$-$T^\mr{log}$-connection on  $\mcV$.
Also, let $Y'^{\mr{log}}$ be a fine log scheme equipped with a {\it log \'{e}tale} morphism $u^\mr{log} : Y'^{\mr{log}} \migi Y^\mr{log}$.
Then, we can  pull-back $\nabla$ to obtain an $\hslash$-$T^\mr{log}$-connection on $u^*(\mcV)$.
Indeed, because of the log \'{e}taleness of $u^\mr{log}$,
the induced  morphism  $u^*(\Omega) \migi \Omega' := \Omega_{Y'^{\mr{log}}/T^\mr{log}}$ is an isomorphism.
Under the identification $u^*(\Omega)=\Omega'$ given by this isomorphism,
we define the {\bf pull-back} $u^* (\nabla)$ 
\index{$\nabla$, $\hslash$-$T^\mr{log}$-connection@--- $u^*(\nabla)$, pull-back of $\nabla$}
of $\nabla$ to be the $f^{-1}(\mcO_T)$-linear morphism
\begin{align} \label{QS8300}
u^*(\nabla) : & \ u^*(\mcV) \left(= \mcO_{Y'} \otimes_{u^{-1}(\mcO_Y)}u^{-1}(\mcV) \right)
 \\ & \ 
 \migi \Omega' \otimes u^*(\mcV)\left(= \mcO_{Y'} \otimes_{u^{-1}(\mcO_Y)}u^{-1}(\Omega \otimes \mcV)\right)
 \notag
\end{align}
given by $a \otimes v \mapsto da \otimes v + a \otimes  \nabla (v)$ for any local sections $a \in \mcO_{Y'}$, $v \in \mcV$.
This morphism  forms an $\hslash$-$T^\mr{log}$-connection on $u^*(\mcV)$  and satisfies  $\psi^{u^*(\nabla)}=u^*(\psi^\nabla)$.
In particular, the flatness of $\nabla$ implies that of $u^*(\nabla)$.

Also, the base-change of $\nabla$ to a fine  $T^\mr{log}$-log scheme   can be constructed.
 To observe this,  suppose that we are given
 a fine log scheme  $T'^{\mr{log}}$ over $k$ equipped with a $k$-morphism $t^\mr{log} : T'^{\mr{log}} \migi T^\mr{log}$.
Write $Y'^{\mr{log}}:= T'^{\mr{log}} \times_{T^\mr{log}} Y^\mr{log}$ and use the notation $t^*(-)$ to denote base-change  to $Y'$.
 Then, $\Omega' := \Omega_{Y'^{\mr{log}}/T'^{\mr{log}}}$ may be identified with $t^*(\Omega)$.
By means of this identification, we can define the {\bf base-change}
\index{$\nabla$, $\hslash$-$T^\mr{log}$-connection@--- $t^*(\nabla)$, base-change}
 \begin{align} \label{QS8301}
t^*(\nabla) : t^*(\mcV) \migi \Omega' \otimes t^*(\mcV)
\end{align}
of $\nabla$ in the same manner as (\ref{QS8300}).
It is immediate that $\psi^{t^* (\nabla)} = t^*(\psi^\nabla)$ and hence, the flatness of $\nabla$ implies that of $t^*(\nabla)$.

\subsection{Determinant}\label{QQ1105}
In what follows, we recall several  constructions of $\hslash$-flat bundles by using  a  given  $\hslash$-flat bundle.
Let us fix  a vector bundle $\mcV$ on $Y$ of rank $n$ and
an $\hslash$-$T^\mr{log}$-connection $\Nab : \mcV \migi \Omega\otimes \mcV$ on $\mcV$.
Then, the  {\bf determinant} of $\nabla$ is the $\hslash$-$T^\mr{log}$-connection
\index{$\nabla$, $\hslash$-$T^\mr{log}$-connection@--- $\mr{det}(\Nab)$, determinant of}
\begin{align} \label{QQ291}
\mr{det}(\Nab) : \mr{det}(\mcV) \migi \Omega \otimes \mr{det}(\mcV)
\end{align}
 on 
 the determinant  bundle $\mr{det}(\mcV) = \bigwedge^n \mcV$ of $\mcV$
   given   by 
\begin{align} \label{QQ289}
\mr{det}(\nabla) (v_1 \wedge v_2 \wedge \cdots \wedge v_n) = \sum_{j=1}^n v_1 \wedge \cdots \wedge \nabla (v_j) \wedge \cdots \wedge v_n
\end{align}
for $v_1, \cdots, v_n  \in \mcV$.
The curvature $\psi^{\mr{det}(\Nab)}$ of $\mr{det}(\Nab)$
satisfies 
\begin{align} \label{QS7006}
\psi^{\mr{det}(\Nab)} = \mr{tr}(\psi^{\Nab})
\end{align}
under the natural identification $\mr{tr}(\mcE nd_{\mcO_Y}(\mcV)) = \mcO_Y$.
In particular,  the flatness of $\nabla$ implies that of $\mr{det}(\nabla)$.

\subsection{Dual} \label{QSQS7001}
Let $\mcV$ and $\nabla$ be as above.
Then, the  {\bf dual} of $\nabla$ is the 
$\hslash$-$T^\mr{log}$-connection
\index{$\nabla$, $\hslash$-$T^\mr{log}$-connection@--- $\Nab^\vee$, dual of $\nabla$}
\begin{align} \label{QQ297}
\Nab^\vee : \mcV^\vee \migi \Omega \otimes \mcV^\vee
\end{align}
on the dual bundle $\mcV^\vee$ of $\mcV$  determined by 
\begin{align} \label{QQ298}
\Nab^\vee (h) (v) = d (h (v))-  (\mr{id}_\Omega \otimes h)  \circ \Nab (v)
\end{align}
for local sections $h \in \mcV^\vee$, $v \in \mcV$.
It is immediate that  $(\mcV^{\vee \vee}, \nabla^{\vee \vee}) \cong  (\mcV, \nabla)$. 
Under the natural identification $\mcE nd_{\mcO_Y} (\mcV^\vee) = \mcE nd_{\mcO_Y}(\mcV)$, the curvature $\psi^{\nabla^\vee}$ of $\nabla^\vee$ satisfies 
\begin{align}
\psi^{\nabla^\vee}  = - {^t}\psi^\nabla, 
\end{align}
where ${^t}\psi^\nabla$ denotes the transpose of $\psi^\nabla$.
In particular, the flatness of $\nabla$ is equivalent to the flatness of $\nabla^\vee$.

\subsection{Tensor product} \label{QS7002}
Next, let 
$\mcV'$ be another vector bundle on $Y$ and $\Nab'$ an $\hslash$-$T^\mr{log}$-connection on $\mcV'$.
The {\bf tensor product} of $\Nab$ and $\Nab'$ is the $\hslash$-$T^\mr{log}$-connection 
\index{$\nabla$, $\hslash$-$T^\mr{log}$-connection@--- $\Nab \otimes \Nab'$, tensor product}
\begin{align} \label{QQ295}
\Nab \otimes \Nab' : \mcV \otimes \mcV' \migi \Omega \otimes (\mcV \otimes \mcV')
\end{align}
 on the tensor product $\mcV \otimes \mcV'$ determined  by
\begin{align} \label{QQ290}
\Nab \otimes \Nab' (v \otimes v') = \Nab (v) \otimes v' + v \otimes \Nab'(v')
\end{align}
for local sections $v \in \mcV$, $v' \in \mcV'$.
In particular, for $m \in \mbZ_{>0}$,  $\Nab$ induces an  $\hslash$-$T^\mr{log}$-connection
\index{$\nabla$, $\hslash$-$T^\mr{log}$-connection@--- $\Nab^{\otimes m}$, $m$-fold tensor product of}
\begin{align} \label{QQ296}
\Nab^{\otimes m} : \mcV^{\otimes m} \migi \Omega\otimes \mcV^{\otimes m}
\end{align}
on 
 the $m$-fold tensor product $\mcV^{\otimes m}$
   of $\mcV$.
 The curvature $\psi^{\Nab \otimes \Nab'}$ (resp., $\psi^{\nabla^{\otimes m}}$)
 of $\nabla \otimes \nabla'$ (resp., $\nabla^{\otimes m}$) can be expressed as
 \begin{align} \label{QS7007}
 & \hspace{12mm}\psi^{\Nab \otimes \Nab'} = \psi^{\nabla} \otimes \mr{id}_{\mcV'} + \mr{id}_\mcV \otimes  \psi^{\nabla'} \\
&\left(\text{resp.,} \  \psi^{\nabla^{\otimes m}} = \sum_{j=1}^m \mr{id}_{\mcV}^{\otimes (j-1)} \otimes \psi^\nabla \otimes \mr{id}_\mcV^{\otimes (m-j)}
\right).\notag
 \end{align}
 In particular, if  $\mcV$ is a line bundle, then we have $\psi^{\nabla^{\otimes m}} = m \cdot \psi^\nabla$. 
 Moreover,  (\ref{QS7007}) implies that $\Nab \otimes \Nab'$ (resp., $\Nab^{\otimes m}$)  is flat if both $\Nab$ and $\Nab'$  are flat (resp., $\Nab$ is flat). 

\subsection{Gauge transformation} \label{QS7100}
Let $\mcV$, $\mcV'$  be  vector bundles on $Y$ and
$\eta : \mcV \isom \mcV'$ an $\mcO_Y$-linear isomorphism.
Then, for an $\hslash$-$T^\mr{log}$-connection $\nabla$ on $\mcV$,
the composite
\begin{align} \label{QW299}
\eta_* (\nabla) : \mcV' \xrightarrow{\eta^{-1}} \mcV \xrightarrow{\nabla} \Omega \otimes \mcV \xrightarrow{\mr{id}_\Omega \otimes \eta} \Omega \otimes \mcV'
\end{align}
specifies 
an $\hslash$-$T^\mr{log}$-connection on $\mcV'$.

\bde \label{QS7109}
Suppose that $\mcV' = \mcV$, i.e., $\eta$ is an automorphism of $\mcV$.
Then, we shall refer to $\eta_* (\nabla)$ as the {\bf gauge transformation} of $\nabla$ by $\eta$. 
\index{$\nabla$, $\hslash$-$T^\mr{log}$-connection@--- $\eta_* (\nabla)$, gauge transformation of $\nabla$ by $\eta$}
\ede

The following assertion may be thought of   as an extension of the discussion in Remark \ref{ppp037} to the case of gauge  transformation by an automorphism that is not necessarily unipotent.

\bpr \label{QS7120}
Let  $\eta$ be  an automorphism of $\mcO_Y^{\oplus n}$ and $\nabla$  an $\hslash$-$T^\mr{log}$-connection on $\mcO_Y^{\oplus n}$.
Then, the gauge transformation $\eta_*(\nabla)$ of $\nabla$ by $\eta$ satisfies the following equality:
\begin{align} \label{QS7121}
\eta_*(\nabla)^{\mr{GL}_n} = - \hslash \cdot \eta^{-1} d \eta + \eta \circ \nabla^{\mr{GL}_n} \circ \eta^{-1}  \ \left(=   \hslash \cdot \eta d (\eta^{-1}) + \eta \circ \nabla^{\mr{GL}_n} \circ \eta^{-1}\right)
\end{align}
(cf. (\ref{QS7144}) for the definition of $(-)^{\mr{GL}_n}$),
where $d \eta$ denotes the image of $\eta$ via 
the morphism $\mcE nd_{\mcO_Y} (\mcO_Y^{\oplus n}) \migi \Omega \otimes \mcE nd_{\mcO_Y} (\mcO_Y^{\oplus n})$ given by $(a_{ij})_{ij} \mapsto (d a_{ij})_{ij}$ under the natural identification $\mcE nd_{\mcO_Y} (\mcO_Y^{\oplus n}) = \mcO_Y^{\oplus n\times n}$.
\epr
\begin{proof}
The assertion follows from the definition of $\eta_*(\nabla)$.
\end{proof}

\section{Logarithmic differential operators with parameter $\hslash$} \label{y0114}

\subsection{The ring of logarithmic differential operators} \label{QS7010}
Let $Y^\mr{log}/T^\mr{log}$ and $\hslash$ be as above.
Let us consider the Lie algebroid  $(\mcT_{}, \hslash \cdot [-, -], \hslash \cdot\mr{id}_{\mcT})$, i.e., the tangent Lie algebroid rescaled by $\hslash$.
Denote by
\begin{equation}  \label{QS9042}
\mcD^{< \infty}_{\hslash, Y^\mr{log}/T^\mr{log}}
\index{$\mcD^{< j}_{\hslash, Y^\mr{log}/T^\mr{log}}$, $\mcD^{<j}_\hslash$}
\end{equation}
the universal  enveloping algebra of $(\mcT_{}, \hslash \cdot [-, -], \hslash \cdot\mr{id}_{\mcT})$ (cf. ~\cite{BeBe}, \S\,1.2.5).
 Explicitly, it is 
  the Zariski sheaf on $Y$
   generated, as a sheaf of (noncommutative) rings, by $\mcO_U$ and  
$\mcT_{}$
 subject to the following relations:
\begin{itemize}
\item
$a_1 \ast a_2 = a_1 \cdot a_2$;
 \item
  $a_1 \ast \partial_1 = a_1 \cdot \partial_1$;
\item
$\partial_1 \ast \partial_2 -\partial_2 \ast \partial_1= \hslash \cdot [\partial_1, \partial_2] $; 
\item
 $\partial_1 \ast a_1-a_1 \ast \partial_1 = \hslash \cdot \partial_1 (a_1)$,
\end{itemize}
for local sections $a_1$, $a_2 \in \mcO_Y$ and $\partial_1$, $\partial_2 \in \mcT_{}$, where ``$\ast$'' denotes the multiplication in $ \mcD^{< \infty}_{\hslash, Y^\mr{log}/T^\mr{log}}$.
We occasionally consider $\mcD^{< \infty}_{\hslash, Y^\mr{log}/T^\mr{log}}$ as a sheaf on $\mfE \mft_{/Y}$ in the usual manner. 
Also, 
 it 
admits  a natural ``Poincare-Birkhoff-Witt'' filtration 
\begin{align}
\mcD_{\hslash, Y^\mr{log}/T^\mr{log}}^{<\infty} = \bigcup_{j \in \mbZ}\mcD_{\hslash, Y^\mr{log}/T^\mr{log}}^{< j},
\end{align}
 where $\mcD_{\hslash, Y^\mr{log}/T^\mr{log}}^{<j} = 0$ for  every $j \leq 0$, $\mcD_{\hslash, Y^\mr{log}/T^\mr{log}}^{<1} = \mcO_U$, and $\mcD_{\hslash, Y^\mr{log}/T^\mr{log}}^{< j+1}= \mcD_{\hslash, Y^\mr{log}/T^\mr{log}}^{<j} + \mcT \cdot \mcD_{\hslash, Y^\mr{log}/T^\mr{log}}^{<j}$ for every $j \geq 1$.
 If there is no fear of confusion, then 
 we  shall write $\mcD_\hslash^{<j} := \mcD^{<j}_{\hslash, Y^\mr{log}/T^\mr{log}}$ ($j \in \mbZ \sqcup \{\infty \}$) for simplicity.
The sheaf $\mcD^{< j}_{\hslash}$  admits two different  structures  of  $\mcO_Y$-module, i.e., one as given by left multiplication, where we denote this $\mcO_Y$-module by ${^L}\mcD^{<  j}_{\hslash}$, and  the other given by right multiplication, where  we denote this $\mcO_Y$-module by ${^R}\mcD^{< j}_{\hslash}$.
  If $\mr{gr}( \mcD^{< \infty}_{\hslash})$ denotes the graded $\mcO_Y$-algebra associated to the filtration  $\{ \mcD_\hslash^{<j} \}_j$, then
  the two $\mcO_Y$-module structures on
  this sheaf 
   arising from ${^L}\mcD_\hslash^{< \infty}$ and ${^R }\mcD_\hslash^{< \infty}$ coincide.
  Also,
   there exists a canonical isomorphism
 \begin{equation} \label{GL6}
  \mr{gr}(\mcD^{< \infty}_{\hslash}) \isom \mbS_{\mcO_Y} (\mcT ) \left(= \bigoplus_{j \in \mbZ_{\geq 0}} \mbS^j_{\mcO_Y}(\mcT) \right)
 \end{equation}
 of graded $\mcO_Y$-algebras, in particular, we have 
 $\mcD_{\hslash}^{< j+1}/\mcD_{\hslash}^{< j} \cong \mbS_{\mcO_Y}^j(\mcT)$  for every $j \in \mbZ_{\geq 0}$.

 Given a  vector bundle $\mcV$ on $Y$  (or more generally, an $\mcO_Y$-module), we equip  the tensor product ${\mcD}^{<j}_{\hslash} \otimes \mcV :=  {^R\mcD}^{<j}_{\hslash} \otimes \mcV$ (resp., $\mcV \otimes {\mcD}^{<j}_{\hslash} := \mcV \otimes {^L \mcD}^{<j}_{\hslash}$) with
   a  structure of  $\mcO_Y$-module
   given by left (resp., right) multiplication.
 
 \begin{rema} \label{QH1000}
 If $\hslash =1$,  then 
the  sheaf $\mcD^{< \infty}_{1}$  is nothing but  (the logarithmic version of) the ring of  {\it crystalline differential operators} in the terminology of ~\cite{BMR}, \S\,1.2, or (in a more general situation) {\it PD differential operators} in the terminology of ~\cite{Og2}, \S\,4.
\end{rema}

 \begin{exa}[The case of ``$\hslash =0$''] \label{QQ300}
 Let us consider the case of $\hslash = 0$.
 Then,
 $\mcD_0^{<\infty}$ is, by definition,  canonically isomorphic to
 the symmetric algebra $\mbS_{\mcO_Y} (\mcT)$ of $\mcT$.
  The two  $\mcO_U$-module structures on $\mcD_0^{<\infty}$ coincide, i.e., ${^L}\mcD_0^{<\infty} = {^R}\mcD_0^{<\infty}$.
 In particular, 
 for  an $\mcO_U$-module $\mcV$, the tensor product $ \mcD_0^{< \infty} \otimes \mcV$ has a decomposition
 \begin{align} \label{QQ304}
 \mcD_0^{< \infty} \otimes \mcV \isom \left(\mcV \otimes \mcD_{0}^{<\infty} =\right) \bigoplus_{j \in \mbZ_{\geq 0}} \mcV \otimes \mbS_{\mcO_Y}^j(\mcT).
 \end{align}
 \end{exa}

\subsection{Flat $\hslash$-connections = left $\mcD_\hslash^{<\infty}$-actions}\label{QS7020}
Let $\mcV$ be an $\mcO_Y$-module  and $\Nab  : \mcV \migi \Omega_{} \otimes \mcV$   a flat $\hslash$-$T^\mr{log}$-connection on $\mcV$.
Then, $\Nab$ induces 
  a  left ${\mcD}^{<\infty}_{\hslash}$-action
 \index{$\nabla$, $\hslash$-$T^\mr{log}$-connection@--- ${^\mcD}\hspace{-1.5mm}\nabla$}
 \begin{equation}
 \label{GL70}
 {^\mcD}\hspace{-1.5mm}\nabla : {\mcD}^{<\infty}_{\hslash} \otimes \mcV \migi \mcV \ \left(\text{or equivalently}, \ {\mcD}^{<\infty}_{\hslash}  \migi \mcE nd_{f^{-1}(\mcO_T)}(\mcV) \right)
  \end{equation}
on $\mcV$ determined uniquely by the condition that 
${^\mcD}\hspace{-1.5mm}\nabla (\partial \otimes v) = \nabla_{\partial} (v)$ for any local sections $v \in \mcV$,  $\partial \in \mcT$.
The assignment $\nabla \mapsto {^\mcD}\hspace{-1.5mm}\nabla$ yields  
 a bijective correspondence
 \begin{align} \label{QS7f040}
\begin{pmatrix}
\text{the set of} \\
\text{flat $\hslash$-$T^\mr{log}$-connections on $\mcV$} 
\end{pmatrix}
\isom \begin{pmatrix}
\text{the set  of } \\
\text{ left ${\mcD}^{<\infty}_{\hslash}$-actions on $\mcV$} \\
\end{pmatrix}.
\end{align}

\begin{rema} \label{QS7048}
The structure sheaf $\mcO_Y$ is equipped with the  left $\mcD_\hslash^{<\infty}$-action ${^\mcD}\hspace{-1.5mm}\nabla_\hslash$ determined uniquely by  $(\partial, a) \mapsto \hslash \cdot \partial (a)$ for local sections $\partial \in \mcT$, $a \in \mcO_Y$;
it corresponds, via (\ref{QS7f040}),   to the $\hslash$-$T^\mr{log}$-connection $\hslash \cdot d$ (cf. (\ref{QS6000})).
Note that, if $\mr{char}(k) >0$,  then
 the  $\mcD_\hslash^{< \infty}$-action ${^\mcD}\hspace{-1.5mm}\nabla_\hslash$  is not faithful,
i.e., the corresponding  morphism $\mcD_\hslash^{< \infty} \migi \mcE nd_{f^{-1}(\mcO_T)}(\mcO_Y)$
 is not injective.
But, by considering the local description of ${^\mcD}\hspace{-1.5mm}\nabla_\hslash$,
we see that
 its restriction 
 $\mcD_{\hslash}^{< j} \migi \mcE nd_{f^{-1}(\mcO_T)}(\mcO_Y)$
      to  the subsheaf $\mcD_{\hslash}^{< j}$ becomes injective when $j \leq \mr{char}'(k)$ and  $\hslash \in \Gamma (T, \mcO_T^\times)$, where   $\mr{char}' (k) := \infty$   (resp., $\mr{char}' (k) := p$) if $\mr{char}(k) =0$ (resp., $\mr{char}(k) = p >0$).
Thus,   $\mcD_{\hslash}^{< j}$ may be regarded as a subsheaf of  $\mcE nd_{f^{-1}(\mcO_T)}(\mcO_Y)$ via this injection.
\end{rema}

\begin{rema} \label{QQ425}

The discussion in the above remark can be extended to 
  differential operators between two line bundles.
For simplicity,  suppose that $Y^\mr{log}/T^\mr{log}= Y/T$, i.e., $f : Y \migi T$ is a smooth scheme over a $k$-scheme $T$.
Let $\mcL_i$ ($i=1,2$) be line bundles on $Y$.
By a {\bf  (linear) differential operator} \index{(linear) differential operator} (over $T$) from $\mcL_1$ to $\mcL_2$, we mean an  $f^{-1}(\mcO_T)$-linear  morphism $D  :\mcL_1 \migi \mcL_2$ 
locally expressed, after fixing identifications $\mcL_1 \cong \mcL_2 \cong \mcO_Y$ and a local coordinate system $\vec{x}:= (x_1, \cdots, x_n)$ in $Y$ relative to $T$, as 
\begin{align} \label{E442}
D 
: v \mapsto 
D(v) =  \sum_{\alpha \in \mbZ_{\geq 0}^n} a_\alpha \cdot \partial^\alpha_{\vec{x}} (v),
\end{align}
 where $a_\alpha$'s are 
  local sections of $\mcO_Y$  with $a_\alpha =0$ for almost all $\alpha$, and 
 for each $\alpha := (\alpha_1, \cdots, \alpha_n) \in \mbZ_{\geq 0}^n$, 
  we write
 $\partial^\alpha_{\vec{x}}(v) := \frac{\partial^{|\alpha|}(v)}{\partial x_1^{\alpha_1} \cdots \partial x_n^{\alpha_n}}$ ($|\alpha| := \alpha_1 + \cdots + \alpha_n$).
Now, suppose that these local sections $a_\alpha$ satisfy   $a_\alpha=0$ for every  $\alpha$ with $|\alpha| \geq \mr{char}'(k)$, where $\mr{char}'(k)$ is as in Remark \ref{QS7048}.
Then,    the nonnegative integer $j_\mr{max} := \mr{max}\left\{ |\alpha| \, | \, a_\alpha \neq 0 \right\}$ is well-defined (i.e., depend only on $D$,  not on  the local expression (\ref{E442})), where $j_\mr{max} := -\infty$ when $D =0$.
In this situation, 
 we say that $D$ is {\bf of order $j_\mr{max}$}. 
 Given a nonnegative integer $j \leq \mr{char}'(k)$, we denote by
\begin{align} \label{QS7100}
\mcD \textit{iff}^{\,< j}(\mcL_1, \mcL_2)
\index{$\mcD \textit{iff}^{\,< j}(\mcL_1, \mcL_2)$}
\end{align} 
the Zariski sheaf on $Y$ consisting of locally defined differential operators from $\mcL_1$ to $\mcL_2$ of order $< j$.
Because of the fact mentioned in Remark \ref{QS7048},
$\mcD \textit{iff}^{\,< j}(\mcL_1, \mcL_2)$ may be regarded as 
a subsheaf of  $\mcH om_{f^{-1}(\mcO_T)} (\mcL_1, \mcL_2)$.
Moreover,  composing with the $f^{-1}(\mcO_T)$-linear morphism $\mcL_2 \otimes \mcD_1^{< j} \migi \mcL_2$ given by
$v \otimes D \mapsto v \otimes D(1)$ yields
  an isomorphism
\begin{align} \label{Edr56}
\mcH om_{\mcO_Y} (\mcL_1, \mcL_2 \otimes \mcD_1^{< j}) \isom \mcD {\it iff}^{< j}(\mcL_1, \mcL_2) \left( \subseteq \mcH om_{f^{-1}(\mcO_T)} (\mcL_1, \mcL_2) \right).
\end{align}
This isomorphism for $\mcL_1 = \mcL_2 = \mcO_Y$ is nothing but the inclusion  $\mcD_\hslash^{j}\migiincl \mcE nd_{f^{-1}(\mcO_T)}(\mcO_Y)$  for  $\hslash =1$ discussed  in Remark \ref{QS7048}.
\end{rema}

\section{Comparison of  log connections} \label{y0106}
 
 \subsection{}
 This section is devoted to proving the following proposition, asserting that 
  giving 
 a rank $n$ $\hslash$-flat bundle  is equivalent to  giving  
  an $\hslash$-flat $\mr{GL}_n$-bundle.

\bpr \label{QQ97}
Let $\mcV$ be a rank  $n$ vector bundle on $Y$ and denote by $\mcE$ the $\mr{GL}_n$-bundle corresponding to $\mcV$.
Then, 
there exists a canonical bijection 
\normalfont{
  \begin{align} \label{QHH34}
\begin{pmatrix}
\text{the set of} \\
\text{$\hslash$-$T^\mr{log}$-connections on $\mcV$} 
\end{pmatrix}
\isom \begin{pmatrix}
\text{the set   of } \\
\text{$\hslash$-$T^\mr{log}$-connections on $\mcE$} \\
\end{pmatrix}.
\end{align}}
Moreover,  an $\hslash$-$T^\mr{log}$-connection   on $\mcV$ is flat if and only if the corresponding  $\hslash$-$T^\mr{log}$-connection on $\mcE$ is flat.
\epr

The following assertion
follows directly from the construction of
(\ref{QHH34}).

 \bco \label{QS8702}
 There is an equivalence of categories
 \normalfont{
  \begin{align} \label{QS7040}
\begin{pmatrix}
\text{the groupoid  of  rank $n$} \\
\text{$\hslash$-flat bundles on $Y^\mr{log}/T^\mr{log}$} 
\end{pmatrix}
\isom \begin{pmatrix}
\text{the groupoid  of } \\
\text{$\hslash$-flat $\mr{GL}_n$-bundles on $Y^\mr{log}/T^\mr{log}$} \\
\end{pmatrix}.
\end{align}}
\eco

 In the following discussions, we  construct  the  bijection desired   in Proposition \ref{QQ97} by examining logarithmic derivations on $\mcE^\mr{log}$ in detail.
 
\subsection{Defining $\mcD er_{\mcV/Y/T}^\mr{log}$} \label{QS7070}
 Let $\mcV$ be as in the statement of Proposition \ref{QQ97}.
 To begin with, we set
 \begin{equation}
 \mcD er_{\mcV/Y/T}
 \end{equation}
to be the $\mcO_Y$-submodule   of $\mcE nd_{f^{-1}(\mcO_T)}(\mcV)$   consisting of local sections $\widetilde{D}$ satisfying that 
 \begin{equation} \label{Der1}
  [\widetilde{D}, \mcE nd_{\mcO_Y}(\mcV)] \subseteq \mcE nd_{\mcO_Y}(\mcV) \ \ \text{and} \ \  [\widetilde{D}, \mr{Im}(\Delta)] \subseteq \mr{Im}(\Delta),
  \end{equation}
where $\Delta$ \index{$\Delta$} denotes the diagonal embedding $\mcO_Y \migiincl \mcE nd_{f^{-1}(\mcO_T)}(\mcV)$.
If $\widetilde{D}$ is a section of $\mcE nd_{f^{-1}(\mcO_T)}(\mcV)$ over an open subscheme $U$ of $Y$, then the  well-defined endomorphism
\begin{align} \label{QW75}
\widetilde{D} |_\Delta : \Delta^{-1} ([\widetilde{D}, \Delta (-)]) : \mcO_U \migi \mcO_U
\end{align}
of $\mcO_U$ specifies a derivation on $\mcO_U$ relative to $T$, i.e., a section of $\mcT_{Y/T}$ over $U$.

Next, we shall  write
\begin{equation} \label{QH1188}
 {\mcD} er_{\mcV/Y/T}^\mr{log} 
 \index{${\mcD} er_{\mcV/Y/T}^\mr{log}$}
\end{equation}
for the Zariski sheaf on $Y$ which,  to any open subscheme $U$ of $Y$, assigns the set of pairs 
\begin{equation} \label{Der4}
\widetilde{\partial}:= (\widetilde{D}, \delta),
\end{equation}
 where  $\widetilde{D}$ is
 a section of $\mcD er_{\mcV/Y/T}$ over $U$ 
   and $\widetilde{\delta}$
   is a monoid homomorphism from $M_Y |_U$ to the underlying additive monoid of
   $\mcO_U$
     satisfying the following two conditions:
 \vspace{1mm}
  \begin{itemize}
  \item
  $\widetilde{D}|_\Delta (\alpha_Y (m))  = \alpha_Y (m) \cdot \delta (m)$ for any local section $m \in M_Y$;
 \vspace{1.5mm}
  \item
    $\delta ({f^*(m)}) =0$ for any local section $m \in M_T$.
    \end{itemize} 
 The sheaf $\mcD er^\mr{log}_{\mcV/Y/T}$ has  a structures of $\mcO_Y$-module, as well as a  Lie bracket operation, obtained  from those of $\mcE nd_{f^{-1}(\mcO_T)}(\mcV)$ via restriction.
 In particular, the Lie bracket operation may be described  as 
 \begin{align}
 [(\widetilde{D}_1, \widetilde{\delta}_1), (\widetilde{D}_2, \widetilde{\delta}_2)] := (\widetilde{D}_1\circ \widetilde{D}_2-   \widetilde{D}_2 \circ \widetilde{D}_1,  \widetilde{D}_1 \circ \widetilde{\delta}_2 - \widetilde{D}_2 \circ    \widetilde{\delta}_1).
 \end{align}
 The assignment $\widetilde{\partial} = (\widetilde{D}, \delta) \mapsto (\widetilde{D} |_\Delta, \delta)$
 determines an $\mcO_Y$-linear morphism 
 \begin{equation}
 d_\mcV^\mr{log} : {\mcD} er_{\mcV/Y/T}^\mr{log} \migi \mcT_{Y^\mr{log}/T^\mr{log}}
 \end{equation}
 compatible with the Lie bracket operations.
 Moreover, 
the pair 
\begin{align} \label{QS7099}
({\mcD} er_{\mcV/Y/T}^\mr{log}, d_\mcV^\mr{log})
\end{align}
 specifies a Lie algebroid on $Y^\mr{log}/T^\mr{log}$ and fits into the following short exact sequence:
 \begin{align} \label{EE334}
0 \longmigi  
\mcE nd_{\mcO_Y}(\mcV)
 \xrightarrow{\widetilde{D} \mapsto (\widetilde{D}, 0)}  {\mcD} er_{\mcV/Y/T}^\mr{log} \xrightarrow{d^\mr{log}_\mcV}  \mcT_{Y^\mr{log}/T^\mr{log}} \longmigi 0. 
\end{align}
(The surjectivity of the anchor map  $d_\mcV^\mr{log}$ can  be  verified immediately. For example, we can  prove  this fact  by considering  the local description of this morphism, as discussed in \S\,\ref{QS7072} later.)


 \subsection{Local construction of  $\widetilde{\mcT}_{\mcE^\mr{log}/T^\mr{log}} \isom \mcD er^\mr{log}_{\mcV/Y/T}$} \label{QS7072}
We shall construct an isomorphism
 $(\widetilde{\mcT}_{\mcE^\mr{log}/T^\mr{log}}, d_\mcE^\mr{log}) \isom (\mcD er^\mr{log}_{\mcV/Y/T}, d_\mcV^\mr{log})$ of Lie algebroids.
 To this end,   we first  consider the case of the trivial vector bundle.
Let  $U$ be an open subscheme of $Y$ and let
 \begin{align} \label{EE333}
  \nabla_0 :  \mcT_{U^\mr{log}/T^\mr{log}} \migi   {\mcD} er_{\mcO_U^{\oplus n}/U/T}^\mr{log} 
 \end{align}
 be the $\mcO_Y$-linear  injection given by  $(D, \delta) \mapsto (D^{\oplus n}, \delta)$.
 Since the equality 
 $d_{\mcO_U^{\oplus n}}^\mr{log} \circ \nabla_0 = \mr{id}_{\mcT_{U^\mr{log}/T^\mr{log}}}$ holds, the morphism
 \begin{align} \label{QQ106}
 (\nabla_0, \mr{incl}) : \mcT_{U^\mr{log}/T^\mr{log}}\oplus  \mcO_U \otimes_k \mfg \mfl_n \migi \mcD er_{\mcO_U^{\oplus n}/U/T}^\mr{log}
 \end{align}
 is an isomorphism, where
 $\mr{incl}$ denotes the inclusion $\mcO_U \otimes_k \mfg \mfl_n \migi \mcD er_{\mcO_U^{\oplus n}/U/T}^\mr{log}$ defined by $\widetilde{D} \mapsto (\widetilde{D}, 0)$
   under the natural identification $\mcO_U \otimes_k \mfg \mfl_n \cong \mcE nd_{\mcO_U}(\mcO_U^{\oplus n})$.
Recall  (cf. (\ref{QQ136})) that the trivial $\mr{GL}_n$-bundle $U \times_k \mr{GL}_n$, i.e., the $\mr{GL}_n$-bundle corresponding to $\mcO_U^{\oplus n}$,
admits a canonical isomorphism
\begin{align}
\tau : \widetilde{\mcT}_{U^\mr{log} \times_k \mr{GL}_n/T^\mr{log}} \isom \mcT_{U^\mr{log}/T^\mr{log}} \oplus \mcO_U \otimes_k \mfg \mfl_n.
\index{$\tau$}
\end{align}


\ble  \label{QQ104}
Let $h$ be an element  of $\mr{GL}_n (U) = \mr{Aut}_{\mcO_U}(\mcO_U^{\oplus n})$.
 Write 
  \begin{align} \label{QS7090}
d_h :   {\mcD} er_{\mcO_U^{\oplus n}/U/T}^\mr{log} \isom  {\mcD} er_{\mcO_U^{\oplus n}/U/T}^\mr{log}
  \end{align}
  for the automorphism  ${\mcD} er_{\mcO_U^{\oplus n}/U/T}^\mr{log}$ given by assigning 
  $(\widetilde{D}, \delta) 
  \mapsto (h \circ \widetilde{D}\circ h^{-1}, \delta)$.
Then, $d_h$ specifies an automorphism of the Lie algebroid $(\mcD er^\mr{log}_{\mcO_U^{\oplus n}/U/T}, d_{\mcO_U^{\oplus n}}^\mr{log})$ and  makes the following diagram commute:
\begin{align} \label{Eh500fg7}
\vcenter{\xymatrix@C=56pt@R=36pt{
\widetilde{\mcT}_{U^\mr{log} \times_k \mr{GL}_n/T^\mr{log}} \ar[r]_-{\sim}^-{ (\nabla_0, \mr{incl}) \circ \tau} \ar[d]^-{\wr}_-{d\mr{L}_h^{\mr{GL}_n}} & \mcD er_{\mcO_U^{\oplus n}/U/T}^\mr{log} \ar[d]_-{\wr}^{d_h}\\
\widetilde{\mcT}_{U^\mr{log} \times_k \mr{GL}_n/T^\mr{log}} \ar[r]^-{\sim}_-{ (\nabla_0, \mr{incl}) \circ \tau}&
\mcD er_{\mcO_U^{\oplus n}/U/T}^\mr{log}.
}}
\end{align}
\ele
\begin{proof}
The former assertion is immediate, so we shall prove  the latter assertion, i.e., the commutativity of (\ref{Eh500fg7}).
Let  $((D, \delta), \widetilde{D})$ be a local section of $\mcT_{U^\mr{log}/T^\mr{log}} \oplus (\mcO_U \otimes_k \mfg \mfl_n)$.
On the one hand, the following sequence of equalities holds:
\begin{align} \label{QS8001}
d_h \circ (\nabla_0, \mr{incl}) (((D, \delta), \widetilde{D})) 
&= 
d_h (\nabla_0 ((D, \delta)) + \mr{incl} (\widetilde{D}))  \\
 & =    
d_h ((D^{\oplus n} + \widetilde{D}, \delta)) \notag \\
 &  =  (h \circ (D^{\oplus n} + \widetilde{D}) \circ h^{-1}, \delta ) \notag \\
 & = (h\circ D^{\oplus n}\circ h^{-1} + h\circ \widetilde{D} \circ h^{-1}, \delta) \notag   \\
 & = (D^{\oplus n} + h \cdot  D (h^{-1})+ h \circ \widetilde{D} \circ h^{-1},  \delta). \notag 
\end{align}
On the other hand, we have
\begin{align} \label{QS8000}
&\ \ \ \  (\nabla_0, \mr{incl}) \circ  \tau \circ d \mr{L}_h^{\mr{GL}_n} \circ\tau^{-1}((D, \delta), \widetilde{D}) \\
&=(\nabla_0, \mr{incl})  ((D, \delta), (h^{-1})^*(\Theta_{\mr{GL}_n})(D) + \mr{Ad}_{\mr{GL}_n}(h)  (\widetilde{D})) \notag \\
& = (\nabla_0, \mr{incl})(((D, \delta), h \cdot  D(h^{-1})  + h \circ  \widetilde{D} \circ h^{-1}) ) \notag \\
& = \nabla_0 ((D, \delta)) + \mr{incl} (h \cdot D (h^{-1}) + h\circ  \widetilde{D} \circ h^{-1}) \notag \\
& = (D^{\oplus n} + h \cdot D (h^{-1}) + h\circ  \widetilde{D} \circ h^{-1}, \delta), \notag \hspace{-9mm}
\end{align}
where the first equality follows from Lemma \ref{QR349} and the second equality follows from the discussion in Example \ref{QQ29}.
By (\ref{QS8001}) and (\ref{QS8000}), 
we obtain  
\begin{align}
d_h \circ (\nabla_0, \mr{incl}) = (\nabla_0, \mr{incl}) \circ\tau \circ d\mr{L}_h^{\mr{GL}_n} \circ \tau^{-1}.
\end{align}
This implies the commutativity of (\ref{Eh500fg7}).
\end{proof}
 
 \subsection{Global construction of $\widetilde{\mcT}_{\mcE^\mr{log}/T^\mr{log}} \isom \mcD er^\mr{log}_{\mcV/Y/T}$}\label{QS8010}
Now, suppose that $\mcV$ may be trivialized on $U$.
 Let us fix a trivialization $\mcV|_U \isom \mcO_U^{\oplus n}$ of $\mcV |_U$, which induces a trivialization $\mcE |_U \isom U \times_k \mr{GL}_n$.
 These trivializations give
 the composite isomorphism
 \begin{align} \label{QH1201}
\Xi_U^{\mcE \rightleftarrows \mcV}  : \widetilde{\mcT}_{\mcE^\mr{log}/T^\mr{log}}|_{U} \isom 
 \widetilde{\mcT}_{U^\mr{log} \times_k \mr{GL}_n/T^\mr{log}}  \xrightarrow{(\nabla_0, \mr{incl})\circ \tau} \mcD er_{\mcO_U^{\oplus n}/U/T}^\mr{log} \isom
 \mcD er_{\mcV/Y/T}^\mr{log} |_U.\hspace{-5mm}
 \end{align}
 Because of the commutativity of (\ref{Eh500fg7}),  this morphism does not  depend on the choice of  $\mcV |_U \isom \mcO_U^{\oplus n}$.
 Hence, the isomorphisms $\Xi_U^{\mcE \rightleftarrows \mcV}$ for various $U$'s (on which $\mcV$ may be trivialized)  may be glued together to obtain an $\mcO_Y$-linear  isomorphism
 \begin{align} \label{QW144}
 \Xi^{\mcE \rightleftarrows \mcV}:  \widetilde{\mcT}_{\mcE^\mr{log}/T^\mr{log}} \isom 
 \mcD er_{\mcV/Y/T}^\mr{log}.
  \index{$\Xi^{\mcE \rightleftarrows \mcV}$}
 \end{align}
It follows from  Lemma \ref{QQ104} that $\Xi^{\mcE \rightleftarrows \mcV}$ forms an isomorphism $(\widetilde{\mcT}_{\mcE^\mr{log}/T^\mr{log}}, d_\mcE^\mr{log}) \isom (\mcD er^\mr{log}_{\mcV/Y/T}, d_\mcV^\mr{log})$ of Lie algebroids.
In particular, $\Xi^{\mcE \rightleftarrows \mcV}$ restricts to an isomorphism $(\mfg \mfl_n)_\mcE \isom \mcE nd_{\mcO_Y}(\mcV)$ of $\mcO_Y$-Lie algebras  fitting  into the following morphism of short exact sequences:
\begin{align} \label{QW30}
\begin{CD}
0 @>>> (\mfg \mfl_n)_\mcE @>>>  \widetilde{\mcT}_{\mcE^\mr{log}/T^\mr{log}}@> d_\mcE^\mr{log}>> \mcT_{Y^\mr{log}/T^\mr{log}} @>>>0
\\
@. @V \wr VV @V \wr V  \Xi^{\mcE \rightleftarrows \mcV} V @V \wr V \mr{id}V @.
\\
0 @>>>  \mcE nd_{\mcO_Y} (\mcV)@>>> \mcD er^\mr{log}_{\mcV/Y/T}@>>d_\mcV^\mr{log}> \mcT_{Y^\mr{log}/T^\mr{log}}@>>>0,
\end{CD}
\end{align}
where the upper  and lower horizontal sequences are  (\ref{Ex0}) and   (\ref{EE334}) respectively.

\begin{rema} \label{QW111}
If $\mr{char} (k) =p>0$, then  $\mcD er^\mr{log}_{\mcV/Y/T}$ admits a structure of restricted Lie algebroid on $Y^\mr{log}/T^\mr{log}$ with $p$-power  operation given by
\begin{align} \label{QW101}
(\widetilde{D}, \delta)^{[p]} := (\widetilde{D}^p, F_Y^* \circ \delta + (\widetilde{D}|_\Delta)^{p-1}\circ \delta)
\end{align}
(cf. ~\cite{Ogu3}, Chap.\,V,\S\,4, Lemma 4.1.5).
One verifies that $\Xi^{\mcE \rightleftarrows \mcV}$ is compatible with the $p$-power operations, i.e., specifies an isomorphism of restricted Lie algebroids.
\end{rema}

\subsection{Correspondence of connections} \label{QS8012}
We now turn to 
finish the proof of Proposition \ref{QQ97}.
We shall construct a bijective correspondence  between  $\hslash$-$T^\mr{log}$-connections  on $\mcV$ and  $\hslash$-$T^\mr{log}$-log connections on $\mcE$.
Let  $\nabla : \mcT \migi  \widetilde{\mcT}_{\mcE^\mr{log}/T^\mr{log}}$  be an $\hslash$-$T^\mr{log}$-connection on $\mcE$.
If we take a local section $\partial  := (D, \delta)\in \mcT$,  then
the local section $\Xi^{\mcE \rightleftarrows \mcV} \circ \nabla (\partial)$ of $\mcD er_{\mcV/Y/T}^\mr{log}$ can be expressed as  $(\widetilde{D},  \delta)$ for some $\widetilde{D} \in \mcE nd_{f^{-1}(\mcO_T)} (\mcV)$.
The morphism $\mcT\migi \mcE nd_{f^{-1}(\mcO_T)} (\mcV)$ given by $\partial \mapsto \widetilde{D}$ induces  an $f^{-1}(\mcO_T)$-linear morphism $\breve{\nabla} : \mcV \migi \Omega \otimes \mcV$, forming   an $\hslash$-$T^\mr{log}$-connection on $\mcV$.

Conversely, let  $\breve{\nabla} : \mcV \migi \Omega \otimes \mcV$ be an $\hslash$-$T^\mr{log}$-connection on $\mcV$.
Then, the morphism $\nabla : \mcT \migi \widetilde{\mcT}_{\mcE^\mr{log}/T^\mr{log}}$ given by 
 $\partial :=(D, \delta) \mapsto (\Xi^{\mcE \rightleftarrows \mcV})^{-1} ((\breve{\nabla}_\partial,  \delta))$ specifies 
 an $\hslash$-$T^\mr{log}$-connection on $\mcE$.

It is immediately verified that the resulting  assignments $\nabla \mapsto \breve{\nabla}$, $\breve{\nabla} \mapsto \nabla$ are inverses of   each other  and  give  the desired  bijective correspondence.
Moreover, since the isomorphism $\Xi^{\mcE \rightleftarrows \mcV}$ is compatible with the Lie bracket operations, 
we have 
\begin{align} \label{QW148}
 \Xi^{\mcE \rightleftarrows \mcV} \circ \psi^{\nabla} = \psi^{\breve{\nabla}}.
\end{align}
It follows that the flatness of $\nabla$ is equivalent to the flatness of $\breve{\nabla}$.
This completes  the proof of Proposition \ref{QQ97}.

\section{$(\mr{GL}_n, \hslash)$-opers and comparison with $(\mfs \mfl_n, \hslash)$-opers}\label{y0107}

Now, we define the notion of a $(\mr{GL}_n, \hslash)$-oper and a certain equivalence relation on the set of $(\mr{GL}_n, \hslash)$-opers. 
As asserted in Proposition \ref{y0112} later, the quotient set of $(\mr{GL}_n, \hslash)$-opers by this equivalence relation corresponds bijectively to the set of    $(\mfs \mfl_n, \hslash)$-opers.

\subsection{} \label{QG800}
Hereafter, 
we suppose that $Y^\mr{log}/T^\mr{log} = U^\mr{log}/S^\mr{log}$, where
$S^\mr{log}$ is 
  an fs  log scheme $S^\mr{log}$ over $k$ and $U^\mr{log}$  is  a log curve over $S^\mr{log}$ (hence $\Omega := \Omega_{U^\mr{log}/S^\mr{log}}$, $\mcT := \mcT_{U^\mr{log}/S^\mr{log}}$, and  $\hslash \in \Gamma (S, \mcO_S)$).
Also, suppose that  $n>1$.

\bde[] \label{y0108} 
\begin{itemize}
 \item[(i)]
 Let us consider   a collection of  data 
 \begin{equation} 
 \label{GL1}
 \msF^\heartsuit := (\mcF, \Nab, \{ \mcF^j \}_{j=0}^n),\end{equation}
where
\begin{itemize}
\item
$\mcF$ is a vector bundle  on $U $ of rank $n$;
\item
$\Nab$ is an $\hslash$-$S^\mr{log}$-connection
 $\mcF \migi \Omega \otimes\mcF$ on $\mcF$;
\item
$\{ \mcF^j \}_{j=0}^n$ is   a decreasing filtration
\begin{equation}
0 = \mcF^n \subseteq \mcF^{n-1} \subseteq \dotsm \subseteq  \mcF^0= \mcF
\end{equation}
 on $\mcF$ consisting of subbundles.
\end{itemize}
Then, we say that $\msF^\heartsuit$ is  a {\bf $(\mr{GL}_n, \hslash)$-oper}  \index{$\mcF^\heartsuit$, $(\mr{GL}_n, \hslash)$-oper}
 on  $U^\mr{log}/S^\mr{log}$ if it    satisfies  the following three  conditions:
\begin{itemize}
\item
For each  $j =0, \cdots,  n-1$, the subquotient $\mcF^j / \mcF^{j+1}$
   is a line bundle;
\item
For each  $j = 1, \cdots,   n-1$, $\nabla (\mcF^j)$ is contained in $\Omega \otimes \mcF^{j-1}$;
\item
For each $j=1, \cdots, n-1$, the well-defined {\it $\mcO_U$-linear} morphism
\begin{equation} \label{GL2}
\mr{KS}^j_{\msF^\heartsuit} :  \mcF^j/\mcF^{j+1} \stackrel{}{\migi}   \Omega \otimes (\mcF^{j-1}/\mcF^j)
\index{$\mr{KS}^j_{\msF^\heartsuit}$}
\end{equation}
defined by $\overline{a} \mapsto \overline{\Nab (a)}$ for any local section $a \in \mcF^j$ (where $\overline{(-)}$'s denote the images in the respective quotients)
  is an isomorphism.
\end{itemize}
If $U^\mr{log}/S^\mr{log}$ arises from
a local pointed stable curve $\msU$ 
  (e.g., $U^\mr{log} = X^\mr{log}$ for a pointed stable curve $\msX := (X/S, \{ \sigma_i \}_i)$), then
we shall refer to any  $(\mr{GL}_n, \hslash)$-oper on $U^\mr{log}/S^\mr{log}$ as a {\bf $(\mr{GL}_n, \hslash)$-oper on $\msU$}.
For simplicity, a $(\mr{GL}_n, 1)$-oper will  be called  a {\bf $\mr{GL}_n$-oper}.

\item[(ii)]
Let $\msF^\heartsuit := (\mcF, \Nab, \{ \mcF^j \}_{j=0}^n)$, $\msF'^\heartsuit := (\mcF', \nabla', \{ \mcF'^j \}_{j=0}^n)$ be $(\mr{GL}_n, \hslash)$-opers on $U^\mr{log}/S^\mr{log}$.
An {\bf  isomorphism of $(\mr{GL}_n, \hslash)$-opers}  from $\msF^\heartsuit$ to $\msF'^\heartsuit$ is an isomorphism
$(\mcF, \nabla) \isom (\mcF', \nabla')$ of $\hslash$-flat bundles   compatible with the respective filtrations $\{ \mcF^j\}_{j}$, $\{ \mcF'^j \}_{j}$.
\end{itemize}  
  \ede

\begin{rema}[Local description] \label{QH366}
Let $\msF^\heartsuit := (\mcF, \nabla, \{ \mcF^j \}_j)$ be a $(\mr{GL}_n, \hslash)$-oper on $U^\mr{log}/S^\mr{log}$.
If we choose a log chart $(V, \partial)$ on $U^\mr{log}/S^\mr{log}$ and a local trivialization $\mcO_V^{\oplus n} \isom \mcF |_V$ of $\mcF$, identifying
$\{ \mcF^j |_V \}_j$ with the standard filtration on  the trivial rank $n$ bundle $\mcO_V^{\oplus n}$, then the element $\nabla^{\mr{GL}_n} (\partial)$ of $\mfg \mfl_n (V)$ can be  expressed as 
\begin{align} \label{QH555}
\nabla^{\mr{GL}_n} (\partial)  =  \begin{pmatrix} * & *& *& \cdots & *& * \\ u_1& * & *& \cdots & * & * \\ 0 & u_2  & *& \cdots & * & * \\ 0 & 0 &u_3& \cdots & *& * \\ \vdots & \vdots & \vdots & \ddots& \vdots & \vdots  \\ 0 & 0 & 0 & \cdots & u_{n-1} & *  \end{pmatrix},
\end{align}
 where the $*$'s indicate arbitrary   functions on $V$ and   $u_1, \cdots, u_{n-1}$ are nowhere vanishing functions, i.e.,  elements in  $\mbG_m (V)$.

\end{rema}

\begin{rema} \label{y01Frh1}
Let $\msF^\heartsuit := (\mcF, \nabla, \{ \mcF^j \}_j)$ be a $(\mr{GL}_n, \hslash)$-oper on $U^\mr{log}/S^\mr{log}$.
 Various  $\mr{KS}_{\msF^\heartsuit}^j$'s together  induce  a  composite isomorphism 
\begin{equation}    \label{QW1909}
\mr{KS}_{\msF^\heartsuit}^{j \Rightarrow n-1} : \mcF^j/\mcF^{j+1} \isom (\mcF^{j+1}/\mcF^{j+2}) \otimes  \mcT_{} \isom \cdots \isom \mcF^{n-1} \otimes  \mcT_{}^{\otimes n-1-j} \end{equation}
($j = 0, \cdots, n-1$). 
Thus, we obtain 
 a composite isomorphism 
\begin{align} \label{QW1906}
 \mr{det}(\mcF) 
 &
 \isom \bigotimes_{j=0}^{n-1}\mcF^j/\mcF^{j+1}
 \\&
  \xrightarrow{\bigotimes_{j=0}^{n-1}\mr{KS}_{\msF^\heartsuit}^{j \Rightarrow n-1}} \bigotimes_{j=0}^{n-1}\mcT_{}^{\otimes (n-1-j)}  \otimes  \mcF^{n-1}
  \notag 
  \\&
   \isom \mcT^{\otimes \frac{n(n-1)}{2}} \otimes (\mcF^{n-1})^{\otimes n}.  \notag
 \end{align}
 In particular, the determinant of the underlying vector bundle  $\mcF$ of $\msF^\heartsuit$ is uniquely  determined, up to isomorphism,  by the $(n-1)$-st step $\mcF^{n-1}$ of the filtration.
\end{rema}

\begin{exa}[$(\mr{GL}_2, \hslash)$-opers]  \label{QH333}
Let us specialize to  the case of $n=2$.
By definition, a $(\mr{GL}_2, \hslash)$-oper is given by a triple $\msF^\heartsuit := (\mcF, \nabla, \mcL)$ consisting of a rank $2$ $\hslash$-flat bundle $(\mcF, \nabla)$ on $U^\mr{log}/S^\mr{log}$ and a line subbundle $\mcL$ of $\mcF$ such that the  {\it $\mcO_U$-linear} composite 
\begin{align} \label{QH334}
\mcL \xrightarrow{\mr{inclusion}} \mcF \xrightarrow{\nabla} \Omega \otimes \mcF \xrightarrow{\mr{quotient}} \Omega \otimes (\mcF/\mcL)
\end{align} 
is an isomorphism.

Now, suppose that $U^\mr{log}/S^\mr{log} = X/k$ for a smooth curve of genus $g$ over $k$ and  $\msF^\heartsuit$ is an {\it $(\mr{SL}_2, \hslash)$-oper}, which means  $(\mr{det}(\mcF), \mr{det}(\nabla)) \cong (\mcO_X, \hslash \cdot d)$. (In general, a $(\mr{GL}_n, \hslash)$-oper with trivial determinant is called an {\bf $(\mr{SL}_n, \hslash)$-oper}.)
Then, 
 we have $\mcO_X \cong \mr{det}(\mcF) \cong \mcL \otimes (\mcF/\mcL) \cong  \mcL^{\otimes 2} \otimes \mcT$, where the last ``$\cong$'' follows from  (\ref{QH334}).
 This implies that $\mcL$  specifies a theta characteristic of $X$ (cf. Example \ref{QQ188}).
 Since $\mcF/\mcL \cong \mr{det}(\mcF) \otimes \mcL^\vee \cong \mcL^\vee$,  the rank $2$ vector bundle $\mcF$ can be obtained as an extension of $\mcL^\vee$ by $\mcL$:
 \begin{align} \label{QH445}
 0 \longmigi \mcL \longmigi \mcF \longmigi \mcL^\vee \longmigi 0.
 \end{align}
One verifies immediately from the construction of (\ref{QH334}) that
the element of $\mr{Ext}^1 (\mcL^\vee, \mcL)$ determined by this extension 
corresponds to the value $u \cdot \hslash \cdot (g-1) \in k$ for some $u \in k^\times$ via the composite of natural isomorphisms
\begin{align}
\mr{Ext}^1 (\mcL^\vee, \mcL) \isom H^1 (X, \Omega) \isom H^0 (X, \mcO_X)^\vee \isom k.
\end{align}  
Note that the factor $u$ comes from the ambiguity in the choice of  the duality isomorphism   $H^1 (X, \Omega) \isom H^0 (X, \mcO_X)^\vee$ in the middle.
In particular,  if $\hslash \cdot (g-1) =0$, then the underlying vector bundle $\mcF$ of the $(\mr{SL}_2, \hslash)$-oper $\msF^\heartsuit$ decomposes into the direct sum $\mcL \oplus \mcL^\vee$.
On the other hand, if $\hslash \cdot (g-1) \neq 0$, then (\ref{QH445}) is a nontrivial (i.e., non-split) extension; such an extension  is uniquely determined up to isomorphism because $\mr{Ext}^1 (\mcL^\vee, \mcL)\cong k$.
\end{exa}

\subsection{Twisting by a flat line bundle}
\label{QS420}
Let $\msF^\heartsuit := (\mcF, \nabla, \{ \mcF^j \}_{j=0}^n)$ be a $(\mr{GL}_n, \hslash)$-oper on $U^\mr{log}/S^\mr{log}$ and $\msL := (\mcL, \nabla_\mcL)$ an $\hslash$-flat line bundle on $U^\mr{log}/S^\mr{log}$.
The tensor product $\mcF \otimes \mcL$  is equipped with 
the decreasing filtration $\{ \mcF^j \otimes \mcL \}_{j =1}^n$ on $\mcF \otimes \mcL$.
Then, the collection of data
\index{$\mcF^\heartsuit$, $(\mr{GL}_n, \hslash)$-oper@--- $\msF^\heartsuit_{\otimes \msL}$}
\begin{equation} \label{QQ111}
  \msF^\heartsuit_{\otimes \msL} := (\mcF\otimes \mcL, \nabla \otimes \nabla_\mcL, \{ \mcF^j \otimes \mcL \}_{j=0}^n)
  \end{equation}
(cf. (\ref{QQ295}) for the definition of $\nabla \otimes \nabla_\mcL$)
forms  a $(\mr{GL}_n, \hslash)$-oper on $U^\mr{log}/S^\mr{log}$.

\bde[] \label{y0109}
 Let 
$\msF^\heartsuit :=(\mcF, \nabla, \{ \mcF^j \}_{j=0}^n)$,
$\mcF'^\heartsuit :=(\mcF', \nabla', \{ \mcF'^j \}_{j=0}^n)$
be $(\mr{GL}_n, \hslash)$-opers on $U^\mr{log}/S^\mr{log}$.
We shall say that $\msF^\heartsuit$ is {\bf equivalent to 
$\msF'^\heartsuit$}
\index{$\mcF^\heartsuit$, $(\mr{GL}_n, \hslash)$-oper@--- $\stackrel{\heartsuit}{\sim}$, equivalent}
if there exists an $\hslash$-flat line bundle $\msL := (\mcL, \nabla_\mcL)$ on $U^\mr{log}/S^\mr{log}$ such that the $(\mr{GL}_n, \hslash)$-oper
$\msF^\heartsuit_{\otimes \msL}$ is isomorphic  to $\msF'^\heartsuit$. 
We use the notation  ``$\msF^\heartsuit \stackrel{\heartsuit}{\sim}\msF'^\heartsuit$''  to indicate the situation that $\msF^\heartsuit$ is  equivalent to $\msF'^\heartsuit$.
(It is immediate that the binary  relation ``$\stackrel{\heartsuit}{\sim}$''  on the set of $(\mr{GL}_n, \hslash)$-opers on $U^\mr{log}/S^\mr{log}$ 
defines  an equivalence relation.)
  \ede

\subsection{Pull-back and base-change}\label{QS421}
Let  $\msF^\heartsuit := (\mcF, \nabla, \{ \mcF^j \}_{j=0}^n)$ be   a $(\mr{GL}_n, \hslash)$-oper on $U^\mr{log}/S^\mr{log}$.
Suppose that we are given 
another  log curve $U'^{\mr{log}}$ over $S^\mr{log}$
 equipped with a {\it log \'{e}tale}   morphism $u^\mr{log} : U'^{\mr{log}} \migi U^\mr{log}$.
Then, one verifies that
 the pull-backs    of  various data
 \index{$\mcF^\heartsuit$, $(\mr{GL}_n, \hslash)$-oper@--- $u^*(\msF^\heartsuit)$, pull-back of}
\begin{equation} \label{QQ114}
u^*(\msF^\heartsuit) := (u^*(\mcF), u^*(\nabla), \{ u^*(\mcF^j) \}_{j=0}^n)
\end{equation}
 defining  $\msF^\heartsuit$   forms a $(\mr{GL}_n, \hslash)$-oper on $U'^{\mr{log}}/S^\mr{log}$.

Also, if $s: S'^{\mr{log}} \migi S^\mr{log}$ is a morphism of fs log schemes over $k$,
then the base-changes 
\index{$\mcF^\heartsuit$, $(\mr{GL}_n, \hslash)$-oper@--- $s^*(\msF^\heartsuit)$, base-change of}
\begin{align} \label{QS8303}
s^*(\msF^\heartsuit) := (s^*(\mcF), s^*(\nabla), \{ s^*(\mcF^j) \}_{j=0}^n)
\end{align}
to  $S'^{\mr{log}}$ 
 specifies 
a $(\mr{GL}_n, s^*(\hslash))$-oper on the log curve $S'^{\mr{log}} \times_{S^\mr{log}} U^\mr{log}/S'^{\mr{log}}$.
The formations  of pull-back and base-change are  compatible with 
 the equivalence relation  defined in Definition \ref{y0109} above.

\subsection{The sheave of $\mfs \mfl_n$-opers and   $(\mr{GL}_n, \hslash)$-opers} \label{QS422}
Until the end of \S\,\ref{y0131}, we suppose that
{\it either ``\,$\mr{char}(k) =0$'' or ``\,$\mr{char}(k) > n$'' 
is fulfilled.}
In particular, under this assumption,  the Lie  algebra $\mfp \mfg \mfl_n$ of $\mr{PGL}_n$ is isomorphic to  $\mfs \mfl_n$ via the composite of natural    morphisms $\mfs \mfl_n \migiincl \mfg \mfl_n \migisurj\mfp \mfg \mfl_n$.
For simplicity, we shall write
 \index{$\mcO p^\spadesuit_{n, \hslash (, \rho)}$, $\mcO p^{\spadesuit, \mr{O}}_{n, \hslash}$, $\mcO p^{\spadesuit, \mr{Sp}}_{n, \hslash}$}
\begin{align} \label{QH999}
\mcO p^\spadesuit_{n, \hslash} 
 := \mcO p_{\mfs \mfl_n, \hslash, U^\mr{log}/S^\mr{log}},
\end{align}
i.e., the set-valued sheaf on  $\mfE \mft_{/U}$ ($:=$ the small \'{e}tale  on $U$) which, to any \'{e}tale $U$-scheme $U'$, assigns the set of isomorphism classes of $(\mfs \mfl_n, \hslash)$ on $U'^{\mr{log}} (= U' \times_{U'}U^\mr{log})/S^\mr{log}$ (cf. (\ref{opersite})).
Also, denote by 
\begin{equation}   \label{QQ113}
\overline{\mcO} p^\heartsuit_{n, \hslash} 
  : \mfE \mft_{/U} \migi \mfS \mfe \mft  \index{$\overline{\mcO} p_{n, \hslash, U^\mr{log}/S^\mr{log}}$}
  \index{$\overline{\mcO} p^\heartsuit_{n, \hslash}$}
\end{equation}
the set-valued sheaf on $\mfE \mft_{/U}$
 associated to the presheaf which, to any \'{e}tale $U$-scheme $U'$, assigns the set of {\it equivalence classes} of $(\mr{GL}_n, \hslash)$-opers on $U'^{\mr{log}}/S^\mr{log}$.

\subsection{Rank $n$ flat  bundles and flat $\mr{PGL}_n$-bundles} \label{QG801}

The following lemma will be elemental but shows  how to apply  the assumption imposed  at the beginning of \S\,\ref{QS422}
 to our discussion.

\bpr \label{y0111}
\begin{itemize}
\item[(i)]
Let  $\mcL$ be a line bundle on $U$.
Assume  that  we are given an $\hslash$-$S^\mr{log}$-connection $\nabla_{\mcL^{\otimes n}}$
on the $n$-fold tensor product $\mcL^{\otimes n}$ of $\mcL$.
Then, there exists a {\it unique}  $\hslash$-$S^\mr{log}$-connection $\nabla$ on $\mcL$ satisfying  $\nabla^{\otimes n} = \nabla_{\mcL^{\otimes n}}$ (cf. (\ref{QQ296})). 
\item[(ii)]
Let $\mcE_{\mr{GL}_n}$ be a $\mr{GL}_n$-bundle on $U$.
Assume  that we are given 
a pair $(\nabla_{\mr{PGL}_n}, \nabla_{\mbG_m})$, where 
\begin{itemize}
\item
$\nabla_{\mr{PGL}_n}$ denotes an $\hslash$-$S^\mr{log}$-connection  on the $\mr{PGL}_n$-bundle $\mcE_{\mr{PGL}_n} := \mcE_{\mr{GL}_n} \times^{\mr{GL}_n} \mr{PGL}_n$ induced by $\mcE_{\mr{GL}_n}$ via the natural quotient $\mu_{\mr{PGL}_n} :\mr{GL}_n \migisurj  \mr{PGL}_n$;
\item
$\nabla_{\mbG_m}$ denotes an $\hslash$-$S^\mr{log}$-connection  on the $\mbG_m$-bundle  $\mcE_{\mbG_m} :=\mcE_{\mr{GL}_n} \times^{\mr{GL}_n} \mbG_m$ induced by $\mcE_{\mr{GL}_n}$ via the determinant map $\mu_{\mbG_m} :\mr{GL}_n \migi \mbG_m$.
\end{itemize}
Then, there exists a unique  $\hslash$-$S^\mr{log}$-connection  $\nabla_{\mr{GL}_n}$ on $\mcE_{\mr{GL}_n}$
satisfying  the equalities
\begin{align}
\nabla_{\mr{GL}_n} \times^{\mr{GL}_n, \mu_{\mr{PGL}_n}} \mr{PGL}_n = \nabla_{\mr{PGL}_n},
\hspace{5mm}
\nabla_{\mr{GL}_n} \times^{\mr{GL}_n, \mu_{\mbG_m}} \mbG_m = \nabla_{\mbG_m}.
\end{align}
\end{itemize}
  \epr
\begin{proof}
The assumption ``$\mr{char}(k) =0$ or $\mr{char}(k)>n$'' implies the \'{e}taleness of 
the $n$-th power morphism $\mcL^\times \migi (\mcL^{\otimes n})^\times$ (resp., the morphism $(\mu_{\mr{PGL}_n}, \mu_{\mbG_m}) : \mr{GL}_n \migi \mr{PGL}_n \times_k \mbG_m$), where $\mcL^\times$ and $(\mcL^{\otimes n})^\times$ denote the $\mbG_m$-bundles associated to $\mcL$ and $\mcL^{\otimes n}$ respectively.
It follows that 
 the morphism
\begin{align} \label{GL3}
& \hspace{15mm} \widetilde{\mcT}_{(\mcL^\times)^\mr{log}/S^\mr{log}} \migi \widetilde{\mcT}_{(\mcL^{\otimes n})^{\times \mr{log}}/S^\mr{log}} \\
& \left(\text{resp.}, \  \widetilde{\mcT}_{\mcE_{\mr{GL}_n}^\mr{log}/S^\mr{log}} \migi \widetilde{\mcT}_{\mcE_{\mr{PGL}_n}^\mr{log}/S^\mr{log}}  \times_{\mcT_{U^\mr{log}/S^\mr{log}}} \widetilde{\mcT}_{\mcE_{\mbG_m}^\mr{log}/S^\mr{log}} \right)\notag
\end{align}
 induced by  this morphism is an isomorphism of Lie algebroids.
Hence, 
 assertion (i) (resp., assertion (ii)) follows from this fact and the definition of an $\hslash$-$S^\mr{log}$-connection.
(In the non-resp'd portion, we use the result of Proposition \ref{QQ97}.)
\end{proof}

\bco \label{QW1800}
Let $\mcV$ be a rank $n$ vector bundle  on $U$ and denote by $\mcE$ the $\mr{PGL}_n$-bundle induced by $\mcV$ via projectivization.
Also,  let   $\nabla$ and $\nabla'$ be  $S^\mr{log}$-connections on $\mcV$
such that the $\hslash$-$S^\mr{log}$-connections on $\mcE$ induced by them coincide.
Then, there exists uniquely 
an $\hslash$-$S^\mr{log}$-connection $\nabla_\mcO$ on $\mcO_U$ with 
$\nabla = \nabla' \otimes \nabla_\mcO$.
\eco
\begin{proof}
We shall regard the tensor product $\mr{det}(\nabla) \otimes \mr{det}(\nabla')^\vee$ as an $\hslash$-$S^\mr{log}$-connection on $\mcO_U$ via the natural identification $\mr{det}(\mcV) \otimes \mr{det}(\mcV)^\vee = \mcO_U$.
By Proposition \ref{y0111}, (i),  there exists uniquely
an $\hslash$-$S^\mr{log}$-connection $\nabla_\mcO$ on $\mcO$
with $\nabla_\mcO^{\otimes n} = \mr{det}(\nabla) \otimes \mr{det}(\nabla')^\vee$.
Then, we have
\begin{align}
\mr{det}(\nabla' \otimes \nabla_\mcO) = \mr{det}(\nabla') \otimes \nabla_{\mcO}^{\otimes n} = \mr{det}(\nabla')  \otimes (\mr{det}(\nabla) \otimes \mr{det}(\nabla')^\vee) = \mr{det}(\nabla).
\end{align}
This shows that  $\nabla$ and $\nabla' \otimes \nabla_\mcO$ coincide after projectivization, as well as  taking determinants.
Hence, it follows from the uniqueness portion of Proposition \ref{y0111}, (ii), that
$\nabla = \nabla' \otimes \nabla_\mcO$, thus completing 
the  proof of the assertion.
\end{proof}

\subsection{The $(\mfs \mfl_n, \hslash)$-oper associated to a $(\mr{GL}_n, \hslash)$-oper} \label{QS8645}
We shall construct an isomorphism of sheaves from 
  $\overline{\mcO} p_{n, \hslash}^\heartsuit$ to 
  $\mcO p^\spadesuit_{n, \hslash}$.
Let $u : U' \migi U$ be an \'{e}tale scheme over $U$  and  $\msF^\heartsuit = (\mcF, \nabla, \{ \mcF^j \}_{j=0}^n)$  a $(\mr{GL}_n, \hslash)$-oper on $U'^{\mr{log}}/S^\mr{log}$.
Denote by $(\mcF_{\mr{PGL}_n}, \nabla_{\mr{PGL}_n})$ the $\hslash$-flat $\mr{PGL}_n$-bundle on $U'^{\mr{log}}/S^\mr{log}$ induced  from $(\mcF, \nabla)$ via projectivization.
The filtration $\{ \mcF^i \}_{i =0}^n$ determines  a $B$-reduction $\mcF_B$ of $\mcF_{\mr{PGL}_n}$ (cf. (\ref{BorelB}) for the definition of $B$).
By the definition of a $(\mr{GL}_n, \hslash)$-oper, the pair 
\index{$\mcF^\heartsuit$, $(\mr{GL}_n, \hslash)$-oper@--- $\mcF^{\heartsuit \Rightarrow \spadesuit}$}
\begin{align} \label{QQ140}
\msF^{\heartsuit \Rightarrow \spadesuit} := (\mcF_B, \nabla_{\mr{PGL}_n})
\end{align}
 constructed  in this way forms an $(\mfs \mfl_n, \hslash)$-oper on $U'^{\mr{log}}/S^\mr{log}$ (cf. (\ref{QQ191}) and (\ref{QH555})).
The isomorphism class $[\msF^{\heartsuit \Rightarrow \spadesuit}]$
 of 
 $\msF^{\heartsuit \Rightarrow \spadesuit}$
  does not depend on the choice of a representative in the equivalence class defining $\msF^\heartsuit$.
Thus, the resulting 
assignment $\msF^\heartsuit \mapsto [\msF^{\heartsuit \Rightarrow \spadesuit}]$ 
is well-defined and gives  a morphism of sheaves
\index{$\Lambda^{\heartsuit \Rightarrow \spadesuit}_{(\rho)}$, $\Lambda^{\heartsuit \Rightarrow \spadesuit, \BL}$, $\Lambda^{\heartsuit \Rightarrow \spadesuit, \BL}_\msX$}
\begin{equation}   \label{QQ118}
\Lambda^{\heartsuit \Rightarrow \spadesuit}_{}  : \overline{\mcO} p^\heartsuit_{n, \hslash}    \migi  \mcO p^\spadesuit_{n, \hslash}.
 \end{equation}
The formation of this  morphism commutes, in an evident sense,  with pull-back to fs log schemes over $S^\mr{log}$.


\bpr \label{y0112} 
$\Lambda^{\heartsuit \Rightarrow \spadesuit}$ is an isomorphism.
\epr
\begin{proof}
By the  functorial formations of 
$\overline{\mcO} p^\heartsuit_{n, \hslash}$, $\mcO p^\spadesuit_{n, \hslash}$, and $\Lambda^{\heartsuit \Rightarrow \spadesuit}_{}$,
it suffices to prove the bijectivity of the induced map of sets $\Lambda^{\heartsuit \Rightarrow \spadesuit} (U)  : \overline{\mcO} p^\heartsuit_{n, \hslash}(U)    \migi  \mcO p^\spadesuit_{n, \hslash}(U)$.

First, let us consider  the surjectivity  of $\Lambda_{}^{\heartsuit \Rightarrow  \spadesuit}(U)$.
Let $\msE^\spadesuit := (\mcE_B, \nabla)$ be an $(\mfs \mfl_n, \hslash)$-oper on $U^\mr{log}/S^\mr{log}$.
Then, the  problem is to find locally a $(\mr{GL}_n, \hslash)$-oper $\msF^\heartsuit$ with $\msF^{\heartsuit \Rightarrow \spadesuit} \cong \msE^\spadesuit$.
Note that $\mcE_B$ comes, \'{e}tale locally on $U$, from a rank $n$ vector bundle $\mcF$ together with a complete flag $\{ \mcF^j \}_{j=0}^n$ on $\mcF$ via change of  structure group $\mr{GL}_n \migisurj \mr{PGL}_n$.
Hence, in order to solve the problem, we can assume, after possibly replacing $U$ with its covering in the \'{e}tale topology,  that both  $\mcF$ and  $\{ \mcF^j \}_{j=0}^n$ can be defined over the entire space $U$ and that the determinant $\mr{det}(\mcF)$ is trivial. 
According to  Corollary  \ref{QS8702} and Proposition  \ref{y0111}, (ii), there exists an $\hslash$-$S^\mr{log}$-connection $\nabla_\mcF$ on $\mcF$ which induces  $\nabla$ on $\mcE_{\mr{PGL}_n}:= \mcE_B \times^B \mr{PGL}_n$ and induces the trivial connection $\hslash \cdot d$ on $\mcO_U \left(=\mr{det}(\mcF)\right)$.
By construction, the collection of data $(\mcF, \nabla_\mcF, \{ \mcF^j \}_{j})$ forms a  desired $(\mr{GL}_n, \hslash)$-oper.
This completes the surjectivity of  $\Lambda^{\heartsuit \Rightarrow  \spadesuit}(U)$.

Next, let us consider the  injectivity of $\Lambda_{}^{\heartsuit \Rightarrow  \spadesuit}(U)$.
Let  $\msF^\heartsuit := (\mcF, \nabla, \{\mcF^j\}_{j})$,  $\msF'^\heartsuit := (\mcF', \nabla', \{ \mcF^j \}_{j})$ be $(\mr{GL}_n, \hslash)$-opers on $U^\mr{log}/S^\mr{log}$ such that  there exists an isomorphism
 \begin{equation} h_{\mr{PGL}} : \msF^{\heartsuit \Rightarrow \spadesuit} \isom \msF'^{\heartsuit \Rightarrow \spadesuit}
 \end{equation}
 of $(\mfs \mfl_n, \hslash)$-opers.
Note that  we can find  a line bundle $\mcL$ on $U$ and an isomorphism $h_\mr{GL} : \mcF \isom \mcF' \otimes \mcL$ of vector bundles inducing $h_{\mr{PGL}}$  via projectivization   (cf. ~\cite{Har}, Chap.\,II, \S\,7, Exercise 7.9).
The determinant of this isomorphism specifies an isomorphism
\begin{equation}
h_{\mbG_m} : \mr{det} (\mcF) \isom \left(\mr{det}(\mcF' \otimes \mcL)=\right) 
\mr{det}(\mcF') \otimes \mcL^{\otimes n}
\end{equation}
between line bundles.
Let  $\nabla_{\mcL^{\otimes n}}$  be 
the $\hslash$-$S^\mr{log}$-connection on $\mcL^{\otimes n}$
such that $\nabla$ corresponds to $\nabla' \otimes \nabla_{\mcL^{\otimes n}}$ via $h_{\mbG_m}$.
By  Proposition \ref{y0111}, (i),  there exists a unique $\hslash$-$S^\mr{log}$-connection $\nabla_\mcL$ on $\mcL$ with $\nabla^{\otimes n}_{\mcL} = \nabla_{\mcL^{\otimes n}}$.
Then, $h_{\mbG_m}$ becomes an isomorphism $(\mr{det}(\mcF), \mr{det}(\nabla))\isom (\mr{det}(\mcF' \otimes \mcL), \mr{det}(\nabla' \otimes \nabla_\mcL))$ between $\hslash$-flat line bundles.
It follows  that $h_\mr{GL}$ is compatible with
the $\hslash$-$S^\mr{log}$-connections $\nabla$, $\nabla' \otimes \nabla_\mcL$ after changing the structure group by the natural quotient $\mr{GL}_n \migisurj \mr{PGL}_n$, as well as the determinant map $\mr{GL}_n \migisurj \mbG_m$.
Because of this fact and the uniqueness portion of Proposition \ref{y0111}, (ii),
$h_\mr{GL}$ is in fact compatible with $\nabla$ and  $\nabla' \otimes \nabla_\mcL$, i.e., defines an isomorphism of 
$\hslash$-flat bundles.
Finally, 
since $h_{\mr{PGL}}$ preserves the  $B$-reductions, 
$h_\mr{GL}$ preserves the  filtrations $\{ \mcF^j \}_{j}$,  $\{ \mcF'^j \otimes \mcL \}_{j}$.
This implies that $h_\mr{GL}$ specifies an isomorphism $\msF^{\heartsuit} \isom  \msF'^{\heartsuit}_{\otimes \msL}$,  so we have  $\msF^\heartsuit \stackrel{\heartsuit}{\sim} \msF'^\heartsuit$.
 This proves  the  injectivity of $\Lambda_{}^{\heartsuit \Rightarrow \spadesuit}(U)$, and hence, completes the proof of this proposition.
\end{proof}

\bpr \label{QW1045}
Let $\msF^\heartsuit :=(\mcF, \nabla, \{ \mcF^j\}_j)$ be a $(\mr{GL}_n, \hslash)$-oper on $U^\mr{log}/S^\mr{log}$.
Then, each automorphism of $\msF^\heartsuit$ may be given by $\Delta (a)$ (cf. (\ref{QQ77})) for a unique global section   $a$ of $\mcO_U^\times \cap \mr{Ker}(\mcO_U \xrightarrow{\hslash \cdot d} \Omega)$.
\epr
\begin{proof}
Let $\eta$ be an automorphism of $\msF^\heartsuit$.
According to  Proposition \ref{y040},  the automorphism of $\msF^{\heartsuit \Rightarrow \spadesuit}$ induced by $\eta$ coincides with the identity morphism.
This implies that $\eta$ may be given by $\Delta (a)$ for some  $a \in \Gamma (U, \mcO_U^\times)$.
Such a section $a$ is uniquely determined because of the local freeness of $\mcF$.
For every local section $v \in \mcF$, 
we have
 \begin{align}
 a \cdot \nabla (v) = (\mr{id}_\Omega \circ \eta) (\nabla (v)) = \nabla (\eta (v)) = \nabla (a \cdot v)= \hslash \cdot  da \otimes v + a \cdot \nabla (v).
 \end{align}
Hence,  the section $a$ satisfies   $\hslash \cdot d a =0$, i.e., lies in $\mr{Ker}(\hslash \cdot d)$.
This completes the proof of the assertion.
\end{proof}

\section{$(n, \hslash, \BB)$-projective connections and $(\mr{GL}_n, \hslash, \BB)$-opers} \label{y0116}

\subsection{$(n, \hslash, \BB)$-projective connections}\label{QS9000}

Let us fix a line bundle $\BB$ on $U$.
For $j \in \mbZ$, 
we shall set
\index{${^\Join}\DD_\BB^{\,<j}$, $\JJ\DD_\BB^{\, <j}$, $\JJ\DD_\BB^{\, \sharp <n}$}
\begin{align}  
& \hspace{11mm} 
 {^\Join}\DD_\BB^{\,<j} := \mcH om_{\mcO_Y} (\BB^\vee,  (\Omega_{}^{\otimes n} \otimes \BB^\vee) \otimes \mcD_{\hslash}^{< j}) \\
&  \left(\text{resp.,} \ 
\JJ\DD_\BB^{\, <j}  := \mcH om_{\mcO_Y} (\mcT_{}^{\otimes n} \otimes \BB,  \mcD_{\hslash}^{< j}\otimes \BB)\right). \notag 
\end{align}
(According to the discussion in Remark \ref{QQ425}, each section of ${^\Join}\DD_\BB^{\,<j}$ may be thought of  as an ``{\it $\hslash$-twisted version}'' of   a  differential operator $\BB^\vee \migi \Omega^{\otimes n} \otimes \BB^\vee$.)
If $j >0$, then the canonical  isomorphism $\mcD^{< j}_{\hslash}/\mcD^{< j-1}_{\hslash} \isom  \mcT_{}^{\otimes j-1}$ (cf. (\ref{GL6}))
yields an isomorphism
\begin{equation}
\label{GL8}
{^\Join}\DD_\BB^{\,<j}/{^\Join}\DD_\BB^{\,<j-1}
 \isom \Omega_{}^{\otimes (n-j+1)} \  
  \left(\text{resp.}, \ 
  \JJ\DD_\BB^{\, <j}
  / \, 
  \JJ\DD_\BB^{\, <j-1}
  \isom \Omega_{}^{\otimes (n-j+1)} \right). \hspace{-5mm}
  \end{equation}
By composing it  with the quotient
$ {^\Join}\DD_\BB^{\, <j}\migisurj  {^\Join}\DD_\BB^{\, <j}/ {^\Join}\DD_\BB^{\, <j-1}$ (resp., $\JJ\DD_\BB^{\, <j} \migisurj \JJ\DD_\BB^{\, <j}/\, \JJ\DD_\BB^{\, <j-1}$),
we obtain a surjection
\begin{equation}
\label{GL9}
{^\Join} \Sigma^{j} : 
 {^\Join}\DD_\BB^{\, <j} 
 \migisurj \Omega_{}^{\otimes (n-j+1)}
  \ \left(\text{resp.,} \ 
  \JJ\Sigma^j : 
   \JJ\DD_\BB^{\, <j}
   \migisurj \Omega_{}^{\otimes (n-j+1)}\right).
\end{equation}
If  $D$ is a local section of ${^\Join}\DD_\BB^{\, <j}$, then ${^\Join} \Sigma^{j}(D)$ is called the {\bf principal symbol} \index{principal symbol} of $D$.

Let us consider the composite isomorphism 
\begin{equation} 
\label{GL12}
{^\circlearrowright \beta}^j :  
{^\Join}\DD^{\, <j}_{\BB} \isom   \Omega_{}^{\otimes n} \otimes \BB^\vee \otimes \mcD_{\hslash}^{< j} \otimes \BB \isom 
\JJ\DD^{\,<j}_{\BB}
\index{${^\circlearrowright \beta}^j$}
\end{equation}
of $\mcO_U$-bimodules, where the first and second isomorphisms follow from the definitions of ${^\Join}\DD^{\, <j}_{\BB}$ and $\JJ\DD^{\,<j}_{\BB}$ respectively.
This composite fits into the following commutative diagram:
\begin{align} \label{QQ240}
\vcenter{\xymatrix@C=36pt@R=36pt{
{^\Join}\DD^{\, <j}_{\BB}  \ar[rr]_-{\sim}^-{{^\circlearrowright \beta}^j} \ar[rd]_-{{^\Join} \Sigma^j}&&\JJ\DD^{\,<j}_{\BB}  \ar[ld]^-{ \JJ\Sigma^j}\\
& \Omega^{\otimes (n-j+1)}.&
}}
\end{align}
Moreover, for each $j' \leq j$,  the equality
   ${^\circlearrowright \beta}^j |_{{^\Join}\DD^{\, <j'}_\BB } = {^\circlearrowright \beta}^{j'}$
 holds.

Now, we shall write
\index{$\DD_{n, 1, \bb}^{\,\clubsuit \mr{full}}$, $\mfD \mfi \mff \mff^{\clubsuit \mr{full}}_{n, 1, \bb  (, \rho)}$}
\begin{equation}   \label{QS9010}
\mcD\text{{\it  iff}}_{n, \hslash, \BB}^{\, \clubsuit}
 : = ({^\Join }\Sigma^{n+1})^{-1}(1) \ \left(\subseteq {^\Join}\DD^{\,<n+1}_\BB\right).   
 \end{equation}
It is clear that this sheaf  is nonempty
   and
admits a natural ($\mr{Ker}({^\Join}\Sigma^{n+1})=$) ${^\Join}\DD^{\,<n}_{\BB}$-torsor structure.
In the natural manner,  we consider  $\mcD\text{{\it  iff}}_{n, \hslash, \BB}^{\, \clubsuit}$
as a $\mfS \mfe \mft$-valued sheaf on $\mfE \mft_{/U}$.
 
\bde \label{QH01}
We shall refer to  any element of $\Gamma (U, \mcD\text{{\it  iff}}_{n, \hslash, \BB}^{\, \clubsuit})$ as an {\bf $(n, \hslash, \BB)$-projective connection} \index{$D^\clubsuit$, $(n, \hslash, \BB)$-projective connection} on $U^\mr{log}/S^\mr{log}$.
\ede

\subsection{$(\mr{GL}_n, \hslash, \BB)$-opers} \label{QG810}
We shall write 
\begin{align} \label{QW47}
\DF_{\BB} := \mcD_\hslash^{<n} \otimes \BB, \hspace{5mm}
\DF_{\BB}^j  := \mcD_\hslash^{<n-j} \otimes \BB \ \ \  (j=0,1, \cdots, n).
\index{$\DF_{\BB}$, $\DF_{\BB}^j$}
\end{align}
Here, recall (cf. the comment preceding Remark \ref{QH1000})
that
an $\mcO_U$-module of the form $\mcD_\hslash^{<j} \otimes \mcV$ for  $j \in \mbZ \sqcup \{ \infty \}$ and a vector bundle $\mcV$
  is, to be exact,
defined as   ${^R}\mcD_\hslash^{<j} \otimes \mcV$ equipped with 
 a  structure of  $\mcO_U$-module
   given by left  multiplication on $\mcD_\hslash^{<j}$.
The subquotient $\mr{gr}^j (\DF_\BB) := \DF_\BB^j /\DF_\BB^{j+1}$ (for each $j = 0, \cdots, n-1$) is naturally isomorphic to $\mcT^{\otimes n-j-1} \otimes \BB$.


\bde[] \label{y0257}
\begin{itemize}
\item[(i)]
A {\bf $(\mr{GL}_n, \hslash, \BB)$-oper} 
\index{$\nabla^\diamondsuit$, $(\mr{GL}_n, \hslash, \BB)$-oper}
 on $U^\mr{log}/S^\mr{log}$ is 
an $\hslash$-$S^\mr{log}$-connection $\nabla^\diamondsuit$ on 
$\DF_\BB$
such that 
the collection of data
\index{$\nabla^\diamondsuit$, $(\mr{GL}_n, \hslash, \BB)$-oper@--- $\nabla^{\diamondsuit \Rightarrow \heartsuit}$}
\begin{equation}  \label{rrr7}
\nabla^{\diamondsuit \Rightarrow \heartsuit}  := (\DF_\BB,
 \nabla^\diamondsuit, \{ \DF^j_\BB
 \}_{j=0}^{n}) 
 \end{equation}
forms a $(\mr{GL}_n, \hslash)$-oper on $U^\mr{log}/S^\mr{log}$.
\item[(ii)]
Let $\nabla^\diamondsuit$, $\nabla'^\diamondsuit$ be $(\mr{GL}_n, \hslash, \BB)$-opers on $U^\mr{log}/S^\mr{log}$.
Then, we shall say that $\nabla^\diamondsuit$ is {\bf  isomorphic to $\nabla'^\diamondsuit$} if the associated $(\mr{GL}_n, \hslash)$-opers $\nabla^{\diamondsuit \Rightarrow\heartsuit}$ and $\nabla'^{\diamondsuit \Rightarrow \heartsuit}$  are isomorphic.
\end{itemize}
  \ede

Let  $u : U' \migi U$ be  an \'{e}tale morphism and $\nabla^\diamondsuit$  a $(\mr{GL}_n, \hslash, \BB)$-oper on $U^\mr{log}/S^\mr{log}$.
Because of the \'{e}taleness of $u$,
the differential of this morphism induces an identification $\mcD^{<j}_{\hslash, U'^{\mr{log}}/S^\mr{log}} \otimes u^*(\BB) \isom  u^{*}(\mcD^{<j}_{\hslash, U^{\mr{log}}/S^\mr{log}} \otimes \BB)$ for every $j$.
Under this identification,   the pull-back  $u^*(\nabla^\diamondsuit)$ (cf. (\ref{QS8300})) specifies  a  $(\mr{GL}_n, \hslash, u^*(\BB))$-oper on $U'^{\mr{log}}/S^\mr{log}$.
Hence, 
we obtain 
 an $\mfS \mfe \mft$-valued  presheaf
 \index{$\mcO p_{n, \hslash, \BB}^\diamondsuit$, $\mcO p_{n, \hslash, \bb (, \rho)}^\diamondsuit$, $\mcO p_{n, \hslash, \BB}^{\diamondsuit, \mr{O}}$, $\mcO p_{n, \hslash, \BB}^{\diamondsuit, \mr{Sp}}$}
\begin{equation}   \label{rrr9}
\mcO p_{n, \hslash, \BB}^\diamondsuit : \mfE \mft_{/U} \migi \mfS \mfe \mft 
\end{equation}
on $ \mfE \mft_{/U} $  which, to any \'{e}tale $U$-scheme $u : U' \migi U$, assigns the set of isomorphism classes of $(\mr{GL}_n, \hslash, u^*(\BB))$-opers on $U'^{\mr{log}}/S^\mr{log}$.  

\subsection{The  $(\mr{GL}_n, \hslash, \BB)$-oper associated to an $(n, \hslash, \BB)$-projective connection} \label{QS9400}


Let 
 $D^\clubsuit : \BB^\vee \migi  \Omega_{}^{\otimes n} \otimes \BB^\vee \otimes \mcD_{\hslash}^{< n+1}$ be an $(n, \hslash, \BB)$-projective connection on $U^\mr{log}/S^\mr{log}$.
It gives a left $\mcD^{< \infty}_{\hslash}$-module
$\mcD^{<\infty}_{\hslash}\otimes \BB/\langle \mr{Im}({^\circlearrowright \beta}^{n+1}  (D^\clubsuit))\rangle$, 
i.e., the quotient 
of $\mcD^{<\infty}_{\hslash}\otimes \BB$ by the $\mcD^{<\infty}_{\hslash}$-submodule generated by the image 
 of the morphism 
\begin{equation}
{^\circlearrowright \beta}^{n+1}  (D^\clubsuit) : \mcT_{}^{\otimes n} \otimes \BB\migi \mcD_{\hslash}^{<n+1}\otimes \BB \ \left(\subseteq \mcD^{<\infty}_{\hslash} \otimes \BB\right).
\end{equation}
Since  ${^\Join} \Sigma^n(D^\clubsuit)=1$,  the composite
\begin{equation}  \label{GL13}
 \DF_\BB
 \xrightarrow{\mr{inclusion}}
  \mcD^{< \infty}_{\hslash} \otimes \BB
   \xrightarrow{\mr{quotient}}\mcD^{< \infty}_{\hslash}\otimes \BB/\langle \mr{Im}({^\circlearrowright \beta}^{n+1}  (D^\clubsuit))
   \rangle\end{equation}
is an isomorphism of  $\mcO_U$-modules.
The  $\mcD^{<\infty}_{\hslash}$-action on  $\mcD^{< \infty}_{\hslash}\otimes \BB/\langle \mr{Im}({^\circlearrowright \beta}^{n+1}  (D^\clubsuit))\rangle$ determines, via this isomorphism, 
an $\hslash$-$S^\mr{log}$-connection
\index{$D^\clubsuit$, $(n, \hslash, \BB)$-projective connection@--- $D^{\clubsuit \Rightarrow\diamondsuit}$, $D^{\clubsuit \Rightarrow\heartsuit}$}
\begin{equation}
\label{GL14}  D^{\clubsuit \Rightarrow\diamondsuit} :  \DF_\BB
  \migi \Omega_{} \otimes \DF_\BB
    \end{equation}
on  $\DF_\BB$.
One verifies that the collection of data
\begin{equation} 
\label{GL61}
 D^{\clubsuit \Rightarrow\heartsuit} := (\DF_\BB,
  D^{\clubsuit \Rightarrow\diamondsuit}, \{ \DF^j_\BB
  \}_{j=0}^n) 
    \end{equation}
forms a $(\mr{GL}_n, \hslash, \BB)$-oper on $U^\mr{log}/S^\mr{log}$.
The assignment $D^\clubsuit \mapsto D^{\clubsuit \Rightarrow \heartsuit}$ is functorial with respect to pull-back  to \'{e}tale $U$-schemes, and hence, determines a morphism
\begin{equation} 
\label{GL62}  \Lambda^{\clubsuit \Rightarrow \diamondsuit}_{\BB} : \DD^{\,\clubsuit}_{n, \hslash, \BB}\migi \mcO p_{n, \hslash, \BB}^\diamondsuit  
\index{$\Lambda_{\BB}^{\clubsuit \Rightarrow \diamondsuit}$, $\Lambda_{\bb (, \rho)}^{\clubsuit \Rightarrow \diamondsuit}$}
\end{equation}
of  presheaves on $\mfE \mft_{/U}$.
The formation of this morphism is  compatible with  pull-back to fs log schemes over $S^\mr{log}$.


\bpr  \label{GL19}
$ \Lambda^{\clubsuit \Rightarrow \diamondsuit}_{\BB}$ is an isomorphism of presheaves. In particular, $\mcO p_{n, \hslash, \BB}^\diamondsuit$ is a sheaf on $\mfE \mft_{/U}$. 
\epr
\begin{proof}
To prove the assertion, 
we shall construct the inverse to $\Lambda^{\clubsuit \Rightarrow \diamondsuit}_{\BB}$.
Let 
$\nabla^\diamondsuit : 
\DF_\BB
 \migi \Omega_{} \otimes  \DF_\BB$
 be 
a $(\mr{GL}_n, \hslash, \BB)$-oper on $U^\mr{log}/S^\mr{log}$.
It corresponds to  a  left $\mcD^{< \infty}_{\hslash}$-action
 ${^\mcD}\hspace{-1.5mm}\nabla^{\diamondsuit} :  \mcD^{< \infty}_{\hslash} \otimes
\DF_\BB
  \migi \DF_\BB$.
For $j \geq 0$,  we denote by   
\begin{equation}  \label{GL25}
 {^\mcD}\hspace{-1.5mm}\nabla^{\diamondsuit <j} 
   :  \mcD_{\hslash}^{< j}  \otimes \BB  \migi 
   \DF_\BB
\end{equation}
 the restriction  of ${^\mcD}\hspace{-1.5mm}\nabla^{\diamondsuit}$ to the $\mcO_U$-submodule
 $\mcD_{\hslash}^{<j}  \otimes \BB  \left(=  \mcD_{\hslash}^{<j} \otimes (\mcD_{\hslash}^{<1} \otimes \BB)\right)$ of  $\mcD^{< \infty}_{\hslash} \otimes (\mcD_{\hslash}^{< n}\otimes \BB)$.
By the definition of a $(\mr{GL}_n, \hslash, \BB)$-oper, 
 ${^\mcD}\hspace{-1.5mm}\nabla^{\diamondsuit <n}$ is an isomorphism and ${^\mcD}\hspace{-1.5mm}\nabla^{\diamondsuit <n+1}$ is surjective.
Hence, 
$({^\mcD}\hspace{-1.5mm}\nabla^{\diamondsuit <n})^{-1} \circ   {^\mcD}\hspace{-1.5mm}\nabla^{\diamondsuit <n+1}$
determines a split surjection of
the short exact sequence
\begin{equation}
\label{GL17}
 0 \migi  \DF_\BB  \migi \mcD^{< n+1}_{\hslash} \otimes \BB \migi \mcT_{}^{\otimes n} \otimes \BB \migi 0
  \end{equation}
arising from 
the natural isomorphism  $\mcD_{\hslash}^{ < n+1}/\mcD_{\hslash}^{ < n} \isom \mcT^{\otimes n}$;
the corresponding split injection
\begin{equation}
{^\circlearrowright}\nabla^{\diamondsuit \Rightarrow \clubsuit} : \mcT_{}^{\otimes n} \otimes \BB \migiincl \mcD^{< n+1}_{\hslash} \otimes \BB,
\end{equation}
of  this sequence  satisfies     $\JJ\Sigma^{n+1}({^\circlearrowright}\nabla^{\diamondsuit \Rightarrow \clubsuit}) =1$.
Since ${^\Join}\Sigma^{n+1}$ is compatible with $\JJ\Sigma^{n+1}$ via  the isomorphism ${^\circlearrowright \beta}^{n+1}$,
$\nabla^{\diamondsuit \Rightarrow \clubsuit} :=({^\circlearrowright \beta}^{n+1})^{-1}({^\circlearrowright}\nabla^{\diamondsuit \Rightarrow\clubsuit})$ belongs to  the set
$\Gamma (U, \DD_{n, \hslash, \mcB}^{\, \clubsuit})$.
The assignment 
$\nabla^\diamondsuit \mapsto  \nabla^{\diamondsuit \Rightarrow \clubsuit}$
 is functorial with respect to pull-back  to  \'{e}tale $U$-schemes,  and 
 the resulting morphism $\mcO p_{n, \hslash, \BB}^\diamondsuit \migi \DD^{\,\clubsuit}_{n, \hslash, \BB}$
determines the inverse to $\Lambda^{\clubsuit \Rightarrow \diamondsuit}_\BB$.
This completes the proof of the assertion.
\end{proof}

\subsection{A torsor structure on $\mcO p_{n, \hslash, \BB}^\diamondsuit$} \label{QS9452}
We shall  study a ${^\Join}\DD_\BB^{\,<n}$-torsor structure on $\mcO p_{n, \hslash, \BB}^\diamondsuit$.
Let us take $\nabla^\diamondsuit \in  \Gamma (U,  \mcO p_{n, \hslash, \BB}^\diamondsuit)$ and $R \in  \Gamma (U, {^\Join}\DD_\BB^{\,<n})$.
As in the proof of Proposition \ref{GL19} (cf. (\ref{GL25})),
$\nabla^\diamondsuit$ induces an automorphism ${^\mcD}\hspace{-1.5mm}\nabla^{\diamondsuit <n}$ of $\DF_\BB$.
Also, let $\gamma$ denote the injection 
\begin{equation}
\label{GL79}
 \JJ\DD^{\,<n}_\BB
  \migiincl \Omega_{} \otimes \mcE nd_{\mcO_U}
  (\DF_\BB) \left(= \mcH om_{\mcO_U} (\mcT \otimes \DF_\BB, \DF_\BB) \right)
    \end{equation}
given by forward composition with the quotient 
$\mcT_{} \otimes\DF_\BB \migisurj \mcT_{} \otimes \mr{gr}^0(\DF_\BB) \left(= \mcT_{}^{\otimes n} \otimes \BB \right)$.
Then,
\index{$\nabla^\diamondsuit$, $(\mr{GL}_n, \hslash, \BB)$-oper@--- $\nabla^\diamondsuit_{+ R}$}
\begin{equation}\label{QS9666}
\nabla^\diamondsuit_{+ R} :=  \nabla^\diamondsuit - ( {^\mcD}\hspace{-1.5mm}\nabla^{\diamondsuit <n})^{-1} \circ \gamma ({^\circlearrowright \beta}^{n}(R)) \circ  {^\mcD}\hspace{-1.5mm}\nabla^{\diamondsuit <n}, 
\end{equation}
defines an  $\hslash$-$S^\mr{log}$-connection on 
  $\DF_\BB$.

\bpr  \label{GL22} 
\begin{itemize}
  \item[(i)]
 Let $\nabla^\diamondsuit$ and  $R$ be as above.
 Then, $\nabla^\diamondsuit_{+ R}$ forms a $(\mr{GL}_n, \hslash, \BB)$-oper on $U^{\mr{log}}/S^\mr{log}$.
 Moreover, the formation of $\nabla^\diamondsuit_{+ R}$ commutes with pull-back to \'{e}tale $U$-schemes.
\item[(ii)]
The ${^\Join}\DD^{\,<n}_{\BB}$-action  
\begin{align}
\mcO p_{n, \hslash, \BB}^\diamondsuit \times  {^\Join}\DD^{\,<n}_{\BB}   &\migi   \mcO p_{n, \hslash, \BB}^\diamondsuit
\end{align}
on $\mcO p_{n, \hslash, \BB}^\diamondsuit$ given by 
$(\nabla^\diamondsuit, R) \mapsto \nabla^\diamondsuit_{+ R}$
is compatible with 
that 
 on
$\DD^{\,\clubsuit}_{n, \hslash, \BB}$ via  the isomorphism $\Lambda_{\BB}^{\clubsuit \Rightarrow \diamondsuit}$.
\end{itemize}
\epr
\begin{proof}
The assertion follows from the various definitions involved.
\end{proof}

\begin{rema}[Local description] \label{QQ198}
In this remark, we consider the local descriptions of  $\Lambda^{\clubsuit \Rightarrow \diamondsuit}_\BB$ and the ${^\Join}\DD^{\,<n}_\BB$-torsor structures  discussed above.

Let us take 
an $(n, \hslash, \BB)$-projective connection $D^\clubsuit$ on $U^\mr{log}/S^\mr{log}$.
Suppose that there exist a global log chart $(U, \partial)$ on $U^\mr{log}/S^\mr{log}$ and  a trivialization  $\BB \cong \mcO_U$ of the line bundle $\BB$.
These data  give rise to a decomposition $\mcD_\hslash^{< n} \otimes \BB = \bigoplus_{j=0}^{n-1} \mcO_U \cdot \partial^j$.
Using this  decomposition,
we can describe $D^\clubsuit$ as 
\begin{align}
D^\clubsuit = \partial^n + q_1 \cdot \partial^{n-1}+ \cdots + q_{n-1}\cdot \partial + q_n
\end{align}
 for some $q_1, \cdots, q_n \in \Gamma (U, \mcO_U)$.
Then, the $(\mr{GL}_n, \hslash, \BB)$-oper  $D^{\clubsuit \Rightarrow \diamondsuit}$ associated to $D^\clubsuit$ may be regarded as an $\hslash$-$S^\mr{log}$-connection on  $\mcO_U^{\oplus n} \left(=\bigoplus_{j=0}^{n-1} \mcO_U \cdot \partial^j\right)$ and 
satisfies  
\begin{align} \label{QQ325}
(D^{\clubsuit \Rightarrow \diamondsuit})^{\mr{GL}_n}(\partial) = 
\begin{pmatrix} 0 & 0& 0& \cdots & 0& - q_n 
\\
 1& 0 & 0& \cdots & 0 & -q_{n-1}
 \\
  0 & 1 & 0& \cdots & 0 & -q_{n-2} \\ 0 & 0 &1& \cdots & 0& -q_{n-3} \\
 \vdots & \vdots & \vdots & \ddots& \vdots & \vdots  \\
  0 & 0 & 0 & \cdots & 1 & -q_1  \end{pmatrix}
\end{align}
(cf. (\ref{QS7144}) for the definition of $(-)^{\mr{GL}_n}$).
Moreover, if $\partial^\vee \left(\in \Gamma (U, \Omega) \right)$ denotes the dual of $\partial$, then the $\hslash$-$S^\mr{log}$-connection $\mr{det}(D^{\clubsuit \Rightarrow \diamondsuit})$ on the determinant bundle $\mr{det}(\mcD_\hslash^{< n} \otimes \BB) \left(= \mr{det}(\bigoplus_{j=0}^{n-1} \mcO_U \cdot \partial^j) =  \mcO_U\right)$ is given by $\hslash \cdot  d - q_1 \cdot \partial^\vee$.

Next, let $R$ be an element of  $\Gamma (U, {^\Join}\DD_\BB^{\, <n})$ defined as $R = a_1 \cdot \partial^{n-1}+a_2 \cdot \partial^{n-2}+ \cdots + a_n$, where $a_1, \cdots, a_n \in \Gamma (U, \mcO_U)$. Then, the $\hslash$-$S^\mr{log}$-connection $D^{\clubsuit \Rightarrow \diamondsuit}_{+ R}$ resulting from  the ${^\Join}\DD^{\,<n}_\BB$-torsor structure on $\mcO p^\diamondsuit_{n, \hslash, \BB}$ satisfies (thanks to Proposition \ref{GL22}, (ii)) the following equality:
\begin{align} \label{QQ324}
(D^{\clubsuit \Rightarrow \diamondsuit}_{+ R})^{\mr{GL}_n}(\partial) = 
\begin{pmatrix} 0 & 0& 0& \cdots & 0& - q_n  -a_{n}
\\
 1& 0 & 0& \cdots & 0 & -q_{n-1} -a_{n-1}
 \\
  0 & 1 & 0& \cdots & 0 & -q_{n-2} - a_{n-2} \\ 0 & 0 &1& \cdots & 0& -q_{n-3} -a_{n-3}\\
 \vdots & \vdots & \vdots & \ddots& \vdots & \vdots  \\
  0 & 0 & 0 & \cdots & 1 & -q_1 -a_{1} \end{pmatrix}.
\end{align}
\end{rema}

\section{$(n, \hslash)$-theta characteristics} \label{y0120}

\subsection{The definition of an $(n, \hslash)$-theta characteristic} \label{QG830}

We provide the definition of an $(n, \hslash)$-theta characteristic, generalizing the usual definition of  a theta characteristic.

\bde \label{y0270}
\begin{itemize}
\item[(i)]
An {\bf $(n, \hslash)$-theta characteristic} \index{$\bb$, $(n, \hslash)$-theta characteristic} of $U^\mr{log}/S^\mr{log}$ is a pair 
\begin{align}
\bb := (\BB, \nabla_\bb) 
\end{align}
 consisting of a line bundle $\BB$ on $U$ and an $\hslash$-$S^\mr{log}$-connection $\nabla_\bb$ on the line bundle 
 $\mcT^{\otimes \frac{n(n-1)}{2}} \otimes \BB^{\otimes n}$.
 For simplicity, each $(n ,1)$-theta characteristic is called  an {\bf $n$-theta characteristic}.
\item[(ii)]
Let $\bb := (\BB, \nabla_\bb)$,  $\bb' = (\BB', \nabla_{\bb'})$ be  $(n, \hslash)$-theta characteristics
of $U^\mr{log}/S^\mr{log}$.
Then, an {\bf isomorphism of $(n, \hslash)$-theta characteristics} from $\bb$ to $\bb'$ is an isomorphism
$\BB \isom \BB'$ of line bundles such that 
the induced isomorphism $\mcT^{\otimes \frac{n(n-1)}{2}} \otimes \BB^{\otimes n} \isom \mcT^{\otimes \frac{n(n-1)}{2}} \otimes \BB'^{\otimes n}$
 is compatible with the respective $\hslash$-$S^\mr{log}$-connections $\nabla_\bb$,  $\nabla'_\bb$.
\end{itemize}
  \ede

\begin{rema} \label{y0121}
Let $\bb := (\BB, \nabla_\bb)$ be an $(n, \hslash)$-theta characteristic of  $U^\mr{log}/S^\mr{log}$.
Then, as discussed  in Remark \ref{y01Frh1}, there exists a canonical isomorphism
\begin{align}  \label{fff0100}
\mr{det} (\DF_\BB)
   \isom \mcT_{}^{\otimes \frac{n (n-1)}{2}} \otimes \BB^{\otimes n}.  
 \end{align}
By  using this,  
we occasionally  regard $\nabla_\bb$ as an $\hslash$-$S^\mr{log}$-connection on the determinant bundle $\mr{det}(\DF_\BB)$.
\end{rema}

\begin{exa}[The case of $\hslash =0$] \label{QQ181}
Suppose that $\hslash =0$.
Then, 
 giving an $(n, 0)$-theta characteristic of   $U^\mr{log}/S^\mr{log}$ is equivalent to giving a line bundle $\BB$ together with a global section of $\Omega_{}$ (cf. Example \ref{QS9900}).
 In particular,  for any $n$, there always exists an $(n, 0)$-theta characteristic on any log curve.
\end{exa}

\subsection{Pull-back and base-change} \label{QQ155}
Let $\bb := (\BB, \nabla_\bb)$ be an $(n, \hslash)$-theta characteristic of $U^\mr{log}/S^\mr{log}$.
Also, let $u: U' \migi U$ be an \'{e}tale $U$-scheme and write $\mcT':= \mcT_{U'^{\mr{log}}/S^\mr{log}}$. 
Then,  there exists a natural isomorphism  $u^*(\mcT^{\otimes \frac{n(n-1)}{2}}\otimes \BB^{\otimes n}) \isom  \mcT'^{\otimes \frac{n(n-1)}{2}} \otimes u^*(\BB)$, by which 
 the pair 
\index{$\bb$, $(n, \hslash)$-theta characteristic@--- ${u}^*(\bb)$, pull-back of $\bb$}
\begin{equation} \label{GL40}
{u}^*(\bb) := ({u}^*(\BB), {u}^*(\nabla_\bb))
\end{equation}
 defines  an $(n, \hslash)$-theta characteristic of  $U'^{\mr{log}}/S^\mr{log}$.
We shall refer to ${u}^*(\bb)$ as  the {\bf pull-back} of $\bb$ to $U'^{\mr{log}}$.

Also. if $s^\mr{log} : S'^{\mr{log}} \migi S^\mr{log}$ is a morphism of fs log schemes over $k$, then 
the pair
\index{$\bb$, $(n, \hslash)$-theta characteristic@--- $s^*(\bb)$, base-change of $\bb$}
\begin{equation}  \label{GL41}
s^*(\bb) := (s^*(\BB), s^*(\nabla_\bb)),
\end{equation}
  where the notation $s^*(-)$ on the right-hand side denotes pull-back by 
  $s^\mr{log}$,
   specifies
an $(n, s^*(\hslash))$-theta characteristic of  the log curve $(S'^{\mr{log}} \times_{S^\mr{log}} U^\mr{log})/{S}'^{\mr{log}}$.
We refer to $s^*(\bb)$ as the {\bf base-change} of $\bb$ to $S'^{\mr{log}}$.

\begin{exa}[Theta characteristic] \label{QQ188}
By a {\bf theta characteristic} \index{theta characteristic} of $U^\mr{log}/S^\mr{log}$,
we means a line bundle $\varTheta$ on $U$ together with an isomorphism $ \mcO_U \isom \mcT \otimes \varTheta^{\otimes 2}$.
The notion of a theta characteristic  is classical and has been  well studied for  the non-logarithmic case.
We know that  
 any proper smooth  curve   over an algebraically closed field always has  a theta characteristic.
For example, if, moreover, the base field is of odd characteristic,
then 
 there are exactly $2^{2g}$ theta characteristics of such a curve  in total, where $g$ denotes its genus (cf. \S\,\ref{QQ156} below).

Now, suppose that we are given a theta characteristic $\varTheta$ of $U^\mr{log}/S^\mr{log}$.
Fix an  isomorphism  $\tau_\BB : \mcO_U \isom \mcT \otimes \varTheta^{\otimes 2}$.
Then, 
$\tau_{\BB*} (\hslash \cdot d) \left(= \hslash \cdot  \tau_{\BB*} (d)\right)$ defines 
an $\hslash$-$S^\mr{log}$-connection on $\mcT \otimes\BB^{\otimes 2}$, and 
 the pair $(\varTheta, \tau_{\BB*} (\hslash \cdot d))$  specifies a $(2, \hslash)$-theta characteristic of $U^\mr{log}/S^\mr{log}$.
 More generally,  
  the line bundle $\mcT^{\otimes \frac{n (n-1)}{2}} \otimes (\varTheta^{\otimes (n-1)})^{\otimes n}$ admits an $\hslash$-$S^\mr{log}$-connection corresponding to $\hslash \cdot d$ via 
 the composite isomorphism 
 \begin{align}
 \mcO_U \xrightarrow{\tau_\BB^{\otimes \frac{n (n-1)}{2}}}
 (\mcT \otimes \varTheta^{\otimes 2})^{\otimes \frac{n(n-1)}{2}} \isom 
 \mcT^{\otimes \frac{n (n-1)}{2}} \otimes (\varTheta^{\otimes (n-1)})^{\otimes n}.
 \end{align}
 The line bundle  $\varTheta^{\otimes (n-1)}$ together with this connection
  forms an $(n, \hslash)$-theta characteristic of $U^\mr{log}/S^\mr{log}$.
 In this way,  there is  a canonical construction of  an $(n, \hslash)$-theta characteristic (for any $n$ and $\hslash$) by using a theta characteristic.
 In other words,  $(n, \hslash)$-theta characteristics  may be thought of as  generalizations of theta characteristics.
\end{exa}

\subsection{$(n, \hslash)$-theta characteristics on pointed smooth curves} \label{QQ156}
Let us consider the case where  $U^\mr{log} = X^\mr{log}$ for an $r$-pointed proper  {\it smooth} curve $\msX := (X/S, \{ \sigma_i \}_{i=1}^r)$  (hence, $S^\mr{log} = S$).
Assume   that 
  either $r$ or $n-1$ is even.
 In what follows, we shall prove the claim that 
 {\it $X^\mr{log}/S$ admits,  at least \'{e}tale locally on $S$,   an $(n, \hslash)$-theta characteristic.}
For each integer $d$,  let  $\mr{Pic}^d_{X/S}$ \index{$\mr{Pic}^d_{X/S}$, relative Picard scheme of $X/S$} denote the relative Picard scheme of $X/S$ classifying the set of (equivalence classes, relative to the equivalence relation determined by tensoring with a line bundle pulled-back from the base $S$, of) degree $d$ line bundles on $X$.
  The assumption implies that  $n(n-1)(2g-2+r)/2$ is divisible by $n$, so
we obtain  the morphism
\begin{align}  
\times n : \mr{Pic}_{X/S}^{\frac{(n-1)(2g-2+r)}{2}}  \migi   \mr{Pic}_{X/S}^{\frac{n(n-1)(2g-2+r)}{2}}
 \end{align}
given by $[\mcL] \mapsto [\mcL^{\otimes n}]$.
 Since $n$ is invertible in $k$, 
 the inverse image via this morphism of the global section $S \migi \mr{Pic}_{X/S}^{\frac{n(n-1)(2g-2+r)}{2}}$ determined by 
$\Omega^{\otimes \frac{n (n-1)}{2}}$ 
forms a  torsor over $S$ modeled on $\mr{Pic}^0_{X/S}[n] := \mr{Ker}\left( \mr{Pic}_{X/S}^0 \xrightarrow{\times n}\mr{Pic}_{X/S}^0\right)$.
In particular,  
after possibly replacing $S$ with its  covering in the \'{e}tale topology, 
there exists a line bundle $\varTheta$ on $X$ with
$\varTheta^{\otimes n} \cong \Omega^{\otimes \frac{n (n-1)}{2}}$.
Under the resulting  identification $\mcT^{\otimes \frac{n (n-1)}{2}} \otimes \varTheta^{\otimes n} \cong \mcO_X$,  $\hslash \cdot d$ may be regarded as   an  $\hslash$-$S$-connection on $\mcT^{\otimes \frac{n (n-1)}{2}} \otimes \varTheta^{\otimes n}$.
Consequently, the pair $(\varTheta, \hslash \cdot d)$ specifies an $(n, \hslash)$-theta characteristic of $X^\mr{log}/S$.
This completes the proof of the claim.

Notice that, unlike 
this situation,
a nodal curve  may not admit  such a line bundle, i.e., a line bundle $\varTheta$
with $\varTheta^{\otimes n} \cong \Omega^{\otimes \frac{n(n-1)}{2}}$
  (cf. the discussion in ~\cite{Chi}, 1.2.4,  Example). 

\subsection{$(n, \hslash)$-theta characteristics in positive characteristic} \label{QQ183}
If   $\mr{char}(k)=:p>0$,
then 
 any log curve $U^\mr{log}/S^\mr{log}$ admits an $(n, \hslash)$-theta characteristic  $\bb := (\BB, \nabla_\bb)$ with ${^p}\psi^{\nabla_\bb}=0$.
 In particular, there exists an $(n, \hslash)$-theta characteristic  on the universal family of pointed stable curves $\msC_{g,r}$ over $\overline{\mfM}_{g,r}$.
In what follows, we shall prove this fact by constructing an explicit example.

Since $(p, n)=1$,
 one may find  a pair of  integers $(a, b)$ with
 $a \cdot n + b \cdot p = \frac{n (n-1)}{2}$.
By letting $\BB^{a,b} := \Omega_{}^{\otimes a}$,
we have a composite isomorphism
 \begin{align}
 \mcT^{\otimes \frac{n (n-1)}{2}} \otimes (\Omega^{\otimes a})^{\otimes n} \isom \mcT^{\otimes (\frac{n(n-1)}{2}-a \cdot n)}\isom \mcT^{\otimes b \cdot p} \isom F^*_{U/S} (\mcT_{U^{(1)\mr{log}}/S^\mr{log}}^{\otimes b}).
 \end{align}
The canonical $\hslash$-$S^\mr{log}$-connection
on  $F^*_{U/S}(\mcT_{U^{(1)\mr{log}}/S^\mr{log}}^{\otimes b})$ (cf. (\ref{QR110})) 
determines, via the above  composite, an $\hslash$-$S^\mr{log}$-connection $\nabla_\bb^{a,b}$ on $\mcT^{\otimes \frac{n (n-1)}{2}} \otimes (\Omega^{\otimes a})^{\otimes n}$ with vanishing $p$-curvature (cf.  Proposition \ref{y073}).
In particular, the resulting pair 
\begin{align} \label{QH2}
\bb^{a,b} :=(\varTheta^{a,b}, \nabla_\bb^{a, b})
\end{align}
specifies  a required $(n, \hslash)$-theta characteristic.

\subsection{Twisting with a flat line bundle}
Let $\bb := (\BB, \nabla_\bb)$ be  an $(n, \hslash)$-theta characteristic of $U^\mr{log}/S^\mr{log}$.
For an $\hslash$-flat  line bundle $\msL = (\mcL, \nabla_\mcL)$, 
the pair
\index{$\bb$, $(n, \hslash)$-theta characteristic@--- $\bb \otimes  \msL$}
 \begin{equation}
 \label{GL34}
 \bb \otimes  \msL := (\BB \otimes \mcL, \nabla_\bb \otimes \nabla_\mcL^{\otimes n}),
 \end{equation}
 forms an $(n, \hslash)$-theta characteristic of $U^\mr{log}/S^\mr{log}$.
 Conversely,  we have the following property:

\ble \label{y0122}
Let $\bb' := (\varTheta', \nabla_{\bb'})$ be another $(n, \hslash)$-theta characteristic of $U^\mr{log}/S^\mr{log}$.
Then, 
 there exists  an $\hslash$-flat  line bundle 
 \index{$\bb$, $(n, \hslash)$-theta characteristic@--- $\bb'/\bb$}
\begin{align} \label{QQ1138}
\bb'/\bb := (\mcN, \nabla_\mcN)
\end{align}
on $U^\mr{log}/S^\mr{log}$ such  that  $\bb \otimes  \bb'/\bb$ is isomorphic to $\bb'$.
Moreover, such an $\hslash$-flat  line bundle   is uniquely determined up to isomorphism,
\ele
\begin{proof}
Denote by
$\nabla_{\mcL^{\otimes n}}$ the $\hslash$-$S^\mr{log}$-connection on the line bundle $\mcL^{\otimes n}$, where $\mcL := \varTheta^\vee \otimes \varTheta'$,  corresponding to $\nabla_\bb^\vee \otimes \nabla_{\bb'}$ via the composite of natural  isomorphisms
\begin{align}
(\mcT^{\otimes \frac{n (n-1)}{2}} \otimes\varTheta^{\otimes n})^\vee \otimes 
(\mcT^{\otimes \frac{n (n-1)}{2}} \otimes\varTheta'^{\otimes n})
&\isom (\varTheta^{\otimes n})^\vee \otimes \varTheta'^{\otimes n} \isom \mcL^{\otimes n}. \notag
\end{align}
Since $n<p$, it follows from  Proposition \ref{y0111}, (i), that there exists an $\hslash$-$S^\mr{log}$-connection $\nabla_\mcL$ on $\mcL$ with $\nabla_{\mcL}^{\otimes n} = \nabla_{\mcL^{\otimes n}}$.
Then,  one verifies  that
$\bb \otimes  (\mcL, \nabla_\mcL) \isom \bb'$, which completes the proof of the existence portion.
The uniqueness follows immediately from the uniqueness portion of  Proposition  \ref{y0111}, (i).
\end{proof}

\subsection{$(n, \hslash, \bb)$-projective connections and $(\mr{GL}_n, \hslash, \bb)$-opers} \label{y0124}

Let $\bb := (\BB, \nabla_\bb)$ be an $(n, \hslash)$-theta characteristic of  $U^\mr{log}/S^\mr{log}$.

\bde \label{y0123}
A {\bf $(\mr{GL}_n, \hslash, \bb)$-oper}  
\index{$\nabla^\diamondsuit$, $(\mr{GL}_n, \hslash, \bb)$-oper} 
on $U^\mr{log}/S^\mr{log}$ is a $(\mr{GL}_n, \hslash, \BB)$-oper $\nabla^\diamondsuit$ on $U^\mr{log}/S^\mr{log}$ satisfying   $\mr{det}(\nabla^\diamondsuit) = \nabla_\bb$ under the identification $\mr{det}(\DF_\BB) = \mcT^{\otimes \frac{n(n-1)}{2}}\otimes \BB^{\otimes n}$ by (\ref{fff0100}).
Each $(\mr{GL}_n, 1, \bb)$-oper is  called a {\bf $(\mr{GL}_n, \bb)$-oper}.
Next, suppose further that
  $U^\mr{log}/S^\mr{log}$ is  the log curve associated to  a local $r (\geq 0)$-pointed stable curve $\msU$ together with a choice of its base.
  Then,  we refer to any $(\mr{GL}_n, \hslash, \bb)$-oper on $U^\mr{log}/S^\mr{log}$ as a {\bf $(\mr{GL}_n, \hslash, \bb)$-oper on $\msU$} if there is no fear of confusion.
\ede

We shall write 
 \index{$\mcO p_{n, \hslash, \BB}^\diamondsuit$, $\mcO p_{n, \hslash, \bb (, \rho)}^\diamondsuit$, $\mcO p_{n, \hslash, \BB}^{\diamondsuit, \mr{O}}$, $\mcO p_{n, \hslash, \BB}^{\diamondsuit, \mr{Sp}}$}
\begin{equation}   \label{QQ219}
\mcO p^\diamondsuit_{n, \hslash, \bb}  
\end{equation}
for  the subsheaf of $\mcO p^\diamondsuit_{n, \hslash, \BB}$ classifying $(\mr{GL}_n, \hslash, \bb)$-opers, and write
\index{$\DD^{\,\clubsuit}_{n, \hslash, \BB}$, $\DD^{\,\clubsuit}_{n, \hslash, \bb (, \rho)}$, $\mfD \mfi \mff \mff^{\,\clubsuit}_{n, \hslash, \bb (, \rho)}$}
\begin{equation}  
\DD^{\,\clubsuit}_{n, \hslash, \bb}  :=  (\Lambda_\BB^{\clubsuit \Rightarrow \diamondsuit})^{-1}(\mcO p^\diamondsuit_{n, \hslash, \bb}) \left(\subseteq \DD^{\,\clubsuit}_{n, \hslash, \BB}  \right).
  \end{equation}
The isomorphism  $\Lambda^{\clubsuit \Rightarrow \diamondsuit}_{\BB}$ restricts to an isomorphism
\begin{equation}  
\label{GL28}
\Lambda_{\bb}^{\clubsuit \Rightarrow \diamondsuit} :  
\DD^{\, \clubsuit}_{n, \hslash,  \bb}
  \isom   \mcO p^\diamondsuit_{n, \hslash, \bb}
\index{$\Lambda_{\BB}^{\clubsuit \Rightarrow \diamondsuit}$, $\Lambda_{\bb (, \rho)}^{\clubsuit \Rightarrow \diamondsuit}$}
   \end{equation}
of sheaves on $\mfE \mft_{/U}$.
Thus,
we have
\begin{align} \label{QQ239}
\mcO p_{n, \hslash, \BB}^\diamondsuit = \coprod_{\bb'} \mcO p_{n, \hslash, \bb'}^\diamondsuit,
\hspace{5mm}
\DD^{\,\clubsuit}_{n, \hslash, \BB} = \coprod_{\bb'} \DD^{\,\clubsuit}_{n, \hslash,  \bb'}, 
\hspace{5mm}
\Lambda^{\clubsuit \Rightarrow \diamondsuit}_\BB = \coprod_{\bb'} \Lambda^{\clubsuit \Rightarrow \diamondsuit}_{\bb'}, 
\end{align}
where all the sums run over 
the set of $(n, \hslash)$-theta characteristics $\bb' :=(\BB', \nabla_{\bb'})$ of $U^\mr{log}/S^\mr{log}$ with $\BB' = \BB$.

\bde  \label{QG3887}
\begin{itemize}
\item[(i)]
Let $D^\clubsuit$ be an $(n, \hslash, \BB)$-projective connection on $U^\mr{log}/S^\mr{log}$.
Then, the {\bf subprincipal symbol} 
\index{$D^\clubsuit$, $(n, \hslash, \BB)$-projective connection@--- subprincipal symbol of}
of $D^\clubsuit$ is defined as a unique $S^\mr{log}$-connection $\nabla_{\bb'}$ on 
the line bundle $\mcT^{\otimes \frac{n(n-1)}{2}} \otimes \BB^{\otimes n}$ such that  $D^\clubsuit \in \Gamma (U, \DD^{\,\clubsuit}_{n, \hslash,  (\BB, \nabla_{\bb'})})$. 
\item[(ii)]
An {\bf $(n, \hslash, \bb)$-projective connection}
\index{$D^\clubsuit$, $(n, \hslash, \bb)$-projective connection}
 on $U^\mr{log}/S^\mr{log}$ is an $(n, \hslash, \BB)$-projective connection on $U^\mr{log}/S^\mr{log}$ whose subprincipal symbol is $\nabla_\bb$, or equivalently, which is classified by 
 $\Gamma (U, \DD^{\,\clubsuit}_{n, \hslash, \bb})$.
Each $(n, 1, \bb)$-projective connection is called an {\bf $(n, \bb)$-projective connection}.
Next, suppose further that  $U^\mr{log}/S^\mr{log}$ is  the log curve associated to  a local $r (\geq 0)$-pointed stable curve $\msU$ together with a choice of its base.
Then, 
we refer to any $(n, \hslash, \bb)$-projective connection on $U^\mr{log}/S^\mr{log}$ as an {\bf $(n, \hslash, \bb)$-projective connection on $\msU$} if there is no fear of confusion.
\end{itemize}
\ede


By taking account of the local description of  $\Lambda_\BB^{\clubsuit \Rightarrow \diamondsuit}$  (cf. Remark \ref{QQ198}, (\ref{QQ325})),
 we see that
$\DD^{\, \clubsuit}_{n, \hslash,  \bb}$, as well as $\mcO p^\diamondsuit_{n, \hslash, \bb}$, has a section  at least \'{e}tale locally on $U$.
Moreover, 
if we write
\begin{align} \label{QW352}
\DV_\BB := {^\Join}\DD^{\,<n-1}_{\BB},
\index{$\DV_\BB$}
\end{align}
then the local description (\ref{QQ324}) tells us  that
the ${^\Join}\DD^{\,<n}_{\BB}$-torsor structures on  $\DD^{\,\clubsuit}_{n, \hslash, \BB}$ and 
$\mcO p_{n, \hslash, \BB}^\diamondsuit$ restrict to   $\DV_\BB$-torsor structures on $\DD^{\,\clubsuit}_{n, \hslash, \bb}$ and  $\mcO p_{n, \hslash, \bb}^\diamondsuit$ respectively.

\begin{rema} \label{QW388}
Let  us keep the situation in Remark \ref{QQ198} and let $\bb := (\BB, \nabla_\bb)$ be as above.
The  $\hslash$-$S^\mr{log}$-connection $\nabla_\bb$ may be expressed as $\nabla_\bb = \hslash \cdot d + w \cdot \partial^\vee$ for some $w \in \Gamma (U, \mcO_U)$.
Then,  for each $D^\clubsuit \in  \DD^{\,\clubsuit}_{n, \hslash, \BB}$,
the $(\mr{GL}_n, \hslash, \BB)$-oper $D^{\clubsuit \Rightarrow \diamondsuit}$ specifies a $(\mr{GL}_n, \hslash, \bb)$-oper if  and only if the matrix $(D^{\clubsuit \Rightarrow \diamondsuit})^{\mr{GL}_n} (\partial)$ may be expressed as (\ref{QQ325}) with $q_1  = -w$.

\end{rema}

\begin{rema} \label{QR50}
Notice  that if there exists  a global log chart $(U, \partial)$ on $U^\mr{log}/S^\mr{log}$, then we have a canonical 
$(n, \hslash)$-theta characteristic defined on the entire log curve.
Indeed, 
under the trivialization $\mcT^{\otimes \frac{n(n-1)}{2}} = \mcO_U \cdot \partial^{\otimes \frac{n(n-1)}{2}}$,
we can regard the  trivial $\hslash$-$S^\mr{log}$-connection $\hslash \cdot d$  on $\mcO_U$ as an $\hslash$-$S^\mr{log}$-connection on $\mcT^{\otimes \frac{n(n-1)}{2}} \left(= \mcT^{\otimes \frac{n(n-1)}{2}} \otimes \mcO_U^{\otimes n} \right)$. 
Thus, the pair 
\begin{align}  \label{QW34}
\bb_{(U, \partial)}^{\mr{triv}} := (\mcO_U, \hslash \cdot d)
\end{align}
 determines an $(n, \hslash)$-theta characteristic of $U^\mr{log}/S^\mr{log}$.
If we  take  an $(n, \hslash, \bb_{(U, \partial)}^{\mr{triv}})$-projective connection  $D^\clubsuit$ on $U^\mr{log}/S^\mr{log}$,
then 
the  matrix    $(D^{\clubsuit \Rightarrow \diamondsuit})^{\mr{GL}_n}(\partial)$ can be expressed as   (\ref{QQ325}) with $q_1=0$.
\end{rema}

\subsection{Comparison} \label{QG860}
At the end of this section, we summarize the correspondences between various sheaves that have been proved so far.
The remaining one is as follows.

\bpr \label{y0127} 
Let us consider the  morphism of sheaves
\begin{equation}  \label{QW1}
 \Lambda^{\diamondsuit \Rightarrow\heartsuit}_{\bb} :  \mcO p_{n, \hslash, \bb}^\diamondsuit  \migi \overline{\mcO} p_{n, \hslash}^\heartsuit 
 \index{$\Lambda_{\bb (, \rho)}^{\diamondsuit \Rightarrow \heartsuit}$, $\Lambda_\bb^{\diamondsuit \Rightarrow \heartsuit, \BL}$, $\Lambda_{\bb, \msX}^{\diamondsuit \Rightarrow \heartsuit, \BL}$}
 \end{equation}
given  by assigning, to each $(\mr{GL}_n, \hslash, \bb)$-oper 
$\nabla^\diamondsuit$, the equivalence class represented by the $(\mr{GL}_n, \hslash)$-oper $\nabla^{\diamondsuit \Rightarrow\heartsuit}$.
Then, it is an isomorphism and compatible, in an evident sense, with  pull-back to  fs log schemes over $S^\mr{log}$.
\epr
\begin{proof}
By the functorial formations of 
$\mcO p_{n, \hslash, \bb}^\diamondsuit$, $\overline{\mcO} p_{n, \hslash}^\heartsuit$, and 
 $\Lambda^{\diamondsuit \Rightarrow\heartsuit}_{\bb}$, it suffices to prove the bijectivity of the induced map of sets 
 $ \Lambda^{\diamondsuit \Rightarrow\heartsuit}_{\bb} (U) :  \mcO p_{n, \hslash, \bb}^\diamondsuit  (U) \migi \overline{\mcO} p_{n, \hslash}^\heartsuit (U)$.

First, let us  consider the injectivity of $  \Lambda^{\diamondsuit \Rightarrow \heartsuit}_{\bb}(U)$.
Let 
 $\nabla_1^\diamondsuit$,  $\nabla_2^\diamondsuit$ be  $(\mr{GL}_n, \hslash, \bb)$-opers on $U^\mr{log}/S^\mr{log}$ 
 with $\nabla_1^{\diamondsuit \Rightarrow \heartsuit} \stackrel{\heartsuit}{\sim} \nabla_2^{\diamondsuit \Rightarrow \heartsuit}$. 
  In the following, we shall prove that these  $(\mr{GL}_n, \hslash, \bb)$-opers define the same section of $\mcO p_{n, \hslash, \bb}^\diamondsuit$.
  To this end, we are always free to replace $U$ by its \'{e}tale covering in the \'{e}tale topology.
In particular,  it can be assumed 
that there exist 
  an $\hslash$-flat line bundle $\msL := (\mcL, \nabla_\mcL)$ on $U^\mr{log}/S^\mr{log}$
  and  an isomorphism $\alpha : (\nabla_1^{\diamondsuit \Rightarrow \heartsuit})_{\otimes \msL} \isom \nabla_2^{\diamondsuit \Rightarrow \heartsuit}$
of $(\mr{GL}_n, \hslash)$-opers.
The isomorphism $\alpha$ restricts to an isomorphism
\begin{align}
\left( \DF_\BB \otimes \mcL\supseteq \right) (\mcD_{\hslash}^{< 1} \otimes \BB) \otimes \mcL \isom \mcD_{\hslash}^{< 1} \otimes \BB \left(\subseteq \DF_\BB
\right)
\end{align}
 between line bundles.
This implies that  $\mcL \cong \mcO_U$, i.e., $\mcL$ can be identified  with $\mcO_U$ by  fixing an isomorphism $\mcL \isom  \mcO_U$.
Hence, $\alpha$  defines an automorphism of  
$\DF_{\BB}$,
so its determinant $\mr{det}(\alpha) : \mr{det}(\DF_\BB)
\isom \mr{det}(\DF_\BB)$
 is given by multiplication by an element $a$ of $\Gamma (U, \mcO_U^\times)$.
Let us replace $U$ by its covering so that there exists $b \in \Gamma (U, \mcO_U^\times)$ with $b^n = a$.
Since   $\mr{det}(\alpha)$ is compatible with the respective $\hslash$-$S^\mr{log}$-connections $\mr{det}(\nabla_1^\diamondsuit \otimes \nabla_\mcL) \left(= \nabla_\bb \otimes \nabla_{\mcL}^{\otimes n}\right)$ and   $\mr{det}(\nabla_2^\diamondsuit) \left(= \nabla_\bb \right)$,
  the automorphism $\beta$ of $\mcO_U$ given by multiplication by $b$ is compatible with   $\nabla_\mcL$ and $d$.
Thus, we obtain an isomorphism 
$\mr{id}_{\DF_\BB} \otimes \beta : (\nabla_1^{\diamondsuit \Rightarrow \heartsuit})_{\otimes \msL} \isom \nabla_1^{\diamondsuit \Rightarrow \heartsuit}$
 of $(\mr{GL}_n, \hslash)$-opers.
The composite   $(\mr{id}_{\DF_\BB} \otimes \beta)\circ\alpha^{-1}$  specifies an isomorphism $\nabla_2^{\diamondsuit \Rightarrow \heartsuit} \isom \nabla_1^{\diamondsuit \Rightarrow \heartsuit}$.
This completes the injectivity of  $\Lambda^{\diamondsuit \Rightarrow \heartsuit}_{\bb}(U)$.

Next,  let us consider  the surjectivity  of $\Lambda^{\diamondsuit \Rightarrow \heartsuit}_{\bb}(U)$.
Let $\msF^\heartsuit :=(\mcF, \nabla, \{ \mcF^j \}_{j=0}^n)$ be an arbitrary 
 $(\mr{GL}_n, \hslash)$-oper  on $U^\mr{log}/S^\mr{log}$.
Because of the functorialities mentioned at the beginning, the problem is reduced  to finding  a $(\mr{GL}_n, \hslash, \bb)$-oper $\nabla^\diamondsuit$ with  $\nabla^{\diamondsuit \Rightarrow \heartsuit} \stackrel{\heartsuit}{\sim} \msF^\heartsuit$.
By letting  $\mcL := \BB \otimes \mcF^{n-1\vee}$, we obtain the following composite isomorphism:
\begin{align} \label{QH4}
\mr{det}(\mcF \otimes \mcL) &\isom   \mr{det}(\mcF) \otimes \mcL^{\otimes n} \\
&\isom 
 \mcT^{\otimes \frac{n(n-1)}{2}}\otimes  (\mcF^{n-1})^{\otimes n}
 \otimes \left(\BB \otimes \mcF^{n-1\vee} \right)^{\otimes n}
  \notag \\
  &\isom
 \mcT^{\otimes \frac{n(n-1)}{2}}\otimes \BB^{\otimes n}, \notag
\end{align}
where the second arrow follows from (\ref{QW1906}).
 By Proposition  \ref{y0111}, (i),  there exists a unique $\hslash$-$S^\mr{log}$-connection $\nabla_\mcL$ on $\mcL$  such that
  the  determinant $\mr{det}(\nabla \otimes \nabla_\mcL)$
  $\left(=
  \mr{det}(\nabla) \otimes \nabla^{\otimes n}_\mcL\right)$  corresponds to  $\nabla_\bb$ via  (\ref{QH4}).
Here, we shall write $\nabla^\diamondsuit := \nabla \otimes \nabla_\mcL$.
By the definition of a $(\mr{GL}_n, \hslash)$-oper,  
the morphism
\begin{equation} \label{GL33}
\DF_\BB\left(=   \mcD_{\hslash}^{<n}\otimes (\mcF^{n-1} \otimes \mcL)\right) \migi \mcF \otimes \mcL 
\end{equation}
obtained by restricting $ {^\mcD}\hspace{-1.5mm}\nabla^\diamondsuit$
 to 
$\mcD_{\hslash}^{<n}\otimes (\mcF^{n-1} \otimes \mcL) \left(\subseteq \mcD^{<\infty}_{\hslash}\otimes (\mcF \otimes \mcL)\right)$ is an isomorphism.
 Under the identification  $\DF_\BB = \mcF \otimes \mcL$  using this isomorphism,
 $\nabla^\diamondsuit$ specifies 
 a $(\mr{GL}_n, \hslash, \bb)$-oper on $U^\mr{log}/S^\mr{log}$, which  moreover satisfies $\nabla^{\diamondsuit \Rightarrow \heartsuit} \stackrel{\heartsuit}{\sim}\msF^\heartsuit$.
This completes  the proof of the desired  surjectivity.

The remaining compatibility assertion  is immediate from the definition of $\Lambda_\bb^{\diamondsuit \Rightarrow \heartsuit}$, so we finish the proof of the assertion. 
\end{proof}
\vspace{3mm}

By combining  Propositions \ref{y0112}, \ref{y0127}, and the discussion preceding  Definition \ref{QG3887}, 
we have obtained  the following theorem.
\bt \label{y0128}
Let $\bb$
be an $(n, \hslash)$-theta characteristic  of $U^\mr{log}/S^\mr{log}$.
Then, there exists a canonical sequence of isomorphisms between  sheaves on $\mfE \mft_{/U}$
 \begin{align} \label{GL42}
 \vcenter{\xymatrix@C=46pt@R=36pt{
\DD^{\, \clubsuit}_{n,\hslash, \bb} \ar[rd]^-{\sim}_-{\Lambda_\bb^{\clubsuit \Rightarrow \diamondsuit}}&& \overline{\mcO}p^\heartsuit_{n,\hslash}\ar[rd]^-{\sim}_-{\Lambda^{\heartsuit \Rightarrow \spadesuit}}&\\
& \mcO p^\diamondsuit_{n, \hslash, \bb}\ar[ur]_-{\sim}^-{\Lambda_\bb^{\diamondsuit \Rightarrow \heartsuit}}&& \mcO p^\spadesuit_{n,\hslash}.
}}
 \end{align}
 \et

\section{Radii of $(\mr{GL}_n, \hslash, \bb)$-opers and $(n, \hslash, \bb)$-projective connections} \label{QW274}

We shall consider the notion of radius defined for each $(\mr{GL}_n, \hslash, \bb)$-oper, as well as each $(n, \hslash, \bb)$-projective connection.
At the end of this section, we provide  the bijective correspondences  in the above theorem restricted  to subsheaves classifying various objects of prescribed  radii.

Suppose that  the log curve $U^\mr{log}/S^\mr{log}$ comes from a local  $r (\geq 0)$-pointed stable curve $\msU := (f : U \migi S, \{ \sigma^U_i : S_i \migi U \}_{i=1}^r)$ over a $k$-scheme $S$ together with  a choice of its  base.

\subsection{Monodromy  and radius of a flat vector bundle} \label{QW260}

Let $(\mcV, \nabla)$ be an $\hslash$-flat bundle on $U^\mr{log}/S^\mr{log}$.
Also, we suppose that $r>0$ and   fix $i \in \{1, \cdots, r \}$.

\bde[] \label{QQ1201}
 The {\bf monodromy (operator)} 
 \index{$\nabla$, $\hslash$-$T^\mr{log}$-connection@--- $\mu^\nabla_i$, monodromy (operator)}
  of $\nabla$ at $\sigma^{U}_i$
 is the element  $\mu^\nabla_i$ of $\mr{End}(\sigma^{U*}_i (\mcV))$ defined to be the composite
 \begin{align}
 \mu^\nabla_i : \sigma_i^{U*}(\mcV) 
 &\xrightarrow{\sigma_i^{U*}(\nabla)}   \sigma^{U*}_i (\Omega \otimes \mcV) \\
& \isom  \sigma^{U*}_i(\Omega) \otimes \sigma^{U*}_i (\mcV) \notag \\
 &\xrightarrow{\mr{Res}_{\msU, i} \otimes \mr{id}_{\sigma^{U*}_i (\mcV)}}
 \left( \mcO_{S_i} \otimes \sigma^{U*}_i (\mcV) = \right)
\sigma^{U*}_i (\mcV). \notag 
 \end{align}
\ede

Next, let $s_j$ ($j=1, \cdots, n$) are the maps $\mr{End}(\sigma^{U*}_i (\mcV)) \migi \Gamma (S_i, \mcO_{S_i})$
 assigning the coefficients of characteristic polynomials, i.e.,  
 determined uniquely by
\begin{align} \label{QW278}
\mr{det}(t \cdot \mr{id}_{\sigma^{U*}_i (\mcV)} - x) = t^n + \sum_{j=1}^n (-1)^j \cdot  s_j (x) \cdot t^{n-j}
\end{align}
for any $x \in \mr{End}(\sigma^{U*}_i (\mcV))$.
That is to say,  we set $s_j(x) := \mr{tr}(\bigwedge^j x)$, which is a  homogenous polynomial map of degree $j$ from $\mr{End}(\sigma^{U*}_i (\mcV))$ to $\Gamma (S_i, \mcO_{S_i})$.
Write 
\begin{align} \label{QR4}
\rho_i^\nabla := (s_2 (\breve{\mu}_i^\nabla), s_3 (\breve{\mu}_i^\nabla), \cdots, s_n (\breve{\mu}_i^\nabla )) \in \Gamma (S_i, \mcO_{S_i})^{\oplus (n-1)} \left(= \mbA^{n-1} (S_i) \right),
\end{align}
where $\breve{\mu}_i^\nabla := \mu_i^\nabla - \frac{\mr{tr}(\mu_i^\nabla)}{n}\cdot \mr{id}_{\sigma^{U*}_i (\mcV)}$. (Note that $s_1 (\breve{\mu}_i^\nabla) =0$.)
We  regard  $\rho_i^\nabla$ as an element of $\mfc_{\mfs \mfl_n} (S_i)$ under the identification $\mbA^{n-1} = \mfc_{\mfs \mfl_n}$ given by  $s_j = s_j^{\mft}$ (cf. (\ref{QW304})).

\bde \label{QR3}
We shall refer to $\rho^\nabla_i$ as the {\bf radius} 
\index{$\nabla$, $\hslash$-$T^\mr{log}$-connection@--- $\rho^\nabla_i$, radius of $\nabla$}
 of $\nabla$ at $\sigma^{U}_i$.
\ede

\begin{rema} 
\label{QW306}
Let $(\mcV, \nabla)$ be as above,
 and denote by $(\mcE, \nabla_\mcE)$  the  $\hslash$-flat $\mr{PGL}_n$-bundle induced by $(\mcV, \nabla)$ via projectivization.
Then,  by taking account of 
 the description of $\chi : \mfs \mfl_n \migisurj \mfc_{\mfs \mfl_n}  \left(=  \mfc_{\mfp \mfg \mfl_n} \right)$ discussed in  Example \ref{QQ1290}, we see that 
\begin{align} \label{QQ1219}
\rho^{\nabla}_i = \rho^{\nabla_\mcE}_i
\end{align}
for every  $i=1, \cdots, r$.
\end{rema}

\begin{rema}
By choosing another base for $\msU$, we obtain another  log curve $U^{\mr{log}'}/S^{\mr{log}'}$ extending $U/S$.
If  $\omega_{U/S}$ denotes  the dualizing sheaf on  $U$ relative to $S$, then the line bundle $\Omega_{U/S} (D_\msU)$ is isomorphic to $\Omega_{U^\mr{log}/S^\mr{log}}$, as well as isomorphic to $\Omega_{U^{\mr{log}'}/S^{\mr{log}'}}$
 (cf. the comment preceding Remark \ref{QG19}).
 In particular,  there exists an isomorphism  $\Omega_{U^\mr{log}/S^\mr{log}} \isom \Omega_{U^{\mr{log}'}/S^{\mr{log}'}}$ compatible with the residue maps.
 This implies  that the definitions of monodromy and radius associated to 
   an $\hslash$-flat bundle do not depend on  the choice of a base for $\msU$.
\end{rema}

\subsection{Invariance under gauge transformations} \label{QW300}
Let  $\mcV'$ be another vector bundle and $\eta : \mcV \isom \mcV'$ an $\mcO_U$-linear isomorphism.
According to (\ref{QW299}), we obtain an $\hslash$-flat bundle $(\mcV', \eta_*(\nabla))$ on $U^\mr{log}/S^\mr{log}$.
Then, for every $i=1, \cdots, r$,  we can verify that 
\begin{align}
\mu_i^{\eta_*(\nabla)} = \sigma_i^{U*}(\eta) \circ \mu^{\nabla}_i \circ \sigma_i^{U*}(\eta)^{-1} \ \left(\Longrightarrow \ \breve{\mu}_i^{\eta_*(\nabla)} = \sigma_i^{U*}(\eta) \circ \breve{\mu}^{\nabla}_i \circ \sigma_i^{U*}(\eta)^{-1}\right).
\end{align}
 This implies  
\begin{align} \label{QW300}
\rho^{\eta_*(\nabla)}_i =\rho^\nabla_i.
\end{align}
In particular, the radii of an $\hslash$-flat bundle depend only on its isomorphism class.

\subsection{Tensor product and dual}\label{QW498}
Let $(\mcL, \nabla_\mcL)$ be an $\hslash$-flat line bundle on $U^\mr{log}/S^\mr{log}$.
Then,  the monodromy $\mu_i^{\nabla \otimes \nabla_\mcL}$ of  the  $\hslash$-$S^\mr{log}$-connection $\nabla \otimes \nabla_\mcL$ 
(cf. (\ref{QQ295})) satisfies 
\begin{align}
\mu_i^{\nabla \otimes \nabla_\mcL} = \mu_i^\nabla \otimes \mr{id}_{\sigma_i^{U*}(\mcL)} + \mr{id}_{\sigma_i^{U*}(\mcV)} \otimes \mu_i^{\nabla_\mcL} \left(\Longrightarrow  \breve{\mu}_i^{\nabla \otimes \nabla_\mcL}  =  \breve{\mu}_i^\nabla \otimes \mr{id}_{\sigma_i^{U*}(\mcL)}\right)
\end{align}
 for every $i =1, \cdots, r$.
Hence, we have  
\begin{align} \label{QW499}
\rho^{\nabla \otimes \nabla_\mcL}_i =  \rho_i^{\nabla}.
\end{align}

On the other hand, the monodromy $\mu_i^{\nabla^\vee}$ of  the dual $\nabla^\vee$  (cf. (\ref{QQ297})) satisfies  $\mu_i^{\nabla^\vee} = - {^t}\mu^\nabla$.
Hence, by  the definitions of radius and the action $\star$ on $\mfc_{\mfs \mfl_n}$, the equality 
\begin{align} \label{QW2001}
\rho_i^{\nabla^\vee} = (-1) \star \rho_i^\nabla
\end{align}  
holds for every $i=1, \cdots, r$.

\subsection{Radius of a $(\mr{GL}_n, \hslash)$-oper} \label{QG902}

Let us fix an element  $\rho := (\rho_i)_{i=1}^r \in \prod_{i=1}^r \mfc_{\mfs \mfl_n} (S_i)$ $\left(=\prod_{i=1}^r \mbA^{n-1}(S_i)\right)$.
If $r =0$, then we set $\prod_{i=1}^r \mfc_{\mfs \mfl_n} (S_i) := \{ \emptyset \}$ and  $\rho := \emptyset$.

\bde \label{QW341}
Let $\BB$ be a line bundle on $U$.
We shall say that a $(\mr{GL}_n, \hslash)$-oper $\msF^\heartsuit:= (\mcF, \Nab, \{ \mcF^j \}_j)$ (resp., a $(\mr{GL}_n, \hslash, \BB)$-oper $\nabla^\diamondsuit$) on $\msU$ is {\bf of radii $\rho$}
\index{$\mcF^\heartsuit$, $(\mr{GL}_n, \hslash)$-oper@--- of radii $\rho$}
\index{$\nabla^\diamondsuit$, $(\mr{GL}_n, \hslash, \BB)$-oper@--- of radii $\rho$}
  if $\rho^\nabla_i = \rho_i$ (resp., $\rho^{\nabla^\diamondsuit}_i = \rho_i$) for every $i =1, \cdots, r$.
\ede

\begin{rema}
Let $\msF^\heartsuit$ be a $(\mr{GL}_n, \hslash)$-oper on $\msU$ and $\msL := (\mcL, \nabla_\mcL)$  an $\hslash$-flat line bundle on $U^\mr{log}/S^\mr{log}$.
Then,  one verifies from (\ref{QW499})
 that,  if $\msF^\heartsuit$ is of radii $\rho$, then so is the $(\mr{GL}_n, \hslash)$-oper $\msF^\heartsuit_{\otimes \msL}$ (cf. (\ref{QQ111})).
In other words, the radii of a $(\mr{GL}_n, \hslash)$-oper  depend only on its  equivalence class.  
\end{rema}

Let us fix an $(n, \hslash)$-theta characteristic $\bb := (\BB, \nabla_\bb)$ of $U^\mr{log}/S^\mr{log}$.
We shall denote by 
 \index{$\mcO p^\spadesuit_{n, \hslash (, \rho)}$, $\mcO p^{\spadesuit, \mr{O}}_{n, \hslash}$, $\mcO p^{\spadesuit, \mr{Sp}}_{n, \hslash}$}
 \index{$\overline{\mcO} p^\heartsuit_{n, \hslash, \rho}$}
  \index{$\mcO p_{n, \hslash, \BB}^\diamondsuit$, $\mcO p_{n, \hslash, \bb (, \rho)}^\diamondsuit$, $\mcO p_{n, \hslash, \BB}^{\diamondsuit, \mr{O}}$, $\mcO p_{n, \hslash, \BB}^{\diamondsuit, \mr{Sp}}$}
\begin{align} \label{QQ1229}
\mcO p^\spadesuit_{n, \hslash, \rho}
 \ \left(\text{resp.,} \ \overline{\mcO} p^\heartsuit_{n, \hslash, \rho}; \text{resp.,} \ \mcO p^\diamondsuit_{n, \hslash, \bb, \rho}\right)
 \end{align}
the set-valued \'{e}tale sheaf  on $U$ which, to any \'{e}tale $U$-scheme $U'$,
assigns the set of isomorphism classes of $(\mfs \mfl_n, \hslash)$-opers (resp., equivalence classes of $(\mr{GL}_n, \hslash)$-opers; resp.,  isomorphism classes of $(\mr{GL}_n, \hslash, \bb)$-opers) on 
$U' \times_U \msU$
 of radii $(\rho_i |_{S_i\times_{\sigma_i^U, U} U'})_{i=1}^r$;
it may be regarded as 
a subsheaf of  $\mcO p^\spadesuit_{n, \hslash}$  (resp., $\overline{\mcO} p^\heartsuit_{n, \hslash}$; resp., $\mcO p^\diamondsuit_{n, \hslash, \bb}$).
Then, the isomorphism  $\Lambda^{\heartsuit \Rightarrow \spadesuit}$ (resp., $\Lambda_\bb^{\diamondsuit \Rightarrow \heartsuit}$) restricts to  an isomorphism of sheaves
\begin{align} \label{Qw501}
\Lambda_\rho^{\heartsuit \Rightarrow \spadesuit} : \mcO p^\spadesuit_{n, \hslash, \rho} \isom \overline{\mcO} p_{n, \hslash, \rho}^\heartsuit
\
\left(\text{resp.,} \ 
\Lambda_{\bb, \rho}^{\diamondsuit \Rightarrow \heartsuit} : \mcO p^\diamondsuit_{n, \hslash, \bb, \rho} \isom \overline{\mcO} p_{n, \hslash, \rho}^\heartsuit \right).
\index{$\Lambda^{\heartsuit \Rightarrow \spadesuit}_{(\rho)}$, $\Lambda^{\heartsuit \Rightarrow \spadesuit, \BL}$, $\Lambda^{\heartsuit \Rightarrow \spadesuit, \BL}_\msX$}
\index{$\Lambda_{\bb (, \rho)}^{\diamondsuit \Rightarrow \heartsuit}$, $\Lambda_\bb^{\diamondsuit \Rightarrow \heartsuit, \BL}$, $\Lambda_{\bb, \msX}^{\diamondsuit \Rightarrow \heartsuit, \BL}$}
\end{align}

\subsection{Radius of an $(n, \hslash, \bb)$-projective connection}\label{QW4118}
Let $i \in \{ 1, \cdots, r \}$.
Note that the two $\mcO_{S_i}$-modules
$\sigma^{U*}_i ({^L}\mcD_\hslash^{<n+1})$ and
$\sigma^{U*}_i ({^R}\mcD_\hslash^{<n+1})$
are naturally identified.
 So we will not distinguish these $\mcO_{S_i}$-modules and denote them by $\sigma^{U*}_{i} (\mcD_\hslash^{<n+1})$ for simplicity.
Let us write $\overline{\partial}$ for  the section of $\sigma^{U*}_i (\mcT)$ corresponding to $1$ via the dual  $\mr{Res}_{\msU, i}^\vee :   \mcO_{S_i} \isom \sigma^{U*}_i(\mcT)$ of $\mr{Res}_{\msU, i}$.
The assignment $(a_0, \cdots, a_{n}) \mapsto \sum_{j=0}^{n} a_{j} \cdot \overline{\partial}^{n-j}$ determines a well-defined  isomorphism
\begin{align} \label{QW390}
\mcO_{S_i}^{\oplus n+1} \isom \sigma^{U*}_i (\mcD_\hslash^{<n+1}).
\end{align}
Then, we have the following composite isomorphism:
\begin{align} \label{QW400}
\sigma^{U*}_i (\DD^{\,<n+1}_{\BB}) &\isom \mcH om_{\mcO_{S_i}} (\sigma^{U*}_i(\BB^\vee), \sigma^{U*}_i (\Omega^{\otimes n} \otimes \BB^\vee) \otimes \sigma^{U*}_i (\mcD_\hslash^{<n+1})) 
\\
& \isom \mcH om_{\mcO_{S_i}} (\sigma^{U*}_i (\BB^\vee), \sigma^{U*}_i (\BB^\vee) \otimes \mcO_{S_i}^{\oplus n+1}) \notag \\
& \isom \mcO_{S_i}^{\oplus n+1}, \notag
\end{align}
where the second arrow arises from $\mr{Res}_{\msU, i}^{\otimes n} : \sigma^{U*}_i (\Omega)^{\otimes n} \isom \mcO_{S_i}$ and (\ref{QW390}).
If $D^\clubsuit$ is an $(n, \hslash, \bb)$-projective connection on $\msU$,
then, for each $i=1, \cdots, r$,  there exists 
an element
\begin{align} \label{QH8}
\rho^{D^\clubsuit}_i := (\rho^{D^\clubsuit}_{i2}, \cdots, \rho^{D^\clubsuit}_{in}) \in \Gamma (S_i, \mcO_{S_i})^{\oplus (n-1)}  \left(=  \mbA^{n-1}(S_i)\right)
\end{align}
 such that the image of $\sigma_i^{U*}(D^\clubsuit)$ via (\ref{QW400}) coincides with $(1, -\mu^{\nabla_\bb}, ( (-1)^j \cdot \rho^{D^\clubsuit}_{ij})_{j=2}^n)$.


\bde \label{QW342}
 Let $D^\clubsuit$ be  an $(n, \hslash, \bb)$-projective connection on $\msU$.
Then, we shall say that $D^\clubsuit$ is {\bf of radii $\rho$}
\index{$D^\clubsuit$, $(n, \hslash, \bb)$-projective connection@--- of radii $\rho$}
if  the equality  $\rho_i^{D^\clubsuit} = \rho_i$ holds  for every $i =1, \cdots, r$ under the identification $\mbA^{n-1} = \mfc_{\mfs \mfl_n}$ given by $s_j = s_j^{\mft}$ (cf. (\ref{QW305})).
\ede

\subsection{Comparison} \label{QG930}
Denote by
\index{$\DD^{\,\clubsuit}_{n, \hslash, \BB}$, $\DD^{\,\clubsuit}_{n, \hslash, \bb (, \rho)}$, $\mfD \mfi \mff \mff^{\,\clubsuit}_{n, \hslash, \bb (, \rho)}$}
\begin{align}
\DD^{\,\clubsuit}_{n, \hslash, \bb, \rho}
\end{align}
the subshear of $\DD^{\,\clubsuit}_{n, \hslash, \bb}$ consisting of $(n, \hslash, \bb)$-projective connections  of radii $\rho$.
By means of  the local description (\ref{QQ325}),
we see that $\DD^{\,\clubsuit}_{n, \hslash, \bb, \rho}$ admits a section locally on $U$, and 
 the $\DV_\BB$-action on $\DD^{\,\clubsuit}_{n, \hslash, \bb}$ restricts to  a structure of $\DV_\BB (-D_\msU)$-torsor on $\DD^{\,\clubsuit}_{n, \hslash, \bb, \rho}$.
Moreover, it follows from the various definitions involved that
the isomorphism $\Lambda_\bb^{\clubsuit \Rightarrow \diamondsuit}$ (cf. (\ref{GL28})) restricts to an isomorphism of sheaves
\index{$\Lambda_{\BB}^{\clubsuit \Rightarrow \diamondsuit}$, $\Lambda_{\bb (, \rho)}^{\clubsuit \Rightarrow \diamondsuit}$}
\begin{align} \label{QW402}
\Lambda_{\bb, \rho}^{\clubsuit \Rightarrow \diamondsuit} : \DD^{\,\clubsuit}_{n, \hslash, \bb, \rho} \isom \mcO p^\diamondsuit_{n, \hslash, \bb, \rho}.
\end{align}
This implies the following assertion.

\bt \label{QG1700}
Let $\bb$ and $\msU$ be as above.
Then, there exists a canonical sequence of isomorphisms
 \begin{align} \label{GL42radius}
 \vcenter{\xymatrix@C=46pt@R=36pt{
\DD^{\, \clubsuit}_{n,\hslash, \bb, \rho} \ar[rd]^-{\sim}_-{\Lambda_{\bb, \rho}^{\clubsuit \Rightarrow \diamondsuit}}&& \overline{\mcO}p^\heartsuit_{n,\hslash, \rho}\ar[rd]^-{\sim}_-{\Lambda_\rho^{\heartsuit \Rightarrow \spadesuit}}&\\
& \mcO p^\diamondsuit_{n, \hslash, \bb, \rho}\ar[ur]_-{\sim}^-{\Lambda_{\bb, \rho}^{\diamondsuit \Rightarrow \heartsuit}}&& \mcO p^\spadesuit_{n,\hslash, \rho}.
}}
 \end{align}
 In particular, the sheaf  $\mcO p^\spadesuit_{n,\hslash, \rho}$ has a $\DV_\BB (-D_\msU)$-torsor structure  transported from  $\DD^{\,\clubsuit}_{n, \hslash, \bb, \rho}$ via this sequence.
 \et

\section{Change of $(n, \hslash)$-theta characteristic} \label{QQ1155}

In this section, we study the transformations of $(n, \hslash, \bb)$-projective connections  and $(\mr{GL}_n, \hslash, \bb)$-opers induced by changing the  $(n, \hslash)$-theta characteristic $\bb$.
Let  $S^\mr{log}$ be an fs log scheme  over $k$,  $U^\mr{log}$  a log curve over $S^\mr{log}$, and $\hslash$ an element of $\Gamma (S, \mcO_S)$.
Also, let  $\bb := (\BB, \nabla_\bb)$ and $\bb' := (\BB', \nabla_{\bb'})$ be  $(n, \hslash)$-theta characteristics of $U^\mr{log}/S^\mr{log}$.
Then, it follows from Lemma \ref{y0122} that there exists   a unique (up to isomorphism) $\hslash$-flat line bundle $\msL := (\mcL, \nabla_\mcL)$ ($= \bb'/\bb$ in the notation of that lemma) on $U^\mr{log}/S^\mr{log}$ with $\bb \otimes \msL \cong  \bb'$.

\subsection{The case of $(n, \hslash, \bb)$-projective connections} \label{QW2}

 Let $\msN := (\mcN, \nabla_{\mcN})$ be an $\hslash$-flat line bundle   on $U^\mr{log}/S^\mr{log}$.
 Then, the tensor product $\mcN \otimes \mcD^{<\infty}_{\hslash}$ admits
 a structure of left $\mcD_\hslash^{< \infty}$-module 
 \begin{align} \label{QQ1110}
 \mcD_\hslash^{<\infty} \otimes (\mcN \otimes \mcD_\hslash^{<\infty}) \migi  \mcN \otimes \mcD_\hslash^{<\infty}
 \end{align}
 arising from both $\nabla_{\mcN}$ and the left $\mcD_\hslash^{< \infty}$-action  on $\mcD_\hslash^{< \infty}$ itself.
 To be precise,  it is 
  determined uniquely by $\partial (v \otimes D) = (\nabla_{\mcN})_\partial (v)\otimes D  + v \otimes (\partial * D)$ for any local sections $v \in \mcN$, $D \in \mcD_\hslash^{<\infty}$,  and $\partial \in \mcT \left(\subseteq \mcD_\hslash^{<\infty} \right)$.
For  a nonnegative integer $j$,   (\ref{QQ1110}) restricts to an $\mcO_U$-linear isomorphism
 \begin{align}  \label{QQ1116}
\mr{sw}_{\msN}^{j} :   \mcD_{\hslash}^{<j} \otimes \mcN \left(= \mcD_{\hslash}^{<j} \otimes (\mcN \otimes \mcD_\hslash^{<1}) \right)  \isom \mcN \otimes \mcD_\hslash^{<j}.
\index{$\mr{sw}_{\msN}^{j}$}
 \end{align}
 If $\msN' := (\mcN', \nabla_{\mcN'})$ is another $\hslash$-flat line bundle on $U^\mr{log}/S^\mr{log}$, then the  following diagram is commutative:
 \begin{align} \label{QQ1130}
 \vcenter{\xymatrix@C=26pt@R=36pt{
 \mcD_\hslash^{<j} \otimes (\mcN \otimes \mcN') \ar[rd]^-{\sim}_-{\mr{sw}^j_{\msN} \otimes \mr{id}_{\mcN'}} \ar[rr]_-{\sim}^-{\mr{sw}^j_{\msN \otimes \msN'}}&& \mcN \otimes \mcN' \otimes \mcD_\hslash^{<j} \\
&  \mcN \otimes \mcD_\hslash^{<j} \otimes \mcN' \ar[ur]^-{\sim}_-{\mr{id}_{\mcN} \otimes \mr{sw}^j_{\msN'}}.&
}}
 \end{align}

Now, let us  write $\msL^\vee :=(\mcL^\vee, \nabla_\mcL^\vee)$.
Given a global section 
$D$ of ${^\Join}\DD^{\,<n+1}_\BB$,
we shall  denote by  $D_{\otimes \msL}$ the section of ${^\Join}\DD^{\,<n+1}_{\BB \otimes \mcL}$ defined to be   the composite
\begin{align}
D_{\otimes \msL} :  (\BB \otimes \mcL)^\vee & \xrightarrow{\sim} \BB^\vee \otimes \mcL^\vee \\
& \xrightarrow{D \otimes \mr{id}_\mcL}
(\Omega^{\otimes n} \otimes \BB^\vee) \otimes \mcD_\hslash^{<n+1} \otimes \mcL^\vee \notag \\
& \xrightarrow{\mr{id}_{\Omega^{\otimes n} \otimes \BB^\vee}\otimes \mr{sw}_{\msL^\vee}^{n+1}} (\Omega^{\otimes n} \otimes \BB^\vee) \otimes \mcL^\vee \otimes \mcD_\hslash^{<n+1} \notag \\
& \xrightarrow{\sim}
(\Omega^{\otimes n} \otimes (\BB \otimes \mcL)^\vee) \otimes \mcD_\hslash^{<n+1}. \notag
\end{align}
The assignment $D \mapsto D_{\otimes \msL}$ commutes with pull-back to \'{e}tale $U$-schemes and hence induces a morphism
\begin{align} \label{QW8}
\clubsuit^+_{\bb \Rightarrow \bb'} :  {^\Join}\DD^{\,<n+1}_\BB \migi {^\Join}\DD^{\,<n+1}_{\BB \otimes \mcL}.
\end{align}
Since $(D_{\otimes \msL})_{\otimes \msL^\vee} = D$ and $(D_{\otimes \msL^\vee})_{\otimes \msL} = D$, this morphism is an isomorphism.
Moreover, it restricts to an isomorphism
\begin{align}
\clubsuit_{\bb \Rightarrow \bb'}
  : \DD^{\,\clubsuit}_{n,  \bb} \isom \DD^{\,\clubsuit}_{n, \bb'}.  
  \index{$\clubsuit_{\bb \Rightarrow \bb' (,\rho)}$}
    \end{align}
The formation of this  isomorphism commutes with base-change to fs log schemes over $S^\mr{log}$.

\subsection{The case of $(\mr{GL}_n, \hslash, \bb)$-opers} \label{QW3}

Next, let $\nabla^\diamondsuit$ be 
  a $(\mr{GL}_n, \hslash, \BB)$-oper  on $U^\mr{log}/S^\mr{log}$.
  We shall denote by 
  \index{$\nabla^\diamondsuit$, $(\mr{GL}_n, \hslash, \BB)$-oper@--- $\nabla^\diamondsuit_{\otimes \msL}$}
  \begin{equation} \label{QQ1131}
 \nabla^\diamondsuit_{\otimes \msL}
 \end{equation}
the  $\hslash$-$S^\mr{log}$-connection on 
$\DF_{\BB \otimes \mcL}$
corresponding to  $\nabla_\mcL \otimes \nabla^\diamondsuit$ via 
the composite isomorphism
\begin{align}
\DF_{\BB \otimes \mcL} \left(=\mcD_\hslash^{<n} \otimes \BB \otimes \mcL\right) &\xrightarrow{a \otimes b \otimes c \mapsto a \otimes c \otimes b} \mcD_\hslash^{<n} \otimes \mcL \otimes \BB \\
& \xrightarrow{\mr{sw}_\msL^n \otimes \mr{id}_\BB}
\left( \mcL \otimes \mcD_\hslash^{<n} \otimes \BB =\right) \mcL \otimes \DF_\BB. \notag
\end{align}
It is immediately verified that  $\nabla^\diamondsuit_{\otimes \msL}$ forms a $(\mr{GL}_n, \hslash, \BB\otimes \mcL)$-oper on $U^\mr{log}/S^\mr{log}$. 
The assignment $\nabla^\diamondsuit \mapsto \nabla^\diamondsuit_{\otimes \msL}$ commutes with pull-back to \'{e}tale $U$-schemes and hence induces a morphism
\begin{align} \label{QW9}
\diamondsuit^+_{\bb \Rightarrow \bb'} : \mcO p^\diamondsuit_{h, \hslash, \BB} \migi \mcO p_{n, \hslash, \BB \otimes \mcL}^\diamondsuit.
\end{align}
Since $(\nabla^\diamondsuit_{\otimes \msL})_{\otimes \msL^\vee} = \nabla^\diamondsuit$ and $(\nabla^\diamondsuit_{\otimes \msL^\vee})_{\otimes \msL} = \nabla^\diamondsuit$,
this  morphism 
is  an isomorphism.
Moreover, it restricts to an isomorphism of sheaves
 \begin{equation}  
 \label{GL38}
 \diamondsuit_{\bb \Rightarrow \bb'} : \mcO p_{n, \hslash, \bb}^\diamondsuit  \isom   \mcO p_{n, \hslash, \bb'}^\diamondsuit,  \index{$\diamondsuit_{\bb \Rightarrow \bb' (, \rho)}$}
  \end{equation}
 whose formation commutes with base-change to fs log schemes over $S^\mr{log}$.
 

\bpr \label{y0118}
Let us keep the above notation.
Then,  for an $(n, \hslash, \bb)$-projective connection  $D^\clubsuit$ on $U^\mr{log}/S^\mr{log}$,
the following  equality holds:
\begin{align} \label{QW31}
(D^{\clubsuit \Rightarrow \diamondsuit})_{\otimes \msL} = 
(D^\clubsuit_{\otimes \msL})^{\Rightarrow \diamondsuit}.
\end{align}
In particular, 
the isomorphisms
 $\clubsuit_{\bb \Rightarrow \bb'}$ and $\diamondsuit_{\bb \Rightarrow \bb'}$  make the following diagram commute:
\begin{align} \label{QQ246}
\vcenter{\xymatrix@C=50pt@R=10pt{
\DD^{\, \clubsuit}_{n,\hslash, \bb} \ar[r]_-{\sim}^-{\Lambda^{\clubsuit \Rightarrow \diamondsuit}_{\bb}} \ar[dd]^-{\wr}_-{\clubsuit_{\bb \Rightarrow \bb'}}& \mcO p^\diamondsuit_{n, \hslash, \bb} \ar[rd]_-{\sim}^-{\Lambda_{\bb}^{\diamondsuit \Rightarrow \heartsuit}} \ar[dd]^-{\wr}_-{\diamondsuit_{\bb \Rightarrow \bb'}}&\\
 && \overline{\mcO} p_{n, \hslash}^\heartsuit.\\
 \DD^{\, \clubsuit}_{n,\hslash, \bb'} \ar[r]^-{\sim}_-{\Lambda^{\clubsuit \Rightarrow \diamondsuit}_{\bb'}}& \mcO p^\diamondsuit_{n, \hslash, \bb'}  \ar[ur]^-{\sim}_-{\Lambda_{\bb'}^{\diamondsuit \Rightarrow \heartsuit}}&
}}
\end{align}
\epr
\begin{proof}
The assertion follows
 from the definitions of various morphisms in (\ref{QQ246}).
\end{proof}

\subsection{The compatibility of torsor structures} \label{QW6}
The $\DV_\BB$-action on  $\DD^{\, \clubsuit}_{n,\hslash, \bb}$ (resp., $\mcO p^\diamondsuit_{n, \hslash, \bb}$) is transformed into the $\DV_{\BB'}$-action on
$\DD^{\, \clubsuit}_{n,\hslash, \bb'}$ (resp., $\mcO p^\diamondsuit_{n, \hslash, \bb'}$) via $\clubsuit_{\bb \Rightarrow \bb'}^+$ and $\clubsuit_{\bb \Rightarrow \bb'}$ (resp., $\diamondsuit_{\bb \Rightarrow \bb'}$).
To be precise, for  local sections $D^\clubsuit \in \DD^{\, \clubsuit}_{n,\hslash, \bb}$ (resp., $\nabla^\diamondsuit \in \mcO p^\diamondsuit_{n, \hslash, \bb}$), $R \in \DV_\BB$,  we have
\begin{align} \label{QQ1171}
\clubsuit_{\bb \Rightarrow \bb'} (D^\clubsuit + R) = \clubsuit_{\bb \Rightarrow \bb'} (D^\clubsuit) + R'
 \ \left(\text{resp.,} \ \diamondsuit_{\bb \Rightarrow \bb'} (\nabla^\diamondsuit_{+R}) = \diamondsuit_{\bb \Rightarrow \bb'} (\nabla^\diamondsuit)_{+R'} 
\right), \hspace{-5mm}
\end{align}
 where $R' := \clubsuit_{\bb \Rightarrow \bb'}^+ (R)$.

Next, let us focus on  the case of $n=2$.
Since $\mcD_\hslash^{<1} = \mcO_X$,
  both  the sheaves
  $\DV_\BB$ $\left(= \mcH om (\BB^\vee, \Omega^{\otimes 2} \otimes \BB^\vee) \right)$ and 
  $\DV_{\BB'} \left(= \mcH om ((\BB \otimes \mcL)^\vee, \Omega^{\otimes 2} \otimes (\BB \otimes \mcL)^\vee) \right)$
    may be naturally identified with
  $\Omega^{\otimes 2}$. 
Moreover, under these identifications, the isomorphism $\clubsuit_{\bb \Rightarrow \bb'}^+$ coincides with the identity morphism of $\Omega^{\otimes 2}$.
It follows (from (\ref{QQ1171})) that  both $\clubsuit_{\bb \Rightarrow \bb'}$ and $\diamondsuit_{\bb \Rightarrow \bb'}$  are  compatible with the respective  $\Omega^{\otimes 2}$-actions  on the domain and codomain.
Now, let us equip $\mcO p_{2, \hslash}^\spadesuit$ with the  $\Omega^{\otimes 2}$-torsor structure
 coming from $\mcO p^\diamondsuit_{2, \hslash, \bb}$ via  $\Lambda^{\heartsuit \Rightarrow \spadesuit} \circ \Lambda^{\diamondsuit \Rightarrow \heartsuit}_{\bb}$.
The above discussion  implies that
this torsor structure does not depend on the choice of $\bb$.
Therefore, 
by gluing together various local constructions, 
one  can obtain an $\Omega^{\otimes 2}$-torsor structure 
\begin{align} \label{QW12}
\mcO p_{2, \hslash}^\spadesuit \times \Omega^{\otimes 2} \migi \mcO p^\spadesuit_{2, \hslash}
\end{align}
on $\mcO p_{n, \hslash}^\spadesuit$
{\it even when there are no  $(2, \hslash)$-theta characteristics defined on the entire space $U$}.


\subsection{Restriction by radii}
Let us keep the notation in  \S\,\ref{QW274}.
Also, let $\rho$ be an element of   $\prod_{i=1}^r \mfc_{\mfs \mfl_n} (S_i)$, where $\prod_{i=1}^r \mfc_{\mfs \mfl_n} (S_i) := \{ \emptyset \}$ and $\rho := \emptyset$ if $r=0$.
Then, 
$\clubsuit_{\bb \Rightarrow \bb'}$ and  $\diamondsuit_{\bb \Rightarrow \bb'}$ \index{$\diamondsuit_{\bb \Rightarrow \bb' (, \rho)}$} restrict to  isomorphisms
\begin{align} \label{QW455}
\clubsuit_{\bb \Rightarrow \bb', \rho} : \DD^{\,\clubsuit}_{n, \hslash, \bb, \rho} \isom \DD^{\,\clubsuit}_{n, \hslash, \bb', \rho} 
\index{$\clubsuit_{\bb \Rightarrow \bb' (,\rho)}$}
 \ \ \text{and} \ \  
 \diamondsuit_{\bb \Rightarrow \bb', \rho} : \mcO p^\diamondsuit_{n, \hslash, \bb, \rho} \isom  \mcO p^\diamondsuit_{n, \hslash, \bb', \rho}
\end{align}
respectively.
Moreover,    (\ref{QQ246}) restricts to the following commutative diagram:
\begin{align} \label{QW458}
\vcenter{\xymatrix@C=50pt@R=10pt{
\DD^{\, \clubsuit}_{n,\hslash, \bb, \rho} \ar[r]_-{\sim}^-{\Lambda^{\clubsuit \Rightarrow \diamondsuit}_{\bb, \rho}} \ar[dd]^-{\wr}_-{\clubsuit_{\bb \Rightarrow \bb', \rho}}& \mcO p^\diamondsuit_{n, \hslash, \bb, \rho} \ar[rd]_-{\sim}^-{\Lambda_{\bb, \rho}^{\diamondsuit \Rightarrow \heartsuit}} \ar[dd]^-{\wr}_-{\diamondsuit_{\bb \Rightarrow \bb', \rho}}&\\
 && \overline{\mcO} p_{n, \hslash, \rho}^\heartsuit.\\
 \DD^{\, \clubsuit}_{n,\hslash, \bb', \rho} \ar[r]^-{\sim}_-{\Lambda^{\clubsuit \Rightarrow \diamondsuit}_{\bb', \rho}}& \mcO p^\diamondsuit_{n, \hslash, \bb', \rho}  \ar[ur]^-{\sim}_-{\Lambda_{\bb', \rho}^{\diamondsuit \Rightarrow \heartsuit}}&
}}
\end{align}

\section{A canonical  $\mr{PGL}_n$-bundle underlying $(\mfs \mfl_n, \hslash)$-opers} \label{y0131}

In this section, we construct a canonical  $\mr{PGL}_n$-bundle  underlying $(\mfs \mfl_n, \hslash)$-opers.
After that, we  define the notion of a {\it normal} $(\mfs \mfl_n, \hslash)$-oper by  using this $\mr{PGL}_n$-bundle.

Let $S^\mr{log}$ be an fs log scheme over $k$ and
 $U^\mr{log}$ a log curve over $S^\mr{log}$.
 
\subsection{Local construction} \label{QW15}
First,  assume that there exists     an $(n, \hslash)$-theta characteristic $\bb := (\BB, \nabla_\bb)$ of  $U^\mr{log}/S^\mr{log}$.
Denote by 
\begin{align} \label{QQ173}
 \DE^{\bb}_{n, \hslash, U^\mr{log}/S^\mr{log}}, \  
 \text{or simply} \ {^\dagger}\mcE^{\bb}_n,
\end{align} 
the $\mr{PGL}_n$-bundle on $U$ induced by $\DF_\BB \left(:=\mcD^{<n}_{\hslash} \otimes \BB\right)$ via projectivization.
Then, the filtration
$\{ \DF^j_\BB\}_{j=0}^n$
 determines a $B$-reduction  of ${^\dagger}\mcE^\bb_n$ (cf. (\ref{BorelB}) for the definition of $B$), which we denote by
\begin{align} \label{QQ174}
 \DE^{\bb}_{B, \hslash, U^\mr{log}/S^\mr{log}},  \ 
 \text{or simply} \ {^\dagger}\mcE^{\bb}_B.
\end{align}
By definition,  the formations of ${^\dagger}\mcE^\bb_n$ and
  ${^\dagger}\mcE_B^\bb$ commutes  with  pull-back  to \'{e}tale  schemes over $U$ and base-change to fs log schemes over $S^\mr{log}$.

 Now, let  $(V, \partial)$ be a log chart on $U^\mr{log}/S^\mr{log}$.
The trivialization  
$\mcO_V \isom \mcT |_V$ of $\mcT |_V$   given by $a \mapsto a \cdot \partial$  yields 
an isomorphism  $(\mcT^{\otimes \frac{n (n-1)}{2}}  \otimes \BB^{\otimes n})|_V\cong \BB |_V^{\otimes n}$.
Then, it follows from Proposition \ref{y0111}, (i), 
 that there exists a unique  $\hslash$-$S^\mr{log}$-connection $\nabla_{\BB |_V}$  on $\BB |_V$ whose  $n$-fold tensor product $\nabla_{\BB |_V}^{\otimes n}$ corresponds to $\nabla_\bb |_V$ via this isomorphism.
We shall write $\msL_{\BB, V}:=(\BB |_V, \nabla_{\BB |_V})$.
Since the $\mr{PGL}_n$-bundle associated to $\bigoplus_{j=0}^{n-1} \BB |_V \otimes (\mcO_V \cdot \partial^j)$ is trivial, 
the  isomorphism
\begin{align} \label{QQ410}
(\mr{sw}^n_{\msL_{\BB, V}})^{-1} : 
\bigoplus_{j=0}^{n-1} \BB |_V \otimes (\mcO_V \cdot \partial^j) \left(= (\BB \otimes \mcD_\hslash^{<n})|_V\right) \isom  
\DF_\BB |_V
\end{align}
(cf. (\ref{QQ1116}) for the definition of $\mr{sw}_{(-)}^n$)
induces, via projectivization, a trivialization
\begin{align} \label{QQw1139}
\mr{triv}^\bb_{n, (V, \partial)} : V \times_k \mr{PGL}_n \isom \DE_n^\bb |_V
\end{align}
of the $\mr{PGL}_n$-bundle $\DE_n^\bb |_V$.
The isomorphism $(\mr{sw}^n_{\msL_{\BB, V}})^{-1}$
 preserves  the  filtrations  $\{ \bigoplus_{j'=0}^{n-1-j} \BB |_V \otimes (\mcO_V \cdot \partial^{j'}) \}_{j=0}^n$ and $\{ \DF_\BB^j\}_{j=0}^n$. 
This implies that $\mr{triv}^\bb_{n, (V, \partial)}$ preserves  the  $B$-reductions, so it  reduces   to 
   a trivialization
\begin{align} \label{QQ1139}
\mr{triv}^\bb_{B, (V, \partial)} : V \times_k B \isom \DE_B^\bb |_V
\end{align}
of the $B$-bundle $\DE_B^\bb |_V$.


\subsection{Gluing isomorphisms} \label{QW25}

Next, let $\bb' := (\BB', \nabla'_0)$ be another $(n, \hslash)$-theta characteristic of $U^\mr{log}/S^\mr{log}$.
According to Lemma \ref{y0122}, there exists an $\hslash$-flat line bundle $\bb' /\bb := (\mcL, \nabla_\mcL)$ on $U^\mr{log}/S^\mr{log}$ with $\bb \otimes (\bb'/\bb) \cong \bb'$.
The composite isomorphism 
\begin{align} \label{QW19}
\DF_{\BB'}\left( = \mcD_\hslash^{<n} \otimes  (\BB \otimes \mcL)  \right) &\xrightarrow{a \otimes b \otimes c \mapsto a \otimes c \otimes b}
\mcD_\hslash^{<n} \otimes \mcL \otimes \BB \\
& \xrightarrow{\mr{sw}_{\bb'/\bb}^n \otimes \mr{id}_{\BB}} \left(\mcL \otimes \mcD_\hslash^{<n} \otimes \BB= \right) \mcL \otimes \DF_\BB\notag
\end{align}
preserves the filtrations $\{ \DF_{\BB'}^j\}_{j}$ and $\{ \mcL \otimes \DF_\BB^j \}_{j}$, so
induces an isomorphism
\begin{align} \label{QW21}
\xi_{\bb'/\bb} : \DE^{\bb'}_B \isom  \DE^{\bb}_B
\end{align}
of $B$-bundles.
The isomorphisms $\xi_{\bb' /\bb}$ defined for various pairs $(\bb', \bb)$ satisfy the cocycle condition with respect to $(n, \hslash)$-theta characteristics, i.e., 
 if $\bb''$ is a third $(n, \hslash)$-theta characteristic, then we have  
 \begin{align} \label{QW20}
 \xi_{\bb'' /\bb} = \xi_{\bb'' /\bb'} \circ \xi_{\bb' /\bb}.
 \end{align}
Also, if $(V, \partial)$ is as above, then 
 the following square diagram is commutative:
\begin{align} \label{QQ401}
\vcenter{\xymatrix@C=96pt@R=36pt{
\displaystyle{\bigoplus_{j=0}^{n-1}}\BB' |_V \otimes (\mcO_V \cdot \partial^j) \ar[r]_-{\sim}^-{(a \otimes b) \otimes c \cdot \partial^j \mapsto b \otimes (a \otimes c \cdot \partial^j)} \ar[d]^-{\wr}_-{(\mr{sw}^n_{\msL_{\BB', V}})^{-1}} & 
\mcL \otimes \left( \displaystyle{\bigoplus_{j=0}^{n-1}}\BB |_V \otimes (\mcO_V \cdot \partial^j)\right)
\ar[d]_-{\wr}^-{\mr{id}_\mcL \otimes (\mr{sw}^n_{\msL_{\BB, V}})^{-1}}
\\
\DF_{\BB'} |_V \ar[r]^-{\sim}_-{(\ref{QW19}) |_V} &  \mcL \otimes 
\DF_\BB
 |_V.
}}
\end{align}
This diagram induces the following commutative  diagram of $B$-bundles:
\begin{align} \label{QQ1140}
\vcenter{\xymatrix@C=46pt@R=36pt{
 & V \times_k B \ar[ld]^-{\sim}_-{\mr{triv}_{B, (V, \partial)}^{\bb'}} \ar[rd]^-{\mr{triv}_{B, (V, \partial)}^\bb}_-{\sim} &   \\
\DE^{\bb'}_B  \ar[rr]^-{\sim}_-{\xi_{\bb' /\bb}} & & \DE^{\bb}_B.
}}
\end{align}

\subsection{Global construction} \label{QW26}

By the above discussion,  a canonical $\mr{PGL}_n$-bundle that we want 
 can be constructed  
without assuming  the existence   of a globally defined $(n, \hslash)$-theta characteristic.
Indeed, it follows from (\ref{QW20}) that the $\mr{PGL}_n$-bundles $\DE_n^{\bb}$ (resp., the $B$-bundles $\DE_B^{\bb}$) for  various locally defined  $(n, \hslash)$-theta characteristics $\bb$ on $U^\mr{log}/S^\mr{log}$ may be glued together to obtain a $\mr{PGL}_n$-bundle  (resp., a $B$-bundle) on the entire space $U$; we shall  denote it  by
\begin{align} \label{QQ0430}
\DE_{n, \hslash, U^\mr{log}/S^\mr{log}} \  \text{or simply} \  \DE_{n}   \left(\text{resp.,} \ \DE_{B, \hslash, U^\mr{log}/S^\mr{log}} \ \text{or simply} \  \DE_{B}    \right).
\end{align}

Also, let  $(V, \partial)$ be a log chart on $U^\mr{log}/S^\mr{log}$.
Then,   by the commutativity of (\ref{QQ1140}),
 the trivializations $\mr{triv}_{n, (V, \partial)}^\bb$  (resp., $\mr{triv}_{B, (V, \partial)}^\bb$) defined  for various $\bb$'s 
   may be glued together to obtain
 a trivialization 
\begin{align} \label{QQ0431}
\mr{triv}_{n, (V, \partial)} : V \times_k \mr{PGL}_n \isom \DE_{n} |_V  \ 
\left(\text{resp.,} \   \mr{triv}_{B, (V, \partial)} : V \times_k B \isom \DE_{B} |_V \right)
\end{align}
of the $\mr{PGL}_n$-bundle $\DE_{n} |_V$  (resp.,  the $B$-bundle  $\DE_B |_V$).

\begin{rema} \label{QW3982}
In this remark, we suppose that $U^\mr{log}/S^\mr{log}$ comes from a local $r (>0)$-pointed prestable curve $\msU := (U/S, \{ \sigma^U_i : S_i \migi U \}_{i=1}^r)$  over $S$ together with  a choice of its base.
Then, for each $i =1, \cdots, r$,  there is  a canonical trivialization of  the $\mr{PGL}_n$-bundle  $\sigma^{U*}_i (\DE_n)$ on  $S_i$.
To construct this trivialization,
we first assume  that there exists an $(n, \theta)$-theta characteristic $\bb := (\BB, \nabla_\bb)$ on $U^\mr{log}/S^\mr{log}$.
Tensoring (\ref{QW390}) with the identify morphism of $\sigma^{U*}_i (\BB)$
gives  an isomorphism $\sigma^{U*}_i (\BB)^{\oplus n} \isom \sigma^{U*}_i (\DF_\BB)$.
This isomorphism preserves the filtrations when we equip 
$\sigma^{U*}_i (\BB)^{\oplus n}$ with 
the standard filtration.
 Hence,  it
gives  an trivialization $S_i \times_k B \isom \sigma^{U*}_i (\DE_B)$
of  the $B$-bundle $\sigma^{U*}_i (\DE_B)$ via projectivization.
By construction, it  
 does not depend on the choice of $\bb$.
Hence, by gluing together  these isomorphisms defined on various open subschemes of $U$ (admitting an $(n, \hslash)$-theta characteristic), 
one can obtain  the trivialization 
\begin{align} \label{QW4562}
\mr{triv}_{B, \sigma^U_i} : S_i \times_k B \isom \sigma^{U*}_i (\DE_B)
\end{align} 
 even when  there does not exist a globally defined  $(n, \hslash)$-theta characteristic.
\end{rema}

\begin{rema} \label{QQ0588}
Let us consider  the $T$-bundle $\DE_{B}\times^B T$ 
(cf. (\ref{BorelB}) for the definition of $T$)
 induced from $\DE_B$ via change of structure group by the quotient $B \migi \left(B /[B, B] =  \right) T$.
This $T$-bundle can be obtained from the rank $n$ vector bundle
 $\mr{gr}(\DF_\BB) := \bigoplus_{j=0}^{n-1}\mr{gr}^j (\DF_\BB)\cong \bigoplus_{j=0}^{n-1} \mcT^{\otimes n-j-1} \otimes \BB$
   via projectivization.
In particular,
the isomorphism class of $\DE_{B}\times^B T$   does not depend on the choice of $\hslash$.
More precisely, if $\hslash'$ is another element of $\Gamma (S, \mcO_S)$, then there is a canonical isomorphism of $T$-bundles
\begin{align} \label{QW29}
\DE_{B, \hslash, U^\mr{log}/S^\mr{log}} \times^B T \isom \DE_{B, \hslash', U^\mr{log}/S^\mr{log}}\times^B T.
\end{align}
\end{rema}

 \subsection{Normal $(\mfs \mfl_n, \hslash)$-opers} \label{QW58}
 We shall introduce  {\it normal} $(\mfs \mfl_n, \hslash)$-opers, which are defined as  $(\mfs \mfl_n, \hslash)$-opers  coming  from $(n, \hslash, \bb)$-projective connections.
  We use those for  $n=2$ to define   the notion of a  {\it $\, \qq$-normal $(\mfg, \hslash)$-opers} for a general $\mfg$, as described   in Definition \ref{QQ0998}.

\bde[]\label{y0132} 
Let $\msE^\spadesuit := (\mcE_B, \nabla)$ be  an $(\mfs \mfl_n, \hslash)$-oper on $U^\mr{log}/S^\mr{log}$.
Then, we shall say that $\msE^\spadesuit$ is {\bf normal} 
\index{$\msE^\spadesuit$, $(\mfg, \hslash)$-oper@--- normal ($(\mfs \mfl_n, \hslash)$-oper)} 
if $\mcE_B = \DE_B$ and, for any \'{e}tale $U$-scheme $V$ and any $(n, \hslash)$-theta characteristic $\bb$ of $V^\mr{log}/S^\mr{log}$,
there exists 
an $(n, \hslash, \bb)$-projective connection 
 $D^\clubsuit$ 
 on $V^\mr{log}/S^\mr{log}$
 satisfying  $\nabla |_V = D^{\clubsuit \Rightarrow \diamondsuit} \times^{\mr{GL}_n}\mr{PGL}_n$.
   \ede

The following assertion  follows immediately  from the definition of normality and   (\ref{QW31}).

\bpr \label{QQ260}
Suppose that there exists an $(n, \hslash)$-theta characteristic  $\bb$ of $U^\mr{log}/S^\mr{log}$.
Then, an $(\mfs \mfl_n, \hslash)$-oper on $U^\mr{log}/S^\mr{log}$ is normal if and only if it is of the form $(\DE_B, D^{\clubsuit \Rightarrow \diamondsuit} \times^{\mr{GL}_n} \mr{PGL}_n)$ for some
$(n, \hslash, \bb)$-projective connection $D^\clubsuit$ on $U^\mr{log}/S^\mr{log}$.
\epr

\bpr \label{QQ0470}
For each $(\mfs\mfl_n, \hslash)$-oper $\msE^\spadesuit$  on $U^\mr{log}/S^\mr{log}$,
there exists  uniquely a pair
\begin{align} \label{QQ0471}
({^\dagger}\msE^\spadesuit, \mr{nor}_{\msE^\spadesuit})
\end{align}
consisting of a normal  $(\mfs \mfl_n, \hslash)$-oper ${^\dagger}\msE^\spadesuit$ on $U^\mr{log}/S^\mr{log}$ and an isomorphism ${^\dagger}\msE^\spadesuit \isom \msE^\spadesuit$ of $(\mfs \mfl_n, \hslash)$-opers.
\epr
\begin{proof}
By  Propositions \ref{y040} and  \ref{y0118},
one may assume, after possibly replacing $U$ with its covering in the \'{e}tale topology,  that there exists an $(n, \hslash)$-theta characteristic of $U^\mr{log}/S^\mr{log}$.
Then,
the assertion follows from (Proposition \ref{y040} again and) Theorem \ref{y0128}.
\end{proof}

\subsection{The filtration on the adjoint bundle of $\DE_n$} \label{QQ172} 
 We shall study  in somewhat detail  the adjoint bundle associated to
 the $\mr{PGL}_n$-bundle $\DE_n$.
 First, we shall assume that there exists an $(n, \hslash)$-theta characteristic $\bb := (\BB, \nabla_\bb)$ of $U^\mr{log}/S^\mr{log}$.
 Denote by $\DE_{\mr{GL}_n}$ the $\mr{GL}_n$-bundle corresponding to the rank $n$ vector bundle $\DF_\BB$.
Then,  
the adjoint bundle  $(\mfg \mfl_n)_{\DE_{\mr{GL}_n}}$  of  $\DE_{\mr{GL}_n}$ is canonically isomorphic to $\mcE nd (\DF_\BB) := \mcE nd_{\mcO_U}(\DF_\BB)$ (cf. (\ref{QW30})). 
 For each $j \in \mbZ$, we shall denote by 
\begin{equation} \label{QQ176}
\mcE nd (\DF_\BB)^j
\end{equation}
 the $\mcO_U$-submodule of 
 $\mcE nd(\DF_\BB)$ consisting of local sections $w$ with  $w (\DF_\BB^{l}) \subseteq \DF_\BB^{l+j}$ for all $l = 0, \cdots, n$.
The collection
\begin{align}
\{ \mcE nd (\DF_\BB)^j\}_{j \in \mbZ}
\end{align}
forms a decreasing filtration on $\mcE nd (\DF_\BB)$ with $\mcE nd (\DF_\BB)^j =0$ for $j \geq n$ and $\mcE nd (\DF_\BB)^j =\mcE nd (\DF_\BB)$ for $j \leq -n$.

 Let us consider the subquotiens of $\mcE nd (\DF_\BB)$ associated to this filtration.
 Fix an integer $j$ with $-n< j <n$ and choose  a local section
 $w$  of $\mcE nd (\DF_\BB)^j$.
For each integer $l$ with $\mr{max}\{0, -j \} \leq l\leq  \mr{min} \{n-1, n-1-j\}$,
 we set 
 \begin{align}
 w_l : \mr{gr}^l (\DF_\BB) \migi \mr{gr}^{l+j} (\DF_\BB)
 \end{align}
 to be an $\mcO_U$-linear morphism induced, via taking  quotients, by the restriction $w |_{\DF^l_\BB} : \DF_\BB^l \migi \DF_\BB^{l+j}$ of $w$.
 Since $\mr{gr}^l (\DF_\BB)$ and $\mr{gr}^{l+j} (\DF_\BB)$ are canonically
 isomorphic to $\mcT^{\otimes (n-l-1)}\otimes \BB$ and $\mcT^{\otimes (n-l-j-1)}\otimes \BB$ respectively, 
 $\mcH om_{\mcO_U} (\mr{gr}^l (\DF_\BB), \mr{gr}^{l+j} (\DF_\BB))$ may be identified with $\Omega^{\otimes j}$.
 Under this identification, 
 the assignment $w \mapsto$ $(w_l)_{l = \mr{max}\{0, -j\}}^{\mr{min}\{n-1, n-1-j \}}$ gives an isomorphism
 \begin{align} \label{QW50}
 \mcE nd (\DF_\BB)^j/\mcE nd (\DF_\BB)^{j+1} \isom  (\Omega^{\otimes j})^{\oplus (n- |j |)}.
 \end{align}

Next, let us study  the adjoint bundle $(\mfs \mfl_n)_{\DE_n}$ of $\DE_n$.
We shall write
\begin{align} \label{QW61}
\overline{\mcE nd}(\DF_\BB) := \mcE nd(\DF_\BB)/\mr{Im}(\Delta),
\end{align}
where $\Delta$ \index{$\Delta$} denotes the diagonal embedding $\mcO_U \migiincl \mcE nd (\DF_\BB)$.
This $\mcO_U$-module has a filtration $\{ \overline{\mcE nd}(\DF_\BB)^j\}_{j \in \mbZ}$ defined by
\begin{equation} \label{QW60}
\overline{\mcE nd}(\DF_\BB)^j := \mcE nd(\DF_\BB)^j/ (\mcE nd(\DF_\BB)^j \cap  \mr{Im} (\Delta)).
\end{equation}
By the definition of $\DE_n$ and the natural identification $\mfs \mfl_n = \mfp \mfg \mfl_n$,
we have  an identification
\begin{align} \label{QW59}
(\mfs \mfl_n)_{\DE_n} =   \overline{\mcE nd}(\DF_\BB), 
\end{align}
which  restricts to  $(\mfs \mfl_n)_{\DE_n}^j  =  \overline{\mcE nd}(\DF_\BB)^j$ for each $j$.
Hence,  (\ref{QW50}) induces   a canonical isomorphism
\begin{align} \label{QW56}
(\mfs \mfl_n)_{\DE_n}^j/(\mfs \mfl_n)_{\DE_n}^{j+1} \left(= \overline{\mcE nd}(\DF_\BB)^j/\overline{\mcE nd}(\DF_\BB)^{j+1}\right) \isom  (\Omega^{\otimes j})^{\oplus (n- |j|-\delta_{j0})}.
\end{align}
(Strictly speaking, the domain of this isomorphism for $j =0$ is canonically isomorphic to $\mcO_U^{\oplus n}/\mr{Im}(\Delta_0)$, where $\Delta_0$ denotes the diagonal morphism $\mcO_U \migiincl \mcO_U^{\oplus n}$.)

\begin{exa}[The case of $\hslash =0$] \label{y0134}
We shall specialize to  the case of $\hslash =0$.
 Since $ \mcD^{< \infty}_{0}$ is canonically isomorphic to $\mbS_{\mcO_X}(\mcT)$ (cf. Example \ref{QQ300}),
 we have an isomorphism
\begin{equation} \label{fff302}
\DF_\BB \isom \bigoplus_{j=0}^{n-1} \mcT^{\otimes j} \otimes \BB
\end{equation} 
compatible with the filtrations $\{ \DF_\BB^j \}_{j}$, $\{ \bigoplus_{j' =0}^{n-1-j}\mcT^{\otimes j'} \otimes \BB \}_j$.
Hence, for any integer $j$ with $-n< j < n$, the identification  (\ref{QW59}) gives
 a canonical decomposition
\begin{equation}  \label{fff301}
 (\mfs \mfl_n)^j_{\DE_{n}} = \bigoplus_{l=j}^{n-1}  (\mfs \mfl_n)^l_{\DE_{n}} / (\mfs \mfl_n)^{l+1}_{\DE_{n}} \left(= \bigoplus_{l=j}^{n-1} (\Omega^{\otimes j})^{\oplus (n-|j| -\delta_{j0})} \right).  \end{equation}
\end{exa}

\section{Dormant $(\mr{GL}_n, \hslash)$-opers} \label{y0135}

We now turn to discuss  the case of positive  characteristic.
First,   we  recall  the $p$-curvature of a logarithmic connection defined on a vector bundle; it is an invariant measuring, at least in the case of $\hslash =1$,   the number of horizontal sections (cf. Proposition \ref{QW201}).
After that, we introduce and study  dormant $(\mr{GL}_n, \hslash)$-opers on a log curve.
In what follows,  suppose that $\mr{char}(k) =:p>0$.

\subsection{$p$-curvature of an $\hslash$-flat vector bundle} \label{QG1000}

Let $T^\mr{log}$ and $Y^\mr{log}$  be fine log schemes over $k$ and  $f^\mr{log} : Y^\mr{log} \migi T^\mr{log}$ a  $k$-morphism of log schemes.
Also, let $\mcV$ be a rank $n$ vector bundle on $Y$ (or more generally,  an $\mcO_Y$-module) and $\nabla : \mcV \migi \Omega \otimes \mcV$ a flat $\hslash$-$T^\mr{log}$-connection on $\mcV$.

\ble  \label{QW91}
For a   local  section $\partial$ of $\mcT_{Y^\mr{log}/T^\mr{log}}$, 
 the endomorphism 
\begin{align}
{^p}\psi^\nabla (\partial) := \nabla_\partial^p- \hslash^{p-1}\cdot \nabla_{\partial^{[p]}}
\end{align}
 of $\mcV$  is $\mcO_Y$-linear,
  where $\nabla^p_{\partial}$ denotes the $p$-th power iterate of $\nabla_{\partial}: \mcV \migi \mcV$.
Moreover, for  local  sections $\partial, \partial_1, \partial_2 \in \mcT_{Y^\mr{log}/T^\mr{log}}$ and a local section  $a \in \mcO_Y$, we have
\begin{align}
{^p}\psi^\nabla (a \cdot \partial) = a^p \cdot {^p}\psi^\nabla (\partial),
\hspace{5mm}
{^p}\psi^\nabla (\partial_1 + \partial_2) = {^p}\psi^\nabla (\partial_1) + {^p}\psi^\nabla (\partial_2).
\end{align}
\ele
\begin{proof}
First, let us consider the former assertion.
Let $a$ and $v$ be arbitrary local sections of $\mcO_Y$ and $\mcV$ respectively.
By the $\hslash$-twisted  Leibniz rule  $\nabla_\partial (a \cdot v) =\hslash \cdot  \partial (a) \cdot v + a \cdot \nabla_\partial (v)$, the equaltiy
\begin{align}
\nabla_\partial^m (a \cdot v) = \sum_{i=0}^m \binom{m}{i} \cdot (\hslash^i \cdot \partial^i (a)) \cdot \nabla_\partial^{m-i} (v) 
\end{align}
holds for eacn $m \in \mbZ_{\geq 0}$. By putting $m =p$, we have
\begin{align} \label{QW107}
\nabla_\partial^p (a \cdot v) =\left( a \cdot \nabla^p_\partial (v) + \hslash^p \cdot \partial^p (a) \cdot v = \right) a \cdot \nabla^p_\partial (v) + \hslash^p \cdot \partial^{[p]} (a) \cdot v.
\end{align} 
On the other hand,  the $\hslash$-twisted  Leibniz rule with  $\partial$ replaced by $\partial^{[p]}$ reads
\begin{align} \label{QW106}
\nabla_{\partial^{[p]}} (a \cdot v) = \hslash \cdot \partial^{[p]} (a) \cdot v + a \cdot \nabla_{\partial^{[p]}} (v). 
\end{align}
By (\ref{QW107}) and (\ref{QW106}), we have  
\begin{align}
& \ 
{^p}\psi^\nabla (\partial) (a \cdot v) \left(=(\nabla_\partial^p -\hslash^{p-1}\cdot \nabla_{\partial^{[p]}})(a \cdot v) \right)
\\
 = & \  
 (a \cdot \nabla^p_\partial (v) +\cancel{ \hslash^p \cdot \partial^{[p]} (a) \cdot v}) 
 -
 \hslash^{p-1} \cdot (\cancel{\hslash \cdot \partial^{[p]} (a) \cdot v} + a \cdot \nabla_{\partial^{[p]}} (v)) \notag \\
 = & \ a \cdot (\nabla^p_\partial (v) - \hslash^{p-1} \cdot  \nabla_{\partial^{[p]}}(v)) \left(=  a\cdot {^p}\psi^\nabla (\partial) (v)\right),\notag
\end{align}
thus completing  the proof of the former assertion.

The latter assertion follows from an argument similar to the proof of Lemma \ref{QW82}   (together with  (\ref{QW101})).
\end{proof}

It follows from  the above lemma that the assignment
$F^{-1}_{Y}(\partial) \mapsto \nabla_\partial^p - \hslash^{p-1} \cdot \nabla_{\partial^{[p]}}\left(={^p}\psi^\nabla (\partial)\right)$
 for any local section $\partial \in \mcT_{Y^\mr{log}/T^\mr{log}}$ defines an {\it $\mcO_Y$-linear} morphism
\begin{align} \label{QQ203}
{^p}\psi^\nabla : F_Y^*(\mcT_{Y^\mr{log}/T^\mr{log}}) \migi \mcE nd_{\mcO_Y}(\mcV).
\end{align}

\bde \label{QQ202}
We shall refer to ${^p}\psi^\nabla$ as 
the {\bf $p$-curvature} of $\nabla$.
 \index{$\nabla$, $\hslash$-$T^\mr{log}$-connection@--- ${^p}\psi^\nabla$, $p$-curvature of}
Also,  we shall say that $\nabla$ is  {\bf $p$-flat} (or {\bf $p$-integrable}) 
\index{$\nabla$, $\hslash$-$T^\mr{log}$-connection@--- $p$-flat (= $p$-integrable)}
if    ${^p}\psi^\nabla =0$.
   \ede

\begin{rema} \label{QQ200}
Let $\mcE$ denote the $\mr{GL}_n$-bundle corresponding to $\mcV$ and  $\nabla_\mcE$ the $\hslash$-$T^\mr{log}$-connection corresponding to $\nabla$ resulting from Proposition \ref{QQ97}.
Then, by  Remark \ref{QW111}, the equality 
\begin{align} \label{QH10}
{^p}\psi^{\nabla_\mcE} = {^p}\psi^\nabla
\end{align}
holds via  the canonical isomorphism  $(\mfg \mfl_n)_\mcE \isom \mcE nd_{\mcO_Y}(\mcV)$ appearing in (\ref{QW30}).
In particular, the $p$-flatness of  $\nabla$  is equivalent to the $p$-flatness of $\nabla_\mcE$ in the sense of Definition \ref{y068}.
\end{rema}

\subsection{Determinant and tensor product} \label{QW203}
We shall recall  several facts on $p$-curvature.
We refer the reader to ~\cite{JP}, \S\,2, Proposition 2.1.2, for the case where  $\hslash =1$ and $Y^\mr{log}/T^\mr{log} = Y/T$.
Let $\nabla$ be a flat $\hslash$-$T^\mr{log}$-connection on $\mcV$.
Then,  the $p$-curvature of the determinant $\mr{det}(\nabla)$ (resp., the dual $\nabla^\vee$) of $\nabla$ satisfies 
\begin{align} \label{QW210}
{^p}\psi^{\mr{det}(\nabla)} = \mr{tr} ({^p}\psi^\nabla) \ \left(\text{resp.,} \ {^p}\psi^{\nabla^\vee} = - {^t}\psi^\nabla \right)
\end{align}
under the natural identification $\mr{tr}(\mcE nd_{\mcO_Y}(\mcV)) = \mcO_Y$ (resp., $\mcE nd_{\mcO_Y}(\mcV^\vee) = \mcE nd_{\mcO_Y}(\mcV)$).
In particular, the $p$-flatness of $\nabla$ implies  that of $\mr{det}(\nabla)$ (resp., $\nabla^\vee$).

Also, if $\mcV'$ is another vector bundle on $U$ and $\nabla'$ a flat $\hslash$-$T^\mr{log}$-connection on $\mcV'$, then the $p$-curvature of the tensor product $\nabla \otimes \nabla'$, which is flat (cf. \S\,\ref{QS7002}),   satisfies the  equality
\begin{align} \label{QW207}
{^p}\psi^{\nabla \otimes \nabla'} = {^p}\psi^\nabla \otimes \mr{id}_{\mcV'} + \mr{id}_{\mcV} \otimes {^p}\psi^{\nabla'}.
\end{align}
In particular, $\nabla \otimes \nabla'$ is $p$-flat if and only if both $\nabla$ and $\nabla'$ are $p$-flat.

\subsection{Cartier descent property} \label{QQ200}
We shall specialize to  the case of $\hslash=1$.
Given 
an 
$\mcO_{Y^{(1)}}$-module 
 $\mcW$,
then there exists a canonical $T^\mr{log}$-connection  
\index{$\nabla$, $\hslash$-$T^\mr{log}$-connection@--- $\nabla^\mr{can}_{\mcG}$, $\nabla^{\mr{can}}_{\mcG, \hslash}$,  $\nabla_\mcW^\mr{can}$, canonical}
\begin{align} \label{QW114}
\nabla^\mr{can}_\mcW : F^*_{Y/T}(\mcW) \migi \Omega_{Y^\mr{log}/T^\mr{log}} \otimes F^*_{Y/T}(\mcW)
\end{align}
 on the pull-back $F^*_{Y/T}(\mcW) \left(= \mcO_Y \otimes_{F^{-1}_{Y/T}(\mcO_{Y^{(1)}})}F^{-1}_{Y/T} (\mcW) \right)$;
it is determined by assigning $a \otimes v \mapsto da \otimes v$ for any local sections $a \in \mcO_{Y}$, $v \in \mcW$.
Since $\nabla_\mcW^\mr{can}$ coincides locally on $Y$ with the trivial connection,  $\nabla^\mr{can}_\mcW$ turns out to be flat and moreover  $p$-flat.
If
 $\mcW$ is a rank $n$ vector bundle and
  $\mcG$ denotes the $\mr{GL}_n$-bundle  on $Y^{(1)}$ associated to $\mcW$, then  
the flat bundle $(F^*_{Y/T}(\mcW), \nabla^\mr{can}_\mcW)$ 
 corresponds to 
 the flat $\mr{GL}_n$-bundle $(F_{Y/T}^*(\mcG), \nabla_\mcG^\mr{can})$
  (cf. (\ref{QR111}) for the definition of $\nabla_\mcG^\mr{can}$) via the equivalence of categories (\ref{QS7040}).

Conversely,  
let 
$\mcV$ be an $\mcO_Y$-module and $\nabla$ a flat $\hslash$-$T^\mr{log}$-connection on $\mcV$.
Then, the kernel $\mr{Ker}(\nabla)$ of $\nabla$  is closed under the $F_{Y/T}^{-1}(\mcO_{Y^{(1)}})$-action, i.e.,  can be regarded as an $O_{Y^{(1)}}$-module via the underlying homeomorphism of $F_{Y/T}$.
In particular,  we obtain the $\mcO_Y$-linear morphism
\begin{align} \label{QW120}
\nu^\sharp_{(\mcV, \nabla)} : F_{Y/T}^*(\mr{Ker}(\nabla)) \migi \mcV
\index{$\nu^\sharp_{(\mcV, \nabla)}$}
\end{align}
corresponding the natural inclusion $\mr{Ker}(\nabla) \migiincl F_{Y/T*}(\mcV)$ via the adjunction relation ``$F_{Y/T}^*(-) \dashv F_{Y/T*}(-)$''.
One verifies that $\nu^\sharp_{(\mcV, \nabla)}$ is compatible with the the respective $T$-connections, i.e.,  $\nabla_{\mr{Ker}(\nabla)}^\mr{can}$ and  $\nabla$.

Now, suppose that  the log structures of $Y^\mr{log}$ and $T^\mr{log}$ are  both trivial (i.e., $Y^\mr{log} = Y$ and $T^\mr{log} = T$) and that $Y/T$ is smooth.
 It is well-known (cf. ~\cite{Kal}, \S\,5, p.\,190, Theorem 5.1) that 
a flat connection $\nabla$ on a quasi-coherent $\mcO_Y$-module $\mcV$ is $p$-flat
 if and only if $\nu^\sharp_{(\mcV, \nabla)}$ is an isomorphism.
In particular,  the assignments $\mcW \mapsto (F^*_{Y/T}(\mcW), \nabla^{\mr{can}}_{\mcW})$ and $(\mcV, \nabla) \mapsto \mr{Ker}(\nabla)$ determine an equivalence of categories
\begin{equation} \label{equicat}
\begin{pmatrix}
\text{the category of} \\
\text{vector bundles on $Y^{(1)}$}
\end{pmatrix}
\isom
\begin{pmatrix}
\text{the category of flat  bundles $(\mcV, \nabla)$} \\
\text{on $Y/T$ with ${^p}\psi^\nabla =0$} 
\end{pmatrix}.
\end{equation}
We refer the reader to ~\cite{Og}, Theorem 1.3.4, for a generalization of this result to log schemes with nontrivial log structures. 

The following assertion will be used in the proof of Proposition \ref{QW200}.

\bpr \label{QW201}
Suppose that $Y^\mr{log} = Y$,  $T^\mr{log} =T$, and $Y/T$ is smooth.
Let $(\mcV, \nabla)$ be a rank $n$  flat bundle  on $Y/T$.
Then, $\nabla$ is $p$-flat if and only if  the $\mcO_{Y^{(1)}}$-module  $\mr{Ker}(\nabla)$  forms a rank $n$ vector bundle on $Y^{(1)}$.
\epr
\begin{proof}
The ``only if'' part of the equivalence can be verified immediately.  
Indeed, the $p$-flatness of $\nabla$ implies that the $\mcO_Y$-module $F_{Y/T}^*(\mr{Ker}(\nabla))$, being isomorphic to $\mcV$ by (\ref{equicat}), is locally free of rank $n$.
Hence,  $\mr{Ker}(\nabla)$ turns out to be locally free of rank $n$ because of  the faithful flatness of $F_{Y/T}$.

Next,  we shall consider the ``if'' part of the equivalence.
Assume  that $\mr{Ker}(\nabla)$ is locally free and of   rank $n$.
To prove the equality ${^p}\psi^\nabla =0$,
one may suppose, after possibly changing the base $T$ to each of  its geometric  points,
that $T = \mr{Spec}(k)$ (and $k$ is algebraically closed).
Moreover, it suffices to consider the case where $Y$ is connected.
We first  prove the injectivity of   $\nu^\sharp_{(\mcV, \nabla)}$.
Since $\nu^\sharp_{(\mcV, \nabla)}$ is compatible with the $k$-connections, 
$\nabla^{\mr{can}}_{\mr{Ker}(\nabla)}$ restricts to a $T$-connection $\nabla_{\mr{Ker}}$ on  the kernel $\mr{Ker}(\nu^\sharp_{(\mcV, \nabla)}) \left(\subseteq F^*_{Y/k}(\mr{Ker}(\nabla)) \right)$.
The short exact sequence
\begin{align}
0 \longmigi \mr{Ker}(\nu^\sharp_{(\mcV, \nabla)}) \xrightarrow{\mr{inclusion}} F^*_{Y/k}(\mr{Ker}(\nabla)) \xrightarrow{\nu^\sharp_{(\mcV, \nabla)}} \mcV
\end{align}
induces, via restriction to 
 sheaves of horizontal sections, the following exact sequence:
 \begin{align}
 0 \longmigi \mr{Ker}(\nabla_{\mr{Ker}})  \longmigi \mr{Ker}(\nabla) \xrightarrow{\mr{id}} \mr{Ker}(\nabla).
 \end{align}
 Hence,  the equality  $\mr{Ker}(\nabla_{\mr{Ker}}) =0$ holds.
 On the other hand, since  the subsheaf $\mr{Ker}(\nu^\sharp_{(\mcV, \nabla)})$  of $F^*_{Y/k}(\mr{Ker}(\nabla))$ is closed under  $\nabla^\mr{can}_{\mr{Ker}(\nabla)}$, it follows from 
the definition of $p$-curvature that the $p$-flatness of $\nabla^{\mr{can}}_{\mr{Ker}(\nabla)}$ implies that of $\nabla_{\mr{Ker}}$.
In particular,  $\mr{Ker}(\nu^\sharp_{(\mcV, \nabla)})$ is  isomorphic to $F^*_{Y/k}(\mr{Ker}(\nabla_{\mr{Ker}}))$ (cf. (\ref{equicat})), which is zero.
We thus obtain  the desired  injectivity of $\nu^\sharp_{(\mcV, \nabla)}$.
By  this injectivity   and the assumption that $\mr{Ker}(\nabla)$ is locally free of rank $n$,  the morphism 
$\nu^\sharp_{(\mcV, \nabla)}$ turns out to be  an isomorphism over the generic point of $Y$.
That is to say, the $p$-curvature ${^p}\psi^{\nabla}$ of $\nabla$, considered as an element of $\Gamma (Y, F^*_{Y/k}(\Omega) \otimes \mcE nd_{\mcO_Y} (\mcV) )$, vanishes at the generic point.
This fact and the local freeness of  $F^*_{Y/k}(\Omega) \otimes \mcE nd_{\mcO_Y} (\mcV)$  together imply   ${^p}\psi^{\nabla} =0$.
This completes the proof of the ``if '' part, and hence completes the proof of the assertion.
\end{proof}

\subsection{The definition of a dormant $(\mr{GL}_n, \hslash, \bb)$-oper} \label{QW200}

Hereafter, 
we suppose   that  $\mr{char}(k)=:p >n$.
Let $S^\mr{log}$ be an fs log scheme over $k$ and $U^\mr{log}$ a log curve over $S^\mr{log}$.
Also, let $\bb := (\BB, \nabla_\bb)$ be an $(n, \hslash)$-theta characteristic of $U^\mr{log}/S^\mr{log}$.

\bde \label{y0136}
  We shall say that a $(\mr{GL}_n, \hslash)$-oper $\msF^\heartsuit := (\mcF, \nabla, \{ \mcF^j \}_j)$ (resp.,  a $(\mr{GL}_n, \hslash, \bb)$-oper $\nabla^\diamondsuit$) $\nabla^\diamondsuit$ on $U^\mr{log}/S^\mr{log}$ is {\bf dormant}
   \index{$\mcF^\heartsuit$, $(\mr{GL}_n, \hslash)$-oper@--- dormant} \index{$\nabla^\diamondsuit$, $(\mr{GL}_n, \hslash, \bb)$-oper@--- dormant} if  
 $\nabla$ (resp., $\nabla^\diamondsuit$) is $p$-flat.
  \ede

\begin{rema} \label{QQ164}
Let $\nabla^\diamondsuit$ be a $(\mr{GL}_n, \hslash, \bb)$-oper on $U^\mr{log}/S^\mr{log}$.
Then, we have  ${^p}\psi^{\nabla_\bb} = {^p}\psi^{\mr{det} (\nabla^\diamondsuit)}$ $\left( \stackrel{(\ref{QW210})}{=} \mr{tr}({^p}\psi^{\nabla^\diamondsuit})\right) =0$ under the identification $\mr{det}(\DF_\BB) = \mcT^{\otimes \frac{n(n-1)}{2}} \otimes \BB^{\otimes n}$ by (\ref{fff0100}).
It follows that  the dormancy condition on   $\nabla^\diamondsuit$ implies the $p$-flatness of $\nabla_\bb$.
\end{rema}

\bpr \label{y0137}
Suppose that $\nabla_\bb$ is $p$-flat.
Also, let $\nabla^\diamondsuit$ be a $(\mr{GL}_n, \hslash, \bb)$-oper on $U^\mr{log}/S^\mr{log}$ and denote by $\msE^\spadesuit := (\mcE_B, \nabla)$ the corresponding  $(\mfs \mfl_n, \hslash)$-oper  via the sequence of isomorphisms in (\ref{GL42}).
Then, $\nabla^\diamondsuit$ is dormant if and only if $\msE^\spadesuit$ is dormant.
\epr
\begin{proof}
Denote by $\DE_{\mr{GL}_n}$ the $\mr{GL}_n$-bundle corresponding to $\DF_\BB$.
Recall from the discussion in the proof of Proposition \ref{y0111} that the morphism
\begin{align} \label{QW159}
\widetilde{\mcT}_{\DE_{\mr{GL}_n}^\mr{log}/S^\mr{log}} \isom  \widetilde{\mcT}_{\DE_n^\mr{log}/S^\mr{log}} \times_{\mcT_{U^\mr{log}/S^\mr{log}}} \widetilde{\mcT}_{\mr{det}(\DF_\BB)^{\times \mr{log}}/S^\mr{log}}
\end{align}
(cf.  the resp'd portion of (\ref{GL3}))  induced by differentiating the quotient $\mr{GL}_n \migisurj \mr{PGL}_n$ and the determinant map $\mr{GL}_n \migisurj \mbG_m$ is an isomorphism.
Moreover,
 (\ref{QW159}) preserves  the $p$-power operations when we equip 
  the codomain  with the $p$-power operation arising from both $\widetilde{\mcT}_{\DE_n^\mr{log}/S^\mr{log}}$ and $\widetilde{\mcT}_{\mr{det}(\DF_\BB)^{\times \mr{log}}/S^\mr{log}}$.
 It follows that the equality ${^p}\psi^{\nabla^\diamondsuit}= ({^p}\psi^{\nabla}, {^p}\psi^{\nabla_\bb})$ ($= ({^p}\psi^{\nabla}, 0)$ by assumption) holds  via  (\ref{fff0100}), (\ref{QH10}), and (\ref{QW159}).
 The required equivalence follows immediately from this equality.
\end{proof}

\subsection{Differentia operators having a full set of root functions}\label{QW160}
In the rest of this section, 
let us consider the relationship between dormant opers and projective connections  having a full set of root functions.

Suppose that the log structures of $U^\mr{log}$ and $S^\mr{log}$ are trivial, i.e., $U^\mr{log} = U$ and $S^\mr{log} = S$.
Let $\mcL_i$ ($i=1,2$) be line bundles on $U$ and
$D$ be an $\mcO_U$-linear morphism $\mcL_1 \migi \mcL_2 \otimes \mcD_1^{< n+1}$.
Denote by $\Dif : \mcL_1 \migi \mcL_2$ the differential operator,  in the sense of Remark \ref{QQ425},  corresponding to $D$ via (\ref{Edr56}).
Then, $\Dif$ may be 
 considered as an {\it $\mcO_{U^{(1)}}$-linear} morphism 
  via the underlying homeomorphism of $F_{U/S}$.
In particular, the kernel $\mr{Ker} (\Dif)$
of this morphism  forms an $\mcO_{U^{(1)}}$-submodule of $F_{U/S*}(\mcL_1)$.

\bde \label{D01fgh}
We shall say that {\bf $D$ has a full set of  root functions} \index{(linear) differential operator@--- has a full set of  root functions} if  $\mr{Ker} (\Dif)$  is a vector bundle on $U^{(1)}$ of  rank $n$, i.e., has locally a basis consisting of $n$ sections as  an $\mcO_{U^{(1)}}$-module. 
\ede

Denote by 
\index{$\DD_{n, 1, \bb}^{\,\clubsuit \mr{full}}$, $\mfD \mfi \mff \mff^{\clubsuit \mr{full}}_{n, 1, \bb  (, \rho)}$}
\begin{align} \label{QQ392}
\DD_{n, 1, \bb}^{\,\clubsuit \mr{full}}
\end{align}
the subsheaf of $\DD_{n, 1, \bb}^{\,\clubsuit}$ consisting of
 $(n, \bb)$-projective connections
 having a full set of root functions.
 On the other hand, denote by
 \index{$\mfO \mfp^\ZZZ_{\mfg, \hslash, (\rho,) \msX}$, $\mfO \mfp^\ZZZ_{\mfg, \hslash, (\rho,) g,r}$, $\mfO \mfp_{n, \hslash, \bb (, \rho)}^{\diamondsuit \ZZZ}$}
 \begin{align}
 \mcO p_{n, 1, \bb}^{\diamondsuit^{_\mr{Zzz...}}}
 \end{align}
 the subsheaf of $\mcO p_{n, 1, \bb}^{\diamondsuit}$ consisting of dormant $(\mr{GL}_n,  \bb)$-opers.

 \bpr \label{QW200}
 The isomorphism $\Lambda_\bb^{\clubsuit \Rightarrow \diamondsuit}$ (cf. (\ref{GL28})) restricts to an isomorphism of sheaves
 \begin{align} \label{QQ379}
 \DD_{n, 1, \bb}^{\,\clubsuit \mr{full}} \isom \mcO p_{n, 1, \bb}^{\diamondsuit ^{_\mr{Zzz...}}}.
 \end{align}
 \epr
 \begin{proof}
 Since the problem is of local nature,  one may assume, after possibly replacing $U$ with its covering in the \'{e}tale topology, that
there are a global log chart $(U, \partial)$ and an identification $\BB  = \mcO_U$.
 Let us take  a local section
  $D^\clubsuit$ of $\DD_{n, 1, \bb}^{\,\clubsuit}$, which 
 may be expressed as $D^\clubsuit = \partial^n + q_1 \cdot  \partial^{n-1}+ \cdots + q_{n-1} \cdot \partial + q_n$ for 
   local sections $q_1, \cdots, q_n \in \mcO_U$.
Then,  the dual  $(D^{\clubsuit \Rightarrow \diamondsuit})^\vee$ of $D^{\clubsuit \Rightarrow \diamondsuit}$
satisfies   (with respect to a suitable local basis of $\DF_\BB$) the following  equality:
\begin{align}
(D^{\clubsuit \Rightarrow \diamondsuit})^{\vee \mr{GL}_n} (\partial) =  \begin{pmatrix} q_1 & q_2& q_3& \cdots & q_{n-1}&  q_n \\ -1& 0 & 0& \cdots & 0 & 0 \\ 0 & -1 & 0& \cdots & 0 & 0 \\ 0 & 0 &-1& \cdots & 0& 0 \\ \vdots & \vdots & \vdots & \ddots& \vdots & \vdots  \\ 0 & 0 & 0 & \cdots & -1 & 0  \end{pmatrix}
\end{align}
(cf. (\ref{QS7144}) for the definition of $(-)^{\mr{GL}_n}$).
Then, $y \mapsto {^t(}\partial^{n-1}(y), \cdots, \partial (y), y)$
gives a bijective correspondence between 
the local sections of $\mr{Ker} (\Dif^\clubsuit)$ 
 and the horizontal local sections of $\DF_\BB^\vee$ with respect to the connection $(D^{\clubsuit \Rightarrow \diamondsuit})^\vee$.
This implies that $D^\clubsuit$ has a full set of root functions if and only if the sheaf of horizontal sections in $\DF_\BB^\vee$
  forms  a rank $n$ vector bundle on $U^{(1)}$.
But, it follows from Proposition \ref{QW201} that the latter  condition is   equivalent to  the $p$-flatness of $(D^{\clubsuit \Rightarrow \diamondsuit})^\vee$, and hence, equivalent to that of $D^{\clubsuit \Rightarrow \diamondsuit}$.
This completes the proof of the assertion.
\end{proof}


\section{Comparison of the moduli functors} \label{QQ189}

In this section, we suppose that {\it  either ``\,$\mr{char}(k) =0$'' or ``\,$\mr{char}(k)>n$''  is fulfilled}. Let $(g,r)$ be a pair of nonnegative integers with $2g-2 +r >0$.

\subsection{} \label{QG1001}

In what follows, 
let us consider the case where  the log curve ``$U^\mr{log}/S^\mr{log}$'' is taken to be $X^\mr{log}/S^\mr{log}$ associated to 
 an $r$-pointed stable curve  $\msX := (f:X \migi S, \{ \sigma_i \}_{i=1}^r)$   of genus $g$ over  a $k$-scheme $S$. 
Let us fix $\hslash \in \Gamma (S, \mcO_S)$ and  $\rho \in \mfc_{\mfs\mfl_n}^{\times r} (S)$, where $\mfc_{\mfs\mfl_n}^{\times r} (S) := \{ \emptyset \}$ and $\rho := \emptyset$ if $r =0$.
According to  Example \ref{QQ1290}, 
the $k$-scheme $\mfc_{\mfs \mfl_n}$
  and  the adjoint quotient $\chi : \mfs \mfl_n \migisurj \mfc_{\mfs \mfl_n}$ are  well-defined.
In exactly the same way as Definition \ref{y057},  we can define the notion of 
an $(\mfs \mfl_n, \hslash)$-oper of radii $\rho$. 
(Recall that Definition \ref{y057} was described  under the assumption that ``\,$\mr{char}(k)=0$'' or ``\,$\mr{char}(k) > 2h_\mbG$'' is fulfilled.)
Also, just as in (\ref{opermoduli}) (resp., (\ref{QQ1025})), we can construct 
the moduli functor on $\mfS \mfc \mfh_{/S}$ 
classifying $(\mfs \mfl_n, \hslash)$-opers (resp., $(\mfs \mfl_n, \hslash)$-opers of radii $\rho$) on $\msX$; we denote this functor by  
\index{$\mfO \mfp_{n, \hslash (, \rho)}^{\spadesuit}$, $\mfO \mfp_{n, \hslash (, \rho)}^{\spadesuit, \mr{O}}$, $\mfO \mfp_{n, \hslash (, \rho)}^{\spadesuit, \mr{Sp}}$}
\begin{align}
\mfO \mfp_{\mfs \mfl_n, \hslash, \msX} \ \text{or simply} \ \mfO \mfp^\spadesuit_{n, \hslash}   \ \left(\text{resp.,} \ \mfO \mfp_{\mfs \mfl_n, \hslash, \rho, \msX} \ \text{or simply} \  \mfO \mfp^\spadesuit_{n, \hslash, \rho}  
  \right).
\end{align}
 
 Next, suppose that  $X^\mr{log}/S^\mr{log}$ admits  an $(n, \hslash)$-theta characteristic $\bb := (\BB, \nabla_\bb)$ of $\msX$.
We shall write
\index{$\mfO \mfp_{n, \hslash, \bb (, \rho)}^\diamondsuit$, $\mfO \mfp^{\diamondsuit, \mr{O}}_{n, \hslash}$, $\mfO \mfp^{\diamondsuit, \mr{Sp}}_{n, \hslash}$}
\index{$\overline{\mfO} \mfp_{n, \hslash (, \rho)}^\heartsuit$}
\index{$\DD^{\,\clubsuit}_{n, \hslash, \BB}$, $\DD^{\,\clubsuit}_{n, \hslash, \bb (, \rho)}$, $\mfD \mfi \mff \mff^{\,\clubsuit}_{n, \hslash, \bb (, \rho)}$}
\begin{align} \label{1093}
\overline{\mfO} \mfp_{n, \hslash}^\heartsuit \left(\text{resp.,}  \ \mfO \mfp_{n, \hslash, \bb}^\diamondsuit ;\text{resp.,} \ \mfD \mfi \mff \mff^\clubsuit_{n,\hslash, \bb} \right)
\end{align}
the $\mfS \mfe \mft$-valued contravariant functor on $\mfS \mfc \mfh_{/S}$ which, to any $S$-scheme $s : S' \migi S$, assigns the set
of global sections of $\overline{\mcO} p^\heartsuit_{n, \hslash}$  (resp., $\mcO p^\diamondsuit_{n, \hslash, \bb}$; resp., $\DD^{\,\clubsuit}_{n ,\hslash, \bb}$) defined for the log curve
$(S' \times_S X^\mr{log})/ (S' \times_S S^\mr{log})$.
Also,  we write
\index{$\overline{\mfO} \mfp_{n, \hslash (, \rho)}^\heartsuit$}
\index{$\mfO \mfp_{n, \hslash, \bb (, \rho)}^\diamondsuit$, $\mfO \mfp^{\diamondsuit, \mr{O}}_{n, \hslash}$, $\mfO \mfp^{\diamondsuit, \mr{Sp}}_{n, \hslash}$}
\index{$\DD^{\,\clubsuit}_{n, \hslash, \BB}$, $\DD^{\,\clubsuit}_{n, \hslash, \bb (, \rho)}$, $\mfD \mfi \mff \mff^{\,\clubsuit}_{n, \hslash, \bb (, \rho)}$}
\begin{align} \label{QW506}
\overline{\mfO} \mfp_{n, \hslash, \rho}^\heartsuit \left(\text{resp.,}  \ \mfO \mfp_{n, \hslash, \bb, \rho}^\diamondsuit ;\text{resp.,} \ \mfD \mfi \mff \mff^\clubsuit_{n,\hslash, \bb, \rho} \right)
\end{align}
for the subfunctor of $\overline{\mfO} \mfp_{n, \hslash}^\heartsuit$ (resp.,   $\mfO \mfp_{n, \hslash, \bb}^\diamondsuit$; resp., $\mfD \mfi \mff \mff^\clubsuit_{n,\hslash, \bb}$) classifying those of radii $\rho$.

\bt[cf. Theorem \ref{y017}] \label{y0129}
Let us keep the above notation.
Then, there exists a canonical sequence of isomorphisms
  \begin{align} \label{QQ360}
 \vcenter{\xymatrix@C=46pt@R=36pt{
\mfD \mfi \mff \mff^{\, \clubsuit}_{n,\hslash, \bb, \Box} \ar[rd]^-{\sim}_-{\Lambda_{\msX, \bb, \Box}^{\clubsuit \Rightarrow \diamondsuit}}&& \overline{\mfO}\mfp^\heartsuit_{n,\hslash, \Box}\ar[rd]^-{\sim}_-{\Lambda_{\msX, \Box}^{\heartsuit \Rightarrow \spadesuit}}&\\
& \mfO \mfp^\diamondsuit_{n, \hslash, \bb, \Box}\ar[ur]_-{\sim}^-{\Lambda_{\msX, \bb, \Box}^{\diamondsuit \Rightarrow \heartsuit}}&& \mfO \mfp^\spadesuit_{n,\hslash, \Box}.
}}
 \end{align}
of functors, where $\Box$ denotes either the absence or presence of ``$\rho$''.
Moreover,  in the case where  $\Box$ denotes the absence (resp., presence) of $\rho$, then
all the functors appearing in this sequence 
  may be represented by 
 relative affine spaces over $S$ modeled on 
$\mbV (f_* (\DV_\BB))$ (resp., $\mbV (f_* (\DV_\BB (-D_\msX)))$),
 in particular, of relative dimension 
\begin{align}
\aleph (\mfs \mfl_n) = \frac{2g-2+r}{2} \cdot (n^2-1) + \frac{r}{2} \cdot (n-1)  \hspace{4mm}   \\
\left(\text{resp.,} \ {^c}\aleph (\mfs \mfl_n) = \frac{2g-2+r}{2} \cdot (n^2-1) - \frac{r}{2} \cdot (n-1)\right). \notag
\end{align}
\et
\begin{proof}
We only consider the non-resp'd assertion because the proof of  the resp'd assertion is entirely similar.
The former assertion follows from Theorem \ref{y0128} with  $U^\mr{log}/S^\mr{log}$ replaced by $(S' \times_S X^\mr{log})/(S' \times_S S^{\mr{log}})$ for various $S$-schemes $S'$.
To prove the latter assertion, 
it suffices to consider  the representability of 
$\mfD \mfi \mff \mff^{\, \clubsuit}_{n,\hslash, \bb}$ because of the former assertion.
Recall from the discussion in  \S\,\ref{y0124} that
$\DD^{\, \clubsuit}_{n ,\hslash, \bb}$ has a  $\DV_\BB$-torsor structure.
Since   $\mbR^1f_*(\DV_\BB) = 0$ (cf. Lemma \ref{QQ230}, (i),  below),
the direst image  $f_*(\DD^{\, \clubsuit}_{n ,\hslash, \bb})$ has  a section locally on $S$  and the formation of this sheaf commutes  with base-change to $S$-schemes.
It follows that  
the $f_*(\DV_\BB)$-torsor structure on $f_*(\DD^{\, \clubsuit}_{n ,\hslash, \bb})$ induces a structure of relative affine space on
 $\mfD \mfi \mff \mff^{\clubsuit}_{n, \hslash, \bb}$  modeled on  $\mbV (f_*(\DV_\BB))$.
Finally, it follows from Lemma \ref{QQ230}, (ii),  that the $S$-scheme  $\mbV (f_*(\DV_\BB))$ is of relative dimension $\aleph (\mfs \mfl_n)$.
This completes the proof of  the assertion.
\end{proof}

The following lemma was used in the proof of the above proposition.

\ble \label{QQ230}
\begin{itemize}
\item[(i)]
The equality $\mbR^1 f_*(\DV_\BB) =0$ (resp., $\mbR^1 f_*(\DV_\BB (-D_\msX)) =0$) holds.
\item[(ii)]
  The $\mcO_S$-module  $f_*(\DV_\BB)$ (resp.,   $f_*(\DV_{\BB}(-D_\msX))$) forms a vector bundle  of rank $\aleph (\mfs \mfl_n) $ (resp., ${^c}\aleph (\mfs \mfl_n)$).
 Moreover, if $s: S' \migi S$ is a morphism of $k$-schemes, then the natural morphism of $S'$-schemes
 \begin{align} \label{QW247}
 S' \times_S \mbV (f_*(\DV_{\BB}))  \migi \mbV ((\mr{id}_{S'}\times f)_* (\DV_{s^*(\BB)})) \hspace{16mm}\\
 \left(\text{resp.,} \  S' \times_S \mbV (f_*(\DV_{\BB} (-D_\msX)))  \migi \mbV ((\mr{id}_{S'}\times f)_* (\DV_{s^*(\BB)}(-D_{S' \times_S \msX})))  \right) \notag 
   \end{align}
 is an isomorphism.
\end{itemize}
\ele
\begin{proof}
We only consider the resp'd assertion because the proof of the non-resp'd assertion is entirely similar.
For each  $j = 1, \cdots, n-1$, 
the relative degree of $\mcT_{}^{\otimes (n-j)}$ is negative.
Hence,
we have
\begin{align}   \label{QQ380}
 \mbR^1f_*({^\Join}\DD^{\,<j}_{\BB} (-D_\msX) /{^\Join}\DD^{\,<j-1}_{\BB}(-D_\msX)) &\cong \mbR^1 f_* (\Omega^{\otimes (n-j+1)} (-D_\msX)) \\
   &\cong  f_*(\mcT_{}^{\otimes (n-j)})^\vee  \notag \\
   &   =0,  \notag
 \end{align}
 where the first ``$\cong$'' is induced by  (\ref{GL8}) and the second 
 isomorphism follows from Grothendieck-Serre duality.
 It follows that  
 \begin{align} \label{QW510}
 \mbR^1f_*(\DV_\BB (-D_\msX)) \left(= \mbR^1f_*({^\Join}\DD^{\,<n-1}_{\BB}(-D_\msX)))\right) =0,
 \end{align}
  which completes the proof of assertion (i).
 Moreover, this fact implies 
 that $f_*(\DV_\BB (-D_\msX))$ is a vector bundle on $S$ and (\ref{QW247}) is an isomorphism (cf. ~\cite{Har}, Chap.\,III, Theorem 12.11 (b)).
  The rank $\mr{rk}(f_*(\DV_\BB (-D_\msX)))$
 of $f_*(\DV_\BB (-D_\msX))$
  may be calculated  by
\begin{align} 
 & \ \ \ \ \mr{rk}(f_*(\DV_\BB (-D_\msX))) \\
& 
 = \sum_{j=1}^{n-1} \mr{rk}(f_*({^\Join}\DD^{\,<j}_{n, \hslash, \BB} (-D_\msX)/{^\Join}\DD^{\,<j-1}_{n, \hslash, \BB}(-D_\msX)))  \notag \\
 & = \sum_{j=1}^{n-1} \mr{rk}(f_*(\Omega_{}^{\otimes (n-j+1)}(-D_\msX))) \notag \\
& =  \sum_{j=1}^{n-1} ((2g-2+r)(n+1-j)-(g-1+r)) \notag \\ 
& =  (2g-2+r)(n^2-1) -  (2g-2+r) \cdot \frac{n(n-1)}{2} - (g-1+r)  (n-1) \notag \\
& = {^c}\aleph (\mfs \mfl_n), \notag
\end{align}
where the first equality follows from (\ref{QQ380}) and  the third equality follows from the Riemann-Roch theorem.
This completes the proof of  the assertion.
\end{proof}

\bpr \label{y0130} 
Let us keep the above notation.
Moreover, let $\bb'$ be another $(n, \hslash)$-theta characteristic of $X^\mr{log}/S^\mr{log}$.
If $\Box$ denotes either the absence or presence of ``$\rho$'', then
there exist  isomorphisms of $S$-schemes
\begin{equation} 
\label{1094}
\clubsuit_{\msX, \bb \Rightarrow \bb', \Box} :  \mfD \mfi \mff \mff_{n, \hslash, \bb, \Box}^{\clubsuit}   \isom  
  \mfD \mfi \mff \mff_{n, \hslash, \bb', \Box}^{\clubsuit}, 
  \hspace{2mm}
\diamondsuit_{\msX,  \bb \Rightarrow \bb', \Box} :  \mfO \mfp_{n, \hslash, \bb, \Box}^\diamondsuit \isom   \mfO \mfp_{n, \hslash, \bb', \Box}^\diamondsuit  
\end{equation}
which fit into the following commutative diagram:
\begin{align} \label{QQ270}
\vcenter{\xymatrix@C=50pt@R=10pt{
   \mfD \mfi \mff \mff_{n, \hslash, \bb, \Box}^{\clubsuit}
    \ar[r]_{\sim}^{\Lambda^{\clubsuit \Rightarrow \diamondsuit}_{\msX, \Box}} \ar[dd]^-{\wr}_-{\clubsuit_{\msX, \bb \Rightarrow \bb', \Box}}&
        \mfO \mfp_{n, \hslash, \bb, \Box}^\diamondsuit  \ar[rd]_-{\sim}^-{\Lambda_{\msX, \Box}^{\diamondsuit \Rightarrow \heartsuit}} \ar[dd]^-{\wr}_-{\diamondsuit_{\msX,  \bb \Rightarrow \bb', \Box}}&\\
 && \overline{\mfO} \mfp_{n, \hslash, \Box}^\heartsuit. \\
 \mfD \mfi \mff \mff_{n, \hslash, \bb', \Box}^{\clubsuit}
  \ar[r]^-{\sim}_-{\Lambda_{\msX, \Box}^{\clubsuit \Rightarrow \diamondsuit}}&  \mfO \mfp_{n, \hslash, \bb', \Box}^\diamondsuit  \ar[ur]^-{\sim}_-{\Lambda_{\msX, \Box}^{\diamondsuit \Rightarrow \heartsuit}}&
}}
\end{align}
\epr
\begin{proof}
The assertion follows from
Proposition \ref{y0118}.
  \end{proof}


In the case of $n=2$, the $\Omega_{X^\mr{log}/S^\mr{log}}^{\otimes 2}$-torsor structure on $\mcO p^\spadesuit_{n, \hslash}$ induces  a  free transitive  $\mbV (f_* (\Omega^{\otimes 2}_{X^\mr{log}/S^\mr{log}}))$-action on $\mfO \mfp^\spadesuit_{2, \hslash}$, which also induces  a  $\mbV (f_* (\Omega^{\otimes 2}_{X^\mr{log}/S^\mr{log}}(-D_\msX)))$-action on $\mfO \mfp^\spadesuit_{2, \hslash, \rho}$.
Thus, we have the following assertion.
 
\bco \label{QQ1195}
(In this corollary, we do not impose  the existence assumption of a globally defined  $(2, \hslash)$-theta characteristic.)
The functor 
$\mfO \mfp^\spadesuit_{2, \hslash}$
 (resp., 
 $\mfO \mfp^\spadesuit_{2, \hslash, \rho}$)
  can be represented by a relative affine space over $S$ modeled on $\mbV (f_* (\Omega^{\otimes 2}_{X^\mr{log}/S^\mr{log}}))$ (resp., $\mbV (f_* (\Omega^{\otimes 2}_{X^\mr{log}/S^\mr{log}}(-D_\msX)))$).
In particular,  the fiber of $\mfO \mfp_{2, \hslash}$ (resp., $\mfO \mfp_{2, \hslash, \Box}$) over any geometric point of $S$ is an affine space of dimension $3  g-3+2  r$ (resp., $3g-3+r$).
\eco

\subsection{Positive characteristic} \label{QG1004}
Let us consider the case where  {\it the characteristic $\mr{char}(k)$ of $k$ is a prime $p$ with $n<p$}.
Let $\bb := (\BB, \nabla_\bb)$ be an $(n, \hslash)$-theta characteristic of $X^\mr{log}/S^\mr{log}$ with ${^p}\psi^{\nabla_\bb}=0$.
(Recall from  the discussion in \S\,\ref{QQ183} that there always exists such an $(n, \hslash)$-theta characteristic.)
We shall denote by 
\index{$\mfO \mfp^\ZZZ_{\mfg, \hslash, (\rho,) \msX}$, $\mfO \mfp^\ZZZ_{\mfg, \hslash, (\rho,) g,r}$, $\mfO \mfp_{n, \hslash, \bb (, \rho)}^{\diamondsuit \ZZZ}$}
\begin{equation} \label{10034}
\mfO \mfp^\ZZZ_{\mfs \mfl_n, \hslash, \msX} \ \left(\text{resp.,} \
\mfO \mfp^\ZZZ_{\mfs \mfl_n, \hslash, \rho, \msX}; \text{resp.,} \ 
\mfO \mfp_{n, \hslash, \bb}^{\diamondsuit \ZZZ}; \text{resp.,} \  \mfO \mfp_{n, \hslash, \bb, \rho}^{\diamondsuit \ZZZ} \right)
\end{equation}
 the subfunctor of 
 $\mfO \mfp_{\mfs \mfl_n, \hslash, \msX}$ (resp., $\mfO \mfp_{\mfs \mfl_n, \hslash, \rho, \msX}$; resp.,
 $\mfO \mfp_{n, \hslash, \bb}^{\diamondsuit}$; resp., $\mfO \mfp_{n, \hslash, \bb, \rho}^{\diamondsuit}$) 
  classifying {\it dormant} $(\mr{GL}_n, \hslash, \bb)$-opers (cf. (\ref{QS239}), (\ref{QS240})).

\bco \label{1090}
Let  $\Box$ denote either the absence or presence of ``$\rho$''.
Then, the  following assertions  hold:
\begin{itemize}
\item[(i)]
The composite isomorphism $\Lambda_{\msX, \Box}^{\heartsuit \Rightarrow \spadesuit}\circ \Lambda_{\msX, \bb,\Box}^{\diamondsuit \Rightarrow \heartsuit}$ (cf. (\ref{QQ360})) restricts to an isomorphism of $S$-schemes
\begin{equation} \label{1091}
{^p}\Lambda_{\msX, \bb, \Box}^{\diamondsuit \Rightarrow \spadesuit} :  \mfO \mfp_{n, \hslash, \bb, \Box}^{\diamondsuit \ZZZ} \isom \mfO \mfp_{\mfs \mfl_n, \hslash, \Box, \msX}^\ZZZ. 
 \end{equation}

\item[(ii)]
Let 
$\bb':= (\BB', \nabla_{\bb'})$ be another $(n, \hslash)$-theta characteristic of $X^\mr{log}/S^\mr{log}$ with ${^p}\psi^{\nabla_{\bb'}} =0$.
Then, 
 the isomorphism $\diamondsuit_{\msX, \bb \Rightarrow \bb', \Box}$ (cf. (\ref{1094})) 
  restricts to  an isomorphism
\begin{equation}
\label{1097}
 {^p}\diamondsuit_{\msX,  \bb \Rightarrow \bb', \Box} : \mfO \mfp_{n, \hslash, \bb, \Box}^{\diamondsuit \ZZZ}  \isom \mfO \mfp_{n, \hslash, \bb', \Box}^{\diamondsuit \ZZZ} 
\end{equation}
of  $S$-schemes  which makes
the following square diagram commute:
\begin{align} \label{1098}
\vcenter{\xymatrix@C=46pt@R=12pt{
\mfO \mfp_{n, \hslash, \bb, \rho}^{\diamondsuit \ZZZ} \ar[rd]_-{\sim}^-{{^p}\Lambda^{\diamondsuit \Rightarrow \spadesuit}_{\msX, \bb, \Box}} \ar[dd]^{\wr}_-{{^p}\diamondsuit_{\msX,  \bb \Rightarrow \bb', \Box}}&\\
& \mfO \mfp_{\mfs \mfl_n, \hslash, \Box, \msX}^{^\mr{Zzz...}}.\\
\mfO \mfp_{n, \hslash, \bb', \Box}^{\diamondsuit \ZZZ}  \ar[ur]^-{\sim}_-{\ {^p}\Lambda^{\diamondsuit \Rightarrow \spadesuit}_{\msX, \bb', \Box}} &
}}
\end{align}
\end{itemize}
\eco
\begin{proof}
The assertion follows from  Propositions \ref{y0137} and  \ref{y0130}.
\end{proof}

\bpr \label{QG1780}
Let 
$\Box$ denote either the presence or absence of ``$\rho$''.
Let  $\mfD \mfi \mff \mff^{\clubsuit \mr{full}}_{n, 1, \bb, \Box}$
\index{$\DD_{n, 1, \bb}^{\,\clubsuit \mr{full}}$, $\mfD \mfi \mff \mff^{\clubsuit \mr{full}}_{n, 1, \bb  (, \rho)}$}
 denote the subfunctor of $\mfD \mfi \mff \mff^{\clubsuit \mr{full}}_{n, 1,\bb, \Box}$
classifying $(n, \bb)$-projective connections   having a full set of root functions.
Then, $\mfD \mfi \mff \mff^{\clubsuit \mr{full}}_{n, 1, \bb, \Box}$ may be represented by a closed subscheme of $\mfD \mfi \mff \mff^{\clubsuit}_{n, 1, \bb, \Box}$  and there exists a canonical  $S$-isomorphism
\begin{align}
\mfD \mfi \mff \mff^{\clubsuit \mr{full}}_{n, 1, \bb, \Box} \isom \mfO \mfp^\ZZZ_{\mfs \mfl_n, 1, \Box,  \msX}.
\end{align}
In particular,  $\mfD \mfi \mff \mff^{\clubsuit \mr{full}}_{n, 1, \bb, \rho}$ is empty  unless $\rho \in \mfc^{\times r} (\mbF_p)$.
\epr
\begin{proof}
The assertion follows from Propositions \ref{QS241}, (i), and  \ref{QW200} and  Theorem \ref{y0129}.
\end{proof}

\section{Hypergeometric differential  operators} \label{QG1901}

In this final section of the present chapter, 
we study $\mfs \mfl_2$-opers on the $3$-pointed projective line $\msP := (\mbP^1, \{[0], [1], [\infty] \})$ (cf. (\ref{1051})) over $k$.
Suppose that the base field $k$ is an algebraically closed field in which $2$ is invertible.
Let us consider 
 a monic, linear, and   second-order differential operators on    $\mbP^1$ having at most three regular singular points.
 The classical theory of Riemann schemes shows that 
any such operator may be transformed 
 into
a Gauss' hypergeometric differential operator
\begin{align}
\Dif_{a,b, c}^\clubsuit := \partial_x^2 + \left(\frac{c}{x} + \frac{1 -c + a + b}{x-1} \right) \cdot \partial_x +\frac{ab}{x (x-1)} 
\end{align}
determined by some  triple  $(a, b,c) \in k^{\times 3}$, where
$\mr{Spec}(k[x]) = \mbA^1 = \mbP^1 \setminus \{ \infty \}$
 and $\partial_x := \frac{d}{dx}$.
It defines  a second-order  differential operator $\mcO_{\mbP^1} \migi \Omega^{\otimes 2}$ with unit principal symbol, or more precisely, 
 $D_{a, b, c}^\clubsuit := dx^{\otimes 2} \otimes \Dif_{a, b, c}^\clubsuit$ specifies  a global section of 
$\mcD {\it iff}^{\clubsuit}_{2, \mcO_{\mbP^1}}$.
If we write $y := x-1$, $z := 1/x$, then the following equalities hold:
\begin{align}
 D_{a,b, c}^\clubsuit &= \left(\frac{dx}{x}\right)^{\otimes 2} \otimes \left((x \partial_x)^2 + \frac{(a+b) \cdot x+ 1-c}{x-1} \cdot x \partial_x + \frac{ab \cdot  x}{x-1} \right) \\
&= \left( \frac{dy}{y}\right)^{\otimes 2} \otimes \left((y\partial_y)^2 + \frac{(a+b) \cdot y -c + a+b}{y+1} \cdot y \partial_y + \frac{ab \cdot y}{y+1} \right) \notag \\
& = \left( \frac{dz}{z}\right)^{\otimes 2} \otimes \left((z \partial_z)^2 + \frac{(1-c) \cdot z  + a+b}{z-1} \cdot z \partial_z - \frac{ab}{z-1} \right). \notag
\end{align}
These equalities show that  the radius  of $D^\clubsuit_{a,b, c}$ at the marked point $[0]$ (resp., $[1]$; resp., $[\infty]$) is given by  
\begin{align} \label{QG2000}
\rho_{[0]}^{D^\clubsuit_{a,b, c}} = \frac{(1-c)^2}{4} \ \left(\text{resp.,} \ \rho_{[1]}^{D^\clubsuit_{a,b, c}} = \frac{(c-a-b)^2}{4}; \text{resp.,} \ 
\rho_{[\infty]}^{D^\clubsuit_{a,b, c}} = \frac{(b-a)^2}{4}  \right).
\end{align} 
Here, recall that  the {\it exponent differences} of $\Dif^\clubsuit_{a,b,c}$ at $0$, $1$, $\infty$ are $1-c$, $c-a-b$, $b-a$ respectively; we shall write 
\begin{align}
\mr{Ex}_{a,b,c} :=(1-c, c-a-b, b-a),
\end{align}
 which will be regarded as a triple of elements in $k/\{ \pm 1 \}$ ($:=$ the quotient set of $k$ by the equivalence relation generated by $v \sim -v$ for every $v \in k$).

Next, let us write $\partial_0 := \frac{dx}{x (x-1)}$, which generates  $\mcT |_{\mbA^1}$.
Then, 
the  $(\mr{GL}_n, \mcO_{\mbP^1})$-oper 
$D_{a,b,c}^{\clubsuit \Rightarrow \diamondsuit}$ associated to $D_{a,b,c}^\clubsuit$  satisfies 
 \begin{align}
(D_{a,b,c}^{\clubsuit \Rightarrow \diamondsuit})^{\mr{GL}_2} (\partial_0) := \begin{pmatrix} 0& -ab \cdot x (x-1)\\  1 & (1-a-b)\cdot x -1+c  \end{pmatrix}.
\end{align}
Denote by  $\nabla_{a,b,c}$ the unique   $k^\mr{log}$-connection on $\mcT$ 
expressed, on $\mbA^1$,  as 
\begin{align}
\nabla_{a,b, c} := d +  dx \otimes \frac{(1-a-b) \cdot x -1 +c}{x (x-1)}
\end{align}
under   the identification $\mcO_{\mbA^1} \isom \mcT |_{\mbA^1}$ given by 
$v \mapsto v \cdot \frac{dx}{x (x-1)}$.
Then, the pair 
  $\bb_{a,b, c} := (\mcO_{\mbP^{1}}, \nabla_{a,b, c})$ specifies   a $2$-theta characteristic of $\mbP^{1\mr{log}}/k^\mr{log}$, and
$D_{a,b,c}^{\clubsuit}$ (resp., $D_{a,b,c}^{\clubsuit \Rightarrow \diamondsuit}$) lies in    $\mcD{\it iff}^{\clubsuit}_{2, \bb_{a,b,c}}$ (resp., $\mcO p^\diamondsuit_{n, \bb_{a,b,c}}$).

Let $\msE_{a,b,c}^\spadesuit$ denote the $\mfs \mfl_2$-oper on $\msP$ corresponding to $D_{a,b,c}^\clubsuit$ via the sequence of isomorphisms resulting from Theorem  \ref{y0129}.
Observe from the same theorem   that $\mfO \mfp_{\mfs \mfl_2,  \rho, \msP} \cong \mr{Spec}(k)$ for each $\rho \in \mfc^{\times r}(k)$.
That is to say,   (the isomorphism class of) an $\mfs \mfl_2$-oper on $\msP$ is
uniquely determined by its radii.
Hence, it follows from  (\ref{QG2000}) that 
 each  $\mfs \mfl_2$-oper on $\msP$ is isomorphic to $\msE_{a,b,c}^\spadesuit$ for some $(a,b,c)\in k^{\times 3}$, and  that
 \begin{align} \label{QH16}
 \msE^\spadesuit_{a, b, c} \cong  \msE^\spadesuit_{a', b', c'} \ \Longleftrightarrow \ \mr{Ex}_{a,b,c} = \mr{Ex}_{a', b', c'} \ \text{in} \ (k/\{ \pm 1 \})^{\times 3}.
 \end{align}

Hereafter, suppose further  that $\mr{char}(k)=:p>2$.
According to ~\cite{Ihara1}, \S\,1.6 (or ~\cite{Katz2}, \S\,6.4), 
the $(2, \bb_{a,b,c})$-projective connection  $D_{a, b, c}^\clubsuit$ has a full set of root functions  if and only if $(a,b,c)$ lies in $ \mbF_p^{\times 3}  \left(\subseteq k^{\times 3} \right)$ and  either $\widetilde{b} \geq  \widetilde{c} > \widetilde{a}$ or $\widetilde{a} \geq  \widetilde{c} > \widetilde{b}$ is satisfied, where $\widetilde{(-)}$ denotes the inverse of the  bijective restriction  $\{1, \cdots, p \} \isom \mbF_p$  of the natural quotient $\mbZ \migisurj \mbF_p$.
In particular,  the following equality can be verified 
after a straightforward calculation:
\begin{equation} \label{QG2003}
\sharp \begin{pmatrix}
\text{the set of hypergeometric operators} \\
\text{having  a full set of root functions}
\end{pmatrix}
=  \frac{p^3-p}{3}.
\end{equation}
If $D_{a, b, c}^\clubsuit$ has a full set of root functions, then
the two elements 
\begin{align}
_2 \overline{F}_1 (a,b; c;  x), \hspace{5mm} x^{1-\widetilde{c}}  \cdot  {_2 \overline{F}}_1 (a-c+1, b-c+1;  2-c; x)
\end{align}
of  $k (x)$ (i.e., the function field of $\mbP^1$)
forms a basis of the root functions.
Here,  $_2 \overline{F}_1 (a, b; c; x)$ denotes  a polynomial of $x$ defined by the following  truncated hypergeometric series:
\begin{align}
_2 \overline{F}_1 (a, b; c; x) := 
 & \  1 + \frac{a \cdot b}{1 \cdot c} \cdot x + \frac{a \cdot  (a+1) \cdot b \cdot  (b+1)}{1 \cdot 2 \cdot c \cdot (c +1)} \cdot x^2 \\
& +  \frac{a \cdot  (a+1) \cdot (a+2) \cdot b \cdot  (b+1) \cdot (b+2)}{1 \cdot 2 \cdot 3\cdot c \cdot (c +1) \cdot (c+2)} \cdot x^3 + \cdots,  \notag
\end{align}
where we stop the series as soon as the numerator vanishes.
Moreover, since
the exponent differences of $D_{a,b,c}^\clubsuit$  are all nonzero if it has a full set of root functions,
  the condition (\ref{QH16}) implies that the assignment $D_{a, b, c}^\clubsuit  \mapsto \msE^\spadesuit_{a, b, c}$ determines  a $2^3$-to-$1$ correspondence 
\begin{equation} \label{QG2010}
\begin{pmatrix}
\text{the set of hypergeometric operators} \\
\text{having  a full set of root functions}
\end{pmatrix}
\stackrel{2^3 : 1}{\longleftrightarrow}
\begin{pmatrix}
\text{the set of isomorphism classes} \\
\text{of dormant $\mfs \mfl_2$-opers on $\msP$}
\end{pmatrix}.
\end{equation}
In particular,  by this correspondence and (\ref{QH16}),
we can compute explicitly the number of the isomorphism classes of dormant $\mfs \mfl_2$-opers on $\msP$ as follows: 
\begin{align} \label{QG2004}
\sharp \mfO \mfp^\ZZZ_{\mfs \mfl_2, 1, \msP} (k) = \frac{p^3-p}{24}.
\end{align}

\begin{rema}
By taking account  of some results proved later, we see
 that  (\ref{QG2004}) also gives a computation of the number of the set $\mfO \mfp^\ZZZ_{\mfs \mfl_2, 1, \msX} (k)$ for a sufficiently  general  unpointed stable  curve  $\msX$ over $k$ of genus $2$.
Indeed,  let us consider  the unpointed stable curve $\msP_{2}$ of genus $2$ obtained by gluing together two copies of $\msP$ at the  respective corresponding marked points.
Then, since any  $\mfs \mfl_2$-oper on $\msP$ is uniquely determined by its radii,  the restrictions  of a dormant  $\mfs \mfl_2$-oper on $\msP_2$ to  the respective   components are  necessarily isomorphic (cf. Theorem \ref{y0176} for $\mfg = \mfs \mfl_2$); conversely, each dormant $\mfs \mfl_2$-oper on $\msP_2$ can be obtained by gluing together two isomorphic dormant $\mfs \mfl_2$-opers on $\msP$ at the points of attachment.
This implies that there exists a canonical bijection between $\mfO \mfp^\ZZZ_{\mfs \mfl_2, 1, \msP} (k)$ and $\mfO \mfp^\ZZZ_{\mfs \mfl_2, 1, \msP_2} (k)$.
On the other hand, as we will prove later (cf. Theorem  \ref{3778567d} and Remark \ref{QG91}), 
the equality  $\mfO \mfp^\ZZZ_{\mfs \mfl_2, 1, \msP_2} (k) = \mfO \mfp^\ZZZ_{\mfs \mfl_2, 1, \msX} (k)$ holds if  $\msX$ is 
 sufficiently general in $\overline{\mfM}_{2, 0}$.
Consequently,   for such a curve  $\msX$, we have $\mfO \mfp^\ZZZ_{\mfs \mfl_2,  1,\msX} (k) = \left(\mfO \mfp^\ZZZ_{\mfs \mfl_2, 1, \msP} (k) =\right) \frac{p^3-p}{24}$. 
This result is consistent with (\ref{n=2}) in Introduction.
\end{rema}

\vspace{10mm}

\chapter{Duality of  opers} \label{p0139} \vspace{3mm}

In this chapter, we discuss $(\mfg, \hslash)$-opers for $\mfg = \mfs \mfo_{2l+1}$ and $\mfs \mfp_{2m}$.
One important concept is   duality of  $(\mr{GL}_n, \hslash, \bb)$-opers, which 
 may be interpreted  as  taking the  transposes of the corresponding  differential operators.
In the cases of $\mfg = \mfs \mfo_{2l+1}$ and  $\mfg = \mfs \mfp_{2m}$,  we describe $(\mfg, \hslash)$-opers as  certain {\it self-dual} $(\mr{GL}_n, \hslash, \bb)$-opers (for a suitable $n$) equipped with a nondegenerate bilinear form.
This description  gives rise to a canonical closed immersion
from the moduli space of  (dormant) $(\mfg, \hslash)$-opers to that of (dormant) $(\mfs \mfl_n, \hslash)$-opers.
Moreover, the duality  yields a bijection  between (dormant) $(\mfg, \hslash)$-opers of radii $\rho$ and (dormant) $(\mfg, \hslash)$-opers of radii $(-1) \star \rho$ for each $\rho$. 
These observations  will be used to
construct  the pseudo-fusion rings associated to  dormant $(\mfg, \hslash)$-opers for such $\mfg$'s, as discussed in Chapter \ref{y0171}.
 
\begin{notation}
Let $n$ be an integer greater than $1$.
Also, let  $k$ 
\index{$k$, field} 
be a perfect field such that either ``\,$\mr{char}(k) =0$'' or ``\,$\mr{char}(k) > 2  n \left(= 2   h_{\mr{PGL}_n} \right)$'' is fulfilled.
\end{notation}
 
\section{$(\mfs \mfo_{2l+1}, \hslash)$-opers and $(\mfs \mfp_{2m}, \hslash)$-opers}\label{p0140p1}

In the following discussion,
the {\it non-resp'd} portion  deals with $(\mfs \mfo_{2l+1}, \hslash)$-opers  and 
the  {\it resp'd} portion deals with   $(\mfs \mfp_{2m}, \hslash)$-opers.
Suppose that there exists a positive integer $l$ (resp., $m$) with   $n = 2l+1$ (resp., $n =2m$).

\subsection{Algebraic groups: $\mr{GO}_{2l+1}$ and $\mr{GSp}_{2m}$} \label{QW530}

Denote  
 by $\mr{GO}_{2l+1}$
 \index{$\mbG$, algebraic group@--- $\mr{GO}_{2l+1}$, group of  orthogonal similitudes}
  (resp., $\mr{GSp}_{2m}$) 
  \index{$\mbG$, algebraic group@--- $\mr{GSp}_{2m}$, group of symplectic similitudes}
  the group of  orthogonal similitudes  of rank $2l+1$  (resp., the group of symplectic similitudes  of rank  $2m$).
 That is to say,
 \begin{align}
 \mr{GO}_{2l+1} := \left\{ A \in \mr{GL}_{2l+1} \, | \, {^t A} K_{2l+1} A =  c \cdot  K_{2l+1}  \ \text{for some $c \in \mbG_m$} \right\} \hspace{-5mm} \\
 \left(\text{resp.,} \  \mr{GSp}_{2m} :=  \left\{ A \in \mr{GL}_{2m} \, | \, {^t A} J_{2m} A = c \cdot  J_{2m} \ \text{for some $c \in \mbG_m$} \right\} \right),  \notag 
 \end{align}
where  $K_a := \begin{pmatrix} 0 & 0  &\cdots  & 0 & 1 \\  0 & 0 & \cdots & 1 & 0  \\ \vdots & \vdots &  \reflectbox{$\ddots$} & \vdots    & \vdots  \\  0  & 1 & \cdots  & 0 & 0 \\1 & 0  & \cdots & 0 & 0  \end{pmatrix}  \in \mr{GL}_a$
and 
$J_{2a} := \begin{pmatrix} 0 &  K_a \\ - K_a & 0\end{pmatrix}$ $\in \mr{GL}_{2a}$ for each $a \geq 1$.
Then, $\mr{GO}_{2l+1}$ (resp., $\mr{GSp}_{2m}$) and its adjoint group
\index{$\mbG$, algebraic group@--- $\mr{PGO}_{2l+1}$, $\mr{PGSp}_{2m}$}
\begin{align} \label{fff060}
\mr{PGO}_{2l+1} := (\mr{GO}_{2l+1})_\mr{ad} \ \left(\text{resp.,} \ \mr{PGSp}_{2m} := (\mr{GSp}_{2m})_\mr{ad}\right)
\end{align}
 are closed subgroups of $\mr{GL}_{n}$ and $\mr{PGL}_n$ respectively.
Note that, if $k$ is of positive characteristic, then 
 the inequality  $\mr{char}(k) > 2  h_{\mr{PGO}_{2l+1}}$ (resp., $\mr{char}(k) > 2  h_{\mr{PGSp}_{2m}}$) holds because of  the assumption on $\mr{char}(k)$.

Next, we write
\index{$\mfg$, Lie algebra@--- $\mfs \mfo_{2l+1}$, $\mfs \mfp_{2m}$}
\begin{align} \label{QH2231}
\mfs \mfo_{2l+1} := \left\{ v \in \mfg \mfl_{2l+1} \, | \, {^t}v K_{2l+1} + K_{2l+1}v = 0\right\} \hspace{-6mm}\\ 
\left(\text{resp.,} \ \mfs \mfp_{2m} := \left\{ v \in \mfg \mfl_{2m} \, | \, {^t}v J_{2m} + J_{2m}v = 0\right\} \right). \notag
\end{align}
Since $2l+1$ (resp., $2m$) is invertible in $k$ by assumption,
$\mfs \mfo_{2l+1}$ (resp., $\mfs \mfp_{2m}$)
may be naturally identified with the Lie algebra of 
$\mr{PGO}_{2l+1}$ (resp., $\mr{PGSp}_{2m}$).
Under this identification and the usual identification $\mfs \mfl_n = \mfp \mfg \mfl_n$, the closed immersion  $\mr{PGO}_{2l+1} \migiincl \mr{PGL}_n$ (resp., $\mr{PGSp}_{2m} \migiincl \mr{PGL}_n$) induces, via differentiation, an inclusion of Lie algebras $\mfs \mfo_{2l+1} \migiincl \mfs \mfl_n$
(resp., $\mfs \mfp_{2m} \migiincl \mfs \mfl_n$).
Also, the subgroups 
\begin{align} \label{QW1400}
T^{\mr{O}} := T \cap \mr{PGO}_{2l+1}
\ \   \text{and} \ \  
 B^{\mr{O}} := B \cap \mr{PGO}_{2l+1} \hspace{1mm} \\
 \left(\text{resp.,} \  T^{\mr{Sp}} := T \cap \mr{PGSp}_{2m} \ \ \text{and} \ \  B^{\mr{Sp}} := B \cap \mr{PGSp}_{2m}\right) \notag
\end{align}
(cf. (\ref{BorelB}) for the definitions of $T$ and $B$) specifies, respectively, a maximal torus and   a Borel subgroup of 
$\mr{GO}_{2q+1}$ (resp., $\mr{GSp}_{2m}$)   (cf., e.g.,  ~\cite{Malle}, Exercise 10.19 (b)).
We suppose that these choices of subgroups  are selected whenever we are dealing with  the algebraic group 
$\mr{PGO}_{2q+1}$ (resp., $\mr{PGSp}_{2m}$). 
The inclusion  $\mfs \mfo_{2l+1} \migiincl \mfs \mfl_n$
(resp., $\mfs \mfp_{2m} \migiincl \mfs \mfl_n$)
  induces   a closed immersion 
 \begin{align} \label{pp88}
 \mfc_{\mfs \mfo_{2l+1}} \migiincl  \mfc_{\mfs \mfl_n} \  \  \left(\text{resp.,}  \ \mfc_{\mfs \mfp_{2m}} \migiincl \mfc_{\mfs \mfl_n} \right)
 \end{align}
  of $k$-schemes.
  In particular, for a positive  integer $r$ and a $k$-scheme $S$,  the set $\mfc_{\mfs \mfo_{2l+1}}^{\times r} (S)$ (resp., $\mfc_{\mfs \mfp_{2m}}^{\times r} (S)$)  may be regarded as a subset of $\mfc_{\mfs \mfl_n}^{\times r} (S)$.

\subsection{The sheaves}

In what follows, let us fix an fs log scheme $S^\mr{log}$,  a log curve $U^\mr{log}$ over $S^\mr{log}$, and an element $\hslash$ of $\Gamma (S, \mcO_S)$.
For simplicity, we write $\Omega := \Omega_{U^\mr{log}/S^\mr{log}}$ and $\mcT := \mcT_{U^\mr{log}/S^\mr{log}}$.
Also,  we write
\index{$\mcO p^\spadesuit_{n, \hslash (, \rho)}$, $\mcO p^{\spadesuit, \mr{O}}_{n, \hslash}$, $\mcO p^{\spadesuit, \mr{Sp}}_{n, \hslash}$}
\begin{align} \label{QW1411}
\mcO p^{\spadesuit, \mr{O}}_{n, \hslash} := \mcO p_{\mfs \mfo_{2l+1}, \hslash, U^\mr{log}/S^\mr{log}} \ \left(\text{resp.,} \ \mcO p^{\spadesuit, \mr{Sp}}_{n, \hslash} := \mcO p_{\mfs \mfp_{2m}, \hslash, U^\mr{log}/S^\mr{log}}\right).
\end{align}
If  $\msE_{\sphericalangle}^\spadesuit := (\mcE_{B^\mr{O}}, \nabla)$ (resp., $\msE_{\sphericalangle}^\spadesuit:= (\mcE_{B^\mr{Sp}}, \nabla)$) is  an $(\mfs \mfo_{2l+1}, \hslash)$-oper (resp., an $(\mfs \mfp_{2m}, \hslash)$-oper) on $U^\mr{log}/S^\mr{log}$, then 
 the pair 
 \index{$\msE^\spadesuit$, $(\mfg, \hslash)$-oper@--- $\msE_{\sphericalangle \Rightarrow \emptyset}^\spadesuit$}
\begin{align} \label{QW1410}
\msE_{\sphericalangle \Rightarrow \emptyset}^\spadesuit := (\mcE_{B^\mr{O}} \times^{B^{\mr{O}}} B,  \nabla \times^{\mr{PGO}_{2l+1}} \mr{PGL}_n) \hspace{3mm}\\
\ \left(\text{resp.,} \ 
  \msE_{\sphericalangle \Rightarrow \emptyset}^\spadesuit:= (\mcE_{B^\mr{Sp}} \times^{B^{\mr{Sp}}} B,  \nabla \times^{\mr{PGSp}_{2m}} \mr{PGL}_n)\right) \notag
\end{align}
turns out to form an $(\mfs \mfl_n, \hslash)$-oper.
The assignment $\mcE_{\sphericalangle}^\spadesuit \mapsto \mcE_{\sphericalangle \Rightarrow \emptyset}^\spadesuit$ commutes with pull-back to \'{e}tale $U$-schemes and hence  defines a morphism of sheaves
\begin{align} \label{QW1154}
\mcO p^{\spadesuit, \mr{O}}_{n, \hslash} \migi \mcO p_{n, \hslash}^\spadesuit  \ 
\left(\text{resp.,} \ \mcO p^{\spadesuit, \mr{Sp}}_{n, \hslash} \migi \mcO p_{n, \hslash}^\spadesuit \right).
\end{align}

\section{$(\mr{GO}_{2l+1}, \hslash, \bb)$-opers and  $(\mr{GSp}_{2m}, \hslash, \bb)$-opers}
 \label{p0fg0p8}

We  shall introduce  $(\mr{GO}_{2l+1}, \hslash, \bb)$-opers (resp.,    $(\mr{GSp}_{2m}, \hslash, \bb)$-opers) defined for an $(n, \hslash)$-theta characteristic, and study the relationship with $(\mfs \mfo_{2l+1}, \hslash)$-opers (resp., $(\mfs \mfp_{2m}, \hslash)$-opers).

\subsection{$\mr{GO}_{n}$/$\mr{GSp}_{n}$-bundles with  a Borel reduction} \label{QW1000}
Let $\mcF$ and $\mcG$   be  vector bundles on $U$.
In what follows, we will not distinguish an $\mcO_U$-bilinear map $\mcF \times \mcF \migi \mcG$ on $\mcF$ valued in $\mcG$ with the corresponding $\mcO_U$-linear morphism $\mcF^{\otimes 2} \migi \mcG$.
If  $\omega : \mcF^{\otimes 2} \migi \mcG$ is an $\mcO_U$-bilinear map on $\mcF$,
then we write
\begin{align} \label{QW1567}
\omega_{(-, \bullet)} 
: \mcF \migi \mcF^\vee \otimes \mcG 
\index{$\omega_{(-, \bullet)}$, $\omega_{(\bullet, -)}$}
\end{align}
for the $\mcO_U$-linear morphism given by $v \mapsto \omega (v \otimes (-))$ for any local section $v \in \mcF$.

Now, let $B^{\mr{GO}}$ (resp., $B^{\mr{GSp}}$) be the subgroup of $\mr{GO}_{n}$ (resp., $\mr{GSp}_n$) consisting of invertible upper triangular matrices. 
In particular, the image of this subgroup via the quotient $\mr{GO}_{2l+1} \migisurj \mr{PGO}_{2l+1}$ (resp., $\mr{GSp}_{2m} \migisurj \mr{PGSp}_{2m}$) coincides with $B^{\mr{O}}$ (resp., $B^{\mr{Sp}}$).
 One verifies that giving  a $B^\mr{GO}$-bundle (resp., a $B^{\mr{GSp}}$-bundle) on $U$ is equivalent to giving a collection of data
\begin{align} \label{QW90gg0}
(\mcF, \{ \mcF^j \}_{j=0}^n, \mcN, \omega),
\end{align}
where 
\begin{itemize}
\item
$\mcF$ is a rank $n$ vector bundle on $U$;
\item
 $\{ \mcF^j \}_{j=0}^n$ is a decreasing filtration on $\mcF$ forming a complete flag, i.e., satisfying that $\mcF^0 =\mcF$, $\mcF^n =0$, and the  subquotient $\mcF^j/\mcF^{j+1}$ for every $j=0, \cdots, n-1$ is a line bundle;
 \item
  $\mcN$ is a line bundle on $U$; 
  \item
  $\omega$ is a nondegenerate {\it symmetric}  (resp., {\it skew-symmetric}) $\mcO_U$-bilinear map  $\mcF^{\otimes 2} \migi \mcN$ on $\mcF$ satisfying that   $\mcF^{n-j} = (\mcF^{j})^\bot \left(:= \mr{Ker}(\omega_{(-, \bullet)}|_{\mcF^{j}})\right)$ for any $j =0, \cdots, n$.
  \end{itemize}
Let us fix such a collection and denote by 
$\mcE_B$  the corresponding   $B^\mr{GO}$-bundle (resp., $B^\mr{GSp}$-bundle).
 Also, denote by $\mcE$ the $\mr{GO}_{2l+1}$-bundle (resp., the $\mr{GSp}_{2m}$-bundle) induced by  this bundle  via change of structure group. 
 Then, 
giving an $\hslash$-$S^\mr{log}$-connection  on $\mcE$ 
 is equivalent to
giving a pair 
\begin{align}
(\nabla, \nabla_\mcN)
\end{align}
consisting of an $\hslash$-$S^\mr{log}$-connection $\nabla$ on $\mcF$ and an $\hslash$-$S^\mr{log}$-connection $\nabla_\mcN$ on $\mcN$ such that $\omega$ is compatible with $\nabla^{\otimes 2}$ and $\nabla_\mcN$.

\subsection{The definitions}
 \label{QG1010}
 
 Let us introduce the definition of a $(\mr{GO}_{2l+1}, \hslash)$-opers (resp., a $(\mr{GSp}_{2m}, \hslash)$-oper) as 
 a $(\mr{GL}_n, \hslash)$-oper equipped with an additional nondegenerate  bilinear map.
 
\bde \label{pp1}
\begin{itemize}
\item[(i)]
Let us consider a collection of data 
\begin{align}
  \msF^{\heartsuit }_\sphericalangle
  := (\mcF, \Nab, \{ \mcF^j \}_{j=0}^n, \mcN, \nabla_\mcN, \omega),
\end{align}
where
\begin{itemize}
\item
$\msF^\heartsuit_{\sphericalangle \Rightarrow \emptyset} := (\mcF, \Nab, \{ \mcF^j \}_{j=0}^n)$
  \index{$\msF^{\heartsuit }_\sphericalangle$, $(\mr{GO}_{2l+1}, \hslash)$-oper, $(\mr{GSp}_{2m}, \hslash)$-oper@--- $\msF^\heartsuit_{\sphericalangle \Rightarrow \emptyset}$}
 is a $(\mr{GL}_n, \hslash)$-oper on $U^\mr{log}/S^\mr{log}$;
\item
$\msN := (\mcN, \nabla_\mcN)$ is an $\hslash$-flat line bundle  on $U^\mr{log}/S^\mr{log}$;
\item
$\omega$ denotes a nondegenerate {\it symmetric} (resp., {\it skew-symmetric}) $\mcO_U$-bilinear map  $\mcF^{\otimes 2} \migi \mcN$ on $\mcF$ valued in $\mcN$.
\end{itemize}
Then, we shall say that $\msF^{\heartsuit }_\sphericalangle$ 
 is a {\bf $(\mr{GO}_{2l+1}, \hslash)$-oper} (resp., {\bf $(\mr{GSp}_{2m}, \hslash)$-oper})  
 \index{$\msF^{\heartsuit }_\sphericalangle$, $(\mr{GO}_{2l+1}, \hslash)$-oper, $(\mr{GSp}_{2m}, \hslash)$-oper}
 on $U^\mr{log}/S^\mr{log}$ if it satisfies
the following two conditions:
\begin{itemize}
\item
$\omega$ is compatible with the $\hslash$-$S^\mr{log}$-connections $\nabla^{\otimes 2}$ and $\nabla_\mcN$;
\item
 For any  $j = 0, \cdots, n$,     the equality $\mcF^j = (\mcF^{n-j})^\bot$ holds with respect to $\omega$. 
\end{itemize}
We refer to   $\msF^\heartsuit_{\sphericalangle \Rightarrow \emptyset}$    as the {\bf underlying $(\mr{GL}_n, \hslash)$-oper} of $\msF^{\heartsuit }_\sphericalangle$.

\item[(ii)]
Let 
$\msF^{\heartsuit }_\sphericalangle  := (\msF^\heartsuit, \msN, \omega)$ and 
$\msF'^{\heartsuit }_\sphericalangle := (\msF'^\heartsuit, \msN', \omega')$ be 
 $(\mr{GO}_{2l+1}, \hslash)$-opers 
  (resp., $(\mr{GSp}_{2m}, \hslash)$-opers) 
    on $U^\mr{log}/S^\mr{log}$.
Then, an {\bf isomorphism of  $(\mr{GO}_{2l+1}, \hslash)$-opers} (resp., {\bf $(\mr{GSp}_{2m}, \hslash)$-opers}) 
 from $\msF^{\heartsuit }_\sphericalangle$
to $\msF'^{\heartsuit }_\sphericalangle$
is a pair 
\begin{align} \label{QW588}
(\eta_\msF, \eta_\msN)
\end{align}
 consisting of an isomorphism $\eta_\msF : \msF^\heartsuit \isom \msF'^\heartsuit$ between the underlying  $(\mr{GL}_n, \hslash)$-opers and an isomorphism $\eta_\msN : \msN \isom \msN'$ of $\hslash$-flat line bundles satisfying  $\eta_\msN \circ \omega = \omega' \circ \eta_\msF^{\otimes 2}$.

\end{itemize}
 \ede

\subsection{Twisting with a flat line bundle} \label{QW600}
Let $\msF^{\heartsuit }_\sphericalangle := (\mcF, \nabla_\mcF, \{ \mcF^j \}_{j}, \mcN, \nabla_\mcN, \omega)$ be a $(\mr{GO}_{2l+1}, \hslash)$-oper (resp., $(\mr{GSp}_{2m}, \hslash)$-oper) on $U^\mr{log}/S^\mr{log}$ and $\msL := (\mcL, \nabla_\mcL)$ an $\hslash$-flat line bundle on $U^\mr{log}/S^\mr{log}$.
Then,  
\begin{align} \label{QW602}
\omega_{\otimes \mcL} := \omega \otimes (\mr{id}_{\mcL^{\otimes 2}}) : (\mcF \otimes \mcL)^{\otimes 2} \left(= \mcF^{\otimes 2} \otimes \mcL^{\otimes 2} \right) \migi \mcN \otimes \mcL^{\otimes 2}
\end{align}
is a symmetric (resp., a skew-symmetric) $\mcO_U$-bilinear map on $\mcF \otimes \mcL$, and  
 the resulting collection 
 \index{$\msF^{\heartsuit }_\sphericalangle$, $(\mr{GO}_{2l+1}, \hslash)$-oper, $(\mr{GSp}_{2m}, \hslash)$-oper@--- $\msF^{\heartsuit }_{\sphericalangle, \otimes \msL}$}
\begin{align} \label{QW601}
\msF^{\heartsuit }_{\sphericalangle, \otimes \msL}
   := (\mcF \otimes \mcL, \nabla \otimes \nabla_\mcL, \{ \mcF^j \otimes \mcL \}_{j}, \mcN \otimes \mcL^{\otimes 2}, \nabla_\mcN \otimes \nabla^{\otimes 2}_\mcL, \omega_{\otimes \mcL}) 
\end{align}
forms a $(\mr{GO}_{2l+1}, \hslash)$-oper (resp., a $(\mr{GSp}_{2m}, \hslash)$-oper) on $U^\mr{log}/S^\mr{log}$.

\bde  \label{pp2}
 Let  $\msF^{\heartsuit }_{\sphericalangle}$ and $\msF'^{\heartsuit }_{\sphericalangle}$  be $(\mr{GO}_{2m}, \hslash)$-opers (resp., $(\mr{GSp}_{2m}, \hslash)$-opers) 
\index{$\msF^{\heartsuit }_\sphericalangle$, $(\mr{GO}_{2l+1}, \hslash)$-oper, $(\mr{GSp}_{2m}, \hslash)$-oper@--- equivalent}
 on $U^\mr{log}/S^\mr{log}$.
 We shall say that 
 $\msF^{\heartsuit }_{\sphericalangle}$  is {\bf equivalent to} 
 $\msF'^{\heartsuit }_{\sphericalangle}$
  if there exists an $\hslash$-flat  line bundle $\msL$ on $U^\mr{log}/S^\mr{log}$ with   $\msF^{\heartsuit }_{\sphericalangle, \otimes\msL} \cong \msF'^{\heartsuit }_{\sphericalangle}$.
  We use the notation $\msF^{\heartsuit }_{\sphericalangle} \stackrel{\heartsuit}{\sim} \msF'^{\heartsuit }_{\sphericalangle}$ to indicate the situation that $\msF^{\heartsuit }_{\sphericalangle}$  is  equivalent to 
 $\msF'^{\heartsuit }_{\sphericalangle}$.
  (It is immediate that the binary relation ``$\stackrel{\heartsuit}{\sim}$'' on the set of $(\mr{GO}_{2m}, \hslash)$-opers (resp., $(\mr{GSp}_{2m}, \hslash)$-opers) on $U^\mr{log}/S^\mr{log}$ defines an equivalence relation.)
  \ede

\subsection{The sheaves}
Denote by 
\begin{align}
\overline{O} p_{n, \hslash}^{\heartsuit, \mr{O}} \ \left(\text{resp.,} \  \overline{O} p_{n, \hslash}^{\heartsuit, \mr{Sp}}\right): \mfE \mft_{/U} \migi \mfS \mfe \mft
\index{$\overline{O} p_{n, \hslash}^{\heartsuit, \mr{O}}$, $\overline{O} p_{n, \hslash}^{\heartsuit, \mr{Sp}}$}
\end{align}
the $\mfS \mfe \mft$-valued sheaf on $\mfE \mft_{/U}$ associated to the presheaf which, to any \'{e}tale $U$-scheme $U'$, assigns the set of {\it equivalence classes} of $(\mr{GO}_{2l+1}, \hslash)$-opers (resp., $(\mr{GSp}_{2m}, \hslash)$-opers) on $U'^{\mr{log}}/S^\mr{log}$.

By taking account of the discussion in  \S\,\ref{QW1000} and carrying out a  change of structure group by the quotient $\mr{GO}_{2l+1} \migisurj \mr{PGO}_{2l+1}$ (resp., $\mr{GSp}_{2m} \migisurj \mr{PGSp}_{2m}$),
we can construct an $(\mfs \mfo_{2l+1})$-oper   (resp., $(\mfs \mfp_{2m}, \hslash)$-oper)
\index{$\msF^{\heartsuit }_\sphericalangle$, $(\mr{GO}_{2l+1}, \hslash)$-oper, $(\mr{GSp}_{2m}, \hslash)$-oper@--- $\msF_\sphericalangle^{\heartsuit \Rightarrow \spadesuit}$}
\begin{align} \label{QH12}
\msF_\sphericalangle^{\heartsuit \Rightarrow \spadesuit}
\end{align}
associated to 
a $(\mr{GO}_{2l+1}, \hslash)$-oper (resp., $(\mr{GSp}_{2m}, \hslash)$-oper) $\msF_\sphericalangle^\heartsuit$.
The isomorphism class of $\msF_\sphericalangle^{\heartsuit \Rightarrow \spadesuit}$ depends only on the equivalence class of $\msF_\sphericalangle^\heartsuit$.

Then, we can prove the following assertion.
(In the proof of this assertion, we will use  the notations in   Lemma \ref{QW1300} and \S\,\ref{p0140p5}.)


\bpr \label{QW910}
Let $\BL \in \{ \mr{O}, \mr{Sp}\}$.
Then, 
the assignment $\msF_\sphericalangle^\heartsuit \mapsto \msF_\sphericalangle^{\heartsuit \Rightarrow \spadesuit}$ defines an isomorphism of sheaves
\begin{align}
\Lambda^{\heartsuit \Rightarrow \spadesuit, \BL} : \overline{O} p_{n, \hslash}^{\heartsuit, \BL} \migi \mcO p^{\spadesuit, \BL}_{n, \hslash}
\index{$\Lambda^{\heartsuit \Rightarrow \spadesuit}_{(\rho)}$, $\Lambda^{\heartsuit \Rightarrow \spadesuit, \BL}$, $\Lambda^{\heartsuit \Rightarrow \spadesuit, \BL}_\msX$},
\end{align}
which makes the following square diagram commute:
\begin{align} \label{QW1101ww}
\vcenter{\xymatrix@C=46pt@R=36pt{
\overline{\mcO} p_{h, \hslash}^{\heartsuit, \BL}  \ar[r]_-{\sim}^-{\Lambda^{\heartsuit \Rightarrow \spadesuit, \BL}} \ar[d]_-{\msF_{\sphericalangle}^\heartsuit \mapsto \msF_{\sphericalangle \Rightarrow \emptyset}^\heartsuit}&  O p^{\spadesuit, \BL}_{n, \hslash} 
 \ar[d]^-{\msE_{\sphericalangle}^\spadesuit \mapsto \msE_{\sphericalangle \Rightarrow \emptyset}^\spadesuit}\\
 \overline{\mcO} p_{h, \hslash}^{\heartsuit} \ar[r]^-{\sim}_-{\Lambda^{\heartsuit \Rightarrow \spadesuit}}&  O p^{\spadesuit}_{n, \hslash}.
}}
\end{align}
Finally, the formation of this isomorphism commutes with base-change to fs log schemes over $S^\mr{log}$.
\epr
\begin{proof}
We only consider the case of $\BL = \mr{O}$ because the proof of the other is entirely similar.
Since any $B^{\mr{O}}$-bundle can be,  \'{e}tale locally,   lifted to a $B^{\mr{GO}}$-bundle,  the surjectivity of $\Lambda^{\heartsuit \Rightarrow \spadesuit, \mr{O}}$   is immediate from  an argument similar to that of  $\Lambda^{\heartsuit \Rightarrow \spadesuit}$ proved in  Proposition \ref{y0112}.

Next, we shall consider the injectivity of $\Lambda^{\heartsuit \Rightarrow \spadesuit, \mr{O}}$.
Let $\msF^{\heartsuit}_{\sphericalangle} := (\mcF, \nabla, \{ \mcF^j \}_j, \mcN, \nabla_\mcN, \omega)$ and 
$\msF'^{\heartsuit}_{\sphericalangle} := (\mcF', \nabla', \{ \mcF'^j \}_j, \mcN', \nabla_{\mcN'}, \omega')$ be  
$(\mr{GO}_n, \hslash)$-opers on $U^\mr{log}/S^\mr{log}$ with $\msF^{\heartsuit \Rightarrow \spadesuit}_{\sphericalangle} \cong \msF'^{\heartsuit \Rightarrow \spadesuit}_{\sphericalangle}$.
The problem is to prove that
$\msF^{\heartsuit}_{\sphericalangle} \stackrel{\heartsuit}{\sim}\msF'^{\heartsuit}_{\sphericalangle}$.
Because of  the injectivity of $\Lambda^{\heartsuit \Rightarrow \spadesuit}$ (cf. Proposition \ref{y0112}),
the underlying $(\mr{GL}_n, \hslash)$-oper  $\msF^{\heartsuit}_{\sphericalangle \Rightarrow \emptyset}$  of  $\msF^{\heartsuit}_{\sphericalangle}$ is  equivalent to  that of $\msF'^{\heartsuit}_{\sphericalangle}$, i.e., $\msF'^{\heartsuit}_{\sphericalangle \Rightarrow \emptyset}$.
Hence, by replacing $U$ with its covering in the \'{e}tale topology and   tensoring $\msF'^{\heartsuit}_{\sphericalangle}$ with an $\hslash$-flat line bundle, if necessary,  we may assume that
there exists an isomorphism  $\eta_\msF : \msF^{\heartsuit}_{\sphericalangle \Rightarrow \emptyset}  \isom \msF'^{\heartsuit}_{\sphericalangle \Rightarrow \emptyset}$ of
$(\mr{GL}_n, \hslash)$-opers.
Let us here apply  Lemma \ref{QW1300}, (i), which will be described later.
Then, we see that  $\eta_\msF$ induces an isomorphism $\eta_\msN : (\mcN, \nabla_\mcN) \isom (\mcN', \nabla_{\mcN'})$ of $\hslash$-flat line bundles.
According to  Proposition \ref{QW1045}, 
there exists an element  $a$  of $\Gamma (U, \mcO_U^{\times} \cap \mr{Ker}(\mcO_U \xrightarrow{\hslash \cdot d} \Omega))$ satisfying the equality
\begin{align} \label{QH13}
\omega'_{(-, \bullet)}\circ \eta_\msF = \Delta (a) \circ (\eta^\dual_\msF  \otimes \eta_\msN) \circ  \omega_{(-, \bullet)}
\end{align}
 of isomorphisms $\msF^\heartsuit_{\sphericalangle \Rightarrow \emptyset} \isom (\msF'^{\heartsuit \dual}_{\sphericalangle \Rightarrow \emptyset})_{\otimes \msN'}$
(cf. (\ref{pp43}) for the definition of ``$(-)^\dual$''), where $\msN':=(\mcN', \nabla_{\mcN'})$,   $\Delta$ \index{$\Delta$} denotes the diagonal embedding $\mcO_U \migiincl \mcE nd_{\mcO_U}(\mcF')$, and  $\eta^\dual_\msF$ denotes the isomorphism $\msF^{\heartsuit \dual}_{\sphericalangle \Rightarrow \emptyset} \isom \msF'^{\heartsuit \dual}_{\sphericalangle \Rightarrow \emptyset}$ induced naturally by $\eta_\msF$. 
If $\eta_a$ denotes the automorphism of $\msN$ given by multiplication by $a$,
then (\ref{QH13}) implies that  the pair $(\eta_\msF, \eta_a \circ \eta_\msN)$ forms an isomorphism
$\msF^\heartsuit_{\sphericalangle} \isom \msF'^\heartsuit_{\sphericalangle}$.
In particular, we have $\msF^{\heartsuit}_{\sphericalangle} \stackrel{\heartsuit}{\sim}\msF'^{\heartsuit}_{\sphericalangle}$, and  this completes 
 the injectivity of $\Lambda^{\heartsuit \Rightarrow \spadesuit, \mr{O}}$.

Finally, the commutativity of (\ref{QW1101ww}) and the last assertion follows immediately  from the definition of $\Lambda^{\heartsuit \Rightarrow \spadesuit, \BL}$.
\end{proof}


\subsection{$(\mr{GO}_{2l+1}, \hslash, \bb)$ and $(\mr{GSp}_{2m}, \hslash, \bb)$-opers}
 \label{p0140p2}

Next, let $\bb := (\BB, \nabla_\bb)$ be an $(n, \hslash)$-theta characteristic of  $U^\mr{log}/S^\mr{log}$.
We shall write
\begin{align} \label{QW923}
{^\dagger}\mcN_\bb := \mcT^{\otimes (n-1)} \otimes \BB^{\otimes 2}.
\index{${^\dagger}\mcN_\bb$}
\end{align}
It follows from  Proposition  \ref{y0111}, (i), that there exists a unique $\hslash$-$S^\mr{log}$-connection
${^\dagger}\nabla_{\mcN_\bb}$
on ${^\dagger}\mcN_\bb$ satisfying   ${^\dagger}\nabla_{\mcN_\bb}^{\otimes n} = \nabla_\bb^{\otimes 2}$ under the natural identification  ${^\dagger}\mcN_\bb^{\otimes n} = (\mcT^{\otimes \frac{n(n-1)}{2}} \otimes \BB^{\otimes n})^{\otimes 2}$.
Thus, we have an $\hslash$-flat line bundle 
\begin{align} \label{QW1302}
{^\dagger}\msN_\bb := ({^\dagger}\mcN_\bb, {^\dagger}\nabla_{\mcN_\bb})
\index{${^\dagger}\msN_\bb$}
\end{align}
on $U^\mr{log}/S^\mr{log}$.

\bde  \label{pp15}
\begin{itemize}
\item[(i)]
By a {\bf $(\mr{GO}_{2l+1}, \hslash, \bb)$-oper} (resp.,  a {\bf $(\mr{GSp}_{2m}, \hslash, \bb)$-oper}) 
\index{$\nabla^\diamondsuit_{\sphericalangle}$, $(\mr{GO}_{2l+1}, \hslash, \bb)$-oper, $(\mr{GSp}_{2m}, \hslash, \bb)$-oper}
on $U^\mr{log}/S^\mr{log}$, we mean   a pair 
\begin{align} \label{QW611}
\nabla^\diamondsuit_{\sphericalangle} := (\nabla^\diamondsuit,  \omega)
\end{align}
consisting of a $(\mr{GL}_{n}, \hslash, \bb)$-oper $\nabla^\diamondsuit$ on $U^\mr{log}/S^\mr{log}$ and a nondegenerate  {\it symmetric} (resp.,   {\it skew-symmetric}) $\mcO_U$-bilinear map  $\omega : \DF_\BB^{\otimes 2} \migi {^\dagger}\mcN_\bb$ on $\DF_\BB$ (cf. (\ref{QW47})) such that the collection
 \index{$\nabla^\diamondsuit_{\sphericalangle}$, $(\mr{GO}_{2l+1}, \hslash, \bb)$-oper, $(\mr{GSp}_{2m}, \hslash, \bb)$-oper@--- $\nabla^{\diamondsuit \Rightarrow \heartsuit}_{\sphericalangle}$}
\begin{align} \label{QW610}
\nabla^{\diamondsuit \Rightarrow \heartsuit}_{\sphericalangle} :=
 (\DF_\BB, \nabla^\diamondsuit, \{ \DF^j_\BB \}_{j=0}^n, {^\dagger}\mcN_\bb,  {^\dagger}\nabla_{\mcN_\bb}, \omega)
\end{align}
forms a $(\mr{GO}_{2l+1}, \hslash)$-oper (resp.,  $(\mr{GSp}_{2m}, \hslash)$-oper) on $U^\mr{log}/S^\mr{log}$.
If $\nabla_{\sphericalangle}^\diamondsuit$ is as above, then we occasionally write
$\nabla_{\sphericalangle \Rightarrow \emptyset}^\diamondsuit := \nabla^\diamondsuit$ and refer to it as the {\bf underlying $(\mr{GL}_n, \hslash, \bb)$-oper} of $\nabla_{\sphericalangle}^\diamondsuit$.
\item[(ii)]
Let $\nabla^\diamondsuit_{\sphericalangle}$ and $\nabla'^{\diamondsuit}_{\sphericalangle}$ be $(\mr{GO}_{2l+1}, \hslash, \bb)$-opers (resp., $(\mr{GSp}_{2m}, \hslash, \bb)$-opers) on $U^\mr{log}/S^\mr{log}$.
Then, an {\bf isomorphism of $(\mr{GO}_{2l+1}, \hslash, \bb)$-opers} (resp., {\bf $(\mr{GSp}_{2m}, \hslash, \bb)$-opers})
from $\nabla^\diamondsuit_{\sphericalangle}$ to  $\nabla'^{\diamondsuit}_{\sphericalangle}$ is defined as  an isomorphism $\nabla^{\diamondsuit \Rightarrow \heartsuit}_{\sphericalangle} \isom \nabla'^{\diamondsuit \Rightarrow \heartsuit}_{\sphericalangle}$ of  $(\mr{GO}_{2l+1}, \hslash)$-opers (resp., $(\mr{GSp}_{2m}, \hslash)$-opers).
\end{itemize}
  \ede

\ble \label{QW1300}
Let $\msF^\heartsuit_{\sphericalangle} := (\mcF, \nabla, \{ \mcF^j \}_j, \mcN, \nabla_\mcN, \omega)$ be a $(\mr{GO}_{2l+1},  \hslash)$-oper (resp., a $(\mr{GSp}_{2m},  \hslash)$-oper) on $U^\mr{log}/S^\mr{log}$.
Then, the following assertions hold:
\begin{itemize}
\item[(i)]
Denote by 
$\nabla_\circledcirc$ the unique $\hslash$-$S^\mr{log}$-connection on $\mcT^{\otimes (n-1)} \otimes (\mcF^{n-1})^{\otimes 2}$  such that    $\mr{det}(\nabla)^{\otimes 2}$ corresponds to $\nabla_\circledcirc^{\otimes n}$ via the composite isomorphism
\begin{align}
\mr{det}(\mcF)^{\otimes 2} \xrightarrow{(\ref{QW1906})^{\otimes 2}} (\mcT^{\otimes \frac{n(n-1)}{2}} \otimes (\mcF^{n-1})^{\otimes n})^{\otimes 2} \isom (\mcT^{\otimes (n-1)} \otimes (\mcF^{n-1})^{\otimes 2})^{\otimes n}
\end{align}
(cf. Proposition \ref{y0111}, (i)).
Then, the $\hslash$-flat line bundle $(\mcN, \nabla_\mcN)$ is isomorphic to
$(\mcT^{\otimes (n-1)} \otimes (\mcF^{n-1})^{\otimes 2}, \nabla_\circledcirc)$.
In particular,  the isomorphism class of  $(\mcN, \nabla_\mcN)$
 does not depend on the choice of $\omega$.
\item[(ii)]
If, moreover,  $\mcF^{\heartsuit}_{\sphericalangle \Rightarrow \emptyset}$ can be described as $\mcF^{\heartsuit}_{\sphericalangle \Rightarrow \emptyset} = (\DF_\BB, \nabla^\diamondsuit, \{ \DF_\BB^j \}_j)$
 for some $(\mr{GL}_n, \hslash, \bb)$-oper $\nabla^\diamondsuit$, then 
 $(\mcN, \nabla_\mcN)$ is isomorphic to $({^\dagger}\mcN_\bb, {^\dagger}\nabla_{\mcN_\bb})$.
 \end{itemize}
\ele
\begin{proof}
We first  consider assertion (i).
The   isomorphism   $\omega_{(-, \bullet)} : \mcF\isom \mcF^\vee \otimes \mcN$   restricts to an isomorphism $\mcF^{n-1} \isom \left( (\mcF/\mcF^1)^{\vee} \otimes \mcN  =\right) (\Omega^{\otimes (n-1)} \otimes \mcF^{n-1\vee}) \otimes \mcN$.
It induces  an isomorphism $\mcN \isom \mcT^{\otimes (n-1)} \otimes (\mcF^{n-1})^{\otimes 2}$, by which  we  identify  $\mcN$ with $\mcT^{\otimes (n-1)} \otimes (\mcF^{n-1})^{\otimes 2}$. 
Notice that the isomorphism $(\mr{det}(\mcF), \mr{det} (\nabla))^{\otimes 2} \isom
(\mcT^{\otimes (n-1)} \otimes (\mcF^{n-1})^{\otimes 2}, \nabla_\mcN)^{\otimes n}$
 arising from the determinant of $\omega_{(-, \bullet)}$ coincides with 
the $2$-fold tensor product  of  (\ref{fff0100}).
Hence,  it follows from  the definition of $\nabla_\circledcirc$ that  $\nabla_\mcN^{\otimes n} = \left(\mr{det}(\nabla)^{\otimes 2}=\right) \nabla_\circledcirc^{\otimes n}$.
 By the uniqueness portion of Proposition  \ref{y0111}, (i),  we have $\nabla_\mcN = \nabla_\circledcirc$.
This completes  the proof of  assertion (i).

Assertion (ii) follows immediately  from (i),  so we finish the proof of this lemma.
\end{proof}

\subsection{The presheaves}
 \label{QW1500}
Denote by 
 \index{$\mcO p_{n, \hslash, \BB}^\diamondsuit$, $\mcO p_{n, \hslash, \bb (, \rho)}^\diamondsuit$, $\mcO p_{n, \hslash, \BB}^{\diamondsuit, \mr{O}}$, $\mcO p_{n, \hslash, \BB}^{\diamondsuit, \mr{Sp}}$}
\begin{align} \label{QW900}
\mcO p_{n, \hslash, \bb}^{\diamondsuit, \mr{O}} \ \left(\text{resp.,} \  \mcO p_{n, \hslash, \bb}^{\diamondsuit, \mr{Sp}}\right): \mfE \mft_{/U} \migi \mfS \mfe \mft
\end{align}
the $\mfS \mfe \mft$-valued presheaf on $\mfE \mft_{/U}$ which, to any \'{e}tale $U$-scheme $u: U' \migi U$, assigns the set of isomorphism classes of $(\mr{GO}_{2l+1}, \hslash, u^*(\bb))$-opers (resp., $(\mr{GSp}_{2m}, \hslash, u^*(\bb))$-opers) on $U'^{\mr{log}}/S^\mr{log}$.

\bpr \label{QW911}
Let $\BL \in \{ \mr{O}, \mr{Sp}\}$.
Then, the assignment $\nabla^\diamondsuit_{\sphericalangle} \mapsto [\nabla^{\diamondsuit \Rightarrow \heartsuit}_{\sphericalangle}]$ defines a well-defined isomorphism 
 \begin{align}
 \Lambda_\bb^{\diamondsuit \Rightarrow \heartsuit, \BL} : \mcO p_{n, \hslash, \bb}^{\diamondsuit, \BL} \isom \overline{O} p^{\heartsuit, \BL}_{n, \hslash}, 
 \index{$\Lambda_{\bb (, \rho)}^{\diamondsuit \Rightarrow \heartsuit}$, $\Lambda_\bb^{\diamondsuit \Rightarrow \heartsuit, \BL}$, $\Lambda_{\bb, \msX}^{\diamondsuit \Rightarrow \heartsuit, \BL}$}
 \end{align}
 which makes the following square diagram commute:
\begin{align} \label{QW1101}
\vcenter{\xymatrix@C=46pt@R=36pt{
\mcO p_{n, \hslash, \bb}^{\diamondsuit, \BL}  \ar[r]_-{\sim}^-{\Lambda_\bb^{\diamondsuit \Rightarrow \heartsuit, \BL}} \ar[d]_-{\nabla_{\sphericalangle}^\diamondsuit \mapsto \nabla_{\sphericalangle \Rightarrow \emptyset}^\diamondsuit}&  \overline{O} p^{\heartsuit, \BL}_{n, \hslash}  \ar[d]^-{\msF_{\sphericalangle}^\heartsuit \mapsto \msF_{\sphericalangle \Rightarrow \emptyset}^\heartsuit}\\
 \mcO p_{h, \hslash, \bb}^{\diamondsuit} \ar[r]^-{\sim}_-{\Lambda_\bb^{\diamondsuit \Rightarrow \heartsuit}}&  \overline{O} p^{\heartsuit}_{n, \hslash}.
}}
\end{align}
In particular, $\mcO p^{\diamondsuit, \BL}_{n, \hslash, \bb}$ is a sheaf.
Finally, the formation of this isomorphism commutes with base-change to fs log schemes over $S^\mr{log}$.
\epr
\begin{proof}
The assertion can be proved   by 
an argument similar to the argument in the proof of Proposition \ref{y0127} together with  Lemma \ref{QW1300}, (ii) and (ii),  so we omit the details.
\end{proof}

\section{Dual opers}\label{p0140p5}

In this section, we shall consider duality  of $(\mr{GL}_n, \hslash)$-opers
and describe $(\mr{GO}_{2l+1}, \hslash)$-opers (resp.,  $(\mr{GSp}_{2m}, \hslash)$-opers) in terms of the duality.

\subsection{Dual $(\mr{GL}_n, \hslash)$-opers} \label{QW1021}
Let $\msF^\heartsuit := (\mcF, \nabla, \{ \mcF^j \}_{j=0}^n)$ be a $(\mr{GL}_n, \hslash)$-oper on $U^\mr{log}/S^\mr{log}$.
For each $j = 0, \cdots, n$, we
regard $\mcF^{\vee j} := (\mcF/\mcF^{n-j})^\vee$ as a subbundle of $\mcF^\vee$.
Then, one verifies that the  collection of data
\begin{align} \label{tt9}
\msF^{\heartsuit \dual} := (\mcF^\vee, \nabla^\vee, \{ \mcF^{\vee j} \}_{j=0}^n)
\end{align}
 forms a $(\mr{GL}_n, \hslash)$-oper on $U^\mr{log}/S^\mr{log}$.
By construction,   $(\msF^{\heartsuit\dual})^\dual$ is naturally isomorphic to $\msF^\heartsuit$.

Next, suppose that we are given    a $(\mr{GO}_{2l+1}, \hslash)$-oper  (resp., a $(\mr{GSp}_{2m}, \hslash)$-oper)
 $\msF^{\heartsuit}_{\sphericalangle}:= (\mcF, \nabla, \{ \mcF^j \}_{j}, \mcN, \nabla_\mcN, \omega)$ on $U^\mr{log}/S^\mr{log}$.
 Then, the composite
\begin{align} \label{QW918}
\omega^\dual 
 : (\mcF^\vee)^{\otimes 2} &\isom (\mcF \otimes \mcN^\vee)^{\otimes 2} \left(= \mcF^{\otimes 2} \otimes (\mcN^\vee)^{\otimes 2} \right) \\
 & \xrightarrow{\omega \otimes \mr{id}_{(\mcN^\vee)^{\otimes 2}}} \left(\mcN \otimes (\mcN^\vee)^{\otimes 2} = \right)  \mcN^\vee,  \notag
 \index{$\omega^\dual$}
\end{align}
where the first arrow arises from the isomorphism $\omega_{(-, \bullet)} : \mcF \isom \mcF^\vee \otimes \mcN$, defines a symmetric (resp., a skew-symmetric)  $\mcO_U$-bilinear map 
on $\mcF^\vee$.
The collection of data
\index{$\msF^{\heartsuit }_\sphericalangle$, $(\mr{GO}_{2l+1}, \hslash)$-oper, $(\mr{GSp}_{2m}, \hslash)$-oper@--- $\msF^{\heartsuit \dual}_{\sphericalangle}$, dual of}
\begin{align} \label{pp43}
 \msF^{\heartsuit \dual}_{\sphericalangle}
  := (\mcF^\vee, \nabla^\vee, \{ \mcF^{\vee j} \}_{j=0}^n, \mcN^\vee, \nabla_\mcN^\vee, \omega^\dual)
\end{align}
forms a $(\mr{GO}_{2l+1}, \hslash)$-oper (resp., a $(\mr{GSp}_{2m}, \hslash)$-oper) on $U^\mr{log}/S^\mr{log}$.
Just as in the case of $(\mr{GL}_n, \hslash)$-opers discussed above, we have 
$(\mcF^{\heartsuit \dual}_{\sphericalangle})^\dual \cong  \mcF^{\heartsuit}_{\sphericalangle}$.

\bde  \label{pp122}
We shall refer to  $\msF^{\heartsuit \dual}$ (resp.,  $\msF^{\heartsuit \dual}_{\sphericalangle}$) defined above  as the {\bf dual} of
   $\msF^\heartsuit$ (resp., $ \msF^{\heartsuit}_{\sphericalangle}$).
 \index{$\mcF^\heartsuit$, $(\mr{GL}_n, \hslash)$-oper@--- $\msF^{\heartsuit \dual}$, dual of}
  \index{$\msF^{\heartsuit }_\sphericalangle$, $(\mr{GO}_{2l+1}, \hslash)$-oper, $(\mr{GSp}_{2m}, \hslash)$-oper@--- $\msF^{\heartsuit \dual}_{\sphericalangle}$, dual of}
  \ede

\begin{rema}
According to  ~\cite{BD2}, \S\,2.2,
taking the dual   of a $\mr{GL}_n$-oper   on a smooth curve $X$ over $\mbC$ may be identified with taking the {\it transpose} ${^t}\Dif$ of the corresponding  differential operator $\Dif$ (cf.  (\ref{Edr56}) and Theorem \ref{y0128} for the related correspondences).
We should note  that this identification was constructed by using pseudo-differential operators on $X/\mbC$, but the situation could be extended to more general cases. 
Here, recall that if $\Dif$ is a differential operator between line bundles $\mcL_1 \migi \mcL_2$   expressed locally as $\Dif = \sum_j a_j \cdot \partial^j$, where $\partial := \frac{d}{dx}$ for a local coordinate function $x$ on $X$, then its  transpose ${^t}\Dif$ is defined as the differential operator  $\Omega_{X/\mbC} \otimes \mcL_2^\vee \migi \Omega_{X/\mbC} \otimes \mcL_1^\vee$ expressed as  ${^t}\Dif := \sum_{j} (-\partial)^j \cdot a_j$.
\end{rema}

\subsection{The dual of an $(n, \hslash)$-theta characteristic}\label{QW1501}
Let $\bb := (\BB, \nabla_\bb)$ be an $(n, \hslash)$-theta characteristic of  $U^\mr{log}/S^\mr{log}$.
Write
\begin{align} \label{QW945663}
\BB^\dual := \Omega^{\otimes (n-1)} \otimes \BB^\vee.
\index{$\BB^\dual$}
\end{align}
Then, we have a composite isomorphism
\begin{align} \label{QW924}
\mcT^{\otimes \frac{n(n-1)}{2}} \otimes (\BB^\dual)^{\otimes n} \isom \mcT^{\otimes \frac{n(n-1)}{2}} \otimes \mcT^{\otimes - n(n-1)} \otimes (\BB^{\vee})^{\otimes n} \isom  (\mcT^{\otimes \frac{n(n-1)}{2}} \otimes \BB^{\otimes n})^\vee. \hspace{-5mm}
\end{align}
By  this composite, the  dual  $\nabla_\bb^\vee$ of $\nabla_\bb$ may be regarded as an $\hslash$-$S^\mr{log}$-connection on $\mcT^{\otimes \frac{n(n-1)}{2}} \otimes (\BB^\dual)^{\otimes n}$.
Thus, we obtain  an $(n, \hslash)$-theta characteristic 
\index{$\bb$, $(n, \hslash)$-theta characteristic@--- $\bb^\dual$}
\begin{align} \label{QW926}
\bb^\dual := (\BB^\dual, \nabla_\bb^\vee)
\end{align}
 of $U^\mr{log}/S^\mr{log}$, which we call  the {\bf dual} of $\bb$.
 It is immediate that $(\bb^\dual)^\dual \cong \bb$.
 Also, there is a canonical isomorphism of $(n, \hslash)$-theta characteristics
 \begin{align} \label{QW1303}
 \bb^\dual \otimes {^\dagger}\msN_\bb \isom \bb.
 \end{align}
 
\subsection{Self-dual $(\mr{GL}_n, \hslash, \bb)$-opers}\label{ppp0140p7}

Let $\nabla^\diamondsuit$ be a $(\mr{GL}_n, \hslash, \bb)$-oper on $U^\mr{log}/S^\mr{log}$.
Since 
\begin{align} \label{QW1010}
\DF_\BB^{\vee n-1}  = (\DF_\BB /\DF_\BB^1)^\vee = (\mcT^{\otimes (n-1)}\otimes \BB)^\vee = \BB^\dual,
\end{align}
we have an inclusion $\BB^\dual \migiincl \DF_\BB^\vee$.
The composite
\begin{align} \label{pp50}
\DF_{\BB^\dual} \left(=\mcD_{\hslash}^{<n} \otimes \BB^\dual \right)
& \migiincl  \mcD_{\hslash}^{< \infty}  \otimes \DF_\BB^\vee  \xrightarrow{{^\mcD}\hspace{-1.25mm}\nabla^{\diamondsuit\vee}}\DF_\BB^\vee
\end{align}
turns out to be an isomorphism, where the first arrow denotes the tensor product of  
the inclusions $\BB^\dual \migiincl \DF_\BB^\vee$, 
$\mcD_{\hslash}^{<n} \migiincl \mcD_{\hslash}^{<\infty}$.
Hence,  the dual $\nabla^{\diamondsuit \vee}$ of $\nabla^\diamondsuit$ corresponds, via this composite isomorphism, to  an $\hslash$-$S^\mr{log}$-connection 
\index{$\nabla^\diamondsuit$, $(\mr{GL}_n, \hslash, \bb)$-oper@--- $\nabla^{\diamondsuit \dual}$}
\begin{align} \label{QW945}
\nabla^{\diamondsuit \dual} :  \DF_{\BB^\dual} \migi \Omega  \otimes \DF_{\BB^\dual}
\end{align}
 on $\DF_{\BB^\dual}$.
 Moreover, it
  forms a $(\mr{GL}_n, \hslash, \bb^\dual)$-oper on $U^\mr{log}/S^\mr{log}$ whose underlying $(\mr{GL}_n, \hslash)$-oper   coincides  with the dual 
 of $\nabla^{\diamondsuit \Rightarrow \heartsuit}$ via (\ref{pp50}), i.e., $(\nabla^{\diamondsuit \dual})^{\Rightarrow \heartsuit} = (\nabla^{\diamondsuit \Rightarrow \heartsuit})^\dual$.
In particular, we obtain a $(\mr{GL}_n, \hslash, \bb)$-oper
\begin{align} \label{QW1305}
(\nabla^{\diamondsuit \dual})_{\otimes {^\dagger}\msN_\bb}
\end{align}
on $U^\mr{log}/S^\mr{log}$ under the identification  $\bb^\dual \otimes {^\dagger}\msN_\bb = \bb$ by (\ref{QW1303}).

\bde \label{QW1072}
Let $\nabla^\diamondsuit$ be a $(\mr{GL}_n, \hslash, \bb)$-oper on $U^\mr{log}/S^\mr{log}$.
Then, we shall say that $\nabla^\diamondsuit$ is {\bf self-dual}
\index{$\nabla^\diamondsuit$, $(\mr{GL}_n, \hslash, \bb)$-oper@--- self-dual}
if it is isomorphic to $(\nabla^{\diamondsuit \dual})_{\otimes {^\dagger}\msN_\bb}$.
\ede

\subsection{Characterization by self-duality} \label{QW1504}
If $\nabla^\diamondsuit_{\sphericalangle} :=(\nabla^\diamondsuit, \omega)$
 is a  $(\mr{GO}_{2l+1}, \hslash, \bb)$-oper  (resp., a $(\mr{GSp}_{2m}, \hslash, \bb)$-oper) on $U^\mr{log}/S^\mr{log}$,  then the pair   
 \index{$\nabla^\diamondsuit_{\sphericalangle}$, $(\mr{GO}_{2l+1}, \hslash, \bb)$-oper, $(\mr{GSp}_{2m}, \hslash, \bb)$-oper@--- $\nabla^{\diamondsuit\dual}_{\sphericalangle}$, dual of $\nabla^{\diamondsuit}_{\sphericalangle}$}
\begin{align} \label{QW957}
 \nabla^{\diamondsuit\dual}_{\sphericalangle} := (\nabla^{\diamondsuit \dual}, \omega^\dual)
\end{align}
 forms a $(\mr{GO}_{2l+1}, \hslash, \bb^\dual)$-oper (resp., a $(\mr{GSp}_{2m}, \hslash, \bb^\dual)$-oper) on $U^\mr{log}/S^\mr{log}$.
One verifies that $(\nabla_{\sphericalangle}^{\diamondsuit \dual})^\dual = \nabla_{\sphericalangle}^{\diamondsuit}$ and $(\nabla^{\diamondsuit \dual}_{\sphericalangle})^{\Rightarrow \heartsuit} \cong (\nabla_{\sphericalangle}^{\diamondsuit \Rightarrow \heartsuit})^\dual$.
Also. the isomorphism $\omega_{(-,\bullet)} : \DF_\BB \isom \DF_{\BB^\dual} \otimes {^\dagger}\mcN_\bb \left(= \DF_\BB^\vee \otimes {^\dagger}\mcN_\bb \right)$ becomes   an isomorphism 
\begin{align} \label{QW1066}
(\nabla^\diamondsuit_{\sphericalangle \Rightarrow \emptyset}) \isom 
((\nabla^{\diamondsuit}_{\sphericalangle \Rightarrow \emptyset})^\dual)_{\otimes {^\dagger}\msN_\bb} \left(= ((\nabla^{\diamondsuit \dual}_{\sphericalangle})_{\Rightarrow \emptyset})_{\otimes {^\dagger}\msN_\bb}\right)
\end{align}
 of  $(\mr{GL}_n, \hslash, \bb)$-opers.
In particular, the underlying $(\mr{GL}_n, \hslash, \bb)$-oper  $\nabla^\diamondsuit_{\sphericalangle \Rightarrow \emptyset}$ is self-dual.

Conversely, let 
$\nabla^\diamondsuit$ be a self-dual $(\mr{GL}_n, \hslash, \bb)$-oper on $U^\mr{log}/S^\mr{log}$.
We shall fix  an isomorphism $\eta : \nabla^{\diamondsuit} \isom (\nabla^{\diamondsuit\dual})_{\otimes {^\dagger}\msN_\bb}$.
This isomorphism  induces  a nondegenerate  $\mcO_U$-bilinear map 
\begin{align}
\omega_{\nabla^\diamondsuit}  : \DF_\BB^{\otimes 2} \migi {^\dagger}\mcN_\bb
\end{align}
on $\DF_\BB$, which 
 is compatible with the $\hslash$-$S^\mr{log}$-connections $\nabla^{\diamondsuit \otimes 2}$, ${^\dagger}\nabla_{\mcN_\bb}$.
Note that $\omega_{\nabla^\diamondsuit}$ depends on the choice of 
$\eta$,
 but
by Proposition \ref{QW1045}
it is uniquely determined up to composition with an automorphism of the $\hslash$-flat line bundle ${^\dagger}\msN_\bb$.
In particular, 
whether $\omega_{\nabla^\diamondsuit}$ is symmetric (resp., skew-symmetric) or not  depends only on the isomorphism class of $\nabla^\diamondsuit$.
If $\omega_{\nabla^\diamondsuit}$ is symmetric (resp., skew-symmetric), then the collection of data
\begin{align} \label{QW1056}
\nabla^\diamondsuit_{\emptyset  \Rightarrow \sphericalangle} :=(\nabla^{\diamondsuit}, \omega_{\nabla^\diamondsuit})
\end{align}
forms a $(\mr{GO}_{2l+1}, \hslash, \bb)$-oper (resp., a $(\mr{GSp}_{2m}, \hslash, \bb)$-oper) on $U^\mr{log}/S^\mr{log}$.
The isomorphism class of $\nabla^\diamondsuit_{\emptyset  \Rightarrow \sphericalangle}$ depends only on that of $\nabla^\diamondsuit$.
Therefore, the following assertion can be verified from the various definitions involved.

\bpr \label{QW1070}
The assignments
$\nabla^\diamondsuit_\sphericalangle \mapsto \nabla^\diamondsuit_{\sphericalangle \Rightarrow\emptyset}$ and $\nabla^\diamondsuit \mapsto \nabla^\diamondsuit_{\emptyset \Rightarrow \sphericalangle}$
define a bijective correspondence
\begin{equation} \label{equicat10}
\begin{pmatrix}
\text{the set of isomorphism  classes} \\
\text{of self-dual $(\mr{GL}_n, \hslash, \bb)$-opers $\nabla^\diamondsuit$} \\
\text{ on $U^\mr{log}/S^\mr{log}$ with $\omega_{\nabla^\diamondsuit}$} \\
\text{symmetric (resp., skew-symmetric)}
\end{pmatrix}
\isom
\begin{pmatrix}
\text{the set of isomorphism classes} \\
\text{of $(\mr{GO}_{2l+1}, \hslash, \bb)$-opers}  \\
\text{(resp., $(\mr{GSp}_{2m}, \hslash, \bb)$-opers)} \\
\text{on $U^\mr{log}/S^\mr{log}$}
\end{pmatrix}. \hspace{-3mm}
\end{equation}
Moreover, the formation of this correspondence commutes with pull-back to  \'{e}tale  schemes over  $U$, as well as  base-change to fs log schemes over $S^\mr{log}$.
\epr

\begin{exa}[$(\mr{GSp}_2, \hslash, \bb)$-opers] \label{Qw1329}
We shall  construct the equivalence between $(\mr{GL}_2, \hslash, \bb)$-opers and $(\mr{GSp}_2, \hslash, \bb)$-opers.
Let $\nabla^\diamondsuit$ be a $(\mr{GL}_2, \hslash, \bb)$-oper on $U^\mr{log}/S^\mr{log}$.
Notice that the $2$-fold tensor product $\DF_\BB^{\otimes 2}$ decomposes canonically as the direct sum of its symmetric and skew-symmetric parts, i.e., we have
\begin{align} \label{QW1331}
\DF_\BB^{\otimes 2} \isom \mbS^2_{\mcO_U} (\DF_\BB) \oplus \bigwedge^2 \DF_\BB.
\end{align}
Let us consider the composite
\begin{align}
\omega_{\nabla^\diamondsuit} :\DF_\BB^{\otimes 2} \xrightarrow{(\ref{QW1331})} \mbS^2_{\mcO_U} (\DF_\BB) \oplus \bigwedge^2 \DF_\BB \xrightarrow{\mr{pr}_2} \bigwedge^2 \DF_\BB  \left(=\mr{det}(\DF_\BB) \right)\xrightarrow{(\ref{fff0100})} {^\dagger}\mcN_\bb.
\end{align}
This is  compatible with  the $\hslash$-$S^\mr{log}$-connections $\nabla^{\diamondsuit \otimes 2}$,  ${^\dagger}\nabla_{\mcN_\bb}$, and  forms a nondegenerate skew-symmetric $\mcO_U$-bilinear map on $\DF_\BB$.
In particular, it induces an isomorphism $\nabla^\diamondsuit \isom (\nabla^{\diamondsuit \dual})_{\otimes {^\dagger}\msN_\bb}$, that is to say,  $\nabla^\diamondsuit$ is self-dual.
The resulting  assignment $\nabla^\diamondsuit \mapsto (\nabla^\diamondsuit, \omega_{\nabla^\diamondsuit})$ commutes with pull-back to \'{e}tale schemes over $U$ and hence defines a morphism $\mcO p^{\diamondsuit}_{2, \hslash, \bb} \migi \mcO p_{n, \hslash, \bb}^{\diamondsuit, \mr{Sp}}$, which gives  an inverse to the natural morphism $\mcO p_{2, \hslash, \bb}^{\diamondsuit, \mr{Sp}} \migi \mcO p_{2, \hslash, \bb}^{\diamondsuit}$.

\end{exa}

\section{Isomorphisms of moduli spaces induced by duality}\label{p0140p3}

\subsection{The moduli spaces and comparison results}\label{QW613}
Let $(g,r)$ be a pair of nonnegative integers with $2g-2+r>0$.
Also, let $S$ be a $k$-scheme and $\msX := (f: X \migi S, \{ \sigma_i \}_{i=1}^r)$  an $r$-pointed  stable curve of genus $g$ over $S$.
For convenience, we write
\index{$\mfO \mfp_{n, \hslash (, \rho)}^{\spadesuit}$, $\mfO \mfp_{n, \hslash (, \rho)}^{\spadesuit, \mr{O}}$, $\mfO \mfp_{n, \hslash (, \rho)}^{\spadesuit, \mr{Sp}}$}
\begin{align}
\mfO \mfp_{n, \hslash}^{\spadesuit, \mr{O}} := \mfO \mfp^{}_{\mfs \mfo_{2l+1}, \hslash, \msX} \ \left(\text{resp.,} \  \mfO \mfp_{n, \hslash}^{\spadesuit, \mr{Sp}} := \mfO \mfp^{}_{\mfs \mfp_{2m}, \hslash, \msX} \right).
\end{align}
Also,  given an $(n, \hslash)$-theta characteristic $\bb$, we write
\index{$\mfO \mfp_{n, \hslash, \bb (, \rho)}^\diamondsuit$, $\mfO \mfp^{\diamondsuit, \mr{O}}_{n, \hslash}$, $\mfO \mfp^{\diamondsuit, \mr{Sp}}_{n, \hslash}$}
\begin{align} \label{pp21}
\overline{\mfO} \mfp_{n, \hslash}^{\heartsuit, \mr{O}}  \  \left(\text{resp.,} \ \mfO \mfp_{n, \hslash, \bb}^{\diamondsuit, \mr{O}}; \text{resp.,} \ \overline{\mfO} \mfp_{n, \hslash}^{\heartsuit, \mr{Sp}}; \text{resp.,}  \ \mfO \mfp_{n, \hslash, \bb}^{\diamondsuit, \mr{Sp}}\right) : \mfS \mfc \mfh_{/S} \migi \mfS \mfe \mft
\index{$\overline{\mfO} \mfp_{n, \hslash}^{\heartsuit, \BL}$},
\end{align}
for the $\mfS \mfe \mft$-valued contravariant functor on $\mfS \mfc \mfh_{/S}$ which,
to any $S$-scheme $s:S' \migi S$, assigns the set of global sections of 
$\overline{\mcO}p_{n, \hslash}^{\heartsuit, \mr{O}}$  (resp., $\mcO p^{\diamondsuit, \mr{O}}_{h, \hslash, \bb}$; resp., $\overline{\mcO}p^{\heartsuit, \mr{Sp}}_{n, \hslash}$; resp., $\mcO p^{\diamondsuit, \mr{Sp}}_{n, \hslash, \bb}$) defined for the log curves $(S' \times_S X^\mr{log})/(S' \times_S S^\mr{log})$.

\bt \label{QW1199}
Let $\BL \in \{ \mr{O}, \mr{Sp}\}$.
Then, there exist canonical isomorphisms of functors
\begin{align} \label{QW1200}
\Lambda_{\bb, \msX}^{\diamondsuit \Rightarrow \heartsuit, \BL} : 
\mfO \mfp_{n, \hslash, \bb}^{\diamondsuit, \BL} \isom
\overline{\mfO} \mfp_{n, \hslash}^{\heartsuit, \BL},
\hspace{10mm} 
\Lambda_{\msX}^{\heartsuit \Rightarrow \spadesuit, \BL} : 
\overline{\mfO} \mfp_{n, \hslash}^{\heartsuit, \BL} \isom
\mfO \mfp_{n, \hslash}^{\spadesuit, \BL}
\index{$\Lambda_{\bb (, \rho)}^{\diamondsuit \Rightarrow \heartsuit}$, $\Lambda_\bb^{\diamondsuit \Rightarrow \heartsuit, \BL}$, $\Lambda_{\bb, \msX}^{\diamondsuit \Rightarrow \heartsuit, \BL}$}
\index{$\Lambda^{\heartsuit \Rightarrow \spadesuit}_{(\rho)}$, $\Lambda^{\heartsuit \Rightarrow \spadesuit, \BL}$, $\Lambda^{\heartsuit \Rightarrow \spadesuit, \BL}_\msX$}
\end{align}
which make the following diagram commute:
\begin{align} \label{QW1201}
\vcenter{\xymatrix@C=56pt@R=36pt{
 \mfO \mfp_{n, \hslash, \bb}^{\diamondsuit, \BL} 
 \ar[r]^-{\Lambda_{\bb, \msX}^{\diamondsuit \Rightarrow \heartsuit, \BL}}_-{\sim} \ar[d]^-{\nabla^{\diamondsuit}_{\sphericalangle} \mapsto \nabla^{\diamondsuit}_{\sphericalangle \Rightarrow \emptyset}} & \overline{\mfO} \mfp_{n, \hslash}^{\heartsuit, \BL} 
  \ar[r]^-{\Lambda_{\msX}^{\heartsuit \Rightarrow \spadesuit, \BL}}_-{\sim} \ar[d]^-{[\msF^\heartsuit_{\sphericalangle}] \mapsto  [\msF^\heartsuit_{\sphericalangle \Rightarrow \emptyset}]}& \mfO \mfp_{n, \hslash}^{\spadesuit, \BL} \ar[d]^-{\mcE^\spadesuit_{\sphericalangle} \mapsto \mcE^\spadesuit_{\sphericalangle \Rightarrow \emptyset}}\\
 \mfO \mfp_{n, \hslash, \bb}^{\diamondsuit} \ar[r]_-{\Lambda_{\bb, \msX}^{\diamondsuit \Rightarrow \heartsuit}}^-{\sim} &  \overline{\mfO} \mfp_{n, \hslash}^{\heartsuit}  \ar[r]^-{\sim}_-{\Lambda_{\msX}^{\heartsuit \Rightarrow \spadesuit}}& \mfO \mfp_{n, \hslash}^{\spadesuit}.
}}
\end{align}
Moreover, $\mfO \mfp_{n, \hslash, \bb}^{\diamondsuit, \BL}$, $\overline{\mfO} \mfp_{n, \hslash}^{\heartsuit, \BL}$, and $\mfO \mfp_{n, \hslash}^{\spadesuit, \BL}$ may be represented by $S$-schemes and 
all the vertical arrows in this diagram are closed immersions.
\et
  \begin{proof}
  The former assertion follows from Propositions \ref{QW910} and \ref{QW911}.
  To prove the latter assertion, it suffices, from the former assertion,  to prove the claim that
  the left-hand vertical arrow $\mfO \mfp_{n, \hslash, \bb}^{\diamondsuit, \BL} \migi \mfO \mfp_{n, \hslash, \bb}^{\diamondsuit}$ in (\ref{QW1201}) is a closed immersion. 
  In what follows, we will only consider the case of ``\,$? = \mr{Sp}$'' because the case of  ``\,$? = \mr{O}$'' is entirely similar.
The case of $m=1$ was already proved in   Example \ref{Qw1329},  so let us  assume that $m >1$, i.e., $n >2$.
  Write $\delta_\bb$ for the automorphism of $\mfO \mfp_{n, \hslash, \bb}^\diamondsuit$ given by $\nabla^\diamondsuit \mapsto (\nabla^{\diamondsuit \dual})_{\otimes {^\dagger}\msN_\bb}$.
  Since $\mfO \mfp^\diamondsuit_{n, \hslash, \bb}$ is separated over $S$ (cf. Theorem \ref{y0129}), 
  the equalizer
  \begin{align}
  \dot{S} := \mr{Eq}(\delta_\bb, \mr{id}) \left(= \varprojlim \left\{ \mfO \mfp^\diamondsuit_{n, \hslash, \bb}  
  \overset{\delta_\bb}{\underset{\mr{id}}{\rightrightarrows}} \mfO \mfp^\diamondsuit_{n, \hslash, \bb}
   \right\}  \right)
  \end{align}
  of $\delta_\bb$ and the identity morphism forms a closed subschme of $\mfO \mfp^\diamondsuit_{n, \hslash, \bb}$.
  By definition, $ \dot{S}$ represents the subfunctor of   $\mfO \mfp^\diamondsuit_{n, \hslash, \bb}$  classifying self-dual $(\mr{GL}_n, \hslash, \bb)$-opers.
  Denote by $\dot{f} : \dot{X} \migi \dot{S}$ the base-change of $f:X \migi S$ to the $S$-scheme $\dot{S}$ and by $\dot{\nabla}^\diamondsuit$ the $(\mr{GL}_n, \hslash, \bb)$-oper on $\dot{X}^\mr{log}/\dot{S}^\mr{log}$ classified by 
  $\dot{S} \left(\subseteq  \mfO \mfp^\diamondsuit_{n, \hslash, \bb} \right)$.
 In the following, the dot ``$\dot{(-)}$'' will be used to denote base-change by  the projection $\dot{S} \migi S$.
  Let us  fix an isomorphism $\eta : \dot{\nabla}^\diamondsuit \isom (\dot{\nabla}^{\diamondsuit \dual})_{\otimes {^\dagger}\dot{\msN}_\bb}$ and denote by 
  $\omega : {^\dagger}\dot{\mcF}_\BB^{\otimes 2} \migi {^\dagger}\dot{\mcN}_\bb$ the induced $\mcO_{\dot{X}}$-bilinear map on ${^\dagger}\dot{\mcF}_\BB$.
  The sum $\overline{\omega} := \omega + \omega \circ \mr{sw}$, where $\mr{sw}$ denotes the automorphism of ${^\dagger}\dot{\mcF}_\BB^{\otimes 2}$ given by $a \otimes b \mapsto b\otimes a$,   determine a global section of ${\dot{f}}_*( {^\dagger}\dot{\mcN}_\bb \otimes {^\dagger}\dot{\mcF}_\BB^{\vee \otimes 2})$.
  Here, note that 
  ${^\dagger}\dot{\mcN}_\bb \otimes {^\dagger}\dot{\mcF}_\BB^{\vee \otimes 2}$ admits 
   a filtration arising  from the filtration $\{ {^\dagger}\dot{\mcF}^j_\BB \}_{j=0}^n$ each of whose  subquotients  is a   line bundle isomorphic to 
   $ {^\dagger}\dot{\mcN}_\bb \otimes ({^\dagger}\dot{\mcF}_\BB^{\vee j}/{^\dagger}\dot{\mcF}_\BB^{\vee j+1}) \otimes ({^\dagger}\dot{\mcF}_\BB^{\vee j'}/{^\dagger}\dot{\mcF}_\BB^{\vee j'+1})$ for some $j, j'$.
  Since ${^\dagger}\dot{\mcF}^j_\BB/ {^\dagger}\dot{\mcF}^{j+1}_\BB \cong \dot{\mcT}^{\otimes n-1-j}\otimes\dot{\BB}$ ($j=0, \cdots, n-1$), each subquotient of this filtration is isomorphic to $\dot{\Omega}^{\otimes (n+j-1)}$ for some integer $j$ with $0\leq j \leq 2n-2$, so vanishes after applying the functor $\mbR^1 \dot{f}_* (-)$ (because of the assumption $n>2$).
  Hence, $\mbR^1{\dot{f}}_*( {^\dagger}\dot{\mcN}_\bb \otimes {^\dagger}\dot{\mcF}_\BB^{\vee\otimes 2}) =0$ and ${\dot{f}}_*( {^\dagger}\dot{\mcN}_\bb \otimes {^\dagger}\dot{\mcF}_\BB^{\vee \otimes 2})$ forms a vector bundle on $S$ (cf. ~\cite{Har}, Chap.\,III, Theorem 12.11, (b)).
  Let us define $\ddot{S}$  to be the equalizer
   \begin{align}
  \ddot{S} := \mr{Eq}(\overline{\omega}, 0) \left(= \varprojlim \left\{ \dot{S}  
  \overset{\overline{\omega}}{\underset{0}{\rightrightarrows}} \mbV ({\dot{f}}_*( {^\dagger}\dot{\mcN}_\bb \otimes {^\dagger}\dot{\mcF}_\BB^{\vee \otimes 2}))
   \right\}  \right).
  \end{align}
of the two sections of ${\dot{f}}_*( {^\dagger}\dot{\mcN}_\bb \otimes {^\dagger}\dot{\mcF}_\BB^{\vee \otimes 2})$ determined by $\overline{\omega}$ and the zero element $0$.
It  forms a closed subscheme of $\dot{S}$ (and hence of $S$) and
does not depend on the choice of $\eta$ (cf. Proposition \ref{QW1045}). 
By definition,  $\ddot{S}$ represents the subfunctor of $\mfO \mfp^{\diamondsuit}_{n, \hslash, \bb}$  classifying  self-dual $(\mr{GL}_n, \hslash, \bb)$-opers $\nabla^\diamondsuit$ with $\omega_{\nabla^\diamondsuit}$ skew-symmetric.
On the other hand,  it follows from Proposition \ref{QW1070} that  the $S$-scheme $\ddot{S}$ is isomorphic to $\mfO \mfp^{\diamondsuit, \mr{Sp}}_{n, \hslash, \bb}$.
Hence, the natural morphism $\mfO \mfp^{\diamondsuit, \mr{Sp}}_{n, \hslash, \bb} \migi \mfO \mfp_{n, \hslash, \bb}^{\diamondsuit}$ turns out to be a closed immersion.
This completes the proof of the assertion.
 \end{proof} 

\subsection{Duality of   opers on a pointed stable curve}\label{p0140p8}

Next, we describe assertions concerning automorphisms of  moduli spaces   arising from duality.

\bt  \label{yp0128}
\begin{itemize}
\item[(i)]
The assignments $\nabla^\diamondsuit \mapsto (\nabla^{\diamondsuit \dual})_{\otimes {^\dagger}\msN_\bb}$ and  $[\msF^\heartsuit] \mapsto [\msF^{\heartsuit \dual}]$ define  automorphisms
 \index{$\Dual_{\bb (, \rho)}^{\, \diamondsuit}$, $\Dual_{(\rho)}^{\,\heartsuit}$, $\Dual_{(\rho)}^{\, \spadesuit}$, ${^p}\Dual_{(\rho)}^{\, \spadesuit}$}
\begin{align} \label{QG3001}
\Dual_\bb^{\, \diamondsuit}
 : \mfO \mfp_{n, \hslash, \bb}^\diamondsuit \isom \mfO \mfp_{n, \hslash, \bb}^\diamondsuit
\hspace{3mm} \text{and} \hspace{3mm} 
\Dual^{\,\heartsuit}
 : \overline{\mfO} \mfp_{n, \hslash}^\heartsuit \isom \overline{\mfO} \mfp_{n, \hslash}^\heartsuit
\end{align}
respectively,  and the following square diagram is commutative:
\begin{align} \label{QW1341}
\vcenter{\xymatrix@C=46pt@R=36pt{
 \mfO \mfp_{n, \hslash, \bb}^\diamondsuit 
 \ar[r]_-{\sim}^-{\Lambda^{\diamondsuit \Rightarrow \heartsuit}_\bb} \ar[d]^-{\wr}_-{\scalebox{0.7}{$\Dual$}^{\, \diamondsuit}_\bb} &
 \overline{\mfO} \mfp_{n, \hslash}^\heartsuit  \ar[d]_-{\wr}^-{\scalebox{0.7}{$\Dual$}^{\, \heartsuit}}\\
\mfO \mfp_{n, \hslash, \bb}^\diamondsuit \ar[r]^-{\sim}_-{\Lambda^{\diamondsuit \Rightarrow \heartsuit}_{\bb}}& \overline{\mfO} \mfp_{n, \hslash}^\heartsuit.
}}
\end{align}
Moreover,  the equalities  $\Dual^{\, \diamondsuit}_\bb \circ \Dual^{\, \diamondsuit}_\bb  = \mr{id}$ and 
$\Dual^{\, \heartsuit} \circ \Dual^{\, \heartsuit} = \mr{id}$ hold.
\item[(ii)]
Let $\BL \in \{ \mr{O}, \mr{Sp}\}$.
Then, the following diagrams are commutative:
\begin{align} \label{QW9001}
\vcenter{\xymatrix@C=10pt@R=36pt{
& \mfO \mfp^{\diamondsuit, \BL}_{n,\hslash, \bb}\ar[dl]_-{\nabla_\sphericalangle^\nabla \mapsto \nabla_{\sphericalangle \Rightarrow \emptyset}^\nabla} \ar[dr]^-{\nabla_\sphericalangle^\nabla \mapsto \nabla_{\sphericalangle \Rightarrow \emptyset}^\nabla}&\\
\mfO \mfp_{n, \hslash, \bb}^{\diamondsuit}\ar[rr]^-{\sim}_-{\scalebox{0.7}{$\Dual$}_\bb^{\, \heartsuit}}&& \mfO \mfp_{n, \hslash, \bb}^{\diamondsuit},
}}
\hspace{8mm}
\vcenter{\xymatrix@C=10pt@R=36pt{
& \overline{\mfO} \mfp_{n, \hslash}^{\heartsuit, \BL}\ar[dl]_-{[\msF^\heartsuit_\sphericalangle] \mapsto [\msF^\heartsuit_{\sphericalangle\Rightarrow \emptyset}]} \ar[dr]^-{[\msF^\heartsuit_\sphericalangle] \mapsto [\msF^\heartsuit_{\sphericalangle\Rightarrow \emptyset}]}&\\
\overline{\mfO} \mfp_{n, \hslash}^{\heartsuit} \ar[rr]^-{\sim}_-{\scalebox{0.7}{$\Dual$}^{\, \heartsuit}}&& \overline{\mfO} \mfp_{n, \hslash}^{\heartsuit}.
}}
\end{align}
\end{itemize}
\et
\begin{proof}
The assertions follows from Proposition \ref{QW1070} and the equality $(\nabla^{\diamondsuit \dual})^{\Rightarrow \heartsuit} = (\nabla^{\diamondsuit \Rightarrow \heartsuit})^\dual$ obtained  for each $(\mr{GL}_n, \hslash, \bb)$-oper $\nabla^\diamondsuit$.
\end{proof}

By the above assertion, we obtain an automorphism
 \index{$\Dual_{\bb (, \rho)}^{\, \diamondsuit}$, $\Dual_{(\rho)}^{\,\heartsuit}$, $\Dual_{(\rho)}^{\, \spadesuit}$, ${^p}\Dual_{(\rho)}^{\, \spadesuit}$}
\begin{align}
\Dual^{\, \spadesuit}
  := \Lambda^{\heartsuit \Rightarrow \spadesuit} \circ  
  \Dual^{\, \heartsuit}
   \circ (\Lambda^{\heartsuit \Rightarrow \spadesuit})^{-1} : \mfO \mfp_{n, \hslash}^{\spadesuit} \isom \mfO \mfp_{n, \hslash}^{\spadesuit}
\end{align}
of $\mfO \mfp_{n, \hslash}^{\spadesuit}$ with $\Dual^{\, \spadesuit} \circ \Dual^{\, \spadesuit} = \mr{id}$.
Moreover, it   restricts to the identity morphism of $\mfO \mfp_{n, \hslash}^{\spadesuit, \BL}$ for $\BL \in \{ \mr{O}, \mr{Sp}\}$.

\bde  \label{pp79}
Given an $(\mfs \mfl_n, \hslash)$-oper $\msE^\spadesuit$, we shall write
$\msE^{\spadesuit \dual}$ for the $(\mfs \mfl_n, \hslash)$-oper defined as the image of $\msE^\spadesuit$ via $\Dual^{\, \spadesuit}$.
We refer to $\msE^{\spadesuit \dual}$ as the {\bf dual} 
\index{$\msE^\spadesuit$, $(\mfg, \hslash)$-oper@--- $\msE^{\spadesuit \dual}$, dual of}
of  $\msE^\spadesuit$.
  \ede


\bt[cf. Theorem \ref{ypp017}]  \label{yp3128}
Let $\rho \in \mfc_{\mfs \mfl_n}^{\times r} (S)$ (where $\mfc_{\mfs \mfl_n}^{\times r} (S) := \{ \emptyset \}$ and $\rho := \{ \emptyset \}$ if $r=0$).
Then, the following assertions hold:
\begin{itemize}
\item[(i)]
The automorphism $\Dual^{\, \diamondsuit}_\bb$ (resp., $\Dual^{\, \heartsuit}$; resp., $\Dual^{\, \spadesuit}$)
restricts to an isomorphism
\begin{align} \label{QW1360}
\Dual^{\, \diamondsuit}_{\bb, \rho}
 : \mfO \mfp^{\diamondsuit}_{n, \hslash, \bb, \rho} \isom  \mfO \mfp^{\diamondsuit}_{n, \hslash, \bb, (-\hslash)\star\rho} \hspace{35mm}\\
\left(\text{resp.,} \  
\Dual_\rho^{\, \heartsuit} : \overline{\mfO} \mfp^{\heartsuit}_{n, \hslash, \rho} \isom  \overline{\mfO} \mfp^{\heartsuit}_{n, \hslash, (-\hslash)\star\rho};\text{resp.,} \
\Dual_\rho^{\, \spadesuit}   : \mfO \mfp^{\spadesuit}_{n, \hslash,  \rho} \isom  \mfO \mfp^{\spadesuit}_{n, \hslash, (-\hslash)\star\rho} \right)\notag
 \index{$\Dual_{\bb (, \rho)}^{\, \diamondsuit}$, $\Dual_{(\rho)}^{\,\heartsuit}$, $\Dual_{\rho}^{\, \spadesuit}$}
\end{align}
\item[(ii)]
 If $\rho$ lies in $\mfc_{\mfs \mfo_{2l+1}}^{\times r} (S)$ (cf. (\ref{pp88})),
 then the isomorphism $\Dual^{\, \diamondsuit}_{\bb, \rho}$ (resp., $\Dual_\rho^{\, \heartsuit}$; resp., $\Dual_\rho^{\, \spadesuit}$)
 restricts to
 the identity morphism of
$\mfO \mfp^{\diamondsuit, \mr{O}}_{n, \hslash, \bb, \rho} := \mfO \mfp^{\diamondsuit, \mr{O}}_{n, \hslash, \bb} \cap \mfO \mfp^{\diamondsuit}_{n, \hslash, \bb, \rho}$ (resp., $\overline{\mfO} \mfp^{\heartsuit, \mr{O}}_{n, \hslash, \rho} := \overline{\mfO} \mfp^{\heartsuit, \mr{O}}_{n, \hslash} \cap \overline{\mfO} \mfp^{\heartsuit}_{n, \hslash, \rho}$; resp., $\mfO \mfp^{\spadesuit, \mr{O}}_{n, \hslash, \rho} := \mfO \mfp^{\spadesuit, \mr{O}}_{n, \hslash}\cap \mfO \mfp^{\spadesuit}_{n, \hslash,  \rho}$). 
Moreover, the same assertion with $\mr{GO}_{2l+1}$ and $\mfs \mfo_{2l+1}$ replaced by $\mr{GSp}_{2m}$ and $\mfs \mfp_{2m}$ respectively holds.
\end{itemize}
\et
\begin{proof}
The assertions  follow from (\ref{QW2001}) and Theorem  \ref{yp0128}.
\end{proof}

\subsection{Positive characteristic}\label{p0140p9}

We end this chapter by proposing  two assertions
 concerning the positive characteristic situation.
These assertions can be deduced   from results obtained so far, via restricting the relevant  moduli schemes of opers  to their closed subschemes.

Let us consider the situation that {\it $k$ is of characteristic $p>0$, i.e., 
the inequality  $p > 2  n$ holds}.
For simplicity, we will  make the statements only for the case of ``$\spadesuit$''.

\bpr  \label{QW1392}
Let $\BL \in \{ \mfs \mfo_{2l+1}, \mfs \mfp_{2m} \}$.  
Also, let us fix an 
element 
  $\rho \in \mfc_{\BL}^{\times r} (S)$,  
     considered as an element of $\mfc_{\mfs \mfl_n}^{\times r} (S)$ via the closed immersion  $\mfc_{\BL} \migiincl \mfc_{\mfs \mfl_n}$ (cf. (\ref{pp88})).
Then,  the closed immersion $\mfO \mfp_{\BL, \hslash, \Box, \msX}^{} \migiincl \mfO \mfp_{\mfs \mfl_n, \hslash, \Box,  \msX}$ given by  $\msE^\spadesuit_{\sphericalangle} \mapsto \msE_{\sphericalangle \Rightarrow \emptyset}^\spadesuit$
restricts to  closed immersions
\begin{align} \label{pp74}
\mfO \mfp_{\BL, \hslash, \Box, \msX}^{^\mr{Zzz...}} \migiincl \mfO \mfp_{\mfs \mfl_{n}, \hslash, \Box, \msX}^{^\mr{Zzz...}} \hspace{3mm}\text{and} \hspace{3mm}  \mfO \mfp_{\BL, \hslash, \Box, \msX}^{^{p \text{-} \mr{nilp}}} \migiincl \mfO \mfp_{\mfs \mfl_{n}, \hslash,  \Box, \msX}^{^{p \text{-} \mr{nilp}}}
\end{align}
over $S$, where $\Box$ denotes either the absence or presence of ``$\rho$''.
\epr
\begin{proof}
The assertions for  ``$\mr{Zzz...}$'' (resp., ``$p$-nilp'') follows from  (\ref{QS61}) (cf. (\ref{QH23})) in the cases where the pair $(\mbG, \mbG')$ is taken to be $(\mr{PGO}_{2l+1}, \mr{PGL}_n)$ and  $(\mr{PGSp}_{2m}, \mr{PGL}_n)$.
\end{proof}

Moreover, since the dual   of any dormant $(\mr{GL}_n, \hslash)$-oper is  dormant (cf. the resp'd portion of (\ref{QW210})),  the following  theorem  holds.

\bt[cf. Theorem \ref{ypp017}] 
 \label{ppp012}
\begin{itemize}
\item[(i)]
The automorphism $\Dual^{\, \spadesuit}$ restricts to $S$-automorphisms
\begin{align} \label{pdrtp81}
\mfO \mfp_{\mfs \mfl_n, \hslash, \msX}^{^\mr{Zzz...}}  \isom  \mfO \mfp_{\mfs \mfl_n, \hslash, \msX}^{^\mr{Zzz...}} 
\hspace{3mm}\text{and}
 \hspace{3mm}  
 \mfO \mfp_{\mfs \mfl_n, \hslash, \msX}^{^{p \text{-} \mr{nilp}}}  \isom  \mfO \mfp_{\mfs \mfl_n, \hslash, \msX}^{^{p \text{-} \mr{nilp}}}.
\end{align}
Moreover,  for an element $\rho$ of $\mfc_{\mfs \mfl_n}^{\times r} (S)$,   these automorphisms restrict to isomorphisms
\begin{align} \label{pp81}
\mfO \mfp_{\mfs \mfl_n, \hslash, \rho, \msX}^{^\mr{Zzz...}}  \isom  \mfO \mfp_{\mfs \mfl_n, \hslash,  (-\hslash)\star \rho, \msX}^{^\mr{Zzz...}}
\hspace{3mm}
\text{and} \hspace{3mm}  
 \mfO \mfp_{\mfs \mfl_n, \hslash, \rho, \msX}^{^{p \text{-} \mr{nilp}}}  \isom  \mfO \mfp_{\mfs \mfl_n, \hslash,  (-\hslash)\star \rho, \msX}^{^{p \text{-} \mr{nilp}}}
\end{align}
respectively.
\item[(ii)]
Let $\BL \in \{ \mfs \mfo_{2l+1}, \mfs \mfp_{2m} \}$, and let us fix an element $\rho$ of $\mfc_{\BL}^{\times r}(S)$, considered as an element of $\mfc_{\mfs \mfl_n}^{\times r}(S)$.
Then, the automorphisms in (\ref{pdrtp81}) (resp., (\ref{pp81})) restrict to
the identity morphisms of 
$\mfO \mfp^\ZZZ_{\BL, \hslash, \msX}$ (resp., $\mfO \mfp^\ZZZ_{\BL, \hslash, \rho, \msX}$) and  $\mfO \mfp^\NILP_{\BL, \hslash, \msX}$ (resp., $\mfO \mfp^\NILP_{\BL, \hslash, \rho, \msX}$) respectively.
\end{itemize}
\et
\begin{proof}
The assertions follows from Theorem \ref{yp3128}, (i) and (ii).
\end{proof}

\vspace{10mm}

\chapter{Local deformation  of opers} \label{y0139} \vspace{3mm}

The main purpose of this chapter is to study the local theory of (dormant) $(\mfg, \hslash)$-opers.
We first review the classical deformation theory of pointed stable curves and flat vector bundles in the context of logarithmic geometry.
After that, we describe the deformation space of a $(\mfg, \hslash)$-oper using de Rham cohomology of the adjoint bundle associated to that oper.
A key observation concerning opers  in characteristic $p>0$ is that the Cartier operator may be interpreted as the deformation of $p$-curvature. 
By taking this into account, 
we achieve a detailed understanding of  the tangent space of the moduli stack classifying  dormant $(\mfg, \hslash)$-opers.
One of the  results in this chapter is to give a criterion for the \'{e}taleness of that stack, which plays a crucial  role in  the proof of Joshi's conjecture.

\begin{notation}
Let $k$
\index{$k$, field} 
 be a field and $\mbG$
\index{$\mbG$, algebraic group}
 a connected smooth affine algebraic group over $k$ with Lie algebra $\mfg$.
Also, let $(g, r)$ be a pair of nonnegative integers with $2g-2+r>0$,  $S$  a scheme over $k$, and $\msX := (f : X \migi S, \{ \sigma_i \}_{i=1}^r)$ an $r$-pointed stable curve over $S$ of genus $g$.
For simplicity, we write $\Omega:= \Omega_{X^\mr{log}/S^\mr{log}}$, $\mcT := \mcT_{X^\mr{log}/S^\mr{log}}$.

Next, let $k_\epsilon := k[\epsilon]/(\epsilon^2)$. 
For each $k$-scheme $Y$, we shall write $Y_\epsilon := \mr{Spec}(k_\epsilon) \times_{\mr{Spec}(k)} Y$; denote by  
 $\mr{pr}_Y : Y_\epsilon \migi Y$ the natural projection and by $\mr{in}_Y : Y \migi Y_\epsilon$ the closed immersion  induced by reduction modulo $\epsilon$.
Finally,  we write   $\msX_\epsilon := S_\epsilon \times_S \msX$, i.e., the base-change of $\msX$ to  $S_\epsilon$ (cf. (\ref{QQ02})).

\end{notation}

\section{Deformation spaces of a curve and an $\hslash$-flat $\mbG$-bundle}\label{QW6000}

\subsection{The Kodaira-Spencer map of a curve} \label{QW6002}
We begin  by recalling the Kodaira-Spencer map associated to a family of pointed stable curves $\msX$.
Since $f^\mr{log} : X^\mr{log} \migi S^\mr{log}$ is log smooth,  its  differential  yields the following short exact sequence:
\begin{align}\label{QH26}
0 \longmigi \mcT \left(=\mcT_{X^\mr{log}/S^\mr{log}}\right) \longmigi \mcT_{X^\mr{log}/k} \migi f^*(\mcT_{S^\mr{log}/k}) \longmigi 0.
\end{align}
Let us consider  the composite
\begin{align} \label{QW6011}
\mr{KS}'_{\msX} : \mcT_{S^\mr{log}/k} \migi  f_*(f^*(\mcT_{S^\mr{log}/k})) 
 \migi  \mbR^1 f_*(\mcT), 
 \index{$\mr{KS}'_{\msX}$, $\mr{KS}_{\msX}$, Kodaira-Spencer map associated to $\msX$}
\end{align}
where 
the first arrow denotes the morphism corresponding to the identity morphism
of $f^*(\mcT_{S^\mr{log}/k})$ via the adjunction relation ``$f^*(-) \dashv f_*(-)$'' and
the second arrow 
denotes  the connecting morphism   in the long exact sequence of higher direct images associated to (\ref{QH26}).
This  morphism also can be constructed in the case of $\msX = \msC_{g,r}$, i.e.,  the universal family of pointed curves   over $\overline{\mfM}_{g,r}$.
It is well-known that the resulting morphism $\mr{KS}'_{\msC_{g,r}} : \mcT_{\overline{\mfM}_{g,r}^\mr{log}} \migi \mbR^1 f_{\mr{univ}*}(\mcT_{\mfC^\mr{log}_{g,r}/\overline{\mfM}_{g,r}})$ is an isomorphism.
Now, let 
\index{$s^\mr{log}_\mfm$, $s^{\mr{log}}_{\mfo \mfp}$, ${^c}s^{\mr{log}}_{\mfo \mfp}$, ${^p}s^{\mr{log}}_{\mfo \mfp}$}
\begin{align} \label{QW6010}
s^\mr{log}_{\mfm}: S^\mr{log} \migi \overline{\mfM}_{g,r}^\mr{log}
\end{align}
 denote the strict morphism whose underlying $S$-rational point of $\overline{\mfM}_{g,r}$ classifies $\msX$.
The pull-back of $\mr{KS}'_{\msC_{g,r}}$ via $s_\mfm$ becomes  an isomorphism
\begin{align}
\mr{KS}_\msX : s_\mfm^*(\mcT_{\overline{\mfM}^\mr{log}_{g,r}})  \isom \mbR^1 f_*(\mcT)
 \index{$\mr{KS}'_{\msX}$, $\mr{KS}_{\msX}$, Kodaira-Spencer map associated to $\msX$}
\end{align}
 under the identification   $s^*_\mfm (\mbR^1 f_{\mr{univ}*}(\mcT_{\mfC^\mr{log}_{g,r}/\overline{\mfM}_{g,r}})) = \mbR^1 f_*(\mcT)$ given by the base-change map.
The morphism $\mr{KS}'_\msX$ coincides with the composite of $\mr{KS}_\msX$ and the differential $\mcT_{S^\mr{log}/k} \migi s^*_\mfm (\mcT_{\overline{\mfM}^\mr{log}_{g,r}})$ of $s^\mr{log}_\mfm$.
Hence, 
$\mr{KS}'_\msX$ is an isomorphism  if  $s^\mr{log}$ is  log \'{e}tale (or equivalently, $s$ is \'{e}tale),  and vice versa when $S^\mr{log}$ is log smooth over $k$ (cf. ~\cite{KaKa}, \S\,3, Proposition (3.12)). 
We shall refer to $\mr{KS}_\msX$ (or $\mr{KS}'_\msX$) as the {\bf Kodaira-Spencer map} associated to $\msX$. 

\subsection{Deformation of curves}\label{QW6801}
Let us study the local description of the  Kodaira-Spencer map.
By a {\bf first-order deformation} of $\msX$, we shall mean  an $r$-pointed stable curve $\msX_{\epsilon'}$  of genus $g$ over $S_\epsilon$ together with a morphism $\msX \migi \msX_{\epsilon'}$ over the  inclusion $\mr{in}_S : S  \migi S_\epsilon$.
(In particular, we have $\msX \cong S \times_{S_\epsilon} \msX_{\epsilon'}$.)
The notion of an isomorphism between  two first-order deformations of $\msX$  can be defined in an evident way. (We leave to the reader the precise definition.)
By the universal property of $\overline{\mfM}_{g,r}$,
the set of isomorphism classes of first-order deformations of $\msX$
corresponds bijectively to the set of isomorphism classes of morphisms $S_\epsilon \migi \overline{\mfM}_{g,r}$ extending $s_\mfm$.
This correspondence may be regarded as a bijection between {\it pointed} sets by  declaring    the trivial deformation $\msX_\epsilon$ of $\msX$ and   the composite  $s_\mfm   \circ \mr{pr}_S : S_\epsilon \migi \overline{\mfM}_{g,r}$ to be  the distinguished points of the former  and latter sets respectively.
Also, by means of (\ref{QG24}),
first-order deformations of $\msX$ may be identified with first-order deformations (in the logarithmic sense) of  the stable log curve $X^\mr{log}/S^\mr{log}$  of type $(g,r)$.

Hereafter, let us suppose that $S = \mr{Spec}(R)$ for a $k$-algebra $R$ (hence, $\Gamma (S, \mbR^1 f_* (\mcT)) = H^1 (X, \mcT)$).
 Then,  by taking account of the facts mentioned above and the isomorphism  $\mr{KS}_\msX$, we obtain  a canonical bijection of pointed sets
\begin{align} \label{QW6005}
\begin{pmatrix}
\text{the set of isomorphism classes} \\
\text{of first-oder deformations of $\msX$} 
\end{pmatrix}
\cong H^1 (X, \mcT).
\end{align}
In what follows, we recall the explicit  construction of this bijection in terms of \v{C}ech cocycles.
Let $v$ be an element of $H^1 (X, \mcT)$.
We can  find  an affine open covering $\mcU := \{ U_\alpha \}_{\alpha \in I}$ (where $I$ is a finite  index set) of $X$ such that
$v$ may be represented  by a \v{C}ech $1$-cocycle $\{ \partial_{\alpha\beta} \}_{(\alpha, \beta) \in I_2} \in \check{C}^1 (\mcU, \mcT)$,
where 
 $I_2$ denotes the set of pairs $(\alpha, \beta) \in I \times I$ with $U_{\alpha \beta} := U_\alpha \cap U_\beta \neq \emptyset$ and $\partial_{\alpha \beta} := (D_{\alpha \beta}, \delta_{\alpha \beta})$ (for each $(\alpha, \beta) \in I_2$) denotes a logarithmic derivation  belonging to  $\Gamma (U_{\alpha \beta}, \mcT)$.
For each $(\alpha, \beta) \in I_2$,  let us  consider the decompositions
\begin{align}
\mcO_{U_{\alpha\beta, \epsilon}} = \mcO_{U_{\alpha\beta}} \oplus \epsilon \cdot \mcO_{U_{\alpha\beta}}, 
\hspace{5mm}
M_{U_{\alpha \beta, \epsilon}} = M_{U_{\alpha \beta}} \oplus \epsilon \cdot \mcO_{U_{\alpha \beta}}
\end{align}
(where $M_{(-)}$ denotes the sheaf of monoids defininig  the log structure on $(-)$ pulled-back  from $X^\mr{log}$)
 arising from both $\mr{in}_{U_{\alpha \beta}}$ and $\mr{pr}_{U_{\alpha \beta}}$.
 Then, 
$\mr{id}_{\mcO_{U_{\alpha \beta, \epsilon}}} +\epsilon \cdot  D_{\alpha \beta}$ and
$\mr{id}_{M_{U_{\alpha \beta, \epsilon}}} + \epsilon \cdot \delta_{\alpha \beta}$ determines an automorphism of $\mcO_{U_{\alpha \beta, \epsilon}}$ and $M_{U_{\alpha \beta, \epsilon}}$ respectively.
 The pair of these automorphisms 
 corresponds to an  $S_\epsilon^\mr{log}$-automorphism 
\begin{align}
(\tau^{v}_{\alpha \beta})^{\mr{log}} : U^\mr{log}_{\alpha \beta, \epsilon} \isom U^\mr{log}_{\alpha \beta, \epsilon}
\end{align}
of $U_{\alpha \beta, \epsilon}^\mr{log}$.
The automorphisms 
$(\tau^{v}_{\alpha \beta})^{\mr{log}}$ defined for various $(\alpha, \beta) \in I_2$ satisfy the cocycle condition in the usual sense.
Hence,
by means of these  automorphisms,
 the $S_\epsilon^\mr{log}$-log schemes $U_{\alpha, \epsilon}^\mr{log}$ for various $\alpha \in I$ may be glued together 
 to obtain a log scheme
  $(X_{\epsilon'}^{v})^\mr{log}$.
  This log  scheme is  log smooth and integral over $S_\epsilon^\mr{log}$, so  the underlying scheme $X_{\epsilon'}^{v}$ is  flat over $S_\epsilon$  (cf. ~\cite{KaKa}, \S\,4, Corollary (4.5)).
 Moreover,  this log scheme determines a first-order deformation of the stable log curve $X^\mr{log}/S^\mr{log}$  of type $(g,r)$, and hence, corresponds, via (\ref{QG24}),  to a first-order deformation 
  \begin{align} \label{QW6048}
  \msX^v_{\epsilon'} := (X^v_{\epsilon'}/S_\epsilon, \{ \sigma_i^v \}_i)
  \end{align}
 of   $\msX$.
  The resulting assignment $v \mapsto \msX^v_{\epsilon'}$ realizes the bijective  correspondence  (\ref{QW6005}).
  (To construct the inverse to this assignment, we will need  the result of ~\cite{KaKa}, \S\,3, Proposition (3.14).)

 \subsection{Deformation of $\mbG$-bundles}\label{QW6052}
 Next, let us consider the deformation space of a $\mbG$-bundle.
 Let $\pi: \mcE \migi X$ be a $\mbG$-bundle on $X$.
 A {\bf first-oder deformation of $\mcE$}  is defined as a $\mbG$-bundle on $X_\epsilon$ inducing $\mcE$ via  reduction modulo $\epsilon$.
Also, by  a {\bf first-order deformation of $(\msX, \mcE)$}, we mean 
  a pair $(\msX_{\epsilon'}, \mcE_{\epsilon'})$ consisting of 
  a first-order deformation $\msX_{\epsilon'}$ of $\msX$ and a $\mbG$-bundle $\mcE_{\epsilon'}$ inducing $\mcE$ via  reduction modulo $\epsilon$.
In particular, each such pair $(\msX_{\epsilon'}, \mcE_{\epsilon'})$ with 
$\msX_{\epsilon'} = \msX_{\epsilon}$ may be identified with a first-order deformation of $\mcE$.

Under the assumption that $S$ is affine, 
we can verify that first-order deformations of $\mcE$ (resp., $(\msX, \mcE)$) can be described by means of  elements in the $1$-st cohomology group $H^1 (X, \mfg_\mcE)$  (resp., $H^1 (X, \widetilde{\mcT}_{\mcE^\mr{log}/S^\mr{log}})$).
In fact, let us  take an element $\widetilde{v}$ of $H^1 (X, \widetilde{\mcT}_{\mcE^\mr{log}/S^\mr{log}})$.
It is  represented by  a \v{C}ech $1$-cocycle 
$\{ \widetilde{\partial}_{\alpha\beta} \}_{(\alpha, \beta) \in I_2} \in \check{C}^1 (\mcU, \widetilde{\mcT}_{\mcE^\mr{log}/S^\mr{log}})$, where $\mcU := \{ U_\alpha \}_{\alpha \in I}$ is an affine open covering  of $X$ and $\widetilde{\partial}_{\alpha \beta}$ denotes an element of 
$\Gamma (U_{\alpha \beta}, \widetilde{\mcT}_{\mcE^\mr{log}/S^\mr{log}})$.
Each  $\widetilde{\partial}_{\alpha \beta}$  corresponds, via
 $\Xi^{\mcE \rightleftharpoons \mcV}$ (cf. (\ref{QW144})), to a pair 
 $(\widetilde{D}_{\alpha \beta}, \widetilde{\delta}_{\alpha \beta})$
  as in (\ref{Der4}).
If we set  $\partial_{\alpha \beta} := d_\mcE^\mr{log} (\widetilde{\partial}_{\alpha \beta})$, then 
the \v{C}ech $1$-cocycle   $\{ \partial_{\alpha \beta} \}_{(\alpha, \beta) \in I_2} \in \check{C}^1 (\mcU, \mcT)$ represents an element $v \in H^1 (X, \mcT)$, and hence, corresponds to a  first-order deformation $\msX_{\epsilon'}^v := (X^v_{\epsilon'}/S_\epsilon, \{ \sigma^v_i \}_i)$ of $\msX$ 
as discussed above.
Also, $\mr{id}_{\mcO_{(\mcE |_{U_{\alpha \beta}})_\epsilon}} + \epsilon \cdot \widetilde{D}_{\alpha \beta}$  and 
$\mr{id}_{M_{(\mcE |_{U_{\alpha \beta}})_\epsilon}} + \epsilon \cdot \widetilde{\delta}_{\alpha \beta}$ determine 
the automorphisms of
 $\mcO_{(\mcE |_{U_{\alpha \beta}})_\epsilon} \left(= \mcO_{\mcE |_{U_{\alpha \beta}}} \oplus \epsilon \cdot \mcO_{\mcE |_{U_{\alpha \beta}}} \right)$
 and $M_{(\mcE |_{U_{\alpha \beta}})_\epsilon} \left(= M_{\mcE |_{U_{\alpha \beta}}} \oplus \epsilon \cdot \mcO_{\mcE |_{U_{\alpha \beta}}} \right)$ respectively.
The pair of these automorphisms
corresponds to  a $\mbG$-equivariant  $S_\epsilon^\mr{log}$-automorphism 
\begin{align}
(\tau^{\widetilde{v}}_{\alpha \beta})^{\mr{log}} : (\mcE^\mr{log} |_{U_{\alpha \beta}})_\epsilon \isom (\mcE^\mr{log} |_{U_{\alpha \beta}})_\epsilon
\end{align}
 of  $(\mcE^\mr{log} |_{U_{\alpha \beta}})_\epsilon$,  which makes  the following square diagram commute:
\begin{align} \label{QW6902}
\vcenter{\xymatrix@C=46pt@R=36pt{
(\mcE^\mr{log} |_{U_{\alpha \beta}})_\epsilon \ar[r]_-{\sim}^-{(\tau^{\widetilde{v}}_{\alpha \beta})^\mr{log} } \ar[d]_-{\mr{id}_{k_\epsilon} \times \pi |_{\mcE |_{U_{\alpha \beta}}}} & (\mcE^\mr{log} |_{U_{\alpha \beta}})_\epsilon \ar[d]^-{\mr{id}_{k_\epsilon} \times \pi |_{\mcE |_{U_{\alpha \beta}}}}\\
 U^\mr{log}_{\alpha \beta, \epsilon}\ar[r]^-{\sim}_-{(\tau_{\alpha \beta}^v)^\mr{log}} & U^\mr{log}_{\alpha \beta, \epsilon}.
}}
\end{align}
By means of  the various  $(\tau^{\widetilde{v}}_{\alpha \beta})^\mr{log}$'s, 
 the $\mbG$-bundles  $(\mcE |_{U_\alpha})_\epsilon$ may be glued together to obtain a $\mbG$-bundle
 $\mcE^{\widetilde{v}}_{\epsilon'}$ on  $X^v_{\epsilon'}$.
The pair $(\msX^{v}_{\epsilon'}, \mcE^{\widetilde{v}}_{\epsilon'})$ specifies  a first-order deformation of $(\msX, \mcE)$ and 
the resulting assignment $v \mapsto (\msX^{v}_{\epsilon'}, \mcE^{\widetilde{v}}_{\epsilon'})$ gives a bijection between pointed sets
       \begin{align} \label{QW6080}
\begin{pmatrix}
\text{the set of isomorphism classes} \\
\text{of first-oder deformations of $(\msX, \mcE)$} 
\end{pmatrix}
\cong H^1 (X,  \widetilde{\mcT}_{\mcE^\mr{log}/S^\mr{log}})
\end{align}
(cf. ~\cite{Chen}, \S\,4, Proposition 4.1.3), where the distinguished point of the left-hand  side is taken  as the trivial deformation.
  (To verify the surjectivity of this assignment, we may need the fact that $H^1 (X,  \widetilde{\mcT}_{\mcE^\mr{log}/S^\mr{log}})$ is canonically isomorphic to the $1$-st cohomology group  of the \'{e}tale sheaf associated to $\widetilde{\mcT}_{\mcE^\mr{log}/S^\mr{log}}$, each of whose elements may be represented by a   \v{C}ech $1$-cocycle with respect to an  \'{e}tale covering; this is because the $\mbG$-bundle $\mcE$ is locally trivial in the \'{e}tale topology, but not necessarily in the Zariski topology.)
  Moreover,  by the above discussion with   $\widetilde{\mcT}_{\mcE^\mr{log}/S^\mr{log}}$ replaced with $\mfg_\mcE \left(= \widetilde{\mcT}_{\mcE^\mr{log}/X^\mr{log}} \right)$, we obtain   a bijection 
    \begin{align} \label{QW6077}
\begin{pmatrix}
\text{the set of isomorphism classes} \\
\text{of first-oder deformations of $\mcE$} 
\end{pmatrix}
\cong H^1 (X,  \mfg_\mcE).
\end{align}
By construction, the natural maps $H^1 (X,  \widetilde{\mcT}_{\mcE^\mr{log}/S^\mr{log}}) \migi H^1 (X, \mcT)$ and  $H^1 (X, \mfg_\mcE) \migi H^1 (X,  \widetilde{\mcT}_{\mcE^\mr{log}/S^\mr{log}})$ induced by the morphisms in (\ref{Ex0}) correspond, via (\ref{QW6005}), (\ref{QW6080}), and (\ref{QW6077}),   to  the assignments of deformations $(\msX_{\epsilon'}, \mcE_{\epsilon'}) \mapsto \msX_{\epsilon'}$ and $\mcE_{\epsilon'} \mapsto (\msX_\epsilon, \mcE_{\epsilon'})$ respectively.

\begin{rema}[Infinitesimal gauge transformation]\label{QW9004}
Let us take a section $v := (\widetilde{D}, \widetilde{\delta})$ of $\widetilde{\mcT}_{\mcE^\mr{log}/S^\mr{log}}$ over an open subscheme $U$ of $X$.
Let $(\eta^v_U)^\mr{log} : U^\mr{log}_\epsilon \isom U^\mr{log}_\epsilon$ and $(\eta^{v}_\mcE)^\mr{log} : \mcE^\mr{log}_\epsilon |_U \isom \mcE^\mr{log}_\epsilon |_U$ be the $S_\epsilon^\mr{log}$-automorphisms corresponding to
the pairs 
\begin{align} \label{QH40}
(\mr{id}_{U_\epsilon} + \epsilon \cdot \widetilde{D} |_\Delta, \mr{id}_{M_{U_\epsilon}} + \epsilon \cdot \widetilde{\delta}) \ \ \text{and} \ \ 
(\mr{id}_{\mcE_\epsilon |_U} + \epsilon \cdot \widetilde{D}, \mr{id}_{M_{\mcE_\epsilon |_U}} + \epsilon \cdot \widetilde{\delta})
\end{align}
 respectively.
Then,  $(\eta_\mcE^v)^\mr{log}$ preserves the $\mbG$-action and the following square diagram is commutative:
\begin{align} \label{QW9003}
\vcenter{\xymatrix@C=46pt@R=36pt{
\mcE^\mr{log}_\epsilon |_U \ar[r]_-{\sim}^{(\eta_\mcE^{v})^\mr{log}} \ar[d]_-{\mr{id}_{k_\epsilon}\times\pi |_U} & \mcE^\mr{log}_\epsilon |_U \ar[d]^-{\mr{id}_{k_\epsilon}\times \pi |_U} \\
U^\mr{log}_\epsilon \ar[r]^-{\sim}_{(\eta^v_U)^\mr{log}}& U^\mr{log}_\epsilon.
}}
\end{align}
By this diagram, we can think of the $\mbG$-bundle $\mcE_\epsilon^\mr{log} |_U$ on the left top corner  as the pull-back of  $\mcE^\mr{log}_\epsilon |_U$ by $(\eta_U^{v*})^\mr{log}$.
One verifies that the automorphism $d (\eta^{v}_\mcE)^\mbG$ (resp., $d (\eta^{v}_U)^\mbG$) of  $\widetilde{\mcT}_{\mcE_\epsilon^\mr{log}|_U/S_\epsilon^\mr{log}}$ (resp., $\mcT_{U_\epsilon^\mr{log}/S^\mr{log}_\epsilon}$) obtained by differentiating $\eta^{v}_\mcE$ (resp., $\eta^{v}_U$)  is given by $w \mapsto w - \epsilon \cdot \mr{ad}(v) (w)$ (resp., $w \mapsto w - \epsilon \cdot \mr{ad}(d_\mcE^\mr{log}(v)) (w)$) under the natural identification
 \begin{align}
 \widetilde{\mcT}_{\mcE_\epsilon^\mr{log}|_U/S_\epsilon^\mr{log}} = \widetilde{\mcT}_{\mcE^\mr{log}|_U/S^\mr{log}} \oplus \epsilon \cdot \widetilde{\mcT}_{\mcE^\mr{log}|_U/S^\mr{log}} \hspace{-1mm} \\
   \left(\text{resp.,}  \ \mcT_{U_\epsilon^\mr{log}/S^\mr{log}_\epsilon} = \mcT_{U^\mr{log}/S^\mr{log}} \oplus \epsilon \cdot \mcT_{U^\mr{log}/S^\mr{log}} \right). \notag
 \end{align}

Now, let $\hslash$ be an element of $\Gamma (S, \mcO_S)$ and 
 $\nabla_{\epsilon'}$
an  $\hslash$-$S^\mr{log}$-connection on $\mcE_\epsilon |_U$.
 We shall set  $\nabla := \mr{in}_S^* (\nabla_{\epsilon'})$.
Then, by the above discussion, 
the pull-back $\eta^{v *}_\mcE (\nabla_{\epsilon'})$
$\left(=(d (\eta^{v}_\mcE)^\mbG)^{-1} \circ \nabla_{\epsilon'} \circ d (\eta^{v}_U)^\mbG \right)$ can be described as 
\begin{align} \label{QW8600}
\eta_\mcE^{v *} (\nabla_{\epsilon'})(\partial)   = \nabla_{\epsilon'} (\partial) -  \epsilon \cdot \left([\nabla (\partial), v] - \nabla ([\partial, d_\mcE^\mr{log}(v)]) \right). 
\end{align}
for any local section $\partial \in \mcT_{U^\mr{log}_\epsilon/S^\mr{log}_\epsilon}$.
In particular, if $v$ lies in $\mfg_\mcE$, then $\eta^{v}_\mcE$ specifies an automorphism of the $\mbG$-bundle $\mcE_\epsilon |_U$ (i.e., $\eta^v_U = \mr{id}_{U_\epsilon}$) and (\ref{QW8600}) becomes
\begin{align} \label{QW8601}
\eta_\mcE^{v *} (\nabla_{\epsilon'})(\partial)    = \nabla_{\epsilon'} (\partial) -  \epsilon \cdot [\nabla (\partial), v].
\end{align}
In the notation of the subsequent discussion, (\ref{QW8600}) (resp., (\ref{QW8601})) reads  
$\eta^{v *}_\mcE (\nabla_{\epsilon'})= \nabla_{\epsilon'}  -  \epsilon \cdot \widetilde{\nabla}^{\mr{ad}}(v)$ 
(resp., 
$\eta^{v *}_\mcE (\nabla_{\epsilon'})   = \nabla_{\epsilon'}  -  \epsilon \cdot \nabla^{\mr{ad}}(v)$).
\end{rema}

 \subsection{Deformation of $\hslash$-flat $\mbG$-bundles} \label{QW6081}
 Next, let us consider deformations of a flat $\mbG$-bundle.
Let  $\hslash$ be an element of $\Gamma (S, \mcO_S)$ and
 $\nabla$ an $\hslash$-$S^\mr{log}$-connection on $\mcE$. 
Denote by 
\index{$\nabla^\mr{ad}$}
\begin{equation} 
\nabla^\mr{ad} : \mfg_{\mcE} \migi \Omega\otimes \mfg_{\mcE}
 \label{Ad} 
\end{equation}
the $\hslash$-$S^\mr{log}$-connection on the adjoint vector  bundle $ \mfg_{\mcE}$ induced by $\nabla$  via the adjoint representation $\mr{Ad}_\mbG : \mbG \migi \mr{GL}(\mfg)$.

A {\bf first-order deformation of $(\mcE, \nabla)$} is an $\hslash$-flat $\mbG$-bundle $(\mcE_{\epsilon'}, \nabla_{\epsilon'})$  on $\msX_\epsilon$ inducing $(\mcE, \nabla)$ via reduction modulo $\epsilon$.
If $S$ is affine, then such deformations may be described by means of elements in  the $1$-st hypercohomology group $\mbH^1 (X, \mcK^\bullet [\nabla^\mr{ad}])$ associated to  the complex $\mcK^\bullet [\nabla^\mr{ad}]$  (cf. ~\cite{O3}, Proposition 3.6).
That is to say,  there exists 
 a canonical bijection of pointed sets
   \begin{align} \label{QW6082}
\begin{pmatrix}
\text{the set of isomorphism classes of} \\
\text{first-oder deformations of $(\mcE, \nabla)$} 
\end{pmatrix}
\cong \mbH^1 (X,  \mcK^\bullet [\nabla^\mr{ad}]).
\end{align}
To construct this, recall that
one may calculate $\mbH^1 (\mcK^\bullet[\nabla^\mr{ad}])$ as the total cohomology   of the \v{C}ech double complex   $\mr{Tot}^\bullet (\check{C}^\bullet (\mcU, \mcK^\bullet [\nabla^\mr{ad}]))$
 associated to $\mcK^\bullet [\nabla^\mr{ad}]$.
Each element $v$ of $\mbH^1 (\mcK^\bullet[\nabla^\mr{ad}])$  may be given by a collection of data
\begin{align}\label{W105}
 v = (\{ \partial_{\alpha \beta} \}_{(\alpha, \beta)}, \{ \delta_\alpha \}_{\alpha})
\end{align}
consisting of a \v{C}ech $1$-cocycle 
$\{ \partial_{\alpha\beta} \}_{(\alpha, \beta) \in I_2} \in \check{C}^1 (\mcU, \mfg_{\mcE})$
 (where $\partial_{\alpha\beta} \in \Gamma (U_{\alpha \beta}, \mfg_{\mcE})$) and a \v{C}ech $0$-cochain $\{ \delta_\alpha \}_{\alpha \in I} \in \check{C}^0 (\mcU, \Omega  \otimes \mfg_\mcE)$
(where $\delta_\alpha \in \Gamma (U_\alpha, \Omega \otimes \mfg_\mcE)$ $= \mr{Hom}_{\mcO_{U_\alpha}}(\mcT |_{U_\alpha},  \mfg_\mcE |_{U_\alpha})$)
which agree under $\nabla^\mr{ad}$ and the \v{C}ech coboundary map.
 The flat $\mbG$-bundles  
 $\msE_{\alpha}^v := (\mcE_\epsilon |_{U_\alpha}, \mr{pr}^*_S(\nabla) |_{\mcE |_{U_\alpha}}+ \epsilon \cdot \delta_\alpha)$
 may be glued together by means of the isomorphisms
 $\msE_{\beta}^v|_{U_{\alpha \beta}}   \isom \msE_{\alpha}^v |_{U_{\alpha \beta}} $  determined by the automorphism 
 \begin{align} \label{QW8400}
\mr{id}_{\mcO_{\mcE_\epsilon |_{U_{\alpha\beta}}}} + \epsilon \cdot \partial_{\alpha \beta}
 :  \Gamma (\mcE_\epsilon |_{U_{\alpha\beta}}, \mcO_{\mcE_\epsilon |_{U_{\alpha\beta}}})
\end{align}
of $\mcO_{\mcE_\epsilon |_{U_{\alpha\beta}}}$.
The resulting $\hslash$-flat $\mbG$-bundle on $X^\mr{log}/S^\mr{log}$, which will be denoted  by $(\mcE_{\epsilon'}^v, \nabla_{\epsilon'}^v)$, specifies a first-order deformation of $(\mcE, \nabla)$.
 One verifies that the assignment $v \mapsto (\mcE_{\epsilon'}^v, \nabla_{\epsilon'}^v)$
obtained in this way gives the desired correspondence.
(To prove the injectivity of that assignment, we will use  (\ref{QW8601}).)

\subsection{Deformation of $(\msX, \mcE, \nabla)$'s} \label{QW7000}
Finally, we shall describe the deformation space of the triple $(\msX, \mcE, \nabla)$ by means of the hypercohomology of a certain complex. 

 For any local sections
$v \in  \widetilde{\mcT}_{\mcE^{\mr{log}}/S^\mr{log}}$ and
$\partial \in \mcT$, 
we shall write $\widetilde{\nabla}^\mr{ad}_\partial (v):= [\nabla (\partial), v] - \nabla ([\partial, d^\mr{log}_{\mcE}(v)]) \in \widetilde{\mcT}_{\mcE^{\mr{log}}/S^\mr{log}}$.

\ble \label{y0145}
The local section $\widetilde{\nabla}^\mr{ad}_\partial (v)$
lies  in 
$\mfg_{\mcE}$.
\ele
\bpf
Note that $d^\mr{log}_{\mcE}$ preserves the Lie bracket operations.
Hence, for any local sections $v \in  \widetilde{\mcT}_{\mcE^{\mr{log}}/S^\mr{log}}$ and  $\partial \in \mcT$,    the following sequence of equalities holds:
 \begin{align}   d^\mr{log}_{\mcE}  ( \widetilde{\nabla}_{\partial}^{\mr{ad}} (v)) 
 & = d^\mr{log}_{\mcE} ( [\nabla (\partial), v] - \nabla ([\partial,  d^\mr{log}_{\mcE}(v)]) )  \\
 & = [d^\mr{log}_{\mcE} \circ \nabla (\partial), d^\mr{log}_{\mcE} (v)]-  d^\mr{log}_{\mcE} \circ \nabla ([\partial, d^\mr{log}_{\mcE}(v)])  \notag \\
&
 =  [\hslash \cdot \partial, d^\mr{log}_{\mcE} (v)] -\hslash \cdot [\partial, d^\mr{log}_{\mcE} (v)] \notag \\
 &  = 0. \notag
 \end{align} 
It follows that the image of $\widetilde{\nabla}^\mr{ad}_\partial (v)$ lies  in 
$\Omega\otimes \mr{Ker}(d^\mr{log}_{\mcE})  \left(= \Omega\otimes  \mfg_{\mcE}\right)$.
\epf

By the above lemma, we can verify that the assignment
$v \mapsto  \widetilde{\nabla}^{\mr{ad}}_{(-)} (v)$ defines 
an  $f^{-1}(\mcO_S)$-linear morphism
\index{$\nabla^\mr{ad}$@--- $\widetilde{\nabla}^{\mr{ad}}$}
\begin{equation} 
 \widetilde{\nabla}^{\mr{ad}} :  \widetilde{\mcT}_{\mcE^{\mr{log}}/S^\mr{log}} \migi  \left(\mcH om_{\mcO_X} (\mcT,   \mfg_\mcE)= \right)\Omega\otimes  \mfg_\mcE.\label{circledast}
\end{equation}
This morphism fits into the following morphism of short exact sequences:
\begin{align} \label{QW7010}
\vcenter{\xymatrix@C=40pt@R=36pt{
0 \ar[r] & \mfg_\mcE \ar[r]^-{\mr{inclusion}}\ar[d]^{\nabla^\mr{ad}}& \widetilde{\mcT}_{\mcE^\mr{log}/S^\mr{log}} \ar[r]^-{d^\mr{log}_{\mcE}}\ar[d]^{\widetilde{\nabla}^\mr{ad}}& \mcT  \ar[r] \ar[d] & 0 \\
0 \ar[r]& \Omega \otimes \mfg_\mcE \ar[r]_-{\mr{id}}& \Omega \otimes \mfg_\mcE \ar[r] & 0  \ar[r] & 0. 
}}
\end{align}
In other words,
we have  the following  short exact sequence of complexes:
\begin{align} \label{QW7034}
0 \longmigi \mcK^\bullet [\nabla^\mr{ad}] \longmigi \mcK^\bullet [\widetilde{\nabla}^\mr{ad}] \longmigi \mcT[0] \longmigi 0. 
\end{align}
Now, suppose that $S$ is affine.
Since $\Gamma (X, \mcT) = H^2 (X, \mcT) = 0$, the above  sequence  gives rise to 
a short exact sequence
\begin{align} \label{QW7011}
0 \longmigi \mbH^1 (\mcK^\bullet [\nabla^\mr{ad}]) \stackrel{}{\longmigi} \mbH^1 (\mcK^\bullet [\widetilde{\nabla}^\mr{ad}]) \stackrel{}{\longmigi} H^1 (X, \mcT) \longmigi 0.
\end{align}
Here,  a {\bf first-order deformation of $(\msX, \mcE, \nabla)$} is defined as   a  triple $(\msX_{\epsilon'}, \mcE_{\epsilon'}, \nabla_{\epsilon'})$ consisting of a first-order deformation $(\msX_{\epsilon'}, \mcE_{\epsilon'})$ of $(\msX, \mcE)$ and 
an $\hslash$-$S_\epsilon^\mr{log}$-connection $\nabla_{\epsilon'}$
on $\mcE_{\epsilon'}$ inducing $\nabla$ via reduction modulo $\epsilon$.
By referring to the constructions of (\ref{QW6080}) and (\ref{QW6082}), 
  we can construct, in an entirely similar manner,  the following bijective correspondence:
   \begin{align} \label{QW6083}
\begin{pmatrix}
\text{the set of isomorphism classes of} \\
\text{first-oder deformations of $(\msX, \mcE, \nabla)$} 
\end{pmatrix}
\cong \mbH^1 (X,  \mcK^\bullet [\widetilde{\nabla}^\mr{ad}])
\end{align}
(cf. ~\cite{Chen}, \S\,4, Proposition 4.3.1).
By construction, the morphisms $\mbH^1 (X, \mcK^\bullet [\widetilde{\nabla}^\mr{ad}]) \migi H^1 (X, \mcT)$ and $\mbH^1 (X, \mcK^\bullet [\nabla^\mr{ad}]) \migi \mbH^1 (X, \mcK^\bullet [\widetilde{\nabla}^\mr{ad}])$ appearing in  (\ref{QW7011}) correspond, via  (\ref{QW6005}), (\ref{QW6082}), and (\ref{QW6083}),   to  the assignments of deformations $(\msX_{\epsilon'}, \mcE_{\epsilon'}, \nabla_{\epsilon'}) \mapsto \msX_{\epsilon'}$ and $(\mcE_{\epsilon'}, \nabla_{\epsilon'}) \mapsto (\msX_\epsilon, \mcE_{\epsilon'}, \nabla_{\epsilon'})$.
 Also, the map $\mbH^1 (X, \mcK^\bullet [\widetilde{\nabla}^\mr{ad}]) \migi H^1 (X, \widetilde{\mcT}_{\mcE^\mr{log}/S^\mr{log}})$ induced by  the natural morphism $\mcK^\bullet [\widetilde{\nabla}^\mr{ad}] \migi \widetilde{\mcT}_{\mcE^\mr{log}/S^\mr{log}}[0]$ is compatible with the assignment of deformations $(\msX_{\epsilon'}, \mcE_{\epsilon'}, \nabla_{\epsilon'}) \mapsto (\msX_{\epsilon'}, \mcE_{\epsilon'})$.
 Finally, it is immediate that the image of the map $\Gamma (X, \Omega \otimes \mfg_\mcE) \migi \mbH^1 (X, \mcK^\bullet [\widetilde{\nabla}^\mr{ad}])$ induced by the natural morphism $\Omega \otimes \mfg_\mcE [-1] \migi \mcK^\bullet [\widetilde{\nabla}^\mr{ad}]$
 classifies the deformations $(\msX_{\epsilon'}, \mcE_{\epsilon'}, \nabla_{\epsilon'})$'s  with $(\msX_{\epsilon'}, \mcE_{\epsilon'}) \cong (\msX_\epsilon, \mcE_\epsilon)$.

\section{Cohomology of complexes  associated to a $(\mfg, \hslash)$-oper}\label{y0140}

In the rest of this chapter,
let 
 $k$
 \index{$k$, field} 
  be   a perfect field,   $(\mbG, \mbT)$ 
 \index{$\mbG$, algebraic group}
 \index{$\mbG$, algebraic group@— $\mbT$, maximal torus of}
  a split  semisimple algebraic group of adjoint type over $k$,
$\mbB$ 
\index{$\mbG$, algebraic group@— $\mbB$, Borel subgroup of}
a Borel subgroup of $\mbG$ defined over $k$ containing the maximal torus $\mbT$, and $\ \qq := \{ x_\alpha \}_{\alpha \in \Gamma}$
\index{$\ \qq\,$}
 a collection  of generators $x_\alpha$ of $\mfg^\alpha$'s (for simple positive roots $\alpha$ in $\mbB$).
Denote by $\mfg$ 
\index{$\mfg$, Lie algebra}
and $\mfb$ the Lie algebras of $\mbG$ and $\mbB$ respectively. We suppose that {\it either ``\,$\mr{char}(k) =0$'' or ``\,$\mr{char}(k) > 2  h_\mbG$'' is fulfilled}.

In this section, we compute  the  cohomology of complexes associated to a $(\mfg, \hslash)$-oper; it is carried out by  using the filtration arising from the Cartan decomposition of $\mfg$.

\subsection{Spectral sequences}\label{QW1003}
First, we recall two spectral sequences  associated to a morphism of sheaves.
Let $\nabla : \mcK^0 \migi \mcK^1$ be a morphism of sheaves of abelian groups on $X$.

Recall that there exists  a spectral sequence of the form
\begin{equation}
{'E}^{p,q}_{1} :=\mbR^qf_*(\mcK^p[\nabla]) \Rightarrow \mbR^{p+q}f_*(\mcK^\bullet[\nabla])  \label{HDsp}
\end{equation}
which is the  {\it Hodge to de Rham spectral sequence} of the complex $\mcK^\bullet [\nabla]$.
This spectral sequence
 yields a short exact sequence
   \index{$'e_\sharp [\nabla]$, ${'e}_{\flat}[\nabla]$, ${''e_\sharp[\nabla]}$, ${''e}_\flat[\nabla]$}
\begin{equation}   
 0 \longmigi \mr{Coker}(\mbR^0 f_*(\nabla)) \xrightarrow{'e_\sharp [\nabla]} \mbR^1f_*(\mcK^\bullet [\nabla]) \xrightarrow{{'e}_{\flat}[\nabla]}   \mr{Ker}(\mbR^1f_*(\nabla))   \longmigi 0, \label{HDeq}
 \end{equation}
where  $\mbR^if_*(\nabla)$ ($i = 0, 1$) denotes the morphism $\mbR^if_*(\mcK^0) \migi \mbR^if_*(\mcK^1)$ 
induced by $\nabla$.

On the other hand,   the {\it conjugate spectral sequence} of $\mcK^\bullet [\nabla]$ is a spectral sequence
\begin{equation}
''E_{2}^{p,q} := \mbR^p(\mcH^q(\mcK^\bullet [\nabla])) \Rightarrow \mbR^{p+q}f_*(\mcK^\bullet [\nabla]),  \label{Csp}
\end{equation}
where $\mcH^q(\mcK^\bullet [\nabla])$ denotes the $q$-th cohomology sheaf of the complex $\mcK^\bullet [\nabla]$.
This spectral sequence 
 induces a  short exact sequence
  \index{$'e_\sharp [\nabla]$, ${'e}_{\flat}[\nabla]$, ${''e_\sharp[\nabla]}$, ${''e}_\flat[\nabla]$}
\begin{equation}    0 \longmigi \mbR^1f_*(\mr{Ker}(\nabla)) \xrightarrow{{''e_\sharp[\nabla]}}
 \mbR^1f_*(\mcK^\bullet[\nabla]) \xrightarrow{{''e}_\flat[\nabla]}
  f_*(\mr{Coker}(\nabla)) \longmigi 0. 
\label{Ceq}
\end{equation}

\subsection{Filtrations of  $\nabla^{\mr{ad}}$ and $\dV_\mfg^\vee$} \label{QW2034}
 Let $\hslash \in k$, and 
 $\msE^{\spadesuit} = (\DE_\mbB, \nabla)$ be a $\ \qq$-normal  $(\mfg, \hslash)$-oper on $\msX$.
Denote by 
$\nabla^\mr{ad} : \mfg_{\DE_{\mbB}} \migi \Omega\otimes \mfg_{\DE_{\mbB}}$
the $\hslash$-$S^\mr{log}$-connection on the adjoint vector  bundle $ \mfg_{\DE_{\mbB}}$ ($=  \mfg_{\DE_{\mbG}} = \mfg_{\DE_{\mbB^\odot}}$) induced by $\nabla$.
If $j \in \mbZ$, then
 the definition of a $(\mfg, \hslash )$-oper (and the fact that $[q_{-1}, \mfg^j] \subseteq \mfg^{j-1}$) implies that
$\nabla^\mr{ad} (\mfg_{\DE_{\mbB}}^j) \subseteq \Omega\otimes \mfg_{\DE_{\mbB}}^{j-1}$.  
Hence, we can obtain a morphism
\index{$\nabla^\mr{ad}$@--- $\nabla^{\mr{ad}(j)}$}
\begin{equation}
\nabla^{\mr{ad}(j)} : \mfg_{\DE_{\mbB}}^j \migi \Omega \otimes \mfg_{\DE_{\mbB}}^{j-1}
\label{Ad(i)}
\end{equation}
  by restricting $\nabla^\mr{ad}$.
This morphism  induces a well-defined {\it $\mcO_X$-linear} morphism
\index{$\nabla^\mr{ad}$@--- $\nabla^{\mr{ad}(j/j+1)}$}
\begin{equation} 
 \nabla^{\mr{ad}(j/j+1)} : \mfg_{\DE_{\mbB}}^j/\mfg_{\DE_{\mbB}}^{j+1}
  \migi
    \Omega \otimes (\mfg_{\DE_{\mbB}}^{j-1}/ \mfg_{\DE_{\mbB}}^{j}).
     \label{g/g}
      \end{equation}

Next, let us recall the vector bundle $\dV_\mfg := \dV_{\mfg,  X^\mr{log}/S^\mr{log}}$ (cf. (\ref{QQ801})).
The dual of this bundle has  a natural  isomorphism
 $(\mfg^{\mr{ad}(q_1)} )^\vee_{\DE_{\mbB^\odot}} \isom \Omega\otimes \dV^\vee_\mfg$, which restricts to an isomorphism $(\mfg_j^{\mr{ad}(q_1)} )^\vee_{\DE_{\mbB^\odot}} \isom \Omega\otimes (\dV_{\mfg, j})^\vee$ for each $j \in \mbZ$.
Let us define  a decreasing  filtration $\{ \dV_\mfg^{\vee, j} \}_{j \in \mbZ}$ on  $\dV^\vee$ by putting
\begin{align}
\dV^{\vee, j}_\mfg := (\dV_\mfg /\dV_\mfg^{-j+1})^\vee \left(\subseteq \dV_\mfg \right).
\end{align}
In particular, we have $\dV_\mfg^{\vee, j}/\dV_\mfg^{\vee, j+1} = (\dV_{\mfg, -j})^{\vee}$.


\subsection{The kernel/cokernel of $\nabla^{\mr{ad}(j/j+1)}$} \label{QW2035}
 In what follows, we  describe the kernel and cokernel of $\nabla^{\mr{ad}(j/j+1)}$ by using $\dV_{\mfg, j}$'s.
 Recall that 
 $\mfg$ is equipped with a $\mbB^\odot$-action given by the composite $\mbB^\odot \xrightarrow{\iota_\mbB} \mbB \xrightarrow{\mr{inclusion}} \mbG \xrightarrow{\mr{Ad}_\mbG} \mr{GL}(\mfg)$.
 This action preserves the filtration $\{ \mfg^j \}_j$, so it induces  a $\mbB^\odot$-action  on $\mfg^j/\mfg^{j+1}$ for each $j \in \mbZ$.
It is compatible  with the natural $\mbG_m$-action on $\mfg_j \left(=\mfg^j/\mfg^{j+1}\right)$ via the quotient
  $\mbB^\odot \migisurj \mbB^\odot /\mbN^\odot \left(= \mr{GL}(\mfn^\odot) = \mbG_m\right)$.
 
  Next, for $\Box \in \{-1,  1\}$, denote by  $\mfg_{\mr{ad}(q_\Box)}$ 
  \index{$\mfg$, Lie algebra@--- $\mfg^{\mr{ad}(q_1)}$, $\mfg^{\mr{ad}(q_{-1})}$, $\mfg_{\mr{ad}(q_1)}$, $\mfg_{\mr{ad}(q_{-1})}$}
    the space of $\mr{ad}(q_\Box)$-coinvariants, i.e., $\mfg_{\mr{ad}(q_\Box)} := \mfg / \mr{Im}(\mr{ad}(q_\Box))$.
The decomposition  $\mfg = \bigoplus_{j \in \mbZ} \mfg_j$  induces a decomposition  $\mfg_{\mr{ad}(\Box)} = \bigoplus_{j \in \mbZ} \mfg_{\mr{ad}(\Box), j}$ 
 via the quotient $\mfg \migisurj \mfg_{\mr{ad}(\Box)}$.
Also, each component  $\mfg_{\mr{ad}(\Box), j}$ is equipped with a well-defined $\mbG_m$-action arising from the $\mbG_m$-action on $\mfg_j$.
Since $\{ q_{-1}, 2 \check{\rho}, q_1 \}$ forms an $\mfs \mfl_2$-triple and either ``\,$\mr{char}(k) =0$'' or ``\,$\mr{char}(k) > 2  h_\mbG$'' is  fulfilled,
 the two composites
\begin{equation}   \mfg^{\mr{ad}(q_1)} \migiincl \mfg \migisurj  \mfg_{\mr{ad}(q_{-1})}, \hspace{5mm}  \mfg^{\mr{ad}(q_{-1})} \migiincl \mfg \migisurj \mfg_{\mr{ad}(q_1)}\label{composite}
\end{equation}
are verified to be  isomorphisms (cf. (\ref{QR02})) and preserve the $\mbG_m$-actions.
 There isomorphisms 
  restrict, for any $j \in \mbZ$,  to   $\mbG_m$-equivariant  isomorphisms 
\begin{equation}
 \mfg_j^{\mr{ad}(q_1)} \migiincl \mfg \migisurj  \mfg_{\mr{ad}(q_{-1}), j}\hspace{3mm} \text{and} \hspace{3mm} 
   \mfg^{\mr{ad}(q_{-1})}_j  \isom \mfg_{\mr{ad}(q_1),j} \label{isom8} \end{equation}
   respectively.
On the other hand, the  isomorphism $\mfg \isom \mfg^\vee$ determined by the Killing form on $\mfg$ (cf. Remark \ref{QR1000})  induces, via taking subquotients, 
an isomorphism  
\begin{equation} 
 \left(\mfg_j^{\mr{ad}(q_{-1})} \stackrel{(\ref{isom8})}{\cong}\right)  \mfg_{\mr{ad}(q_1), j} \isom ( \mfg_{-j}^{\mr{ad}(q_1)} )^\vee.
\label{isom3}
\end{equation}
By
taking account of  (\ref{isom8}) and  (\ref{isom3}) (and the fact that $q_{-1} = \iota_\mfg (q_{-1}^\odot)\in \iota_{\mfg}(\mfg^\odot_{-1})$), we obtain
 an exact sequence
\begin{align} \label{composite5}
\vcenter{\xymatrix@C=6pt@R=36pt{
 0 
 \ar[rd]^-{}_-{}
 &&\mfg_j \left(=\mfg^j/\mfg^{j+1}\right)
  \ar[rd]^{\mr{ad}(q_{-1})}_-{}& &  \mfg^{\mr{ad}(q_1)}_{j-1}(1)  \ar[dr]&
\\
&   ( \mfg^{\mr{ad}(q_1)}_{-j} )^\vee \ar[ur]_-{}^-{}& &
 \mfg_{j-1}(1) \left(= \mfg^{j-1}/\mfg^j(1)\right) \ar[ur] &&0,
}}
\end{align}
   all of whose arrows are $\mbG_m$-equivariant.
Since $\msE^\spadesuit$ is $\ \qq$-normal, 
 the morphism $\nabla^{\mr{ad}(j/j+1)}$ fits into the following exact sequence obtained by twisting   this sequence  by $\DE_{\mbB^\odot} \times^{\mbB^\odot} \mbG_m \left(= \Omega^{\times} \right)$:
\begin{align} \label{QH31}
\vcenter{\xymatrix@C=10pt@R=36pt{
 0 
 \ar[rd]^-{}_-{}
 && \mfg_{\DE_{\mbB}}^j /  \mfg_{\DE_{\mbB}}^{j+1}
  \ar[rd]^{\nabla^{\mr{ad}(j/j+1)}}_-{}& & \dV_{\mfg, j-1}  \ar[dr]&
\\
&  \Omega \otimes (\dV_{\mfg, -j})^\vee \ar[ur]_-{}^-{}& &
\Omega \otimes (\mfg_{\DE_{\mbB}}^{j-1} /  \mfg_{\DE_{\mbB}}^{j})  \ar[ur] &&0.
}}
\end{align}


Now, let $j$ be an integer with
 $j \leq 1$ (resp., $j \geq 0$).
 Since $\mfg^{\mr{ad}(q_1)}_{j-1} =0$ (resp., $( \mfg^{\mr{ad}(q_1)}_{-j-1} )^\vee =0$),  (\ref{composite5}) becomes a  short exact sequence 
\begin{align} \label{composite52} 
  0 \longmigi  ( \mfg^{\mr{ad}(q_1)}_{-j} )^\vee \longmigi \mfg^j/\mfg^{j+1} \xrightarrow{\mr{ad}(q_{-1})} \mfg^{j-1}/\mfg^j (1)\longmigi 0
  \\
  \left(\text{resp.,} \   0  \migi \mfg^j/\mfg^{j+1} \xrightarrow{\mr{ad}(q_{-1})} \mfg^{j-1}/\mfg^j  (1)\migi  \mfg^{\mr{ad}(q_1)}_{j-1} (1)\migi 0 \right). \notag
  \end{align}
 This sequence  is split by the $\mbG_m$-equivariant  composite surjection $\mfg^j/\mfg^{j+1} \left(=\mfg_j\right) \migisurj \mfg_{\mr{ad}(q_1), j} \xrightarrow{(\ref{isom3})} (\mfg_{-j}^{\mr{ad}(q_1)})^\vee$ (resp.,  the $\mbG_m$-equivariant  inclusion  $\mfg^{\mr{ad}(q_1)}_{j-1}(1) \migiincl \left(\mfg_{j-1} (1)= \right)$ $\mfg^{j-1}/\mfg^j (1)$). 
Hence, (\ref{QH31}) becomes  a {\it split} short exact sequence
\begin{align}  \label{Ex5} 0  \longmigi  \Omega\otimes (\dV_{\mfg, -j})^\vee   \migi   \mfg^{j}_{\DE_{\mbB}}/\mfg^{j+1}_{\DE_{\mbB}}
\xrightarrow{\nabla^{\mr{ad}(j/j+1)}}  \Omega\otimes (\mfg^{j-1}_{\DE_{\mbB}}/ \mfg^{j}_{\DE_{\mbB}})  \longmigi 0
\\
\left(\text{resp.,} \  0 \longmigi    \mfg^{j}_{\DE_{\mbB}}/\mfg^{j+1}_{\DE_{\mbB}} \xrightarrow{\nabla^{\mr{ad}(j/j+1)}}  \Omega\otimes (\mfg^{j-1}_{\DE_{\mbB}}/ \mfg^{j}_{\DE_{\mbB}}) \longmigi \dV_{\mfg, j-1} \longmigi 0\right). \notag
\end{align}

\subsection{Cohomology of the complex $\mcK^\bullet [\nabla^{\mr{ad}}]$}
The following Lemma \ref{y0142} 
can be proved by applying  the above discussion   and
 will be used in the proof of Proposition \ref{y0143} described below.

\ble \label{y0142}
Let $j$ be an integer.
Then, the following assertions hold:
\begin{itemize}
\item[(i)]
Let   $\varsigma^j :  \dV^j_\mfg \migi \Omega\otimes \mfg^j_{\DE_{\mbB}}$ be  
 the restriction 
 of  $\varsigma$ 
(cf. (\ref{inclV})) to the respective $j$-th subsheaves in the filtrations.
Then, the composite
\begin{equation}   f_*( \dV^j_\mfg)  \xrightarrow{f_*(\varsigma^j)} f_*(\Omega\otimes \mfg^j_{\DE_{\mbB}}) \migisurj  \mr{Coker}(f_*(\nabla^{\mr{ad}(j+1)}))  \label{Lemma1}
\end{equation}
 is an isomorphism of $\mcO_S$-modules.
  In particular, we have
\begin{equation}
 f_*( \dV_\mfg)  \cong  \mr{Coker}(f_*(\nabla^{\mr{ad}})).
\label{Lemma2}\end{equation}
\item[(ii)]
Let  $\varsigma^{\veebar}$ be  the filtered morphism
$\mfg_{\DE_{\mbB}} \migi \Omega\otimes \dV_\mfg^{\vee}$ 
defined to be the composite  of  the isomorphism $\mfg_{\DE_{\mbB}} \isom \mfg_{\DE_{\mbB}}^\vee$ induced by the Killing form on $\mfg$ and
the surjection  $\mfg_{\DE_{\mbB}}^\vee \migisurj \Omega\otimes \dV^{\vee}_\mfg$ induced by the dual of $\varsigma$.
Also, let 
$\varsigma^{\veebar, j}: \mfg_{\DE_{\mbB}}^j \migi  \Omega\otimes \dV_\mfg^{\vee, j}$ be  
the restriction  of $\varsigma^\veebar$ to the respective $j$-th subsheaves in the filtrations.
Then, the composite
\begin{align}  \label{Lemma3} 
  \mr{Ker}(\mbR^1f_*(\nabla^{\mr{ad}(j)}))  & \migiincl   \mbR^1 f_*(\mfg^j_{\DE_{\mbB}}) 
\xrightarrow{\mbR^1f_*( \varsigma^{\veebar, j})} \mbR^1f_*(\Omega\otimes \dV_\mfg^{\vee, j})
\end{align}
is an isomorphism of $\mcO_S$-modules.
In particular, we have
\begin{equation}
\label{Lemma4}
 \mr{Ker}(\mbR^1f_*(\nabla^{\mr{ad}})) \cong  \mbR^1f_*(\Omega\otimes \dV_\mfg^\vee) \left(\cong  f_*(\dV_\mfg (-D_\msX))^\vee\right). 
\end{equation}
\end{itemize}
\ele
\bpf
We only consider assertion (ii)  since assertion (i) follows from   a similar argument.

If $j \geq 0$, then $\dV_\mfg^{\vee, j} \left(:= (\dV_\mfg/\dV_\mfg^{-j+1})^\vee \right) = (\dV_\mfg/\dV_\mfg)^\vee=0$,   and the equality $\mr{Ker}(\mbR^1f_*(\nabla^{\mr{ad}(j/j+1)})) = 0$ holds (which implies that $\mr{Ker}(\mbR^1f_*(\nabla^{\mr{ad}(j)})) = 0$) because the  resp'd portion of (\ref{Ex5})  is split.
Hence, it suffices to prove  the assertion for $j \leq -1$.
We shall consider the morphism of sequences
\begin{align} \label{MorSeq}
\vcenter{\xymatrix@C=14pt@R=36pt{
0 \ar[r]& 
 \mbR^1 f_*(\mfg^{j+1}_{\DE_{\mbB}})
\ar[r] \ar[d]^-{ \mbR^1f_*(\nabla^{\mr{ad}(j+1)})}
&\mbR^1 f_*(\mfg^j_{\DE_{\mbB}})  \ar[r] \ar[d]^-{\mbR^1f_*(\nabla^{\mr{ad}(j)})}
& \mbR^1 f_*(\mfg^{j}_{\DE_{\mbB}}/\mfg^{j+1}_{\DE_{\mbB}})
\ar[r] \ar[d]^-{\mbR^1f_*(\nabla^{\mr{ad}(j/j+1)})} & 0
\\
0 \ar[r]& \mbR^1f_*(\Omega\otimes \mfg^j_{\DE_{\mbB}}) \ar[r] &
\mbR^1f_*(\Omega\otimes \mfg^{j-1}_{\DE_{\mbB}}) 
\ar[r]
& \mbR^1f_*(\Omega\otimes (\mfg^{j-1}_{\DE_{\mbB}}/\mfg^j_{\DE_{\mbB}})) \ar[r]& 0.
}}
\end{align}
Since $\mfg^{j}_{\DE_{\mbB}}/\mfg^{j+1}_{\mcE^\dagger_{\mbB}}$ and $\Omega \otimes (\mfg^{j-1}_{\DE_{\mbB}}/\mfg^j_{\DE_{\mbB}})$ are direct sums of finite  copies of $\Omega^{\otimes j}$ (cf. (\ref{QQ0104})), 
we have 
\begin{equation} 
\label{f=0}
f_*(\mfg^{j}_{\DE_{\mbB}}/\mfg^{j+1}_{\DE_{\mbB}}) = f_*(\Omega\otimes (\mfg^{j-1}_{\DE_{\mbB}}/\mfg^j_{\DE_{\mbB}}))   = 0.
\end{equation}
By 
these equalities
together with  the fact that $X/S$ is of relative dimension $1$ (which implies that  $\mbR^2f_*(-) =0$), 
both the upper and lower horizontal sequences in  (\ref{MorSeq}) turns out to be  exact.
Also,  
the right-hand vertical  arrow in (\ref{MorSeq}) 
 is surjective because of   the long exact sequence arising from the non-resp'd portion of (\ref{Ex5}).
By descending induction on $j$, 
one  can verify the surjectivity  of  both the left-hand  and middle  vertical  arrows  in (\ref{MorSeq}).
Thus,   
the snake lemma  applied to (\ref{MorSeq}) shows that
 the natural   sequence
\begin{equation} 0 \migi   \mr{Ker}(\mbR^1f_*(\nabla^{\mr{ad}(j+1)}))\migi  \mr{Ker}(\mbR^1f_*(\nabla^{\mr{ad}(j)})) \migi  \mr{Ker}(\mbR^1f_*(\nabla^{\mr{ad}(j/j+1)})) \migi 0  \label{Ker} \hspace{-3mm}
\end{equation}
is  exact.
Let us consider  the  morphism of short exact sequences
\begin{align} \label{MorSeq2}
\vcenter{\xymatrix@C=9pt@R=36pt{
0 \ar[r]& \mr{Ker}(\mbR^1f_*(\nabla^{\mr{ad}(j+1)})) \ar[d] \ar[r] 
&\mr{Ker}(\mbR^1f_*(\nabla^{\mr{ad}(j)}))  \ar[r] \ar[d]^-{(\ref{Lemma3})} &  \mr{Ker}(\mbR^1f_*(\nabla^{\mr{ad}(j/j+1)}))
\ar[d] \ar[r] & 0
\\
0 \ar[r]& \mbR^1f_*(\Omega\otimes \dV_\mfg^{\vee, j+1})  \ar[r]& \mbR^1f_*(\Omega\otimes \dV_\mfg^{\vee, j})   \ar[r]
& \mbR^1f_*(\Omega\otimes (\dV_{\mfg, -j})^\vee) \ar[r]& 0
}}\hspace{-4mm}
\end{align}
induced by (\ref{Lemma3}),  where  the left-hand vertical arrow is  (\ref{Lemma3}) with $j$ replaced by $j+1$ and the exactness of the lower horizontal sequence follows from  the natural  decomposition $\dV_\mfg^{\vee, j} = \dV_\mfg^{\vee, j+1} \oplus (\dV_{\mfg, -j})^\vee$.
 The right-hand vertical  arrow in 
this  diagram is an isomorphism because it coincides with the inverse of the isomorphism $\mbR^1 f_*(\Omega \otimes (\dV_{\mfg, -j})^\vee) \isom \mr{Ker}(\mbR^1 f_*(\nabla^{\mr{ad}(j/j+1)}))$ arising from the splitting  of 
the non-resp'd portion of  (\ref{Ex5}).
Hence, by descending induction on $j$, we see that  the middle vertical  arrow (\ref{Lemma3})   is   an isomorphism.
This completes the proof of the assertion.
\epf

Now, let us prove 
the following  proposition; the case  
  of 
 $\mfg = \mfs \mfl_2$ 
 can be found  in  ~\cite{Mzk1}, Chap.\,I, Theorem 2.8.

\bpr \label{y0143} 
\begin{itemize}
\item[(i)]
$\mbR^0f_*(\mcK^\bullet[\nabla^\mr{ad}]) \left(\cong f_*(\mr{Ker}(\nabla^\mr{ad}))\right)  =0$.
\item[(ii)]
There exists a natural short exact sequence
\begin{equation}     0 \longmigi f_*(\dV) \longmigi \mbR^1f_*(\mcK^\bullet[\nabla^\mr{ad}]) \longmigi \mbR^1f_*(\Omega\otimes \dV^\vee) \longmigi 0  \label{coh}
\end{equation}
of $\mcO_S$-modules.
In particular, $\mbR^1f_*(\mcK^\bullet[\nabla^\mr{ad}])$ is a vector bundle on $S$ of rank $ (2g-2+r) \cdot \mr{dim}(\mfg)$.
\item[(iii)]
$\mbR^2f_*(\mcK^\bullet[\nabla^\mr{ad}])  \left( \cong \mbR^1f_*(\mr{Coker}(\nabla^\mr{ad}))\right)  =0$.
\end{itemize}
\epr
\bpf
First, we shall prove assertion (i). Consider the following morphism of short exact sequences defined for each $j \in \mbZ$: 
\begin{align} \label{diag4}
\vcenter{\xymatrix@C=16pt@R=36pt{
0 \ar[r]  &\mfg^{j+1}_{\DE_{\mbB}} \ar[r] \ar[d]^-{\nabla^{\mr{ad}(j+1)}}
& \mfg^{j}_{\DE_{\mbB}} \ar[r] \ar[d]^-{\nabla^{\mr{ad}(j)}} 
&\mfg^{j}_{\DE_{\mbB}}/\mfg^{j+1}_{\DE_{\mbB}} \ar[r] \ar[d]^-{\nabla^{\mr{ad}(j/j+1)}}
&0\\
0 \ar[r]& \Omega\otimes  \mfg^{j}_{\DE_{\mbB}}\ar[r] 
&\Omega\otimes  \mfg^{j-1}_{\DE_{\mbB}} \ar[r]& \Omega\otimes (\mfg^{j-1}_{\DE_{\mbB}}/ \mfg^{j}_{\DE_{\mbB}})  \ar[r] & 0.
}}
\end{align}
 The exactness of 
 the resp'd portion of (\ref{Ex5}) implies
 the injectivity of 
   the right-hand vertical arrow in this  diagram for $j\geq 0$.
Hence, by descending induction on $j$, 
we can verify  
the injectivity  of 
$f_*(\nabla^{\mr{ad}(0)}) : f_*(\mfg^{0}_{\DE_{\mbB}})\migi f_*(\Omega\otimes  \mfg^{-1}_{\DE_{\mbB}})$.
Since  both $\mfg^{j}_{\DE_{\mbB}}/\mfg^{j+1}_{\DE_{\mbB}} $ and $ \Omega\otimes (\mfg^{j-1}_{\DE_{\mbB}}/ \mfg^{j}_{\DE_{\mbB}})$ are direct sums of finite copies of $\Omega^{\otimes j}$ (cf.  (\ref{QQ0104})),
we have 
\begin{align}
f_*(\mfg^{0}_{\DE_{\mbB}}) &= f_*(\mfg^{-1}_{\DE_{\mbB}}) = \cdots = f_*(\mfg_{\DE_{\mbB}}), \\
f_*(\Omega\otimes  \mfg^{-1}_{\DE_{\mbB}})  &= f_*(\Omega\otimes  \mfg^{-2}_{\DE_{\mbB}})  =      \cdots  
 = f_*(\Omega\otimes  \mfg_{\DE_{\mbB}}). \notag
\label{eq7}
\end{align}
It follows that
$f_*(\nabla^{\mr{ad}(0)})$ is nothing but $f_*(\nabla^{\mr{ad}})$.
In particular, 
  $f_*(\nabla^\mr{ad})$ turns out to be  injective, i.e.,  $\mr{Ker}(f_*(\nabla^\mr{ad})) = \mbR^0f_*(\mcK^\bullet[\nabla^\mr{ad}]) = 0$.
This completes the proof of assertion (i).

Assertion (ii) follows from  the exactness of  (\ref{HDeq}),   Lemma \ref{y0142},  (i), (ii), and Proposition \ref{y046}, (ii).
Finally, assertion (iii) follows from an argument similar to the proof of  (i) (hence we omit the details). 
\epf


\section{Infinitesimal deformations of a $(\mfg, \hslash)$-oper} \label{y0147}

In this section, we describe the deformation space of a $(\mfg, \hslash)$-oper  using the cohomology of a certain complex associated to a $(\mfg, \hslash)$-oper.

\subsection{The complex $\widetilde{\nabla}^{\mr{ad} (0)}$ associated to a $(\mfg, \hslash)$-oper} \label{y0144}

The $\hslash$-$S^\mr{log}$-connection $\nabla$ induces 
an $f^{-1}(\mcO_S)$-linear morphism $\widetilde{\nabla}^\mr{ad} : \widetilde{\mcT}_{\DE_\mbG^\mr{log}/S^\mr{log}} \migi \Omega \otimes \mfg_{\DE_\mbG}$ (cf. (\ref{circledast})).
By the definition of a $(\mfg, \hslash)$-oper,  the image of the composite  $ \widetilde{\nabla}^{\mr{ad}} \circ \widetilde{\iota}_{\mfg/\mfb} : \widetilde{\mcT}_{\DE_\mbB^\mr{log}/S^\mr{log}} \migi \Omega \otimes \mfg_{\DE_\mbG}$  (cf. (\ref{extension}) for the definition of $\widetilde{\iota}_{\mfg/\mfb}$) is contained 
  in 
$\Omega\otimes  \mfg^{-1}_{\DE_{\mbB}}$.
Hence, $\widetilde{\nabla}^\mr{ad}$ induces 
an $f^{-1}(\mcO_S)$-linear morphism
\index{$\nabla^\mr{ad}$@--- $\widetilde{\nabla}^{\mr{ad} (0)}$}
\begin{equation}  
 \widetilde{\nabla}^{\mr{ad} (0)} :  
\widetilde{\mcT}_{\DE_{\mbB}^{\mr{log}}/S^\mr{log}} \migi \Omega\otimes  \mfg^{-1}_{\DE_{\mbB}}.
\label{Adcircledast}
\end{equation}
The following lemma will be used in the subsequent discussion.

\ble \label{y0146}
The following equalities hold:
\begin{align}
\mbR^2f_*(\mcK^\bullet [\nabla^{\mr{ad}(0)}] ) = \mbR^2f_*(\mcK^\bullet [\widetilde{\nabla}^{\mr{ad}(0)}] ) = 0.
\end{align}
 \ele
\begin{proof}
By  applying   the Hodge to de Rham spectral sequence (cf. (\ref{HDsp}))
to the complex $\mcK^\bullet [\nabla^{\mr{ad}(0)}] $ (resp., $\mcK^\bullet [\widetilde{\nabla}^{\mr{ad}(0)}]$), we see that 
$\mbR^2f_*(\mcK^\bullet [\nabla^{\mr{ad}(0)}] )$ (resp.,  $\mbR^2f_*(\mcK^\bullet [\widetilde{\nabla}^{\mr{ad}(0)}] )$) is isomorphic to  the cokernel of the morphism
\begin{align} \label{Mor66}
 \mbR^1f_*(\nabla^{\mr{ad}(0)}) : \mbR^1f_*(\mfg_{\DE_{\mbB}}^{0}) \migi \mbR^1f_*(\Omega\otimes \mfg_{\DE_{\mbB}}^{-1}) \hspace{7mm} \\
 \left(\text{resp.,} \   \mbR^1f_*(\widetilde{\nabla}^{\mr{ad}(0)}) : \mbR^1f_*(\widetilde{\mcT}_{\DE^{\mr{log}}_{\mbB}}) \migi \mbR^1f_*(\Omega\otimes \mfg_{\DE_{\mbB}}^{-1}) \right). \notag
\end{align}
Since $ \mbR^1f_*(\nabla^{\mr{ad}(0)})$ coincides with the compose of  
$\mbR^1f_*(\widetilde{\nabla}^{\mr{ad}(0)})$ and  the morphism 
$\mbR^1f_*(\mfg_{\DE_{\mbB}}^{0}) \migi \mbR^1f_*(\widetilde{\mcT}_{\DE^{\mr{log}}_{\mbB}})$ induced by the inclusion $\mfg_{\DE_{\mbB}}^{0} \migiincl \widetilde{\mcT}_{\DE^{\mr{log}}_{\mbB}}$,
the problem can be  reduced to proving the surjectivity of $\mbR^1f_*(\nabla^{\mr{ad}(0)})$.
To prove this,  we observe  that 
$\mbR^1f_*(\Omega\otimes \mfg_{\DE_{\mbB}}^{1}) = 0$, which follows   
 from (\ref{QQ0104}).
 This implies the surjectivity of  
$\mbR^1f_*(\nabla^{\mr{ad}(2)})$.
On the other hand, 
the exactness of   (\ref{Ex5}) for  $j =0, 1$  shows 
\begin{equation}
\mr{Coker}(\mbR^1f_*(\nabla^{\mr{ad}(0/1)})) = \mr{Coker}(\mbR^1f_*(\nabla^{\mr{ad}(1/2)})) = 0.
\end{equation}
By these facts together with  
 the morphism of short exact sequences   (\ref{MorSeq}) for $j=0, 1$,  we can verify   the desired surjectivity.
 This  completes
 the proof of this lemma.
\end{proof}


Since  $\widetilde{\nabla}^{\mr{ad}(0)}$ restricts to $\nabla^{\mr{ad}(0)}$, 
there exists a short exact sequence of complexes
\begin{equation} \label{short7}
0 \longmigi \mcK^\bullet [\nabla^{\mr{ad}(0)}] \longmigi \mcK^\bullet [\widetilde{\nabla}^{\mr{ad}(0)}] \longmigi  \mcT [0] \longmigi 0
\end{equation}
 obtained   in the same way as (\ref{QW7034}).
By applying   the functor $\mbR^1f_*(-)$ to 
this sequence, we
obtain  a sequence of $\mcO_S$-modules
\begin{equation} \label{short8}
0 \longmigi \mbR^1f_*(\mcK^\bullet [\nabla^{\mr{ad}(0)}]) \longmigi \mbR^1f_*(\mcK^\bullet [\widetilde{\nabla}^{\mr{ad}(0)}]) \longmigi  \mbR^1f_*(\mcT) \longmigi 0.
\end{equation}

This sequence  turns out to be  exact because of   Lemma \ref{y0146} and the equality $f_*(\mcT) = 0$.
On the other hand, 
consider the short exact sequence
\begin{equation}
0 \migi \mr{Coker}(f_*(\nabla^{\mr{ad}(0)})) \xrightarrow{{' e}_\sharp[\nabla^{\mr{ad}(0)}]} \mbR^1f_*(\mcK^\bullet [\nabla^{\mr{ad}(0)}]) \xrightarrow{{' e}_\flat [\nabla^{\mr{ad}(0)}]} \mr{Ker}(\mbR^1f_*(\nabla^{\mr{ad}(0)})) \migi 0
\end{equation}
i.e., the sequence  (\ref{HDeq}) for  $\nabla = \nabla^{\mr{ad}(0)}$.
The equality $\mr{Ker}(\mbR^1f_*(\nabla^{\mr{ad}(0)})) =0$ holds by Lemma \ref{y0142}, (ii), 
so the composite
\begin{equation} \label{isom9}
 f_*(\dV_\mfg) \xrightarrow{f_*(\varsigma)} \mr{Coker}(f_*(\nabla^{\mr{ad}(0)})) \xrightarrow{{' e}_\sharp[\nabla^{\mr{ad}(0)}]} \mbR^1f_*(\mcK^\bullet [\nabla^{\mr{ad}(0)}]) 
\end{equation}
is an isomorphism.
By combining (\ref{short8}) and  (\ref{isom9}), we obtain a  sequence
\begin{equation}   0 \longmigi f_*(\dV_\mfg) \longmigi \mbR^1f_*(\mcK^\bullet [\widetilde{\nabla}^{\mr{ad}(0)}])  \longmigi \mbR^1f_*(\mcT) \longmigi 0\label{tangent}
\end{equation}
of $\mcO_S$-modules.

\subsection{Deformation of $(\mfg, \hslash)$-opers}

Now, by  a {\bf first-order deformation of $\msE^\spadesuit$} (resp., {\bf $(\msX, \msE^\spadesuit)$}), we shall mean  
a $(\mfg, \hslash)$-oper $\msE_{\epsilon'}^\spadesuit$ on $\msX_\epsilon$ (resp., a pair $(\msX_{\epsilon'}, \msE^\spadesuit_{\epsilon'})$ of a first-order deformation $\msX_{\epsilon'}$ of $\msX$ and a $(\mfg, \hslash)$-oper 
$\msE^\spadesuit_{\epsilon'}$ on $\msX_{\epsilon'}$)
 inducing  $\msE^\spadesuit$ (resp., $(\msX, \msE^\spadesuit)$) via  reduction modulo $\epsilon$.
 Suppose that $S$ is affine.
Then, since  $\mfO \mfp_{\mfg, \hslash, \msX}$ forms a relative affine space modeled on  $\mbV (f_*({^\dagger}\mcV_{\mfg}))$  (cf. (\ref{QR094})),
there is a canonical bijective correspondence
   \begin{align} \label{QW7083}
\begin{pmatrix}
\text{the set of isomorphism classes of} \\
\text{first-oder deformations of $\msE^\spadesuit$} 
\end{pmatrix}
\cong \Gamma (X,  \dV_\mfg).
\end{align}
To be precise, this bijection is constructed  in such a way that each element $R \in \Gamma (X, \dV_\mfg)$ corresponds to the $(\mfg, \hslash)$-oper $(\DE_{\mbB, \epsilon}, \mr{pr}_S^*(\nabla) + \epsilon \cdot R)$.
Moreover, 
by considering  the explicit description of  $\mbH^1 (X, \mcK^\bullet [\widetilde{\nabla}^{\mr{ad} (0)}])$ in terms of the \v{C}ech double complex associated to  $\mcK^\bullet [\widetilde{\nabla}^{\mr{ad}}]$,  
we obtain,  in the same manner as (\ref{QW6083}), the following bijection:
   \begin{align} \label{QW7101}
\begin{pmatrix}
\text{the set of isomorphism classes of} \\
\text{first-oder deformations of $(\msX, \msE^\spadesuit)$} 
\end{pmatrix}
\cong \mbH^1 (X, \mcK^\bullet [\widetilde{\nabla}^{\mr{ad} (0)}]).
\end{align}
It is verified that the maps $\Gamma (X, \dV_\mfg) \migi \mbH^1 (X, \mcK^\bullet [\widetilde{\nabla}^{\mr{ad} (0)}])$ and $\mbH^1 (X, \mcK^\bullet [\widetilde{\nabla}^{\mr{ad} (0)}]) \migi H^1 (X, \mcT)$ correspond, via  (\ref{QW6005}), (\ref{QW7083}),  and (\ref{QW7101}), to the assignments of deformations $\msE_{\epsilon'}^\spadesuit \mapsto (\msX_\epsilon, \msE_{\epsilon'}^\spadesuit)$ and $(\msX_{\epsilon'}, \msE_{\epsilon'}^\spadesuit) \mapsto \msX_{\epsilon'}$ respectively.
 Therefore, we obtain the following  assertion (even when $S$ is not necessarily affine).

\bpr \label{y0148}
Suppose that $\hslash \in k$.
Let 
$s^{\mr{log}}_{\mfo \mfp} : S^\mr{log} \migi \mfO \mfp_{\mfg, \hslash, g,r}^\mr{log}$ (cf. (\ref{univ1}) 
\index{$s^\mr{log}_\mfm$, $s^{\mr{log}}_{\mfo \mfp}$, ${^c}s^{\mr{log}}_{\mfo \mfp}$, ${^p}s^{\mr{log}}_{\mfo \mfp}$}
for the definition of $\mfO \mfp_{\mfg, \hslash, g,r}$ and (\ref{univ1p}) for its log structure) be the  strict morphism
 whose underlying $S$-rational point  of  $\mfO \mfp_{\mfg, \hslash, g,r}$ classifyies    the pair  $(\msX, \msE^\spadesuit)$.
Denote by 
\begin{equation}   \label{isom2}
 \mr{KS}_{\msE^\spadesuit /\msX} :
\left(s^{*}_{\mfo \mfp}(\mcT_{\mfO \mfp_{\mfg, \hslash,g,r}/\overline{\mfM}_{g,r}}) = \right) \ s_{\mfo \mfp}^{*}(\mcT_{\mfO \mfp^\mr{log}_{\mfg, \hslash,g,r}/\overline{\mfM}^\mr{log}_{g,r}})   \isom f_*(\dV_\mfg)   
\index{${^{(c)}}\mr{KS}_{\msE^\spadesuit /\msX}$, ${^{(c)}}\mr{KS}_{\msE^\spadesuit}$, ${^p}\mr{KS}_{\msE^\spadesuit}$}
\end{equation}
 the canonical isomorphism arising from the structure of relative affine space on $\mfO \mfp_{\mfg, \hslash,\msX}$ modeled on $\mbV (f_*(\dV_\mfg))$ (cf. Proposition \ref{y051}).
Then, there exists a canonical isomorphism
\begin{equation}    \mr{KS}_{\msE^\spadesuit} : 
s^{*}_{\mfo \mfp}(\mcT_{\mfO \mfp^\mr{log}_{\mfg, \hslash, g,r}/k}) \isom \mbR^1f_*(\mcK^\bullet [\widetilde{\nabla}^{\mr{ad}(0)}])    \label{isom3333}
\end{equation}
which fits into  the following  diagram of $\mcO_S$-modules:
\begin{align} \label{OpCoh}
\vcenter{\xymatrix@C=16pt@R=36pt{
0 \ar[r]
& s^{*}_{\mfo \mfp}(\mcT_{\mfO \mfp^\mr{log}_{\mfg, \hslash,g,r}/\overline{\mfM}^\mr{log}_{g,r}}) \ar[r] \ar[d]_-{\wr}^-{\mr{KS}_{\msE^\spadesuit /\msX}}
& s^{ *}_{\mfo \mfp}(\mcT_{\mfO \mfp^\mr{log}_{\mfg, \hslash, g,r}/k}) \ar[r] \ar[d]_-{\wr}^-{\mr{KS}_{\msE^\spadesuit} }
&s^{*}_{\mfm}(\mcT_{\overline{\mfM}^\mr{log}_{g,r}/k}) \ar[r] \ar[d]_-{\wr}^-{\mr{KS}_{\msX}} &0\\
0 \ar[r]&  f_*(\dV_\mfg) \ar[r] &\mbR^1f_*(\mcK^\bullet [\widetilde{\nabla}^{\mr{ad}(0)}]) \ar[r]
& \mbR^1f_*(\mcT) \ar[r] & 0, 
}}
\end{align}
where the lower horizontal sequence is (\ref{tangent}) and the upper horizontal  sequence is
obtained by differentiating  the projection $\mfO \mfp^\mr{log}_{\mfg, \hslash, g,r} \migisurj \overline{\mfM}_{g,r}^\mr{log}$.
In particular,  (\ref{tangent}) is exact.
 \epr

\subsection{The differential of the adjoint quotient $\chi$} \label{y0149}

In what follows, we consider deformations of a $(\mfg, \hslash)$-oper of prescribed radii.
Suppose further that $r>0$ and  the  $(\mfg, \hslash)$-oper  $\msE^{\spadesuit}$ is of radii $\rho$ for some $\rho \in \mfc^{\times r} (k)$.
Write $\mcT\mfg$ and  $\mcT\mfc$  for  the total space of the tangent bundles over $k$
   of $\mfg$ and  $\mfc$  respectively.
   Also, write 
 $d\chi : \mcT\mfg \migi  \mcT\mfc$ for the differential of the morphism $\chi : \mfg \migi \mfc$ (cf. (\ref{chi})).
Let us fix $i \in \{ 1, \cdots, r \}$.
By using the differential $d\mr{triv}_{\mbG, \sigma_i}: \mcO_S \otimes_k \mfg \isom \sigma^{*}_i (\mfg_{\DE_\mbG})$  of $\mr{triv}_{\mbG, \sigma_i} : S \times_k \mbG \isom \sigma^*_i (\DE_{\mbG})$ (cf. (\ref{QW4421})), we shall regard the monodromy $\mu_i := \mu_i^\nabla$ at each marked point $\sigma_i$ as a global section  of $\mcO_S \otimes_k \mfg$, or equivalently, an $S$-rational point of $\mfg$. 
The pull-back of $d \chi$ by the composite $\rho_i := \chi \circ \mu_i : S \migi \mfc$ determines
a morphism
 \begin{equation}
 d_{\mu_i}\chi : \left(\mbV (\mcO_S \otimes_k \mfg) = \right) S \times_{\mu_i, \mfg} \mcT \mfg  \migi  S \times_{\rho_i, \mfc}\mcT \mfc 
 \label{diffmu}
 \end{equation}
 of relative affine spaces over $S$.
 We shall define 
 \begin{equation}  \mcK^1[{^c}\widetilde{\nabla}^{\mr{ad}(0)}] \label{cK}
 \index{$\mcK^1[{^c}\widetilde{\nabla}^{\mr{ad}(0)}]$}
 \end{equation}
  to be the $\mcO_X$-submodule 
 of $\Omega\otimes \mfg^{-1}_{\DE_{\mbB}}$  consisting of  local sections $v$ such that, for every $i =1, \cdots, r$,  its restriction  $\sigma^*_i(v)$ to $\mbV (\sigma_i^*(\Omega\otimes \mfg^{-1}_{\DE_{\mbB}})) $ 
 is mapped to the zero section via the composite
 \begin{align}  \label{cLmap}
 \mbV (\sigma_i^*(\Omega\otimes \mfg^{-1}_{\DE_{\mbB}})) 
 & \migiincl
 \mbV(\sigma_i^*(\Omega\otimes \mfg_{\DE_{\mbG}}))  \left(= \mbV(\sigma_i^*(\Omega) \otimes  \sigma_i^*(\mfg_{\DE_{\mbG}})) \right)\\
 &\xrightarrow{\mr{Res}_{\msX, i} \otimes \mr{id}} \mbV (\sigma^*_i(\mfg_{\DE_{\mbB}})) 
 \notag \\
 & \xrightarrow{d\mr{triv}_{\mbG, \sigma_i}^{-1}} \mbV(\mcO_S \otimes_k \mfg)  \notag \\
 &\xrightarrow{d_{\mu_i}\chi} S \times_{\rho_i, \mfc}  \mcT \mfc. \notag
 \end{align}

\ble \label{y0150} 
 The image of $\widetilde{\nabla}^{\mr{ad}(0)}$ (cf. (\ref{Adcircledast})) lies  in $ \mcK^1[{^c}\widetilde{\nabla}^{\mr{ad}(0)}]$ $\left(\subseteq \Omega\otimes \mfg^{-1}_{\DE_{\mbB}}\right)$.
 \ele
\begin{proof}
For simplicity, we write $\mcE := \DE_\mbG$ and $\nabla_\epsilon := \mr{pr}_S^*(\nabla)$.
Let us fix $i \in \{1, \cdots, r\}$.
Also, let $U$  an open subscheme  of $X$ with $U \cap \mr{Im}(\sigma_i) \neq \emptyset$ and 
  $v$ an element of $\Gamma (U,  \widetilde{\mcT}_{\DE_{\mbB}^{\mr{log}}/S^\mr{log}})$.
The problem is to  prove that $\widetilde{\nabla}^{\mr{ad}(0)} (v)$ belongs to $\Gamma (U,  \mcK^1[{^c} \widetilde{\nabla}^{\mr{ad}(0)}])$.
To  this end, we may suppose, after changing the base  $S$ to its open subscheme if necessary,
that $\mr{Im}(\sigma_i) \subseteq U$.
If $\eta^v_\mcE$ denotes  the automorphism of $\mcE_{\epsilon} |_U$
obtained from $v$ as in Remark  \ref{QW9004} (cf. (\ref{QH40})), then
we have $\eta^{v*}_\mcE (\nabla_\epsilon) = \nabla_\epsilon  -\epsilon \cdot  \widetilde{\nabla}^{\mr{ad}(0)}(v)$ (cf. (\ref{QW8600})).
 This  implies  
\begin{align}
\mu_i^{\eta^{v*}_\mcE (\nabla_\epsilon)} = \mr{pr}_S^*(\mu_i) - \epsilon \cdot  \sigma_i^*(\widetilde{\nabla}^{\mr{ad}(0)}(v)),
\end{align}
where $\sigma_i^*(\widetilde{\nabla}^{\mr{ad}(0)}(v))$ is considered as an element of $\Gamma (S, \sigma_i^*(\mfg_\mcE))$ via $\mr{Res}_{\msX, i} \otimes \mr{id} : \sigma_i^*(\Omega \otimes \mfg_\mcE) \isom \sigma_i^*(\mfg_\mcE)$.
This means that  the tangent vector $[\mu_i^{\eta^{v*}_\mcE (\nabla_\epsilon)}]$ in $S \times_{\mu_i, \mfg} \mcT \mfg$
determined by $\mu_i^{\eta^{v*}_\mcE (\nabla_\epsilon)}$ is given by $\sigma_i^*(\widetilde{\nabla}^{\mr{ad}(0)}(v))$.
On the other hand, the automorphism $\eta_\mcE^{v}$ restricts to   an isomorphism $(\sigma^{*}_{i}(\mcE_\epsilon), \mu_i^{\eta^{v*}_\mcE (\nabla_\epsilon)}) \isom (\sigma^{*}_{i}(\mcE_\epsilon), \mr{pr}_S^*(\mu_i))$ of objects in $[\mfg/\mbG](S_\epsilon)$.
This isomorphism induces  an equality  $\rho_i^{\eta^{v*}_\mcE (\nabla_\epsilon)} = \rho_i$,  and hence, $[\mu_i^{\eta^{v*}_\mcE (\nabla_\epsilon)}]$ is mapped, via $d_{\mu_i} \chi$,  to the zero element of $S \times_{\rho_i, \mfc}\mcT \mfc$.
It follows that 
$\sigma_i^*(\widetilde{\nabla}^{\mr{ad}(0)}(v))$ is contained in the kernel of (\ref{cLmap}), i.e., 
$\widetilde{\nabla}^{\mr{ad}(0)}(v) \in \Gamma (U,  \mcK^1[{^c} \widetilde{\nabla}^{\mr{ad}(0)}])$.
This completes the proof of the assertion.
\end{proof}

\subsection{Deformation of $(\mfg, \hslash)$-opers of prescribed radii} \label{QG1028}

By Lemma \ref{y0150} above,
 $\widetilde{\nabla}^{\mr{ad}(0)}$ induces  an $f^{-1}(\mcO_S)$-linear morphism
  \index{$\nabla^\mr{ad}$@--- ${^c}\widetilde{\nabla}^{\mr{ad}(0)}$}
\begin{equation} \label{903}  
{^c}\widetilde{\nabla}^{\mr{ad}(0)} :
  \widetilde{\mcT}_{\DE_{\mbB}^{\mr{log}}/S^\mr{log}} \migi \mcK^1[{^c}\widetilde{\nabla}^{\mr{ad}(0)}].
\end{equation}
The sheaf  $\dV_\mfg (-D_\msX)$ lies in $\mcK^1[{^c}\widetilde{\nabla}^{\mr{ad}(0)}]$,
 so
(\ref{short7}) restricts to a sequence  of complexes
\begin{equation}
0 \longmigi  \dV_\mfg (-D_\msX) [-1] \longmigi \mcK^\bullet [{^c}\widetilde{\nabla}^{\mr{ad}(0)} ] \longmigi \mcT [0] \longmigi 0.
\end{equation}
By applying  the functor $\mbR^1f_*(-)$ to this sequence,  we obtain a sequence of $\mcO_S$-modules
\begin{equation}   0 \migi f_*(\dV_\mfg (-D_\msX)) \migi \mbR^1f_*(\mcK^\bullet [{^c}\widetilde{\nabla}^{\mr{ad}(0)}])  \migi \mbR^1f_*(\mcT) \migi 0,
\label{ctangent}
\end{equation}
which is compatible with (\ref{tangent}) via the inclusions 
$\dV_\mfg (-D_\msX) \migiincl \dV_\mfg$,  $\mcK^1[{^c}\widetilde{\nabla}^{\mr{ad}(0)}] \migiincl \Omega \otimes \mfg_{\DE_\mbB}^{-1}$.

\bpr \label{y0151}
 Let 
 ${^c s}_{\mfo \mfp}^{\mr{log}}: S^\mr{log} \migi \mfO \mfp^\mr{log}_{\mfg,\hslash,  \rho, g,r}$
\index{$s^\mr{log}_\mfm$, $s^{\mr{log}}_{\mfo \mfp}$, ${^c}s^{\mr{log}}_{\mfo \mfp}$, ${^p}s^{\mr{log}}_{\mfo \mfp}$}
  be  the
 strict morphism whose underlying $S$-rational point of $\mfO \mfp_{\mfg,\hslash,  \rho, g,r}$
   classifies  the pair  $(\msX, \msE^\spadesuit)$.
Denote by 
\begin{equation}  
 {^c}\mr{KS}_{\msE^\spadesuit /\msX} :
\left( {^c s}_{\mfo \mfp}^{*}(\mcT_{\mfO \mfp_{\mfg, \hslash,\rho, g,r}/\overline{\mfM}_{g,r}})  =\right)    {^c s}_{\mfo \mfp}^{*}(\mcT_{\mfO \mfp^\mr{log}_{\mfg, \hslash,\rho, g,r}/\overline{\mfM}^\mr{log}_{g,r}}) 
  \isom  f_*(\dV_\mfg (-D_\msX)) 
   \end{equation}
    the canonical   isomorphism arising from the structure of relative affine space on $\mfO \mfp_{\mfg, \hslash, \rho,\msX}$ modeled on $\mbV (f_*( \dV_\mfg (-D_\msX)))$ (cf. Theorem \ref{QQ1060}).
Then, there exists a canonical isomorphism
\begin{equation}
\label{904}
    {^c}\mr{KS}_{\msE^\spadesuit} : 
    {^c s}_{\mfo \mfp}^{*}(\mcT_{\mfO \mfp^\mr{log}_{\mfg, \hslash, \rho, g,r}/k}) 
    \isom \mbR^1f_*(\mcK^\bullet [{^c}\widetilde{\nabla}^{\mr{ad}(0)} ])
    \index{${^{(c)}}\mr{KS}_{\msE^\spadesuit /\msX}$, ${^{(c)}}\mr{KS}_{\msE^\spadesuit}$, ${^p}\mr{KS}_{\msE^\spadesuit}$}
     \end{equation}
which fits into  the following  diagram of $\mcO_S$-modules:
\begin{align} \label{QW6089}
\vcenter{\xymatrix@C=16pt@R=36pt{
0 \ar[r]
& {^c s}^{*}_{\mfo \mfp}(\mcT_{\mfO \mfp^\mr{log}_{\mfg, \hslash, \rho,g,r}/\overline{\mfM}^\mr{log}_{g,r}})  \ar[r] \ar[d]_-{\wr}^-{{^c}\mr{KS}_{\mcE^\spadesuit /\msX}}
& {^c s}_{\mfo \mfp}^{*}(\mcT_{\mfO \mfp^\mr{log}_{\mfg, \hslash, \rho, g,r}/k})    \ar[r] \ar[d]_-{\wr}^-{{^c}\mr{KS}_{\mcE^\spadesuit}}
& s^{*}_{\mfm}(\mcT_{\overline{\mfM}^\mr{log}_{g,r}/k}) \ar[r] \ar[d]_-{\wr}^-{\mr{KS}_{\msX}}
& 0\\
0 \ar[r]
& f_*(\dV_\mfg (-D_\msX)) \ar[r]
& \mbR^1f_*(\mcK^\bullet [{^c} \widetilde{\nabla}^{\mr{ad}(0)}]) \ar[r]
&\mbR^1f_*(\mcT) \ar[r]& 0, 
}}
\end{align}
where the lower horizontal sequence is (\ref{ctangent}) and the upper horizontal    sequence is 
obtained by differentiating the projection $\mfO \mfp^\mr{log}_{\mfg, \hslash, g,r, \rho} \migi \overline{\mfM}_{g,r}^\mr{log}$.
In particular,    (\ref{ctangent}) is exact. 
\epr
\begin{proof}
The assertion follows from the construction of (\ref{isom3333}) resulting from  Proposition \ref{y0148}  and the definition of $\mcK^1 [{^c}\widetilde{\nabla}^{\mr{ad}(0)}]$.
\end{proof}

\section{Flat vector bundles in positive characteristic} \label{y0155}

We shall study basic properties on flat bundles in positive characteristic  to   compute  the deformation space of a dormant $(\mfg, \hslash)$-oper.
One of the results in this section is  the nondegeneracy of the pairing on the de Rham cohomology of a flat bundle  with vanishing $p$-curvature.
Hereafter, 
 suppose that $k$ is of characteristic $p > 0$, i.e., $p > 2  h_\mbG$ is satisfied.
(But, the results in this section hold even for a general smooth algebraic group $\mbG$.)

\subsection{The Cartier operator of  a flat  vector bundle} \label{QG1034}

First, we shall recall  the  Cartier operator associated to a flat  bundle.
Let us  write $\Omega_{X^{(1)\mr{log}}/S^\mr{log}} := \Omega^{(1)}$ for simplicity.
Let $(\mcV, \nabla)$ be a flat  bundle on $X^\mr{log}/S^\mr{log}$.
 Since $\nabla$ preserves the $F^{-1}_{X/S}(\mcO_{X^{(1)}})$-actions, 
 both $\mr{Ker}(\nabla)$ and $\mr{Coker}(\nabla)$ may  be regarded as $\mcO_{X^{(1)}}$-modules  via the underlying homeomorphism of $F_{X/S}$.

Now, recall from ~\cite{Og}, Chap.\,I, Proposition 1.2.4, that the {\bf Cartier operator}   associated to 
$(\mcV, \nabla)$ is, by definition, a unique $\mcO_{X^{(1)}}$-linear morphism
\begin{equation}
C^{(\mcV, \nabla)} : F_{X/S*}(\Omega \otimes \mcV) \migi \Omega^{(1)} \otimes F_{X/S*}(\mcV)
\index{$C^{(\mcV, \nabla)}$, Cartier operator associated to $(\mcV, \nabla)$}
\label{CartierF}
\end{equation}
satisfying the following condition:
 for any local sections  $\partial \in \mcT$ and $a \in \Omega\otimes \mcV$, the equality 
\begin{equation}
\langle (\mr{id}_X \times F_S)^{-1}(\partial), C^{(\mcV, \nabla)} (a)\rangle = \langle \partial^{[p]}, a\rangle - \nabla_{\partial}^{p-1}(\langle \partial, a\rangle)
\label{Cartier2F}
\end{equation}
holds, where 
$\langle -, -\rangle$'s
  denote  the pairings induced by the natural pairing  $\mcT \times \Omega \migi \mcO_X$.
Let us  consider   the  composite  isomorphism
\begin{align} \label{projformula1}
\Omega^{(1)}\otimes F_{X/S*}(\mcV)  \isom 
F_{X/S*}(F_{X/S}^*(\Omega^{(1)}) \otimes \mcV) \isom F_{X/S*}(F_X^*(\Omega) \otimes \mcV), 
\end{align}
where  the first isomorphism follows from the projection formula.
By using this composite, 
we occasionally regard $C^{(\mcV, \nabla)}$ as a morphism of sheaves $\Omega \otimes \mcV \migi F_X^*(\Omega) \otimes \mcV$.

Moreover, since $(\mcV, \nabla)$ has vanishing $p$-curvature,
it follows from ~\cite{Og3}, Theorem 3.1.1 (and the discussion following ~\cite{Og},   Chap.\,I, Proposition 1.2.4), that 
there exists an isomorphism 
\begin{align} \label{Qw8345}
\overline{C}^{(\mcV, \nabla)} : \mr{Coker} (\nabla) \isom \Omega^{(1)} \otimes \mr{Ker}(\nabla)
\end{align}
which makes the following square diagram commute:
\begin{align}\label{Cartier3F}
\vcenter{\xymatrix@C=46pt@R=36pt{
F_{X/S*}(\Omega \otimes \mcV) \ar[r]^-{C^{(\mcV, \nabla)}} 
\ar[d]_-{\mr{quotient}}
 & \Omega_{X^{(1)}}\otimes F_{X/S*}(\mcV) \\
\mr{Coker}(\nabla) \ar[r]_-{\overline{C}^{(\mcV, \nabla)}}& \Omega^{(1)} \otimes \mr{Ker}(\nabla). \ar[u]_-{\mr{inclusion}}
}}
\end{align}

\begin{rema} \label{QH45}
Recall from the comment following ~\cite{KaKa}, \S\,4,  Definition (4.8), and ~\cite{Og}, Chap.\,I, Remark 1.2.3,   that  $X^\mr{log}/S^\mr{log}$
  is {\it of Cartier type}, so the notion of the exact relative Frobenius map in the statement of {\it loc.\,cit.} coincides with the usual relative Frobenius morphism $F_{X/S}$.
\end{rema}

\subsection{Deformation of flat  bundles with ${^p}\psi=0$} \label{y0156}

We shall  prove the following proposition. (For the case where $\mbG = \mbG_m$, i.e., the case of line bundles, we refer to ~\cite{Katz2}, (7.2.2), Proposition.)

\bpr \label{y0157}
Let $(\mcE, \nabla)$ be a flat  $\mbG$-bundle over $X^\mr{log}/S^\mr{log}$  with ${^p}\psi^\nabla =0$, and write $\nabla_\epsilon := \mr{pr}_S^*(\nabla)$.
Also, let $w$ be an element of  
$\Gamma (X, \Omega \otimes \mfg_\mcE)$, considered as 
an $\mcO_X$-linear morphism
$\mcT\migi \mfg_\mcE \left( \subseteq \widetilde{\mcT}_{\mcE^\mr{log}/S^\mr{log}}\right)$.
Then, the $p$-curvature 
 of  the $S_\epsilon^\mr{log}$-connection $\nabla_\epsilon + \epsilon \cdot w$ on $\mcE_\epsilon$
satisfies 
\begin{equation} \label{910}
 {^p \psi}^{\nabla_\epsilon + \epsilon \cdot w} 
  = - \epsilon \cdot C^{(\mfg_\mcE, \nabla^\mr{ad})} (w) \left(\in \Gamma (X_\epsilon, \Omega^{\otimes p}_{X^\mr{log}_\epsilon/S^\mr{log}_\epsilon} \otimes \mfg_{\mcE_\epsilon}) \right).
\end{equation}
In particular, if the $p$-curvature of $\nabla_\epsilon + \epsilon \cdot w$  vanishes, then $w$ may be expressed, Zariski locally on $X$, as
$w = \nabla^\mr{ad}(v)$ for some local section $v$ of $\mfg_\mcE$. 
  \epr
\begin{proof}
Let us take an open subscheme $U$ of $X$ on which 
there exists 
a section $x \in \Gamma (U, M_U)$ with  $\mcO_U \cdot d \mr{log} (x) = \Omega |_U$.
The dual $\partial  \in \Gamma (U, \mcT)$
 of  $d \mr{log}(x)$ satisfies 
$\partial^{[p]} = \partial$ (cf. ~\cite{Og},  Remark 1.2.2).
 Let us observe 
 the following  sequence of equalities   in the  enveloping algebra of the Lie algebroid $(\widetilde{\mcT}_{\mcE^{\mr{log}}/S^\mr{log}}, d^\mr{log}_{\mcE})$ (cf. ~\cite{BF}, Appendix, A.3):
\begin{align}  \label{xi1}
 & \  \ \  \  \langle F_{U_\epsilon}^{-1}(\mr{pr}^{-1}_U(\partial)), -\epsilon \cdot  C^{(\mfg_{\mcE}, \nabla^\mr{ad})}(w) \rangle \\
  & =  -\epsilon \cdot  \langle  (\mr{id}_U \times F_S)^{-1}(\partial), C^{(\mfg_{\mcE}, \nabla^\mr{ad})}(w)\rangle \notag  \\
 & =   - \epsilon \cdot  \langle \partial^{[p]}, w\rangle + \epsilon \cdot  (\nabla_{\partial}^\mr{ad})^{p-1}(\langle \partial, w  \rangle)  \notag \\
&=   - \epsilon \cdot  \langle \partial, w \rangle +  \epsilon \cdot  \mr{ad}(\nabla (\partial))^{p-1}(\langle \partial, w \rangle) \notag \\
& = - \epsilon \cdot  \langle \partial, w\rangle + \epsilon \cdot   \sum_{i=0}^{p-1} (-1)^i \cdot \begin{pmatrix} p-1 \\ i \end{pmatrix} \cdot \nabla (\partial)^{p-1 -i} \cdot \langle \partial, w \rangle  \cdot \nabla (\partial)^{i} \notag \\ 
& =    - \epsilon \cdot  \langle \partial, w\rangle +  \epsilon \cdot   \sum_{i=0}^{p-1} \nabla (\partial)^{p-1 -i} \cdot \langle  \partial, w \rangle  \cdot \nabla (\partial)^{i}, \notag
\end{align}
where the second  equality follows from  (\ref{Cartier2F}).
On the other hand,  we have 
\begin{align}  \label{xi2}
 & \  \langle F_{U_\epsilon}^{-1}(\mr{pr}^{-1}_U(\partial)), \psi^{\nabla_\epsilon + \epsilon \cdot w} \rangle  \\
= & \  (\nabla_\epsilon + \epsilon \cdot w)(\mr{pr}_U^{-1}(\partial))^{[p]}-(\nabla_\epsilon + \epsilon \cdot w)(\mr{pr}_U^{-1}(\partial)^{[p]}) \notag \\
= &  \  (\mr{pr}_U^{-1}(\nabla (\partial)) + \epsilon \cdot \langle \partial, w \rangle)^{[p]} -  \mr{pr}_U^{-1}(\nabla (\partial^{[p]})) - \epsilon \cdot \langle \partial^{[p]}, w \rangle\notag \\
=  & \ \epsilon \cdot \sum_{i=0}^{p-1} \nabla (\partial)^{ p-1-i} \cdot  \langle \partial, w \rangle \cdot  \nabla (\partial)^{i}  
 + \mr{pr}_U^{-1}(\nabla (\partial)^{[p]} - \nabla  (\partial^{[p]}))  - \epsilon \cdot \langle \partial, w \rangle. \notag \\
= & \ -  \epsilon \cdot  \langle \partial, w  \rangle + \epsilon \cdot   \sum_{i=0}^{p-1} \nabla (\partial)^{ p-1-i} \cdot  \langle \partial, w \rangle \cdot  \nabla (\partial)^{i}, \notag
\end{align} 
where  the last equality follows from   the assumption $\psi^{\nabla}(\partial) \left(= \nabla (\partial)^{[p]} -  \nabla (\partial^{[p]})\right) =0$.
Thus,    (\ref{xi1}) and  (\ref{xi2}) together show that
\begin{align}
\langle F_{U_\epsilon}^{-1}(\mr{pr}^{-1}_U(\partial)), \psi^{\nabla_\epsilon + \epsilon \cdot w} \rangle =  \langle F_{U_\epsilon}^{-1}(\mr{pr}^{-1}_U(\partial)), -\epsilon \cdot  C^{(\mfg_{\mcE}, \nabla^\mr{ad})}(w) \rangle.
\end{align}
By considering this equality for various open subschemes $U$, 
 we obtain  (\ref{910}), which   completes  the proof of the former assertion.

The latter assertion follows from
 the former assertion  and the commutativity of (\ref{Cartier3F}).
\end{proof}

\bco \label{y0158}
 Let $\mcE$ be a $\mbG$-bundle on $X$ and   $\nabla_{\epsilon'}$ an $S^\mr{log}_\epsilon$-connection on $\mcE_\epsilon$  with ${^p \psi}^{\nabla_{\epsilon'}} = 0$.
 Also, write $\nabla := \mr{in}^*_S (\nabla_{\epsilon'})$ and $\nabla_{\epsilon} := \mr{pr}_S^*(\nabla)$.
  Then, 
  the flat $\mbG$-bundle $(\mcE_\epsilon, \nabla_{\epsilon'})$ on $X_\epsilon^\mr{log}/S^\mr{log}_\epsilon$
   is, Zariski locally on $X$,
   isomorphic to 
   $(\mcE_\epsilon, \nabla_{\epsilon})$.
  \eco
\begin{proof}
Let 
$w$
 be the unique  element of  $\Gamma (X, \Omega\otimes \mfg_{\mcE})$
 satisfying   $\nabla_{\epsilon'} = \nabla_\epsilon + \epsilon \cdot w$.
It follows from  the latter assertion of Proposition \ref{y0157} that, for each point of $X$,  there exist
an open neighborhood  $U$ of this point and an element $v$ of $\Gamma (U, \mfg_\mcE)$ with $\nabla^{\mr{ad}}(v) =w |_U$.
Let $\eta_\mcE^v$ be the automorphism
 of $\mcE_\epsilon |_U$ obtained from $v$ as in Remark \ref{QW9004}.
 Then, we have $\eta_\mcE^{v*}(\nabla_{\epsilon'}|_U) = \nabla_{\epsilon'} |_U - \epsilon \cdot \nabla^\mr{ad} |_U (v) = \nabla_\epsilon |_U$ (cf. (\ref{QW8601})).
Hence,  $\eta_\mcE^v$ determines an isomorphism $(\mcE_\epsilon |_U, \nabla_{\epsilon'}|_U) \isom (\mcE_\epsilon |_U, \nabla_{\epsilon}|_U)$.
This completes the proof of the assertion.
\end{proof}

\bpr \label{y0159} 
Let $(\mcV, \nabla)$ be a flat  vector bundle on $X^\mr{log}/S^\mr{log}$ with ${^p}\psi^\nabla =0$.
Then,
both $\mr{Ker}(\nabla)$ and $\mr{Coker}(\nabla)$ are relatively torsion-free sheaves (cf. Definition \ref{QS1406}) on $X^{(1)}$ of rank  $\mr{rk}(\mcV)$.
 \epr
\begin{proof}
To complete the proof, it suffices to consider  the assertion for 
 $\mr{Ker}(\nabla)$ because $\mr{Coker}(\nabla)$ is isomorphic to $\Omega^{(1)} \otimes \mr{Ker}(\nabla)$ by   
 $\overline{C}^{(\mcV, \nabla)}$.
Let us take an $S$-scheme $s: S' \migi S$ and write $X' := S' \times_S X$.
Then, the following  square diagram is commutative:
 \begin{align} 
\vcenter{\xymatrix@C=46pt@R=36pt{
s^*(F_{X/S*}(\mcV))
\ar[r]^-{\sim}\ar[d]_-{s^*(F_{X/S*}(\nabla))}& F_{X'/S'*}(s^*(\mcV))
  \ar[d]^-{F_{X'/S'*}(s^*(\nabla))}\\
s^*(F_{X/S*}(\Omega \otimes \mcV))
 \ar[r]_-{\sim}
 &F_{X'/S'*}(s^*(\Omega \otimes \mcV)),
}}
\end{align}
where the upper and lower horizontal arrows are base-change maps and we use the notation  $s^*(-)$ to    denote  base-change by $s$.
This square diagram 
 induces an isomorphism of $\mcO_{X'^{(1)}}$-modules
$s^*(\mr{Coker}(\nabla)) \isom \mr{Coker}(s^*(\nabla))$.
Hence, since  the formations of $C^{(\mcV, \nabla)}$ and     (\ref{Cartier3F}) commute with
 the base-change by $s$, 
  the natural morphism 
$s^*(\mr{Ker}(\nabla)) \migi  \mr{Ker}(s^*(\nabla))$
 turns out be  an isomorphism.
 In particular,   we obtain the injectivity of  the pull-back
$s^*(\mr{Ker}(\nabla)) \migi  s^*( F_{X/S*}(\mcV))$
 of the inclusion $\mr{Ker}(\nabla) \migiincl F_{X/S*}(\mcV)$.
 By applying this argument to various  $S$-schemes  $s : S' \migi S$,
we conclude that $\mr{Ker}(\nabla) \migi F_{X/S*}(\mcV)$ is universally injective with respect to base-change to $S$-schemes. 
This implies the flatness of $ F_{X/S*}(\mcV)/\mr{Ker}(\nabla)$  over $S$
(cf. ~\cite{MAT}, p.\,17, Theorem 1).
By ~\cite{Ja}, Proposition 1.1.1, (1), $\mr{Ker}(\nabla)$ is verified  to be  relatively torsion-free, as desired.

Moreover, the computation of its rank $\mr{rk}(\mr{Ker}(\nabla))$ follows from  the equivalence of categories (\ref{equicat}).
In fact,  the rank of $\mr{Ker}(\nabla)$ coincides with that of $\mcV$ (as a vector bundle on $X$) because the natural morphism $F_{X/S}^*(\mr{Ker}(\nabla)) \migi \mcV$ is an isomorphism over the scheme-theoretically dense open subscheme  $X^\mr{sm}$ of $X$.
This  completes the proof of Proposition \ref{y0159}.
\end{proof}

The following corollary follows immediately from  Corollary \ref{y0158} and Proposition \ref{y0159}.

\bco \label{y0160} 
Let $n$ be a positive integer and $(\mcV_{\epsilon'}, \nabla_{\epsilon'})$  a  rank $n$ flat bundle  on $X_\epsilon^\mr{log}/S_\epsilon^\mr{log}$ with ${^p}\psi^{\nabla_{\epsilon'}} =0$.
Denote by $(\mcV, \nabla)$ the base-change of $(\mcV_{\epsilon'}, \nabla_{\epsilon'})$ via  $\mr{in}_S : S \migi S_\epsilon$.
Then,  the $\mcO_{{X_{\epsilon}^{(1)}}}$-module 
 $\mr{Ker}(\nabla_{\epsilon'})$ is 
 relatively torsion-free sheaf of rank $n$
 on ${X_{\epsilon}^{(1)}}$,  and the canonical $\mcO_{X^{(1)}}$-linear morphism 
 $\mr{in}_{X^{(1)}}^*(\mr{Ker}(\nabla_{\epsilon'})) \isom \mr{Ker}(\nabla)$ 
is an isomorphism.
 \eco

\subsection{The dual connection and the Frobenius morphism} \label{y0161}

We study the  dual of a logarithmic connection  viewed as $\mcO_{X^{(1)}}$-linear morphism.

Denote by $\omega_{X/S}$ and $\omega_{X^{(1)}/S}$
the dualizing sheaves
    of $X/S$ and $X^{(1)}/S$ respectively (cf. \S\,\ref{QQ03}).
 Then,  there exist isomorphisms   
$\omega_{X/S} \isom \Omega (-D_\msX)$,
  $F_{X/S}^!(\omega_{X^{(1)}/S}) \isom   \omega_{X/S}$.
For each $\mcO_{X^{(1)}}$-module $\mcG$, we write
\begin{equation}
\label{1003}
 \mcG^\VV \hspace{-1mm}  := \mcH om_{\mcO_{X^{(1)}}} (\mcG, \omega_{X^{(1)}/S}).
 \index{$\mcG^\VV \hspace{-1mm}$}
\end{equation}
It is clear that $(\mcG^\VV \hspace{-1mm})^\VV \cong \mcG$ if $\mcG$ is a vector bundle.
Given a vector bundle $\mcV$ on $X$, we obtain a composite isomorphism  
\begin{align} \label{1002}
\Psi^\mcV : F_{X/S*}(\mcV)^\VV  \hspace{-1mm} \isom F_{X/S*}(\mcH om_{\mcO_{X}}(\mcV, \omega_{X/S})) \isom F_{X/S*}( \Omega \otimes \mcV^\vee (-D_\msX)),
 \index{$\Psi^\mcV$}
 \end{align}
where the
first isomorphism follows from  Grothendieck-Serre duality (cf.   ~\cite{Har}, Chap III, Exercise 6.10 (b)).

Now, let $(\mcV, \nabla)$ be a rank $n (>0)$ flat  bundle on $X^\mr{log}/S^\mr{log}$   with ${^p}\psi^\nabla =0$.
 By applying the contravariant functor $(-)^{\ooalign{$\vee$ \cr $\wedge$}}$ 
  to the  morphism $F_{X/S*}(\nabla) : F_{X/S*}(\mcF) \migi F_{X/S*}(\Omega \otimes \mcF)$ induced by $\nabla$,
  we obtain a morphism of $\mcO_{X^{(1)}}$-modules
\begin{align}
F_{X/S*}(\nabla)^\VV \hspace{-1mm}  : F_{X/S*}(\Omega  \otimes \mcV)^\VV \hspace{-1mm} 
\migi F_{X/S*}(\mcV)^\VV.
\end{align}
On the other hand,  
observe that  the dual $\nabla^\vee$ of $\nabla$ restricts to a morphism
\begin{equation}
{^c}\nabla^\vee
: \mcV^\vee(-D_\msX) \migi \Omega\otimes \mcV^\vee(-D_\msX),
\label{rest}
\end{equation}
 which forms an $S^\mr{log}$-connection on $\mcV^\vee (-D_\msX)$ with vanishing $p$-curvature.
 In particular, it induces 
 an $\mcO_{X^{(1)}}$-linear morphism 
 \begin{align}
 F_{X/S*}({^c}\nabla^\vee): F_{X/S*}(\mcV^\vee (-D_\msX)) \migi F_{X/S*}(\Omega \otimes \mcV^\vee (-D_\msX)).
 \end{align}

\bpr \label{y0162}
\begin{itemize}
\item[(i)]
 The square diagram
\begin{align} \label{Sq3}
\vcenter{\xymatrix@C=46pt@R=36pt{
F_{X/S*}(\Omega\otimes \mcV)^\VV \hspace{-1mm}
\ar[r]_-{\sim}^-{\Psi^{\Omega\otimes \mcV}} \ar[d]_-{F_{X/S*}(\nabla)^{\scalebox{0.7}{$\VV$}} \hspace{-1mm}}& F_{X/S*}(\mcV^\vee(-D_\msX))  \ar[d]^-{F_{X/S*}({^c}\nabla^\vee)}\\
F_{X/S*}(\mcV)^\VV \hspace{-1mm} \ar[r]^-{\sim}_-{ \Psi^{\mcV}}& F_{X/S*}(\Omega\otimes \mcV^\vee(-D_\msX) )
}}
\end{align}
 is anti-commutative, i.e., the following equality holds:
 \begin{equation} \label{1266789f}
 F_{X/S*}({^c}\nabla^\vee) \circ \Psi^{\Omega\otimes \mcV} = -   \Psi^{\mcV} \circ F_{X/S*}(\nabla)^\VV \hspace{-1mm}. 
 \end{equation}

 \item[(ii)]
 There exist   canonical isomorphisms
\begin{equation}
\mr{Ker}({^c}\nabla^\vee)^\VV  \hspace{-1mm} \isom \mr{Coker}(\nabla),
\hspace{5mm}
\mr{Coker}({^c}\nabla^\vee)^\VV \hspace{-1mm} \isom  \mr{Ker}(\nabla).
 \label{firstisom}
\end{equation}
 \end{itemize}
 \epr
\begin{proof}
We first consider assertion (i).
By  Proposition \ref{torsionE} and  ~\cite{Ja}, Proposition 1.1.3,  the $\mcO_{X^{(1)}}$-module on the left-hand upper corner of (\ref{Sq3}) is relatively torsion-free (of rank $n \cdot p$).
Hence, it suffices to verify the anti-commutativity of   (\ref{Sq3}) over a  scheme-theoretically dense open subscheme of $X$.
Moreover, since  $\nabla$ has vanishing $p$-curvature, $(\mcV, \nabla)$ is, Zariski locally on $X^\mr{sm} \left( \subseteq X\right)$, isomorphic to the trivial flat  bundle $(\mcO_X^{\oplus n}, d^{\oplus n})$ (cf. (\ref{equicat})).
It follows that the problem is reduced to
  the case where $(\mcV, \nabla) = (\mcO_X, d)$.

We shall prove the assertion of this case.
Let $U$ be an open subscheme of $X^\mr{sm}$ (hence $D_\msX \cap U = \emptyset$ and $\omega_{X/S} |_U = \Omega |_U= \Omega_{U/S}$).
Recall   that the trace morphism 
\begin{equation}
\left(F_{X/S*}(F^!_{X/S}(\Omega^{(1)}))|_{U^{(1)}}=\right)   F_{X/S*}(\Omega)|_{U^{(1)}} \migi \Omega^{(1)}|_{U^{(1)}}
\end{equation}
 coincides (after composing with the inclusion $\Omega^{(1)} \migiincl \Omega^{(1)} \otimes F_{X/S*}(\mcO_X)$) with  the restriction of the Cartier operator $C^{(\mcO_X, d)}|_{U^{(1)}}$ associated to $(\mcO_X, d)$.
Hence, for any local sections $a$, $b \in F_{X/S*}(\mcO_X)|_{U^{(1)}}$ (regarded  as local sections of $\mcO_U$), we have a sequence of equalities
\begin{align}
(F_{X/S*}(d)^{\ooalign{$\vee$ \cr $\wedge$}} \circ (\Psi^{\Omega})^{-1} (a))(b)  & = (\Psi^{\Omega})^{-1} (a)(db)  \\
&  = C^{(\mcO_X, d)} (a \cdot d(b))  \notag \\
& = C^{(\mcO_X, d)} (d(a \cdot b)) - C^{(\mcO_X, d)} (da \cdot b) \notag \\
& = - C^{(\mcO_X, d)} (da \cdot b) \notag \\
&=  - ( \Psi^{\mcO_X})^{-1} (d(a))(b) \notag \\
& =   (-( \Psi^{\mcO_X})^{-1} \circ F_{X/S*}(d) (a))(b), \notag
\end{align}
where the fourth  equality follows from the commutativity of   (\ref{Cartier3F}).
This implies the required anti-commutativity.

Next, let us  consider assertion (ii).
By Propositions \ref{torsionE} and  \ref{y0159},   both $\mr{Coker}(\nabla)$ and $F_{X/S*}(\Omega\otimes \mcV)$ are    relatively torsion-free.
It follows that $F_{X/S*}(\mr{Im}(\nabla))$ is relatively torsion-free (cf. ~\cite{Ja}, Proposition 1.1.1,  (i)).
Hence, by Proposition 1.1.2 in {\it loc.\,cit.,}  the sequence
\begin{align} 
\vcenter{\xymatrix@C=6pt@R=36pt{
 0 
 \ar[rd]^-{}_-{}
 &&F_{X/S*}(\Omega\otimes \mcV)^\VV 
  \ar[rd]^{F_{X/S*}(\nabla)^{\scalebox{0.7}{$\VV$}}}_-{}& &  \mr{Ker}(\nabla)^\VV  \ar[dr]&
\\
&  \mr{Coker}(\nabla)^\VV  \ar[ur]_-{}^-{}& &
 F_{X/S*}(\mcV)^\VV  \ar[ur] &&0
}}
\end{align}
arising  naturally  from  $F_{X/S*}(\nabla)$
 is  exact.
 This sequence induces  isomorphisms
 \begin{align} \label{QH36}
 \mr{Coker}(\nabla)^\VV \hspace{-1mm}  \isom  \mr{Ker}(F_{X/S*}(\nabla)^\VV \hspace{-1mm}),\hspace{5mm} 
   \mr{Ker}(\nabla)^\VV \hspace{-1mm} \isom \mr{Coker}(F_{X/S*}(\nabla)^\VV \hspace{-1mm}).
 \end{align}
On the other hand,  taking the kernels and  cokernels of 
the vertical arrows in (\ref{Sq3}) yields 
  isomorphisms
\begin{align} \label{QH37}
\mr{Ker}(F_{X/S*}(\nabla)^\VV \hspace{-1mm}) \isom \mr{Ker}({^c}\nabla^\vee)),
\hspace{5mm}
 \mr{Coker}(F_{X/S*}(\nabla)^\VV \hspace{-1mm}) \isom \mr{Coker}({^c}\nabla^\vee). 
\end{align} 
 By composing  (\ref{QH36}) with  (\ref{QH37})  and successively taking the duals (cf. ~\cite{Ja}, Proposition 1.1.4), we obtain the desired isomorphisms.
This completes the  proof of the assertion.
\end{proof}


\subsection{The nondegenerate pairing on  de Rham cohomology} \label{QG1054}
Next, consider the tensor product
\begin{equation}  
 (\mcV^\vee(-D_\msX)\otimes\mcV,  {^c}\nabla^\vee \otimes \nabla)  \end{equation}
(cf. (\ref{QQ295})) 
 of  $(\mcV^\vee(-D_\msX),  {^c}\nabla^\vee )$ and $(\mcV, \nabla)$.
Denote by ${^c d}$  the restriction $\mcO_X (-D_\msX) \migi \mcO_X (-D_\msX) \otimes \Omega \left(= \omega_{X/S} \right)$ of the universal derivation $d$.
Then, the natural pairing 
$\mcV^\vee(-D_\msX)\otimes\mcV \migi \mcO_X(-D_\msX)$
 is compatible with the  $S^\mr{log}$-connections, i.e.,  ${^c\nabla}^\vee \otimes \nabla$ and ${^c d}$.
By composing this pairing and the cup product in the de Rham cohomology, we obtain 
a  composite $\mcO_S$-bilinear pairing 
\begin{align}
 \oint : \mbR^1f_*(\mcK^\bullet[{^c}\nabla^\vee]) \times \mbR^1f_*(\mcK^\bullet [\nabla]) 
 \migi \mbR^2f_*(\mcK^\bullet[{^c} d])  \isom \mbR^1f_*(\omega_{X/S})  \isom \mcO_S. \hspace{-3mm}
 \index{$\oint$, $\oint'$}
\end{align}
In particular, this morphism yields a morphism of $\mcO_S$-modules
\begin{equation} \label{QW9700}
 \oint' : \mbR^1f_*(\mcK^\bullet[{^c}\nabla^\vee]) \migi  \mbR^1f_*(\mcK^\bullet [\nabla])^\vee.
\end{equation}
 By the definition of $\oint$ (and by, e.g.,  the explicit description of the cup product in the de Rham cohomology in terms of the \v{C}ech double complex), 
 the restriction of $\oint$ to the image of ${'' e}_\sharp[{^c}\nabla^\vee] \times {'' e}_\sharp [\nabla]$ (cf.  (\ref{Ceq})) vanishes  identically.
 Thus, the morphism $\oint'$ fits into the following diagram:
\begin{align} \label{oint}
\vcenter{\xymatrix@C=30pt@R=36pt{
0 \ar[r]
&\mbR^1f_*(\mr{Ker}({^c}\nabla^\vee )) \ar[r]^-{ {'' e}_\sharp[{^c}\nabla^\vee]} \ar[d]^-{\oint''}
&\mbR^1f_*(\mcK^\bullet[{^c}\nabla^\vee]) \ar[r]^-{{'' e}_\flat [{^c}\nabla^\vee ] } \ar[d]^-{\oint'}
&f_* (\mr{Coker}({^c}\nabla^\vee )) \ar[r] \ar[d]^-{\oint'''} & 0\\
0 \ar[r]&f_*(\mr{Coker}(\nabla))^\vee \ar[r]_-{{'' e}_\flat[\nabla]^\vee}
&\mbR^1f_*(\mcK^\bullet [\nabla])^\vee \ar[r]_-{{'' e}_\sharp[\nabla]^\vee}
&\mbR^1f_*(\mr{Ker}(\nabla))^\vee \ar[r]& 0,  \hspace{-3mm}
}}
\end{align}
where we recall that the upper horizontal  sequence is exact (cf. (\ref{Ceq})) and the lower sequence is left-exact.
\bco \label{y0163}
The morphisms $\oint'$, $\oint''$, and $\oint'''$ are all isomorphisms.
In particular, the lower horizontal  sequence of (\ref{oint}) is exact, and the $\mcO_S$-bilinear pairing $\oint$ is nondegenerate.
 \eco
\begin{proof}
Since $X/S$ is of relative dimension $1$, the equality 
$\mbR^2f_*(\mr{Ker}({^c}\nabla^\vee)) =0$ holds.
Moreover, according to  Proposition \ref{y0159},  
$\mr{Ker}({^c}\nabla^\vee)$ is  flat over $S$.
Hence,  $\mbR^1f_*(\mr{Ker}({^c}\nabla^\vee))$ is a vector bundle 
(cf. ~\cite{Har}, Chap.\,III, Theorem 12.11 (b)), and  we have 
$\mbR^1f_*(\mr{Ker}({^c}\nabla^\vee))^{\vee \vee} = \mbR^1f_*(\mr{Ker}({^c}\nabla^\vee))$.
 One verifies from the definition of  $\oint$ that the dual of the composite isomorphism
\begin{align}
 f_*(\mr{Coker}(\nabla)) \isom  f_*^{(1)}(\mr{Ker}({^c}\nabla^\vee)^\VV  \hspace{-1mm}) \isom \mbR^1f_*(\mr{Ker}({^c}\nabla^\vee))^\vee
\end{align}
coincides with
$\oint''$, where the first  arrow  arises  from  the first isomorphism in   (\ref{firstisom}) and  the second  arrow follows from    Grothendieck-Serre duality.
In particular,  $\oint''$ turns out to be  an isomorphism.
A similar argument (using  the dual of the second  isomorphism in  (\ref{firstisom})) shows  that $\oint'''$ is an isomorphism, which implies the exactness of 
 the lower horizontal sequence  of (\ref{oint}).
By   the five lemma  
applied to   (\ref{oint}),  
we see that $\oint'$ is  an isomorphism, i.e., $\oint$ is nondegenerate. 
\end{proof}

\section{Infinitesimal deformations of a dormant $(\mfg, \hslash)$-oper} \label{y0164}

We shall  go back to the discussion on  opers (restricted to the case of $\hslash =1$).
Let us keep the assumption in the previous section, i.e., that $\mr{char}(k) =:p > 2h_\mbG$.
Also, let  
 $\msE^\spadesuit = (\DE_{\mbB}, \nabla)$ be  a $\ \qq$-normal  dormant  $\mfg$-oper on $\msX$.

\subsection{The adjoint operator determined by  monodromy} \label{y0152}

Consider the adjoint operator
\begin{equation}
\label{908}   \mr{ad}(\mu_i^{\nabla}) : \sigma^*_i (\mfg_{\DE_\mbG}) \migi \sigma^*_i (\mfg_{\DE_\mbG})  \end{equation}
($i \in \{ 1, \cdots, r \}$) on $\sigma^*_i (\mfg_{\DE_\mbG}) $  determined by the monodromy $\mu_i^{\nabla}$ of $\nabla$;
this element of  $\mr{End}_{\mcO_S}(\sigma_i^*(\mfg_{\DE_\mbG}))$ coincides, by definition, with
the monodromy
$\mu_{i}^{\nabla^\mr{ad}}$ of $\nabla^\mr{ad}$ in  the sense of Definition  \ref{QQ1201}.
  The following proposition concerning $\mr{ad} (\mu_i^\nabla)$ will be used in the proof of Proposition \ref{9011} below.

\bpr \label{y0153} 
Both 
$\mr{Ker}( \mr{ad}(\mu_i^\nabla))$ and $\mr{Coker}( \mr{ad}(\mu_i^\nabla)$ are vector bundles on $S$ of rank $\mr{rk}(\mfg)$.
 \epr
\begin{proof}
The $\mcO_S$-modules  $\mr{Ker}( \mr{ad}(\mu_i^\nabla))$,  $\mr{Coker}( \mr{ad}(\mu_i^{\nabla})$ may be obtained by pulling-back the corresponding sheaves for  
the universal family of  $\mfg$-opers on $\mfC_{g,r} \times_{\overline{\mfM}_{g,r}} \mfO \mfp_{\mfg,  g,r}$ via the  classifying morphism $s_{\mfo \mfp} : S \migi  \mfO \mfp_{\mfg, g,r}$ of $\msE^\spadesuit$.
Since $\mfO \mfp_{\mfg, g,r}$ is reduced,  we may assume, without loss of generality, that  $S$ is reduced.
Under this assumption,  a given  coherent $\mcO_S$-module $\mcH$ is locally free of rank $l \geq 0$ if and only if for  any closed point $s \in S$, the fiber $\mcH \otimes_{\mcO_S} k(s)$ of $\mcH$, where $k(s)$ denotes the residue field of $s$,  is a $k(s)$-vector space of rank $l$.
Hence, it suffices to prove the assertion in the case where   $S = \mr{Spec}(k)$ and $k$ is algebraically closed.
Then, 
since $\msE^\spadesuit$ is $\, \qq$-normal,
 the monodromy $\mu_i^{\nabla} \in  \mfg$ may be expressed as   $q_{-1} + R$ for some $R \in \mfg^{\mr{ad}(q_1)}$.
The element  $q_{-1} + R$   is known to be regular (cf. ~\cite{Ngo}, the comment  preceding Lemma 1.2.3), so both the kernel and  cokernel of the adjoint operator $\mr{ad}(q_{-1} + R)$ are of rank $\mr{rk}(\mfg)$.
This completes the proof of  Proposition \ref{y0153}.
\end{proof}

\bpr \label{9011} 
 The $\mcO_S$-module  $f_*(\mr{Coker}(\nabla^\mr{ad}))$ (resp., $\mbR^1f_*(\mr{Ker}(\nabla^\mr{ad}))$) is a vector bundle 
 of rank $\aleph (\mfg)$ (resp., ${^c \aleph} (\mfg)$) (cf. (\ref{aleph})).
 \epr
\begin{proof}
According to  Proposition \ref{y0159},  both $\mr{Coker}(\nabla^\mr{ad})$ and 
$\mr{Ker}(\nabla^\mr{ad})$ are flat over $S$.
Hence, since  $\mbR^1f_*(\mr{Coker}(\nabla^\mr{ad})) =0$  (by Proposition \ref{y0143}, (iii)) and  $\mbR^2f_*(\mr{Ker}(\nabla^\mr{ad})) =0$ (by the fact that  $X/S$ is of relative dimension $1$), 
 both $f_*(\mr{Coker}(\nabla^\mr{ad}))$  and $\mbR^1f_*(\mr{Ker}(\nabla^\mr{ad}))$ are vector bundles on $S$ (cf. ~\cite{Har}, Chap.\,III, Theorem 12.11 (b)).

In the following, we shall compute the ranks of these vector bundles.
First,  by  Corollary \ref{y0163} in  the case where  $(\mcF, \nabla)$ is taken to be $(\mfg_{\DE_{\mbG}}, \nabla^\mr{ad})$,
there exists an isomorphism  of $\mcO_S$-modules
\begin{equation} \label{1010}
\oint'' : \mbR^1f_*(\mr{Ker}({^c}\nabla^{\mr{ad}\vee})) \isom f_*(\mr{Coker}(\nabla^\mr{ad}))^\vee.
\end{equation}
In particular, $\mbR^1f_*(\mr{Ker}({^c}\nabla^{\mr{ad}\vee}))$ is a vector bundle on $S$ of rank equal to the rank of $f_*(\mr{Coker}(\nabla^\mr{ad}))$.
On the other hand, the Killing form on $\mfg$
induces an isomorphism of flat bundles
$(\mfg_{\DE_{\mbG}}, \nabla^\mr{ad}) \isom (\mfg_{\DE_{\mbG}}^\vee,  \nabla^{\mr{ad}\vee})$.
If ${^c}\nabla^\mr{ad}$ (resp., ${^c}\nabla^{\mr{ad}\vee}$) denotes the $S^\mr{log}$-connection on 
$\mfg_{\DE_{\mbG}}(-D_\msX)$ (resp., $\mfg^\vee_{\DE_{\mbG}}(-D_\msX)$) obtained by restricting $\nabla^\mr{ad}$ (resp., $\nabla^{\mr{ad}\vee}$), then
this isomorphism  restricts to  an isomorphism 
\begin{equation}
(\mfg_{\DE_{\mbG}}(-D_\msX), {^c}\nabla^\mr{ad}) \isom (\mfg^\vee_{\DE_{\mbG}}(-D_\msX), {^c}\nabla^{\mr{ad}\vee}).
\end{equation}
It  induces an isomorphism of $\mcO_S$-modules
\begin{equation} \label{1011}
\mbR^1f_*(\mr{Ker}({^c}\nabla^{\mr{ad}})) \isom \mbR^1f_*(\mr{Ker}({^c}\nabla^{\mr{ad}\vee})). 
\end{equation}
By composing (\ref{1010}) with (\ref{1011}), we obtain an isomorphism of $\mcO_S$-modules
\begin{equation} \label{1012}
\mbR^1f_*(\mr{Ker}({^c}\nabla^{\mr{ad}})) \isom f_*(\mr{Coker}(\nabla^\mr{ad}))^\vee,
\end{equation}
which  implies  the following equality: 
\begin{equation}
\mr{rk}(\mbR^1f_*(\mr{Ker}({^c}\nabla^{\mr{ad}}))) = \mr{rk}(f_*(\mr{Coker}(\nabla^\mr{ad}))).
\label{KC9}
\end{equation}

Next, consider the natural  injection of short exact sequences
\begin{align} \label{Fex}
\vcenter{\xymatrix@C=16pt@R=36pt{
0 \ar[r]
&\mfg_{\DE_{\mbG}}(-D_\msX)  \ar[r] \ar[d]^-{{^c}\nabla^\mr{ad}}
&\mfg_{\DE_{\mbG}} \ar[r] \ar[d]^-{\nabla^\mr{ad}}
&\bigoplus_{i=1}^r \sigma_{i*}(\sigma_i^*(\mfg_{\DE_{\mbG}})) \ar[r] \ar[d]^-{\bigoplus_{i=1}^r \sigma_{i*}(\mr{ad}(\mu_i^{\nabla}))}
& 0 \\
0 \ar[r]&\Omega\otimes \mfg_{\DE_{\mbG}}(-D_\msX) \ar[r]
&\Omega \otimes \mfg_{\DE_{\mbG}} \ar[r]
&\bigoplus_{i=1}^r \sigma_{i*}(\sigma_i^*(\Omega \otimes \mfg_{\DE_{\mbG}})) \ar[r] & 0.
}}
\end{align}
By the snake lemma,  this diagram gives rise to  an exact sequence
\begin{align}  \label{Lex}
\vcenter{\xymatrix@C=-3pt@R=36pt{
 0 
 \ar[rd]^-{}_-{}
 && \mr{Ker}(\nabla^\mr{ad}) 
  \ar[rd]^{}_-{}& & \mr{Coker}({^c}\nabla^\mr{ad})  \ar[dr]^-{\delta}&
\\
&  \mr{Ker}({^c}\nabla^\mr{ad})   \ar[ur]_-{}^-{}& &
\bigoplus_{i=1}^r \sigma_{i*}(\mr{Ker}(\mr{ad}(\mu_i^\nabla)))  \ar[ur] && \mr{Coker}(\nabla^\mr{ad}).
}}
\end{align}
Here,  let us consider the square diagram
\begin{align}  \label{Comm3}
\vcenter{\xymatrix@C=46pt@R=36pt{
\mr{Coker}({^c}\nabla^\mr{ad}) \ar[r]^{\delta} \ar[d]_-{\wr}& \mr{Coker}(\nabla^\mr{ad}) \ar[d]^-{\wr}\\
\Omega^{(1)}\otimes \mr{Ker}({^c}\nabla^\mr{ad}) \ar[r]_-{\mr{inclusion}}& \Omega^{(1)}\otimes \mr{Ker}(\nabla^\mr{ad}),
}}
\end{align}
where 
  the left-hand and right-hand vertical arrows denote the isomorphisms $\overline{C}^{(\mfg_{\DE_\mbG}(-D_\msX), {^c}\nabla^\mr{ad})}$ and 
  $\overline{C}^{(\mfg_{\DE_\mbG}, \nabla^\mr{ad})}$
  respectively.
By  the construction of  Cartier operators,   this diagram is commutative.
This implies the injectivity of $\delta$, and 
(\ref{Lex}) yields a short  exact sequence
\begin{equation} \label{ex15}
0  \longmigi  \mr{Ker}({^c}\nabla^\mr{ad}) \longmigi \mr{Ker}(\nabla^\mr{ad}) \longmigi \bigoplus_{i=1}^r \sigma_{i*}(\mr{Ker}(\mr{ad}(\mu_i^\nabla)))  \longmigi 0.
\end{equation}
Moreover,  observe that $\mbR^1f_*(\bigoplus_{i=1}^r \sigma_{i*}(\mr{Ker}(\mr{ad}(\mu_i^\nabla)))) = 0$ and $f_*( \mr{Ker}(\nabla^\mr{ad})) =0$
(cf. Proposition \ref{y0143}, (i)).
Hence, by applying the functor $\mbR^1f_*(-)$ to (\ref{ex15}), we obtain a short exact sequence
\begin{align} \label{1014}
0 \longmigi \bigoplus_{i=1}^r \mr{Ker}(\mr{ad}(\mu_i^\nabla))   \longmigi \mbR^1f_*(\mr{Ker}({^c}\nabla^\mr{ad}))  \longmigi \mbR^1f_*(\mr{Ker}(\nabla^\mr{ad})) \longmigi 0. 
\end{align}
Since $\mbR^1f_*(\mr{Ker}(\nabla^\mr{ad}))$ is a vector bundle as proved above,
 the exactness of   (\ref{1014}) shows that $\mbR^1f_*(\mr{Ker}({^c}\nabla^\mr{ad}))$ is a vector bundle on $S$.
The rank of this vector bundle satisfies 
\begin{align} \label{un}
\mr{rk}(\mbR^1f_*(\mr{Ker}({^c}\nabla^\mr{ad}))) 
& = \mr{rk}(\mbR^1f_*(\mr{Ker}(\nabla^\mr{ad})) ) + \mr{rk}(\bigoplus_{i=1}^r \mr{Ker}(\mr{ad}(\mu_i^\nabla))   )  \\
& = \mr{rk}(\mbR^1f_*(\mr{Ker}(\nabla^\mr{ad})) ) + r \cdot \mr{rk}(\mfg), \notag 
\end{align}
where the second equality follows from Proposition \ref{y0153}.

Finally, by Proposition \ref{y0143}, (ii),  and  the exactness of  (\ref{Ceq}) for  $\nabla = \nabla^\mr{ad}$, 
  the following equalities hold:
\begin{align} \label{eq+}
 \mr{rk}(f_*(\mr{Coker}(\nabla^\mr{ad}))) + \mr{rk}(\mbR^1f_*(\mr{Ker}(\nabla^\mr{ad})))  & = \mr{rk}(\mbR^1f_*(\mcK^\bullet [\nabla^\mr{ad}]))  \\
 & = (2g-2+r) \cdot \mr{dim}(\mfg).  \notag
\end{align}

Consequently,    (\ref{KC9}), (\ref{un}), and (\ref{eq+}) induce  the desired  equalities
\begin{equation}
\mr{rk}(f_*(\mr{Coker}(\nabla^\mr{ad}))) =\aleph (\mfg), \ \ \ \mr{rk}(\mbR^1f_*(\mr{Ker}(\nabla^\mr{ad}))) = {^c \aleph} (\mfg).
\end{equation}
This completes the proof of Proposition \ref{9011}.
\end{proof}

\subsection{Deformation of a dormant $\mfg$-oper} \label{QG1057}
The $\mbG$-bundle ``$\DE_{\mbG, 1, U^\mr{log}/S^\mr{log}}$''  introduced in (\ref{canonicalEg2}) can be defined in the case where the log curve  $U^\mr{log}/S^\mr{log}$ 
is taken to be $\mfC_{g,r}^\mr{log}/\overline{\mfM}_{g,r}^\mr{log}$.
If $\DE_{\mbG, \msC_{g,r}}$ denotes the resulting $\mbG$-bundle, then we obtain 
the relative scheme 
\begin{align}
\ooalign{$\bigoplus$ \cr $\bigotimes^{\mfg}_{\msC_{g,r}}$} := f_{\mr{univ}*} (\Omega_{\mfC_{g,r}^\mr{log}/\overline{\mfM}_{g,r}}^{\otimes p} \otimes \mfg_{\DE_{\mbG, \msC_{g,r}}})
\end{align}
over $\overline{\mfM}_{g,r}$.
We equip this stack with the log structure pulled-back from $\overline{\mfM}_{g,r}^\mr{log}$; the resulting log stack will be denoted by $\ooalign{$\bigoplus$ \cr $\bigotimes^{\mfg \mr{log}}_{\msC_{g,r}}$}$.
The assignment sending   each $\, \qq$-normal  $\mfg$-oper 
to its $p$-curvature
 determines a strict morphism
\begin{equation}
\rotatebox[origin=c]{58}{{\large $\kappa$}}^{\mfg \mr{log}}_{}  :  \mfO \mfp^\mr{log}_{\mfg, g,r} \migi \ooalign{$\bigoplus$ \cr $\bigotimes^{\mfg \mr{log}}_{\msC_{g,r}}$}.
\end{equation}
of log stacks over $\overline{\mfM}_{g,r}^\mr{log}$.

Denote by 
$\text{\Pisymbol{pzd}{42}}^\mfg_{\msX} :  S \migi \ooalign{$\bigoplus$ \cr $\bigotimes^{\mfg}_{\msC_{g,r}}$}$
the zero section of  $\ooalign{$\bigoplus$ \cr $\bigotimes^{\mfg}_{\msC_{g,r}}$}$ over $S$ and by $s_{\mfo \mfp}^\mr{log} : S^\mr{log} \migi \mfO \mfp_{\mfg, g,r}^\mr{log}$ the strict morphism whose underlying $S$-rational point of $\mfO \mfp_{\mfg, g,r}$ classifies $(\msX, \msE^\spadesuit)$.
 Then, we have a canonical isomorphism of $\mcO_S$-modules
\begin{equation} \label{isom60}
  \left(\text{\Pisymbol{pzd}{42}}^{\mfg*}_\msX(\mcT_{\tiny{\ooalign{$\bigoplus$ \cr $\bigotimes^{\mfg \mr{log}}_{\msC_{g,r}}$}}/\overline{\mfM}^\mr{log}_{g,r}}) = \right)  \text{\Pisymbol{pzd}{42}}^{\mfg*}_\msX(\mcT_{\tiny{\ooalign{$\bigoplus$ \cr $\bigotimes^{\mfg}_{\msC_{g,r}}$}}/\overline{\mfM}_{g,r}}) \isom f_*(\Omega^{\otimes p} \otimes \mfg_{\DE_{\mbG}}).
\end{equation}
Also, the differential of   $\rotatebox[origin=c]{58}{{\large $\kappa$}}^{\mfg \mr{log}}_{}$ gives a composite 
\begin{align} \label{QH66}
d \rotatebox[origin=c]{58}{{\large $\kappa$}}^{\mfg \mr{log}} : s^{*}_{\mfo \mfp}(\mcT_{\mfO \mfp_{\mfg, g,r}^\mr{log}/k}) \migi \text{\Pisymbol{pzd}{42}}^{\mfg*}_\msX (\mcT_{\tiny{\ooalign{$\bigoplus$ \cr $\bigotimes^{\mfg \mr{log}}_{\msC_{g,r}}$}}/k}) \migisurj  
\text{\Pisymbol{pzd}{42}}^{\mfg*}_\msX (\mcT_{\tiny{\ooalign{$\bigoplus$ \cr $\bigotimes^\mfg_{\msC_{g,r}}$}}/\overline{\mfM}_{g,r}}),
\end{align}
where the second arrow denotes the projection to the second factor with respect to the decomposition $\text{\Pisymbol{pzd}{42}}^{\mfg*}_\msX(\mcT_{\tiny{\ooalign{$\bigoplus$ \cr $\bigotimes^{\mfg \mr{log}}_{\msC_{g,r}}$}}/k})  =  s^{*}_\mfm(\mcT_{\overline{\mfM}^\mr{log}_{g,r}/k})  \oplus  \text{\Pisymbol{pzd}{42}}^{\mfg*}_\msX(\mcT_{\tiny{\ooalign{$\bigoplus$ \cr $\bigotimes^{\mfg}_{\msC_{g,r}}$}}/\overline{\mfM}_{g,r}})$ induced by  differentiating the zero section of 
$\ooalign{$\bigoplus$ \cr $\bigotimes^{\mfg}_{\msC_{g,r}}$}/\overline{\mfM}_{g,r}$.

\ble \label{y0166}
The subsheaf  $\mr{Im}(\widetilde{\nabla}^{\mr{ad}})$ of  $\Omega \otimes \mfg_{\DE_{\mbG}}$ coincides with $\mr{Im}(\nabla^\mr{ad})$.   
 \ele
\begin{proof}
The inclusion relation  $\mr{Im}(\nabla^\mr{ad}) \subseteq \mr{Im}(\widetilde{\nabla}^\mr{ad})$ is clear because of  the fact $\nabla^\mr{ad} = \widetilde{\nabla}^\mr{ad} |_{\mfg_{\DE_\mbG}}$.
 We shall prove the inverse of this inclusion relation.
Denote by $\gamma$
 the $\mcO_X$-linear endomorphism of $\widetilde{\mcT}_{\DE_{\mbG}^{\mr{log}}/S^\mr{log}}$ given by assigning $v \mapsto v - \nabla  \circ d^\mr{log}_{\DE_\mbG}(v)$ for any local section $v\in \widetilde{\mcT}_{\DE_{\mbG}^{\mr{log}}/S^\mr{log}}$. Then, given  each $\partial \in \mcT$,   we have
\begin{align}
&  \ \ \  \  \langle \partial,  \widetilde{\nabla}^{\mr{ad}}(\gamma (v)) \rangle \\
&= \langle  \partial, \widetilde{\nabla}^{\mr{ad}} (v - \nabla \circ d^\mr{log}_{\DE_\mbG}(v)) \rangle \notag  \\
& = [ \nabla (\partial),  v - \nabla \circ d^\mr{log}_{\DE_\mbG}(s)] - \nabla ([\partial, d^\mr{log}_{\DE_\mbG}(v- \nabla \circ d^\mr{log}_{\DE_\mbG}(v))]) \notag \\
& = [\nabla (\partial),v] - \nabla ([\partial, d^\mr{log}_{\DE_\mbG}(s)]) 
  - \nabla  (  [\partial, d^\mr{log}_{\DE_\mbG}(v)]- [\partial, (d^\mr{log}_{\DE_\mbG} \circ \nabla) \circ d^\mr{log}_{\DE_\mbG} (v)]   ) \notag\\
& = [\nabla (\partial), v] - \nabla ([\partial, d^\mr{log}_{\DE_\mbG}(v)]) \notag \\
& =  \langle \partial, \widetilde{\nabla}^{\mr{ad}} (v)\rangle,
\notag \end{align}
where $\langle -,-\rangle$ denotes the pairing $\mcT \times (\Omega  \otimes \widetilde{\mcT}_{\DE_\mbG^\mr{log}/S^\mr{log}}) \migi \widetilde{\mcT}_{\DE_\mbG^\mr{log}/S^\mr{log}}$ induced by the natural pairing $\mcT \times \Omega \migi \mcO_X$.
This implies  $\mr{Im}(\widetilde{\nabla}^{\mr{ad}}) = \mr{Im}(\widetilde{\nabla}^{\mr{ad}}\circ \gamma)$. 
But,  
\begin{align}
 d^\mr{log}_{\DE_\mbG} (\gamma (v))= d^\mr{log}_{\DE_\mbG}(v) - (d^\mr{log}_{\DE_\mbG} \circ \nabla) \circ d^\mr{log}_{\DE_\mbG} (v) = d^\mr{log}_{\DE_\mbG} (v) - \mr{id}_\mcT \circ d^\mr{log}_{\DE_\mbG} (v)  = 0, 
\end{align}
so the image of $\gamma$ is contained in  $\mfg_{\DE_\mbG} \left(= \mr{Ker}(d^\mr{log}_{\DE_\mbG} )\right)$.
It follows that
 \begin{align}
 \mr{Im}(\widetilde{\nabla}^{\mr{ad}}) = \widetilde{\nabla}^\mr{ad} (\mr{Im}(\gamma)) \subseteq \mr{Im} (\widetilde{\nabla}^\mr{ad} |_{\mfg_{\DE_\mbG}}) = \mr{Im}(\nabla^\mr{ad}),
 \end{align}
thus completing  the proof of Lemma \ref{y0166}.
\end{proof}

By Lemma \ref{y0166} above, the square diagram
\begin{align} \label{QW9806}
\vcenter{\xymatrix@C=46pt@R=36pt{
\widetilde{\mcT}_{\DE_{\mbG}^{\mr{log}}/S^\mr{log}} \ar[r]_-{} \ar[d]_-{\widetilde{\nabla}^{\mr{ad}}}& 0  \ar[d]\\
 \Omega\otimes \mfg_{\DE_{\mbG}} \ar[r]_-{\mr{quotient}}& \mr{Coker}(\nabla^\mr{ad})
}}
\end{align}
is commutative.
This diagram defines
 a  morphism  of complexes  $\mcK^\bullet[\widetilde{\nabla}^{\mr{ad}}] \migi \mr{Coker}(\nabla^\mr{ad})[-1]$.
By applying the functor $\mbR^1f_*(-)$ to the composite
\begin{equation}
\mcK^\bullet [\widetilde{\nabla}^{\mr{ad}(0)}] \migiincl \mcK^\bullet [\widetilde{\nabla}^{\mr{ad}}] \migi \mr{Coker}(\nabla^\mr{ad})[-1],\end{equation}
 we obtain a morphism of $\mcO_S$-modules
\begin{align} \label{1016}
\widetilde{e}^{(0)}_\flat :   \mbR^1f_*(\mcK^\bullet [\widetilde{\nabla}^{\mr{ad}(0)}])  \migi f_*( \mr{Coker}(\nabla^\mr{ad})).
\end{align}

On the other hand, let us consider the composite 
\begin{align} \label{comp998}
f_*(\mr{Coker}(\nabla^\mr{ad})) &\isom f_*(\Omega^{(1)}\otimes \mr{Ker}(\nabla^\mr{ad}))  \\
& \xrightarrow{\mr{inclusion}}
f_*^{(1)}(\Omega^{(1)}\otimes F_{X/S*}(\mfg_{\DE_\mbG})) \notag \\
& \isom
f_*^{(1)}(F_{X/S*}(F_{X/S}^*(\Omega^{(1)}) \otimes \mfg_{\DE_{\mbG}})) \notag \\
&  \isom  f_*(\Omega^{\otimes p} \otimes \mfg_{\DE_{\mbG}}),  \notag
\end{align}
where the first arrow is the direct image of  $\overline{C}^{(\mfg_{\DE_\mbG}, \nabla^\mr{ad})}$ by $f$ and the third  arrow  follows from the projection formula.
The composite of $\widetilde{e}^{(0)}_\flat$ and (\ref{comp998}) specifies a morphism
\begin{align} \label{inj44}
\mbR^1f_*(\mcK^\bullet [\widetilde{\nabla}^{\mr{ad}(0)}]) \migi f_*(\Omega^{\otimes p} \otimes \mfg_{\DE_{\mbG}}).
\end{align} 
Since (\ref{comp998}) is injective, the kernel of (\ref{inj44}) coincides with that of $\widetilde{e}^{(0)}_\flat$.

\bpr \label{y0167}
 Let us keep the above notation.
 Then, the following square diagram is anti-commutative:
\begin{align} \label{diag56}
\vcenter{\xymatrix@C=46pt@R=36pt{
s_{\mfo \mfp}^{*}(\mcT_{\mfO \mfp^\mr{log}_{\mfg, g,r}/k}) 
\ar[r]^-{d \rotatebox[origin=c]{58}{{\large $\kappa$}}^{\mfg \mr{log}}} \ar[d]^-{\wr}_-{\mr{KS}_{\msE^\spadesuit}} &  
 \text{\Pisymbol{pzd}{42}}^{\mfg*}_\msX(\mcT_{\tiny{\ooalign{$\bigoplus$ \cr $\bigotimes^\mfg_{\msC_{g,r}}$}}/\overline{\mfM}_{g,r}})
\ar[d]_-{\wr}^-{(\ref{isom60})} \\
\mbR^1f_*(\mcK^\bullet(\widetilde{\nabla}^{\mr{ad}(0)})) \ar[r]_-{(\ref{inj44})}   & f_*(\Omega^{\otimes p} \otimes \mfg_{\DE_{\mbG}})
}}
\end{align} 
 (cf. (\ref{isom3333}) for the definition of $\mr{KS}_{\msE^\spadesuit}$).
\epr
\begin{proof}
The assertion follows from the definitions of various morphisms involved  together with  Proposition \ref{y0157} and the commutativity of (\ref{Cartier3F}).
\end{proof}

\bco \label{y0169}
Let ${^p}s_{\mfo \mfp}^\mr{log} : S^\mr{log} \migi \mfO \mfp^{^\mr{Zzz...} \mr{log}}_{\mfg, g,r}$
\index{$s^\mr{log}_\mfm$, $s^{\mr{log}}_{\mfo \mfp}$, ${^c}s^{\mr{log}}_{\mfo \mfp}$, ${^p}s^{\mr{log}}_{\mfo \mfp}$}
 be  the 
strict morphism whose underlying $S$-rational point of 
$\mfO \mfp^{^\mr{Zzz...}}_{\mfg, g,r}$ classifies 
   the pair $(\msX, \msE^\spadesuit)$.
Then, the isomorphism $\mr{KS}_{\msE^\spadesuit}$  restricts to 
an  isomorphisms
\begin{align} \label{ggg03}
 {^p}\mr{KS}_{\msE^\spadesuit}: {^p}s_{\mfo \mfp}^{*}(\mcT_{\mfO \mfp^{^\mr{Zzz...} \mr{log}}_{\mfg, g,r}/k}) \isom  
 \mr{Ker} (\widetilde{e}^{(0)}_\flat).
 \index{${^{(c)}}\mr{KS}_{\msE^\spadesuit /\msX}$, ${^{(c)}}\mr{KS}_{\msE^\spadesuit}$, ${^p}\mr{KS}_{\msE^\spadesuit}$}
  \end{align}
 \eco
\begin{proof} 
The assertion follows from Proposition \ref{y0167} because $\mfO \mfp_{\mfg, g,r}^\ZZZ$ is defined as the scheme-theoretic inverse image of the zero section in $\ooalign{$\bigoplus$ \cr $\bigotimes^\mfg_{\msC_{g,r}}$}$ via $\rotatebox[origin=c]{58}{{\large $\kappa$}}^\mfg$.
\end{proof}

\subsection{Deformation of a dormant $\mfg$-oper of prescribed radii} \label{y0168}

In what follows, we  perform a computation of the deformation space of $\msE^\spadesuit$  and  give  a  criterion for  the \'{e}taleness of the universal  moduli stack classifying  dormant $\mfg$-opers.

Let ${^p}s_{\mfo \mfp}^\mr{log}$ be as in Corollary \ref{y0169}, and 
 suppose further  that $S = \mr{Spec}(k_v)$ for a field $k_v$ over $k$.
 We write
$v:= {^p}s_{\mfo \mfp}$.
  According to  Lemma \ref{y0166}, 
  the image of $\widetilde{\nabla}^{\mr{ad} (0)}$ is contained in the kernel of 
  the composite 
\begin{equation} 
\varPsi  : \Omega\otimes \mfg_{\DE_{\mbB}}^{-1} \xrightarrow{\mr{inclusion}}  \Omega \otimes \mfg_{\DE_{\mbG}}  \xrightarrow{\mr{quotient}} \mr{Coker}(\nabla^\mr{ad}).
\index{$\varPsi$}
\end{equation}
Hence,  by restricting    the codomain of  $\widetilde{\nabla}^{\mr{ad} (0)}$  to $\mr{Ker}(\varPsi)$, we obtain  an $f^{-1}(\mcO_S)$-linear  morphism
\index{$\nabla^\mr{ad}$@--- $\widetilde{\nabla}^{\mr{ad}(0)}_\flat$}
\begin{equation}
\widetilde{\nabla}^{\mr{ad}(0)}_\flat : \widetilde{\mcT}_{\DE_\mbB^\mr{log}/S^\mr{log}} \migi  \mr{Ker}(\varPsi).
\end{equation}
Since 
 the equality $\mbH^2(X, \mcK^\bullet [\widetilde{\nabla}^{\mr{ad}(0)}] ) =0$ holds (cf.  Lemma \ref{y0146}),
 the natural short exact sequence
\begin{equation}
0 \longmigi \mcK^\bullet [\widetilde{\nabla}^{\mr{ad}(0)}_\flat] \longmigi \mcK^\bullet [\widetilde{\nabla}^{\mr{ad}(0)}] \longmigi \mr{Im}(\varPsi)[-1] \longmigi 0
\end{equation}
 induces the exact sequence
\begin{align} \label{QW9999}
\vcenter{\xymatrix@C=-3pt@R=26pt{
0  \ar[rd]&& \mbH^1 (X, \mcK^\bullet [\widetilde{\nabla}^{\mr{ad}(0)}])  \ar[rd]_-{\widetilde{e}^{(0)}_\flat}&& \mbH^2 (X, \mcK^\bullet [\widetilde{\nabla}_\flat^{\mr{ad}(0)}]) \ar[rd] &  \\
& \mbH^1 (X, \mcK^\bullet [\widetilde{\nabla}_\flat^{\mr{ad}(0)}]) \ar[ru]&&\Gamma (X, \mr{Im}(\varPsi)) \ar[ru]  & &0.
}}
\end{align}
The second arrow in this sequence yields an isomorphism of $k_v$-vector spaces
\begin{align} \label{QW10000}
\mbH^1 (X, \mcK^\bullet [\widetilde{\nabla}_\flat^{\mr{ad}(0)}]) \isom \left(\mr{Ker}(\widetilde{e}^{(0)}_\flat) \cong \right) \mcT_v \mfO \mfp^{\ZZZ \mr{log}}_{\mfg, g,r}
\end{align}
(cf. Corollary \ref{y0169}).
 Here, let  $\mfA \mfr \mft_{/k_v}$ denote the category of local $k_v$-algebras with residue field $k_v$.
By means of the description of deformation  spaces in terms of 
 \v{C}ech cocycles, we see that 
 the deformation functor 
 \begin{align} \label{QW10002}
\mr{Def}_v : \mfA \mfr \mft_{/k_v} \migi \mfS \mfe \mft; R \mapsto  \begin{pmatrix}
\text{the set of isomorphism classes} \\
\text{of objects $(\msX_R, \msE_R^\spadesuit)$ in $\mfO \mfp^\ZZZ_{\mfg, g,r}(R)$} \\
 \text{with special fiber $\cong (\msX, \msE^\spadesuit)$}
\end{pmatrix}
 \end{align}
  has a tangent-obstruction theory by putting $\mbH^1 (X, \mcK^\bullet [\widetilde{\nabla}_\flat^{\mr{ad}(0)}])$ as the tangent space  and $\mbH^2 (X, \mcK^\bullet [\widetilde{\nabla}_\flat^{\mr{ad}(0)}])$ as the obstruction  space.
  (We refer the reader to ~\cite{FGA}, Chap.\,6, Definition 6.1.21, for the definition of a tangent-obstruction theory.)
  Moreover, the following assertion holds:

\bco \label{y0170}
Let  us keep the above notation.
Then, 
the logarithmic tangent space $\mcT_v \mfO \mfp^{\ZZZ \mr{log}}_{\mfg, g,r}$
of $\mfO \mfp^{^\text{Zzz..} \mr{log}}_{\mfg, g,r}$ at $v$ is a $k_v$-vector space 
  of dimension $\geq 3g-3 +r$.
If, moreover,  the natural projection  $\mfO \mfp^\ZZZ_{\mfg,  g,r} \migi \overline{\mfM}_{g,r}$ is unramified at $v$, then
 this projection is also flat (hence,  \'{e}tale)  at $v$.
 \eco
\begin{proof}
Recall  from   Corollary \ref{y0169} that
$\mcT_v \mfO \mfp_{\mfg, g,r}^{\ZZZ \mr{log}}$
 is isomorphic to  the kernel $\mr{Ker}(\widetilde{e}_\flat^{(0)})$ of the morphism
\begin{equation}   \label{1021}
\widetilde{e}_\flat^{(0)} : \mbH^1(X, \mcK^\bullet[\widetilde{\nabla}^{\mr{ad}(0)}])) \migi \Gamma (X, \mr{Coker}(\nabla^\mr{ad})). \end{equation}
Hence, by Propositions \ref{y0148} and  \ref{9011}, we have 
\begin{align} \label{1022}
\mr{dim}(\mr{Ker}(\widetilde{e}_\flat^{(0)})) & \geq \mr{dim}(\mbH^1(X, \mcK^\bullet[\widetilde{\nabla}^{\mr{ad}(0)}]))) - \mr{dim}(\Gamma (X, \mr{Coker}(\nabla^\mr{ad})))  \\
& = (\aleph (\mfg) + 3g-3+r) - \aleph (\mfg) \notag  \\
&  = 3g-3+r.  \notag
\end{align} 
This completes the proof of the former assertion.

Next, we shall consider the latter assertion.
Suppose that $\mfO \mfp^{^\mr{Zzz...}}_{\mfg, g,r}/\overline{\mfM}_{g,r}$ is unramified at $v$.
Then, by Corollary \ref{y0169} 
(and the fact that the log structure of $\mfO \mfp^{^\mr{Zzz...} \mr{log}}_{\mfg, g,r}$ is obtained by pull-back from  $\overline{\mfM}^\mr{log}_{g,r}$), 
the composite 
\begin{equation}
\label{isom77}
 \mr{Ker}(\widetilde{e}_\flat^{(0)}) \xrightarrow{\mr{inclusion}}\mbH^1(X, \mcK^\bullet[\widetilde{\nabla}^{\mr{ad}(0)}])) \migi H^1(X, \mcT)\end{equation}
is verified to be  injective, where the second arrow denotes the penultimate morphism in (\ref{tangent}).
This composite  is in fact  an isomorphism because of 
(\ref{1022}) together with the equality $\mr{dim}(H^1 (X, \mcT)) =3g-3+r$.
 In particular,  the equality   $\mr{dim}(\mr{Ker}(\widetilde{e}_\flat^{(0)})) = 3g-3+r$ holds.
We shall  observe the following sequence:
  \begin{align}
\mr{dim} (\Gamma (X, \mr{Coker}(\nabla^\mr{ad}))) \geq   &
  \ \mr{dim}(\Gamma (X, \mr{Im}(\varPsi))) \\
 \geq  & \  \mr{dim}(\mbH^1 (X, \mcK^\bullet [\widetilde{\nabla}^{\mr{ad}(0)}])) - \mr{dim}(\mr{Ker} (\widetilde{e}^{(0)}_\flat)) \notag \\
= &  \ (\aleph (\mfg) + 3g-3+r) - (3g-3+r) \notag \\
 = & \ \mr{dim} (\Gamma (X, \mr{Coker}(\nabla^\mr{ad})), \notag
 \end{align}
 where the last equality follows from  Proposition \ref{9011}.
This implies that the second  inequality in this  sequence must be an equality, so 
 the third arrow in (\ref{QW9999})  turns out to be  surjective.
Hence, by the exactness of (\ref{QW9999}), we have  $\mbH^2 (\mcK^\bullet [\widetilde{\nabla}_\flat^{\mr{ad}(0)}]) = 0$.
This  means that 
 any infinitesimal  deformation of $(\msX, \msE^\spadesuit)$ is unobstructed (cf. the discussion preceding this corollary).
That is to say, $\mfO \mfp^{\ZZZ \mr{log}}_{\mfg, g,r}$ is log  smooth at $v$.
Moreover, as is proved above, 
 the
 differential of the projection  $\mfO \mfp^{\ZZZ \mr{log}}_{\mfg, g,r} \migi \overline{\mfM}_{g,r}^\mr{log}$ is bijective at $v$.
 It follows that $\mfO \mfp^{^\mr{Zzz...} \mr{log}}_{\mfg,  g,r} \migi \overline{\mfM}^\mr{log}_{g,r}$ is log \'{e}tale at $v$ (cf. ~\cite{KaKa}, \S\,3, Proposition (3.12)), so the underlying morphism $\mfO \mfp^{^\mr{Zzz...}}_{\mfg,  g,r} \migi \overline{\mfM}_{g,r}$ is \'{e}tale  at the same  point  (cf. ~\cite{KaKa}, \S\,3, Proposition (3.8)).
This completes the proof of the assertion.
\end{proof}

\vspace{10mm}

\chapter{The pseudo-fusion ring for dormant opers} \label{y0171}\vspace{3mm}

In this chapter, we examine  the behavior of the generic degree
$\mr{deg}(\mfO \mfp^\ZZZ_{\mfg, \hslash, \hslash \star \rho, g,r}/\overline{\mfM}_{g,r})$
  of  $\mfO \mfp^\ZZZ_{\mfg, \hslash, \hslash \star \rho, g,r}/\overline{\mfM}_{g,r}$  according to decompositions arising from  various   clutching morphisms.
  We observe that the $\mbZ$-valued  function 
  given by $\rho \mapsto \mr{deg}(\mfO \mfp^\ZZZ_{\mfg, \hslash, \hslash \star \rho, g,r}/\overline{\mfM}_{g,r})$
   satisfies a certain factorization property, which will be called  a ``{\it pseudo-fusion rule}'' on $\mfc (\mbF_p)$.
  This notion may be  thought of  as a somewhat weaker variant of a fusion rule.
   The latter half of this chapter is mainly devoted to developing a general theory of pseudo-fusion rules.
  By applying results proved  there, we obtain a way to calculate 
  the value $\mr{deg}(\mfO \mfp^\ZZZ_{\mfg, \hslash, \hslash \star \rho, g,r}/\overline{\mfM}_{g,r})$
    combinatorially from the case of $(g,r) = (0, 3)$.
  Moreover,  the resulting pseudo-fusion  rule induces a commutative ring, which we call the {\it pseudo-fusion ring for dormant $\mfg$-opers}.
  The explicit understanding  of  its ring  structure  allows us to perform a computation of the generic degrees.

\section{Semi-graphs and clutching data}\label{y0172}

\subsection{Semi-graphs} \label{QG1059}

First, recall the definition of a semi-graph  (cf. ~\cite{Mzk3}, Appendix, or ~\cite{Mzk4}, \S\,1).

\bde \label{QW9930}
 \begin{itemize}
 \item[(i)]
 A  {\bf  semi-graph} 
 \index{$\GRRR$, semi-graph} 
  is a triple
\begin{align}
\GR= (V, E, \zeta),
\end{align}
consisting of the following data:
 \begin{itemize}
 \item
   a set  $V$, whose elements are
 called {\bf vertices}; 
 \index{$\GRRR$, semi-graph@--- vertex of}
 \item
a set  $E$,   whose elements are
  called {\bf edges},  
  \index{$\GRRR$, semi-graph@--- edge} 
  consisting of   sets  with cardinality $2$ such   that $e  \neq e' \in E$ implies $e \cap e' = \emptyset$;
 \item
 a map $\zeta : \bigsqcup_{e \in E} e \migi V \sqcup \{  \circledast \}$ (where
  ``$\circledast$'' denotes an abstract symbol with 
  $\circledast \notin V$),
    called a {\bf coincidence map}, 
    \index{$\GRRR$, semi-graph@--- coincidence map of}
    such that $\zeta (e) \neq \{ \circledast\}$ for any $e \in E$.
 \end{itemize}
Each edge $e \in E$ with $\circledast \in \zeta (e)$ (resp., $\circledast \notin \zeta (e)$) is called {\bf open} (resp., {\bf closed}),  
\index{$\GRRR$, semi-graph@--- $E^{\mr{op}}$, set of open edges}
\index{$\GRRR$, semi-graph@--- $E^{\mr{cl}}$, set of closed edges} 
and
write $E^{\mr{op}}$ (resp., $E^{\mr{cl}}$) for the set of open  (resp., closed) edges in $E$.
(In particular, we have $E = E^{\mr{op}} \sqcup E^{\mr{cl}}$.)
Also, for each edge $e  \in E$,  we shall refer to any element $b$ of $e$   as a {\bf branch} \index{$\GRRR$, semi-graph@--- branch}  of  $e$.
For each 
 $v\in V \sqcup \{ \circledast\}$, we write
 \begin{align}
 B_v  :=  \zeta^{-1} (\{ v\}).
  \index{$B_v$}
 \end{align}
  \item[(ii)]
A {\bf graph} \index{$\GRRR$, semi-graph@--- graph} is a semi-graph with $B_{\circledast}  \left(= \zeta^{-1} (\{ \circledast  \})\right)= \emptyset$, or equivalently, with $E^{\mr{op}} = \emptyset$.
  \end{itemize}
  \ede

\bde  \label{QH1055}
Let  $\GR= (V, E, \zeta)$ be   a  semi-graph. 
\begin{itemize}
\item[(i)]
We shall say that $\GR$ is {\bf finite} \index{$\GRRR$, semi-graph@--- finite} if both $V$ and $E$ are finite.
\item[(ii)]
We shall say that  $\GR$ is {\bf connected}  \index{$\GRRR$, semi-graph@--- connected} if for any two distinct  vertices $u$, $v \in V$,
 there exists  a sequence of edges $e_1, e_2, \cdots, e_l \in E$ ($l \geq 1$)  such that  $u \in \zeta (e_1)$, $v \in  \zeta (e_l)$, and $\zeta (e_j) \cap \zeta (e_{j+1}) \neq \emptyset$ for any $j = 1, \cdots, l-1$.
\end{itemize}
  \ede

\begin{rema}[Semi-graph as a topological space] \label{QH1001}
We have defined the notion of a semi-graph
  in a purely set-theoretic formulation, but
each semi-graph $\GR= (V, E, \zeta)$ can be regarded as a {\it topological space} in the usual way.
Indeed, we realize each closed edge $e := \{ b, b' \}\in E$ as a copy of the  closed unit interval $[0,1] \left( \subseteq \mbR \right)$ whose endpoints are indexed by $b$ and $b'$, and realize each open edge  $e := \{ b, b' \}$ (with $\zeta(b) \neq \circledast$) as a copy of $[0, 1) \left( \subseteq \mbR\right)$, where the point $0 \in \mbR$ is indexed by $b$.
Then, we can construct a topological space obtained from these constituents  in such a way that, for each $v \in V$, the points indexed by elements of $\zeta^{-1}(\{v \})$ are attached to form a single point.
The connectedness for semi-graphs defined in Definition \ref{QH1055}, (ii), is equivalent to the usual definition of connectedness for the corresponding   topological spaces.
Also, it makes sense to speak of the {\it first Betti number} $b_1 (\GR) := \mr{dim}_\mbQ (H_1 (\GR, \mbQ))$ of the semi-graph $\GR$.
Let us draw two examples of  semi-graphs as follows.
\vspace{-2mm}
\begin{center}
\begin{picture}(370,80)

\put(60,40){\circle*{5}}

\put(90,58){\circle{5}}
\put(60,40){\line(5,3){27.9}}
\put(66.5,73){\circle{5}} 
\put(60,40){\line(1,5){6}}
\put(40,65){\circle{5}}
\put(60,40){\line(-4,5){18.5}}
\put(27,45.8){\circle{5}}
\put(60,40){\line(-6,1){30.8}}
\put(90,22){\circle{5}}
\put(60,40){\line(5,-3){28}}
\put(66,10){\circle{5}}
\put(60,40){\line(1,-5){5.5}}
\put(42,13){\circle{5}}
\put(60,40){\line(-2,-3){16.4}}

\put(38,31){\rotatebox[origin=c]{116}{$\cdots$}}

\put(60,40){\line(1,0){80}}

\put(140,40){\circle*{5}}
\put(110,58){\circle{5}}\put(140,40){\line(-5,3){27.9}}
\put(134,72){\circle{5}}\put(140,40){\line(-1,5){5.9}}
\put(165,15){\circle{5}}\put(140,40){\line(1,-1){23.3}}
\put(175,34){\circle{5}}\put(140,40){\line(6,-1){32.1}}
\put(110,22){\circle{5}}\put(140,40){\line(-5,-3){27.9}}
\put(134,10){\circle{5}}\put(140,40){\line(-1,-5){5.5}}
\put(170,60){\circle{5}}\put(140,40){\line(3,2){27.5}}

\put(141,55){\rotatebox[origin=c]{-20}{$\cdots$}}

\put(260,40){\circle*{5}}
\put(275,70){\circle{5}}\put(260,40){\line(1,2){13.9}}
\put(247,72){\circle{5}}\put(260,40){\line(-2,5){11.9}}
\put(228,53){\circle{5}}\put(260,40){\line(-5,2){29.8}}
\put(247,8.5){\circle{5}}\put(260,40){\line(-2,-5){11.7}}
\put(276,8){\circle{5}}\put(260,40){\line(1,-2){14.8}}
\put(231,21){\circle{5}}\put(260,40){\line(-3,-2){26.5}} 

\put(80,-5){Figure 1}
\put(260,-5){Figure 2}

\qbezier(260,40)(320, 70)(340,50)
\qbezier(260,40)(320, 10)(340,30)
\qbezier(340,50)(350, 40)(340,30)

\put(235,36){\rotatebox[origin=c]{95}{$\cdots$}}

\end{picture}

\vspace{2mm}
\end{center}
(Here, ``$\circ\hspace{-1.5mm}-\hspace{-1.5mm}-\hspace{-1.5mm}\bullet$'' represents an open edge and 
``$\bullet\hspace{-1.5mm}-\hspace{-1.5mm}-\hspace{-1.5mm}\bullet$'' represents a closed edge.)
A finite connected semi-graph  having  only one closed edge is visualized as either of these two figures.
\end{rema}


\subsection{The dual semi-graph of a pointed (pre)stable curve} \label{R001}
 Let $k$ be an algebraically closed field, $(g,r)$ a pair of nonnegative integers with $2g-2+r>0$,  and
  $\msX := (X/k, \{ \sigma_i \}_{i=1}^r)$  an $r$-pointed prestable curve   of genus $g$  over $k$.
Then,  one can associate to $\msX$   a  semi-graph 
\begin{align}
\GR_{\msX} := (V_{\msX}, E_{\msX}, \zeta_{\msX})
\end{align}
defined as follows.
$V_{\msX}$ is taken as   the set of irreducible components of $X$ and  $E_{\msX}$ is  the disjoint union $N_{\msX} \sqcup  \{ \sigma_i \}_{i=1}^r$
 of  the set of nodal points  $N_{\msX}$ ($\subseteq X (k)$)  and the set of marked points $\{ \sigma_i \}_{i=1}^r$.
 Here,  note that any  nodal point  $e \in X(k)$ has two distinct  branches $b$ and $b'$, each of which lies on some well-defined irreducible component of $X$.
(Even when the two branches $b$, $b'$ lie on the same component, they are still treated differently from each other.)
 We identify  $e$ 
  with  the set  $\{ b, b'\}$  and,   moreover,  identify each marked point $\sigma_i$ with the set $\{ \sigma_i,  \sigma_i^\circledast \}$ (where ``$\sigma_i^\circledast$'' denotes an abstract symbol not equal to $\sigma_i$).
 In this way, we consider the  elements of $E_{\msX}$ as  sets  with cardinality $2$.
Finally, let $\zeta_{\msX}$ denote  the  map  $\bigsqcup_{e \in E_{\msX}} e  \migi V_{\msX} \sqcup \{ \circledast \}$ determined in the following manner:
\begin{itemize}
\item
If $b \in e$ (for some $e \in N_{\msX}$) or $b = \sigma_i$ (for some $i \in \{1, \cdots, r \}$), then
$\zeta_{\msX} (b)$ is the irreducible component  (i.e., an element of $V_{\msX}$) in which $b$ lies;
\item
 $\zeta_{\msX} (\sigma_i^\circledast) = \circledast$ for any $i =1, \cdots, r$.
\end{itemize}
(In particular, we have $E_\msX^{\mr{op}} = \{ \sigma_i \}_{i=1}^r$ and $E_\msX^{\mr{cl}} = N_\msX$.)
It is immediate that $\GR_\msX$ is finite and connected.
Moreover, $\GR_{\msX}$ becomes a {\it graph} if  and only if $r=0$, i.e.,  $\msX$ is unpointed.
We shall refer to $\GR_{\msX}$ as the {\bf dual semi-graph} associated to $\msX$. \index{$\GRRR$, semi-graph@--- $\GRRR_\msX$, dual semi-graph associated to $\msX$}

The semi-graph  $\GR_{\msX}$ is equipped with a weight function $h_\msX : V_\msX \migi \mbZ_{\geq 0}$ given by assigning,    to each irreducible component of $X$, the  genus of its normalization.
Then, $\msX$ is stable if and only if the weighted semi-graph $(\GR_\msX, h_\msX)$ is stable in the sense of ~\cite{ACP}, \S\,3,  i.e., the inequality $2 h_\msX (v)-2 + \sharp (B_v)>0$ holds for each $v \in V_\msX$.

\subsection{Clutching data} \label{QH70}
Let $n$ be a positive integer, and
 suppose that we are given a collection of data:
\begin{equation}   \msG :=  ( \GR, \{ (g_j, r_j) \}_{j=1}^n, \{ \lambda_j \}_{j=1}^n),   
\index{$\msG$, clutching data}
\end{equation}
where
\begin{itemize}
\item
$\GR := (V, E, \zeta)$ is  a finite connected  semi-graph with  $n$ vertices $V = \{ v_1, \cdots, v_n \}$, numbered $1$ through $n$;
\item
$(g_j, r_j)$ ($j = 1, \cdots, n $) is  a pair of nonnegative integers with $2g_j -2 + r_j >0$;
\item
$\lambda_j$ ($j = 1, \cdots, n$) denotes  a bijection of sets  $\lambda_j : B_{v_j} \isom  \{ 1, \cdots, r_j \}$.
\end{itemize} 
\vspace{3mm}
 
Let $k$
\index{$k$, field} 
 be a field and $\overline{k}$ an algebraically closed field over $k$.
Also,   let   $\{ {\msX_j} \}_{j=1}^n$  be   an ordered set of pointed stable curves over $\overline{k}$,
where  the $j$-th  curve 
${\msX_j} := (X_j/\overline{k}, \{ \sigma_{j, i} \}_{i})$ is $r_j$-pointed and of genus $g_j$ 
  and moreover the underlying  curve $X_j$ is {\it irreducible}.
Once  a collection of data $\msG$ described above has been specified, we may use it to glue together 
$\{ {\msX_j} \}_{j =1}^n$ and  obtain a new pointed stable curve $\msX := (X/\overline{k}, \{\sigma_i  \}_{i})$
in such a way that 
\begin{itemize}
\item
the dual semi-graph associated to $\msX$ can be identified with  $\GR$;
\item
the $j$-th vertex $v_j$ corresponds, as an irreducible component of $X$, to $X_j$; 
\item
if $e = \{ b, b'\}  \in E^{\mr{cl}}$ is a closed  edge with 
$\zeta (b) =v _{j}$, $\zeta (b') = v_{j'}$ (for some $j$, $j' \in \{1, \cdots, n \}$),
then $e$ corresponds to a node of $X$ obtained by gluing together $X_{j}$ at the $\lambda_{j} (b)$-th marked point $\sigma_{j, \lambda_{j}(b)}$ to $X_{j'}$ at the $\lambda_{j'} (b')$-th  marked point $\sigma_{j', \lambda_{j'}(b')}$;
\item
for each $j \in \{1, \cdots, n \}$,  the set of open edges $e \in E^{\mr{op}}$ with $\zeta (e) = \{ v_j,  \circledast \}$ is in bijection with the set of marked points  of $\msX$ lying on $X_j$;
\item
we order the marked points $\{\sigma_i \}_{i}$ of $\msX$ {\it lexicographically}, i.e., a marked point lying on $X_j$ comes before a marked point lying on $X_{j'}$ (if $j <j'$), and among marked points lying on $X_j$, we take the ordering induced by the original ordering of marked points on $X_j$. 
\end{itemize}
\vspace{3mm}
The genus ``$g$",  and the number of marked points ``$r$", of the resulting pointed curve $\msX$ may be computed combinatorially from the collection of data $\msG$.
(To be precise, the genus $g$ is given by the formula $g= b_1 (\GR)+ \sum_{j=1}^n g_i$ with the notation in Remark \ref{QH1001}.)
One may extend, in an evident way, this construction to the case where 
 ${\msX_i}$'s  are  (possibly {\it reducible}) pointed  stable curves  over an {\it arbitrary scheme} $S$ over $k$.
Moreover, by carrying out the above construction in the universal case, we obtain a morphism between moduli stacks
\index{$\mr{Clut}_\msG$, clutching morphism}
\begin{equation}
\label{1030}
  \mr{Clut}_\msG : \prod_{j=1}^n \overline{\mfM}_{g_j, r_j} \migi \overline{\mfM}_{g,r}
  \end{equation}
 given by assigning $\{ \msX_{j} \}_{j=1}^n \mapsto \msX$, where the product ``$\prod$''  is  taken over $k$.

\bde \label{y0173} 
We shall refer  (cf. ~\cite{Mzk2}, Chap.\,I, Definition 3.11) to such a collection  of data $\msG$ as  {\bf clutching data} \index{$\msG$, clutching data} of type $(g,r)$
  and to $\mr{Clut}_\msG$ as the {\bf clutching morphism} \index{$\mr{Clut}_\msG$, clutching morphism}  associated to $\msG$.
 \ede

\section{Gluing $\hslash$-flat $\mbG$-bundles}\label{y0174}

In  what follows, we shall discuss  the  procedure for gluing together $\hslash$-flat  $\mbG$-bundles by means of clutching data.

Let $k$
\index{$k$, field} 
 be a field,  $\mbG$
\index{$\mbG$, algebraic group}
 a connected smooth affine algebraic group over $k$, $S$ a $k$-scheme, and  $\hslash$ an element of $\Gamma (S, \mcO_S)$.
Denote by $\mfg$
\index{$\mfg$, Lie algebra}
 the Lie algebra of $\mbG$.
Also, let $n$ be a positive integer and  $\msG := (\GR, \{ (g_j, r_j)\}_{j=1}^n, \{ \lambda_j\}_{j=1}^n)$  clutching data  of type $(g,r)$.

\subsection{} \label{QG1070}
Suppose that, for each $j = 1, \cdots, n$, we are given an $r_j$-pointed stable curve ${\msX_j} := (X_j /S_j, \{ \sigma_{j, i} \}^{r_j}_{i=1})$   of genus $g_j$  over $S_j = S$.
Denote by $\msX := (X/S, \{ \sigma_i \}_{i=1}^r)$  the pointed stable curve obtained by gluing together $\{ {\msX_j} \}_{j=1}^n$ by means of $\msG$ in the  way explained above.
The clutching morphism $\mr{Clut}_\msG$ associated to $\msG$ gives  natural morphisms 
\begin{align}
\mr{Clut}_j : X_j \migi X
\end{align}
 ($j =1, \cdots, n$)  of $S \left(=S_j\right)$-schemes.
Here, recall  from \S\,\ref{QR781} that, for each $j$,   
the $k$-scheme $S_j$ (resp., $X_j$) is equipped with  the   log structure  pulled-back from  $\overline{\mfM}_{g_j,r_j}^\mr{log}$ (resp., $\mfC_{g_j,r_j}^\mr{log}$)  via the classifying morphism of ${\msX_j}$.
Denote by  $S_j^{\mr{log}}$ (resp., $X^{\mr{log}}_j$) and  the resulting log scheme.
On the other hand, 
we shall write, as usual,  $S^{\mr{log}}$ (resp.,  $X^\mr{log}$) for  the log scheme obtained by equipping $S$ (resp., $X$)  with the log structure pulled-back from $\overline{\mfM}_{g,r}^\mr{log}$ (resp., $\mfC_{g,r}^\mr{log}$) via the classifying morphism of $\msX$.
For simplicity, we shall write $\mcT := \mcT_{X^\mr{log}/S^\mr{log}}$ and $\mcT_j := \mcT_{X_j^\mr{log}/S_j^\mr{log}}$ ($j =1, \cdots, n$).
Also, write $X^{\mr{log} |_X}_j$ for the log scheme obtained by equipping $X_j$ with the log structure pulled-back from  $X^\mr{log}$ via $\mr{Clut}_j$.
The structure morphism $X_j \migi S_j \left(=S \right)$ of  $X_j/S_j$ extends to a morphism $X_j^{\mr{log} |_X} \migi S^{\mr{log}}$ of log schemes (but we note that it may not be a log curve).
Moreover, the natural morphisms $S^\mr{log} \migi S \left(=S_j \right)$ and  $X_j^{\mr{log} |_X} \migi X_j$ 
extend to
 morphisms $S^\mr{log} \migi S_j^{\mr{log}}$ and   $X_j^{\mr{log} |_X} \migi X_j^{\mr{log}}$ respectively, which make the following square diagram commute:
\begin{align}  \label{fff01}
\vcenter{\xymatrix@C=46pt@R=36pt{
 X_j^{\mr{log} |_X} \ar[r] \ar[d]_-{\mr{projection}} &X_j^{\mr{log}}  \ar[d]^-{\mr{projection}}\\
S^{\mr{log}} \ar[r]&S_j^{\mr{log}}.
}}
\end{align}
This diagram induces  a morphism 
$X_j^{\mr{log} |_X} \migi X_j^{\mr{log}} \times_{S_j^{\mr{log}}} S^{\mr{log}}$
   over $S^{\mr{log}}$ 
   whose
underlying morphism of  $S$-schemes  coincides with the identity morphism of $X_j$.
The differential of this morphism  yields  an isomorphism of $\mcO_{X_j}$-modules
 \begin{align} \label{fff02}
  \delta_{X, j} : \mcT_{X_j^{\mr{log} |_X}/S^\mr{log}} \left(= \mr{Clut}_{j}^*(\mcT)\right)
  \isom  \mcT_j.
 \end{align}

 \subsection{} \label{QG1071}
 Now, let  $\mcE$ be a $\mbG$-bundle  on $X$, and for each $j = 1, \cdots, n$, denote by $\mcE_j$  the pull-back  of $\mcE$ to $X_j$.
Also, denote by $\mcE_j^{\mr{log}}$ (resp., $\mcE_j^{\mr{log} |_X}$) the log scheme defined to be $\mcE_j$ equipped with the log structure pulled-back from $X_j^{\mr{log}}$ (resp., $X_j^{\mr{log} |_X}$), i.e., $\mcE_j^{\mr{log}} := \mcE_j \times_{X_j} X_j^{\mr{log}}$ (resp., $\mcE_j^{\mr{log} |_X} := \mcE_j \times_{X_j} X_j^{\mr{log} |_X}$).
The morphism $\mcE_j^{\mr{log}|_X} \migi \mcE_j^\mr{log}$ induced by the upper horizontal arrow in (\ref{fff01}) gives  
 an isomorphism
\begin{align} \label{QW9901}
 \delta_{\mcE, j} :  \widetilde{\mcT}_{\mcE_j^{\mr{log}|_X}/S^\mr{log}} \left(= \mr{Clut}_j^*(\widetilde{\mcT}_{\mcE^\mr{log}/S^\mr{log}}) \right)\isom \widetilde{\mcT}_{\mcE_j^\mr{log}/S_j^\mr{log}},
\end{align}
that is compatible with (\ref{fff02}) via the respective anchor maps.

 If 
 $e := \{ b, b' \} \in E^{\mr{cl}}$ is  a closed  edge with   $\zeta (b) = v_{j}$,  $\zeta (b') = v_{j'}$ (for some $j$, $j' \in \{ 1, \cdots, n \}$),
 then 
 the equality $\mr{Clut}_{j} \circ \sigma_{j, i} = \mr{Clut}_{j'} \circ \sigma_{j',  i'}$ of morphisms $S \migi X$ holds, where $i:= \lambda_j (b)$ and $i':= \lambda_{j'} (b')$.
Hence, 
 \begin{align}
 \gamma_{X, e}:= \sigma_{j', i'}^*(\delta_{X, j'}) \circ \sigma_{j,  i}^*(\delta_{X, j})^{-1} \hspace{2mm}
 \text{and} \hspace{2mm}
  \widetilde{\gamma}_{\mcE, e}:= \sigma_{j', i'}^*(\delta_{\mcE, j'}) \circ \sigma_{j,  i}^*(\delta_{\mcE, j})^{-1}
 \end{align}
 define isomorphisms fitting into the following isomorphism of short exact sequences:
 \begin{align} \label{QW9907}
\vcenter{\xymatrix@C=16pt@R=36pt{
0 \ar[r] & \sigma_{j, i}^*(\mfg_{\mcE_{j}})  \ar[r]\ar[d]_-{\wr}^{\gamma_{\mcE, e}}
& \sigma_{j, i}^*(\widetilde{\mcT}_{\mcE_{j}^{\mr{log}}/S_j^{\mr{log}}})  \ar[r]\ar[d]_{\wr}^-{\widetilde{\gamma}_{\mcE, e}}
& \sigma_{j, i}^*(\mcT_j) \ar[r] \ar[d]_{\wr}^-{\gamma_{X, e}} & 0 \\
0 \ar[r]&  \sigma_{j', \lambda_{j'}(b')}^*(\mfg_{\mcE_{j'}})   \ar[r]& \sigma_{j', i'}^*(\widetilde{\mcT}_{\mcE_{j'}^{\mr{log}}/S_{j'}^{\mr{log}}})  \ar[r] 
& \sigma_{j', i'}^*(\mcT_{j'}) \ar[r] & 0, 
}}
\end{align}
where $\gamma_{\mcE, e}$ is  obtained by restricting   $\widetilde{\gamma}_{\mcE, e}$ and the upper and lower horizontal sequences are the pull-backs of (\ref{Ex0}) for the cases    $\mcE = \mcE_j$ and $\mcE_{j'}$ respectively.

 \bde \label{QG2}
 A  {\bf $\msG$-$\hslash$-$S^\mr{log}$-connection} \index{$\msG$-$\hslash$-$S^\mr{log}$-connection} on $\mcE$ is a collection
\begin{align}
\{ \nabla_j \}_{j=1}^n,
\end{align}
where
    each $\nabla_j$ denotes  an $\hslash$-$S_j^\mr{log}$-connection on $\mcE_j$,
such that, for  any closed edge  $e = \{ b, b' \} \in  E^{\mr{cl}}$ with $\zeta (b) = v_{j}$,  $\zeta (b') = v_{j'}$ ($j, j \in \{1, \cdots, n \}$), 
 the equality 
 $\mu_{\lambda_j (b)}^{\nabla_{j}} = -\mu_{\lambda_{j'}(b')}^{\nabla_{j'}}$
  holds under the identification $\Gamma (S, \sigma_{j, \lambda_j (b)}^*(\mfg_{\mcE_{j}})) = \Gamma (S, \sigma_{j', \lambda_{j'} (b')}^*(\mfg_{\mcE_{j'}}))$ given by $\gamma_{\mcE, e}$.
 \ede

 \subsection{} \label{QG1074}
 Let  $\nabla : \mcT \migi \widetilde{\mcT}_{\mcE^\mr{log}/S^\mr{log}}$   an $\hslash$-$S^\mr{log}$-connection  on $\mcE$.
Then, for each $j =1, \cdots, n$, 
the composite
\begin{align}
\nabla |_{j} : \mcT_j \xrightarrow{\delta_{X, j}^{-1}} \mcT_{X_j^{\mr{log}|_X}/S^\mr{log}} \xrightarrow{\mr{Clut}^*_j (\nabla)}\widetilde{\mcT}_{\mcE_j^{\mr{log} |_X} /S^\mr{log}} \xrightarrow{\delta_{\mcE, j}} \widetilde{\mcT}_{\mcE_j^{\mr{log}}/S_j^{\mr{log}}}
\end{align}
   specifies an $\hslash$-$S_j^\mr{log}$-connection on $\mcE_j^{\mr{log}}$. We shall refer to $\nabla |_j$ as the {\bf restriction} of $\nabla$ to $\mcE_j$.
 
\bpr \label{y0166pfp}
For each $\hslash$-$S^\mr{log}$-connection $\nabla$ on $\mcE$,
the collection  of its restrictions $\{ \nabla |_j \}_{j=1}^n$ specifies
a $\msG$-$\hslash$-$S^\mr{log}$-connection on $\mcE$.
 Moreover, the resulting assignment $\nabla \mapsto \{ \nabla |_j \}_{j=1}^n$ determines a bijection of sets
\begin{align}
 \begin{pmatrix}
\text{the set of } \\
\text{$\hslash$-$S^\mr{log}$-connections on $\mcE$} 
\end{pmatrix}
\isom 
\begin{pmatrix}
\text{the set of } \\
\text{$\msG$-$\hslash$-$S^\mr{log}$-connections on $\mcE$} 
\end{pmatrix}.
\end{align}
 \epr
\begin{proof}
We shall  consider the former assertion.
Let us take an arbitrary closed edge $e := \{ b_1, b_2 \}$ with $\zeta  (b_1) = v_{j_1}$, $\zeta (b_2) = v_{j_2}$.
Write $i_1 := \lambda_{j_1} (b_1)$ and $i_2 := \lambda_{j_2} (b_2)$.
To verify that  $\{ \nabla |_j \}_{j=1}^n$  forms a $\msG$-$\hslash$-$S^\mr{log}$-connection, we may assume, after possibly replacing $S$ with its open covering,   that 
for each $l =1, 2$,  there exists  a pair $(U_{l}, x_l)$, where
\begin{itemize}
\item
$U_{l}$  denotes  an affine open subscheme  of $X_{j_l}$ with  $\mr{Im} (\sigma_{j_l, i_l}) \subseteq U_{l}$;
\item
$x_l$ denotes   a global  section  in $\mcO_{U_{l}}$ defining    the closed subscheme  $\mr{Im}(\sigma_{j_l, i_l})$  and satisfying  that $\mcO_{U_{l}} \cdot d \mr{log} (x_l) = \Omega_{X_{j_l}^{\mr{log}}/S_{j_l}^{\mr{log}}} |_{U_{l}}$.
\end{itemize}
The closed immersion  $\sigma_{j_l, i_l}$  ($l =1, 2$)
 corresponds to a surjection  $\pi_l : \Gamma (U_{l}, \mcO_{U_{l}}) \migisurj \Gamma (S, \mcO_S)$.
   If $U$ denotes a unique  open subscheme of $X$ satisfying that $U \times_X X_{j_1} = U_{1}$ and $U \times_X X_{j_2} = U_{2}$, then we have  a natural   identification
\begin{align} \label{fff020}
\Gamma (U, \mcO_U) = \left\{ (a_1, a_2) \in \Gamma (U_{1}, \mcO_{U_{1}}) \times \Gamma (U_{2}, \mcO_{U_{2}})  \ | \ \pi_1 (a_1) = \pi_2 (a_2)   \right\}.
\end{align}
For  $l =1, 2$,  the composite of the   inclusion $\Gamma (S, \mcO_S) \migiincl  \Gamma (U_{l}, \mcO_{U_{l}})$ followed by  $\pi_l$    coincides with the identity morphism of $\Gamma (S, \mcO_S)$. 
Hence, the above identification allows us to consider $\Gamma (U_{1}, \mcO_{U_{1}})$
 (resp., $\Gamma (U_{2}, \mcO_{U_{2}})$)
   as a subring of  $\Gamma (U, \mcO_U)$ via
the inclusion $\Gamma (U_{1}, \mcO_{U_{1}}) \migiincl \Gamma (U, \mcO_U)$  (resp., $\Gamma (U_{2}, \mcO_{U_{2}}) \migiincl \Gamma (U, \mcO_U)$)
 given by assigning $a \mapsto (a, \pi_1 (a))$ (resp., $a \mapsto (\pi_2 (a), a)$).
If we regard $x_l$ (for each $l=1,2$)  as  an element of  $\Gamma (U, \mcO_U)$, then  
 the induced section $d \mr{log} (x_{l}) \in \Gamma (U, \Omega)$
     forms a generator  of $\Omega |_U$.
Write  $\partial_{l} \in \Gamma (U, \mcT)$
  for  its  dual; 
  we abuse notation to use $\partial_l$ for writing 
  the dual of $d \mr{log} (x_{l})$ as an element of $\Gamma (U_{l}, \mcT_{j_l})$.
Since the equality  $x_{1} x_{2} = 0$ holds in $\Gamma (U, \mcO_U)$,  we have 
$d \mr{log} (x_{1}) + d \mr{log} (x_{2}) = 0$, and hence,  $\partial_{1} = - \partial_{2}$.
It follows that
\begin{align}
(\sigma_{j_1, i_1})^*(\nabla |_{j_1} (\partial_{1})) =  
& \ 
(\mr{Clut}_{j_1} \circ \sigma_{j_1, i_1})^*(\nabla (\partial_{1})) \\
= & \ - (\mr{Clut}_{j_1} \circ \sigma_{j_1, i_1})^*(\nabla (\partial_{2})) \notag \\
= & \ - (\mr{Clut}_{j_2} \circ \sigma_{j_2, i_2})^*(\nabla (\partial_{2})) \notag \\
= & \ - (\sigma_{j_2,  i_2})^*(\nabla |_{j_2} (\partial_{2})). \notag
\end{align}
Hence, the collection $\{ \nabla |_j \}_j$ turns out to be a $\msG$-$\hslash$-$S^\mr{log}$-connection on $\mcE$.

The latter assertion follows from  the discussion in the proof of the former assertion.
In fact, in order to determine the inverse to the assignment $\nabla \mapsto \{ \nabla |_j \}_j$,  we may reverse the steps in the construction
  by taking  account of 
  canonical isomorphisms
\begin{align} \label{fff07}
 \Gamma  (U, \mcT)  &  
\isom   \left\{ (a_l)_{l=1}^2 
  \,  \big| \,  a_l \in \Gamma (U_{l}, \mcT_{j_l}),  \gamma_{X, e} (\sigma_{j_1, i_1}^* (a_1)) =    \sigma_{j_2,  i_2}^* (a_2)\right\},   \\
  \Gamma (U, \widetilde{\mcT}_{\mcE^\mr{log}/S^\mr{log}}) 
  &\isom    \left\{ (a_l)_{l=1}^2 \, \Big| \,
 a_l \in \Gamma (U_{l}, \widetilde{\mcT}_{\mcE_{j_l}^{\mr{log}}/S_{j_l}^{\mr{log}}}), 
 \widetilde{\gamma}_{\mcE, e}(\sigma_{j_1, i_1}^* (a_1)) = \sigma_{j_2, i_2}^* (a_2) \right\}.   \notag
\end{align}
This completes the proof of the assertion.
\end{proof}

\subsection{} \label{QG10}
Let us suppose that $k$ is of  characteristic $p>0$.
Let $\nabla$ be an $\hslash$-$S^\mr{log}$-connection on $\mcE$.
Note that, for each $j =1, \cdots, n$,
 both $\delta_{X, j}$ and $\delta_{\mcE, j}$ preserve  the $p$-power operations.
Hence, the $p$-curvature ${^p}\psi^{\nabla |_j}$ of the restriction $\nabla |_j$ fits into the following commutative diagram:
\begin{align} \label{Eh5007}
\vcenter{\xymatrix@C=56pt@R=36pt{
F_{X_j}^*(\mcT_{X_j^{\mr{log}|_X}/S^\mr{log}})\ar[r]^-{\mr{Clut}_j^*({^p}\psi^\nabla)} \ar[d]^-{\wr}_-{F_{X_j}^*(\delta_{X, j})}& \widetilde{\mcT}_{\mcE_j^{\mr{log}|_X}/S^\mr{log}} \ar[d]_-{\wr}^-{\delta_{\mcE, j}}\\
F_{X_j}^*(\mcT_j) \ar[r]_-{\psi^{\nabla |_j}} & \widetilde{\mcT}_{\mcE_j^\mr{log}/S_j^\mr{log}}.
}}
\end{align}
In particular, $\nabla$ is $p$-flat if and only if $\nabla |_j$ is $p$-flat for every  $j=1, \cdots, n$.

\section{Factorization of  $(\mfg, \hslash)$-opers}\label{y0174s}

In this section, we shall discuss the gluing procedure for  $(\mfg, \hslash)$-opers by means of clutching data.
We shall keep  the notations introduced at the beginning of 
the previous section.
Suppose further that $k$ is perfect and $(\mbG, \mbT)$
\index{$\mbG$, algebraic group}
\index{$\mbG$, algebraic group@— $\mbT$, maximal torus of}
 specifies  a split semisimple algebraic group over $k$ of adjoint type, where $\mbT$ denotes a maximal torus of $\mbG$.
Denote by $\mfg$
\index{$\mfg$, Lie algebra}
 the Lie algebra of $\mbG$.
Also, suppose that  {\it either   ``\,$\mr{char}(k) =0$''  or ``\,$\mr{char}(k) > 2  h_\mbG$'' is fulfilled}.
Next, let $\mbB$
\index{$\mbG$, algebraic group@— $\mbB$, Borel subgroup of}
 be a Borel subgroup defined over $k$ containing $\mbT$ and 
 $\ \qq := \{ x_\alpha \}_{\alpha \in \Gamma}$ 
 \index{$\ \qq\,$}
  a collection of generators $x_\alpha$ of $\mfg^\alpha$'s as in (\ref{QG5697}).


\subsection{The defnition of a  $\msG$-$(\mfg, \hslash)$-oper} \label{QG1090}

Let $S$, $\msX_j$ ($j=1, \cdots, n$), and $\msX$ be as in \S\,\ref{QG1070}.

\bde \label{y0175}
A {\bf  set of  radii for $\msG$} over $S$ is an ordered  set 
 \index{$\rho_\msG$, set of  radii for $\msG$}
\begin{equation}
 \rho_\msG := \{  \rho^j \}_{j=1}^n,
 \end{equation}
where each $\rho^j$ is an element $\rho^j := (\rho^j_i )_{i=1}^{r_j}$ of $\mfc^{\times r_j} (S)$ such   that for 
every closed edge $e = \{ b, b' \} \in E^{\mr{cl}}$ with  $\zeta (b) = v_{j}$, $\zeta (b') = v_{j'}$ (for some $j$, $j' \in \{ 1, \cdots, n \}$),
the equality  $\rho^{j}_{\lambda_{j} (b)} = (-1)\star \rho^{j'}_{\lambda_{j'} (b')}$ holds.
 \ede

\bde \label{QG8}
\begin{itemize}
\item[(i)]
A {\bf $\msG$-$(\mfg, \hslash)$-oper} \index{$\msG$-$(\mfg, \hslash)$-oper} on $\msX$ is a collection
\begin{align}
\msE_\msG^\spadesuit := (\mcE_\mbB, \{ \nabla_j \}_{j=1}^n),
\end{align}
where $\mcE_\mbB$ denotes a $\mbB$-bundle on $X$ and each $\nabla_j$ ($j=1, \cdots, n$) denotes an $\hslash$-$S_j^\mr{log}$-connection on the $\mbB$-bundle $(\mcE_\mbB \times^\mbB \mbG)_j \left(:= (\mcE_{\mbB})_j \times^\mbB \mbG\right)$, satisfying the following conditions:
\begin{itemize}
\item
For each $j =1, \cdots, n$,  the pair $\msE^\spadesuit_{\msG, j} := ((\mcE_\mbB)_j, \nabla_j)$ forms a $(\mfg, \hslash)$-oper on $\msX_j$.
\item
For each closed edge $e := \{ b, b' \} \in E^{\mr{cl}}$ with $\zeta (b) = v_j$, $\zeta (b') = v_{j'}$ (for some $j, j' \in \{1, \cdots, n \}$), the equality $\rho^{\nabla_j}_{\lambda_j (b)} = (-1)\star \rho^{\nabla_{j'}}_{\lambda_{j'}(b')}$ holds.
\end{itemize}
If, moreover, $\mr{char}(k) >0$, then we say that $\msE^\spadesuit_\msG$ is {\bf dormant} if the $(\mfg, \hslash)$-oper  $\msE^\spadesuit_{\msG, j}$ is dormant for every  $j$.
\item[(ii)]
Let $\msE_\msG^\spadesuit := (\mcE_\mbB, \{ \nabla_j \}_{j=1}^n)$,
${\msE_\msG'}^{\hspace{-0.8mm}\spadesuit} := (\mcE'_\mbB, \{ \nabla'_j \}_{j=1}^n)$
be $\msG$-$(\mfg, \hslash)$-opers on $\msX$.
Then, an {\bf isomorphism of $\msG$-$(\mfg, \hslash)$-opers}  from 
$\msE_\msG^\spadesuit$ to ${\msE_\msG'}^{\hspace{-0.8mm}\spadesuit}$ is defined as an isomorphism $\mcE_\mbB \isom \mcE'_\mbB$ of $\mbB$-bundles such  that, for every $j=1, \cdots, n$, the induced isomorphism $(\mcE_\mbB)_j \isom (\mcE'_\mbB)_j$ defines an isomorphism of $(\mfg, \hslash)$-opers $\msE^\spadesuit_{\msG, j} \isom {\msE_{\msG, j}'}^{\hspace{-2.8mm}\spadesuit}\,$.
\end{itemize}
\ede

\bde
Let $\rho_\msG := \{\rho^j \}_{j=1}^n$ be a set of radii for $\msG$ and 
$\msE_\msG^\spadesuit := (\mcE_\mbB, \{ \nabla_j \}_{j=1}^n)$ a $\msG$-$(\mfg, \hslash)$-oper on $\msX$.
Then, we shall say that $\msE_\msG^\spadesuit$ is {\bf of radii $\rho_\msG$} \index{$\msG$-$(\mfg, \hslash)$-oper@--- of radii $\rho_\msG$} if
for any $j=1, \cdots, n$, the $(\mfg, \hslash)$-oper  $\msE^\spadesuit_{\msG, j}$ is of radii $\rho^j$.
\ede

\subsection{$\ \qq$-normal $\msG$-$(\mfg, \hslash)$-opers} \label{QH90}

For each $j =1, \cdots, n$ and $\Box \in \{ \mbB, \mbG \}$, we shall write
\begin{align} \label{fff0111}
\DE_\Box := \DE_{\Box, \hslash, X^\mr{log}/S^\mr{log}},
\hspace{3mm} 
\DE_{\Box}|_j  := X_j \times_{\mr{Clut}_j, X}\DE_\Box, 
\hspace{3mm} 
\DE_{\Box, j}  := \DE_{\Box, \hslash, X_j^\mr{log}/S^\mr{log}_j}.
\end{align}

\ble \label{y0166pffp}
  There exists
 a unique  isomorphism 
\begin{align} \label{fff010}
\wp_{\mbB, j}  :\DE_{\mbB}|_j  \isom \DE_{\mbB, j} 
\end{align}
  of  $\mbB$-bundles satisfying the following property: for any   log charts $(U, \partial)$ and $(U_{j}, \partial_j)$ on $X^\mr{log}/S^\mr{log}$  and  $X_j^{\mr{log}}/S_j^{\mr{log}}$ respectively such that $U \times_X X_j = U_j$ and 
$\delta_{X, j}(\mr{Clut}_j^{-1}(\partial)) =\delta_j$,
 the following square diagram is commutative:
\begin{align} \label{QG4}
\vcenter{\xymatrix@C=60pt@R=36pt{
U_j \times_U (U \times_k \mbB)
\ar[r]_-{\sim}^-{\mr{id}_{U_j} \times \mr{triv}_{\mbB, (U, \partial)}}  \ar[d]_-{\wr}& U_j \times_U \DE_{\mbB} |_U \left(= \DE_{\mbB} |_{U_j} \right)   \ar[d]_-{\wr}^-{\wp_{\mbB, j}|_{U_j}}\\
U_j \times_k \mbB 
\ar[r]^-{\sim}_-{\mr{triv}_{\mbB, (U_j, \partial_j)}}
 &
\DE_{\mbB, j}|_{U_j}.
}}
\end{align}
  \ele
\begin{proof}
The assertion follows from  the constructions of  $\DE_{\mbB}|_j$,   $\DE_{\mbB, j}$, and  $\delta_{X, j}$.
\end{proof}

By applying  the  change of structure group to $\wp_{\mbB, j}$, we obtain  an isomorphism of $\mbG$-bundles
 \begin{equation} \label{fff025}
\wp_{\mbG, j}  : \DE_{\mbG} |_j \isom \DE_{\mbG, j}.
\end{equation}
Given an $\hslash$-$S^\mr{log}$-connection $\nabla$ on $\DE_\mbG$,
we obtain
an $\hslash$-$S_j^\mr{log}$-connection
\begin{align} \label{QG11}
{^\dagger}\nabla |_j := \wp_{\mbG, j*}(\nabla |_j) : \mcT_{j} \migi \widetilde{\mcT}_{\DE_{\mbG, j}^\mr{log}/S_j^\mr{log}}
\end{align}
on $\DE_{\mbG, j}$ for each $j =1, \cdots, n$.

 \bde \label{QH88}
 We shall say that  a $\msG$-$(\mfg, \hslash)$-oper  $\msE^\spadesuit_\msG := (\mcE_\mbB, \{ \nabla_j \}_{j=1}^n)$ is {\bf $\, \qq$-normal} \index{$\msG$-$(\mfg, \hslash)$-oper@--- $\, \qq$-normal}
  if $\mcE_\mbB = \DE_\mbB$ and  the $(\mfg, \hslash)$-opers $(\DE_{\mbB, j}, {^\dagger}\nabla |_j)$ ($j=1, \cdots, n$) are all $\ \qq$-normal.
 \ede

\bpr \label{y0166pp}
If $\msE^\spadesuit := (\DE_\mbB, \nabla)$ is a $\ \qq$-normal $(\mfg, \hslash)$-oper on $\msX$, then
 the collection 
\begin{align}
\msE^\spadesuit_{\emptyset \Rightarrow \msG} := (\DE_\mbB, \{ {^\dagger}\nabla |_j \}_{i=1}^n)
\end{align}
 forms a $\ \qq$-normal  $\msG$-$(\mfg, \hslash)$-oper on $\msX$.
 Moreover, 
   the resulting assignment $\msE^\spadesuit \mapsto \msE^\spadesuit_{\emptyset \Rightarrow \msG}$ determines an equivalence of categories 
 \begin{align}
 \begin{pmatrix}
\text{the  groupoid of  $\ \qq$-normal} \\
\text{$(\mfg, \hslash)$-opers on $\msX$} 
\end{pmatrix}
\isom 
\begin{pmatrix}
\text{the groupoid of $\ \qq$-normal} \\
\text{$\msG$-$(\mfg, \hslash)$-opers on $\msX$} 
\end{pmatrix}.
\end{align}
Finally, if $\mr{char}(k) >0$, then this equivalence  restricts to an equivalence of categories
\begin{align}
 \begin{pmatrix}
\text{the groupoid  of  $\ \qq$-normal} \\
\text{dormant $(\mfg, \hslash)$-opers on $\msX$} 
\end{pmatrix}
\isom 
\begin{pmatrix}
\text{the groupoid of  $\ \qq$-normal} \\
\text{dormant $\msG$-$(\mfg, \hslash)$-opers on $\msX$} 
\end{pmatrix}.
\end{align}
 \epr
\begin{proof}
We shall  consider the first assertion.
One verifies immediately  from  Lemma \ref{y0166pffp} that   the pair $(\DE_{\mbB, j}, {^\dagger}\nabla |_j)$ defined for each $j$ forms a $\ \qq$-normal $(\mfg, \hslash)$-oper.
The remaining portion is to 
 verify the consistency   of radii; it   can be verified 
  by the following sequence of equalities, where    $e := \{ b, b' \} \in E^{\mr{cl}}$ is 
a closed edge  with $\zeta (b) = v_j$, $\zeta (b') = v_{j'}$  and 
 we write $i := \lambda_j (b)$,   $i' := \lambda_{j'}(b')$:
\begin{align} \label{fff040}
\rho^{{^\dagger}\nabla |_j}_{i} &=  \rho^{\nabla |_j}_{i}   \\
&=  [\chi] ((\sigma^*_{j, i}(\mfg_{\DE_{\mbG}|_j}),  \mu_{i}^{\nabla |_j})) \notag  \\
& =  [\chi]
((\sigma^*_{j', i'}(\mfg_{\DE_{\mbG}|_{j'}}), - \mu_{i'}^{\nabla |_{j'}}))
 \notag \\
 & =   (-1) \star  [\chi]
((\sigma^*_{j', i'}(\mfg_{\DE_{\mbG}|_{j'}}),  \mu_{i'}^{\nabla |_{j'}}))
   \notag \\
& =  (-1) \star  \rho^{\nabla |_{j'}}_{i'}  \notag \\
&=   (-1) \star \rho^{{^\dagger}\nabla |_{j'}}_{i'},  \notag
\end{align}
  where
the third equality follows from 
 Proposition \ref{y0166pfp}.
Thus, we complete the proof of the former assertion.

Next,  we shall   prove  the second assertion by 
constructing  the inverse to the assignment   $\msE^\spadesuit \mapsto \msE^\spadesuit_{\emptyset \Rightarrow \msG}$.
Suppose that we are given a  $\ \qq$-normal  $\msG$-$(\mfg, \hslash)$-oper   $\msE^\spadesuit_{\msG}  :=(\DE_\mbB, \{ \nabla_j \}_{j=1}^n)$  on $\msX$.
Let $e$, $U$,  $j_l$, $i_l$, 
 $(U_{l}, x_l)$, and $\partial_l$ ($l =1, 2$)   be as in the proof of Proposition \ref{y0166pfp}.
 Since the $(\mfg, \hslash)$-opers  $(\DE_{\mbB, j_l}, \nabla_{j_l})$ ($l=1,2$) are $\ \qq$-normal, 
$\mu_{i_l}^{\mr{triv}_{\mbG, (U_l, \partial_l)}^*(\nabla_{j_l})} - 1 \otimes q_{-1}$ is contained in $\mfg^{\mr{ad}(q_1)} (S)$.
Moreover, if $\alpha$ denotes  the automorphism of $\mcO_S \otimes_k \mfg$ determined by $ \mr{Ad}_\mbG \left(\iota_\mbB \left(\overline{\begin{pmatrix} 1 & 0 \\ 0 & -1\end{pmatrix}}\right) \right)$, 
then  the element $\alpha(\mu_{i_1}^{\mr{triv}_{\mbG, (U_1, \partial_1)}^*(\nabla_{j_1})}))) + 1 \otimes q_{-1}$ is contained in $\mfg^{\mr{ad}(q_1)} (S)$.
Hence, we obtain the following sequence of equalities:
\begin{align} \label{fff032}
\mr{Kos}_\mfg  (\mu_{i_2}^{\mr{triv}_{\mbG, (U_2, \partial_2)}^*(\nabla_{j_2})} - 1 \otimes q_{-1} ) 
&= \rho_{i_2}^{\mr{triv}_{\mbG, (U_2, \partial_2)}^*(\nabla_{j_2})}    \\
&=   \rho_{i_2}^{\nabla_{j_2}}  \notag   \\
& =   (-1) \star \rho_{i_1}^{\nabla_{j_1}} \notag \\
&=(-1)\star \rho_{i_1}^{\mr{triv}_{\mbG, (U_1, \partial_1)}^*(\nabla_{j_1})} \notag \\
& = (-1)\star \chi (\mu_{i_1}^{\mr{triv}_{\mbG, (U_1, \partial_1)}^*(\nabla_{j_1})}) \notag \\
& = \chi (-\mu_{i_1}^{\mr{triv}_{\mbG, (U_1, \partial_1)}^*(\nabla_{j_1})}) \notag \\
&= \chi  (-\alpha(\mu_{i_1}^{\mr{triv}_{\mbG, (U_1, \partial_1)}^*(\nabla_{j_1})})) \notag \\
& = \mr{Kos}_\mfg ( -\alpha(\mu_{i_1}^{\mr{triv}_{\mbG, (U_1, \partial_1)}^*(\nabla_{j_1})}))) - 1 \otimes q_{-1}). \notag 
\end{align}
The bijectivity of  $\mr{Kos}_\mfg$  implies 
\begin{align} \label{fff061}
\mu_{i_2}^{\mr{triv}_{\mbG, (U_2, \partial_2)}^*(\nabla_{j_2})}  =  -\alpha(\mu_{i_1}^{\mr{triv}_{\mbG, (U_1, \partial_1)}^*(\nabla_{j_1})}).
\end{align}
For each $l = 1, 2$, denote by 
\begin{align}
\mr{triv}_l : \mcO_S \otimes \mfg \isom \sigma_{j_l, i_l}^* (\mfg_{\DE_{\mbG}|_{j_l}})  \left(= (\mr{Clut}_{j_l} \circ \sigma_{j_l, i_l})^*(\mfg_{\DE_{\mbG}})\right)
\end{align}
the $\mcO_S$-linear isomorphism induced, via restriction by $\sigma_{j_l, i_l}$, 
from  the differential of $\mr{triv}_{\mbG, (U, \partial_l)}$.
Since $\partial_1 = - \partial_2$ in $U$, 
 it follows from Proposition \ref{y0130hh2} that the equality $\mr{triv}_1 =  \mr{triv}_2 \circ \alpha$
holds under the identification $\sigma_{j_1, i_1}^* (\mfg_{\DE_{\mbG}|_{j_1}}) = \sigma_{j_2, i_2}^* (\mfg_{\DE_{\mbG}|_{j_2}})$ given by 
$\mr{Clut}_{j_1} \circ \sigma_{j_1, i_1} = \mr{Clut}_{j_2} \circ \sigma_{j_2, i_2}$.
Therefore,  by letting
\begin{align}
\nabla'_j :=  (\wp_{\mbG, j}^{-1})_*(\nabla_{j}) : \mcT_j \migi \widetilde{\mcT}_{\DE_\mbG |_j}
\end{align}
($j = 1, \cdots, n$), 
we obtain the following  sequence of equalities:
\begin{align} \label{fff033}
 \mu_{i_1}^{\nabla'_{j_1}} &=  \mr{triv}_1
  ( \mu_{i_1}^{\mr{triv}_{\mbG, (U_1, \partial_1)}^*(\nabla_{j_1})})   \\
 & = 
(\mr{triv}_2 \circ \alpha)   ( \mu_{i_1}^{\mr{triv}_{\mbG, (U_1, \partial_1)}^*(\nabla_{j_1})}) \notag \\
& \hspace{-2mm}\stackrel{(\ref{fff061})}{=} - \mr{triv}_2 (\mu_{i_2}^{\mr{triv}_{\mbG, (U_2, \partial_2)}^*(\nabla_{j_1})}) \notag \\
& = -\mu_{i_2}^{\nabla'_{j_2}}. \notag
\end{align}
This fact implies, by Proposition \ref{y0166pfp}, that there exists  a unique  $\hslash$-$S^\mr{log}$-connection $\nabla$ on $\DE_{\mbG}$ whose restriction to  $\DE_{\mbG} |_j$  coincides with $\nabla'_j$ for any $j$.
By construction, the pair $\msE^\spadesuit_{\msG \Rightarrow \emptyset} := (\DE_{\mbB}, \nabla)$ forms a  $\ \qq$-normal  $(\mfg, \hslash)$-oper 
with
$(\msE^\spadesuit_{\msG \Rightarrow \emptyset})_{\emptyset \Rightarrow \msG} =  \msE^\spadesuit_{\msG}$.
The resulting assignment $\msE^\spadesuit_{\msG} \mapsto \msE^\spadesuit_{\msG \Rightarrow \emptyset}$ defines  the inverse to $\msE^\spadesuit \mapsto \msE^\spadesuit_{\emptyset \Rightarrow \msG}$.
This completes the proof of the second assertion.

Finally, the last assertion follows from the second assertion together with the result  in \S\,\ref{QG10}.
\end{proof}

\subsection{Clutching morphisms of (dormant)  $(\mfg, \hslash)$-opers}\label{y0174sf}

Let $\rho_\msG := \{ \rho^j \}_{j=1}^n$ be a set of radii for $\msG$ over $k$ and 
$\msE^\spadesuit_\msG$ be a $\ \qq$-normal  $\msG$-$(\mfg, \hslash)$-oper on $\msX$.
According to the above proposition,
there exists a unique $\ \qq$-normal  $(\mfg, \hslash)$-oper $\msE^\spadesuit$  with
$\msE^\spadesuit_{\emptyset \Rightarrow \msG} = \msE^\spadesuit_\msG$.
 The radii  of $\msE^\spadesuit$ are constructed from $\rho^j$ ($j = 1, \cdots, n$) by means of the clutching data $\msG$ in an evident way; we denote the resulting radii of $\msE^\spadesuit$ by 
\index{$\rho_\msG$, set of  radii for $\msG$@--- $\rho_{\msG \Rightarrow \emptyset}$}
\begin{equation} \label{fff036}
\rho_{\msG \Rightarrow \emptyset} \in \mfc^{\times r} (S).
 \end{equation}
Since  any point of $\prod_{j =1}^n \mfO \mfp_{\mfg, \hslash, \rho^j, g_j,r_j}$ classifies  a  unique  $\ \qq$-normal 
$\msG$-$(\mfg, \hslash)$-oper of radii $\rho_\msG$   (cf. Proposition \ref{y049}),
the assignment $(\{\msX_j\}_{j=1}^n, \msE^\spadesuit_{\msG}) \mapsto (\msX, \msE^\spadesuit)$  determines
a well-defined morphism
\index{$\mr{Clut}_\msG$, clutching morphism@--- $\mr{Clut}_{\msG, \rho_\msG}$, ${^p}\mr{Clut}_{\msG, \rho_\msG}$}
\begin{align}
\mr{Clut}_{\msG, \rho_\msG} : \prod_{j =1}^n \mfO \mfp_{\mfg, \hslash, \rho^j, g_j,r_j} \migi  
\mfO \mfp_{\mfg, \hslash, \rho_{\msG \Rightarrow \emptyset}, g,r}
\end{align}
fitting into the following square diagram of moduli stacks:
\index{$\mr{Clut}_\msG$, clutching morphism@--- $\mr{Clut}_{\msG, \rho_\msG}$, ${^p}\mr{Clut}_{\msG, \rho_\msG}$}
\begin{align} \label{QG26}
\vcenter{\xymatrix@C=46pt@R=36pt{
\prod_{j =1}^n \mfO \mfp_{\mfg, \hslash, \rho^j, g_j,r_j}  
\ar[r]^-{\mr{Clut}_{\msG, \rho_\msG}} \ar[d]_-{\mr{projection}}& \mfO \mfp_{\mfg, \hslash, \rho_{\msG \Rightarrow \emptyset}, g,r}  \ar[d]^-{\mr{projection}}\\
 \prod_{j=1}^n \overline{\mfM}_{g_j, r_j} \ar[r]_-{\mr{Clut}_\msG} & \overline{\mfM}_{g,r},
}} 
\end{align}
where the products ``$\prod$'' are taken over $k$.
Moreover,  if $k$ is of characteristic $p>0$,   this diagram restricts to  a commutative square diagram
\begin{align} \label{1032}
\vcenter{\xymatrix@C=46pt@R=36pt{
\prod_{j =1}^n \mfO \mfp^\ZZZ_{\mfg, \hslash, \rho^j, g_j,r_j}
\ar[r]^-{{^p}\mr{Clut}_{\msG, \rho_\msG}} \ar[d]_-{\mr{projection}}& \mfO \mfp^\ZZZ_{\mfg, \hslash, \rho_{\msG \Rightarrow \emptyset}, g,r}
 \ar[d]^-{\mr{projection}}\\
 \prod_{j=1}^n \overline{\mfM}_{g_j, r_j}  \ar[r]_-{\mr{Clut}_\msG}& \overline{\mfM}_{g,r}.
}}
\end{align}
Here, recall (cf.   Theorem \ref{y0102}, (i)) that 
if $\hslash \in k^\times$, then $\mfO \mfp^\ZZZ_{\mfg, \hslash, \rho_{\msG \Rightarrow \emptyset}, g,r}$ is empty unless the radii $\rho_{\msG \Rightarrow \emptyset}$ is of the form $\rho_{\msG \Rightarrow \emptyset} = \hslash \star \rho$ for some $\rho \in \mfc^{\times r}(\mbF_p)$.

\bt \label{y0176}
Suppose that $\mr{char}(k) =:p > 2  h_\mbG$.
Let  $\msG :=  ( \GR, \{ (g_j, r_j) \}_{j=1}^n, \{ \lambda_j \}_{j=1}^n)$ be  clutching data  of type $(g,r)$,  $\hslash$ an element of $k^\times$,  and $\rho$ an element of $\mfc^{\times r}(\mbF_p)$.
Then, the following commutative square diagram  is cartesian:
\begin{align} \label{1050}
\vcenter{\xymatrix@C=76pt@R=36pt{
\coprod_{ \rho_\msG := \{ \rho^j\}_{j=1}^n}\prod_{j =1}^n \mfO \mfp^\ZZZ_{\mfg, \hslash, \hslash \star \rho^j, g_j,r_j} 
\ar[r]^-{\coprod_{\rho_\msG} {^p}\mr{Clut}_{\msG, \rho_\msG}} \ar[d]_-{\mr{projection}}& \mfO \mfp^\ZZZ_{\mfg, \hslash, \hslash \star \rho, g,r} 
\ar[d]^-{\mr{projection}}\\
 \prod_{j=1}^n \overline{\mfM}_{g_j, r_j}  \ar[r]_-{\mr{Clut}_\msG}& \overline{\mfM}_{g,r},
}}
\end{align}
 where the disjoint union on the upper-left corner runs  over  the  sets of radii $\rho_\msG$ for $\msG$ over $\mbF_p$ with $\rho_{\msG \Rightarrow \emptyset} = \rho$. 
 \et
\begin{proof}
The assertion follows from Propositions 
\ref{y0102}, (i), 
 and  \ref{y0166pp}.
\end{proof}

\begin{exa}[Classical clutching morphisms]
Let us describe the cartesian diagram (\ref{1050}) in two particular cases corresponding to clutching morphisms, in the classical sense,  between moduli spaces of pointed stable curves. 

First, given  two pairs of nonnegative integers $(g_j, r_j)$ ($j=1,2$) with $2g_j-1+r_j >0$, we shall consider 
clutching data 
\begin{align} \label{QH1009}
\msG_\mr{tree} := (\GR_\mr{tree}, \{(g_j, r_j+1)\}_{j=1}^2, \{ \lambda_j \}_{j=1}^2)
\end{align}
  determined by the following conditions:
\begin{itemize}
\item
$\msG_\mr{tree}$ is of type $(g,r)$, where  $r= r_1 +r_2$
 and   $g = g_1 + g_2$;
\item
$\GR_\mr{tree} := (\{ v_1, v_2 \},  \{e_{1i}\}_{i=1}^{r_1} \sqcup \{e_{2i}\}_{i=1}^{r_2} \sqcup \{e \}, \zeta)$, where $\zeta (e_{ji}) = \{ \circledast, v_j \}$ ($j=1,2$, $1 \leq i \leq r_j$)
 and  $\zeta (e) = \{v_1, v_2 \}$;
\item
For each $j=1,2$, the bijection $\lambda_j : B_{v_j} \isom \{ 1, \cdots, r_j+1 \}$ is defined in such a way that
$\lambda_j (b) = i$ if $b \in e_{ji}$ ($1 \leq i \leq r_i$) and $\lambda (b) = r_i+1$ if $b \in e$.
\end{itemize}
In particular,  we can visualize the semi-graph $\GR_\mr{tree}$  as Figure 1 in Remark \ref{QH1001}.
The clutching morphism  associated to $\msG_{\mr{tree}}$ is the morphism
\begin{align} \label{QH1008}
\mr{Clut}_{\msG_\mr{tree}} : \overline{\mfM}_{g_1, r_1+1} \times_k \overline{\mfM}_{g_2, r_2+1} \migi \overline{\mfM}_{g, r} 
\end{align}
obtained by attaching the respective last marked points of curves classified by
$\overline{\mfM}_{g_1, r_1+1}$ and $\overline{\mfM}_{g_2, r_2+1}$  to form a nodal point.
Now, we shall take an element 
\begin{align} \label{QH1033}
\rho :=
 (\rho^1, \rho^2)
 \in \mfc^{\times r_1}(\mbF_p) \times \mfc^{\times r_2} (\mbF_p) \left(= \mfc^{\times r}(\mbF_p) \right). 
\end{align}
 Note that each set   of radii $\rho_{\msG_\mr{tree}}$ for $\msG_{\mr{tree}}$ over $\mbF_p$ with $\rho_{\msG_\mr{tree} \Rightarrow \emptyset} = \rho$ can be described as 
$\rho_{\msG_\mr{tree}} := \{ \rho^j_\lambda \}_{j=1}^2$ for a unique $\lambda \in \mfc (\mbF_p)$, where 
$\rho_\lambda^1 := (\rho^1, \lambda) \in \mfc^{\times (r_1+1)}(\mbF_p)$ and 
$\rho_\lambda^2 := (\rho^2, \lambda^\dual) \in \mfc^{\times (r_2+1)} (\mbF_p)$.
Hence,  the cartesian square diagram (\ref{1050}) for this pair  $(\msG_\mr{tree}, \rho)$
becomes the following diagram: 
\begin{align} \label{QH1004}
\vcenter{\xymatrix@C=74pt@R=36pt{
\coprod_{ \lambda \in \mfc (\mbF_p)} \mfO \mfp^\ZZZ_{\mfg, \hslash, \hslash \star \rho^1_\lambda, g_1,r_1+1}  \times_k \mfO \mfp^\ZZZ_{\mfg, \hslash, \hslash \star \rho^2_\lambda, g_2,r_2+1}  
\ar[r]^-{\coprod_{\lambda} {^p}\mr{Clut}_{\msG_\mr{tree}, \{\rho_\lambda^j\}_j}} \ar[d]_-{\mr{proj.}}& \mfO \mfp^\ZZZ_{\mfg, \hslash, \hslash \star \rho, g,r} 
\ar[d]^-{\mr{proj.}}\\
 \overline{\mfM}_{g_1, r_1+1} \times_k \overline{\mfM}_{g_2, r_2+1}  \ar[r]_-{\mr{Clut}_{\msG_\mr{tree}}}& \overline{\mfM}_{g,r}.
}}
\hspace{-3mm}
\end{align}

Next, given a pair of nonnegative integers $(g,r)$ with $2g-2+r >0$ and $g>0$, we shall consider 
clutching data 
\begin{align} \label{QH1023}
\msG_\mr{loop} := (\GR_\mr{loop}, \{(g-1, r+2)\}, \{ \lambda\})
\end{align}
of type  $(g, r)$  determined by the following conditions:
\begin{itemize}
\item
$\GR_\mr{loop} := (\{ v\},  \{e_{i}\}_{i=1}^{r} \sqcup \{e\}, \zeta)$, where $\zeta (e_{i}) = \{ \circledast, v \}$ ($1 \leq i \leq r$), $e := \{ b_{r+1}, b_{r+2} \}$,  and  $\zeta (b_{r+1}) = \zeta (b_{r+2}) = v$;
\item
 $\lambda : B_{v} \isom \{ 1, \cdots, r+2 \}$ is the bijection defined in such a way that
$\lambda (b) = i$ if $b \in e_{i}$ ($1 \leq i \leq r$) and $\lambda (b_{r+1}) = r+1$, $\lambda (b_{r+2}) =r +2$. 
\end{itemize}
We can visualize the semi-graph $\msG_{\mr{loop}}$ as Figure 2 in Remark \ref{QH1001}.
The clutching morphism associated to $\msG_{\mr{loop}}$ is the morphism
\begin{align} \label{QH1020}
\mr{Clut}_{\msG_{\mr{loop}}} : \overline{\mfM}_{g-1,r+2} \migi \overline{\mfM}_{g,r}
\end{align}
obtained by attaching the last two marked points of each curve classified by $\overline{\mfM}_{g-1,r+2}$ to form a nodal point.
Let us take an element $\rho \in \mfc^{\times r}(\mbF_p)$.
Then, each set of radii $\rho_{\msG_\mr{loop}}$ for $\msG_{\mr{loop}}$ over $\mbF_p$ with $\breve{\msG}_{\mr{loop}} = \rho$ can be described as $\rho_{\msG_\mr{loop}}:= \{(\rho, \lambda, \lambda^\dual)\}$ for a unique $\lambda \in \mfc (\mbF_p)$.
It follows that   (\ref{1050}) for this pair $(\msG_{\mr{loop}}, \rho)$ becomes the following diagram:
\begin{align} \label{QH1021}
\vcenter{\xymatrix@C=76pt@R=36pt{
\coprod_{ \lambda \in \mfc (\mbF_p)} \mfO \mfp^\ZZZ_{\mfg, \hslash, \hslash \star (\rho, \lambda, \lambda^\dual), g-1,r+2}  
\ar[r]^-{\coprod_{\lambda} {^p}\mr{Clut}_{\msG_\mr{loop}, (\rho, \lambda, \lambda^\dual))}} \ar[d]_-{\mr{proj.}}& \mfO \mfp^\ZZZ_{\mfg, \hslash, \hslash \star \rho, g,r} 
\ar[d]^-{\mr{proj.}}\\
 \overline{\mfM}_{g-1, r+2}   \ar[r]_-{\mr{Clut}_{\msG_\mr{loop}}}& \overline{\mfM}_{g,r}.
}}
\end{align}

\end{exa}

\subsection{Totally degenerate curves}\label{y0177}

We shall recall a certain kind of pointed stable  curves,  referred  to as  {\it totally degenerate} curves.

In this subsection, assume   that the base  field $k$   is algebraically closed.
Also, let $(g,r)$ be a pair of nonnegative integers with $2g-2+r >0$,
$\msX: = (X/k, \{ \sigma_i \}_{i =1}^r)$   an $r$-pointed stable curve over $k$ of genus $g$.
 Write
 $\nu_{\msX} : \coprod_{l =1}^{m} X_l \migi X$ for  the normalization of $X$, where  $m$ denotes some positive integer and  each $X_l$ ($l = 1, \cdots, m$) denotes a  {\it smooth}  curve over $k$.

\bde \label{y0178} 
We shall say that  $\msX$ is {\bf totally degenerate} 
\index{$\msX$, pointed (stable) curve@--- totally degenerate}
 if, for any $l =1, \cdots, m$, the pointed smooth curve
\begin{equation}
{\msX_l}:= (X_l/k, \nu_{\msX}^{-1}(E_{\msX}) \cap X_l (k))
\end{equation}
(defined after ordering elements in $\nu_{\msX}^{-1}(E_{\msX}) \cap X_l (k)$)  is isomorphic to $\msP$ (cf. (\ref{1051})), where   $E_{\msX}$ denotes the set of edges in the dual semi-graph $\GR_\msX$ associated to $\msX$ (cf. \S\,\ref{R001}), regarded  as a subset of $X (k)$.
 \ede

\begin{rema} \label{RRR001} 
One may associate an $r$-pointed totally degenerate curve of genus $g$ to each clutching data $\msG = (\GR, \{ (g_j, r_j)\}_{j=1}^n, \{ \lambda_j \}_{j=1}^n )$  of type $(g,r)$ with   $(g_j, r_j) = (0,3)$ for every $j = 1, \cdots, n$.
Indeed, if $\msG$ is  such a collection, then
the pointed stable curve classified by 
the clutching morphism 
\begin{equation}
\mr{Clut}_\msG : \left(\mr{Spec}(k) \cong \right)    \prod_{j=1}^n \overline{\mfM}_{0,3} \migi \overline{\mfM}_{g,r}\end{equation}
 turns out to be totally degenerate.
 Conversely, each totally degenerate curve may be obtained in this way for some clutching data.
The locus in $\overline{\mfM}_{g,r}$ classifying totally degenerate curves is of zero-dimensional and  forms a direct sum of finite copies of $\mr{Spec}(k)$.
 \end{rema}

\subsection{The generic  degree of $\mfO\mfp^\ZZZ_{\mfg, \hslash, \hslash \star \rho, g, r}/\overline{\mfM}_{g,r}$} \label{QG1100}

Now, we shall consider   the following two conditions on $\mbG$:

\begin{itemize}
\item[$(*)$  ]
For any $\rho \in \mfc^{\times 3}(\mbF_p)$, the finite $k$-scheme  $\mfO\mfp^{^\mr{Zzz...}}_{\mfg, \rho, 0, 3}$ is, if it is nonempty,  unramified (or equivalently,  \'{e}tale)  over $k$.
\item[$(**)$]
For any $\rho  \in \mfc^{\times 3}(\mbF_p)$, 
the $k$-scheme $\mfO\mfp^{^\mr{Zzz...}}_{\mfg, \rho, 0, 3}$ is isomorphic to 
$\mfO\mfp^{^\mr{Zzz...}}_{\mfg, (-1)\star\rho, 0, 3}$.
\end{itemize}

\begin{rema} \label{QH100}
Let $\hslash$ be an element of $k^\times$.
Then, because of the resp'd portion of (\ref{QS1100}),
  the condition  $(*)$ (resp., $(**)$) is  equivalent to the condition that $\mfO\mfp^{^\mr{Zzz...}}_{\mfg, \hslash, \hslash \star \rho, 0, 3}$ is unramified over $k$ (resp., $\mfO\mfp^{^\mr{Zzz...}}_{\mfg, \hslash, \hslash \star \rho, 0, 3}$ is isomorphic to 
$\mfO\mfp^{^\mr{Zzz...}}_{\mfg,  \hslash, (-\hslash)\star\rho, 0, 3}$)  for every $\rho \in \mfc^{\times 3}(\mbF_p)$.
\end{rema}
\begin{rema}
By Theorem \ref{ppp012},
the condition $(**)$ is satisfied if $\mbG = \mr{PGL}_n$ (with $2 n < p$), $\mr{GO}_{2l+1}$ (with $4  l+2<p$), or $\mr{GSp}_{2m}$ (with $4 m <p$).
Also, as will be shown later (cf.  Theorem  \ref{3778567d}), the condition $(*)$ is satisfied for  the same   types  of $\mbG$.
The case of $\mbG = \mr{PGL}_2$ can be verified  immediately.
In fact, the condition $(**)$ follows from the equality  $(-1) \star \rho = \rho$ for $\rho \in \mfc_{\mfs \mfl_2}(\mbF_p)$, and the condition $(*)$ follows from Theorem \ref{QR90}, which say that $\mfO\mfp^{^\mr{Zzz...}}_{\mfg, \rho, 0, 3}$ is, if it is nonempty, isomorphic to $\mr{Spec}(k)$.
\end{rema}

\bpr \label{c43} 
Let   $(g,r)$ a pair of nonnegative integers with $2g-2+r>0$,
$\hslash$ an element of $k^\times$,
and  $\rho$ an element of $ \mfc^{\times r} (\mbF_p)$ (where  $\rho = \hslash \star \rho := \emptyset$ if $r =0$).
Suppose   that $\mbG$ satisfies  the condition  $(*)$.
Then,
 $\mfO\mfp^\ZZZ_{\mfg, \hslash, \hslash \star \rho, g, r}$ 
 is (finite,  by virtue of Theorem \ref{y0102}, and) \'{e}tale over the points of  $\overline{\mfM}_{g,r }$ classifying totally degenerate curves.
In particular, 
 $\mfO\mfp^\ZZZ_{\mfg, \hslash, \hslash \star \rho, g, r}$ is generically \'{e}tale over $\overline{\mfM}_{g,r }$, i.e., any irreducible component that dominates  $\overline{\mfM}_{g,r}$ admits a dense open subscheme which is \'{e}tale over $\overline{\mfM}_{g,r}$.

\epr
\begin{proof}
Let us take an arbitrary $r$-pointed  totally degenerate curve $\msX$
of genus $g$ over $k$.
We can find clutching data
  $\msG = (\GR, \{ (g_j, r_j)\}_{j=1}^n, \{ \lambda_j\}_{j=1}^n )$  defining   $\msX$   in the manner  of Remark \ref{RRR001} (hence $(g_j, r_j) = (0,3)$ for every  $j$).
By   Theorem \ref{y0176} for  this clutching data,
the condition $(*)$ (cf. Remark \ref{QH100})
implies that
the left-hand vertical arrow  in   (\ref{1050}) is unramified.
Hence,  the projection  $\mfO\mfp^\ZZZ_{\mfg, \hslash, \hslash \star \rho, g, r} \migi \overline{\mfM}_{g,r }$, i.e., the right-hand vertical arrow in  (\ref{1050}),  is unramified over the point of $\overline{\mfM}_{g,r }$ classifying $\msX$.
But, it follows from Corollary \ref{y0170} that $\mfO\mfp^\ZZZ_{\mfg, \hslash, \hslash \star \rho, g, r} \migi \overline{\mfM}_{g,r }$  is also \'{e}tale over this point.
This completes the proof of the  former assertion.

The latter assertion, i.e., the generic \'{e}taleness of $\mfO\mfp^\ZZZ_{\mfg, \hslash, \hslash \star \rho, g, r}/\overline{\mfM}_{g,r }$,  follows from the former assertion.
Indeed, since $\mfO\mfp^\ZZZ_{\mfg, \hslash, \hslash \star \rho, g, r}/\overline{\mfM}_{g,r }$ is finite (cf. Proposition \ref{y0102}, (i)) and $\overline{\mfM}_{g,r }$ is irreducible,
any irreducible component  $\mfN$ of $\mfO\mfp^\ZZZ_{\mfg, \hslash, \hslash \star \rho, g, r}$ that dominates  $\overline{\mfM}_{g,r }$ surjects onto $\overline{\mfM}_{g,r }$.
In particular, the fiber of the projection $\mfN  \migi \overline{\mfM}_{g,r }$
over the point classifying a totally degenerate curve is nonempty.
By the former assertion and the open nature of flatness, 
the \'{e}tale locus in the component $\mfN$
  (relative to  $\overline{\mfM}_{g,r }$) forms a {\it nonempty} (hence, dense) open subscheme.
\end{proof}

 For $\rho := (\rho_i)_{i=1}^r \in \mfc^{\times r}(\mbF_p)$, 
 we shall write 
 \begin{equation}
 \label{1052}  \mr{deg}(\mfO\mfp^\ZZZ_{\mfg, \hslash, \hslash \star \rho, g, r}/\overline{\mfM}_{g,r })  \ \left(\geq 0 \right)
 \index{$\mr{deg}(\mfO\mfp^\ZZZ_{\mfg, \hslash, \hslash \star \rho, g, r}/\overline{\mfM}_{g,r })$}
 \end{equation}
for  the {\it generic degree}   of $\mfO\mfp^\ZZZ_{\mfg, \hslash, \hslash \star \rho, g, r}$ over $\overline{\mfM}_{g,r }$, i.e., the degree of the fiber of the projection $\mfO\mfp^\ZZZ_{\mfg, \hslash, \hslash \star \rho, g, r} \migi \overline{\mfM}_{g,r }$ over the generic point of $\overline{\mfM}_{g,r }$. 
 This value depends neither on 
   the choice of $\hslash \in k^\times$ (cf. (\ref{QS1100}))
nor on the ordering of $\hslash \star \rho_1, \cdots, \hslash \star \rho_r$.
Here, if  $\varrho := (\varrho_i)_{i=1}^r$ is an element of $\mfc^{\times r} (\mbF_p)$, then we shall write $\overline{\varrho}$ for the multiset $[\varrho_1, \cdots, \varrho_r]$ (with cardinality $r$) over the set $\mfc (\mbF_p)$.
(For the definition of a {\it multiset}, we refer the reader  to ~\cite{SIYS}.)
By the above observation, it makes sense to abbreviate 
the value $ \mr{deg}(\mfO\mfp^\ZZZ_{\mfg, \hslash, \hslash \star \rho, g, r}/\overline{\mfM}_{g,r })$  for an arbitrary $\hslash \in k^\times$
 to 
\index{${^p}N_{\mfg,  0}$, ${^p}N_{\mfg,  g}$}
\begin{equation} \label{d6}
{^p}N_{\mfg, g} (\overline{\rho}).
\end{equation}


\section{Pseudo-fusion rules and pseudo-fusion rings}\label{y0179}

\subsection{Pseudo-fusion rules} \label{QG1101}

In what follows, we shall develop a general theory of {\it pseudo-fusion rules}.
The notion of a pseudo-fusion rule  may be thought of as a somewhat weaker variant of the notion of a {\it fusion rule}.
First, we shall recall  from ~\cite{Beau}, \S\,5 (or  ~\cite{Fuchs}), the definition of a fusion rule as follows.
Let $I$ be a nonempty finite set equipped with an involution $\lambda \mapsto \lambda^\dual$ ($\lambda \in I$), i.e., $(-)^{\dual \dual}=\mr{id}_I$.
Denote by $\mbN^{I}$ the free commutative monoid  with basis $I$,
i.e., the set of sums $\sum_{\alpha \in I}n_\alpha \cdot \alpha$ with $n_\alpha \in \mbN$;
we shall always identify $I$ with a subset of $\mbN^{I}$ in the natural fashion.
The involution of $I$ extends by linearity to an involution $x \mapsto x^\dual$ of $\mbN^{I}$.
Recall that  a {\bf fusion rule}\index{fusion rule} on $I$ (cf.  ~\cite{Beau}, \S\,5, Definition) is, by definition,  a map of sets  $N : \mbN^{I} \migi \mbZ$  satisfying the following three conditions:
\vspace{1mm}
\begin{itemize}
\item[(i)]
$N (0) = 1$, and $N(\alpha) >0$ for some $\alpha \in I$;
\vspace{1mm}
\item[(ii)]
$N(x^\dual) = N (x)$ for any $x \in \mbN^{I}$; 
\vspace{1mm}
\item[(iii)]
$N (x +y) = \sum_{\lambda \in I}N (x +\lambda) \cdot N (y + \lambda^\dual)$ for any $x$, $y \in \mbN^{I}$.
\end{itemize}
\vspace{3mm}

\subsection{Pseudo-fusion rings} \label{QG1103}

As described below,  the notion of a pseudo-fusion rule is defined as a map  satisfying conditions  somewhat weaker than the above  conditions (i)-(iii).
To describe the definition, let us consider the augmentation map
\begin{align}
\mr{aug}_I : \mbN^{I} \migi  \mbN;  \sum_{\alpha \in I}n_\alpha \cdot \alpha \mapsto  \sum_{\alpha \in I} n_\alpha
\end{align}
which is a monoid homomorphism.
In particular, we have  $\mr{aug}_I^{-1} (1) = I$.
 For each  integer $l$, we shall write
\begin{align}
\mbN^{I}_{\geq l} := \left\{ x \in \mbN^{I} \, | \, \mr{aug}_I (x) \geq l \right\}
\end{align}
(hence $\mbN^{I}_{\geq l} = \mbN^{I}$ if $l \leq 0$).
$\mbN^{I}_{\geq l}$  may be naturally  identified with the set of multisets with cardinality $\geq l$ over $I$.

\bde \label{c34} 
 A {\bf pseudo-fusion rule} \index{pseudo-fusion rule} on $I$ is a map $N : \mbN^{I}_{\geq 3} \migi \mbZ$ satisfying the following two conditions:
\begin{itemize}
\item[(ii)']
$N(x^\dual) = N(x)$ for any $x \in \mbN^{I}_{\geq 3}$; 
\vspace{1mm}
\item[(iii)']
$N (x +y) = \sum_{\lambda \in I}N (x +\lambda) \cdot N (y + \lambda^\dual)$ for any $x$, $y \in \mbN^{I}_{\geq 2}$ .
\end{itemize}
 \ede

\begin{rema} \label{QH123}
Because of (iii)', the condition (ii)' can be replaced with the following weaker  condition: 
\begin{itemize}
\item[(ii)'']
$N (\alpha^\dual + \beta^\dual + \gamma^\dual) = N (\alpha + \beta + \gamma)$ for any $\alpha, \beta, \gamma \in I$.
\end{itemize}
\end{rema}

Let us  fix a pseudo-fusion rule $N$ on $I$.
By means of this map,
one may  define a multiplication law 
$\ast' : \mbZ^{I} \times \mbZ^{I} \migi \mbZ^{I}$ 
on the free abelian group $\mbZ^{I}$ with basis $I$ by putting 
\begin{equation} \label{M12}
  \alpha \ast' \beta := \sum_{\lambda \in I} N (\alpha + \beta + \lambda^\dual)  \cdot \lambda \end{equation}
for any $\alpha$, $\beta \in I$ and extending by bilinearity.
This multiplication law is verified to be  commutative and  associative.
Moreover, by induction on $l \geq 2$, one verifies from  the condition (iii)' that 
\begin{equation} \label{c56}
\alpha_1 \ast' \alpha_2 \ast' \cdots \ast' \alpha_l= \sum_{\lambda \in I} N ((\sum_{i =1}^l \alpha_i) + \lambda^\dual) \cdot \lambda
\end{equation}
for any $\alpha_1, \cdots, \alpha_l \in I$.

Denote by $\Fus_N'$ the  ring $\mbZ^{I}$ endowed with the multiplication ``\,$\ast'$\," defined above.
Note that  $\Fus'_N$ \index{$\Fus'_N$} is  not necessarily unital.
Denote by 
\begin{equation} \label{c58}
\Fus_N
\index{$\Fus_N$, pseudo-fusion ring associated to $N$}
\end{equation}
 the {\it unitization} 
   of $\Fus'_N$.
Namely, $\Fus_N$ is 
the unital ring defined as follows:
\begin{itemize}
\item
the  underlying abelian group is the direct sum $\mbZ \cdot \varepsilon \oplus \mbZ^{I}$ (each of whose element will be  denoted by $a \cdot \varepsilon + x$ for some $a \in \mbZ$ and $x \in \mbZ^{I}$)
  of $\mbZ^{I}$ and the  free abelian group $\mbZ \cdot  \varepsilon$ generated by a  symbol $\epsilon$;
\item
the multiplication law $\ast : \Fus_N \times \Fus_N \migi \Fus_N$
on $\Fus_N$  is defined by 
\begin{equation} \label{c75}
(a \cdot \varepsilon +x) \ast (b \cdot \varepsilon + y) := a b \cdot  \varepsilon + a  \cdot y + b \cdot   x + x \ast' y
\end{equation}
for any $a$, $b \in \mbZ$ and $x$, $y \in \mbZ^{I}$.
\end{itemize}
\vspace{2mm}
Thus, $\Fus_N$ is, by construction, a commutative, associative, and moreover,  {\it unital} ring  with identity $\varepsilon$.
The involution of $\mbN^{I}$ extends, by assigning $\varepsilon$ identically,  to an involution $z \mapsto z^\dual$
of $\Fus_N$ that 
  is compatible with the multiplication of $\Fus_N$.
Moreover, for any field $K$ over $\mbQ$,   the $K$-algebra $K \otimes_\mbZ \Fus_N$ is equipped with a $K$-linear  involution induced naturally by the involution of $\Fus_N$. 
The nonunital $\Fus'_N$ may be considered as the ideal defined to be the kernel of the ring homomorphism 
\begin{align} \label{f123}
\mr{aug}_{\scalebox{0.7}{$\Fus_N$}} :   \Fus_N  \migi   \mbZ;  a \cdot \varepsilon + x  \mapsto  a.
\index{$\mr{aug}_{\scalebox{0.7}{$\Fus_N$}}$}
\end{align}

\bde \label{c4fert}
We shall refer to the ring $\Fus_N$ as the {\bf pseudo-fusion ring} associated to the pseudo-fusion rule $N$.
\ede

\begin{rema}[Fusion rule vs. pseudo-fusion rule] \label{f91}  
Note that each fusion rule  on $I$ 
   determines  a pseudo-fusion rule on $I$.
More precisely,  if $\breve{N}$ is a fusion rule on $I$, then 
the restriction $\breve{N}_{\geq 3}$ of $\breve{N}$ to $\mbN_{\geq 3}^{I}$ ($\subseteq \mbN^I$)
specifies  evidently  a pseudo-fusion rule on $I$. 
Let $\Fus_{\breve{N}_{\geq 3}}$ be the pseudo-fusion ring
   associated to the pseudo-fusion rule $\breve{N}_{\geq 3}$.
Also,  
denote by $\Fus^\circledcirc_{\breve{N}}$ the fusion ring  (in the sense of ~\cite{Beau}, \S\,5, Definition) associated to the fusion rule $\breve{N}$.
The underlying abelian group of $\Fus^\circledcirc_{\breve{N}}$ may be identified with that of $\Fus'_{\breve{N}_{\geq 3}} \left(\subseteq \Fus_{\breve{N}_{\geq 3}}\right)$.
According to the discussion in ~\cite{Beau}, \S\,5, the ring $\Fus^\circledcirc_{\breve{N}}$ is  commutative, associative, and moreover, unital; we shall denote  its identity element by   $\varepsilon^\circledcirc$.
The map
\begin{align} \label{QH102}
\Fus_{\breve{N}_{\geq 3}} \migi \mbZ \times \Fus_{\breve{N}}^\circledcirc; 
a \cdot \varepsilon + b  \mapsto (a, a \cdot  \varepsilon^\circledcirc + b)
\end{align}
 defines
 an isomorphism of rings. 
Consequently, {\it the pseudo-fusion ring $\Fus_{\breve{N}_{\geq 3}}$ is isomorphic to the direct product of the fusion ring $\Fus_{\breve{N}}^\circledcirc$ and the ring of integers $\mbZ$}.
\end{rema}

\section{The ring structure of the fusion ring}\label{y0180}

In this section, we examine  the ring structure of the pseudo-fusion ring $\Fus_N$.
Write  
$\sqrt{(0)}$ for the nilradical of $\Fus_N$ and 
\begin{equation}
\Fus_{N \mr{red}}
\index{$\Fus_{N \mr{red}}$, reduced ring associated to $\Fus_N$}
\end{equation}
 for  the reduced ring associated to $\Fus_{N}$, i.e., $\Fus_{N \mr{red}} := \Fus_{N}/\sqrt{(0)}$.
  The underlying abelian group of $\Fus_{N \mr{red}}$  is   torsion-free,  so $\Fus_{N \mr{red}}$ may be 
 considered   as a subring  of $ \Fus_{N \mr{red}} \otimes_\mbZ \mbQ$.
 The ring homomorphism $\mr{aug}_{\Fus_N}$  factors through the natural surjection $\Fus_N \migisurj \Fus_{N  \mr{red}}$; we shall write
\begin{equation} \label{f124}
\mr{aug}_{\scalebox{0.7}{$\Fus_{N \mr{red}}$}} : \Fus_{N \mr{red}} \migi \mbZ
\end{equation} 
for the resulting homomorphism.

\begin{rema}[The kernel of a  pseudo-fusion rule] \label{f141} 
 We shall write
\begin{equation} \label{f393}
 \mr{Ker} (N) \ \left(\subseteq I \right)\index{$\mr{Ker} (N)$, kernel of $N$}
 \end{equation}
  for 
 the set  of elements $\alpha$ in $I$ such that $N (\alpha +x) = 0$ for all $x \in \mbN^{I}_{\geq 2}$;
  we shall refer to $\mr{Ker} (N)$ as the {\bf kernel} of  $N$.
The subset $\mr{Ker}(N)$ (as well as its complement $I \setminus \mr{Ker}(K)$) of $I$ is closed under the involution.

Now, let $I'$ be a subset of $I$ that is closed under the involution
 and satisfies the inclusion relation  $I \setminus  I' \subseteq \mr{Ker}(N)$.
The involution $(-)^\dual$ restricts to  an involution on $I'$.
Then, the  restriction $N'$ of $N$ to $\mbN^{I'} \left(\subseteq \mbN^I\right)$ forms a pseudo-fusion rule on $I'$, and hence, gives rise to  the pseudo-fusion ring $\Fus_{N'}$ associated to it.
Since $N$ may be  obtained  as the extension  of  $N'$
 by zero  to $\mbN^{I}$, 
 we see that $\Fus_N$ is a retract of $\Fus_{N'}$, in   particular, the projection  $\pi : \Fus_{N} \migiincl \Fus_{N'}$ given by $a \cdot \varepsilon + \sum_{\alpha \in I} n_{\alpha} \cdot \alpha \mapsto a \cdot \varepsilon + \sum_{\alpha \in I'} n_{\alpha} \cdot \alpha$ defines a ring homomorphism.
 By the assumption ``$I \setminus I' \subseteq \mr{Ker}(N)$'', 
 this projection induces 
 an isomorphism   of rings
\begin{equation} \label{w33344}
\Fus_{N \mr{red}} \isom \Fus_{N' \mr{red}}.
\end{equation}

We shall  write $\overline{\mr{Cas}}_{N}$ (resp.,  $\overline{\mr{Cas}}_{N'}$)  for the  element of $\Fus_{N \mr{red}}$ (resp., $\Fus_{N' \mr{red}}$) defined as   the image of $\mr{Cas}_N := \sum_{\lambda \in I} \lambda \ast \lambda^\dual \in \Fus_N$ (resp., $\mr{Cas}_{N'} := \sum_{\lambda \in I'} \lambda \ast \lambda^\dual \in \Fus_{N'}$) (cf. (\ref{w33333}))  via the surjection $\Fus_N \migisurj \Fus_{N \mr{red}}$ (resp., $\Fus_{N'} \migisurj \Fus_{N' \mr{red}}$).
Since $\lambda \ast \lambda^\dual =0$ for any $\lambda \in I \setminus I'$,   (\ref{w33344}) sends $\overline{\mr{Cas}}_{N'}$ to $\overline{\mr{Cas}}_{N}$.
Hence,  in many cases, $I$ can be replaced with $I'$ in performing some calculations that we want to derive from the structure of the ring $\Fus_N$ (cf. Propositions \ref{ccc35} and  \ref{c35}, (ii)).
 \end{rema}


Let $K$ be a field over $\mbQ$, $V$ a finite-dimensional  $K$-vector space, and $h$ a $K$-linear endomorphism of $V$.
Then,
we shall write
\begin{equation}\label{QH1300}
\mr{Tr}_K(h \,  | \, H)   \left(\in K \right)
\index{$\mr{Tr}_K(h  \mid  H)$}
\end{equation}
 for the trace of $h$.
If, moreover,  $V$ admits a structure of  $K$-algebra and  $x$ is an element of $V$,
then the multiplication by $x$ defines a $K$-linear  endomorphism  of the underlying $K$-vector space $V$, which we denote by 
\begin{equation} \label{QH1302}
\mr{mult} (x) \left(\in \mr{End}_K (V)\right).
\index{$\mr{mult} (x)$}
\end{equation}  
The assertions in the following proposition can be  proved immediately from  (\ref{c56}).

\bpr \label{f81}
\begin{itemize}
\item[(i)]
For any $\alpha \in I$,  the following equalities hold:
\begin{equation}
\mr{Tr}_\mbQ (\mr{mult} (\alpha ) \, | \, \mbQ \otimes_\mbZ\Fus_N) = \mr{Tr}_\mbQ (\mr{mult} (\alpha^\dual) \, | \,\mbQ \otimes_\mbZ\Fus_N) = \sum_{\lambda \in I} N (\alpha + \lambda + \lambda^\dual).
\end{equation}
\item[(ii)]
 For any $\alpha$, $\beta \in I$,  the following equality holds:
\begin{align} \label{f11}
\mr{Tr}_\mbQ (\mr{mult} (\alpha \ast \beta^\dual) \, | \, \mbQ \otimes_\mbZ \Fus_N) = \sum_{\lambda_1, \lambda_2 \in I} N (\alpha + \lambda_1 + \lambda_2) \cdot N (\beta + \lambda_1 + \lambda_2).
\end{align}
\item[(iii)]
For any positive integer $l$ and any $\alpha_1, \cdots, \alpha_l \in I$,  the following equality holds:
\begin{equation} \label{f40}
\mr{Tr}_\mbQ(\mr{mult} (\prod_{j =1}^l\alpha_j) \, | \, \mbQ \otimes_\mbZ \Fus_N) = \sum_{\lambda \in I} N ((\sum_{j =1}^l \alpha_j) +\lambda + \lambda^\dual).
\end{equation}
\end{itemize}
\epr

The following lemma will be used in the proof of Propositions \ref{f9} and  \ref{c33}.

\ble \label{f5} 
Let $x$ be an element of $\mbR \otimes_\mbZ \Fus_N$, where $\mbR$ denotes the field of real numbers.
Then, the following assertions hold:
\begin{itemize}
\item[(i)]
The following inequality holds:
\begin{equation}
\mr{Tr}_\mbR (\mr{mult} ( x \ast x^\dual) \, | \, \mbR \otimes_\mbZ\Fus_N)\geq 0.
\end{equation}
\item[(ii)]
Assume  that  
$\mr{Tr}_\mbR (\mr{mult} ( x \ast  x^\dual)\, | \, \mbR \otimes_\mbZ\Fus_N)=0$.
Then, 
 the following inclusion relation holds:
\begin{equation}
\mbR \otimes_\mbZ\Fus'_N \subseteq \mr{Ann}(x),
\end{equation}
 where $\mr{Ann}(x)$ denotes the annihilator of $x$ in $\mbR \otimes_\mbZ\Fus_N$.
\end{itemize}
\ele
\begin{proof}
First, we  consider assertion (i).
Let us   write 
$x = u \cdot \varepsilon + \sum_{j =1}^J u_i  \cdot \alpha_i$,
 where  $J$ is a nonnegative integer and  $u$, $u_1, \cdots, u_J \in \mbR$.
Then, we have a sequence of equalities
\begin{align} \label{f8}
&  \ \ \ \ \mr{Tr}_\mbR (\mr{mult} (x \ast x^\dual) \, | \, \mbR \otimes_\mbZ\Fus_{N}) \\
& = u^2 \cdot \mr{Tr}_\mbR (\mr{id}_{\mbR \otimes_\mbZ \scalebox{0.70}{$\Fus_{N}$}} \, | \, \mbR \otimes_\mbZ\Fus_N) 
+  u \cdot \sum_{j=1}^J u_j \cdot \mr{Tr}_\mbR ( \mr{mult} (\alpha_j) \, | \, \mbR\otimes_\mbZ\Fus_N) \notag \\
& \ \ \ + u \cdot \sum_{j=1}^J u_j \cdot \mr{Tr}_\mbR ( \mr{mult} (\alpha^\dual_j) \,  |  \, \mbR \otimes_\mbZ\Fus_N)  \notag \\
& \ \ \  +  \sum_{j, j' =1}^J u_j \cdot u_{j'} \cdot \mr{Tr}_\mbR ( \mr{mult} (\alpha_j \ast \alpha^\dual_{j'}) \, | \, \mbR \otimes_\mbZ\Fus_N) \notag  \\
 & =u^2 \cdot \sum_{\lambda \in I} 1 +2 \cdot u \cdot \sum_{\lambda \in I}\sum_{j =1}^J u_j \cdot N (\alpha_j + \lambda + \lambda^\dual) \notag \\
 &  \ \ \ + \sum_{\lambda_1, \lambda_2 \in I}\sum_{j, j' =1}^J u_j \cdot u_{j'}  \cdot N (\alpha_j + \lambda_1 + \lambda_2) \cdot  N (\alpha_{j'} + \lambda_1 + \lambda_2) \notag  \\
& = \sum_{\lambda \in I}(u + \sum_{j =1}^J u_j \cdot N (\alpha_j + \lambda + \lambda^\dual))^2  + \hspace{-2mm} \sum_{\lambda_1, \lambda_2 \in I, \lambda_1 \neq \lambda_2^\dual}(\sum_{j=1}^J u_j \cdot N (\alpha_j + \lambda_1 + \lambda_2))^2, \hspace{-3mm}
 \notag
\end{align}
where the second equality follows from Proposition \ref{f81}, (i) and (ii).
Since the rightmost  of this sequence is nonnegative, 
we finish  the proof of assertion (i).

Next, we consider assertion (ii).
The assumption  and  (\ref{f8}) together  imply that
\begin{equation} \label{f94}
u + \sum_{j =1}^J u_j \cdot N (\alpha_j +\lambda + \lambda^\dual) = 0 \ \ 
\text{and} \ \  \sum_{j=1}^J N (\alpha_j + \lambda_1 + \lambda_2) = 0
\end{equation}
for any $\lambda, \lambda_1, \lambda_2 \in I$   with $\lambda_1 \neq \lambda_2^\dual$.
Hence, for any $\alpha \in I$, we have
\begin{align}
& \ \ \ \ x \cdot \alpha \\
&   = u \cdot \alpha + \sum_{j =1}^J u_j \cdot \alpha_j \cdot \alpha \notag \\
&  = u  \cdot \alpha + \sum_{j =1}^J u_j \cdot \sum_{\lambda \in I} N (\alpha_j + \alpha + \lambda^\dual) \cdot  \lambda  \notag   \\
&= (u  + \sum_{j =1}^J u_j \cdot N (\alpha_j + \alpha + \alpha^\dual)) \cdot  \alpha +  \sum_{\lambda \in I, \lambda \neq \alpha} \sum_{j =1}^J u_j \cdot N (\alpha_j + \alpha + \lambda^\dual)  \cdot \lambda \notag \\
  & \hspace{-1.5mm}\stackrel{(\ref{f94})}{=} 0.\notag
\end{align}
Since the $\mbR$-vector space $\mbR \otimes_\mbZ\Fus'_N$ is generated by the  set  $I$, 
this implies that $x$ annihilates any element of $\mbR \otimes_\mbZ \Fus'_N$.
This completes the proof of assertion (ii).
\end{proof}

\bpr \label{f9} 
Let $y$ be an element of $\mbR \otimes_\mbZ \Fus'_N$.
Then, the following assertions hold:
\begin{itemize}
\item[(i)]
The following 
 inclusion relation holds:
\begin{equation}
\mbR \otimes_\mbZ \sqrt{(0)} \subseteq \mr{Ann}(y)
\end{equation}
 (cf. Lemma \ref{f5} for the definition of $\mr{Ann} (-)$).
In particular, $y \in \mbR \otimes_\mbZ \sqrt{(0)} \left(\subseteq \mbR \otimes_\mbZ \Fus'_N\right)$  if and only if $y^2 =0$.
\item[(ii)]
If $\overline{y}$ denotes the image of $y$ via the surjection $\mbR \otimes_\mbZ \Fus_N \migi \mbR \otimes_\mbZ\Fus_{N \mr{red}}$, then the following equality holds:
\begin{equation}
\mr{Tr}_\mbR (\mr{mult} (y) \, | \, \mbR \otimes_\mbZ \Fus_N) = \mr{Tr}_\mbR ( \mr{mult} (\overline{y}) \,  | \, \mbR \otimes_\mbZ \Fus_{N \mr{red}}).
\end{equation}
\end{itemize}
\epr
\begin{proof}
Let us take an element $x$ of $\mbR \otimes_\mbZ \sqrt{(0)}$.
Since the endomorphism $\mr{mult} (x \ast x^\dual)$ of $\mbR \otimes_\mbZ \Fus_N$
 is nilpotent,  the equality 
  $\mr{Tr}_\mbR (\mr{mult} (x \ast x^\dual) \, | \, \mbR \otimes_\mbZ \Fus_N) =0$ holds.
  Hence, 
 by Lemma \ref{f5} (ii),
we have $y \ast x =0$.
This completes the proof of  assertion (i).

Assertion (ii) follows from the former assertion and an equality
\begin{equation} \label{f190}
\mr{Tr}_\mbR (\mr{mult} (y) \, | \, \mbR \otimes_\mbZ \Fus_N) = \mr{Tr}_\mbR ( \mr{mult} (\overline{y}) \,  |  \, \mbR \otimes_\mbZ\Fus_{N \mr{red}}) + \mr{Tr}_\mbR (\mr{mult} (y) \, | \, \mbR \otimes_\mbZ\sqrt{(0)})
\end{equation}
arising from  the short exact sequence 
\begin{equation}
0 \longmigi \mbR \otimes_\mbZ \sqrt{(0)} \longmigi \mbR \otimes_\mbZ \Fus_N  \longmigi \mbR \otimes_\mbZ \Fus_{N \mr{red}} \longmigi 0 
\end{equation}
of $\mbR$-vector spaces.
\end{proof}

Note that the nilradical  $\sqrt{(0)}$ of $\Fus_N$ is closed under  the involution $(-)^\dual$.
Hence, this involution induces  an involution on 
  the reduced ring $\Fus_{N \mr{red}}$, as well as on 
 the $K$-algebra $K \otimes_\mbZ \Fus_{N \mr{red}}$ for any field $K$ over $\mbQ$.
  Since $\mbQ \otimes_\mbZ \Fus_{N \mr{red}}$ is a reduced finite-dimensional $\mbQ$-algebra, it admits a decomposition
\begin{equation} \label{f70}
\xi_\mbQ :  \mbQ \otimes_\mbZ \Fus_{N \mr{red}} \isom \prod_{j =1}^J K_j
\end{equation}
 into  the direct  product of a finite number of  finite extensions $K_1, \cdots, K_J$ ($J \geq 2$) of $\mbQ$.
One may choose $\xi_\mbQ$  so that $K_1 = \mbQ$ and 
 the square diagram 
\begin{align} \label{Eh5007}
\vcenter{\xymatrix@C=46pt@R=36pt{
\mbQ \otimes_\mbZ \Fus_{N\mr{red}}\ar[r]_-{\sim}^-{\xi_\mbQ} \ar[d]_-{\mr{id}_\mbQ \otimes \mr{aug}_{\scalebox{0.5}{$\Fus_{N \mr{red}}$}}}& \prod_{j =1}^J K_j  \ar[d]^-{(a_1, \cdots, a_J) \mapsto a_1}\\
\mbQ \ar[r]^-{\sim}_-{\mr{id}_\mbQ}& K_1
}}
\end{align}
is commutative.
 Each factor of the product  $\prod_{j =1}^J K_j$ corresponds to an indecomposable idempotent of $\mbQ \otimes_\mbZ \Fus_{N \mr{red}}$.
  Hence, the construction of such a  decomposition is canonical in the following  sense:
    the collection of (the isomorphism classes of)   fields $K_1, \cdots, K_J$ is uniquely determined up to  ordering, and (after fixing an ordering in $K_1, \cdots, K_J$) $\xi_\mbQ$ is  uniquely determined  up to composition with a product $\prod_{j=1}^J \delta_j$ of  automorphisms $\delta_j : K_j \isom K_j$ ($j = 1, \cdots, J$).
 
\bpr \label{c33} 
Let $K_1, \cdots, K_J$ and $\xi_\mbQ$ be as in (\ref{f70}).
 Then, the following assertions hold:
\begin{itemize}
\item[(i)]
Each $K_j$ ($j =1, \cdots, J$)
   satisfies exactly one of the  two conditions (a),  (b) described as follows:
\begin{itemize}
\vspace{1mm}
\item[(a)]
$K_j$ is a totally real extension of $\mbQ$;
\vspace{1mm}
\item[(b)]
$K_j$ is a totally imaginably extension of $\mbQ$, which is a quadratic extension of a totally real extension $K'_j$ of $\mbQ$.
\end{itemize}
\vspace{1mm}

\item[(ii)]
For each $j \in \{ 1, \cdots, J \}$, we shall  consider the involution $\eta_j$ of $K_j$ defined as follows:
\begin{itemize}
\vspace{1mm}
\item
If $K_j$ satisfies the condition (a), then $\eta_j := \mr{id}_{K_j}$;
\vspace{1mm}
\item
If $K_j$ satisfies the condition (b), then $\eta_j$ is defined to be  the nontrivial automorphism over the quadratic extension $K'_j$.
\end{itemize}  
\vspace{1mm}
Moreover, we shall equip  $\prod_{j =1}^J K_j$ with  the involution given by the product $\prod_{j =1}^J \eta_j$ of  the involutions $\eta_j$.
Then, the isomorphism $\xi_\mbQ$ is compatible with the  respective involutions on $\mbQ \otimes_\mbZ \Fus_{N_{\mr{red}}}$ and $\prod_{j=1}^J K_j$.
\end{itemize}
\epr
\begin{proof}
We consider simultaneously assertions (i) and (ii). 
For each $j = 1, \cdots, J$, let us fix a decomposition  of the $\mbR$-algebra $\mbR \otimes_\mbQ K_j$: 
\begin{equation} \label{f558}
\mbR \otimes_\mbQ K_j  \isom \prod_{l =1}^{J_j} L_{(j, l)},
\end{equation}
 where $L_{(j, l)} = \mbR$ or $\mbC$ ($l =  1, \cdots, J_j $).
Here,  we set 
\begin{equation}
\mcJ := \left\{ (j, l) \, | \,  \text{$1 \ \leq j \leq J$ and  $1 \leq l \leq J_j$}\right\}.
\end{equation}
 The isomorphism $\xi_\mbQ$  and   (\ref{f558}) for all $j$'s give an isomorphism of $\mbR$-algebras
\begin{equation} \label{f111}
\xi_\mbR :  \mbR \otimes_\mbZ \Fus_{N \mr{red}} \isom \prod_{a \in \mcJ} L_{a}.
\end{equation}

Now, let us consider 
the involution $\theta$ of $\prod_{a \in \mcJ} L_{a}$ corresponding to  the involution of 
$\mbR \otimes_\mbZ \Fus_{N \mr{red}}$
  via   $\xi_\mbR$.
We can find
a collection of data
\begin{equation}
\varTheta, \ \ \left\{ \theta_a \, | \, a \in \mcJ \right\},
\end{equation}
where 
\begin{itemize}
\item
$\varTheta$ denotes  an involution  of the set $\mcJ$;
\item
$\theta_a$ (for each $a \in \mcJ$) denotes an isomorphism of fields $L_a \isom L_{\varTheta (a)}$ 
which makes 
the following square diagram commute:
\begin{align} 
\vcenter{\xymatrix@C=46pt@R=36pt{
L_a \ar[r] \ar[d]_-{\theta_{a}}^-{\wr}&
 \prod_{a \in \mcJ} L_a \ar[d]_-{\wr}^-{ \theta}\\
L_{\varTheta (a)} \ar[r]& 
 \prod_{a \in \mcJ} L_a,
}}
\end{align}
where the upper and lower horizontal  arrows denote the inclusions into the $a$-th and $\varTheta (a)$-th factors respectively. 
\end{itemize}
(Indeed, each factor $L_{a}$ corresponds  to an indecomposable idempotent  in $\prod_{a \in \mcJ} L_a$, and such an  idempotent is mapped by the involution $\theta$ to an indecomposable idempotent, which corresponds to a factor $L_{a'}$ for a unique $a'$; the resulting assignment $a \mapsto a'$ defines an involution of $\mcJ$, which we denote by $\varTheta$.)
In particular,  the equality $\theta_{\varTheta (a)} \cdot \theta_{a} = \mr{id}_{L_a}$ holds.

In the following, we shall prove the claim that {\it $\varTheta$ equals  the identity map of $\mcJ$}.
Suppose, on the contrary,  that  there exists  an element $a_0 \in \mcJ$ with  $a_0 \neq \varTheta (a_0)$.
Let us identify $L_{a_0} \times L_{\varTheta (a_0)}$ as an $\mbR$-subvector space of $\prod_{a \in \mcJ} L_a$ via the inclusion 
$L_{a_0} \times L_{\varTheta (a_0)} \migiincl 
\prod_{a \in \mcJ} L_a$
into the $(a_0, \varTheta (a_0))$-th factor.
The involution $\theta$ of $\prod_{a \in L} L_a$ induces, by restricting, an involution of  $L_{a_0} \times L_{\varTheta (a_0)}$; it is   given by  $(y, z) \mapsto (\theta_{\varTheta (a_0)}(z),\theta_{a_0}(y))$.
Now, let us consider the composite 
\begin{equation} \label{f20}
L_{a_0} \times L_{\varTheta (a_0)} \xrightarrow{\mr{inclusion}} \prod_{a \in \mcJ} L_a \xrightarrow{\xi_\mbR^{-1}} 
\mbR \otimes_\mbZ \Fus_{N \mr{red}} \xrightarrow{ \mr{id}_\mbR \otimes  \mr{aug}_{\scalebox{0.5}{$\Fus_{N \mr{red}}$}}} \mbR.
\end{equation}
The third arrow $ \mr{id}_\mbR \otimes \mr{aug}_{\scalebox{0.7}{$\Fus_{N \mr{red}}$}}$ is, by construction, invariant under composition with the involution of $\mbR \otimes_\mbZ \Fus_{N \mr{red}}$, 
so  (\ref{f20})  is invariant under composition with the involution of $L_{a_0} \times L_{\varTheta (a_0)}$.
This composite   preserves  the respective multiplications of the domain and codomain. 
It follows that  (\ref{f20}) must  be
 the zero map.
Hence, since $\mbR \otimes_\mbZ (\Fus'_N / \sqrt{(0)})$ coincides with the kernel of $\mr{id}_\mbR \otimes \mr{aug}_{\scalebox{0.7}{$\Fus_{N \mr{red}}$}}$ we have the inclusion relation
 \begin{equation} \label{122}
 L_{a_0} \times L_{\varTheta (a_0)}  \subseteq \xi_\mbR (\mbR \otimes_\mbZ(\Fus'_N / \sqrt{(0)})).
 \end{equation} 
By  Lemma \ref{f5}, (i),  and Proposition \ref{f9}, (ii), we have the following sequence of equalities  
for each element  $x \in L_{a_0} \times L_{\varTheta (a_0)}$:
 \begin{align} \label{f128}
& \  \ \  \ 
\mr{Tr}_\mbR (\mr{mult} (x \cdot  \theta (x)) \, | \, L_{a_0} \times L_{\varTheta (a_0)}) \\
&  = \mr{Tr}_\mbR (\mr{mult} (x \cdot \theta (x)) \, | \, \prod_{a\in \mcJ} L_a) \notag \\
 & = \mr{Tr}_\mbR ( \mr{mult} (\xi_\mbR^{-1}(x) \ast \xi_\mbR^{-1}(\theta (x))) \,  | \, \mbR\otimes_\mbZ \Fus_{N \mr{red}})  \notag \\
  &= \mr{Tr}_\mbR (\mr{mult} (w \ast w^\dual) \, | \, \mbR \otimes_\mbZ\Fus_N) \notag \\
  & \geq 0,\notag
 \end{align}
where $w$ denotes an element of $\mbR \otimes_\mbZ \Fus_N$  mapped to  $\xi_\mbR^{-1}(x)$  by the surjection $\mbR \otimes_\mbZ\Fus_N \migisurj  \mbR \otimes_\mbZ \Fus_{N \mr{red}}$.
On the other hand, any  element in $L_{a_0} \times L_{\varTheta (a_0)}$
of the form 
  $x \cdot \theta (x)$ (for some $x \in L_{a_0} \times L_{\varTheta (a_0)}$)
 lies in the image of 
the injection
\begin{align} \label{f100}
\kappa : L_{a_0} \migiincl L_{a_0} \times L_{\varTheta (a_0)};
z   \mapsto   (z,   \theta_{a_0} (z)),
\end{align}
and  the map 
$L_{a_0} \migi \mbR$ given by assigning 
$z \mapsto \mr{Tr}_\mbR ( \mr{mult} (\kappa (z)) \, | \, L_{a_0} \times L_{\varTheta (a_0)})$
 is, in fact,  an $\mbR$-linear injection.
 In particular, there exists an element $x'$ of $L_{a_0} \times L_{\varTheta (a_0)}$ satisfying that
 \begin{equation} \label{r457}
  \mr{Tr}_\mbR (\mr{mult} (x' \cdot \theta (x'))\, | \, L_{a_0} \times L_{\varTheta (a_0)}) < 0.
  \end{equation}
By  (\ref{f128}) and (\ref{r457}), we obtain  a  contradiction.
This completes the proof of  the claim.  

Therefore,  the involution of $\prod_{j=1}^J K_j$ induced, via  $\xi_\mbQ$,
 by the involution of $\mbQ \otimes_\mbZ \Fus_N$ may be described  as the product $\prod_{j =1}^J \eta_j$, where each $\eta_j$ is an involution of $K_j$. 
Let us fix   $j_0 \in \{ 1, \cdots, J\}$, and then,  write
$r_1$ (resp., $r_2$) for the cardinality of the subset 
\begin{align}
\mcJ_0^{\cong \mbR} := \left\{ (j_0, l) \in \mcJ \, | \, L_{(j_0, l)} \cong  \mbR\right\} \  \left(\text{resp.}, \  \mcJ_0^{\cong\mbC} := \left\{ (j_0, l) \in \mcJ \, | \, L_{(j_0, l)} \cong  \mbC\right\}\right)
\end{align}
 of $\mcJ$.
In particular,  $\mbR \otimes_\mbQ K_{j_0} \cong \mbR^{\times r_1}\times \mbC^{\times r_2}$.
Here, we shall determine  the involution $\theta_a : L_a \isom L_{\varTheta (a)}$ ($= L_a$) of each $L_a$ in two cases (i.e., the cases ``$a \in \mcJ_0^{\cong \mbR}$" and ``$a \in \mcJ_0^{\cong \mbC}$").
On the one hand, since $\theta_a$ is  an $\mbR$-algebra automorphism, 
 {\it the involution $\theta_{a} : L_a \isom L_a$ coincides with  the identity morphism if $a \in \mcJ_0^{\cong \mbR}$}.
On the other hand,  we can verify  that 
{\it the involution $\theta_{a}$  coincides with the complex conjugation if $a \in \mcJ_0^{\cong \mbC}$}.
Indeed, the $\mbR$-linear  composite 
\begin{equation} \label{f132}
L_a \xrightarrow{\mr{inclusion}} \prod_{a \in \mcJ} L_a \xrightarrow{\xi_\mbR^{-1}}  \mbR \otimes_\mbZ \Fus_{N \mr{red}}  \xrightarrow{\mr{id}_\mbR \otimes \mr{aug}_{\scalebox{0.5}{$\Fus_{N \mr{red}}$}}} \mbR
\end{equation}
 is compatible with the respective multiplications on $L_a$ and $\mbR$.
 Hence,  if $a \in \mcJ_0^{\cong \mbC}$, then this composite must be the zero map since
 $[L_a : \mbR] >1$.
This implies the inclusion relation 
\begin{equation}
L_a \subseteq \xi_\mbR (\mbR \otimes_\mbZ(\Fus'_N /\sqrt{(0)}) )  \left(= 
\mr{Ker} ((\mr{id}_\mbR \otimes \mr{aug}_{\scalebox{0.7}{$\Fus_{N \mr{red}}$}})\circ \xi_\mbR^{-1})\right).
\end{equation}
By considering a sequence of equalities exactly like  (\ref{f128}), 
we obtain the inequality
\begin{align}
\mr{Tr}_\mbR (\mr{mult} (x \cdot \theta_a (x)) \, | \, L_a) \geq 0
\end{align}
 for any $x \in L_a$.
It follows that the involution $\theta_a : \mbC \isom \mbC$ cannot be the identity morphism, i.e.,  it is nothing but the complex conjugation.
This proves  the second italicized assertion described above.

Now, let us suppose that $r_2 =0$.
The $\mbR$-algebra $\mbR \otimes_\mbQ K_{j_0}$ is isomorphic to the direct product of finite copies of $\mbR$, and 
the involution $\eta_{j_0}$, which corresponds to $\prod_{l =1}^{J_{j_0}} \theta_{(j_0, l)}$ via  (\ref{f558}),  is trivial (by the first italicized assertion in the previous paragraph).
Hence,  $K_{j_0}$ satisfies the condition (a) and $\eta_{j_0} = \mr{id}_{K_{j_0}}$.

On the other hand, suppose that  $r_2 > 0$. 
If $K_{j_0}^{\eta_{j_0}}$ denotes the subfield of $K_{j_0}$ fixed  by  $\eta_{j_0}$, then $\mbR \otimes_\mbQ K_{j_0}^{\eta_{j_0}}$ is isomorphic to $\mbR^{\times (r_1 + r_2)}$ (by the second italicized assertion in the previous paragraph).
In particular, the field $K_{j_0}^{\eta_{j_0}}$
 is totally real and strictly smaller than $K_{j_0}$, i.e, $[K_{j_0} : K_{j_0}^{\eta_{j_0}}]\geq 2$.
 Also, we have 
\begin{equation}
r_1 + 2  r_2=  \mr{dim}_\mbR (\mbR \otimes_\mbQ K_{j_0}) = [K_{j_0} : \mbQ] \geq 2  [K_{j_0}^{\eta_{j_0}} : \mbQ] = 2   (r_1 + r_2).
\end{equation}
It follows  that  $r_1 =0$.
Hence,  $K_{j_0}$  satisfies the condition (b) and the  involution $\eta_{j_0}$ coincides with  the nontrivial automorphism over $K_{j_0}^{\eta_{j_0}}$.
This completes the proof of Proposition \ref{c33}.
 \end{proof}

\bco \label{c335} 
Let $K_1, \cdots, K_J$ and $\xi_\mbQ$ be
  as in (\ref{f70}).
Then, the following assertions hold:
\begin{itemize}
\item[(i)]
 There exist  a local $\mbQ$-algebra $\widetilde{K}_1$ with residue field 
 $K_1 \left(= \mbQ\right)$ and 
an isomorphism of $\mbQ$-algebras
\begin{equation}
\widetilde{\xi}_\mbQ : \mbQ \otimes_\mbZ\Fus_N  \isom \widetilde{K}_1 \times \prod_{j =2}^J K_j
\end{equation}
 fitting into 
the following  square diagram:
\begin{align} 
\vcenter{\xymatrix@C=46pt@R=36pt{
\mbQ \otimes_\mbZ\Fus_N  
\ar[r]_-{\sim}^-{\widetilde{\xi}_\mbQ} \ar[d]_-{\mr{quotient}}&\widetilde{K}_1 \times \prod_{j =2}^J K_j \ar[d]^-{\mr{quotient}}\\
\mbQ \otimes_\mbZ \Fus_{N \mr{red}} \ar[r]^-{\sim}_-{ \xi_\mbQ }& \prod_{j =1}^J K_0.
}}
\end{align}
  In particular, we have an isomorphism of $\mbQ$-vector spaces
  \begin{equation} \label{f598}
  \widetilde{K}_0 \isom \mbQ \oplus (\mbQ \otimes_\mbZ\sqrt{(0)}).
  \end{equation}

 \item[(ii)]
Let $\eta_1, \cdots, \eta_J$
 be as in Proposition \ref{c33} (hence $\eta_1 = \mr{id}_{K_1}$).
 Then, there exists a unique involution $\widetilde{\eta}_1$ of $\widetilde{K}_1$ such that 
 the involution $\widetilde{\eta}_1 \times \prod_{j =2}^J \eta_j$ of the product  $\widetilde{K}_1 \times \prod_{j =2}^J K_j$ 
 corresponds, via $\widetilde{\xi}_\mbQ$, to   the  involution of $\mbQ \otimes_\mbZ\Fus_N$. 
\end{itemize}
\eco
\begin{proof}
First, we shall consider assertion (i).
Since $\mbQ \otimes_\mbZ \Fus_N$ is a finite-dimensional $\mbQ$-algebra, there exists a collection  of data $(\{\widetilde{K}_j\}_{j=1}^J, \widetilde{\xi}_\mbQ)$, where  
each $\widetilde{K}_j$ denotes a  local $\mbQ$-algebra with residue field  $K_j$
 and $\widetilde{\xi}_\mbQ$ denotes an isomorphism of $\mbQ$-algebras
$\mbQ \otimes_\mbZ\Fus_N \isom   \prod_{j =1}^J \widetilde{K}_j$
 fitting into  
the  following  square diagram
\begin{align} \label{Eh5007}
\vcenter{\xymatrix@C=46pt@R=36pt{
\mbQ \otimes_\mbZ \Fus_N  \ar[r]_-{\sim}^-{\widetilde{\xi}_\mbQ} \ar[d]_-{\mr{quotient}}&\prod_{j =1}^J  \widetilde{K}_j  \ar[d]^-{\mr{quotient}}\\
\mbQ \otimes_\mbZ\Fus_{N \mr{red}} \ar[r]^-{\sim}_-{\xi_\mbQ}& \prod_{j =1}^J K_j.
}}
\end{align}
Then, the problem is to verify  that {\it $\widetilde{K}_j$ is isomorphic to $K_j$ for every $j > 1$}.
Let us fix $j\in \{2, \cdots, J \}$, and suppose that we are given an element $x_j$  in  the kernel of the surjection $\widetilde{K}_j \migisurj K_j$.
If  $\iota_j : \widetilde{K}_j \migiincl  \prod_{j =1}^J  \widetilde{K}_j$ denotes the inclusion into the $j$-th factor,
then $\widetilde{\xi}_\mbQ^{-1} (\iota_j (x_j)) \in \mbQ \otimes_\mbZ \sqrt{(0)}$ and $\widetilde{\xi}_\mbQ^{-1} (\iota_j (1)) \in \mbQ \otimes_\mbZ\Fus'_N$.
It follows from Proposition \ref{f9}, (i),  that 
\begin{equation}
\widetilde{\xi}_\mbQ^{-1} (\iota_j (x_j)) = \widetilde{\xi}_\mbQ^{-1} (\iota_j (x_j)) \cdot \widetilde{\xi}_\mbQ^{-1} (\iota_j (1)) =0,
\end{equation}
 so  the equality $x_j =0$ holds.
Thus,  the surjection $\widetilde{K}_j \migisurj K_j$ turns out to be  injective, i.e.,  bijective.
This completes the proof of assertion (i).

Assertion (ii) follows immediately  from assertion (i).
Indeed,  $\mbQ \otimes_\mbZ\sqrt{(0)}$ is closed  under the involution of $\mbQ \otimes_\mbZ \Fus_N$, and  the resulting involution of $\mbQ \otimes_\mbZ \sqrt{(0)}$ induces, via  (\ref{f598}), an involution of $\widetilde{K}_0$   satisfying  the required condition.
 \end{proof}

\section{Factorization property of  a fusion rule}\label{QG1109}

We shall write 
\index{$\mfS_N$, $\mfS_N^{\epsilon}$, $\mfS_N^+$, $\widetilde{\mfS}_N$}
\begin{equation} \label{f160}
\mfS_N  \left(:=\mr{Hom} (\Fus_N, \mbC)\right)
\end{equation}
 for the set of ring homomorphisms $\Fus_N \migi \mbC$, i.e., $\mfS_N$ is the underlying set of $\mr{Spec}(\Fus_N) (\mbC)$.
 One may think of $\mfS_N$ as the set of $\mbC$-algebra homomorphisms  $\mbC \otimes_\mbZ\Fus_{N \mr{red}} \migi \mbC$.  
Also, the product  $\mbC^{\mfS_N}$ of copies (indexed by the set $\mfS_N$) of $\mbC$ admits a structure of  $\mbC$-algebra in the natural fashion and  an involution induced by  the complex conjugation $\overline{(-)}$ of each factor $\mbC$  in $\mbC^{\mfS_N}$.

\bpr \label{ccc35} 
The map 
\begin{equation} \label{c39}
\mbC \otimes_\mbZ \Fus_{N \mr{red}}  \migi \mbC^{\mfS_N}
\end{equation}
 given by  $x \mapsto (\chi (x))_{\chi \in \mfS_N}$ for any $x \in \mbC \otimes_\mbZ \Fus_{N \mr{red}}$ is an isomorphism of $\mbC$-algebras  compatible with the  involutions, i.e., $(\chi (x^\dual))_{\chi \in \mfS_N} = (\overline{\chi (x)})_{\chi \in \mfS_N}$.
\epr
\begin{proof}
 The  assertion follows  immediately from Proposition \ref{c33}.
 \end{proof}

Let us write
\begin{equation} \label{w33333}
\mr{Cas}_N : = \sum_{\lambda \in I} \lambda \ast \lambda^\dual \left(\in \Fus_N\right).
\index{$\mr{Cas}_N$}
\end{equation}
Then, an explicit knowledge 
 of the ring homomorphisms $\Fus_N \migi \mbC$
(i.e., of the isomorphism $\mbC \otimes_\mbZ \Fus_N \isom \mbC^{\mfS_N}$ asserted in Proposition \ref{ccc35}) and the element $\mr{Cas}_N$ of $\Fus_N$
 will allow us to perform some  computation that we need in the ring $\Fus_N$ as follows.

\bpr \label{c35} 
Suppose that   we are given a collection of   map of sets
 \begin{equation}
 N_g : \mbN^{I}_{\geq 3-2g} \migi \mbZ
 \end{equation}
($g =0, 1, 2, \cdots$)
 such   that $N_0 = N$ and
\begin{equation} \label{f43}
N_g (x) = \sum_{\lambda \in I} N_{g-1} (x + \lambda + \lambda^\dual)
\end{equation}
for any $g >0$ and $x \in \mbN^{I}_{\geq 3-2g}$. Then, the following assertions hold:
\begin{itemize}
\item[(i)]
For any nonnegative integers $g_1$, $g_2$ and any elements $x \in \mbN^{I}_{\geq 3-2g_1}$, $y \in \mbN^{I}_{\geq 3 -2 g_2}$,  the following equality holds:
\begin{equation} \label{c61}
N_{g_1 +g_2} (x + y) = \sum_{\lambda \in I} N_{g_1} (x + \lambda) \cdot N_{g_2} (y + \lambda^\dual).
\end{equation}
\item[(ii)]
For any pair of nonnegative integers $(g, r)$ with  $2g-2+r >0$ and any $\rho_1, \cdots, \rho_r \in I$,  the following equality holds:
\begin{equation} \label{c60}
N_g  (\sum_{i =1}^r \rho_i) = \sum_{\chi \in \mfS_N} \chi (\mr{Cas}_N)^{g-1} \cdot \prod_{i =1}^r \chi (\rho_i).
\end{equation} 
Here, if $r =0$ (hence $g>1$), then this equality means that
\begin{equation} \label{f5c60}
N_g  (0) = \sum_{\chi \in \mfS_N} \chi (\mr{Cas}_N)^{g-1}.
\end{equation} 
\end{itemize}
\epr
\begin{proof}
We shall consider assertion (i).
By induction on $g$, one may verify  from  (\ref{f43}) that,  for any $g \geq 1$ and $z \in \mbN^{I}_{3-2g}$, the equality
\begin{equation} \label{f66}
N_g (z) = \sum_{\lambda_1, \cdots, \lambda_g \in I } N (z + \sum_{i =1}^g(\lambda_i + \lambda_i^\dual))
\end{equation}
holds. Thus, the following sequence of equalities holds:
\begin{align}
& \ \ \ \ \sum_{\lambda \in I} N_{g_1} (x + \lambda) \cdot N_{g_2} (y + \lambda^\dual) \\
& = \sum_{\lambda \in I}  (\sum_{\lambda_1, \cdots, \lambda_{g_1} \in I } \hspace{-3mm} N (x + \lambda + \sum_{i =1}^{g_1}(\lambda_i + \lambda_i^\dual))) 
\cdot  (\hspace{-3mm}\sum_{\lambda_1, \cdots, \lambda_{g_2} \in I } \hspace{-3mm} N (y + \lambda^\dual + \sum_{i =1}^{g_2}(\lambda_i + \lambda_i^\dual)) ) \hspace{-4mm}\notag \\
& =\sum_{\lambda_1, \cdots, \lambda_{g_1 +g_2}\in I}   \sum_{\lambda \in I}  N (x + \lambda + \sum_{i=1}^{g_1} (\lambda_i + \lambda^\dual_i)) \cdot N (y + \lambda^\dual + \sum_{i=g_1+1}^{g_1 +g_2} (\lambda_i + \lambda^\dual_i))\notag\\
& = \sum_{\lambda_1, \cdots, \lambda_{g_1 +g_2}\in I}  N (x + y+ \sum_{i =1}^{g_1 + g_2} (\lambda_i + \lambda_i^\dual)) \notag \\
& =  N_{g_1 + g_2} (x +y). \notag
\end{align}
This completes the proof of assertion (i).

Next, we consider assertion (ii).
Let $x$ be an element of $\Fus'_N \left(\subseteq \mbC \otimes_\mbZ\Fus_N\right)$.
 If  $\overline{x}$ denotes  the image of $x$ via the surjection $\Fus_N \migisurj \Fus_{N \mr{red}}$, then 
it  corresponds, via  (\ref{c39}), to $(\chi (\overline{x}))_{\chi \in\mfS_N} \in \mbC^{\mfS_N}$.
The matrix corresponding to the endomorphism $\mr{mult} (\overline{x})$ of $\mbC \otimes_\mbZ \Fus_{N \mr{red}}$ with respect to the standard basis of $\mbC^{\mfS_N}$ coincides with the diagonal matrix with entries $(\chi (\overline{x}))_{\chi \in \mfS_N}$.
Hence, it follows from Proposition \ref{f9}, (ii), that 
\begin{align} \label{c40}
\mr{Tr}_\mbQ(\mr{mult} (x)\, | \, \mbQ \otimes_\mbZ\Fus_N ) &  = \mr{Tr}_\mbC(\mr{mult} (\overline{x}) \, | \, \mbC \otimes_\mbZ\Fus_{N \mr{red}})  \\
& =   \mr{Tr}_\mbC( \mr{mult} ((\chi (\overline{x}))_{\chi \in \mfS_N}) \, | \, \mbC^{\mfS_N})  \notag\\
&  = \sum_{\chi \in \mfS_N} \chi (x). \notag
\end{align}
 
 On the other hand, observe the sequence of equalities
\begin{align} \label{f41}
N_g (\sum_{i=1}^r\rho_i) &=\sum_{\lambda_1, \cdots,  \lambda_g \in I } N ((\sum_{i =1}^r \rho_i )+ \sum_{i=1}^g(\lambda_i + \lambda_i^\dual)) \\
& = \sum_{\lambda_1, \cdots,  \lambda_{g-1} \in I} \sum_{\lambda \in I}
 N ((\sum_{i =1}^r \rho_i +  \sum_{i=1}^{g-1}(\lambda_i + \lambda_i^\dual))+ \lambda + \lambda^\dual) \notag \\
& =  \sum_{\lambda_1, \cdots,  \lambda_{g-1} \in I} \mr{Tr}_\mbQ ( \mr{mult} ( (\prod_{i=1}^r \rho_i) \ast \prod_{i =1}^{g-1} \lambda_i \ast \lambda_i^\dual) \, | \,  \mbQ \otimes_\mbZ\Fus_N) \notag \\
& = \mr{Tr}_\mbQ (\mr{mult} ( \sum_{\lambda_1, \cdots,  \lambda_{g-1} \in I}  (\prod_{i=1}^r \rho_i) \ast \prod_{i =1}^{g-1} \lambda_i \ast \lambda_i^\dual) \, | \,  \mbQ \otimes_\mbZ \Fus_N) \notag \\
& = \mr{Tr}_\mbQ ( \mr{mult} (\mr{Cas}_N^{g-1} \ast \prod_{i=1}^r \rho_i) \, | \,  \mbQ \otimes_\mbZ \Fus_N), \notag
\end{align}
where  the first arrow follows from (\ref{f66}) and  the third  equality follows from Proposition \ref{f81}, (iii).
Thus, by combining  (\ref{c40}) (in the case where $x$ is taken to be the element $\mr{Can}_N^{g-1} \cdot \prod_{i=1}^r \rho_i$) and (\ref{f41}), we obtain a sequence of equalities
 \begin{align}
N_g (\sum_{i =1}^r \rho_i)  &= \mr{Tr}_\mbQ ( \mr{mult} (\mr{Cas}_N^{g-1} \cdot \prod_{i=1}^r \rho_i) \, | \, \mbQ \otimes_\mbZ \Fus_N) \\
&  = \sum_{\chi \in\mfS_N} \chi (\mr{Cas}_N^{g-1} \cdot \prod_{i=1}^r \rho_i) \notag \\
&   =   \sum_{\chi \in \mfS_N} \chi (\mr{Cas}_N)^{g-1} \cdot \prod_{i =1}^r \chi (\rho_i). \notag
 \end{align}
This completes the proof of assertion (ii).
 \end{proof}


In what follows, we shall examine the relationship
 among 
the structure constants 
\begin{equation} \label{c67}
N_{\alpha, \gamma}^\beta := N (\alpha + \beta^\dual + \gamma)
\index{$N_{\alpha, \gamma}^\beta$, $N_\alpha$}
\end{equation}
 (where $\alpha$, $\beta$, $\gamma \in I$) of the pseudo-fusion ring $\Fus_N$.
 Fix an element $\alpha$ of $I$, and
consider 
the restriction $\mr{mult} (\alpha) |_{\mbC \otimes_\mbZ \scalebox{0.7}{$\Fus'_N$}}$ to $\mbC \otimes_\mbZ \Fus'_N$ of the $\mbC$-linear endomorphism $\mr{mult} (\alpha) : \mbC \otimes_\mbZ \Fus_N \migi \mbC \otimes_\mbZ \Fus_N$.
The matrix $N_\alpha$  corresponding to this endomorphism 
 with respect to the basis $I$   is given by 
\begin{equation} \label{QH6004}
N_\alpha = (N^\beta_{\alpha, \gamma})_{(\beta, \gamma) \in I \times I}.
\index{$N_{\alpha, \gamma}^\beta$, $N_\alpha$}
\end{equation}

Next, let us consider the matrix representation of $\mr{mult} (\alpha) |_{\mbC \otimes_\mbZ \scalebox{0.7}{$\Fus'_N$}}$ with respect to another basis.
We shall write 
 \index{$\mfS_N$, $\mfS_N^{\epsilon}$, $\mfS_N^+$, $\widetilde{\mfS}_N$}
\begin{equation} \label{QH6003}
\mfS_N^{\epsilon} \left(:= \mr{Hom} (\Fus_N, \mbC [\epsilon]/(\epsilon^2))\right)
\end{equation}
for the set of ring homomorphisms $\Fus_N \migi \mbC [\epsilon]/(\epsilon^2)$, i.e., $\mfS_N^{\epsilon}  $ is the underlying set of $\mr{Spec}(\Fus_N) (\mbC [\epsilon]/(\epsilon^2))$.
 Write   
 \begin{equation} \label{QH6002}
 \chi_0 : \Fus_N \migi \mbC
 \end{equation}
  for   the composite homomorphism $\Fus_N \xrightarrow{\mr{aug}_{\scalebox{0.5}{$\Fus_N$}}} \mbZ \migiincl \mbC$ (cf. (\ref{f123}) for the definition of $\mr{aug}_{\scalebox{0.7}{$\Fus_N$}}$), and
 \index{$\mfS_N$, $\mfS_N^{\epsilon}$, $\mfS_N^+$, $\widetilde{\mfS}_N$}
 \begin{equation} \label{QH6001}
 \mfS^+_N := \mfS_N \setminus \{ \chi_0 \}.
 \end{equation}
By Corollary \ref{c335},  the natural surjection
\begin{equation} \label{QH6000}
\pi_N : \mfS_N^{\epsilon} \migisurj \mfS_N
\end{equation}
restricts to a bijection $\pi_N |_{\pi_N^{-1} (\mfS^+_N)} : \pi_N^{-1} (\mfS^+_N) \isom  \mfS^+_N$.
 Since the equality $(\mbC \otimes_\mbZ \sqrt{(0)})^2 =0$ holds by virtue of Proposition \ref{f9}, (i),
  the inverse image $\pi_N^{-1} (\chi_0)$ of $\chi_0$  admits a structure of $\mbC$-vector space canonically isomorphic to $(\mbC \otimes_\mbZ \sqrt{(0)})^\vee$.  Namely, 
 \begin{equation} \label{f390}
 \mcT_{\chi_0} \mr{Spec}(\mbC \otimes_\mbZ\Fus_N) 
\left(= \pi_N^{-1} (\chi_0)\right)
  =  (\mbC \otimes_\mbZ \sqrt{(0)})^\vee.
 \end{equation}

 Now, let us fix a basis  $\mfs_{N} \left(\subseteq \pi_N^{-1} (\chi_0)\right)$
  for the $\mbC$-vector space 
 $\pi_N^{-1} (\chi_0) \left(= (\mbC \otimes_\mbZ\sqrt{(0)})^\vee \right)$, 
  and write
  \index{$\mfS_N$, $\mfS_N^{\epsilon}$, $\mfS_N^+$, $\widetilde{\mfS}_N$}
\begin{equation}
\widetilde{\mfS}_N := \mfs_{N} \sqcup \mfS^+_N.
\end{equation}
The  $\mbC$-linear morphisms
\begin{equation} \label{f355}
\mbC \otimes_\mbZ\sqrt{(0)} \migi \mbC^{\mfs_N} \ \ \text{and} \ \  \mbC \otimes_\mbZ (\Fus'_{N }/\sqrt{(0)}) \migi \mbC^{\mfS^+_N} 
\end{equation}
determined uniquely  by assigning $1 \otimes x  \mapsto (\chi (x))_{\chi \in \mfs_N}$ ($x \in \sqrt{(0)}$) and  $1 \otimes \overline{x} \mapsto (\chi (x))_{\chi \in \mfS^+_N}$ ($x \in \Fus'_N$) respectively are isomorphisms.
 
 Since $\mbC \otimes_\mbZ \Fus_N $ is a finite-dimensional $\mbC$-algebra, 
 the $\mbC$-algebra structures on the respective  factors of its  direct product decomposition determine
 a canonical split injection $\mbC \otimes_\mbZ\Fus_{N \mr{red}} \migiincl  \mbC \otimes_\mbZ\Fus_N$ of 
 the natural short exact sequence  of $\mbC$-vector spaces
\begin{equation}
0 \longmigi \mbC \otimes_\mbZ\sqrt{(0)} \longmigi  \mbC \otimes_\mbZ \Fus_N  \longmigi \mbC \otimes_\mbZ \Fus_{N \mr{red}} \longmigi 0.
\end{equation}
By composing with the corresponding  split surjection $\mbC \otimes_\mbZ \Fus_N\migisurj \mbC \otimes_\mbZ \sqrt{(0)}$, we shall regard each element of $\mfs_{N} \left(\subseteq (\sqrt{(0)}\otimes_\mbZ \mbC)^\vee \right)$
  as a $\mbC$-linear morphism  $\mbC \otimes_\mbZ \Fus_N  \migi \mbC$.
Thus,  by taking account of the two  isomorphisms in (\ref{f355}),  we obtain a $\mbC$-linear isomorphism
\begin{equation} \label{f366}
\mbC\otimes_\mbZ\Fus'_N \isom \mbC^{\widetilde{\mfS}_N}
\end{equation}
determined uniquely by  $1 \otimes x  \mapsto (\chi (x))_{\chi \in \widetilde{\mfS}_N}$.

For each $\chi \in \widetilde{\mfS}_N$, 
the sum 
\begin{equation} \label{f367}
\mr{Cas}_{\chi} :=  \sum_{\lambda \in I} |\chi (\lambda)|^2
\index{$\mr{Cas}_{\chi}$}
\end{equation}
(hence $\mr{Cas}_{\chi} = \chi (\mr{Cas}_N)$ if $\chi \in \mfS^+_N$) is necessarily  positive since the set $I$ generates  $\mbC \otimes_\mbZ \Fus'_N$.
It follows that we have a well-defined positive real number $\sqrt{\mr{Cas}_{\chi}}$.
Let $\Sigma_{\mfs_N}$ be the  base-change matrix of $\mbC \otimes_\mbZ \Fus'_N$ from the basis $I$ to the basis $\{ \sqrt{\mr{Cas}_{\chi}} \cdot \chi  \, | \,  \chi \in \widetilde{\mfS}_N \}$ under the identification given by  (\ref{f366}), i.e.,
\begin{equation} \label{f368}
\Sigma_{\mfs_N} := 
\left(\frac{\chi (\lambda)}{ \sqrt{\mr{Cas}_{\chi}} }\right)_{(\chi, \lambda)  \in \widetilde{\mfS}_N \times  I}.
\end{equation}

For each $\alpha \in I$, denote by $D_\alpha$ the diagonal matrix  in   $\mbC^{\widetilde{\mfS}_N \times \widetilde{\mfS}_N} \left( = \mr{End}_\mbC (\mbC^{\widetilde{\mfS}_N})\right)$
 whose entry indexed by $\chi \in \widetilde{\mfS}_N$ equals $0$ if $\chi \in \mfs_N$, and  $\chi (\alpha)$ if $\chi \in \mfS^+_N$.
 That is to say, we set
 \begin{align}
&  \hspace{5mm}\mbC^{\mfs_N} \hspace{25mm} \mbC^{\mfS^+_N} \notag \\
 D_\alpha :=& \left( \begin{array}{@{\,}cc | ccccc@{\,}}
\multicolumn{2}{c |}{\raisebox{-8pt}[0pt][0pt]{\Large $O$}}   &  \multicolumn{5}{c }{\raisebox{-8pt}[0pt][0pt]{\Large $O$}} \\ 
& & & & & & \\ \hline
 & & \chi (\alpha) & 0 & 0 & \cdots & 0 \\
 & & 0 & \chi'(\alpha) & 0 & \cdots & 0 \\
\multicolumn{2}{c |}{\raisebox{-3pt}[0pt][0pt]{\Large $O$}}   & 0 & 0 & \chi''(\alpha) & \cdots & 0 \\
 &   & \vdots &  \vdots & \vdots & \ddots & \vdots \\
 &  & 0 & 0 & 0 & \cdots & \chi'''(\alpha)
 \end{array} \right) \hspace{1mm}\begin{array}{@{\,}c @{\,}} \mbC^{\mfs_N} \\ \\ \\ \mbC^{\mfS^+_N}  \\  \\ \\  \end{array}
 \end{align}
where $O$'s denote zero matrices.
This matrix  coincides with the matrix corresponding to the endomorphism $\mr{mult} (\alpha) |_{\mbC \otimes_\mbZ \scalebox{0.7}{$\Fus'_N$}}$ with respect to the basis $\{ \sqrt{\mr{Cas}_{\chi}} \cdot  \chi  \, | \, \chi \in \widetilde{\mfS}_N \}$  via  (\ref{f366}).
Thus, for any $\alpha \in I$, 
  the following equality   holds:
\begin{equation} \label{c69}
 N_\alpha = \Sigma_{\mfs_N}^{-1} \cdot D_\alpha \cdot \Sigma_{\mfs_N}.
\end{equation}

\bpr \label{f200} 
There exists a basis $\mfs_N$ for the $\mbC$-vector space  $(\mbC \otimes_\mbZ \sqrt{(0)})^\vee$
$\left(=  \mcT_{\chi_0} \mr{Spec}(\mbC \otimes_\mbZ \Fus_N)\right)$
  such that  the associated matrix $\Sigma_{\mfs_N}$ is unitary.
In particular, there exists a unitary matrix $\Sigma$ (i.e., $\Sigma_{\mfs_N}$ for such a basis $\mfs_N$) which  diagonalizes  simultaneously all matrices $N_\alpha \in \mr{End}_\mbQ (\mbQ \otimes_\mbZ\Fus_N)$ ($\alpha \in  I$).
\epr
\begin{proof}
 Let us fix a  basis $\mfs'_N$ for $(\mbC \otimes_\mbZ \sqrt{(0)})^\vee$, which induces  isomorphisms 
 \begin{equation}
 \label{f377}
 \mbC \otimes_\mbC \sqrt{(0)} \isom  \mbC^{\mfs'_N}, \   \text{and} \ \ \mbC \otimes_\mbZ\Fus'_N \isom \mbC^{\widetilde{\mfS}_N}
 \end{equation}
  as in (the first isomorphism of) (\ref{f355}) and (\ref{f366}).
Since  $N_{\alpha^\dual} = N_\alpha^*$ (i.e., the conjugate transpose of $N_\alpha$), 
the commutativity of the ring $\Fus_N$ implies that $[N_\alpha, N_\alpha^*] \ (=[N_\alpha, N_{\alpha^\dual}]) =0$
and $[N_\alpha, N_\beta] =0$ for any $\alpha$, $\beta \in I$.
Namely, the matrices $N_\alpha$ ($\alpha \in I$) are normal and commute with each other.
Hence, by a well-known fact in linear algebra (cf. ~\cite{New}, Theorem 2),   there exists a unitary matrix $\Sigma_0$ which  diagonalizes  simultaneously all matrices $N_\alpha$ ($\alpha \in  I$).
By means of the second  isomorphism in (\ref{f377}), we regard $\Sigma_0$ as a matrix with rows indexed by $\widetilde{\mfS}_N$ and columns indexed by $I$.
Also, after possibly permuting the ordering  of basis, we may assume that 
 \begin{equation} \label{f346}
 N_\alpha = \Sigma_0^{-1} \cdot D_\alpha \cdot \Sigma_0
 \end{equation}
for all $\alpha \in I$.
If we write $\Sigma_1 : =  \Sigma_{\mfs'_N}\cdot \Sigma_0^{-1} \left(= \Sigma_{\mfs'_N}\cdot \Sigma_0^*\right)$, then 
  (\ref{c69}) and (\ref{f346}) yield  an equality $D_\alpha \cdot \Sigma_1 = \Sigma_1 \cdot D_\alpha$ for every  $\alpha$.
This implies that $\Sigma_1$ is the direct sum of a regular matrix $R \in \mr{GL} (\mbC^{\mfs'_N})$ and a {\it unitary} diagonal matrix $D \in \mr{GL}(\mbC^{\mfS^+_N})$, i.e., of the form
 \begin{align}
&  \hspace{4mm}\mbC^{\mfs'_N} \hspace{2mm} \mbC^{\mfS^+_N} \notag \\
\Sigma_1 :=& \left( \begin{array}{@{\,}cc | cc@{\,}}
\multicolumn{2}{c |}{\raisebox{-7pt}[0pt][0pt]{\large $R$}}   &  \multicolumn{2}{c }{\raisebox{-7pt}[0pt][0pt]{\large  $O$}} \\ 
 & & & \\ \hline
\multicolumn{2}{c |}{\raisebox{-7pt}[0pt][0pt]{\large $O$}}   &  \multicolumn{2}{c }{\raisebox{-7pt}[0pt][0pt]{\large $D$}} \\ 
&& &  
 \end{array} \right) \hspace{1mm}\begin{array}{@{\,}c @{\,}} \mbC^{\mfs'_N} \\ \\  \mbC^{\mfS^+_N}.  \\   \end{array}
 \end{align}
In particular, the matrix
\begin{equation} \label{f347}
 \left( \begin{array}{@{\,}c | c@{\,}} 
   R^{-1} &  O  \\ \hline  O  & E
 \end{array} \right) 
\cdot \Sigma_{\mfs'_N}  \left(=  \left( \begin{array}{@{\,}c | c@{\,}}  E &  O  \\  \hline O  & D
 \end{array} \right)  \cdot \Sigma_0
\right)
\end{equation}
is unitary, where $E$
    denote the identity matrices.
Let $\mfs_N$ be  the  basis for $(\mbC \otimes_\mbZ\sqrt{(0)})^\vee$ obtained by transforming the basis $\mfs'_N$ by means of   the automorphism of $\mbC \otimes_\mbZ\sqrt{(0)}$ corresponding to the matrix $R^{-1}$ via  the first  isomorphism $\mbC \otimes_\mbZ \sqrt{(0)} \isom  \mbC^{\mfs'_N}$ in  (\ref{f377}).
Then, the matrix $\Sigma_{\mfs_N}$ associated to $\mfs_N$ coincides with    (\ref{f347}), and hence, the basis $\mfs_N$ satisfies the required condition. 
This completes the proof of the  assertion.
 \end{proof}

\section{The pseudo-fusion ring associated to dormant opers}\label{c8}

Now, let us go back to our  study of 
dormant $\mfg$-opers.
Let us keep the assumption imposed  at the beginning of \S\,\ref{y0174s}.
 We fix an involution of 
 the finite set $\mfc (\mbF_p)$
defined by   $\varrho \mapsto \varrho^\dual :=   (-1) \star \varrho$.
 Then,  the following proposition  holds.

\bpr \label{c2} 
Suppose that the algebraic group $\mbG$ satisfies the conditions $(*)$ and $(**)$  described in \S\,\ref{QG1100}. 
Then, the   map of sets
\index{${^p}N_{\mfg,  0}$, ${^p}N_{\mfg,  g}$}
\begin{align} \label{f186}
{^p}N_{\mfg, 0}  : \mbN^{\mfc(\mbF_p)}_{\geq 3} \migi \mbZ
\end{align}
 given by assigning $\overline{\rho} \mapsto {^p}N_{\mfg, 0}(\overline{\rho})$ (cf. (\ref{d6})) for any 
 $\overline{\rho} \in  \mbN^{\mfc(\mbF_p)}_{\geq 3}$
  defines 
  a  pseudo-fusion rule on $\mfc(\mbF_p)$.
\epr
\begin{proof}
The map ${^p}N_{\mfg, 0}$ satisfies  (ii)' in Definition \ref{c34} because of the assumed condition $(**)$ (cf. Remark \ref{QH123}).
To verify that (iii)' is satisfied, let us consider
 the clutching data $\msG_\mr{tree}$
 defined in (\ref{QH1009}) with $g_j=0$ and $r_j \geq 2$ ($j=1,2$).
Also, we shall take an element $\rho$ as in (\ref{QH1033}).
As mentioned there,  each set   of radii $\rho_{\msG_\mr{tree}}$ for $\msG_\mr{tree}$ over $\mbF_p$ with $\rho_{\msG_\mr{tree}\Rightarrow \emptyset} = \rho$ can be described as 
$\rho_\msG := \{ \rho^j_\lambda \}_{j=1}^2$ for some  $\lambda \in \mfc (\mbF_p)$.
According to   Proposition \ref{c43}, the product $\prod_{j=1}^2 \mfO \mfp_{\mfg, \rho^j_\lambda, 0, r_j+1}^\ZZZ$ is generically \'{e}tale over $\prod_{j=1}^2 \overline{\mfM}_{0,r_j+1}$ because of the condition $(*)$.
Hence, since the square diagram   (\ref{QH1004}) for $\hslash =1$ is cartesian by Theorem \ref{y0176}, 
$ \mfO \mfp_{\mfg, \rho, 0, r}^\ZZZ/\overline{\mfM}_{0,r}$ is unramified at the points that lie in  the image of generic points of $\coprod_{\lambda \in \mfc (\mbF_p)} \prod_{j=1}^2 \mfO \mfp_{\mfg, \rho^j_\lambda, 0, r_j+1}^\ZZZ$ via the upper horizontal arrow in (\ref{QH1004}).
It follows from Corollary \ref{y0170} that $ \mfO \mfp_{\mfg, \rho, 0, r}^\ZZZ/\overline{\mfM}_{0,r}$ is also  \'{e}tale at the same points.
Moreover, 
 since $\mfO \mfp_{\mfg, \rho, 0, r}^\ZZZ/\overline{\mfM}_{0, r}$ is finite, each generic point of $ \mfO \mfp_{\mfg, \rho, 0, r}^\ZZZ$ specializes to one of these points.
This implies the second equality in the following sequence:
\begin{align}
{^p}N_{\mfg, 0}(\overline{\rho}) &= \mr{deg}(\mfO \mfp_{\mfg, \rho, 0, r}^\ZZZ/\overline{\mfM}_{0, r}) \\
& =
\mr{deg}(\coprod_{\lambda \in \mfc (\mbF_p)} \prod_{j=1}^2 \mfO \mfp_{\mfg, \rho^j_\lambda, 0, r_j+1}^\ZZZ/\prod_{j=1}^2 \overline{\mfM}_{0, r_j+1}) \notag \\
& = \sum_{\lambda \in \mfc (\mbF_p)} \mr{deg}(\prod_{j=1}^2\mfO \mfp_{\mfg, \rho^j_\lambda, 0, r_j+1}^\ZZZ/\prod_{j=1}^2 \overline{\mfM}_{0, r_j+1}) \notag \\
& =  \sum_{\lambda \in \mfc (\mbF_p)}  \mr{deg} (\mfO \mfp_{\mfg, \rho^1_\lambda, 0, r_1+1}^\ZZZ/\overline{\mfM}_{0, r_1+1}) \cdot \mr{deg}(\mfO \mfp_{\mfg, \rho^2_\lambda, 0, r_2+1}^\ZZZ/\overline{\mfM}_{0, r_2+1})\notag \\
& = \sum_{\lambda \in \mfc (\mbF_p)} {^p}N_{\mfg, 0}(\overline{\rho}^1 + \lambda) \cdot  {^p}N_{\mfg, 0}(\overline{\rho}^2 + \lambda^\dual). \notag
\end{align}
It follows  that ${^p}N_{\mfg, 0}$ satisfies  (iii)', and hence,  we finish the proof.
 \end{proof}

By applying 
 the discussion in \S\,\ref{y0179} to ${^p}N_{\mfg, 0}$ (for  $\mbG$ satisfying  both $(*)$ and $(**)$), one obtains  the pseudo-fusion ring 
\begin{equation} \label{1054}
 {^p}\Fus_{\mfg}  \left(:= \Fus_{{^p}N_{\mfg, 0}}\right)
 \end{equation}
associated to  the pseudo-fusion rule
${^p}N_{\mfg, 0}$.

\bde \label{c4} 
We shall refer to the ring ${^p}\Fus_{\mfg}$
 as the {\bf  pseudo-fusion ring  for   dormant $\mfg$-opers}.
  \index{${^p}\Fus_{\mfg}$, pseudo-fusion ring for dormant $\mfg$-opers}
\ede

\bpr \label{f176} 
With the above notation, we have
\begin{equation}\label{f170}
[0] \in \mr{Ker}({^p}N_{\mfg, 0}),  \hspace{5mm}  \mr{Ker}({^p}N_{\mfg, 0}) \neq \mfc (\mbF_p)
\end{equation}
 (cf. Remark \ref{f141} for the definition of $\mr{Ker}(-)$).
 Moreover,  the reduced ring 
 $({^p}\Fus_{\mfg})_{\mr{red}}$
 $\left(:= \Fus_{{^p}N_{\mfg, 0}\mr{red}}\right)$
  associated to ${^p}\Fus_{\mfg}$   satisfies the inequalities
 \begin{equation} \label{f171}
 1 < \mr{dim}_\mbQ (\mbQ\otimes_\mbZ({^p}\Fus_{\mfg})_{\mr{red}}) \leq  p^{\mr{rk}(\mfg)}.
 \end{equation}
\epr
\begin{proof}
As proved in Remark \ref{QG37},  
$\mfO\mfp^\ZZZ_{\mfg,  \rho, 0, r}$  is empty unless $\rho \in (\mfc (\mbF_p) \setminus \{[0] \})^{\times r}$.
That is to say, the kernel $\mr{Ker}({^p}N_{\mfg, 0})$ of ${^p}N_{\mfg, 0}$ contains the element $[0]$. 
In particular,   the equality $[0]^2 =0$ holds  in ${^p}\Fus_{\mfg}$.
This implies that 
  the $\mbQ$-subvector space $\mbQ \cdot  [0]$ of  $\mbQ \otimes_\mbZ {^p}\Fus_{\mfg}$ generated by $[0]$ is contained in $\mbQ \otimes_\mbZ\sqrt{(0)}$.
 Hence,  we have
\begin{align}
\mr{dim}_\mbQ ( \mbQ \otimes_\mbZ ({^p}\Fus_{\mfg})_{\mr{red}} )  &
\leq \mr{dim}_\mbQ(\mbQ \otimes_\mbZ {^p}\Fus_{\mfg}) -\mr{dim}_\mbQ (\mbQ \cdot  [0]) \\
 & = \mr{dim}_\mbQ (\mbQ  \cdot \varepsilon \oplus (\mbQ \otimes_\mbZ \Fus'_{{^p}N_{\mfg, 0}})) -1  \notag \\
 &  = 1 + \sharp (\mfc (\mbF_p)) -1 \notag \\
 & = p^{\mr{rk}(\mfg)}. \notag
\end{align}
We complete  the proofs of the first assertion in (\ref{f170}) and the second inequality in (\ref{f171}).

To complete the proof, it suffices to consider  the first inequality in (\ref{f171})
because the second assertion in (\ref{f170}) is a  direct consequence of that inequality.
 By Proposition \ref{f9}, (i), we  have 
 \begin{equation}
 \mbR \otimes_\mbZ \sqrt{(0)} = \big\{ y \in \mbR \otimes_\mbZ \Fus'_{{^p}N_{\mfg, 0}} \, \big| \, y^2 =0\big\}.
 \end{equation}
Since   $\alpha^2 = \sum_{\lambda \in \mfc (\mbF_p)} {^p}N_{\mfg, 0} (2 \alpha + \lambda^\dual) \cdot  \lambda$  for any $\alpha \in \mfc (\mbF_p)$,
the problem is reduced to finding 
   elements $\alpha$, $\beta \in \mfc (\mbF_p)$ with ${^p}N_{\mfg, 0} (2 \alpha + \beta) \neq 0$.
For arbitrary elements $\alpha^\odot$, $\beta^\odot \in \mfc_{\mfg^\odot} (\mbF_p)$,
the restriction  $\mfO \mfp^\ZZZ_{\mfg^\odot, \rho, 0,3} \migiincl \mfO \mfp^\ZZZ_{\mfg,  \iota_{\mfc*}(\rho), 0,3}$  of the closed immersion  (\ref{QS4567})
 with $\rho = (\alpha^\odot, \alpha^\odot, \beta^\odot)$ 
 yields  the inequality
\begin{equation}
{^p}N_{\mfg^\odot, 0} (2  \alpha^\odot + \beta^\odot) \leq {^p}N_{\mfg, 0} (2   \iota_{\mfc*}(\alpha^\odot) +  \iota_{\mfc*}(\beta^\odot)).
\end{equation}
Thus, we may assume that $\mfg = \mfg^\odot$.
But, 
this case may be easily verified from 
a combinatorial description  of dormant $\mfg^\odot$-opers on the $3$-pointed projective line, asserted in 
~\cite{Mzk2}, Introduction, \S\,1.2, Theorem 1.3 (or the discussion in  \S\,\ref{y0772} below).
This completes the proof of the remaining portion.
\end{proof}

In addition to ${^p}N_{\mfg,  0} $, let us consider   a collection of  maps 
\index{${^p}N_{\mfg,  0}$, ${^p}N_{\mfg,  g}$}
\begin{equation} \label{f180}
  {^p}N_{\mfg, g} : \mbN^{ \mfc (\mbF_p)}_{\geq 3-2g} \migi \mbZ
  \end{equation}
($g = 1, 2, \cdots $) given by  $\overline{\rho} \mapsto   {^p}N_{\mfg,  g} (\overline{\rho})$  (cf. (\ref{d6})) for any $\overline{\rho} \in \mbN^{{\mfc (\mbF_p)}}_{\geq 3-2g}$.
By applying  Theorem \ref{y0176} to  clutching data whose underlying semi-graph has   precisely  one vertex and one closed edge (which forms a loop), 
we have 
\begin{equation} \label{f181}
  {^p}N_{\mfg,  g} (x) = \sum_{\lambda \in {\mfc}(\mbF_p)} {^p}N_{\mfg,  g-1} (x + \lambda + \lambda^\dual)
  \end{equation}
for  any $g \geq 1$ and   $x \in \mbN_{\geq 3-2g}^{\mfc(\mbF_p)}$.
Thus,    the following theorem holds by  Proposition \ref{c35}.

\bt[cf. Theorem \ref{y018}] 
\label{c6}
Suppose that $\mbG$ satisfies the conditions $(*)$ and $(**)$.
Then, the following assertions hold:
\begin{itemize}
\item[(i)]
Let 
 $g_1$, $g_2$ be  nonnegative integers and let  $x \in \mbN^{\mfc (\mbF_p)}_{\geq 3-2g_1}$, $y \in \mbN^{\mfc (\mbF_p)}_{\geq 3 - 2 g_2}$.
 Then,  the collection of functions $\{  {^p}N_{\mfg, g}\}_{g \geq 0}$  satisfies the following rule:
\begin{equation}
{^p}N_{\mfg, g_1 +g_2} (x + y) 
= \sum_{\lambda \in \mfc (\mbF_p)} {^p}N_{\mfg, g_1} (x + \lambda) \cdot {^p}N_{\mfg, g_2} (y + \lambda^\dual).
\end{equation}
\item[(ii)]
Let $(g, r)$ be a pair of nonnegative integers with  $2g-2+r >0$ and let $\rho := (\rho_i)_{i=1}^r \in \mfc^{\times r} (\mbF_p)$.
Write ${^p}\mfS_{\mfg} \left(:= \mr{Hom} ({^p}\Fus_{\mfg}, \mbC)\right)$ for the set of ring homomorphisms ${^p}\Fus_{\mfg} \migi \mbC$
  and write ${^p}\mr{Cas}_\mfg := \sum_{\lambda \in \mfc (\mbF_p)} \lambda \ast \lambda^\dual \left(\in {^p}\Fus_{\mfg}\right)$.
Then,  the following equality holds:
\begin{align}   \label{1060}
{^p}N_{\mfg, g} (\sum_{i=1}^r \rho_i) 
  \left( =  \mr{deg}(\mfO\mfp^\ZZZ_{\mfg,   \rho,  g, r}/\overline{\mfM}_{g,r})
\right)
 = \sum_{\chi \in {^p}\mfS_{\mfg}} \chi ({^p}\mr{Cas}_\mfg)^{g-1}\cdot \prod_{i=1}^r \chi(\rho_i).  
  \end{align}
In particular, if  $r=0$ (and  $g>1$), then this equality reads
\begin{align}   \label{10ddfg60}
{^p}N_{\mfg, g} (0)
  \left(= \mr{deg}(\mfO\mfp^\ZZZ_{\mfg,   g, 0}/\overline{\mfM}_{g,0})
 \right)
 = \sum_{\chi \in {^p}\mfS_{\mfg}} \chi ({^p}\mr{Cas}_\mfg)^{g-1}.    \end{align}
\end{itemize}
\et

Also, by Proposition \ref{f200}, the structure constants of  ${^p}\Fus_{\mfg}$ satisfy the following Corollary \ref{c3445}.
Here, we shall write  $\chi_0 : {^p}\Fus_{\mfg}  \migi \mbC$ for  the composite of the natural augmentation map ${^p}\Fus_{\mfg}  \migi \mbZ$ (cf. (\ref{f123})) and the inclusion $\mbZ \migiincl \mbC$.
 By abuse of notation, we also write $\chi_0$ for the corresponding 
 $\mbC$-rational point of  $\mr{Spec}(\mbC \otimes_\mbZ {^p}\Fus_{\mfg})$.
If   $\mfs$ is a basis  for the $\mbC$-vector space $\mcT_{\chi_0} \mr{Spec}(\mbC \otimes_\mbZ {^p}\Fus_{\mfg}) \left(= (\mbC \otimes_\mbZ \sqrt{(0)})^\vee \right)$, then
we
 set
\begin{equation}
{^p}\widetilde{\mfS}_{\mfg, \mfs} :=  \mfs \sqcup ({^p}\mfS_{\mfg} \setminus \{ \chi_0\}) \ \ \text{and} \ \ {^p}\mr{Cas}_{\mfg, \chi} := \sum_{\lambda \in \mfc (\mbF_p)} |\chi (\lambda)|^2 \ \ \left(\chi \in {^p}\widetilde{\mfS}_{\mfg, \mfs} \right).
\end{equation}

\bco \label{c3445}
Suppose that  $\mbG$ satisfies the conditions $(*)$ and $(**)$. 
Then, there exists a basis $\mfs$ for  $\mcT_{\chi_0} \mr{Spec}(\mbC \otimes_\mbZ {^p}\Fus_{\mfg})$
such that  the matrix
\begin{equation}
\Sigma_{\mfs} := \left(  \frac{\chi (\lambda)}{\sqrt{{^p}\mr{C as}_{\mfg, \chi}}}\right)_{(\chi, \lambda) \in {^p}\widetilde{\mfS}_{\mfg, \mfs}  \times \mfc (\mbF_p)}
\end{equation}
is unitary  and diagonalizes simultaneously all matrices
\begin{equation}
{^p}N_{\mfg, \alpha} := \left(\mr{deg}(\mfO\mfp^\ZZZ_{\mfg,    (\alpha, (-1)\star \beta, \gamma),  0, 3}/k) \right)_{(\beta, \gamma)  \in  \mfc (\mbF_p) \times \mfc (\mbF_p)} \ \ \left(\alpha \in \mfc (\mbF_p) \right).
\end{equation}
\eco

\section{The pseudo-fusion rule for $\mfs \mfl_2$}\label{y0772}

At the time of writing this manuscript, 
the author does not know much about the pseudo-fusion ring ${^p}\Fus_\mfg$ and the ring homomorphisms   ${^p}\Fus_\mfg \migi \mbC$ for a general $\mfg$.
However, the case of $\mfg = \mfg^\odot \left(= \mfs \mfl_2 \right)$
can be concretely understood by applying related previous studies and comparing with another type of fusion rule.
We shall end this section with  discussing this topic.
Note that the involution
of $\mfc_{\mfg^\odot}$ given by $\lambda \mapsto (-1) \star \lambda$ 
 coincides with the identity map.

\subsection{The explicit description of ${^p}\Fus_{\mfg^\odot}$}
We first give an explicit   description of ${^p}N_{\mfg^\odot, 0}$ following S. Mochizuki.
Let $\overline{k}$ be an algebraic closure of $k$, and denote by  $\overline{k}/ \{ \pm 1\}$ the set of equivalence classes of elements $\lambda \in \overline{k}$, in which $\lambda$ and $- \lambda$ are identified.
By letting    $ \sqrt{\mbF_p}:= \{ a \in \overline{k} \, | \, a^2 \in \mbF_p \} \left(\subseteq \overline{k}\right)$,  we  obtain 
a subset $ \sqrt{\mbF_p} / \{ \pm 1\}$ of $\overline{k}/ \{ \pm 1\}$.
The composite
\begin{equation} \label{QH150}
\overline{k} \migi \mfg^\odot (\overline{k}) \xrightarrow{\chi (\overline{k})} \mfc_{\mfg^\odot} (\overline{k})
\end{equation}
(cf. (\ref{chi}) for the definition of $\chi$),  where the first arrow is given by  assigning $a \mapsto \begin{pmatrix} a & 0 \\ 0 & -a\end{pmatrix}$,  factors through the natural quotient  $\overline{k} \migisurj \overline{k}/ \{ \pm 1\}$.
The resulting map  
$\overline{k}/ \{ \pm 1\} \migi  \mfc_{\mfg^\odot} (\overline{k})$
 is bijective and restricts to a bijection 
 \begin{align} \label{QH801}
 \sqrt{\mbF_p} / \{ \pm 1\} \isom  \mfc_{\mfg^\odot} (\mbF_p).
 \end{align}
 Using  this bijection, we shall identify $\mfc_{\mfg^\odot} (\mbF_p)$ with $ \sqrt{\mbF_p}/ \{ \pm 1\}$.
 Next, let us define a subset $V$ of $\mbZ^{\times 3}$ to be
  \begin{equation} \label{QH802}
 V :=
 \left\{ (a_1, a_2, a_3) \in \mbZ^{\times 3} \,  
  \Biggm| 
  \text{  $\sum_{i=1}^3 a_i$ is odd $< 2 p$,  and } \   \begin{matrix}
 0 < a_1 < a_2 +a_3,  \\
 0 < a_2 < a_3 +a_1, \\ 
  0 < a_3 < a_1 + a_2.
 \end{matrix} \right\}.
 \end{equation}
Also,  define an equivalence relation $\sim$ on $V$ given as follows:
\begin{equation}
(a_1, a_2, a_3) \sim (a'_1, a'_2, a'_3) \stackrel{\mr{def}}{\Longleftrightarrow}  a_i^2 \equiv (a'_i)^2   \ (\mr{mod} \  p)  \ \text{for any $i \in \{ 1, 2,3 \}$}.
\end{equation}
The  map of sets 
\begin{align} \label{QH791}
w :V/\sim  \ \migi ( \sqrt{\mbF_p} / \{ \pm 1\})^{\times 3}; \overline{(a_1, a_2, a_3)} \mapsto \left(\frac{\overline{a_1}}{2}, \frac{\overline{a_2}}{2}, \frac{\overline{a_3}}{2}\right)
\end{align}
 is verified to be  injective, where we recall that $p$ is odd.
 According to ~\cite{Mzk2}, Introduction, \S\,1.2, Theorem 1.3, the structure constants $N^{\beta}_{\alpha, \gamma} := {^p}N_{\mfg^\odot, 0} (\alpha + \beta + \gamma)$ (for  $(\alpha, \beta, \gamma) \in ( \sqrt{\mbF_p} / \{ \pm 1\})^{\times 3}$) of the pseudo-fusion ring ${^p}\Fus_{\mfg^\odot}$ are given as follows:
 \begin{equation}
N^{\beta}_{\alpha, \gamma} = \begin{cases} 1 & \text{if $(\alpha, \beta, \gamma) \in \mr{Im}(w)$,} \\
0 & \text{if $(\alpha, \beta, \gamma) \notin \mr{Im}(w)$.}
\end{cases}
\end{equation}
This result is consistent with the result in \S\,\ref{QG1901} via  a suitable correspondence between $V$ and  the  set of exponent differences of hypergeometric differential operators.

In what follows,  we consider  a slightly simpler description of  the subset  $\mr{Im}(w) \subseteq ( \sqrt{\mbF_p} / \{ \pm 1\})^{\times 3}$ and the pseudo-fusion rule for dormant $\mfg^\odot$-opers.
Let us write
\begin{equation} \label{QH188}
P_\mbZ := \left\{ a \in \mbZ \, \bigg| \, 0 \leq a \leq \frac{p-2}{2} \right\}
\end{equation}
 and 
set   $W_\mbZ$ to be the subset of $P_\mbZ^{\times 3}$  defined as  
 \begin{equation} \label{QH790}
 W_\mbZ := \left\{ (b_1, b_2, b_3) \in P_\mbZ^{\times 3} \,   \Biggm| 
  \text{  $b_1+ b_2+b_3 \leq p-2$   and } \   \begin{matrix}
  b_1 \leq b_2 + b_3,  \\
  b_2 \leq b_3+b_1,  \\ 
  b_3\leq b_1 +b_2.
 \end{matrix}
   \right\}.
 \end{equation}
Then, the map of sets 
\begin{align} \label{QH789}
w_P: P_\mbZ \migi \sqrt{\mbF_p} / \{ \pm 1\}; b \mapsto \frac{\overline{2b +1}}{2} \  \left(= \overline{b + \frac{p+1}{2}}\right) 
 \end{align}
  is injective; we shall regard $P_\mbZ$ as a subset of $\sqrt{\mbF_p} / \{ \pm 1\}$ by this injection.
Moreover,   the resulting inclusion
  $w_P^{\times 3} : P_\mbZ^{\times 3} \migiincl (\sqrt{\mbF_p} / \{ \pm 1\})^{\times 3}$
gives  an identification  $W_\mbZ = \mr{Im}(w)$.
The restriction $N_{P_\mbZ}$ of ${^p}N_{\mfg^\odot, 0} : \mbN^{\sqrt{\mbF_p} / \{ \pm 1\}}_{\geq 3} \migi \mbZ$ to $\mbN^{P_\mbZ}$ forms a pseudo-fusion rule on $P_\mbZ$ (with  trivial involution), and hence, yields  the associated pseudo-fusion ring $\Fus_{N'}$.
The above discussion shows that the structure constants $N_{P_\mbZ} (a + b+ c)$ (for $(a, b, c) \in P_\mbZ^{\times 3}$) 
of $\Fus_{P_\mbZ} := \Fus_{N_{P_\mbZ}}$
 is given as follows:
 \begin{equation} \label{w010w}
N_{P_\mbZ} (a + b+ c) := \begin{cases} 1 & \text{if $(a,b,c) \in W_\mbZ$,} \\
0 & \text{if $(a,b,c) \notin W_\mbZ$.}
\end{cases}
\end{equation}
According to the discussion in  Remark \ref{f141}, 
   the projection  ${^p}\Fus_{\mfg^\odot} \migisurj \Fus_{P_\mbZ}$ induces  an isomorphism  of rings
 $({^p}\Fus_{\mfg^\odot})_{\mr{red}} \isom (\Fus_{P_\mbZ})_{\mr{red}}$.
In summary, we have proved  the following proposition.

\bpr \label{c3445ee}
The ring $({^p}\Fus_{\mfg^\odot})_{\mr{red}}$  is isomorphic to the reduced ring $(\Fus_{P_\mbZ})_{\mr{red}}$ associated to the ring $\Fus_{P_\mbZ}$ defined as follows:
\begin{itemize}
\item
The underlying abelian group of $\Fus_{P_\mbZ}$ is the free abelian group  $\mbZ \cdot \varepsilon \oplus \bigoplus_{a \in P_\mbZ}  \mbZ \cdot e_a$
with  basis $\{ \varepsilon \} \sqcup \{ e_a \}_{a \in P_\mbZ}$.
\item
The multiplication law ``\,$\ast$'' on $\Fus_{P_\mbZ}$ is defined by putting
$\varepsilon \ast e_a = e_a \ast \varepsilon := e_a$, 
$e_a \ast  e_b := \sum_{c \in P_\mbZ} N_{P_\mbZ} (a+ b+c) \cdot e_c$ 
for any $a$, $b$, $c \in P_\mbZ$ and extending by bilinearity,
where $N_{P_\mbZ} (a+ b+c)$ denotes the integer defined by (\ref{w010w}).
\end{itemize}
\epr
\begin{rema} \label{ee111}
Since $(0, 0, 0) \in W_\mbZ$, 
   the moduli stack $\mfO \mfp^{^\mr{Zzz...}}_{\mfg^\odot, 0, 3}$, as well as $\mfO \mfp^{^\mr{Zzz...}}_{\mfg^\odot, \hslash, 0, 3}$ for $\hslash \in k^\times$ (cf. (\ref{QS1100})),  is nonempty.
 Hence, by applying Theorem \ref{y0176} to suitable clutching data, one may  conclude that  {\it for each  pair of nonnegative integers $(g,r)$ with $2g-2+r >0$ and each $\hslash \in k^\times$, the moduli stack     $\mfO \mfp^{^\mr{Zzz...}}_{\mfg^\odot, \hslash, g, r}$ is nonempty.}
\end{rema}



\begin{exa}[The case of $p=7$]
In this remark, we perform  a few computations to indicate how the pseudo-fusion rules for dormant $\mfg^\odot$-opers work for low genus.
We shall deal with the case of  $p =7$ and, for simplicity, we write $N_g := {^7}N_{\mfg^\odot, g}$.
For each $a \in P_\mbZ$,  $[a]$ denotes  the element  $w_{P} (a) \in \sqrt{\mbF_7}/ \{ \pm 1\} \left(\subseteq \mbN^{\sqrt{\mbF_7}/ \{ \pm 1\}}\right)$.
First, recall from  the discussion preceding Proposition \ref{c3445ee}  that  
 $N_0 ([a] + [b] + [c]) = 1$ (resp., $N_0 ([a] + [b] + [c]) = 0$) if and only if $(a, b,c) \in W_\mbZ$ (resp.,  $(a, b,c) \notin W_\mbZ$).
That is to say,  the set of  isomorphism classes of dormant $\mfg^\odot$-opers on $\msP$ (cf. (\ref{1051}))
  is in bijection with the set 
\begin{align}
W = & \ \big\{ (0, 0, 0), (0, 1,1), (1, 0,1), (1,1,0), (1,1,1), (0, 2,2), (2, 0,2), \\
& \ \ \  (2,2,0),
 (1,1,2), (1,2,1), (2,1,1),(1,2,2), (2,1,2), (2,2,1)\big\}.\notag
\end{align}
Then, by  (\ref{f181}) and Theorem \ref{c6}, (i), we have the following computations:
\vspace{2mm}
\begin{itemize}
\item
$(g,r) = (1, 1)$:
\begin{align}
N_1 ([0]) &= {\sum}_{i =0}^2 N_0 ([0] + 2[i]) =  \sharp \big\{ (0, 0, 0), (0, 1, 1), (0, 2,2) \big\} = 3;\\
N_1 ([1]) &= {\sum}_{i =0}^2 N_0 ([1] + 2[i])  =  \sharp \big\{ (1, 1,1), (1, 2,2)\big\} =2; \notag \\
N_1 ([2]) &=  {\sum}_{i =0}^2 N_0 ([2] + 2[i])  = \sharp \big\{ (2, 1,1)\big\} =1. \notag  
\end{align}
\vspace{1mm}
\item
$(g,r) = (1,2)$:
\begin{align}
N_1 (2[0]) &=  {\sum}_{i=0}^2 N_1 ([i]) \cdot   N_0 ([i] + 2 [0]) = 3 \cdot 1+ 2 \cdot 0 + 1 \cdot 0 =3; \\
N_1 ([0]+ [1]) &=  {\sum}_{i=0}^2 N_1 ([i]) \cdot   N_0 ([i] + [0] + [1]) = 3 \cdot 0 + 2 \cdot 1 + 1 \cdot 0 = 2; \notag \\
N_1 (2[1]) &=  {\sum}_{i=0}^2 N_1 ([i]) \cdot   N_0 ([i] + 2 [1]) = 3 \cdot 1 + 2 \cdot 1 + 1 \cdot 1 =6; \notag \\
N_1 ([0]+ [2]) &=  {\sum}_{i=0}^2 N_1 ([i]) \cdot   N_0 ([i] + [0] + [2]) = 3 \cdot 0 + 2 \cdot 0 + 1 \cdot 1 = 1; \notag \\
N_1 ([1]+ [2]) &=  {\sum}_{i=0}^2 N_1 ([i]) \cdot   N_0 ([i] + [1] + [2]) = 3\cdot 0 + 2 \cdot 1 + 1\cdot 1 =  3; \notag \\
N_1 (2[2]) &=  {\sum}_{i=0}^2 N_1 ([i]) \cdot   N_0 ([i] + 2 [2]) =  3 \cdot 1 + 2 \cdot 1 =5. \notag
\end{align}
\vspace{1mm}
\item
$(g,r) = (2,0)$:
\begin{align}
N_2 (0) = {\sum}_{i=0}^2 N_1 ([i]) \cdot N_1 ([i]) = 3 \cdot 3 + 2 \cdot 2 + 1 \cdot 1 = 14.
\end{align}
\vspace{1mm}
\item
$(g,r) = (2, 1)$:
\begin{align}
N_2 ([0]) &= {\sum}_{i=0}^2 N_1 ([i]) \cdot N_1 ([i] + [0]) = 3 \cdot 3 + 2 \cdot 2 + 1 \cdot 1 = 14; \\
N_2 ([1]) &= {\sum}_{i=0}^2 N_1 ([i]) \cdot N_1 ([i] + [1]) = 3 \cdot 2 + 2 \cdot 6 + 1 \cdot 3 = 21; \notag  \\
N_2 ([2]) &= {\sum}_{i=0}^2 N_1 ([i]) \cdot N_1 ([i] + [2]) = 3 \cdot 1 + 2 \cdot 3 + 1 \cdot 5 = 14. \notag
\end{align}
\vspace{1mm}
\item
$(g,r) = (3, 0)$:
\begin{align} \label{fff078}
N_3 (0) = {\sum}_{i=0}^2 N_2 ([i]) \cdot N_1 ([i]) = 14 \cdot 3 + 21 \cdot 2 + 14 \cdot 1 = 98.
\end{align}
\end{itemize}
Here, $N_2 (0)$ may  be  calculated also by $N_2 (0) = \sum_{0 \leq a, b,c \leq 2}N_0 ([a] + [b] + [c])^2 = \sharp (W_\mbZ) = 14$.
The equality  $N_2 (0) = 14$ 
 obtained here  is  consistent with    (\ref{n=2}) in Introduction.
\end{exa}

\subsection{The relationship with the $\mfs \mfl_2 (\mbC)$-WZW model} \label{QH541}

Next, let us study the relationship with
 the fusion ring of the $\mfs \mfl_2 (\mbC)$-WZW (= Wess-Zumino-Witten) conformal field theory at  level $p-2$.
 By means of this relationship, we describe (cf. Theorem  \ref{TT11} below) explicitly the various values of ${^p}N_{\mfg^\odot}$.

   Denote  by $(-, -)$ the Killing form on $\mfs \mfl_2 (\mbC)$ (= the Lie algebra $\mfg^\odot$ over $\mbC$) normalized as $(\alpha^\odot, \alpha^\odot) = 2$, where we regard the unique   root $\alpha^\odot$  (cf. \S,\ref{QR40}) as  an element of $(\mft^\odot)^\vee$  via differentiation.
Let us set
\begin{equation} \label{QH700}
P := \left\{ a \in \frac{1}{2}\mbZ \, \bigg| \, 0 \leq a \leq \frac{p-2}{2} \right\} \left(\subseteq \frac{1}{2} \mbZ := \left\{0, \pm \frac{1}{2}, \pm 1, \pm \frac{3}{2}, \pm 2, \cdots  \right\} \right).
\end{equation}
In particular, we have  $P \cap \mbZ = P_\mbZ$.
The set of dominant weights $\lambda$ of $\mfs \mfl_2 (\mbC)$ with $0 \leq (\alpha^\odot, \lambda) \leq p-2$
coincides with 
$\left\{ 0, \frac{1}{2} \alpha^\odot, \alpha^\odot, \cdots, \frac{p-2}{2} \alpha^\odot \right\}$.
By  the correspondence $j \alpha^\odot \leftrightarrow j$, 
we shall identify this set with $P$.
Next,
we shall 
set
 \begin{equation}
 W := \left\{ (b_1, b_2, b_3) \in P^{\times 3} \,   \Biggm| 
  \text{  $b_1+ b_2+b_3 \leq p-2$   and } \   \begin{matrix}
  b_1 \leq b_2 + b_3,  \\
  b_2 \leq b_3+b_1,  \\ 
  b_3\leq b_1 +b_2.
 \end{matrix}
   \right\}
 \end{equation}
(hence $W \cap \mbZ = W_\mbZ$).
Then, it is known that 
the fusion ring $\Fus_{P}^\circledcirc$ associated to the $\mfs \mfl_2 (\mbC)$-WZW conformal field theory  at level $p-2$ (i.e., ``$\mcR_{p-2}(\mfs \mfl_2)$'' in the notation of ~\cite{Beau}) is defined as follows:
\begin{itemize}
\item
The underlying abelian group of $\Fus_{P}^\circledcirc$ is the  free abelian group $\bigoplus_{a \in P}\mbZ \cdot e_a$ with basis $\{ e_a \}_{a \in P}$;
\item
The multiplication law ``$\ast$'' on $\Fus_{P}^\circledcirc$  is given  by putting $e_a \ast e_b = \sum_{c \in P} N_P (a + b+c) \cdot e_c$ for any $a, b, c \in P$ and extending bilinearity, where $N_P (a + b + c)$ denotes the integer defined by
 \begin{equation} \label{wghki}
N_{P} (a + b+ c) := \begin{cases} 1 & \text{if $(a,b,c) \in W$,} \\
0 & \text{if $(a,b,c) \notin W$.}
\end{cases}
\end{equation}
\end{itemize}
One verifies that
there exists a ring homomorphism 
\begin{align} \label{QH744}
\varkappa : \Fus_{P_\mbZ} \migi \Fus_{P}^\circledcirc
\end{align}
 given by $\varepsilon \mapsto e_0$ and $e_a \mapsto e_a$ for any $a \in P_\mbZ$.
The kernel of $\varkappa$ is the square nilpotent  ideal generated by the  element $\varepsilon - e_0$, and  if we consider $\Fus_{P}^\circledcirc$ as 
 a $\Fus_{P_\mbZ}$-module  by using $\varkappa$, then it is finitely generated.
 It follows that the map of sets $\varkappa^* : \mr{Spec}(\Fus_{P}^\circledcirc)(\mbC) \migi \mr{Spec}(\Fus_{P_\mbZ}) (\mbC)$ induced by $\varkappa$  is surjective.

Next, let us describe the ring homomorphisms $\Fus_{P_\mbZ} \migi \mbC$ of $\Fus_{P_\mbZ}$.
 By  ~\cite{Beau}, Lemma 9.3, 
 the set of ring homomorphisms $\Fus_{P}^\circledcirc \migi \mbC$ coincides with the set $\{ \mr{Tr}_* (j) \}_{j=1}^{p-1}$, where each $\mr{Tr}_* (j)$ is given by 
 \begin{align} \label{ee772}
\mr{Tr}_* (j)  : e_a \mapsto  \frac{\sin \left(\frac{(2a+1)j \pi}{p}\right)}{\sin \left(\frac{j \pi}{p} \right)} \hspace{5mm} \left(a \in P \right).
 \end{align}
The equality  $\mr{Tr}_* (j_1) \circ \varkappa= \mr{Tr}_* (j_2) \circ \varkappa$ holds  if and only if $j_1 + j_2 = p$.
   It follows that the set $\mfS_{P_\mbZ}$ 
  of ring homomorphisms $\Fus_{P_\mbZ} \migi \mbC$  coincides with $\{ \mr{Tr}_* (j) \circ \varkappa\}_{j =1, 2, \cdots, \frac{p-1}{2}}$.
For each $j \in \{1, 2, \cdots, \frac{p-1}{2} \}$ and $a \in P$, the equality 
  $| \mr{Tr}_* (j) (e_a) | = | \mr{Tr}_* (j) (e_{\frac{p-2}{2} -a}) |$ holds.
  Hence, 
  we have 
  \begin{align} \label{ee770}
   (\mr{Tr}_* (j) \circ \varkappa) (\mr{Cas}_{N_{P_\mbZ}}) 
  & = \sum_{a \in P_\mbZ} | (\mr{Tr}_* (j) \circ \varkappa)  (e_{a}) |^2  \\
  & =   \sum_{a \in P_\mbZ} \frac{1}{2}  \cdot \left( | \mr{Tr}_* (j) (e_{a}) |^2 +  | \mr{Tr}_* (j) (e_{\frac{p-2}{2}-a}) |^2 \right) \notag \\
 & =   \frac{1}{2} \cdot  \sum_{a \in P} | \mr{Tr}_* (j) (e_a) |^2   \notag \\
  & =  \frac{p}{ \left(2 \sin \left( \frac{j \pi}{p}\right) \right)^2},   \notag
  \end{align}
where  the fourth equality follows from ~\cite{Beau}, Lemma 9.7.
 It follows that,  for each  pair of nonnegative integers $(g,r)$ with $2g-2+r >0$ and  each $(\rho_i)_{i=1}^r \in P_\mbZ^{\times r}$ (considered as an element  of $\mfc_{\mfg^\odot}(\mbF_p)$ via (\ref{QH801}) and (\ref{QH789}) together), the following sequence of equalities  holds:
  \begin{align} \label{QH808}
{^p}N_{\mfg^\odot, g} (\sum_{i=1}^r \rho_i) 
  &   \stackrel{(\ref{1060})}{=} \sum_{\chi \in {^p}\mfS_{\mfg^\odot}} \chi ({^p}\mr{C as}_{\mfg^\odot})^{g-1} \cdot \prod_{i=1}^r \chi (\rho_i)  \\
   &  \stackrel{}{=} \sum_{\chi \in \mfS_{P_\mbZ}} \chi (\mr{C as}_{N_{P_\mbZ}})^{g-1} \cdot \prod_{i=1}^r \chi (\rho_i) \notag  \\
 &  = \sum_{j\in \left\{1, \cdots, \frac{p-1}{2}\right\}} (\mr{Tr}_* (j) \circ \varkappa) (\mr{C as}_{N_{P_\mbZ}})^{g-1} \cdot \prod_{i=1}^r (\mr{Tr}_* (j) \circ \varkappa)(\rho_i) \notag  \\
  &   = \frac{p^{g-1}}{2^{2g-2}} \cdot  \sum_{j =1}^{\frac{p-1}{2}} \frac{\prod_{i=1}^r \sin \left( \frac{(2 \rho_i +1) j \pi}{p}\right)}{\sin^{2g-2+r} \left(\frac{j \pi}{p} \right)} \notag \\
  &  = \frac{p^{g-1}}{2^{2g-1}} \cdot \sum_{j =1}^{p-1} \frac{\prod_{i=1}^r \sin \left( \frac{(2 \rho_i +1) j \pi}{p}\right)}{\sin^{2g-2+r} \left(\frac{j \pi}{p} \right)},
   \notag
  \end{align}
 where the fourth equality follows from 
  (\ref{ee770}).
 Consequently, we obtain the following assertion, generalizing (\ref{n=2g}).

\vspace{3mm}
\bt \label{TT11} 
Let $(g,r)$ and $\rho$ be as above.
Then,
 the following equality holds:
 \begin{align}
 {^p}N_{\mfg^\odot, g}(\sum_{i=1}^r \rho_i)  \left(= \mr{deg} (\mfO \mfp^{^\mr{Zzz...}}_{\mfg^\odot, \rho, g,r}/\overline{\mfM}_{g,r}) \right)
 =  \frac{p^{g-1}}{2^{2g-1}} \cdot \sum_{j =1}^{p-1} \frac{\prod_{i=1}^r \sin \left( \frac{(2\rho_i +1) j \pi}{p}\right)}{\sin^{2g-2+r} \left(\frac{j \pi}{p} \right)}. \hspace{-3mm}
 \end{align}
 In particular, the generic degree $\mr{deg} (\mfO \mfp^{^\mr{Zzz...}}_{\mfg^\odot,  g,r}/\overline{\mfM}_{g,r})$ of  $\mfO \mfp^{^\mr{Zzz...}}_{\mfg^\odot,  g,r}/\overline{\mfM}_{g,r}$ is given by the following formula:
 \begin{align} \label{QH882}
 \mr{deg} (\mfO \mfp^{^\mr{Zzz...}}_{\mfg^\odot,  g,r}/\overline{\mfM}_{g,r}) =  \frac{p^{g-1}}{2^{2g-1+r}} \cdot \sum_{j =1}^{p-1} 
\frac{\left(1-(-1)^j \cdot \cos \left(\frac{j\pi}{p} \right)\right)^r}{\sin^{2(g-1+r)} \left(\frac{j \pi}{p} \right)}.
 \end{align}
  \et
\begin{proof}
The former assertion was already proved in the above discussion.
The letter assertion follows directly from the former one.
Indeed, we have
\begin{align}
& \ \ \ \ \mr{deg} (\mfO \mfp^{^\mr{Zzz...}}_{\mfg^\odot,  g,r}/\overline{\mfM}_{g,r}) \\
 &= \sum_{(\rho_i)_{i=1}^r \in P_\mbZ^{\times r}} {^p}N_{\mfg^\odot, g}(\sum_{i=1}^r \rho_i) 
\notag \\
& =  \sum_{(\rho_i)_{i=1}^r \in P_\mbZ^{\times r}} \frac{p^{g-1}}{2^{2g-1}} \cdot \sum_{j =1}^{p-1} \frac{\prod_{i=1}^r \sin \left( \frac{(2\rho_i +1) j \pi}{p}\right)}{\sin^{2g-2+r} \left(\frac{j \pi}{p} \right)} \notag \\
& =  \frac{p^{g-1}}{2^{2g-1}} \cdot \sum_{j =1}^{p-1} 
\frac{1}{\sin^{2g-2+r} \left(\frac{j \pi}{p} \right)}\cdot \left(\sum_{(\rho_i)_{i=1}^r \in P_\mbZ^{\times r}}\prod_{i=1}^r \sin \left( \frac{(2\rho_i +1) j \pi}{p}\right)\right) \notag \\
& =  \frac{p^{g-1}}{2^{2g-1}} \cdot \sum_{j =1}^{p-1} 
\frac{1}{\sin^{2g-2+r} \left(\frac{j \pi}{p} \right)} \cdot \left(\sum_{\rho\in P_\mbZ} \sin \left( \frac{(2\rho +1) j \pi}{p}\right)\right)^r \notag \\
& =  \frac{p^{g-1}}{2^{2g-1}} \cdot \sum_{j =1}^{p-1} 
\frac{1}{\sin^{2g-2+r} \left(\frac{j \pi}{p} \right)} \cdot
\left( \frac{1-(-1)^j \cdot \cos \left(\frac{j\pi}{p} \right)}{2 \sin \left(\frac{j \pi}{p}\right)}\right)^r \notag \\
& =  \frac{p^{g-1}}{2^{2g-1+r}} \cdot \sum_{j =1}^{p-1} 
\frac{\left(1-(-1)^j \cdot \cos \left(\frac{j\pi}{p} \right)\right)^r}{\sin^{2(g-1+r)} \left(\frac{j \pi}{p} \right)},  \notag
\end{align}
thus completing the proof of the assertion.
\end{proof}

\vspace{10mm}
\chapter{Generic \'{e}taleness of the moduli space of dormant opers
} \label{y0186}

This chapter is devoted to proving the generic \'{e}taleness of $\mfO \mfp^\ZZZ_{\mfs \mfl_n, \rho, g,r}$ over $\overline{\mfM}_{g,r}$.
By taking account of results proved so far, we see that  this result  deduces moreover  the generic \'{e}taleness of $\mfO \mfp^\ZZZ_{\mfg, \hslash,  \hslash \star \rho, g,r}$ for $\hslash \in k^\times$ and 
$\mfg = \mfs \mfl_n,  \mfs \mfo_{2l+1},  \mfs \mfp_{2m}$.
The simplest  case, i.e., $\mfg = \mfs \mfl_2$, was already  proved by S. Mochizuki in ~\cite{Mzk2}.
Just as in that case, we give a proof of this claim based on the idea that infinitesimal deformations of a dormant oper on the $3$-pointed projective line    may be determined completely by their radii and   the sheaves of their horizontal sections.

We here 
describe a rough sketch of the proof.
To prove the generic \'{e}taleness,
it suffices to verify 
that there are nontrivial deformations of 
 a fixed dormant $(\mr{GL}_n, \bb)$-oper $\nabla^\diamondsuit$ on a 
 genus-$0$ curve $X$  equipped with three marked points
  for an $n$-theta characteristic $\bb$.
 Let us denote by  ${^p}\mr{Def}^{\mr{Op}}_{\nabla^\diamondsuit}$ the deformation space of $\nabla^\diamondsuit$.
In the discussion of the present chapter, we will construct, for a suitable $\bb$,   
 the following  sequence of inclusions between certain  pointed sets:
\begin{align} \label{QQ699}
{^p}\mr{Def}^{\mr{Op}}_{\nabla^\diamondsuit}
 &\subseteq  \varpi (\Gamma (X, \JJ\DDD^{\, <n}_\BB))  \cap {^p}\mr{Def}^\sharp_{\nabla^\diamondsuit}  \\
&\subseteq
\varpi (\Gamma (X, \JJ\DDD^{\, <n}_\BB))  \cap \Upsilon^{-1}(0)  \notag  \\
&\subseteq \varpi (\Gamma (X, \JJ\DD^{\,<n}_{\BB} \cap \Omega^\sharp  \otimes\mcE nd (\DF_\BB^\sharp)))  \notag \\
&\subseteq \varpi (\Gamma (X, \JJ\DDD^{\,\sharp <n}_{\BB})).\notag
\end{align}
On the other hand, we will also prove  the equality $\Gamma (X, \JJ\DDD^{\, \sharp <n}_{\BB}) =0$,   which implies that  the rightmost of this sequence is a singleton.
By these facts,  
  ${^p}\mr{Def}^{\mr{Op}}_{\nabla^\diamondsuit}$  turns out to consist exactly  of one element corresponding to the trivial deformation, as desired.
 
\begin{notation}

Let  $p$ be a prime, $k$ 
\index{$k$, field} 
 a  field  of characteristic $p$, $n$ an integer with $1 < n < p$,  and $(g,r)$ a pair of nonnegative integers with $2g-2+r>0$.
Unless otherwise stated,  $k$ is assumed to be  algebraically closed .
Also, let 
$\msX := (X/k, \{ \sigma_i \}_{i=1}^r)$  an $r$-pointed ``smooth'' curve of genus $g$ over $k$.
For convenience, we write $\Omega := \Omega_{X^\mr{log}/k^\mr{log}}$ and $\Om := \Omega_{X/k}$.
(Hence, there exists a natural inclusion $\Om \migiincl \Omega$, inducing an isomorphism $\Om\isom \Omega (-D_\msX)$.)
Write $k_\epsilon := k[\epsilon]/(\epsilon^2)$, and denote by
\begin{align} \label{QG56}
\msX_\epsilon := (X_\epsilon/k_\epsilon, \{ \sigma_{i\epsilon } \}_{i=1}^r)
\end{align}
the base-change of $\msX$ over $k_\epsilon$.
Also, write $\mr{pr}_X : X_\epsilon \migi X$ and $\mr{in}_X : X \migi X_\epsilon$ for the natural projection and  closed immersion.

\end{notation}


\section{The deformation space of a flat bundle with ${^p}\psi =0$} \label{QG1170}

\subsection{The map $\varpi$} \label{y0187}

First,  we construct the map ``$\varpi$'' 
 appearing in (\ref{QQ699}).
Let $(\mcF, \nabla)$ be a  rank $n$ flat bundle on $X^\mr{log}/k^\mr{log}$
 with ${^p }\psi^\nabla =0$.
 The sheaf $\mcE nd (\mcF) := \mcE nd_{\mcO_X}(\mcF)$ 
 may be naturally  identified with the adjoint vector bundle associated to the $\mr{GL}_n$-bundle corresponding to  $\mcF$; it is
  equipped with  a $k^\mr{log}$-connection $\nabla^\mr{ad}$ induced by $\nabla$ via the adjoint representation $\mr{Ad}_{\mr{GL}_n}  : \mr{GL}_n \migi \mr{GL}(\mfg \mfl_n)$ (cf. (\ref{adjrep})).
One may identify  this connection $\nabla^\mr{ad}$ with $\nabla^\vee \otimes \nabla$ via  the natural  isomorphism $\mcF^\vee \otimes \mcF \isom \mcE nd (\mcF)$.   
Denote by 
\begin{align} \label{QG40}
\mr{Def}_{(\mcF, \nabla)}
\index{$\mr{Def}_{(\mcF, \nabla)}$, $\mr{Der}_\nabla$}
\end{align}
the set of isomorphism classes   of
first-order  deformations
  of $(\mcF,\nabla)$, i.e.,
the set of isomorphism classes of  flat bundles $(\mcF_{\epsilon'}, \nabla_{\epsilon'})$ on $X_\epsilon^\mr{log}/k^\mr{log}_\epsilon$ 
together with an isomorphism $\mr{in}_X^*(\mcF_{\epsilon'}, \nabla_{\epsilon'}) \isom (\mcF, \nabla_{\mcF})$.
As discussed in \S\,\ref{QW6081}, 
there exists a canonical bijection
\begin{equation} \label{10602}
\mr{Def}_{(\mcF, \nabla)} \isom \mbH^1 (X, \mcK^\bullet [\nabla^\mr{ad}])
\end{equation}
between $\mr{Def}_{(\mcF, \nabla)}$
and (the underlying  set of) the $1$-st hypercohomology group $\mbH^1 (X, \mcK^\bullet [\nabla^\mr{ad}])$ of the complex $\mcK^\bullet [\nabla^\mr{ad}]$.
Let us consider the subset 
\begin{align} \label{QQ700}
\mr{Def}_{\nabla}  
\index{$\mr{Def}_{(\mcF, \nabla)}$, $\mr{Der}_\nabla$}
\end{align}
 of $\mr{Def}_{(\mcF, \nabla)}$ consisting of deformations $(\mcF_{\epsilon'}, \nabla_{\epsilon'})$ that  admits  an $\mcO_{X_\epsilon}$-linear  isomorphism $\mcF_{\epsilon'} \isom \mr{pr}_X^*(\mcF)$ inducing the identity  morphism of  $\mcF$ via reduction modulo $\epsilon$.

\bpr \label{QQ513}
 The canonical bijection $\mr{Def}_{(\mcF, \nabla)} \isom \mbH^1 (X, \mcK^\bullet [\nabla^\mr{ad}])$ (cf. (\ref{10602})) restricts to a bijection
\begin{align} \label{QQ502}
\mr{Def}_{\nabla} \isom  \mr{Ker}({'e}_\flat [\nabla^\mr{ad}])
\end{align}
(cf. (\ref{HDeq}) for the definition of ${'e}_\flat [-]$). 
\epr
\begin{proof}
The set of isomorphism classes of first-order deformations 
 of the vector bundle $\mcF$ corresponds bijectively to $H^1 (X, \mcE nd (\mcF))$ (cf. (\ref{QW6077})).
The assignment $(\mcF_{\epsilon'}, \nabla_{\epsilon'}) \mapsto \mcF_{\epsilon'}$ determines, via this correspondence and (\ref{10602}),  a morphism 
$ \mbH^1 (X, \mcK^\bullet [\nabla^\mr{ad}]) \migi H^1 (X, \mcE nd (\mcF))$ of $k$-vector spaces.
By considering  the description of the relevant deformation spaces in terms of \v{C}ech cocycles, 
we see  that this morphism coincides with ${'e}_\flat [\nabla^\mr{ad}]$.
Hence,  the subspace
 \begin{equation}
\mr{Ker}({'e}_\flat [\nabla^\mr{ad}]) \subseteq \mbH^1 (X, \mcK^\bullet [\nabla^\mr{ad}])\end{equation}
 corresponds  to the deformations  $(\mcF_{\epsilon'}, \nabla_{\epsilon'})$  of $(\mcF, \nabla_\mcF)$
 with $\mcF_{\epsilon'}$ trivial.
This completes the proof of this proposition.
\end{proof}

Then, we define $\varpi$ to be the composite surjection 
\begin{align}
\label{QQ701}
\varpi :  \Gamma (X, \Omega\otimes \mcE nd(\mcF))  &\xrightarrow{\mr{quotient}} \mr{Coker}(\Gamma(X, \nabla^\mr{ad})) \\
&  \xrightarrow
{{'e}_\sharp [\nabla^\mr{ad}]} \mr{Ker}({'e}_\flat [\nabla^\mr{ad}]) \notag \\
& \xrightarrow
{(\ref{QQ502})^{-1}}
\mr{Def}_\nabla. \notag
\index{$\varpi$}
\end{align}

\subsection{The sets  ${^p}\mr{Def}_{(\mcF, \nabla)}$ and $ {^p}\mr{Def}_{\nabla}$} \label{QG1230}

Next, let us consider the subset 
 \begin{align} \label{QG47}
 {^p}\mr{Def}_{(\mcF, \nabla)} 
 \index{${^p}\mr{Def}_{(\mcF, \nabla)}$, ${^p}\mr{Def}_{\nabla}$, ${^p}\mr{Def}^\sharp_{\nabla}$, ${^p}\mr{Def}^{\mr{Op}}_{\nabla^\diamondsuit}$}
 \end{align}
 of $\mr{Def}_{(\mcF, \nabla)}$ consisting of deformations $(\mcF_{\epsilon'}, \nabla_{\epsilon'})$ with ${^p}\psi^{\nabla_{\epsilon'}} =0$.
Here, notice that
the $k$-linear map
\begin{align} \label{QQ702}
H^1 (X, \mr{Ker}(\nabla^\mr{ad})) \migi \mbH^1 (X, \mcK^\bullet [\nabla^\mr{ad}]) 
\end{align}
induced by the natural  morphism $\mr{Ker}(\nabla^\mr{ad})[0] \migi \mcK^\bullet [\nabla^\mr{ad}]$ is verified to be injective;
by using  this injection, we shall regard $H^1 (X, \mr{Ker}(\nabla^\mr{ad}))$ as a 
a subset of $\mbH^1 (X, \mcK^\bullet [\nabla^\mr{ad}])$ for convenience.

\bpr \label{QQ512}
The canonical bijection $\mr{Def}_{(\mcF, \nabla)} \isom \mbH^1 (X, \mcK^\bullet [\nabla^\mr{ad}])$ (cf. (\ref{10602})) 
 restrict to a bijection
\begin{align} \label{QQ501}
{^p}\mr{Def}_{(\mcF, \nabla)} \isom  
H^1 (X, \mr{Ker}(\nabla^\mr{ad}))
\end{align}
\epr
\begin{proof}
 Denote by
\begin{equation}
\label{10600}
 \overline{\nabla}^\mr{ad} : \mcE nd (\mcF) \migi \mr{Im} (\nabla^\mr{ad}) \left(= \mr{Ker}(C^{(\mcF, \nabla)}) \right)\end{equation}
  the morphism obtained by restricting the codomain of $\nabla^\mr{ad}$.
Since 
 the natural morphism $\mr{Ker}(\nabla^\mr{ad}) [0] \migi \mcK^\bullet [\overline{\nabla}^\mr{ad}]$ is a quasi-isomorphism, we have  an isomorphism
\begin{equation} \label{10601}
H^1 (X, \mr{Ker}(\nabla^\mr{ad})) \isom \mbH^1 (X, \mcK^\bullet [\overline{\nabla}^\mr{ad}]).
\end{equation}
Hence, it follows from Proposition \ref{y0157} that   the image of (\ref{QQ702}) corresponds, via   (\ref{10602}), to the deformations 
 of $(\mcF, \nabla)$  with vanishing $p$-curvature.
\end{proof}

Let us write
\index{${^p}\mr{Def}_{(\mcF, \nabla)}$, ${^p}\mr{Def}_{\nabla}$, ${^p}\mr{Def}^\sharp_{\nabla}$, ${^p}\mr{Def}^{\mr{Op}}_{\nabla^\diamondsuit}$}
\begin{align} \label{QQ641}
{^p}\mr{Def}_\nabla := \mr{Def}_\nabla \cap {^p}\mr{Def}_{(\mcF, \nabla)}.
\end{align}

\subsection{The map $\Upsilon$} \label{QG1110}

In what follows, we shall construct a map  ``$\Upsilon$''  from ${^p}\mr{Def}_\nabla$ to the tangent space of a certain Quot scheme.
We shall write 
\begin{equation} 
\label{1070}
\A
 := F^*_{X/k}(\mr{Ker}(\nabla)).
 \index{$\A$}
 \end{equation}
The natural inclusion $\mr{Ker}(\nabla) \migiincl F_{X/k*}(\mcF)$ corresponds,
via  the adjunction relation 
 ``$F^*_{X/k} (-) \dashv F_{X/k*}(-)$" ,
to an $\mcO_{X}$-linear  morphism 
\begin{equation} \label{inj77}
 \nu^\sharp : \mcF^\sharp  \migi \mcF.
 \index{$\nu^\sharp$}
 \end{equation}
The morphism $ \nu^\sharp$  is injective and (by (\ref{equicat})) becomes an isomorphism when restricted to the open subscheme $X^\mr{sm} = X \setminus \bigcup_{i=1}^r \mr{Im}(\sigma_i)$.
Also, it  specifies a $k$-rational point $[\nu^\sharp]$ of the Quot-scheme $\mr{Quot}_{\mcF/X/k}^{n, \mr{deg}(\A)}$  (cf. \S\,\ref{QG1409}),  in particular,  we obtain 
   the tangent space
\begin{align} \label{QQ506}
\mcT_{[\nu^\sharp]} \mr{Quot}_{\mcF/X/k}^{n, \mr{deg}(\mcF^\sharp)}
\end{align} 
 of  $\mr{Quot}_{\mcF/X/k}^{n, \mr{deg}(\A)}$ at the point $[\nu^\sharp]$.

Now, let us take an arbitrary  deformation $(\mcF_{\epsilon'}, \nabla_{\epsilon'})$ classified by ${^p}\mr{Def}_{\nabla}$.
The $\mcO_{X_\epsilon}$-module
\begin{align} \label{QQ601}
\A_{\epsilon'} := F_{X_\epsilon/k_\epsilon}^*(\mr{Ker}(\nabla_{\epsilon'}))
\end{align}
specifies  a deformation of $\A$ (cf. Corollary \ref{y0160}). 
The inclusion $\mr{Ker}(\nabla_{\epsilon'}) \migi F_{X_\epsilon/ k_\epsilon*} (\mcF_{\epsilon'})$ induces,  in the same manner as $\nu^\sharp$,   an  injection
\begin{align} \label{QQ603}
\nu_{\epsilon'}^\sharp  : \A_{\epsilon'} \migiincl \mcF_{\epsilon'}, 
\end{align}
forming  a deformation of $\nu^\sharp$.
Here, we shall use  the symbol $\EE$ to denote either the presence or absence of ``\,$\epsilon'$\,".
Then, there exists  a short exact sequence 
\begin{equation} \label{short4}
0 \longmigi \A_{\EE}
\xrightarrow{\nu^\sharp_{\EE}} \mcF_{\EE} \xrightarrow{\nu^\flat_{\EE}}  \bigoplus_{i=1}^r \Lambda_{i, {\EE}} \longmigi 0,
\index{$\Lambda_{i, {\EE}}$}
\end{equation} 
 where
each  $\Lambda_{i, {\EE}}$ ($i=1, \cdots, r$) is  supported on
 $\mr{Im}(\sigma_{i})$.

\ble \label{15232} 
For each $i=1, \cdots, r$, 
  the sheaf  $\Lambda_{i, \epsilon'}$ 
    is  flat over $k_{\epsilon}$.
 \ele
\begin{proof}
It follows from  Corollary \ref{y0158} that $(\mcF_{\epsilon'}, \nabla_{\epsilon'})$ can be, Zariski locally on $X_{\epsilon}$, obtained as the base-change of $(\mcF, \nabla)$.
Hence, by the definitions of $\nu^\sharp_{\epsilon'}$ and $\nu^\flat_{\epsilon'}$,
each $\Lambda_{i, \epsilon'}$ is, Zariski locally on $X_{\epsilon}$,  isomorphic to $\mr{pr}_X^*(\Lambda_i)$,  so it is  flat over $k_\epsilon$.
\end{proof}

Since $(\mcF_{\epsilon'}, \nabla_{\epsilon'})$ is classified by $\mr{Def}_\nabla$,
there exists an isomorphism $\mcF_{\epsilon'} \isom \mr{pr}_X^*(\mcF)$.
By Lemma \ref{15232} above,  the cokernel of the composite
\begin{align}
[\nu_{\epsilon'}^\sharp] : \A_{\epsilon'} \migiincl \mcF_{\epsilon'} \isom \mr{pr}_X^*(\mcF)
\end{align}
 is flat over $k_\epsilon$.
 Hence, it  specifies an element of $\mcT_{[\nu^\sharp]} \mr{Quot}_{\mcF/X/k}^{n, \mr{deg}(\mcF^\sharp)}$.
This element does not depend on the choice of the isomorphism $\mcF_{\epsilon'} \isom \mr{pr}_X^*(\mcF)$.
Thus, the assignment $(\mcF_{\epsilon'}, \nabla_{\epsilon'}) \mapsto [\nu^\sharp_{\epsilon'}]$ determines a map of sets
\begin{align} \label{QQ507}
\Upsilon :
{^p}\mr{Def}_\nabla \migi \mcT_{[\nu^\sharp]} \mr{Quot}_{\mcF/X/k}^{n, \mr{deg}(\A)}.
\index{$\Upsilon$}
\end{align}

\subsection{The set ${^p}\mr{Def}_{\nabla}^\sharp$ and the second ``$\subseteq$''} \label{y0191}
We prove  an inclusion relation between two sets inducing the second  ``$\subseteq$''   in (\ref{QQ699}) under the assumption that $g=0$.
Denote by
\index{${^p}\mr{Def}_{(\mcF, \nabla)}$, ${^p}\mr{Def}_{\nabla}$, ${^p}\mr{Def}^\sharp_{\nabla}$, ${^p}\mr{Def}^{\mr{Op}}_{\nabla^\diamondsuit}$}
\begin{align} \label{QG50}
 {^p}\mr{Def}_{\nabla}^\sharp
\end{align}
the subset of ${^p}\mr{Def}_\nabla$ consisting of deformations $(\mcF_{\epsilon'}, \nabla_{\epsilon'})$
with $\A_{\epsilon'} \cong \mr{pr}_X^*(\A)$.

\bpr \label{QG71}
Suppose that $g=0$.
Then, we have an inclusion relation
\begin{align}
{^p}\mr{Def}_{\nabla}^\sharp \subseteq \Upsilon^{-1}(0)
\end{align}
between subsets of ${^p}\mr{Def}_\nabla$.
\epr
\begin{proof}
Since $g=0$, i.e., $X \cong  \mbP^1$,  it  follows from  the Birkhoff-Grothendieck  theorem that both the vector bundles $\mcF$ and $\A$ can be decomposed into direct sums of line bundles;
let us fix decompositions $\mcF = \bigoplus_{j=1}^n \mcL_i$, $\A = \bigoplus_{j=1}^n \mcL^\sharp_i$, where  $\mcL_i$'s and $\mcL^\sharp_i$'s  are line bundles.
Now, let 
$v$ be an arbitrary element of ${^p}\mr{Def}_{\nabla}^\sharp$, and choose a deformation  
$(\mcF_{\epsilon'}, \nabla_{\epsilon'})$ representing $v$.
By the definition of ${^p}\mr{Def}_{\nabla}^\sharp$, we can identify $\mcF_\epsilon$ and $\A_\epsilon$ with $\bigoplus_{j=1}^n \mr{pr}_X^*(\mcL_i)$ and $\bigoplus_{j=1}^n \mr{pr}_X^*(\mcL^\sharp_i)$ respectively.
The morphism
$\nu^\sharp_{\epsilon'} : \A_{\epsilon'} \migi \mcF_{\epsilon'}$ determines an element of  
\begin{align}
\mr{Hom}_{\mcO_X}(\A_{\epsilon'}, \mcF_{\epsilon'}) &= \mr{Hom}_{\mcO_X}(\bigoplus_{j=1}^n \mr{pr}_X^*(\mcL^\sharp_i), \bigoplus_{j=1}^n \mr{pr}_X^*(\mcL_i)) \\
& = \bigoplus_{i, i^\sharp}\mr{Hom}_{\mcO_X}(\mr{pr}^*_X(\mcL^\sharp_{i^\sharp}), \mr{pr}^*_X(\mcL_i)). \notag
\end{align}
Let  $(\nu^\sharp_{\epsilon'})_{i, i^\sharp}$ denote the $(i, i^\sharp)$-th component of $\nu^\sharp_{\epsilon'}$ with respect to this decomposition.
It follows from the $k_\epsilon$-flatness of $\mr{Coker} (\nu^\sharp_{\epsilon'})$  that the cokernel $\mr{Coker}((\nu^\sharp_{\epsilon'})_{i, i^\sharp})$ for every pair  $(i, i^\sharp)$ is flat over $k_\epsilon$.
Since both  $\mr{pr}^*_X(\mcL^\sharp_{i^\sharp})$ and $\mr{pr}^*_X(\mcL_i)$ are line bundles, 
$\mr{Im}((\nu^\sharp_{\epsilon'})_{i, i^\sharp})$ is the trivial deformation of the image of  the morphism $\mr{in}_X^*((\nu^\sharp_{\epsilon'})_{i, i^\sharp}) : \mcL^\sharp_{i^\sharp} \migi \mcL_i$ between the special fibers.
This means that $\mr{Im}(\nu^\sharp_{\epsilon'}) = \mr{pr}^*_X (\nu^\sharp)$, which implies  $\varpi (v) =0$.
This completes the proof of this proposition.
\end{proof}

\subsection{The third ``$\subseteq$''} \label{QG1299}
Next, we prove an inclusion relation inducing the  third  ``$\subseteq$''   in (\ref{QQ699}).
Let us consider the  morphism
\begin{equation}
\label{10018}
\mcE nd (\mcF) \migi\mcH om_{\mcO_X} (\A, \mcF), 
\hspace{5mm}
\mcE nd (\A) \migi \mcH om_{\mcO_X}(\A, \mcF)
\end{equation}
given by, respectively,   $h \mapsto h \circ \nu^\sharp$ and  $h' \mapsto \nu^\sharp \circ h'$ for any local sections $h \in \mcE nd (\mcF)$, 
$h' \in \mcE nd (\A)$.
These morphisms are injective, so we can  regard both 
$\mcE nd (\mcF)$ and $\mcE nd  (\A)$ as 
$\mcO_X$-submodules of $\mcH om_{\mcO_X} (\A, \mcF)$ via these injections.
In particular, we have a subspace
\begin{align} \label{QG58}
\Gamma (X,  \Omega\otimes \mcE nd (\mcF) \cap \Om\otimes \mcE nd (\mcF^\sharp))
\end{align}
of $\Gamma (X, \Omega\otimes \mcE nd (\mcF))$.

\bpr \label{QG70}
We have an inclusion relation
\begin{align} \label{QQ708}
\Upsilon^{-1}(0) \subseteq  \varpi(\Gamma (X, \Omega\otimes \mcE nd (\mcF) \cap \Om\otimes \mcE nd (\mcF^\sharp)))
\end{align}
between subsets of $\mr{Def}_\nabla  \left(\supseteq {^p}\mr{Def}_\nabla\right)$.
\epr
\begin{proof}
Let $v$ be an element of 
${^p}\mr{Def}_\nabla$
   and choose an element of the preimage   $v' \in \varpi^{-1}(v) \subseteq \Gamma (X, \Omega \otimes \mcE nd (\mcF))$.
Since $v$ corresponds, via (\ref{QQ501}) and (\ref{10601}), to an element of 
$\mbH^1 (X, \mcK^\bullet [\overline{\nabla}^\mr{ad}])$,
there exists
a Zariski open covering $\{U_i \}_{i \in I}$ of $X$ and sections $v''_i \in \Gamma (U_i, \mcE nd (\mcF))$ ($i \in I$) with 
$v' |_{U_i}= \nabla^\mr{ad}(v''_i)$.
Then, the sections
\begin{align}
\nu^\flat \circ v''_i \circ \nu^\sharp \in \Gamma (U_i, \mcH om_{\mcO_X} (\A,  \bigoplus_{i=1}^r \Lambda_i))
\end{align}
(cf. (\ref{short4}) for the definitions of $\nu^\sharp$ and $\nu^\flat$)
  may be glued  together to obtain  a well-defined global  section
$v''' \in \mr{Hom}_{\mcO_X}(\A,  \bigoplus_{i=1}^r \Lambda_i)$.
This section satisfies  the  equality $v''' = \Upsilon (v)$ 
under the canonical bijection
\begin{align} \label{QQ510}
\mcT_{[\nu^\sharp]} \mr{Quot}_{\mcF/X/k}^{n, \mr{deg}(\A)}
\isom \mr{Hom}_{\mcO_X} (\A, \bigoplus_{i=1}^r \Lambda_i)
\end{align}
(cf. (\ref{10101})).

Now, we shall suppose that $v \in \Upsilon^{-1}(0)$, i.e.,  $v''' =0$.
By the above discussion,  the equality  $\nu^\flat \circ v''_i \circ \nu^\sharp =0$ holds,  so
$v''_i $ belongs to  $\mcE nd (\A) \subseteq \mcH om_{\mcO_X}(\A, \mcF)$.
Note that  
 $\mcE nd (\A)$ and  $\mcH om_{\mcO_X} (\A, \mcF)$  
admit
$k^\mr{log}$-connections 
  $\nabla^{\mr{can}\vee}_{\mcV} \otimes \nabla^\mr{can}_{\mcV}$ and $\nabla^{\mr{can}\vee}_{\mcV} \otimes  \nabla$ respectively, where $\mcV := \mr{Ker}(\nabla)$.
  Moreover, 
  the  inclusions displayed in (\ref{10018})  
 are  compatible with the  $k^\mr{log}$-connections.
 It follows that
 the element  $v'$ belongs to $\Gamma (X, \mr{Im}(\nabla_\mcV^{\mr{can}\vee} \otimes \nabla_\mcV^\mr{can}))$.
 Since the image of $\nabla_\mcV^{\mr{can}\vee} \otimes \nabla_\mcV^\mr{can}$
 lies in $\Om\otimes \mcE nd (\A)$,
 we see that $v' \in \Gamma (X, \Om\otimes \mcE nd (\A))$, i.e., $v \in \varpi (\Gamma (X, \Om\otimes \mcE nd (\A)))$.
 Consequently,  we obtain  (\ref{QQ708}), as desired.
\end{proof}

\section{Formally local description of a flat bundle} \label{y0189}

In this section, we study the local description of a given $\mr{GL}_n$-oper at the  marked points of $\msX$.
Suppose that $r >0$ and fix
 $i \in \{1, \cdots, r \}$.
Denote  by $[\sigma_i]$  the reduced effective divisor on $X$ defined by the image of  $\sigma_i$.
For $l \in \mbZ_{\geq 0}$, there exists a unique $k^\mr{log}$-connection
\index{$\nabla$, $\hslash$-$T^\mr{log}$-connection@--- $\nabla_{l \cdot [\sigma_i]}$}
\begin{equation}
\nabla_{l \cdot [\sigma_i]} : \mcO_X (l \cdot [\sigma_i]) \migi \Omega  \otimes \mcO_X (l \cdot [\sigma_i]) 
\end{equation} 
on the line bundle $\mcO_X (l \cdot [\sigma_i]) $
whose restriction to $\mcO_X \left(\subseteq \mcO_X (l\cdot [\sigma_i]) \right)$ coincides with the universal  derivation $d : \mcO_X \migi  \Omega$.
For $l \in \mbZ_{\geq 0}$, let
$\breve{\Lambda}^i_l$ denote 
the $\mcO_X$-module fitting into the following short exact sequence:
\begin{align} \label{QQ440}
0 \longmigi \mcO_X \xrightarrow{\nu^{\sharp i}_{l}} 
\mcO_X (l  \cdot [\sigma_i]) 
\xrightarrow{\nu^{\flat i}_{l}} \breve{\Lambda}^i_l \longmigi 0,
\index{$\breve{\Lambda}^i_l$}
\end{align}
where $\nu^{\sharp i}_{l}$ denotes the natural inclusion.

\ble \label{55232}
\begin{itemize}
\item[(i)]
For $l_1$, $l_2 \in  \mbZ_{\geq 0}$, we have  
\begin{equation} 
\label{1080}
\nabla_{l_1 \cdot [\sigma_i]} \otimes \nabla_{l_2 \cdot [\sigma_i]} = \nabla_{(l_1 +l_2) \cdot [\sigma_i]}\end{equation}
via  the natural identification 
 $\mcO_X (l_1 \cdot [\sigma_i]) \otimes \mcO_X (l_2 \cdot [\sigma_i]) = \mcO_X ((l_1 +l_2) \cdot [\sigma_i])$.
\item[(ii)]
Let $l \in \mbZ_{\geq 0}$.
Then, the monodromy of  $(\mcO_X (l\cdot [\sigma_i]), \nabla_{l \cdot [\sigma_i]})$ at $\sigma_i$, considered as an element of $k$, 
satisfies 
\begin{equation} \label{86695}
\mu_i^{\nabla_{l \cdot [\sigma_i]}}  = - \overline{l},
\end{equation}
where $\overline{(-)}$ denotes the image  via the quotient $\mbZ \migisurj \mbF_p \left(\subseteq k \right)$.
Moreover, we have an equality 
\begin{align}
\mr{Ker}(\nabla_{l \cdot [\sigma_i]}) = \mcO_{X^{(1)}} (\lfloor  l/p\rfloor  \cdot [\sigma_i^{(1)}])
\end{align}
of $\mcO_{X^{(1)}}$-submodules of $F_{X/k*}(\mcO_X)$, where $\lfloor  - \rfloor$ denotes the floor function.
 \end{itemize}
 \ele
\begin{proof} 
The assertions follow immediately from the definition of $\nabla_{l \cdot [\sigma_i]}$.
\end{proof}

Let us fix a local function $t \in \mcO_X$ defining $[\sigma_i]$.
Then, we can construct  the formal completion $( \mcO_X (l \cdot [\sigma_i]), \nabla_{l \cdot [\sigma_i]})^\wedge$ of $( \mcO_X (l \cdot [\sigma_i]), \nabla_{l \cdot [\sigma_i]})$ with respect to the $t$-adic topology;
 it is isomorphic to 
 the pair of $\widehat{\mcO} := \mcO_{\mr{Spf}(k[[t]])}$  and  the  regular singular connection  on $\widehat{\mcO}$ of the form
  \index{$\nabla$, $\hslash$-$T^\mr{log}$-connection@--- $\widehat{\nabla}_l$}
 \begin{equation}
 \widehat{\nabla}_l : = d  - l  \cdot \frac{dt}{t}.
 \end{equation}

Now,
 recall that $(\mcF, \nabla)$ is a rank $n$ flat vector bundle on $X^\mr{log}/k^\mr{log}$ with ${^p}\psi^\nabla =0$.
According to  ~\cite{O2}, Corollary 2.10,   the $t$-adic completion $(\mcF, \nabla)^\wedge$ of $(\mcF, \nabla)$ is isomorphic to a direct sum  of various $ (\widehat{\mcO}, \widehat{\nabla}_l )$'s with $l \in \mbZ_{\geq 0}$.
Since  $ \widehat{\nabla}_{l'} =  \widehat{\nabla}_{l''}$ if and only if $l' \equiv l''$ mod $p$,
 all such integers $l$ can be assumed to be less than $p$.
Hence, there exists   an isomorphism
\begin{equation} \label{decom9}
\eta_i : (\mcF, \nabla)^\wedge \isom \bigoplus_{l =1}^n ( \widehat{\mcO}, \widehat{\nabla}_{m_{i, l}}),
\end{equation}
where $m_{i, 1}, \cdots, m_{i, l}$ are nonnegative integers with
$m_{i, 1} \leq \cdots \leq m_{i, n} < p$.
If $\widehat{\nu}^\sharp$ and $\widehat{\nu}^{\sharp i}_{m_{i, j}}$ ($j=1, \cdots, n$) denote the $t$-adic completions of
$\nu^\sharp$ and $\nu^{\sharp i}_{m_{i, j}}$ respectively, then
 $\eta_i$ 
   induces a composite  isomorphism
\begin{equation} \label{18006}
\Lambda_i \isom 
 \mr{Coker}(\widehat{\nu}^\sharp) \isom \mr{Coker}(\bigoplus_{j=1}^n \widehat{\nu}^{\sharp i}_{m_{i, j}}) \isom \bigoplus_{j=1}^n \mr{Coker}(\widehat{\nu}^{\sharp i}_{m_{i, j}}) \isom \bigoplus_{j=1}^n \breve{\Lambda}^i_{m_{i, j}}.
\end{equation}

\bde \label{QQ716}
Keeping the above situation, 
let us write 
\begin{equation} \label{1082}
\mr{Res}_{\sigma_i}( \nabla) := (- m_{i, l})_{l=1}^n \in \mbZ^n,
\index{$\mr{Res}_{\sigma_i}( \nabla)$, residue of $\nabla$}
\end{equation}
and refer to it as the {\bf residue} of $\nabla$ at $\sigma_i$.
(Notice that this 
$n$-tuple of integers  depends neither on the choice of  $t$ nor on the choice of    $\eta_i$.)
\ede


The relationship between residue and monodromy is as follows.

\bpr[cf. ~\cite{O2}, Corollary 2.10] \label{QG78}
Let $(\mcF, \nabla)$ be a  rank $n$ flat vector bundle on  $X^\mr{log}/k^\mr{log}$ with ${^p}\psi^\nabla =0$  and $i \in \{1, \cdots, r \}$.
Also,   write $(-m_{i, l})_{l=1}^n :=  \mr{Res}_{\sigma_i}(\nabla)$, where $0 \leq m_{i, 1} \leq \cdots \leq m_{i, n} < p$.
Then, the monodromy $\mu_i^\nabla$ of $\nabla$ at $\sigma_i$ (considered as an element of $\mfg \mfl_n$ by means of an  identification $\sigma_i^*(\mcF) \cong k^{\oplus n}$) is conjugate to the diagonal matrix whose diagonal entries are $m_{i,1}, \cdots, m_{i, n}$.
\epr
\begin{proof}
The assertion follows from  Lemma \ref{55232}, (ii).
\end{proof}

\bpr \label{10002}
Let us keep the above notation.
Moreover, suppose   that there exists a decreasing  filtration $\{ \mcF^j \}_{j =0}^n$  on $\mcF$ such that  the triple 
$(\mcF, \nabla, \{ \mcF^j \}_{j =0}^n)$ forms a (dormant) $\mr{GL}_n$-oper on $\msX$.
Then, the following assertions hold:
\begin{itemize}
\item[(i)]
 Let $v$ be a global section of 
 the $t$-adic completion $\widehat{\mcF}$ of $\mcF$
  generating   (formally) 
 the line subbundle 
 $\widehat{\mcF}^{n-1} \left(\subseteq \widehat{\mcF} \right)$.
 We shall write
\begin{equation}
    (u_l(t))_{l=1}^n  := \eta_i (v)   \in k[[t]]^{\oplus n}\left(= \Gamma (\mr{Spf}(k[[t]]), \widehat{\mcO})^{\oplus n}\right). \end{equation}
  Then, each  element  $u_l(t)$ belongs to  $k[[t]]^\times$, or equivalently,  
  $u_l(0) \neq 0$.
 \item[(ii)]
 The factors $-m_{i, 1}, \cdots,  -m_{i, n}$ in the residue  $\mr{Res}_{\sigma_i} (\nabla)$ are mutually distinct.
 \end{itemize}
 \epr
\begin{proof}
First, we shall consider assertion (i).
Using  the assignment $w \mapsto \frac{dt}{t} \otimes w$,  we shall identify  $\widehat{\mcF}$ with 
$(\Omega\otimes \mcF)^\wedge$.
Then, the $t$-adic completion $\widehat{\nabla}$ of $\nabla$
may be regarded as an endomorphism of $\widehat{\mcF}$.
We shall write $\widehat{\nabla}^{j}$ for the $j$-th iterate 
of this endomorphism,
where $\widehat{\nabla}^{0} := \mr{id}_{\widehat{\mcF}}$.
It follows from the definitions of a $\mr{GL}_n$-oper and the section $v$ together  that the set $\{ \widehat{\nabla}^{j}(v )\}_{j=0}^{n-1}$ forms a basis of $\widehat{\mcF}$.
On the other hand, 
since $\eta_i$ is compatible with the $k^\mr{log}$-connections, 
the equality 
$\eta_i (\widehat{\nabla}^{j}(v )) = (\widehat{\nabla}_{m_{i, l}}^{j} (u_l (t)))_{l=1}^n$
holds for every $j = 0, \cdots, n-1$.
In particular, 
for each $l$,  the set $\{  \widehat{\nabla}_{m_{i, l}}^{ j} (u_l (t))\}_{j=0}^{n-1}$ generates the trivial $\widehat{\mcO}$-module  $\widehat{\mcO}$.
If we  write 
\begin{equation}
u_l (t) = \sum_{s =0}^\infty u_{l, s} \cdot t^s
\end{equation}
($u_{l,s} \in k$), 
then  $\widehat{\nabla}_{m_{i, l}}^{ j} (u_l (t))$ may be expressed as
\begin{equation} \label{10001}
\widehat{\nabla}_{m_{i, l}}^{ j} (u_l (t)) = \sum_{s =0}^\infty (\overline{s} - \overline{m}_{i, l})^j \cdot u_{l, s} \cdot t^s.
\end{equation}
Thus,   $\{  \widehat{\nabla}_{m_{i, l}}^{ j} (u_l (t))\}_{j=0}^{n-1}$ generates $\widehat{\mcO}$ if and only if $(-\overline{m}_{i,l})^j \cdot u_{l, 0} \neq 0$ for some $j \geq 0$, equivalently,
$u_{l, 0} \ (= u_l (0)) \neq 0$. 
This completes the proof of assertion (i).

Next, we  shall consider assertion (ii).
Since $\{  (\widehat{\nabla}_{m_{i, l}}^{ j} (u_l (t)))_{l =1}^n\}_{j=0}^{n-1}$ forms a basis for $\widehat{\mcO}^{\oplus n}$,
the $n \times n$ matrix
\begin{equation} (\widehat{\nabla}_{m_{i, a}}^{b-1} (u_a (t)) |_{t=0})_{a, b =1}^n,
 \end{equation} 
i.e., the matrix whose $(a, b)$-entry equals $\widehat{\nabla}_{m_{i, a}}^{b-1} (u_a (t)) |_{t=0}$,  is regular.
But, by  (\ref{10001}), 
 the following sequence of  equalities holds:
\begin{align} \label{177002}
\mr{det}( (\widehat{\nabla}_{m_{i, a}}^{b-1} (u_a (t)) |_{t=0})_{a, b =1}^n) 
&= \mr{det} (((-\overline{m}_{i, a})^{b-1}\cdot u_{a, 0})_{a, b =1}^n) \\
& = \mr{det} (((-\overline{m}_{i, a})^{b-1})_{a, b =1}^n)  \cdot \prod_{a =1}^n \cdot (u_{a, 0})^n \notag \\
& = \prod_{a =1}^n  (u_{a, 0})^n  \cdot \prod_{1 \leq a_2 < a_1 \leq n} (-\overline{m}_{i, a_1} + \overline{m}_{i, a_2}), \notag
\end{align}
where the last equality follows from Vandermonde's determinant and the factor $ \prod_{a =1}^n  (u_{a, 0})^n$  in the rightmost is invertible  by virtue of assertion (i).
Hence, the elements  $\overline{m}_{i, 1}, \cdots, \overline{m}_{i, n}$ (in particular,   the elements  $-m_{i, 1}, \cdots, -m_{i, n}$) must be  mutually distinct.
This completes   the proof of assertion (ii).
\end{proof}

\section{Deformation of dormant $\mr{GL}_n$-opers} \label{QQ710}

\subsection{} \label{QG1400}
Let $\msF^\heartsuit := (\mcF, \nabla, \{ \mcF^j \}_{j=0}^n)$
be a dormant $\mr{GL}_n$-oper on $\msX$.
By applying  the previous discussions  to the flat bundle $(\mcF, \nabla)$,   
 we obtain  the various objects associated to this flat bundle, e.g.,  $\A$, $\nu^\sharp$, and $\nu^\flat$, etc..
Let us define  a decreasing filtration $\{ \mcF^{\sharp j}\}_{j=0}^n$ on $\A$ (resp., $\{ \Lambda_{i}^j \}_{j =0}^n$ on each $\Lambda_{i}$)  as  
\begin{equation} \label{1074}
\mcF^{\sharp j}:=  (\nu^\sharp)^{-1}(\mcF^j)  \ \left(\text{resp.,} \
\Lambda^j_{i} :=  \Lambda_{i} \cap  \nu^\flat (\mcF^j)   \right).
\end{equation}
The collection  $\{ \bigoplus_{i=1}^r \Lambda^j_{i} \}_{j=0}^n$ defines a filtration on $\bigoplus_{i=1}^r \Lambda_{i}$. 
Both  $\nu^\sharp$ and $\nu^\flat$ preserves  the filtrations, so  (\ref{short4}) yields  short exact sequences of $\mcO_{X}$-modules
\begin{equation} \label{short2}
  0 \longmigi  \mcF^{\sharp j}/\mcF^{\sharp j+1}\xrightarrow{\nu^{\sharp j}} \mcF^j /\mcF^{j+1}  \xrightarrow{\nu^{\flat j}} \bigoplus_{i=1}^r \Lambda_{ i}^j/\Lambda^{j+1}_{i} \longmigi 0
\end{equation}
($j =0, \cdots, n-1$).

\bpr \label{11003}
Let us fix $j \in \{0, \cdots, n-1 \}$, and write 
 $(- m_{i, l})_{l=1}^n := \mr{Res}_{\sigma_i}(\nabla)$ ($i=1, \cdots, r$).
Then, the following assertion holds:
\begin{itemize}
\item[(i)]
For each $i =1, \cdots, r$, the following equality holds:
\begin{equation}
\mr{length}(\Lambda^j_i /\Lambda^{j+1}_i) =  m_{i, j+1}.
\end{equation}
\item[(ii)]
The following equality holds:
 \begin{align}
 \mr{deg}(\mcF^{\sharp j}/\mcF^{\sharp j+1}) - \mr{deg}(\mcF^{\sharp}/\mcF^{\sharp 1})& = j \cdot (2g-2+r)  -  \sum_{i =1}^r (m_{i, j+1} -m_{i, 1}).
 \end{align}
\end{itemize}
 \epr
\begin{proof}
First, we shall consider assertion (i).
Let  $v$ and $u_l (t) = \sum_{s =0}^\infty u_{l,s} \cdot t^s$ ($l =1, \cdots, n$) be as in (the proof of) Proposition \ref{10002}.
By a discussion similar to (\ref{177002}), 
the $(n-j-1)\times (n-j-1)$ matrix 
\begin{equation}
(\widehat{\nabla}_{m_{i, n+1-a}}^{b-1}(u_{n+1-a} (t))|_{t=0})_{a,b=1}^{n-j-1}
\end{equation}
 over $k$  turns out to be regular.
Hence, since $k[[t]]^\times = k^\times \oplus t \cdot k[[t]]$, 
the  $(n-j-1)\times (n-j-1)$ matrix 
\begin{equation}
(\widehat{\nabla}_{m_{i, n+1-a}}^{b-1}(u_{n+1-a} (t)))_{a,b=1}^{n-j-1}
\end{equation}
 over $k[[t]]$ is  regular.
This 
means that
the image of  $\{\eta_i (\widehat{\nabla}^{j'}(v)) \}_{0 \leq j' \leq n-j-2}$
    via the projection $k[[t]]^{\oplus n} \migisurj k[[t]]^{\oplus (n-j-1)}$ onto the last $(n-j-1)$ factors 
  forms a basis.
Also, $\widehat{\mcF}^{j+1}$ is generated formally by the elements 
$\{ \widehat{\nabla}^{j'} (v)\}_{0 \leq j' \leq n-j-2}$.
It follows that  any element $w$ of $\widehat{\mcF}/\widehat{\mcF}^{j+1}$ admits a representative $\widetilde{w} \in \widehat{\mcF}$ such that
the last    $(n-j-1)$
factors of $\eta_i(\widetilde{w})$ equal zero.
The image of such a representative via the composite $\widehat{\mcF} \xrightarrow{\eta_i} \widehat{\mcO}^{\oplus n} \migisurj \bigoplus_{j=1}^n \mr{Coker}(\widehat{\nu}_{m_{i, j}}^{\sharp i})$  is annihilated by $t^{m_{i, j+1}}$.
This observation  implies that 
  any element of $\Lambda_i /\Lambda_{i}^{j+1}$ is verified to be  annihilated by $t^{m_{i, j+1}}$.
Since   $\Lambda_i^j/\Lambda^{j+1}_i$
  can be obtained as a quotient of  the line bundle $\mcF^j /\mcF^{j+1}$,  
 $\Lambda_i^j/\Lambda^{j+1}_i$  is   generated by a single section  annihilated by $t^{m_{i, j+1}}$.
Thus, we have
 \begin{equation}\label{18181}
 \mr{length}(\Lambda_i^j/\Lambda^{j+1}_i)  \leq m_{i, j+1}.
 \end{equation}
On the other hand,  (\ref{18006}) induces the following sequence of equalities:
\begin{equation} \label{19191}
\sum_{j=0}^{n-1} \mr{length}(\Lambda_i^j/\Lambda^{j+1}_i) = \mr{length}(\Lambda_i) = \mr{length}(\bigoplus_{l=1}^n \breve{\Lambda}_{i, m_{i,l}}) = \sum_{l=1}^n m_{i, l}.
\end{equation}
By (\ref{18181}) and (\ref{19191}) together, we obtain an equality  $\mr{length}(\Lambda^j_i /\Lambda^{j+1}_i) =  m_{i, j+1}$, as desired.

Next, we shall consider assertion (ii).
Note that $\nu^\sharp$ is compatible with the $k^\mr{log}$-connections $\nabla^\mr{can}_{\mr{Ker}(\nabla)}$ and $\nabla$.
Hence,
the isomorphism  $\mr{KS}_{\mcF^\heartsuit}^j$ 
(cf. (\ref{GL2}))
induces  morphisms
\begin{align} \label{QQ733}
\mr{KS}^j_{\sharp} : \mcF^{\sharp j}/\mcF^{\sharp j+1} \migi \Omega \otimes (\mcF^{\sharp j-1}/\mcF^{\sharp j}),  \hspace{5mm}
\mr{KS}^j_{\Lambda_i} : \Lambda_i^j/\Lambda^{j+1}_i \migi \Omega \otimes (\Lambda_i^{j-1}/\Lambda^{j}_i) 
\end{align}
($i=1, \cdots, r$)
which  fit into the following morphism of short exact sequecnes:
\begin{align} \label{Eh5007}
\vcenter{\xymatrix@C=14pt@R=36pt{
0 \ar[r] & \mcF^{\sharp j}/\mcF^{\sharp j+1}  \ar[rr]^-{\nu^{\sharp j}}\ar[d]^-{\mr{KS}^j_{\sharp}}&& 
 \mcF^j/\mcF^{j+1} \ar[rr]^-{\nu^{\flat j}}\ar[d]^-{\mr{KS}^j_{\mcF^\heartsuit}}&&
 \displaystyle{\bigoplus_{i=1}^r \Lambda^j_i /\Lambda_i^{j+1}}  \ar[r]^-{} \ar[d]^-{\bigoplus_i \mr{KS}^j_{\Lambda_i}} & 0 \\
0 \ar[r]& \Omega\otimes (\mcF^{\sharp j-1}/\mcF^{\sharp j}) \ar[rr]_-{\mr{id} \otimes \nu^{\sharp (j-1)}}
&&  \Omega \otimes (\mcF^{j-1}/\mcF^{j})\ar[rr]_-{\mr{id} \otimes \nu^{\flat (j-1)}}& & \displaystyle{\bigoplus_{i=1}^r\Omega\otimes  (\Lambda^{j-1}_i /\Lambda_i^{j})} \ar[r] & 0. 
}}
\end{align}
By applying  the snake lemma to this morphism, we have
$\mr{Ker} (\bigoplus_{i=1}^r\mr{KS}_{\Lambda_i}^j) \cong \mr{Coker}(\mr{KS}^j_{\sharp})$.
Therefore, the surjectivity of $\bigoplus_{i=1}^r \mr{KS}_{\Lambda_i}^j$
 implies  the following  sequence of equalities:
\begin{align} \label{QQ744}
& \hspace{4mm} \sum_{i=1}^r (\mr{length}(\Lambda^j_i/\Lambda^{j+1}_i) - \mr{length}(\Lambda^{j-1}_i/\Lambda^{j}_i)) \\
& =
\sum_{i=1}^r \mr{length}(\Lambda^j_i/\Lambda^{j+1}_i) -
\sum_{i=1}^r \mr{length}(\Omega \otimes\Lambda^{j-1}_i/\Lambda^{j}_i) \notag \\
&=\mr{length}(\bigoplus_{i=1}^r \Lambda_i^j/\Lambda^{j+1}_i) - \mr{length}(\bigoplus_{i=1}^r \Omega \otimes (\Lambda^{j-1}_i/\Lambda_i^{j})) \notag \\
& = \mr{length}(\mr{Ker}(\bigoplus_{i=1}^r \mr{KS}^j_{\Lambda_i})) \notag\\
& = \mr{length}(\mr{Coker}(\mr{KS}^j_{\sharp})) \notag \\
&=\mr{deg}(\Omega\otimes (\mcF^{\sharp j-1}/\mcF^{\sharp j})) -\mr{deg}(\mcF^{\sharp j}/\mcF^{\sharp j+1}) \notag \\
& =  2g-2+r + \mr{deg}(\mcF^{\sharp j-1}/\mcF^{\sharp j}) - \mr{deg}(\mcF^{\sharp j}/\mcF^{\sharp j+1}). \notag
\end{align}
Using  this equality inductively on $j$, we have
\begin{align}
& \ \mr{deg}(\mcF^{\sharp j}/\mcF^{\sharp j+1}) - \mr{deg}(\mcF^{\sharp}/\mcF^{\sharp 1})  \\
= & \ j \cdot (2g-2+r) - \sum_{i=1}^r (\mr{length}(\Lambda_i^j/\Lambda_i^{j+1}) - \mr{length}(\Lambda_i^0/\Lambda_i^{1})).  \notag
\end{align}
Thus, the required equality   follows from this equality and assertion (i).
\end{proof}

\subsection{The fourth ``$\subseteq$''} \label{y0190}

 In what follows, we prove both  the fourth inclusion relation   ``$\subseteq$'' in (\ref{QQ699}) and the equality
  $\Gamma (X, \JJ\DDD^{\, \sharp <n}_{\BB})=0$ (under the assumption that $g=0$).

 Let  $\BB$ be a line bundle on $X$
 and let $\nabla^\diamondsuit$ be a  $(\mr{GL}_n, 1, \BB)$-oper on $X^\mr{log}/k^\mr{log}$ with ${^p}\psi^{\nabla^\diamondsuit} =0$.
The flat bundle $(\DF_\BB, \nabla^\diamondsuit)$ and the filtration $\{ \DF_\BB^j \}_{j=0}^n$  (cf. (\ref{QW47})) associate  a vector bundle $\DF^{\sharp}_\BB$ and    its filtration $\{ \DF^{\sharp j}_\BB \}_{j=0}^n$ as defined  in (\ref{1070}) and  (\ref{1074})  respectively.
Also, write 
\index{${^\Join}\DD_\BB^{\,<j}$, $\JJ\DD_\BB^{\, <j}$, $\JJ\DD_\BB^{\, \sharp <n}$}
\begin{equation}\label{10041}
\JJ\DDD_{\BB}^{\,\sharp <n}
  := \Omega^\sharp  \otimes (\DF_\BB^{\sharp}/\DF_\BB^{\sharp 1})^\vee \otimes \DF^\sharp_\BB.
\end{equation}
We regard $\JJ\DDD_{\BB}^{\,\sharp <n}$ as an $\mcO_X$-submodule of 
 $\Omega\otimes \mcH om_{\mcO_X} (\DF_\BB^\sharp, \DF_\BB)$
  via the composite injection
\begin{align} \label{10040}
\Omega^\sharp \otimes (\DF_\BB^{\sharp}/\DF_\BB^{\sharp 1})^\vee \otimes \A
& \migiincl \Omega^\sharp\otimes  \DF_\BB^{\sharp \vee} \otimes \DF_\BB^\sharp  \\ 
& \isom  \Omega^\sharp \otimes \mcE nd_{\mcO_X} (\DF_\BB^\sharp) \notag \\
&  \xrightarrow{\mr{inclusion}} \Omega \otimes \mcH om_{\mcO_X} (\DF_\BB^\sharp, \DF_\BB),  \notag 
\end{align}
where the first  arrow is given by forward composition with the quotient   $\DF^\sharp_\BB \migisurj \DF_\BB^{\sharp}/\DF_\BB^{\sharp 1}$. 
In a similar manner, we also regard 
\begin{align}
\JJ\DD_\BB^{\,<n}  \left(:= \mcH om_{\mcO_X}(\Omega^{\otimes (-n)}\otimes\BB, \DF_\BB) \stackrel{}{\cong} \Omega \otimes (\DF_\BB/\DF_\BB^1)^\vee\otimes \DF_\BB\right)
\end{align}
  as an $\mcO_X$-sbubmodule of  $\Omega \otimes \mcH om_{\mcO_X} (\DF^\sharp_\BB, \DF_\BB)$.
\bpr \label{11000}
 We have an inclusion relation
 \begin{equation} \label{10500}
\JJ\DD^{\,<n}_\BB \cap (\Omega^\sharp  \otimes \mcE nd_{\mcO_X}(\DF_\BB^\sharp) ) \subseteq \JJ\DDD_\BB^{\, \sharp <n}
 \end{equation}
 between  $\mcO_X$-submodules of  $\Omega \otimes \mcH om_{\mcO_X} (\DF_\BB^\sharp, \DF_\BB)$.
 In particular, we have an inclusion relation
 \begin{align}
 \varpi (\Gamma (X, \JJ\DD^{\,<n}_\BB \cap \Omega^\sharp  \otimes \mcE nd (\DF_\BB^\sharp))) \subseteq \varpi (\Gamma (X, \JJ\DDD^{\, \sharp <n}_\BB))
 \end{align}
 between (pointed) subsets of $\mr{Def}_{\nabla}$.
 \epr
\begin{proof}
Let $v$ be an arbitrary local section of 
$\JJ\DD^{\,<n}_\BB \cap (\Omega^\sharp  \otimes \mcE nd_{\mcO_X}(\DF_\BB^\sharp) )$.
Since  $v \in \JJ\DD^{\,<n}_\BB$, it may be given by an $\mcO_X$-linear morphism $\DF_\BB/\DF_\BB^1 \migi \Omega \otimes \DF_\BB$.
Denote by $v'$ the composite
\begin{align}
v' : \DF_\BB^\sharp \xrightarrow{\nu^\sharp} \DF_\BB \xrightarrow{\text{quotient}} \DF_\BB/\DF_\BB^1 \xrightarrow{v} \Omega \otimes \DF_\BB.
\end{align}
Then, we have $v' (\DF_\BB^{\sharp 1}) =0$,  so $v'$ factors through the surjection $\DF_\BB^\sharp  \migisurj \DF_\BB^\sharp /\DF_\BB^{\sharp 1}$.
Moreover,   the image $\mr{Im}(v')$ is  contained in 
$\Omega^\sharp \otimes \DF_\BB^\sharp$ because of the assumption   $v \in \Omega^\sharp  \otimes \mcE nd_{\mcO_X}(\DF_\BB^\sharp)$.
It follows that $v'$ induces an $\mcO_X$-linear morphism $v'' : \DF^\sharp / \DF_\BB^{\sharp 1} \migi \Omega^\sharp  \otimes \DF_\BB^\sharp$, which means 
 that  $v \in \JJ\DDD_\BB^{\, \sharp <n}$.
 This completes the proof of the proposition.
 \end{proof}

\bpr \label{12000}
 Assume  further that $g =0$.
 Then,    the following equality holds:
 \begin{equation}\label{10551}
\Gamma (X, \JJ\DDD^{\, \sharp <n}_\BB) =0.
 \end{equation}
 \epr
\begin{proof}
By Proposition  \ref{11003}, (ii),   and the condition $g =0$ (which implies    $\mr{deg}(\Omega^\sharp) = -2$),  the following equality holds:
\begin{equation} \label{11009}
\mr{deg}(\Omega^\sharp \otimes (\DF_\BB^\sharp/\DF_\BB^{\sharp 1})^\vee \otimes  (\DF_\BB^{\sharp j}/\DF_\BB^{\sharp j+1})) = -2 + j \cdot (-2+r) - \sum_{i=1}^r (m_{i,j+1}- m_{i,1}).
\end{equation} 
On the other hand, since the integers $m_{i, 1}, \cdots, m_{i, n}$ are mutually distinct (by Proposition \ref{10002}, (ii)),
we have 
\begin{equation} \label{11010}
\sum_{i=1}^r (m_{i, j+1} -m_{i, 1}) = \sum_{i=1}^r \sum_{l =1}^j (m_{i, l+1} -m_{i, l}) \geq \sum_{i=1}^r\sum_{l =1}^j  1 = j \cdot r. 
\end{equation}
It follows from (\ref{11009}) and (\ref{11010}) that 
\begin{align}
&  \ \ \ \ \mr{deg}(\Omega^\sharp \otimes (\DF_\BB^\sharp/\DF_\BB^{\sharp 1})^\vee \otimes  (\DF_\BB^{\sharp j}/\DF_\BB^{\sharp j+1})) \\
&  \leq -2 + j \cdot (-2 +r) - j\cdot r = -2 -2  j \notag \\
& <0. \notag
\end{align}
Hence,  the equality $\Gamma (X, \Omega^\sharp \otimes (\DF_\BB^\sharp/\DF_\BB^{\sharp 1})^\vee \otimes  (\DF_\BB^{\sharp j}/\DF_\BB^{\sharp j+1})) = 0$ holds.
By induction on $j$, one verifies that $\Gamma (X,  \Omega^\sharp\otimes (\DF_\BB^\sharp/\DF_\BB^{\sharp 1})^\vee \otimes \DF_\BB^{\sharp j}) =0$ for all $j \geq 0$, in particular, 
\begin{equation}
\Gamma (X, \JJ\DDD_\BB^{\, \sharp <n})
 \left(= \Gamma (X, \Omega^\sharp\otimes (\DF_\BB^\sharp/\DF_\BB^{\sharp 1})^\vee \otimes \DF_\BB^{\sharp 0})\right) =0.\end{equation}
This completes the proof of the assertion.
\end{proof}

\subsection{The set ${^p}\mr{Def}^{\mr{Op}}_{\nabla^\diamondsuit}$} \label{QG1402}
Let 
$\bb := (\BB, \nabla_\bb)$ be an $n$-theta characteristic of $X^\mr{log}/k^\mr{log}$ with ${^p}\psi^{\nabla_\bb} =0$.
(Notice that such an $n$-theta characteristic always exists because of the result in \S\,\ref{QQ183}.)
Also, let $\nabla^\diamondsuit$ be a dormant  $(\mr{GL}_n, \bb)$-oper on $\msX$.  
We shall write
\index{${^p}\mr{Def}_{(\mcF, \nabla)}$, ${^p}\mr{Def}_{\nabla}$, ${^p}\mr{Def}^\sharp_{\nabla}$, ${^p}\mr{Def}^{\mr{Op}}_{\nabla^\diamondsuit}$}
\begin{align} \label{QQ751}
{^p}\mr{Def}^{\mr{Op}}_{\nabla^\diamondsuit}
\end{align}
for the set of isomorphism classes of dormant  $(\mr{GL}_n, \bb)$-opers $\nabla_{\epsilon'}^\diamondsuit$ on $\msX_\epsilon$ together with an isomorphism   $(\mr{in}_X^* (\nabla_{\epsilon'}^\diamondsuit))^{\Rightarrow \heartsuit} \isom  \nabla^{\diamondsuit \Rightarrow \heartsuit}$.
Then, by means of the following assertion, 
${^p}\mr{Def}^{\mr{Op}}_{\nabla^\diamondsuit}$ can be regarded as a subset of ${^p}\mr{Def}_{\nabla^\diamondsuit}$.
The following result is an intermediate step for proving the first ``$\subseteq$'' in (\ref{QQ699}).

\bpr \label{QG51}
\begin{itemize}
\item[(i)]
The map
\begin{align} \label{QG52}
{^p}\mr{Def}^{\mr{Op}}_{\nabla^\diamondsuit} \migi {^p}\mr{Def}_{\nabla^\diamondsuit}
\end{align}
given by  assigning  $\nabla_{\epsilon'}^\diamondsuit \mapsto (\DF_{\mr{pr}_X^*(\BB)}, \nabla_{\epsilon'}^\diamondsuit)$ is injective.
(By using this injection, we regard ${^p}\mr{Def}^{\mr{Op}}_{\nabla^\diamondsuit}$ as a (pointed) subset  of 
${^p}\mr{Def}_{\nabla^\diamondsuit}$.)
\item[(ii)]
Suppose further that $\nabla^\diamondsuit = D^{\clubsuit \Rightarrow \diamondsuit}$ for an $(n,  \bb)$-projective connection $D^\clubsuit$ on $X^\mr{log}/k^\mr{log}$.
Then, we have an inclusion relation
\begin{align}
{^p}\mr{Def}^{\mr{Op}}_{\nabla^\diamondsuit} \subseteq \varpi (\Gamma (X, \JJ\DDD^{\, <n}_\BB)) \cap   {^p}\mr{Def}_{\nabla^\diamondsuit}.
\end{align} 
\end{itemize}
\epr
\begin{proof}
First, we shall consider assertion (i).
 Denote by 
$\msE^\spadesuit :=(\DE_B, \nabla)$ the dormant $\mfs \mfl_n$-oper corresponding to $\nabla^\diamondsuit$ via the sequence of isomorphisms resulting from Theorem \ref{y0129}.
By using that theorem again,
we see that ${^p}\mr{Def}^{\mr{Op}}_{\nabla^\diamondsuit}$
 corresponds bijectively to the 
 space of first-order deformations of   $\msE^\spadesuit$ preserving  the dormancy condition.
By  (\ref{isom9}) and Propositions  \ref{y0143}, (ii), and  \ref{y0148}, this set  may be regarded as a subset of $\mbH^1 (X, \mcK^\bullet [\nabla^\mr{ad}])$.
On the other hand,  
let us consider  the  map $\alpha : {^p}\mr{Def}_{\nabla^\diamondsuit} \migi \mbH^1 (X, \mcK^\bullet [\nabla^{\mr{ad}}])$
obtained by composing  the canonical bijection $\left({^p}\mr{Def}_{\nabla^\diamondsuit} \subseteq\right)\mr{Def}_{(\DF_\BB, \nabla^\diamondsuit)} \isom \mbH^1 (X, \mcK^\bullet [\nabla^{\diamondsuit \mr{ad}}])$ (cf. (\ref{10602}))
and the map $\mbH^1 (X, \mcK^\bullet [\nabla^{\diamondsuit \mr{ad}}]) \migi \mbH^1 (X, \mcK^\bullet [\nabla^{\mr{ad}}])$
 arising  from  the quotient  $\mcE nd (\DF_\BB) \migisurj \overline{\mcE nd} (\DF_\BB) \left(\cong (\mfs \mfl_n)_{\DE_B} \right)$ (cf. (\ref{QW61})).
One verifies that the inclusion ${^p}\mr{Def}^{\mr{Op}}_{\nabla^\diamondsuit} \migiincl \mbH^1 (X, \mcK^\bullet [\nabla^{\mr{ad}}])$ factors through  $\alpha$ and
the resulting map ${^p}\mr{Def}^{\mr{Op}}_{\nabla^\diamondsuit}  \migi {^p}\mr{Def}_{\nabla^\diamondsuit}$ coincides with (\ref{QG52}).
This implies the injectivity of (\ref{QG52}).

Next, we shall consider assertion (ii).
Recall  (cf. Theorem \ref{y0129}) 
that $\mfO \mfp^\diamondsuit_{n, 1, \bb}$ admits a structure of 
  affine space on  modeled on $\Gamma (X, \DV_\bb)$ with respect to the $\Gamma (X, \DV_\bb)$-action  given by restricting  (\ref{QS9666}).
Since the isomorphism 
$ {^\mcD}\hspace{-1.5mm}\nabla^{\diamondsuit <n}$ (cf. (\ref{GL25})) coincides with  the identity morphism $\mr{id}_{\DF_\BB}$ by assumption, 
 each element of  ${^p}\mr{Def}^{\mr{Op}}_{\nabla^\diamondsuit}$  is represented by a first-order  deformation of $\nabla^\diamondsuit$  described as $\mr{pr}_X^*(\nabla^\diamondsuit) + \epsilon \cdot {^\circlearrowright \beta}^n (R)$ for some $R \in  \Gamma (U,\DV_\bb)$.
 The sheaf  $\JJ\DDD^{\, <n}_\BB$ contains the image of 
  $\DV_\bb$
  via ${^\circlearrowright \beta}^n$, so 
 the set 
 ${^p}\mr{Def}^{\mr{Op}}_{\nabla^\diamondsuit}$ is contained in $\varpi (\Gamma (X, \JJ\DDD^{\, <n}_\BB))$.
 This completes the proof of assertion (ii).
\end{proof}

\section{Dormant opers on the $3$-pointed projective line}\label{y0194}

In Proposition \ref{002} described later, we will prove  an inclusion relation between two deformation spaces.
By combining this result and Proposition \ref{QG51}, (ii), proved already,
 we conclude  the first inclusion relation  ``$\subseteq$'' appearing in (\ref{QQ699}) under a certain assumption.

Let $\msP := (\mbP^1/k, \{ [0], [1], [\infty] \})$ be as in (\ref{1051}).
Since the underlying curve $\mbP^1$ of $\msP$, i.e., the projective line over $k$, may be defined over $\mbF_p$, the Frobenius twist $\mbP^{1(1)}$   is isomorphic to  $\mbP^1$ itself.
Here, recall the Birkhoff-Grothendieck theorem,  asserting   that any vector bundle $\mcV$ on $\mbP^{1}$ of rank $n$ is isomorphic to a direct sum of $n$ line bundles:
\begin{equation}
\mcV \cong \bigoplus_{j= 1}^n \mcO_{\mbP^1} (w_j),
\end{equation}
where
$w_1, \cdots, w_n$ are integers with $w_1 \leq w_2 \leq \cdots \leq w_n$.
The ordered set of integers $(w_j)_{j=1}^n \in \mbZ^n$  depends only on the isomorphism class of  $\mcV$.
In particular, we have $\mr{deg}(\mcV) =\sum_{j=1}^n w_j$.

\bde \label{QQ580}
In the above situation, we shall say that $\mcV$ is {\bf of type $(w_j)_{j=1}^n$}. \index{vector bundle@--- of type $(w_j)_{j=1}^n$}
\ede

\ble\label{311789}
For $\square = 1,2$, 
let $\mcV_\square$ be 
 a rank $n$ vector bundle   on $\mbP^{1}$ of type $(w_{\square, j})_{j =1}^n$.
 If there exists  an $\mcO_{\mbP^{1}}$-linear   injection $\mcV_1 \migiincl \mcV_2$, then
 the inequality   $w_{1, j} \leq w_{2, j}$ holds  for any $j \in \{ 1, \cdots, n\}$. 
\ele
\begin{proof}
The assertion follows immediately from, e.g.,  an inductive argument with respect to $n$.
\end{proof}
The above lemma deduces the following lemma, which  will be used in the proof of Proposition \ref{001} described later.

\ble \label{6234} 
 Let $s$ be an integer and $\{ \mcV_l\}_{l\in \mbZ_{\geq 0}}$  a set of rank $n$ vector bundles on $\mbP^{1}$ such that each $\mcV_l$ is of degree $s + l$ and type $(w_{l, j})_{j=1}^n$ (hence $\sum_{j=1}^n w_{l, j} = s+l$) with   $w_{l, n} - w_{l, 1} \leq 2$.
Also, assume that we are given  a sequence of $\mcO_{\mbP^1}$-linear injections 
 \begin{equation} \mcV_0 \migiincl \mcV_1 \migiincl \mcV_2 \migiincl \cdots. \end{equation}
 Then, there exists $l_0 \in \mbZ_{\geq 0}$ such that $\mcV_{l_0}$ is of type $(w_{l_0, j})_{j=1}^n$ with   $w_{l_0, n} - w_{l_0, 1} \leq 1$.
 \ele
\begin{proof}
It suffices to consider the case where $w_{0, n} - w_{0, 1} = 2$.
By the assumption and Lemma \ref{311789}, the integer 
\begin{align}
l_1 : = \mr{min} \left\{ l \in \mbZ_{\geq 0} \,  | \,  w_{l, n} \neq  w_{0,n}\right\}
\end{align}
is well-defined and satisfies that 
 $l_1 >0$ and $w_{l_1-1, n} = w_{0, n}$.
Since  $w_{0, n} - w_{0, 1} = 2$ and $w_{l_1, n} - w_{l_1, 1} \leq 2$, 
the inequality $w_{l_1-1, 1} \geq  w_{0, 1} +1$ must hold because of  Lemma \ref{311789} again.
It follows that
\begin{equation}
w_{l_1-1, n} - w_{l_1-1, 1} \leq  w_{0, n} - (w_{0, 1} +1) = 1.
\end{equation}
Thus,  the integer $l_0 := l_1-1$ turns out to be  a desired integer.
\end{proof}

Hereafter, let us fix  an $n$-theta characteristic $\bb := (\BB, \nabla_\bb)$  of   $\mbP^{1\mr{log}}/k^\mr{log}$ with  ${^p \psi}^{\nabla_\bb}=0$.

\bpr \label{500987}
Let $\nabla$
 be a dormant $(\mr{GL}_n, \bb)$-oper on $\msP$.
Assume that the rank $n$ vector bundle $\mr{Ker}(\nabla)$  on $\mbP^{1(1)}$  is of type $(w_l)_{l=1}^n$, i.e., 
  decomposes into a direct sum of  $n$ line bundles 
$\mr{Ker}(\nabla) \cong \bigoplus_{l = 1}^n \mcO_{\mbP^{1(1)}} (w_l)$.
Then, the  inequality $w_n -w_1 \leq2$ holds.
 \epr
\begin{proof}
By assumption, the rank $n$ vector bundle $\DF_\BB^\sharp := F_{\mbP^1/k}^*(\mr{Ker}(\nabla))$ on $\mbP^1$ decomposes into the direct sum
\begin{equation} \label{61298}
\DF_\BB^\sharp \cong  \bigoplus_{l = 1}^n \mcO_{\mbP^{1}} (p \cdot w_l).
\end{equation}
It follows from Propositions \ref{10002}, (ii),  and  \ref{11003}, (ii)  (in the case of  $(g,r) = (0,3)$),  that
\begin{equation} \label{429996}
\mr{deg}(\DF_\BB^{\sharp n-1}/\DF_\BB^{\sharp n}) < \mr{deg}(\DF_\BB^{\sharp n-2}/\DF_\BB^{\sharp n-1}) < \cdots < \mr{deg}(\DF_\BB^{\sharp 0}/\DF_\BB^{\sharp 1}).
\end{equation}
Here, we shall write 
\begin{equation}
\xi' : \mcO_{\mbP^{1}} (p \cdot w_n) \migiincl \DF_\BB^\sharp \  \left(\text{resp.},  \ \xi'' : \DF_\BB^\sharp \migisurj \mcO_{\mbP^{1}} (p \cdot w_1) \right)
\end{equation}
 for the inclusion into the $n$-th factor (resp., the projection onto the $1$-st factor) with respect to the decomposition     (\ref{61298}).
Also,  write
\begin{align}
  n' := \mr{max} \big\{ j \, \big| \,  \mr{Im}(\xi') \subseteq \DF_\BB^{\sharp j} \big\}  \   \left(\text{resp.},  \   n'' := \mr{max} \big\{ j \, \big| \,  \xi''(\DF_\BB^{\sharp j} ) \neq 0\big\} \right). 
 \end{align} 
Then, $\xi'$ (resp., $\xi''$) induces a {\it nonzero} morphism between line bundles
\begin{align}
 \overline{\xi}' : \mcO_{\mbP^{1}} (p \cdot w_n) \migi 
  \DF_\BB^{\sharp n'}\hspace{-1mm}/\DF_\BB^{\sharp n'+1}   \ 
 \left(\text{resp.,} \   \overline{\xi}'' :
  \DF_\BB^{\sharp n''}\hspace{-1mm}/\DF_\BB^{\sharp n''+1} \migi \mcO_{\mbP^{1}} (p \cdot w_1)\right). \hspace{-3mm}
\end{align}
 In particular,  $\overline{\xi}'$ and $ \overline{\xi}''$ are injective,  so we have
\begin{align} \label{722465}
p \cdot w_n &= \mr{deg} (\mcO_{\mbP^{1}} (p \cdot w_n)) \leq \mr{deg}(\DF_\BB^{\sharp n'}\hspace{-1mm}/\DF_\BB^{\sharp n'+1}) \leq  \mr{deg}(\DF_\BB^{\sharp}\hspace{0mm}/\DF_\BB^{\sharp 1}), \\
p \cdot w_1 &=  \mr{deg} (\mcO_{\mbP^{1}} (p \cdot w_1)) \geq \mr{deg}(\DF_\BB^{\sharp n''}\hspace{-1mm}/\DF_\BB^{\sharp n''+1}) \geq \mr{deg}(\DF_\BB^{\sharp n-1}\hspace{0mm}/\DF_\BB^{\sharp n}), \notag
\end{align}
where the   last inequalities in the respective  sequences    follow from (\ref{429996}).
But,  Proposition  \ref{11003}, (ii),  implies 
\begin{align} \label{1}
& \hspace{5mm}\mr{deg}(\DF_\BB^{\sharp}/\DF_\BB^{\sharp 1}) - \mr{deg}(\DF_\BB^{\sharp n-1}/\DF_\BB^{\sharp n}) \\
 &= - (n-1)  (2 \cdot 0 - 2 + 3) + \sum_{\sigma \in \{ [0], [1], [\infty]\}} (m_{\sigma,n}-m_{\sigma,1})  \notag \\
 &< -n+1 + 3  p.
 \notag
\end{align}
where $(-m_{\sigma, l})_{l=1}^n := \mr{Res}_\sigma (\nabla)$.
By combining  (\ref{722465}) and (\ref{1}), we obtain
\begin{equation}
w_n -w_1 \leq \frac{1}{p} \cdot (\mr{deg}(\DF_\BB^{\sharp}/\DF_\BB^{\sharp 1}) - \mr{deg}(\DF_\BB^{\sharp n-1}/\DF_\BB^{\sharp n})) < 3- \frac{n-1}{p}.
\end{equation}
Hence, since both $w_n$ and $w_1$ are integers,  the equality  $w_n -w_1 \leq 2$ holds.
This completes the proof of the assertion. 
\end{proof}

\bpr \label{001}
Let $\nabla^\diamondsuit$ be a dormant $(\mr{GL}_n, \bb)$-oper on $\msP$.
 Then, there exists a flat  line bundle $\msL := (\mcL, \nabla_\mcL)$ on $\mbP^{1\mr{log}}/k^\mr{log}$ with ${^p}\psi^{\nabla_\mcL}=0$ such that 
 the rank $n$ vector bundle  $\mr{Ker}(\nabla^\diamondsuit_{\otimes \msL})$  
 on $\mbP^{1(1)}$ (cf. (\ref{QQ1131}))
   is of type $(w_l)_{l=1}^n$ with $w_n -w_1 \leq 1$.
 \epr
\begin{proof}
For convenience, we write  $\bb_0 := \bb$ and  $\nabla_{0}^\diamondsuit := \nabla^\diamondsuit$.
By means of   this pair $(\bb_0, \nabla_{0}^\diamondsuit)$,
 we shall construct a new
 $n$-theta characteristic  $\bb_1$, 
  a dormant $(\mr{GL}_n, \bb_1)$-oper $\nabla_{1}^\diamondsuit$ on $\msP$, and an  injection $\iota_1 : \mr{Ker}(\nabla_0^\diamondsuit) \migiincl \mr{Ker}(\nabla_1^\diamondsuit)$, as follows. 
First,  by letting  $( -m_l )_{l=1}^n := \mr{Res}_{[0]}(\nabla_0^\diamondsuit)$, 
 we  obtain  the flat line bundle 
\begin{equation}
\msL_1 = (\mcL_1, \nabla_{\mcL, 1}) := (\mcO_{\mbP^1}((p-m_{n})\cdot [0]), \nabla_{(p-m_{n})\cdot [0]})
\end{equation}
on $\mbP^{1 \mr{log}}/k^\mr{log}$, which has vanishing $p$-curvature.
Then, 
twisting   $\bb_0$ by  $\msL_1$ induces
 an  $n$-theta characteristic 
\begin{equation}
\bb_1: = \bb \otimes\msL_1  = (\BB \otimes \mcL_1, \nabla_\bb \otimes \nabla_{\mcL, 1}^{\otimes n}), \end{equation}
 which  satisfies ${^p}\psi^{\nabla_0 \otimes \nabla_{\mcL, 1}^{\otimes n}}=0$.
The $k^\mr{log}$-log connection
 \begin{equation}
 \nabla_{1}^\diamondsuit := \nabla_{0, \otimes \msL_1}^\diamondsuit 
 \end{equation}
 forms a dormant $(\mr{GL}_n, \bb_1)$-oper on $\msP$.
 It follows from the definition of this connection that 
 the injection
 $\DF_\BB \migiincl \DF_{\BB \otimes \mcL_1}$ arising from the inclusion 
 $\mcO_X \migiincl \mcL_1 \left(= \mcO_{\mbP^1}((p-m_n)\cdot [0]) \right)$
 is compatible with the $k^\mr{log}$-connections $\nabla_0^\diamondsuit$, $\nabla_1^\diamondsuit$.
 Hence, it restricts to
   an $\mcO_{\mbP^{(1)}}$-linear  injection
 \begin{equation}
 \iota_1 : \mr{Ker}(\nabla^\diamondsuit_{0}) \migiincl \mr{Ker}(\nabla^\diamondsuit_{1})
 \end{equation}
between rank $n$ vector bundles on $\mbP^{1(1)}$.
In this way, we have obtained  the  triple of data  $(\bb_1, \nabla_1^\diamondsuit, \iota_1)$.

In what follows, we shall compare the degrees of the vector bundles $\mr{Ker}(\nabla^\diamondsuit_{0})$ and  $\mr{Ker}(\nabla^\diamondsuit_{1})$.
 Since the restriction of the flat line bundle $\msL_1$ to  $\mbA^1 := \mbP^1 \setminus \{[0]\} \left(\subseteq \mbP^1 \right)$  becomes  trivial,
    the  natural inclusion
\begin{equation} \label{15}
\mr{Ker}(\nabla^\diamondsuit_{0}) \otimes \mr{Ker}(\nabla_{(p-m_{n})\cdot [0]}) \migiincl \mr{Ker}(\nabla^\diamondsuit_{1})
\end{equation}
 induces an isomorphism
 \begin{equation} \label{15a}
\mr{Ker}(\nabla^\diamondsuit_{0})  |_{\mbA^1} \isom \mr{Ker}(\nabla^\diamondsuit_{1}) |_{\mbA^1}
\end{equation}
when restricted to $\mbA^1$.
 Now, let us fix a local function $t \in \mcO_{\mbP^1}$ defining $[0]$, and choose    a decomposition
 \begin{equation} \label{11}
  (\DF_\BB, \nabla^\diamondsuit_{0})^\wedge \isom \bigoplus_{l=1}^n (\widehat{\mcO}, \widehat{\nabla}_{m_l})
 \end{equation}
   as in  (\ref{decom9}).
By taking the respective sheaves of horizontal sections, one obtains from this decomposition  an isomorphism of $\mcO_{\mr{Spf}(k[[t^p]])}$-modules
\begin{equation} \label{16}
 \mr{Ker}(\nabla^\diamondsuit_{0})^\wedge \isom \bigoplus_{l=1}^n t^{m_l} \cdot \mcO_{\mr{Spf}(k[[t^p]])},
\end{equation}
where $\mr{Ker}(\nabla^\diamondsuit_{0})^\wedge$ denotes the $t^p$-adic completion of $\mr{Ker}(\nabla^\diamondsuit_{0})$.
Also,  tensoring  (\ref{11})   with $\msL_1^\wedge$ gives an isomorphism
  \begin{equation} \label{12}
 (\DF_{\BB \otimes \mcL_1}, \nabla_1^\diamondsuit)^\wedge \isom \bigoplus_{l=1}^n (\widehat{\mcO}, \widehat{\nabla}_{p-m_n+ m_l}).
 \end{equation}
Since  $\mr{Ker}(\widehat{\nabla}_p) = \mcO_{\mr{Spf}(k[[t^p]])}$ by Lemma \ref{55232}, (ii),
 (\ref{12})  induces an isomorphism between the respective sheaves of horizontal sections
\begin{equation} \label{17}
  \mr{Ker}(\nabla^\diamondsuit_{1})^\wedge \isom \bigoplus_{l=1}^n t^{p-m_n+ m'_l} \cdot \mcO_{\mr{Spf}(k[[t^p]])},
\end{equation}
where $m_l' := m_l$ if $l = 1, \cdots, n-1$ and $m_l' :=m_l-p$ if $l =n$.
Under   the identifications given by  (\ref{16}), (\ref{17}), and  $ \mr{Ker}(\nabla_{(p-m_{n})\cdot [0]})^\wedge = t^{p-m_n}\cdot  \mcO_{\mr{Spf}(k[[t^p]])}$,
the $t^p$-adic completion of   (\ref{15})
may be identified with the morphism
\begin{align}
 \big( \bigoplus_{l=1}^n t^{m_l} \cdot \mcO_{\mr{Spf}(k[[t^p]])}\big)\widehat{\otimes}  t^{p-m_n}\cdot  \mcO_{\mr{Spf}(k[[t^p]])} & \migi \bigoplus_{l=1}^n t^{p-m_n+ m'_l} \cdot \mcO_{\mr{Spf}(k[[t^p]])} 
\end{align}
given by $(\sum_{l=1}^n t^{m_l} \cdot a_{l}) \otimes (t^{p-m_m} \cdot b) \mapsto  \sum_{l=1}^n t^{p-m_n+m_l} \cdot a_{l} \cdot b$ ($a_1, \cdots, a_n, b \in  \mcO_{\mr{Spf}(k[[t^p]])}$).
The cokernel of this morphism is of length $1$, so
it follows from  this fact and  (\ref{15a}) together  that
\begin{equation}
\label{36}
\mr{deg}(\mr{Ker}(\nabla^\diamondsuit_{1})) = \mr{deg}(\mr{Ker}(\nabla^\diamondsuit_{0})) +1.
\end{equation}

By applying the above construction of $(\bb_1, \nabla_1^\diamondsuit, \iota_1)$  using $(\bb_0, \nabla_0^\diamondsuit)$  successively,
 we obtain a collection  of  pairs 
 \begin{align}
 \left\{ (\bb_m, \nabla^\diamondsuit_{m}) \right\}_{m \geq 0},
 \end{align}
  where each $\bb_m$ denotes an $n$-theta characteristic of $\mbP^{1\mr{log}}/k^\mr{log}$ (with vanishing $p$-curvature) and $\nabla^\diamondsuit_{m}$ denotes  a dormant $(\mr{GL}_n, \bb_m)$-oper on $\msP$.
Moreover,   the equality $\mr{deg}(\mr{Ker}(\nabla_{m+1}^\diamondsuit)) = \mr{deg}(\mr{Ker}(\nabla_{m}^\diamondsuit)) +1$ holds for any $m \geq 0$, and 
there exists a sequence of $\mcO_{\mbP^{1(1)}/k}$-linear injections
\begin{equation}
\mr{Ker}(\nabla^\diamondsuit_{0}) \stackrel{\iota_1}{\migiincl} \mr{Ker}(\nabla^\diamondsuit_{1})  \stackrel{\iota_2}{\migiincl} \mr{Ker}(\nabla^\diamondsuit_{2}) \stackrel{\iota_3}{\migiincl} \cdots.
\end{equation}
Hence, by  Lemma \ref{6234} and Proposition \ref{500987},  we can find  $m\geq 0$ such that 
$\mr{Ker}(\nabla_{m}^\diamondsuit)$ 
 is of type $(w_l)_{l=1}^n$  with  $w_n -w_1 \leq 1$.
It follows that 
 $\bb_m /\bb_{0}$  (cf. (\ref{QQ1138})) 
 specifies  a desired flat line bundle.
This completes the proof of the proposition.
\end{proof}

\bpr \label{002}
Let $\nabla^\diamondsuit$ be a dormant $(\mr{GL}_n, \bb)$-oper on $\msP$.
Then, there exists a flat line bundle $\msL$ on $\mbP^{1 \mr{log}}/k^\mr{log}$ with vanishing $p$-curvature for which  we have the inclusion relation 
\begin{align}
{^p}\mr{Def}^{\mr{Op}}_{\nabla^\diamondsuit_{\otimes \msL}} \subseteq {^p}\mr{Def}_{\nabla^\diamondsuit_{\otimes \msL}}^\sharp
\end{align}
between subsets of ${^p}\mr{Def}_{\nabla^\diamondsuit_{\otimes \msL}}$.
 \epr
\begin{proof}
To prove the assertion, 
we shall  take ``$\msL$'' to be a flat line bundle  resulting from Proposition \ref{001}, i.e., satisfying that the rank $n$ vector bundle $\mr{Ker}(\nabla^\diamondsuit_{\otimes \msL})$ on $\mbP^{1(1)}$ is of type $(w_l)_{l=1}^n$ with   $w_n-w_1 \leq 1$.
Since $\mr{Ker}(\nabla_{\otimes \msL}^\diamondsuit) \cong \bigoplus_{l=1}^n \mcO_{\mbP^{1(1)}_k} (w_l)$,
  we obtain the following sequence of isomorphisms: 
\begin{align} \label{QHH567}
H^1(\mbP^{1(1)}, \mcE nd (\mr{Ker}(\nabla^\diamondsuit_{\otimes \msL}))) & \isom
H^1(\mbP^{1(1)}, \mcE nd (\bigoplus_{l=1}^n \mcO_{\mbP^{1(1)}} (w_l)))   \\
& \isom H^1(\mbP^{1(1)}, \bigoplus_{l, l' =1}^n\mcO_{\mbP^{1(1)}} (w_l -w_{l'}) )  \notag \\
&\isom \bigoplus_{l, l' =1}^n H^0 (\mbP^{1(1)}, \mcO_{\mbP^{1(1)}} (-w_l +w_{l'}-2))^\vee  \notag\\
& \isom 0, \notag
\end{align} 
where the third isomorphism arises from Serre duality and $\Omega_{\mbP^{1(1)}/k} \cong \mcO_{\mbP^{1(1)}}(-2)$ together, and  the fourth isomorphism follows from the fact that   $\mr{max} \left\{ w_{l'} -w_{l} \, | \, 1  \leq l, l' \leq n \right\}$ $= w_n-w_1 \leq 1$.
Recall    that there exists  a canonical  bijection  between  the underlying  set  of $H^1(\mbP^{1(1)}, \mcE nd (\mr{Ker}(\nabla_{\otimes \msL}^\diamondsuit)))$ and the space of 
  first-order deformations  of the vector bundle $\mr{Ker}(\nabla_{\otimes \msL}^\diamondsuit)$ (cf. (\ref{QW6077})).
Hence,  by (\ref{QHH567}),
any first-order deformation of $\mr{Ker}(\nabla^\diamondsuit_{\otimes \msL})$ is trivial.
This implies that, for each 
deformation $\nabla^\diamondsuit_{\otimes \msL, \epsilon'}$ 
 classified by
${^p}\mr{Def}^{\mr{Op}}_{\nabla^\diamondsuit_{\otimes \msL}}$,
the deformation  $F_{\mbP^1/k}^*(\mr{Ker}(\nabla^\diamondsuit_{\otimes \msL, \epsilon'}))$ of $\DF^\sharp_{\BB \otimes \mcL}$ must be trivial.
This completes the proof of the assertion.
\end{proof}

As a consequence of the results obtained so far, we obtain the following assertion:

\bco \label{QG55}
For each   dormant $(\mfs \mfl_n, \hslash)$-oper $\msE^\spadesuit$ on $\msP$ with $\hslash \in k^\times$, 
there are no nontrivial first-order deformations of $\msE^\spadesuit$.
In other words, the  finite $k$-scheme $\mfO \mfp_{\mfs \mfl_n,  \hslash, \msP}^\ZZZ \left(= \mfO \mfp_{\mfs \mfl_n,  \hslash, 0, 3}^\ZZZ \right)$ is a (possibly empty)  disjoint union of finite copies of $\mr{Spec}(k)$.
\eco 
\begin{proof}
By (\ref{QS1100}),  it suffices to consider the case of $\hslash =1$.
Let us take  an $n$-theta characteristic $\bb := (\BB, \nabla_\bb)$ of $\mbP^{1 \mr{log}}/k^\mr{log}$,   and denote by $\nabla^\diamondsuit$  the dormant $(\mr{GL}_n, \bb)$-oper corresponding to   $\msE^\spadesuit$ 
via
${^p} \Lambda_{\msP, \bb}^{\diamondsuit \Rightarrow \spadesuit}$ 
  resulting from Corollary \ref{1090}, (i).
The space of first-order deformations of $\msE^\spadesuit$ (preserving the dormancy condition) may be identified with ${^p}\mr{Def}^{\mr{Op}}_{\nabla^\diamondsuit}$.
By Theorem \ref{y0129}, we may assume  that
$\nabla^\diamondsuit = D^{\clubsuit \Rightarrow \diamondsuit}$ for an  
$(n, \hslash)$-projective connection $D^\clubsuit$ on $X^\mr{log}/k^\mr{log}$. 
Hence, it follows from Proposition \ref{QG51}, (ii),  that
\begin{align} \label{QHH588}
{^p}\mr{Def}^{\mr{Op}}_{\nabla^\diamondsuit} \subseteq \varpi (\Gamma (\mbP^1, \JJ\DDD^{\, <n}_\BB)) \cap   {^p}\mr{Def}_{\nabla^\diamondsuit}.
\end{align}
Moreover, by  Proposition \ref{002}, 
we may assume, after possibly twisting $\bb$ by some flat line bundle with vanishing $p$-curvature, that
${^p}\mr{Def}^{\mr{Op}}_{\nabla^\diamondsuit} \subseteq {^p}\mr{Def}_{\nabla^\diamondsuit}^\sharp$.
(Because of (\ref{QW31}), the assumption that $\nabla^\diamondsuit$ comes from an $(n, \bb)$-projective connection is  still satisfied after twisting $\bb$.
In particular,   (\ref{QHH588}) still holds.)
Thus, we have ${^p}\mr{Def}^{\mr{Op}}_{\nabla^\diamondsuit}  \subseteq \varpi (\Gamma (\mbP^1, \JJ\DDD^{\, <n}_\BB)) \cap {^p}\mr{Def}_{\nabla^\diamondsuit}^\sharp$.
Also, it follows from Propositions \ref{QG71}, \ref{QG70},  and  \ref{11000},
 that 
 $\varpi (\Gamma (\mbP^1, \JJ\DDD^{\, <n}_\BB)) \cap {^p}\mr{Def}_{\nabla^\diamondsuit}^\sharp$
 is contained in $\Gamma (\mbP^1, \JJ\DDD^{\, \sharp <n}_\BB)$.
 By combining these results, we have ${^p}\mr{Def}^{\mr{Op}}_{\nabla^\diamondsuit} \subseteq \Gamma (\mbP^1, \JJ\DDD^{\, \sharp <n}_\BB)$.
On the other hand, we have already obtained the equality  $\Gamma (\mbP^1, \JJ\DDD^{\, \sharp <n}_\BB) =0$ in Proposition \ref{12000}.
It follows that ${^p}\mr{Def}^{\mr{Op}}_{\nabla^\diamondsuit}$ consists exactly of one element.
By Corollary \ref{1090}, (i), applied to the case where $\msX = \msP_\epsilon$, 
this is equivalent to the fact that there are no nontrivial first-order deformations of $\msE^\spadesuit$.
This completes the proof of the assertion.
\end{proof}

\section{The generic \'{e}taleness of the moduli stack of dormant opers}\label{y0250}

\bpr \label{355495}
Let $\msX := (X/k, \{ \sigma_i\}_{i=1}^r)$ be   a totally degenerate curve over $k$ (cf. Definition \ref{y0178}) and  let $\hslash \in k^\times$, $\rho \in \mfc_{\mfs \mfl_n}^{\times r}(\mbF_p)$.
Then, any dormant $(\mfs \mfl_n, \hslash)$-oper   on
$\msX$ does not have any nontrivial first-order deformation.
Equivalently, the finite $k$-scheme $\mfO \mfp^\ZZZ_{\mfs  \mfl_n, \hslash, \hslash \star \rho, \msX}$  is a (possibly empty) disjoint union of finite copies of $\mr{Spec}(k)$. 
 \epr
\begin{proof}
Let us  apply  Theorem \ref{y0176} to  the clutching data $\msG$
 specifying   $\msX$  in the manner of the discussion in Remark \ref{RRR001}.
Then, we see that each component of $\mfO \mfp^\ZZZ_{\mfs  \mfl_n, \hslash, \hslash \star \rho, \msX}$ (i.e., the fiber of the projection $\mfO \mfp^\ZZZ_{\mfs \mfl_n,  \hslash, \hslash \star \rho, g,r} \migi \overline{\mfM}_{g,r}$ over the point classifying $\msX$) is isomorphic to a product of $\mfO \mfp^\ZZZ_{\mfs \mfl_n, \hslash, \hslash \star \varrho, \msP}$'s  for various radii $\varrho$.
Hence,  the problem is reduced to  the case where $\msX = \msP$, but   it  was already proved in  Corollary \ref{QG55}.
\end{proof}

By applying the above proposition,  we obtain the following theorem, asserting   the generic \'{e}taleness of  $\mfO \mfp_{\mfg, \hslash, \hslash \star \rho, g,r}^\ZZZ/\overline{\mfM}_{g,r}$ for  $\mfg = \mfs \mfl_n,  \mfs \mfo_{2l+1}, \mfs \mfp_{2m}$.
The case of $(\mfg, \hslash) = (\mfs \mfl_2, 1)$ was already proved by S. Mochizuki (cf. Theorem \ref{QS1067}).

\bt[cf. Theorem \ref{QG3457}] 
\label{3778567d} 
Let  $k$ be a perfect of characteristic $p$,  $(g,r)$  a pair of nonnegative integers with $2g-2+r>0$,  and  $\hslash$ an element of $k^\times$.
Also, let $\mfg$ be a Lie algebra which is 
 either $\mfs \mfl_n$ with $2  n<p$,  $\mfs \mfo_{2l+1}$ with   $4  l+2 <  p$,  or $\mfs \mfp_{2m}$ with    $4  m < p$.
Then,  for each $\rho \in \mfc^{\times r}_{\mfg} (\mbF_p)$, 
 $\mfO \mfp_{\mfg, \hslash, \hslash \star \rho, g,r}^\ZZZ$  is
\'{e}tale over the points of $\overline{\mfM}_{g,r}$ classifying totally degenerate curves.
In particular, (because of the irreducibility of  $\overline{\mfM}_{g,r}$ and the  finiteness of  $\mfO \mfp_{\mfg, \hslash, \hslash \star \rho, g,r}^\ZZZ/\overline{\mfM}_{g,r}$)  $\mfO \mfp_{\mfg, \hslash, \hslash \star \rho, g,r}^\ZZZ$ is  
 generically  \'{e}tale over $\overline{\mfM}_{g,r}$,  i.e., any irreducible component that dominates  $\overline{\mfM}_{g,r}$ admits a dense open subscheme \'{e}tale over  $\overline{\mfM}_{g,r}$.
\et
\begin{proof}
The assertion for $\mfs \mfl_n$ follows from Propositions \ref{c43} and \ref{355495}.
For the remaining cases,  let us apply Proposition \ref{QW1392} to the universal family of curves.
Then,  there exists  a closed immersion $\mfO \mfp_{\mfg, \hslash, \hslash \star \rho, g,r}^\ZZZ \migiincl \mfO \mfp_{\mfs \mfl_n, \hslash, \hslash \star \rho, g,r}^\ZZZ$ over $\overline{\mfM}_{g,r}$ (for a suitable $n$ with $2 n <p$), by which we may  consider 
 $\mfO \mfp_{\mfg, \hslash, \hslash \star \rho, g,r}^\ZZZ$  as a closed substack of $\mfO \mfp_{\mfs \mfl_n, \hslash, \hslash \star \rho, g,r}^\ZZZ$.
The assertion in the case of $\mfs \mfl_n$  implies 
that $ \mfO \mfp_{\mfg, \hslash, \hslash \star \rho, g,r}^\ZZZ$ is unramified  over the points of $\overline{\mfM}_{g,r}$ classifying totally degenerate curves.
Hence, by Corollary \ref{y0170}, $ \mfO \mfp_{\mfg, \hslash, \hslash \star \rho, g,r}^\ZZZ$ is also \'{e}tale over these points.
\end{proof}

\bco \label{QG72}
Let us keep the notation 
 in Theorem \ref{3778567d}.
Then, the map of sets ${^p}N_{\mfg, 0} : \mbZ^{\mfc (\mbF_p)}_{\geq 3} \migi \mbZ$ defined in  (\ref{f186}) for this $\mfg$ forms a pseudo-fusion rule on $\mfc_\mfg (\mbF_p)$.
In particular,  
we can compute  the generic degrees $\mr{deg} (\mfO \mfp_{\mfg, \hslash, \hslash \star \rho, g,r}^\ZZZ/\overline{\mfM}_{g,r})$ ($\hslash \in k^\times$, $\rho \in \mfc^{\times r}(\mbF_p)$, $2g-2+r>0$) by using  the ring structure of 
the  resulting pseudo-fusion ring ${^p}\Fus_{\mfg}$ for dormant  $\mfg$-opers (cf. Theorem \ref{c6} and Corollary \ref{c3445}). 
\eco
\begin{proof}
The assertion follows from  Proposition \ref{c2}  and Theorems \ref{ppp012}, 
  \ref{3778567d}.
\end{proof}

\begin{rema} \label{QG91}
Let us keep the notation 
 in Theorem \ref{3778567d}.
Then, one may find
a dense open substack   
\begin{align} \label{EEE222}
\overline{\mfM}_{g, r}^\divideontimes
\index{$\overline{\mfM}_{g, r}^\divideontimes$}
\end{align}
 of $\overline{\mfM}_{g,r}$ over which the projection $\mfO \mfp^\ZZZ_{\mfg, \hslash, g, r}/\overline{\mfM}_{g,r}$ is finite, surjective,  and \'{e}tale (not only generically \'{e}tale).
Indeed, if  $\mfV$  denotes the \'{e}tale locus of $\mfO \mfp^\ZZZ_{\mfg, \hslash, g, r}$ over $\overline{\mfM}_{g,r}$,  being  an open substack of $\mfO \mfp^\ZZZ_{\mfg, \hslash, g, r}$, 
then the complement  (in $\overline{\mfM}_{g,r}$) of the image of $\mfO \mfp^\ZZZ_{\mfg, \hslash, g, r}  \setminus \mfV$ via the projection $\mfO \mfp^\ZZZ_{\mfg, \hslash, g, r} \migi \overline{\mfM}_{g,r}$  defines  a desired open substack.

If $\overline{k}$ is an algebraically closed field containing  $k$ and 
$\msX$ is 
a pointed stable curve  over $\overline{k}$ classified by $\overline{\mfM}_{g,r}^\divideontimes$,  then the number of (isomorphism classes of)  dormant $(\mfg, \hslash)$-opers of radii  $\hslash \star \rho$ (where $\rho \in \mfc (\mbF_p)$) on $\msX$ is precisely  ${^p}N_{\mfg, g} (\overline{\rho})$ (cf. (\ref{f180})).
In particular,  sufficiently general $r$-pointed genus-$g$  stable curves over an algebraically closed field  have the same number of dormant $(\mfg, \hslash)$-opers.

\end{rema}

\chapter{Comparison with Quot-schemes} \label{y0198} \vspace{3mm}

In this chapter, we give 
a proof of Joshi's conjecture by investigating the relationship between the moduli space of dormant opers and  Quot-schemes.
Our discussion in that chapter is parallel to the discussion  in ~\cite{Wak}, where the case of $\mfg = \mfs \mfl_n$ was dealt with.
Just as in {\it loc.\,cit.,} we construct 
a bijection between the set of dormant $\mfs \mfl_n$-opers and the points of a  certain Quot-scheme in positive characteristc.
The generic \'{e}taleness proved in the previous chapter enables us to  relate these sets  with the points of a Quot-scheme  over $\mbC$.
Then, by applying   a kind of the Vafa-intriligator formula computing the Gromov-Witten invariant of this Quot-scheme,
we conclude that Joshi's conjecture holds for sufficiently general curves.

\section{Quot schemes} \label{y0199}

\subsection{} \label{QG1409}
To prepare for the following  discussion, 
we introduce some notation  for Quot-schemes
in {\it arbitrary characteristic}.
Let $n$ be a positive integer and $d$ an integer.
Let $T$ be a noetherian  scheme,
$Y$  a projective smooth curve of genus $g>1$ over $T$,
  and $\mcV$ a vector bundle on $Y$. 
Denote by 
\index{$\mr{Quot}_{\mcV/Y/T}^{n, d}$, Quot-scheme}
\begin{equation} \label{QQ770}
 \mr{Quot}_{\mcV/Y/T}^{n, d} : \mfS \mfc \mfh_{/T} \longmigi \mfS \mfe \mft
 \end{equation}
the $\mfS \mfe \mft$-valued contravariant  functor on $\mfS \mfc \mfh_{/T}$ which to any $T$-scheme $t :T' \migi T$ associates the set of isomorphism classes of $\mcO_{T' \times_T Y}$-linear 
 injections 
 \begin{equation} \label{10100}
   \eta  :  \mcF \migiincl  (t \times \mr{id}_Y)^*(\mcV)
     \end{equation}
such that  the cokernel $\mr{Coker}(\eta)$ 
  is flat over $T'$ (which, since $Y/T$ is smooth of relative dimension $1$, implies that $\mcF$ is {\it locally free}) and $\mcF$ is of rank $n$ and degree $d$. 
It  is known (cf. ~\cite{FGA}, Chap.\,5, \S\,5.5, Theorem 5.14)  that $ \mr{Quot}_{\mcV/Y/T}^{n, d} $ may be  represented by a proper  scheme over $T$.

\subsection{} \label{QG1455}
Suppose  further that $T = \mr{Spec}(k)$ for a field $k$  and we are given an $\mcO_Y$-linear injection $\eta_0 : \mcF \migiincl \mcV$
classified by $\mr{Quot}_{\mcV/Y/T}^{n, d}$.
Denote by 
\begin{align} \label{QQ771}
\mcT_{[\eta_0]}  \mr{Quot}_{\mcV/Y/k}^{n, d}
\end{align}
 the tangent space  at the point $[\eta_0] \in \mr{Quot}_{\mcE/Y/k}^{n, d} (k)$  corresponding to $\eta_0$.
Recall from ~\cite{FGA}, Chap.\,6, \S\,6.4, Theorem 6.4.9,
that the deformation functor classifying deformations of $\eta_0$ has a tangent-obstruction theory by putting $\mr{Hom}_{\mcO_Y}(\mcF, \mr{Coker}(\eta_0))$ as the tangent space  and $\mr{Ext}_{\mcO_Y} (\mcF, \mr{Coker}(\eta_0))$ as the obstruction space.
In particular, 
 there exists a canonical bijection 
\begin{equation} \label{10101}
\mcT_{[\eta_0]}  \mr{Quot}_{\mcE/Y/T}^{n, d} \isom \mr{Hom}_{\mcO_Y}(\mcF, \mr{Coker}(\eta_0))
\end{equation}
of $k$-vector  spaces.

\subsection{}\label{y0200}

Let $S$ be a noetherian scheme of characteristic $p>0$ and $n$ a positive  integer with $1 <  n <\frac{p}{2}$.
Fix  an unpointed projective  smooth curve $\msX := (f: X \migi S, \emptyset)$
of genus $g>1$ over $S$.
Hence, both the   log structures of  $X^\mr{log}$ and $S^\mr{log}$ obtained in the manner of \S\,\ref{QR781} are trivial.
For convenience, we write $\Omega := \Omega_{X/S}$ and $\mcT := \mcT_{X/S}$.
Suppose that there exists a line bundle $\BB$ on $X$
together with an isomorphism $\tau_\BB : \mcO_X \isom \mcT^{\otimes \frac{n(n-1)}{2}} \otimes \BB^{\otimes n}$.
(As discussed in \S\,\ref{QQ156}, such a line bundle necessarily exists  at least \'{e}tale locally on $S$.)
In particular,  the relative degree  $\mr{det}(\BB)$ of $\BB$ is given by $(g-1)(n-1)$, and 
the pair $\bb := (\BB, \tau_{\BB*}(d))$ (cf. (\ref{QW299})) forms an $n$-theta characteristic of $X/S$.
Also, the line bundle $\BB$ associates the rank $n$ vector bundle  $\DF_\BB := \mcD^{< n}_{1, X/S} \otimes \BB$ equipped with a decreasing filtration $\{ \DF_\BB^j \}_{j=0}^n$ (cf. (\ref{QW47})).
If we set
\begin{equation} \label{104h}
\BB^\TT := \mcT^{\otimes (n-1)} \otimes \BB \left(= (\BB^\dual)^\vee \right),
\index{$\BB^\TT$}
\end{equation}
then  the quotient $\DF_\BB/\DF_\BB^1$ is naturally isomorphic to $\BB^\TT$.
Also, by using $\BB^\TT$,  we obtain  the Quot-scheme 
\index{$\mr{Quot}_{\mcV/Y/T}^{n, d}$, Quot-scheme@--- $\mr{Quot}_{\BB}^{n, 0}$, $\mr{Quot}_{\BB}^{n, \mcO}$}
\begin{equation}  \label{QG80}
\mr{Quot}_{\BB}^{n, 0} := \mr{Quot}_{F_{X/S*}(\BB^\TT)/X^{(1)}/S}^{n, 0}.
 \end{equation}
It has a subfunctor 
\begin{align} \label{QG81}
\mr{Quot}_{\BB}^{n, \mcO}
\end{align}
which, to any $T$-scheme $t: T' \migi T$, assigns the set of isomorphism classes of injections $\mcF \migiincl (t \times \mr{id}_Y)^*(\mcV)$ with $\mr{det}(\mcF) \cong \mcO_{T' \times_T Y}$.
This subfunctor may be represented by 
  a closed
subscheme of $\mr{Quot}_{\BB}^{n, 0}$.

\section{Comparison with the moduli space of dormant $(\mr{GL}_n, \bb)$-opers}\label{y0201}

Let us keep the notation in \S\,\ref{y0200}.
In this section, we shall examine   the  relationship between 
  the Quot-scheme $\mr{Quot}_{\BB}^{n, \mcO}$ and  the moduli scheme $\mfO \mfp^{\diamondsuit ^\text{Zzz...}}_{n,  \bb} := \mfO \mfp^{\diamondsuit ^\text{Zzz...}}_{n,  1,  \bb}$ classifying dormant $(\mr{GL}_n, \bb)$-opers on  $\msX$ (cf. (\ref{10034})).

\subsection{The dormant $\mr{GL}_p$-oper $\msA_\BB$} \label{QG1500}
Let us consider the $\mcO_X$-module
\begin{equation} \label{QG82}
 \mcA_{\BB} := F_{X/S}^*(F_{X/S*}(\BB^\TT)).
\index{$\mcA_{\BB}$}
 \end{equation}
 Since $F_{X/S}$ is finite and faithfully flat of degree $p$,
  $\mcA_{\BB}$ forms a rank $p$ vector bundle.
Denote by 
 \begin{align} \label{QG83}
 q : \mcA_{\BB}
  \stackrel{}{\migi} \BB^\TT
  \index{$q$}
 \end{align}
  the morphism corresponding to the identity morphism of $F_{X/S*}(\BB^\TT)$ via  
the adjunction relation ``$F_{X/S}^*(-) \dashv F_{X/S*}(-)$".
 As discussed  in \S\,\ref{QQ200}, 
 there exists a canonical $S$-connection
$\nabla_\mcA^\mr{can} := \nabla^{\mr{can}}_{F_{X/S*}(\BB^\TT)} : \mcA_{\BB}  \migi \Omega \otimes \mcA_{\BB}$
on $\mcA_\BB$ with vanishing $p$-curvature; it
 gives rise to  a {\it $p$-step decreasing filtration} 
$\{ \mcA_{\BB}^j \}_{j=0}^p$
 on
$\mcA_{\BB}$
 as follows.
\begin{align} \label{QG125}
\mcA_{\BB}^0  &:= \mcA_{\BB};\\
\mcA_{\BB}^1 &:= \mr{Ker}(\mcA_{\BB} \stackrel{q}{\migisurj} \BB^\TT); \notag \\
\mcA_{\BB}^j &:= \mr{Ker}(\mcA_{\BB}^{j-1} \xrightarrow{\nabla^{\mr{can}}_\mcA |_{\mcA_{\BB}^{j-1}}} \Omega\otimes \mcA_{\BB} \xrightarrow{\mr{quotient}} \Omega\otimes (\mcA_{\BB}/\mcA_{\BB}^{j-1})) \notag
\end{align}
($j= 2, \cdots, p$).
In particular,  $q$ induces an isomorphism $\overline{q} : \mcA_\BB /\mcA^1_\BB \isom \BB^\TT$.

\bpr \label{y0202}
 The collection of data 
\begin{equation}
 \msA_{\BB}^\heartsuit : = ( \mcA_{\BB}, \nabla^{\mr{can}}_\mcA, \{  \mcA_{\BB}^j \}_{j=0}^p)
 \index{$\msA_{\BB}^\heartsuit$}
\end{equation}
 forms a dormant $\mr{GL}_p$-oper on $\msX$.
 In particular, we have  $\mcA_{\BB}^j/\mcA_{\BB}^{j+1} \cong \Omega^{\otimes (-n+1+j)} \otimes \BB$ for $j=0, \cdots, n-1$ and $\mr{det}( \mcA_{\BB}/\mcA_\BB^n) \cong \mcO_X$.
 \epr
\begin{proof}
The former assertion  can be proved, even when $S$ is a general scheme, by  
 an argument 
  similar to
the argument 
(in the case where $S= \mr{Spec}(k)$ for an algebraically closed field $k$)
 given in the proofs of ~\cite{JRXY}, \S\,5.3, Theorem, and ~\cite{SUN}, \S\,2, Lemma 2.1.
 The latter assertion follows from (\ref{QW1909}) and (\ref{QW1906}).
\end{proof}

\subsection{Constructing dormant $\mr{GL}_n$-opers} \label{QG204}

Let $\eta : \mcV \migiincl F_{X/S*}(\BB^\TT)$ be the  injection classified by an $S$-rational point of $\mr{Quot}_{\BB}^{n, 0}$ and write $\mcF_\eta := F_{X/S}^*(\mcV)$.
By using the pull-back $F_{X/S}^*(\eta) : \mcF_\eta \migiincl \mcA_\BB$ of $\eta$, we define an $n$-step  filtration   $\{ \mcF^j_\eta \}_{j=0}^n$ on $\mcF_\eta$ by putting  $\mcF^j_\eta := F_{X/S}^*(\eta)^{-1} (\mcA^j_\BB)$.

\bpr \label{QG180}
\begin{itemize}
\item[(i)]
 The composite
 \begin{equation}\label{QG100}
\breve{\eta} :   \mcF_\eta \xrightarrow{F^*_{X/S}(\eta)} \mcA_{\BB} \xrightarrow{\mr{quotient}} \mcA_{\BB} /\mcA_{\BB}^n \end{equation}
   is an isomorphism of $\mcO_X$-modules.
\item[(ii)]
The triple 
\begin{align} \label{QG190}
\msF_\eta^\heartsuit := (\mcF_\eta, \nabla^\mr{can}_\mcV, \{ \mcF_\eta^j \}_{j=0}^n)
\end{align}
 forms a dormant $\mr{GL}_n$-oper on $\msX$.
\end{itemize}
\epr
\begin{proof}
First, let us  consider  assertion (i).
Since both $\mcF_\eta$ and $\mcA_{\BB} /\mcA_{\BB}^n$ are flat over $S$,
we may assume, by  
considering the various fibers
over $S$, that $S = \mr{Spec}(k)$ for  a field $k$ containing $\mbF_p$.  
If we write $\mcN^j := \mcF_\eta^j/\mcF_\eta^{j+1}$
 ($j= 0, \cdots, p-1$),
 then  $F_{X/S}^*(\eta):\mcF_\eta \migiincl \mcA_\BB$ induces 
    an injection $\mr{gr}^j (F_{X/S}^*(\eta)) : \mcN^j \migiincl \mcA^j_{\BB} /\mcA^{j+1}_{\BB}$.
 Note that 
 $\mcN^j$ is either trivial or a line bundle because  $\mcA^j_{\BB} /\mcA^{j+1}_{\BB}$ is a line bundle on the smooth curve $X$ over $k$.
The rank of  $\mcF_\eta$ equals  $n$, so
the cardinality of the set $I := \big\{ j  \, \big|  \, \mcN^j \neq 0\big\}$ is exactly $n$.
Since  $\msA^\heartsuit_\BB$ defines  a $\mr{GL}_p$-oper and $F_{X/S}^*(\eta)$   is compatible with the $S$-connections $\nabla^{\mr{can}}_\mcV$ and  $\nabla^{\mr{can}}_\mcA$,
  we see that  $\mcN^{j+1} \neq 0$ implies $\mcN^j \neq 0$.
  By this fact,     $I$ coincides with $\{ 0,1, \cdots, n-1\}$  and 
 the composite $\breve{\eta}$ 
is  an isomorphism at the generic point of $X$.
On the other hand, the following equalities hold:
\begin{align} \label{QG104}
\mr{deg}(\mcF_\eta) &=  p \cdot \mr{deg}(\mcV) = p \cdot 0 = 0,  \\
\label{QH888}\mr{deg}(\mcA_{\BB}/\mcA_{\BB}^n)  
  &= \sum_{j=0}^{n-1}  \mr{deg}(\mcA_{\BB}^j/\mcA_{\BB}^{j+1}) \\
  &= \sum_{j=0}^{n-1}  \mr{deg}(\Omega^{\otimes (-n+1+j)}\otimes\BB) \notag \\
  &  = \sum_{j=0}^{n-1}  (g-1)(-n+1 +2j) \notag  \\
  & = 0, \notag
 \end{align}
where the second equality of (\ref{QH888}) follows from Proposition \ref{y0202}.
Thus,  by comparing the  respective degrees of $\mcF_\eta$ and $\mcA_{\BB}/\mcA_{\BB}^n$,
we conclude that  $\breve{\eta}$ is  an isomorphism.  (In particular,  the morphism $\mr{gr}^j (F^*_{X/S}(\eta))$ for each $j =0, \cdots, n-1$  is an isomorphism.) This completes the proof of  assertion (i).

Next, we shall consider assertion (ii).
Since 
the  $S$-connection $\nabla^\mr{can}_\mcV$ is compatible with  $\nabla^\mr{can}_\mcA$ via 
$F_{X/S}^*(\eta)$,  
the following diagram defined  for each $j=1, \cdots, n-1$ is  commutative:
\begin{align} \label{QG106}
\vcenter{\xymatrix@C=66pt@R=36pt{
\mcN^j \ar[r]^-{\mr{gr}^j (F_{X/S}^*(\eta))}_-{\sim} \ar[d]_-{\mr{KS}_{\msF^\heartsuit_\eta}^j}  & \mcA^j_{\BB} /\mcA^{j+1}_{\BB} \ar[d]^-{\mr{KS}_{\msA_{\BB}^\heartsuit}^j}\\
\Omega\otimes \mcN^{j-1} \ar[r]^-{\sim}_-{\mr{id} \otimes \mr{gr}^{j-1} (F_{X/S}^*(\eta))} & \Omega \otimes (\mcA^{j-1}_{\BB} /\mcA^{j}_{\BB}), 
}}
\end{align}
where the right-hand vertical arrow $\mr{KS}_{\msA_{\BB}^\heartsuit}^j$ is an isomorphism because $\msA_\BB^\heartsuit$ forms a $\mr{GL}_p$-oper.
Hence,   $\mr{KS}_{\msF^\heartsuit_\eta}^j$'s   are  isomorphisms, i.e.,  $\msF^\heartsuit_\eta$ specifies a  $\mr{GL}_n$-oper.
The dormancy condition on $\msF^\heartsuit_\eta$ is clear  from  the definition of $\nabla_\mcV^\mr{can}$, so
  we finish the proof of the assertion.
\end{proof}

\subsection{Constructing dormant $(\mr{GL}_n, \bb)$-opers} \label{QG1504}
By Proposition \ref{QG180}, (i), the isomorphism $\breve{\eta}$ preserves  the filtrations $\{ \mcF_\eta^j \}_j$,  $\{ \mcA_\BB^j/\mcA_\BB^n \}_j$.
In particular, $\breve{\eta}$ restricts to an isomorphism $\breve{\eta}|_{\mcF^{n-1}_\eta} : \mcF^{n-1}_\eta \isom \mcA^{n-1}_\BB/\mcA^{n}_\BB$.
We shall identify $\mcF^{n-1}_\eta$ with $\BB$ via the composite isomorphism
\begin{align} \label{QG199}
 \mcF^{n-1}_\eta \xrightarrow{\breve{\eta}|_{\mcF^{n-1}_\eta}} \mcA^{n-1}_\BB/\mcA^{n}_\BB \isom \Omega^{\otimes (n-1)} \otimes (\mcA_\BB/\mcA_\BB^1) \xrightarrow{\mr{id} \otimes \overline{q}} \Omega^{\otimes (n-1)} \otimes \BB^\TT \isom \BB,
\end{align}
where the second arrow is  induced from (\ref{QW1909}) in the case where ``$\msF^\heartsuit$'' is taken to be $\msA_\BB^\heartsuit$ and the last arrow arises from  the natural pairing $\Omega \times \mcT \migi \mcO_X$.
Since $\msF_\eta^\heartsuit$ forms a $\mr{GL}_n$-oper,
the composite
\begin{align} \label{QG198}
\widehat{\eta} : \DF_\BB \left(:= \mcD^{<n}_{1, X/S} \otimes \BB\right) \xrightarrow{\mr{inclusion}} \mcD^{<\infty}_{1, X/S} \otimes \mcF_\eta \xrightarrow{{^\mcD}\hspace{-1.5mm}\nabla^\mr{can}_\mcV} \mcF_\eta
\end{align}
is verified to be an isomorphism.
If  
\begin{align} \label{QG201}
\nabla^\diamondsuit_\eta : \DF_\BB \migi \Omega \otimes \DF_\BB
\end{align}
 denotes  the $S$-connection on $\DF_\BB$ corresponding to $\nabla_\mcV^\mr{can}$ via this composite,
then Proposition \ref{QG180}, (ii), shows that $\nabla^\diamondsuit_\eta$ forms a dormant $(\mr{GL}_n, \BB)$-oper.

\bpr \label{QG187}
Suppose that $\eta$ lies in the subset $\mr{Quot}_\BB^{n, \triv}(S)$, i.e., $\mr{det}(\mcV) \cong \mcO_{X^{(1)}}$.
Then,  $\nabla_\eta^\diamondsuit$ specifies a dormant $(\mr{GL}_n, \bb)$-oper on $\msX$.
\epr
\begin{proof}
By assumption, there exists 
 an isomorphism $\nu : \mr{det}(\mcV) \isom \mcO_{X^{(1)}}$.
Since $\Gamma (X, \mcO^\times_X) = \Gamma (S, \mcO_S^\times)$,  the composite 
\begin{align}
\mcO_X &\xrightarrow{\tau_\BB} \mcT^{\otimes \frac{n(n-1)}{2}} \otimes \BB^{\otimes n} \\
& \xrightarrow{(\ref{QW1906})^{-1}} \mr{det}(\DF_\BB) \notag \\
&  \xrightarrow{\mr{det}(\widehat{\eta})} \mr{det} (\mcF_\eta) 
\left(= F^*_{X/S}(\mr{det}(\mcV))\right) \notag \\
&  \xrightarrow{F^*(\nu)} \mcO_X \notag 
\end{align}
coincides with  the automorphism $\mr{mult}_u$ of $\mcO_X$ given by multiplication by some element $u \in \Gamma (S, \mcO_S^\times)$.
In particular,   this automorphism commutes with  the universal derivation $d$.
Hence, the isomorphism $\mr{det}(\DF_\BB) \isom \mcT^{\otimes \frac{n(n-1)}{2}} \otimes \BB^{\otimes n}$ defined in (\ref{QW1906}) can be described as the composite of isomorphisms between  flat line bundles
\begin{align}
(\mr{det}(\DF_\BB), \mr{det}(\nabla_\eta^\diamondsuit)) & \xrightarrow{\mr{det}(\widehat{\eta})}
(\mr{det}(\mcF_\eta), \mr{det}(\nabla_\mcV^\mr{can})) \\
& \xrightarrow{F^*(\nu)} 
(\mcO_X, d) \left(= (\mcO_X, \nabla_{\mcO_{X^{(1)}}}^\mr{can})\right) \notag \\
& \xrightarrow{\mr{mult}_u^{-1}}
(\mcO_X, d) \notag \\
& \xrightarrow{\tau_\BB}
(\mcT^{\otimes \frac{n(n-1)}{2}} \otimes \BB^{\otimes n}, \tau_{\BB*}(d)). \notag
\end{align}
That is to say, the $S$-connection $\mr{det} (\nabla_\eta^\diamondsuit)$ corresponds to $\tau_{\BB*}(d)$ via this isomorphism.
This means that $\nabla_\eta^\diamondsuit$ turns out to be  a dormant $(\mr{GL}_n, \bb)$-oper on $\msX$, as desired.
\end{proof}


\subsection{A bijective correspondence} \label{QG1560}
By applying the above proposition,   we can construct an isomorphism between the moduli space $\mfO \mfp_{n, \bb}^{\diamondsuit ^\text{Zzz...}}\hspace{0mm}$
  and the Quot-scheme $\mr{Quot}_\BB^{n, \triv}$, as follows.

\bpr  \label{y0204}
Let $\msX$ and $\bb$ be as above.
Then,
there exists an isomorphism  of $S$-schemes
\begin{equation} \label{QQ790}
\mr{Quot}_\BB^{n, \triv}  \isom \mfO \mfp_{n,  \bb}^{\diamondsuit ^\text{Zzz...}}.  
 \end{equation}
 In particular, $\mr{Quot}_\BB^{n, \triv}$ is isomorphic to 
 $\mfO \mfp^\ZZZ_{\mfs \mfl_n, \msX} \hspace{-1mm} := \mfO \mfp^\ZZZ_{\mfs \mfl_n, 1, \msX}$, which implies that  $\mr{Quot}_\BB^{n, \triv}$  is  
 nonempty and finite over $S$, and moreover,  \'{e}tale over $S$ if $\msX$ is classified by $\overline{\mfM}_{g, 0}^\divideontimes$ (cf. (\ref{EEE222})).
 \epr
\begin{proof}
We shall consider the former assertion.
The assignment  $\eta \mapsto \nabla^\diamondsuit_\eta$
resulting from  Proposition \ref{QG187}
determines   a map of sets
\begin{equation} \label{QHH1212}
\mr{Quot}_\BB^{n, \triv} (S) \migi  \mfO \mfp_{n,  \bb}^{\diamondsuit ^\text{Zzz...}} (S). 
\end{equation}
Since the formation of this map is functorial 
 with respect to $S$, the problem is reduced to 
 proving the bijectivity of this map.

Let   us take 
  a dormant $(\mr{GL}_n, \bb)$-oper $\nabla^\diamondsuit$ on $\msX$.
We shall  identify $\DF_\BB$ with $F^*_{X/S}(\mr{Ker}(\nabla^\diamondsuit))$ via the isomorphism 
$\nu^\sharp_{(\DF_\BB, \nabla^\diamondsuit)} : F^*_{X/S}(\mr{Ker}(\nabla^\diamondsuit)) \isom \DF_\BB$ (cf. (\ref{QW120})).
If  $q_\mcF$  denotes the natural quotient $\left(F^*_{X/S}(\mr{Ker}(\nabla^\diamondsuit)) = \right) \DF_\BB \migisurj \BB^\TT$, then
it corresponds to a morphism 
\begin{align}
\eta_{\nabla^\diamondsuit}:\mr{Ker}(\nabla^\diamondsuit) \migi F_{X/S*}(\BB^\TT)
\end{align}
 via  the adjunction relation ``$F^*_{X/S}(-) \dashv F_{X/S*}(-)$".
In what follows, we shall prove the claim that $\eta_{\nabla^\diamondsuit}$ is injective.
Write 
 \begin{equation}
  \eta_{\nabla^\diamondsuit}^F : 
  \DF_\BB \left(=  F^*_{X/S}(\mr{Ker}(\nabla^\diamondsuit)) \right)
   \migi F_{X/S}^*(F_{X/S*}(\BB^\TT))   \left(=\mcA_{\BB}\right)
  \end{equation} 
for the pull-back of this morphism 
 by $F_{X/S}$.
By construction, $\eta_{\nabla^\diamondsuit}^F$ commutes  with the respective  $S$-connections (i.e., $\nabla^\diamondsuit$ and $\nabla_\mcA^\mr{can}$), as well as  with   the respective surjections onto $\BB^\TT$ (i.e., $q_\mcF$ and  $q$).
Here, notice that, since $(\DF_\BB, \nabla^\diamondsuit, \{ \DF_\BB^j \}_j)$ forms a $\mr{GL}_n$-oper, the filtration $\{ \DF_\BB^j \}_j$ can be characterized inductively  by means of $\nabla^\diamondsuit$ and $q_\mcF$ as follows:
\begin{align} \label{QG122}
\DF_{\BB}^0  &:= \DF_{\BB};\\
\DF_{\BB}^1 &:= \mr{Ker}(\DF_{\BB} \stackrel{q_\mcF}{\migisurj} \BB^\TT); \notag \\
\DF_{\BB}^j &:= \mr{Ker}(\DF_{\BB}^{j-1} \xrightarrow{\nabla^\diamondsuit |_{\DF_{\BB}^{j-1}}} \Omega \otimes \DF_{\BB} \xrightarrow{\mr{quotient}} \Omega\otimes (\DF_{\BB}/\DF_{\BB}^{j-1})) \notag
\end{align}
($j= 2, \cdots, n$).
By taking account of this characterization of $\{\DF_\BB^j \}_j$ and  the definition  of $\{ \mcA_\BB^j \}_j$ described in  (\ref{QG125}),
we see that $\eta^F_{\nabla^\diamondsuit}$ preserves the filtrations.
Hence, $\eta_{\nabla^\diamondsuit}^F$ induces, via taking subquotients,  morphisms  $\mr{gr}^j (\eta_{\nabla^\diamondsuit}^F) : \DF_\BB^j/\DF_\BB^{j+1} \migi \mcA^j_\BB/\mcA^{j+1}_\BB$ ($j=0, \cdots, n-1$) which make 
the following square diagrams  for various $j$'s commute:
\begin{align} \label{QG130}
\vcenter{\xymatrix@C=66pt@R=36pt{
\DF_\BB^{j}/\DF_\BB^{j+1} \ar[r]^-{\mr{gr}^j (\eta_{\nabla^\diamondsuit}^F)} \ar[d]_-{\mr{KS}_{\nabla^{\diamondsuit \Rightarrow \heartsuit}}^j}^-{\wr}& \mcA^j_{\BB}/\mcA^{j+1}_{\BB} \ar[d]_-{\wr}^-{\mr{KS}^j_{\msA^\heartsuit_{\BB}}}\\
\Omega \otimes (\DF_\BB^{j-1}/\DF_\BB^{j}) \ar[r]_-{\mr{id} \otimes \mr{gr}^{j-1} (\eta_{\nabla^\diamondsuit}^F)}&\Omega\otimes (\mcA^{j-1}_{\BB}/\mcA^{j}_{\BB}).
}}
\end{align}
Since $\mr{gr}^{0} (\eta_{\nabla^\diamondsuit}^F) = \mr{id}_{\BB^\TT}$, 
the commutativity of these diagrams shows, by induction on $j$, that
 the morphisms $\mr{gr}^{j} (\eta_{\nabla^\diamondsuit}^F)$   are  isomorphisms.
In particular, $\eta_{\nabla^\diamondsuit}^F \left(= F_{X/S}^*(\eta_{\nabla^\diamondsuit})\right)$ turns out to be injective, and we obtain the injectivity of $\eta_{\nabla^\diamondsuit}$ because of 
the faithful flatness of $F_{X/S}$.
This completes the proof of the claim.
Moreover, by applying the same argument  to the base-changes  of $\eta_{\nabla^\diamondsuit}$ to   $S$-schemes, one concludes that $\eta_{\nabla^\diamondsuit}$ is  universally injective with respect to base-change over $S$.
This implies  that $\mr{Coker}(\eta_{\nabla^\diamondsuit})$ is flat over $S$ (cf.~\cite{MAT}, p.\,17, Theorem 1).
Finally, the composite of (\ref{QW1906}) in our case  and  $\tau_\BB^{-1}$   determines an isomorphism of flat line bundles $(\mr{det}(\DF_\BB), \mr{det}(\nabla^\diamondsuit)) \isom (\mcO_X, d)$, which induces an isomorphism $\mr{det}(\mr{Ker}(\nabla^\diamondsuit)) \isom \mcO_{X^{(1)}}$.
As a consequence,  $\eta_{\nabla^\diamondsuit}$ specifies an  $S$-rational point of $\mr{Quot}_\BB^{n, \mcO}$.
One verifies that the resulting assignment $\nabla^\diamondsuit \mapsto \eta_{\nabla^\diamondsuit}$ gives the inverse to (\ref{QHH1212}), so this completes the proof of the former  assertion.

The latter assertion follows from the former assertion, Theorem \ref{QS1040},   Corollary \ref{1090}, and the discussion in Remark \ref{QG91}.
\end{proof}

\section{The relationship between the Quot-schemes} \label{y0206}

Next, we shall consider the relationship between  $\mr{Quot}_\BB^{n, \triv}$  and  
$\mr{Quot}_\BB^{n, 0}$.
By pulling back line bundles on $X^{(1)}$
via  $F_{X/S} : X \migi X^{(1)}$, we obtain   an $S$-morphism 
\begin{align}
 \mr{Pic}^0_{X^{(1)}/S}  & \migi \mr{Pic}_{X/S}^0;  [\mcN] \mapsto [F_{X/S}^*(\mcN)] 
   \end{align}
between relative Picard schemes classifying degree $0$ line bundles.
We shall denote by
\begin{equation} 
\mr{Ver}_{X/S}
\index{$\mr{Ver}_{X/S}$}
 \end{equation}
the $S$-scheme defined to be the scheme-theoretic inverse image, via this morphism,  of  the $S$-rational point $[\mcO_X] \in \mr{Pic}_{X/S}^0 (S)$.
Let 
\begin{align} \label{QHHH22245}
\mfM_{g, 0}^{\mr{ord}}
\index{$\mfM_{g, 0}^{\mr{ord}}$}
\end{align}
 denote   the locus of $\mfM_{g, 0}$ classifying ordinary smooth curves, which  is open and dense.
It is well-known  (cf. ~\cite{DG}, Expose VII, \S\,4.3; ~\cite{MIL}, \S\,8, Proposition 8.1 and Theorem 8.2; ~\cite{MES}, Appendix, Lemma (1.0))
 that $\mr{Ver}_{X/S}$ is finite and faithfully flat over $S$ of degree $p^g$ and, moreover, \'{e}tale 
 over  the points whose fibers are classified by $\mfM_{g, 0}^{\mr{ord}}$.

\bpr \label{y0207}
There exists an isomorphism of $S$-schemes
\begin{equation} \label{00w}
 \mr{Quot}_\BB^{n, \triv} \times_{S}  \mr{Ver}_{X/S} \isom \mr{Quot}_\BB^{n, 0}. 
  \end{equation}
 \epr
\begin{proof}
It  suffices to prove  that 
there is a bijection between 
the respective sets of $S$-rational points that is functorial with respect to $S$.

Let $(\eta :\mcV \migiincl F_{X/S*}(\BB^\TT), \mcN)$ be 
the  pair of data corresponding to an $S$-rational point  of
$\mr{Quot}_\BB^{n, \triv} \times_{S}  \mr{Ver}_{X/S}$. 
 Then, the  composite 
\begin{align} \label{91134}
   \eta_\mcN : \mcV \otimes \mcN 
   & \xrightarrow{\eta \otimes \mr{id}_\mcN} F_{X/S*}(\BB^\TT) \otimes \mcN \\
   &  \isom F_{X/S*}(\BB^\TT \otimes F^*_{X/S}(\mcN)) \notag \\
   & \isom F_{X/S*}(\BB^\TT \otimes \mcO_X) \notag \\
    & \isom F_{X/S*}(\BB^\TT)   \notag
   \end{align}
 determines an element of $ \mr{Quot}^{n, 0}_\BB (S)$, where the second arrow  arises from the projection formula.
Thus, 
the assignment $(\eta, \mcN) \mapsto \eta_\mcN$ determines 
 a map  of sets
\begin{equation} \label{QH3221}
(\mr{Quot}_\BB^{n, \triv} \times_{S}  \mr{Ver}_{X/S})(S) \isom \mr{Quot}_\BB^{n, 0}(S),
\end{equation}
which is immediately verified to be   functorial with respect to $S$.

Conversely, let $\eta :\mcV \migi F_{X/S*}(\BB^\TT)$ be an injection  classified by the set $\mr{Quot}_\BB^{n, 0} (S)$.
 Since $(p, n) =1$, we can find   a pair of positive integers $(a, b)$ with  $1 + n \cdot a = p \cdot b$. 
Let us consider the injection 
\begin{equation}
\eta_{\mr{det}(\mcV)^{\otimes a}} : \mcV \otimes \mr{det}(\mcV)^{\otimes a} \migi F_{X/S*}(\BB^\TT), 
\end{equation}
 i.e., the morphism defined in the same  fashion as the above construction of $\eta_\mcN$,  where  ``$\mcN$'' is taken to be  $\mr{det}(\mcV)^{\otimes a}$.
It follows from Propositions \ref{y0202} and \ref{QG180} 
 that
\begin{align}
F^*_{X/S}(\mr{det}(\mcV)) \cong  \mr{det}(F_{X/S}^*(\mcV)) \cong
 \mr{det}(\mcA_\BB/\mcA_\BB^n) \cong \mcO_X.
\end{align}
This implies $F^*_{X/S}(\mr{det}(\mcV)^{\otimes (-a)}) \left(\cong F_{X/S}^*(\mr{det}(\mcV))^{\otimes (-a)} \right)\cong \mcO_X$
and
\begin{align}
  \mr{det}(\mcV \otimes \mr{det}(\mcV)^{\otimes a})  
  & \cong  \mr{det}(\mcV) \otimes  \mr{det}(\mcV)^{\otimes n \cdot a}  \\
  &  \cong \mr{det}(\mcV)^{\otimes p \cdot b}  \notag \\
  & \cong  (\mr{id}_X \times F_S)^*(F^*_{X/S}(\mr{det}(\mcV))^{\otimes b}) \notag \\
  & \cong   (\mr{id}_X \times F_S)^* (\mcO_X^{\otimes b}) \notag \\
  & \cong \mcO_X. \notag
  \end{align}
Consequently, the pair of data 
\begin{equation}
(\eta_{\mr{det}(\mcV)^{\otimes a}}, \mr{det}(\mcV)^{\otimes (-a)})
\end{equation}
 defines an 
$S$-rational point of $\mr{Quot}_\BB^{n, \triv} \times_{S}  \mr{Ver}_{X/S}$.
One verifies that the resulting  assignment $\eta \mapsto (\eta_{\mr{det}(\mcV)^{\otimes a}}, \mr{det}(\mcV)^{\otimes (-a)})$  defines   the  inverse to (\ref{QH3221}), so 
this completes the proof of Proposition \ref{y0207}. 
\end{proof}

\bco \label{QG93}
 $\mr{Quot}_\BB^{n, 0}$ are nonempty and finite over $S$.
If, moreover, $\msX$ is classified by $\mfM_{g, 0}^{\mr{ord}} \cap \overline{\mfM}_{g, 0}^\divideontimes$ (cf. (\ref{QHHH22245}), (\ref{EEE222})),
then the structure morphism $\mr{Quot}_\BB^{n, 0} \migi S$  is surjective and \'{e}tale,  and 
its  degree  
 satisfies the following equalities:
\begin{align}
\frac{1}{p^g} \cdot \mr{deg}(\mr{Quot}_\BB^{n, 0}/S) =  \mr{deg}(\mr{Quot}_\BB^{n, \mcO}/S) = \mr{deg}(\mfO \mfp^\ZZZ_{\mfs \mfl_n, g, 0}/\overline{\mfM}_{g,0}).
\end{align}
\eco
\begin{proof}
The assertion follows from Propositions \ref{y0204} and \ref{y0207} together with  the fact mentioned in  Remark \ref{QG91}.
\end{proof}

\section{The Vafa-Intriligator formula} \label{y0210}

In the following, we review a numerical formula concerning  the degree of a certain Quot-scheme over  $\mbC$ and relate it to the degree of the Quot-scheme $\mr{Quot}_\BB^{n, 0}$ discussed above.
Let $C$ be a   smooth  curve over $\mbC$ of genus $g>1$.
If $n$  is a positive  integer, and $\mcV$ is
 a vector bundle  on $C$ of rank $m$ ($\geq n$) and degree $d$,
 then we define invariants
 \begin{align}
  e_{\mr{max}}(\mcV,n) & := \mr{max} \big\{  \mr{deg}(\mcW) \in \mbZ\, \big|\,  \text{$\mcW$ is a subbundle of $\mcV$ of rank $n$}   \big\},   \\
  s_n(\mcV) & := d \cdot n - m \cdot e_{\mr{max}}(\mcV,n).\notag  
  \index{$e_{\mr{max}}(\mcV,n)$}
  \index{$s_n(\mcV)$}
  \end{align}
(One verifies immediately, for instance, by considering an embedding of $\mcV$ into a direct sum of $n$ line bundles, that $e_{\mr{max}}(\mcV,n) $ is well-defined.)

We review some facts concerning these invariants (cf. ~\cite{HIR};  ~\cite{LN};  ~\cite{H}).
Denote by $\mfU_C^{m,d}$ the moduli space of stable bundles on $C$ of rank $m$ and degree $d$ (cf. ~\cite{LN}, \S\,1).
It is known that  $\mfU_C^{m,d}$ is irreducible  (cf.  the discussion at the beginning of ~\cite{LN}, \S\,2).
Thus, it makes sense to speak of a ``sufficiently general" stable bundle in $\mfU_C^{m,d}$, i.e., a stable bundle that corresponds to a point of the scheme $\mfU_C^{m,d}$ that lies outside some fixed  closed subscheme.
If $\mcV$ is  a sufficiently general stable bundle in $\mfU_C^{m,d}$,
then we have 
$ s_n (\mcV) = n (m-n)(g-1) + \epsilon$ (cf. ~\cite{LN}, \S\,1), where
$\epsilon$ is the unique integer such that $0 \leq \epsilon < m$ and $s_n(\mcV) \equiv n \cdot d$ mod $m$.
Also, 
the number $\epsilon$
coincides  with the dimension of every irreducible component of the Quot-scheme 
$\mr{Quot}^{n, e_\mr{max}(\mcV, n)}_{\mcV/C/\mbC}$ (cf. ~\cite{H}, \S\,1).
If, moreover, the equality $s_n(\mcV) =n(m-n)(g-1)$ holds (i.e., $\mr{dim}(\mr{Quot}^{n, e_\mr{max}(\mcV,n)}_{\mcV/C/\mbC})=0$),
then $\mr{Quot}^{n, e_\mr{max}(\mcV,n)}_{\mcV/C/\mbC}$ is \'{e}tale over $\mr{Spec}(\mbC)$ (cf. ~\cite{H}, \S\,1).
Finally,  under this particular assumption, a formula for the degree of this Quot-scheme was given by Y. Holla
as follows. 
\bt  \label{y0230} 
Let $C$ be a  smooth curve over $\mbC$ of genus $g>1$, $\mcV$  a sufficiently general stable bundle in $\mfU_C^{m,d}$.
Write $(a,b)$ for the unique pair of integers such that $d =a \cdot m-b$ with $0 \leq b < m$.
Also, we assume  that the equality 
\begin{equation}
s_n(\mcV) = n (m-n)(g-1),
\end{equation}
holds, or equivalently, the equality
 \begin{equation} e_\mr{max}(\mcE,n) =  \frac{1}{m} \cdot (d \cdot n - n(m-n)(g-1))
\end{equation}
 holds (cf. the above discussion).
 Then, the degree  $  \mr{deg}(\mr{Quot}^{n,e_\mr{max}(\mcV,n)}_{\mcV/C/\mbC}/\mbC)$ of  the finite $\mbC$-scheme  $\mr{Quot}^{n,e_\mr{max}(\mcV,n)}_{\mcV/C/\mbC}$  is given by the following formula:
 \begin{align}
 & \ \ \ \   \mr{deg}(\mr{Quot}^{n,e_\mr{max}(\mcV,n)}_{\mcV/C/\mbC}/\mbC) \\
  &  = \frac{(-1)^{(n-1)(b \cdot n-(g-1)n^2)/m}\cdot m^{n(g-1)}}{n!}\cdot \sum_{\zeta_1 , \cdots , \zeta_n}\frac{(\prod_{i=1}^n \zeta_i)^{b-g+1}}{\prod_{i \neq j} (\zeta_i - \zeta_j)^{g-1}},\notag \end{align}
 where the sum on 
 the right-hand side is taken over the set of $n$-tuples 
  $(\zeta_1 , \cdots , \zeta_n) \in \mbC^{\times n}$ of mutually distinct $m$-th roots of unity in $\mbC$.
 \et
\begin{proof}
The assertion follows from  ~\cite{H}, \S\,4, Theorem 4.2, where ``$k$'' (respectively, ``$r$'') corresponds to our $n$ (respectively, $m$).
\end{proof}


By applying this formula,
we  conclude the same kind of formula 
for certain vector bundles 
in
positive characteristic,
as follows.
\bt \label{y0231} 
 Let  $p$ be a prime, 
 $n$  a positive  integer with $2 n <p$, 
 $k$ an algebraically closed field of characteristic $p$, $X$ a   smooth curve over $k$ of genus $g>1$, and $\BB$ a line bundle on $X$ with  $\BB^{\otimes n} \cong \Omega_{X/k}^{\otimes \frac{n (n-1)}{2}}$.
Assume  that $X/k$ is sufficiently general in $\mfM_{g,0}$.
(Here, we recall  from ~\cite{DM}, \S\,5, that $\overline{\mfM}_{g,0}$ is irreducible; thus, it makes sense to speak of a ``sufficiently general'' $X/k$, i.e., an $X/k$ specifying  a point of $\overline{\mfM}_{g,0}$ that lies outside some fixed closed substack.)
Then,  the Quot-scheme 
\begin{align} \label{QG88}
\mr{Quot}_\BB^{n, 0}:=\mr{Quot}_{F_{X/k*}(\BB^\TT)/X^{(1)}/k}^{n, 0}
\end{align}
where $\BB^\TT := \mcT_{X/k}^{\otimes (n-1)}\otimes \BB$, 
  is (nonempty and)  finite and \'{e}tale  over $\mr{Spec}(k)$.
If, moreover, the inequality $p > n \cdot  (g-1)$ holds, then
 the degree $\mr{deg} (\mr{Quot}_\BB^{n, 0}/k)$
  of $\mr{Quot}_\BB^{n, 0}$ over $\mr{Spec}(k)$ is given by the following formula:
\begin{align}  \label{QG90}
\mr{deg} (\mr{Quot}_\BB^{n, 0}/k)
  =    \frac{ p^{n(g-1)}}{n!} \cdot  
 \sum_{\genfrac{.}{.}{0pt}{}{(\zeta_1, \cdots, \zeta_n) \in \mbC^{\times n} }{ \zeta_i^p=1, \ \zeta_i \neq \zeta_j (i\neq j)}}
\frac{(\prod_{i=1}^n\zeta_i)^{(n-1)(g-1)}}{\prod_{i\neq j}(\zeta_i -\zeta_j)^{g-1}}.  
   \end{align}
 \et
\begin{proof}
One may assume that 
 the curve $X$ is  classified by 
 $\mfM_{g,0}^{\mr{ord}} \cap \overline{\mfM}_{g, 0}^\divideontimes$.
Hence, by Corollary \ref{QG93}, 
$\mr{Quot}_\BB^{n, 0}$ is (nonempty and) finite and \'{e}tale  over $\mr{Spec}(k)$.
In the following, we shall determine the value  $\mr{deg} (\mr{Quot}_\BB^{n, 0}/k)$ under the assumption $p > n \cdot (g-1)$.
Denote by $W$ the ring of Witt vectors with coefficients in $k$ and
$K$ the fraction field of $W$.
Since  
$\mr{dim}(X^{(1)}) =1$ (which implies that $H^2(X^{(1)}, \mcT_{X^{(1)}/k})=0$),
 it follows from  well-known generalities concerning  deformation theory  that 
  $X^{(1)}$ may be lifted to a smooth curve $X^{(1)}_W$ over $W$ of genus $g$.
In a similar vein, the equality  $H^2(X^{(1)}, \mcE nd_{\mcO_{X^{(1)}}} (F_{X/k*}(\BB^\TT)))=0$ implies that  $F_{X/k*}(\BB^\TT)$ may be lifted to a vector bundle $\mcG_W$ on $X^{(1)}_W$.

Now let  $[\eta]$ be the $k$-rational point of  $\mr{Quot}_\BB^{n, 0}$ classifying
an injection  $\eta : \mcV \migiincl F_{X/k*}(\BB^\TT)$.
The tangent space $\mcT_{[\eta]}\mr{Quot}_\BB^{n, 0}$ at $[\eta]$ may be naturally identified  with
the $k$-vector space $\mr{Hom}_{\mcO_{X^{(1)}}}(\mcV, \mr{Coker}(\eta))$
and the obstruction to lifting $\eta$ to any first-order thickening of $\mr{Spec}(k)$ is given by an element of $\mr{Ext}^1_{\mcO_{X^{(1)}}}(\mcV, \mr{Coker}(\eta))$ (cf. \S\,\ref{y0199}).
On the other hand, since  $\mr{Quot}_\BB^{n, 0}$ is   \'{e}tale over $\mr{Spec}(k)$,
we have  $\mr{Hom}_{\mcO_{X^{(1)}}}(\mcV, \mr{Coker}(\eta))$ $=0$, and hence,  $\mr{Ext}^1_{\mcO_{X^{(1)}}}(\mcV, \mr{Coker}(\eta))=0$ by Lemma \ref{y0233}  below.
This implies that
 $\eta$ may be  lifted to a $W$-rational point of the Quot-scheme $\mr{Quot}_{\mcG_W/X^{(1)}_W/W}^{n, 0}$,  so
  $\mr{Quot}_{\mcG_W/X^{(1)}_W/W}^{n, 0}$ is finite and \'{e}tale over $W$.
  A routine argument  shows
 that  $K$ may be supposed   to be a subfield of $\mbC$.
Denote by $X^{(1)}_{\mbC}$  the base-change of $X^{(1)}_{W}$ via the morphism $\mr{Spec}(\mbC) \migi  \mr{Spec}(W)$ induced by the composite embedding $W \migiincl K \migiincl \mbC$.
Also. denote 
by $\mcG_\mbC$ the pull-back of $\mcG_W$ via the natural morphism $X^{(1)}_{\mbC} \migi X^{(1)}_{W}$. 
Then, the following  equalities hold:
\begin{align}
 \mr{deg}(\mr{Quot}_\BB^{n, 0}/k) 
 & = \mr{deg}(\mr{Quot}_{\mcG_W/X^{(1)}_W/W}^{n, 0}/W)  =  \mr{deg}(\mr{Quot}^{n, 0}_{\mcG_\mbC/X^{(1)}_\mbC/C}/\mbC).  \end{align}
It follows that the problem is reduced to 
 calculating  the degree  $ \mr{deg}(\mr{Quot}^{n, 0}_{\mcG_\mbC/X^{(1)}_\mbC/C}/\mbC)$.

By  ~\cite{SUN}, \S\,2, Theorem 2.2, $F_{X/k*}(\BB^\TT)$ turns out to be  stable.
Since the degree of $\mcG_\mbC$ coincides with the degree of   $F_{X/k*}(\BB^\TT)$, 
$\mcG_\mbC$ is a  vector bundle of degree $\mr{deg}(\mcG_\mbC ) = (p-n)(g-1)$ (cf. the proof of Lemma \ref{y0233}). 
On the other hand, 
one verifies easily from the definition of stability and the properness of Quot schemes (cf.  ~\cite{FGA}, \S\,5.5, Theorem 5.14) that $\mcG_\mbC$ is a stable vector bundle.
Next, let us observe that $\mr{Quot}^{n, 0}_{\mcG_\mbC/X^{(1)}_\mbC/C}$ is zero-dimensional (cf. the discussion above), which, by the discussion preceding Theorem \ref{y0230}, implies 
$s_n(\mcG_\mbC) = n(p-n)(g-1)$.  
By choosing the deformation $\mcG_W$ of $F_{X/k*}(\BB^\TT)$ appropriately, 
we may assume, without loss of generality, that $\mcG_\mbC$ is sufficiently general in
$\mfU_{X^{(1)}_\mbC}^{p, (p-n)(g-1)}$
 so that Theorem \ref{y0230} is available.  
Then,  we have 
 \begin{align} 
 e_{\mr{max}}(\mcG_\mbC, n) &=  \frac{1}{p} \cdot (\mr{deg}_\mbC(\mcG_\mbC ) \cdot n - s_n(\mcG_\mbC))  
 \\ 
  &=  \frac{1}{p} \cdot ((p-n)(g-1) \cdot n - n \cdot (p-n)(g-1)) \notag \\
  & = 0 
  \notag
  \end{align} 
  (cf. the discussion preceding Theorem \ref{y0230}). 
If, moreover, $(a,b)$ denotes  the unique pair of integers  such that $\mr{deg}_\mbC(\mcG_\mbC ) = p \cdot a -b$ with $0 \leq b < p$,
 then it follows from the hypothesis $p>n \cdot (g-1)$ that $a = g-1$ and $b= n \cdot (g-1)$.
Thus,
by   Theorem \ref{y0230} in the case where the data 
\begin{align}
 ``(C, \mcV, m, d, n, a, b, e_{\mr{max}}(\mcV, n))"
 \end{align}
 is taken to be
  \begin{align}
  (X^{(1)}_{\mbC}, \mcG_\mbC, p, (g-1)(p-n), n, g-1, n \cdot (g-1), 0),
  \end{align}
  the following equalities hold:
\begin{align}  &  \ \  \ \  \mr{deg}(\mr{Quot}^{n, 0}_{\mcG_\mbC/X^{(1)}_\mbC/\mbC}/\mbC)  \\
 &= \frac{(-1)^{(n-1)(n (g-1) \cdot n-(g-1)n^2)/p} \cdot p^{n(g-1)}}{n!} \cdot \sum_{\zeta_1, \cdots, \zeta_n}\frac{(\prod_{i=1}^n \zeta_i)^{n \cdot (g-1)-g+1}}{\prod_{i \neq j} (\zeta_i - \zeta_j)^{g-1}} \notag \\
&=  \frac{p^{n  (g-1)}}{n!} \cdot \sum_{\zeta_1, \cdots, \zeta_n}\frac{(\prod_{i=1}^n \zeta_i)^{(n-1) (g-1)}}{\prod_{i \neq j} (\zeta_i - \zeta_j)^{g-1}}.  \notag  
 \end{align} 
This completes the proof of the required equality.
\end{proof}
The following lemma was used in the proof of Theorem \ref{y0231}.

\ble \label{y0233}
Let $k$, $X$, and $\BB$ be as in  Theorem  \ref{y0231}, and
  $\eta :\mcV \migiincl F_{X/k*}(\BB^\TT)$ an injection classified by  a $k$-rational point of $\mr{Quot}_\BB^{n, 0}$. 
 Then, 
 $\mr{Coker}(\eta)$ is a vector bundle on $X^{(1)}$  and
   the following equality holds: 
 \begin{equation} \label{QH200}
  \mr{dim}_k(\mr{Hom}_{\mcO_{X^{(1)}}}(\mcV, \mr{Coker}(\eta)) )= \mr{dim}_k(\mr{Ext}^1_{\mcO_{X^{(1)}}}(\mcV, \mr{Coker}(\eta))). 
  \end{equation}
 \ele
\begin{proof}
Recall (cf. Proposition \ref{QG180}, (i)) that 
the composite $F^*_{X/k}(\mcV) \migi  \mcA_\BB/\mcA_\BB^n$ of the pull-back
$F_{X/k}^*(\eta) : F^*_{X/k}(\mcV) \migi \mcA_\BB$ of $\eta$  and the quotient $\mcA_\BB \migisurj \mcA_\BB/\mcA_\BB^n$ is an isomorphism.
This implies that 
the natural composite $\mcA_{\BB}^n \migiincl  \mcA_{\BB} \migisurj F^*_{X/k}(\mr{Coker}(\eta))$ is an isomorphism, and hence that $F^*_{X/k}(\mr{Coker}(\eta))$ is a vector bundle.
Because of the faithful flatness of $F_{X/k}$, $\mr{Coker}(\eta)$ turns out to be a vector bundle.

Next, we shall prove (\ref{QH200}).
Since the morphism $F_{X/k}:X \migi X^{(1)}$ is  affine,   we have an  equality of Euler characteristics $\chi(F_{X/k*}(\BB^\TT)) = \chi (\BB^\TT)$.
Thus, by 
the Riemann-Roch theorem, the following equalities hold:
\begin{align}
\mr{deg}(F_{X/k*}(\BB^\TT)) &= \chi (F_{X/k*}(\BB^\TT)) - \mr{rk}(F_{X/k*}(\BB^\TT))\cdot (1-g) \\
& =  \chi (\BB^\TT) - p (1-g) \notag \\
&=  (p-n)(g-1). \notag \end{align} 
Also, since $ \mr{rk} (\mcH om_{\mcO_{X^{(1)}}}(\mcV, \mr{Coker}(\eta))) = n (p-n)$, the following equalities hold:
\begin{align}
& \hspace{4mm}  
\mr{deg}( \mcH om_{\mcO_{X^{(1)}}}(\mcV, \mr{Coker}(\eta)))  \\
&=  n \cdot \mr{deg}(\mr{Coker}(\eta)) - (p-n) \cdot \mr{deg}(\mcV) \notag \\
  & =  n \cdot \mr{deg}(F_{X/k*}(\BB^\TT)) - 0 \notag \\
  & =  n (p-n) (g-1).  \notag
  \end{align} 
Finally,
by applying the Riemann-Roch theorem again, we have
\begin{align} & \ \mr{dim}_k(\mr{Hom}_{\mcO_{X^{(1)}}}(\mcV, \mr{Coker}(\eta))) - \mr{dim}_k(\mr{Ext}_{\mcO_{X^{(1)}}}^1(\mcV, \mr{Coker}(\eta))) \\
 = &  \ \mr{deg}(\mcH om_{\mcO_{X^{(1)}}} (\mcV, \mr{Coker}(\eta))) + \mr{rk} (\mcH om_{\mcO_{X^{(1)}}} (\mcV, \mr{Coker}(\eta))) \cdot (1-g) \notag \\
 =  & \ n(p-n)(g-1) + n(p-n)  (1-g) =   0. \notag \end{align} 
 This completes the proof of the lemma.
\end{proof}

\section{Joshi's conjecture} \label{y0234}

Finally, we shall conclude the manuscript with the following  theorem (and  two remarks), which gives an affirmative answer to Joshi's conjecture for  sufficiently general curves.

\bt[cf. Theorem \ref{y020}]  \label{y0251}
Let $n$  and $g$ be  integers with $n,g>1$,  $p$ a prime with $p >n \cdot \mr{max} \{g-1, 2 \}$, $k$ a perfect field of characteristic $p$, and $\hslash$ an element of $k^\times$.
  Then, the generic degree $\mr{deg}(\mfO \mfp_{\mfs \mfl_n, \hslash, g,0}^\ZZZ/\overline{\mfM}_{g, 0}) $
 of $\mfO \mfp_{\mfs \mfl_n, \hslash, g,0}^\ZZZ$ over $\overline{\mfM}_{g,0}$
 is given by the following formula:
 \vspace{3mm}
 \begin{equation} \label{QH234}
  \mr{deg}(\mfO \mfp_{\mfs \mfl_n, \hslash, g,0}^\ZZZ/\overline{\mfM}_{g, 0})    =    \frac{ p^{(n-1)(g-1)-1}}{n!} \cdot  
 \sum_{\genfrac{.}{.}{0pt}{}{(\zeta_1, \cdots, \zeta_n) \in \mbC^{\times n} }{ \zeta_i^p=1, \ \zeta_i \neq \zeta_j (i\neq j)}}
\frac{(\prod_{i=1}^n\zeta_i)^{(n-1)(g-1)}}{\prod_{i\neq j}(\zeta_i -\zeta_j)^{g-1}}.  
 \end{equation}
 \et
\begin{proof}
One may assume that $k$ is algebraically closed.
Let us fix  a sufficiently general smooth curve $X/k$ and a line bundle $\BB$ 
as in the statement of Theorem \ref{y0231}.
Then,  it follows from (\ref{QS1100}) (with $\msX$ replaced by the universal family of curves $\msC_{g, 0}$),   Corollary \ref{QG93}, and  Theorem \ref{y0231}  that
 \begin{align} 
 \mr{deg}(\mfO \mfp_{\mfs \mfl_n, \hslash, g,0}^\ZZZ/\overline{\mfM}_{g, 0})   
 & = \mr{deg}(\mfO \mfp_{\mfs \mfl_n, g,0}^\ZZZ/\overline{\mfM}_{g, 0})  \\
 &    =  \frac{1}{p^g} \cdot  \mr{deg}(\mr{Quot}^{n,0}_\BB/k) \notag  \\
 & = \frac{ p^{(n-1)(g-1)-1}}{n!} \cdot  
 \sum_{\genfrac{.}{.}{0pt}{}{(\zeta_1, \cdots, \zeta_n) \in \mbC^{\times n} }{ \zeta_i^p=1, \ \zeta_i \neq \zeta_j (i\neq j)}}
\frac{(\prod_{i=1}^n\zeta_i)^{(n-1)(g-1)}}{\prod_{i\neq j}(\zeta_i -\zeta_j)^{g-1}}.  \notag
 \end{align}
This completes the proof of the assertion.
\end{proof}

\begin{rema}
Let $n$, $g$, and $p$ be as in Theorem \ref{y0251}.
Also, let $k$ be an algebraically closed field of characteristic $p$,  $X$  a {\it sufficiently general} smooth curve over $k$ of genus $g$, and $\BB$ a line bundle on $X$ admitting  an isomorphism $\tau_\BB : \mcO_X \isom \mcT_{X/k}^{\otimes \frac{n(n-1)}{2}} \otimes \BB^{\otimes n}$
 (e.g., the $(n-1)$-st tensor product  $\BB_0^{\otimes (n-1)}$ of a theta characteristic $\BB_0$ of $X$).
In particular, we obtain an $(n, \hslash)$-theta characteristic  $\bb := (\BB, \tau_{\BB*} (d))$ of $X/k$.
Denote by $\mr{Diff}_{n, \BB}^{\clubsuit \mr{full}}$ ($=\mfD \mfi \mff \mff_{n, \bb}^{\, \clubsuit \mr{full}}(k)$ in the notation of Proposition \ref{QG1780})
the set of $(n, \bb)$-projective connections on $X/k$ having a full set of root functions (cf. Definition \ref{D01fgh}).
That is to say,  
this set  classifies 
$n$-th order differential operators $\BB^\vee \migi \Omega_{X/k}^{\otimes n} \otimes \BB^\vee$ (in the sense of Remark \ref{QQ425}) with unit principal symbol and trivial subprincipal symbol, and moreover, having a full set of root functions. 
Then, by Corollary \ref{1090}, Proposition \ref{QG1780},  and  Theorem \ref{y0251} together, 
the cardinality  of the set  $\mr{Diff}_{n, \BB}^{\clubsuit \mr{full}}$ is given  by the right-hand side of   (\ref{QH234}).
\end{rema}

\begin{rema} \label{ghj3}
 In ~\cite{Wak}, we performed   explicit computations of the values  $\mr{deg}(\mfO \mfp_{\mfs \mfl_2, g,0}^\ZZZ/\overline{\mfM}_{g, 0})$ for small genus $g$ (cf. ~\cite{Wak}, \S\,6.2 (2)).
 Recall from  ~\cite{LO}, \S\,2, Theorem 2.1, that  these values may be expressed as a polynomial   of degree $3g-3$ with respect to  $p$ (cf.  (\ref{n=2}) in Introduction).
One may also give a polynomial  description of 
 the   value 
     $d_{3, g} := \mr{deg}(\mfO \mfp_{\mfs \mfl_3, g,0}^\ZZZ/\overline{\mfM}_{g, 0})$.
Indeed, Theorem \ref{y0251} implies 
\begin{align} \label{fff093}
& \  d_{3,g} \\
    =  & \    \frac{ p^{2g-3} \cdot (-1)^{g-1}}{6}  \cdot 
 \sum_{\genfrac{.}{.}{0pt}{}{(\zeta_1, \zeta_2, \zeta_3) \in \mbC^{\times n} }{ \zeta_i^p=1, \ \zeta_i \neq \zeta_j (i\neq j)}}
\frac{(\zeta_1 \zeta_2 \zeta_3)^{2(g-1)}}{((\zeta_1 -\zeta_2) (\zeta_2 - \zeta_3) (\zeta_3 -\zeta_1))^{2 (g-1)}} \notag  \\
= & \   \frac{ p^{2g-3} \cdot (-1)^{g-1}}{6} 
 \cdot   \hspace{-3mm}
 \sum_{\genfrac{.}{.}{0pt}{}{0 \leq \eta_1, \eta_2, \eta_3 < p-1}{\eta_i \neq \eta_j (i \neq j)}} 
  \displaystyle{\prod_{\scalebox{0.6}{$(i,j)=(2,1), (3,2), (1,3)$}}}\left(1 - e^{\frac{2 \pi \sqrt{-1} (\eta_i - \eta_j)}{p}}\right)^{-2(g-1)}
  \notag
  \\
  = &  \  \frac{ p^{2g-3} \cdot (-1)^{g-1} \cdot (2p)}{6}  \cdot \hspace{-5mm}
 \sum_{\genfrac{.}{.}{0pt}{}{0 < \theta_1, \theta_2, p -\theta_1 -\theta_2}{}} 
 \displaystyle{\prod_{\scalebox{0.6}{$\theta = \theta_1, \theta_2, -\theta_1 - \theta_2$}}}\left(1 - e^{\frac{2 \pi \sqrt{-1} \theta}{p}}\right)^{-2(g-1)}
  \notag
   \\
=  & \    \frac{ p^{2g-2} \cdot (-1)^{g-1}}{3}  \cdot \hspace{0mm} \sum_{\genfrac{.}{.}{0pt}{}{0 < \theta_1, \theta_2,}{\theta_1 +\theta_2<p}} 
\hspace{0mm}
\displaystyle{\prod_{\scalebox{0.6}{$\theta = \theta_1, \theta_2, -\theta_1-\theta_2$}}}\left(2  \sin\left( \frac{\pi \theta}{p} \right) e^{\pi \sqrt{-1} \left(\frac{\theta}{p}-\frac{1}{2}\right)}\right)^{-2(g-1)}
 \hspace{-5mm}   \notag 
 \end{align}
    \begin{align}
= & \ \frac{ p^{2g-2}}{3}  \cdot  \hspace{0mm}\sum_{\genfrac{.}{.}{0pt}{}{0 < \theta_1, \theta_2, }{\theta_1 +\theta_2<p}} 
\left(8 \sin\left( \frac{\pi \theta_1}{p} \right) \sin\left( \frac{\pi \theta_2}{p} \right) \sin\left( \frac{\pi (\theta_1 + \theta_2)}{p} \right)\right)^{-2(g-1)}. \notag
 \end{align}
On the other hand, it follows from the discussion in  ~\cite{Sz},  \S\,3, that
\begin{align} \label{fff091}
& \  \sum_{\genfrac{.}{.}{0pt}{}{0 < \theta_1, \theta_2,}{\theta_1 +\theta_2<p}} 
\left(8 \sin\left( \frac{\pi \theta_1}{p} \right) \sin\left( \frac{\pi \theta_2}{p} \right) \sin\left( \frac{\pi (\theta_1 + \theta_2)}{p} \right)\right)^{-2(g-1)} \\
= & \  \mr{Res}_{y =0} \mr{Res}_{x =0} \left( \frac{k^2 \mr{cot} (kx) \mr{cot} (ky) dx dy}{(8 \sin (x) \sin (y) \sin (x+ y))^{2 (g-1)}} \right), \notag 
\end{align}
where  $\mr{Res}_{x=0}(-)$  and  $\mr{Res}_{y=0}(-)$ denote the ordinary $1$-dimensional residues at $x =0$ and $y =0$ respectively taken by assuming  all the other variables to be constants.
 Thus, the right-hand side of  (\ref{fff091})
 may be computed by considering 
 the relation
 $\frac{1}{\mr{sin}^2(t)} = 1 + \mr{cot}^2(t)$
 and  the coefficients of the Laurent expansion
  $\mr{cot}(t) = \frac{1}{t}+\sum_{j=1}^{\infty}\frac{(-1)^{j} \cdot 2^{2j} \cdot B_{2j}}{(2j)!} \cdot t^{2j-1}$,
where $B_{2j}$ denotes the $(2j)$-th Bernoulli number, i.e.,
$\frac{w}{e^w-1} = 1 - \frac{w}{2} + \sum_{j=1}^\infty \frac{B_{2j}}{(2j)!} \cdot w^{2j}$. 
Consequently,  by  (\ref{fff093}),  (\ref{fff091}) and   an explicit computation,  the value $d_{3, g}$ may be expressed as a polynomial with respect to $p$.
Some examples for small $g$ are obtained as follows:
\begin{align}
d_{3,2}  & = \frac{p^8}{181440} + \frac{p^6}{4320} - \frac{11p^4}{8640} + \frac{47p^2}{45360}, \\
 d_{3,3} & = 
  \frac{19p^{16}}{1120863744000} + \frac{p^{14}}{1437004800} + \frac{p^{12}}{62208000}   +  \frac{139 p^{10}}{391910400}  + \frac{7583 p^{8}}{2612736000}  - \frac{269 p^6}{11404800}    \notag \\
  & \ \ \ 
   + \frac{1422991 p^4}{70053984000},
   \notag
   \end{align}
    \begin{align}
   d_{3,4} & = \frac{1031 p^{24}}{15327282651512832000}  + \frac{p^{22}}{250092722880000}  
     +  \frac{31 p^{20}}{251073478656000}+ \frac{105359 p^{18}}{39544072888320000}  
     \notag \\
         & \ \ \    + \frac{2689 p^{16}}{57940033536000} + \frac{19303 p^{14}}{24141680640000} +  \frac{141011551 p^{12}}{15817629155328000}  + \frac{115715659 p^{10}}{1883051089920000}\notag \\
   & \ \ \     - \frac{518113  p^8}{889218570240}  + \frac{1224856468291 p^6}{2394887914298880000}, \notag 
   \end{align}
   \begin{align}
  d_{3,5} & = \frac{32293 p^{32}}{118449240330891165696000000} + \frac{911 p^{30}}{42303300118175416320000} + \frac{1763 p^{28}}{2021180129869824000000}\notag \\
    & \ \ \   + \frac{1400983 p^{26}}{57937128422718504960000} + \frac{2002933  p^{24}}{3872155617224294400000} + \frac{677 p^{22}}{73952551895040000} \notag \\
    & \ \ \   + \frac{339010759 p^{20}}{2391625528285593600000} + \frac{3036911 p^{18}}{1423586623979520000} + \frac{2254841161 p^{16}}{85102321257676800000} \notag \\
   & \ \ \  
    + \frac{4688714632351  p^{14}}{19312376140906168320000}  + \frac{4013074350869 p^{12}}{2554547108585472000000}  - \frac{539250373 p^{10}}{33307587016704000}
    \notag \\
    & \ \ \ 
    + \frac{59746303954645337 p^8}{4164231105382892544000000}, \notag 
\end{align}
\begin{align}
d_{3,6} & = \frac{27739097 p^{40}}{25048031592726595660148637696000000} + \frac{669077 p^{38}}{6121219841819793660837888000000}
\notag \\
& \ \ \ + \frac{17351 p^{36}}{3153337759397169856512000000} + \frac{63632321 p^{34}}{337722474031436891632435200000}
\notag \\
& \ \ \ + \frac{278001151 p^{32}}{56287079005239481938739200000}+ \frac{84750829 p^{30}}{801187718759878754304000000}
\notag \\
& \ \ \ + \frac{85048301147 p^{28}}{44050757482361333691187200000} + \frac{3990465223 p^{26}}{128803384451348928921600000}
\notag \\
& \ \ \ 
+ \frac{1380198109 p^{24}}{3066747248841641164800000}+ \frac{8628232370623 p^{22}}{1376586171323791677849600000}
\notag \\
& \ \ \ 
+ \frac{7568113936878661 p^{20}}{95011437707053856980992000000}+ \frac{22634647835562647 p^{18}}{26803370954875943780352000000}
\notag \\
& \ \ \ 
+ \frac{55113411366175237 p^{16}}{7817649861838816935936000000}+ \frac{11845004996053 p^{14}}{268094988403251609600000}
\notag \\
& \ \ \ 
- \frac{35245261 p^{12}}{73457621539061760} + \frac{11954875855784102515639 p^{10}}{27955392402596646942130176000000}.
\notag
\end{align}

\end{rema}







\printindex
\end{document}